\documentclass{amsart}

\usepackage{amsmath,amsthm,amsfonts}
\usepackage{latexsym}
\usepackage{amssymb}

\usepackage[dvips]{epsfig}

\def\Re{{\rm Re}}
\def\Im{{\rm Im}}
\def\bfx{\mathbf{x}}
\def\bfy{\mathbf{y}}
\def\Leins{L^1_t L^1_x}
\def\Lzweins{L^2_t L^1_x}
\def\calC{{\mathcal C}}
\def\ener{L^\infty_t L^2_x}
\def\enerN{L^1_t L^2_x}
\def\Linf{L^\infty_t L^\infty_x}
\def\term{{\mathit{term}}}
\def\wt{\widetilde}
\def\wh{\widehat}
\def\wht{\wh}
\def\nn{\nonumber}
\def\lhalb{\frac{\ell}{2}}
\def\Ltwotx{L^2_t L^2_x}
\def\trip{|\!|\!|}

\def\calN{{\mathcal Q}}
\def\vecv{\vec{v}}
\def\vecw{\vec{w}}

\def\disk{{\mathbb D}}
\def\eps{\varepsilon}
\def\calM{{\mathcal M}}
\def\calE{{\mathcal E}}
\def\calT{{\mathcal T}}
\def\calL{{\mathcal L}}
\def\del{\partial}
\def\R{{\mathbb R}}
\def\les{\lesssim}
\def\dist{{\rm dist}}

\def\diam{{\rm diam}}
\def\la{\langle}
\def\ra{\rangle}
\renewcommand\angle{\sphericalangle}
\def\NFA{{\rm NFA}}
\def\PWA{{\rm PWA}}
\def\calQ{{\mathcal Q}}
\def\NF{{\rm NF}}
\def\PW{{\rm PW}}
\def\C{{\mathbb C}}
\def\Z{{\mathbb Z}}
\def\supp{{\rm supp}}
\def\calS{{\mathcal S}}
\def\sign{{\rm sign}}
\def\calC{{\mathcal C}}
\def\calF{{\mathcal F}}
\def\enr{{L^\infty_t L^2_x}}
\def\enromega{{L^\infty_{t_\omega} L^2_{x_\omega}}}

\def\caps{{\mathcal C}}
\def\Ltwo{{L^2_t L^2_x}}
\def\Felltwo{{L^2_{\xi_\omega} L^2_{\tau_\omega}}}
\def\sangle{\sphericalangle}
\def\calR{{\mathcal R}}
\def\calD{{\mathcal D}}

\def\NP{{\rm NP}}
\def\med{{\rm med}}
\def\be{{\mathbf e}}
\def\Hyp{{\mathbb H}}

\def\bmo{\mathrm{BMO}}
\def\calH{{\mathcal H}}
\def\const{{\mathrm const}}
\def\vphi{\varphi}
\def\calK{{\mathcal K}}
\def\calR{{\mathcal R}}
\def\Hyp{\mathbb{H}}
\def\slashint{-\!\!\!\!\!\!\int}
\def\uf{\underbar{f}}
\def\ueps{\underline{\varepsilon}}
\def\bfu{\mathbf{u}}
\def\bfv{\mathbf{v}}
\def\Ecrit{E_{\mathrm crit}\,}
\def\weg{\!\!\!\!}

\newtheorem{theorem}{Theorem}
\newtheorem{lemma}[theorem]{Lemma}
\newtheorem{prop}[theorem]{Proposition}
\newtheorem{cor}[theorem]{Corollary}
\newtheorem{exse}[theorem]{Remark}
\newtheorem{defi}[theorem]{Definition}
\theoremstyle{remark}
\newtheorem{remark}[theorem]{Remark}

\numberwithin{equation}{section} \numberwithin{theorem}{section}

\begin{document}

\title{Concentration compactness for critical wave maps}
\author{Joachim Krieger, Wilhelm Schlag}
\address{JK: Department of Mathematics, The University of Pennsylvania, 209 South 33rd Street, Philadelphia, PA 19104, U.S.A.}
\address{WS: Department of Mathematics, The University of Chicago, 5734 South University Avenue, Chicago, IL 60615, U.S.A.}
\thanks{The authors wish to thank Sergiu Klainerman, Carlos Kenig,
and Frank Merle for their interest in this work, as well as their
encouragement. The first author is also indebted to Sergiu
Klainerman for introducing him to the topic of wave maps in general,
and to those into the hyperbolic plane in particular.  Further, they
thank Carlos Kenig, Daniel Tataru, and Jacob Sterbenz for helpful
discussions. The authors were supported in part by the National
Science Foundation, DMS-0757278 (JK) and DMS-0617854 (WS). The first
author was supported by a Sloan fellowship.}

\begin{abstract}
   By means of the concentrated compactness method of
   Bahouri-Gerard~\cite{BG} and Kenig-Merle~\cite{KeM1}, we prove global existence and
   regularity for wave maps with smooth data and large energy from
   $\R^{2+1}\to\Hyp^2$. The argument yields an apriori bound of the Coulomb
   gauged derivative components of our wave map relative to a suitable norm~$\|\cdot\|_S$
   (which holds the solution) in terms of the energy alone. As a
   by-product of our argument, we obtain a phase-space decomposition
   of the gauged derivative components
   analogous to the one of Bahouri-Gerard.
\end{abstract}

\maketitle

\section{Introduction and Overview}\label{sec:intro}

\subsection{The main result and its history}

Formally speaking, wave maps are the analogue of harmonic maps where
the Minkowski metric is imposed on the independent variables. More
precisely, for a smooth $\bfu:\R^{n+1}\to\calM$ with $(\calM,g)$
Riemannian,  define the Lagrangian
\[
\calL(\bfu):= \int_{\R^{n+1}} \big(|\del_t \bfu|_g^2 -
|\nabla\bfu|_g^2\big)\, dtdx
\]
Then the critical points are defined as $\calL'(\bfu)=0$ which means
that $\Box\bfu\perp T_u\calM$ in case $\calM$ is imbedded in some
Euclidean space. This is called the {\em extrinsic} formulation,
which can also be written as
\[
\Box\bfu + A(\bfu)(\del_\alpha \bfu, \del^\alpha\bfu)=0
\]
where $A(\bfu)$ is the second fundamental form.  In view of this, it
is clear that $\gamma\circ\phi$ is a wave map for any geodesic
$\gamma$ in~$\calM$ and any free scalar wave~$\phi$. Moreover, any
harmonic map is a stationary wave map. The {\em intrinsic}
formulation is $D^\alpha \del_\alpha u=0$, where
\[
D_\alpha X^j:= \del_\alpha X^j + \Gamma^j_{ik} \circ \bfu \, X^i
\del_\alpha u^k
\]
is the covariant derivative induced by~$\bfu$ on the pull-back
bundle of $T\calM$ under~$\bfu$ (with the summation convention in
force). Thus, in local coordinates $\bfu=(u^1,\ldots,u^d)$ one has
\begin{equation}\label{eq:intrins}
\Box u^j + \Gamma^j_{ik}\circ\bfu \, \del_\alpha u^i \, \del^\alpha
u^k =0
\end{equation}
The central problem for wave maps is to answer the following
question:

\medskip
\noindent  {\em For which $\calM$ does the Cauchy problem for the
wave map $\bfu:\R^{n+1}\to \calM$ with smooth data
$(\bfu,\dot\bfu)|_{t=0}=(\bfu_0,\bfu_1)$ have global smooth
solutions?}

\medskip
In view of finite propagation speed, one may assume that the data
$(\bfu_0,\bfu_1)$ are trivial outside of some compact set (i.e.,
$\bfu_0$ is constant outside of some compact set, whereas $\bfu_1$
vanishes outside of that set). Let us briefly describe what is known
about this problem.

First, recall that the wave map equation is invariant under the
scaling $\bfu\mapsto \bfu(\lambda\cdot)$ which is critical relative
to $\dot H^{\frac{n}{2}}(\R^n)$, whereas the conserved energy
\[
\calE(\bfu) =  \frac12 \sum_{\alpha=0}^n \int_{\R^n} |\del_\alpha
\bfu(t,\cdot)|^2\, dx
\]
is critical relative to $\dot H^1(\R^n)$. In the supercritical case
$n\ge3$ it was observed by Shatah~\cite{Shatah} that there are
self-similar blowup solutions of finite energy. In the critical case
$n=2$, it is known that there can be no self-similar blowup,
see~\cite{SStruwe}. Moreover, Struwe~\cite{Struwe1} observed that in
the equivariant setting, blowup in this dimension has to result from
a strictly slower than self-similar rescaling of a harmonic sphere
of finite energy. His arguments were based on the very detailed
wellposedness of equivariant wave maps by Christodoulou,
Tavildar-Zadeh~\cite{CT1}, \cite{CT2}, and Shatah,
Tahvildar-Zadeh~\cite{ShaT1}, \cite{ShaT2} in the energy class for
equivariant wave maps into manifolds that are invariant under the
action of~$SO(2,\R)$. Finally,  Rodnianski,
Sterbenz~\cite{RodSterb}, as well as the authors together with
Daniel Tataru~\cite{KST} exhibited finite energy wave maps from
$\R^2\to S^2$ that blow up in finite time by suitable rescaling of
harmonic maps.

Let us now briefly recall some well-posedness results. The
nonlinearity in~\eqref{eq:intrins} displays a {\em nullform
structure}, which was the essential feature in the subcritical
theory of Klainerman-Machedon~\cite{KM1}--\cite{KM3}, and
Klainerman-Selberg~\cite{KS1}, \cite{KS2}. These authors proved
strong  local  well-posedness  for data in~$H^s(\R^n)$ when
$s>\frac{n}{2}$. The important critical theory $s=\frac{n}{2}$ was
begun by Tataru~\cite{Tat3}, \cite{Tat2}. These seminal papers
proved global well-posedness for smooth data satisfying a smallness
condition in $\dot B^{\frac{n}{2}}_{2,1}(\R^n)\times \dot
B^{\frac{n}{2}-1}_{2,1}(\R^n)$. In a breakthrough work,
Tao~\cite{T1}, \cite{T2} was able to prove well-posedness for data
with small $\dot{H}^{\frac{n}{2}}\times \dot H^{\frac{n}{2}}$ norm
and the sphere as target. For this purpose, he introduced the
important {\em microlocal gauge} in order to remove some ``bad''
interaction terms from the nonlinearity. Later results by
Klainerman, Rodnianski~\cite{KlRod}, Nahmod, Stephanov,
Uhlenbeck~\cite{NSU}, Tataru~\cite{Tat1}, \cite{Tat}, and
Krieger~\cite{Krieger}, \cite{Krieger2}, \cite{Krieger3} considered
other cases of targets by using similar methods as in Tao's work.

Recently, Sterbenz and Tataru~\cite{SterbTat1}, \cite{SterbTat2}
have given the following very satisfactory answer to the above
question: {\em If the energy of the initial data is smaller than the
energy of any nontrivial harmonic map $\R^n\to\calM$, then one has
global existence and regularity.}

\noindent Notice in particular that if there are no harmonic maps
other than constants, then one has global existence for all
energies. A particular case of this are the hyperbolic
spaces~$\Hyp^n$ for which Tao~\cite{T3}--\cite{T7} has achieved the
same result (with some apriori global norm control).

The purpose of this paper is to apply the method of compensated
compactness as in Bahouri, Gerard~\cite{BG} and Kenig,
Merle~\cite{KeM1}, \cite{KeM2} to the large data wave map problem
with the hyperbolic plane~$\Hyp^2$ as target. We emphasize that this
gives more than global existence and regularity as already in the
semilinear case considered by the aforementioned authors. The fact
that in the critical case the large data problem should be decided
by the geometry of the target is a conjecture going back to Sergiu
Klainerman.

Let us now describe our result in more detail. Let $\Hyp^2$ be the
upper half-plane model of the hyperbolic plane equipped with the
metric $ds^2=\frac{d\bfx^2+d\bfy^2}{\bfy^2}$. Let
$\bfu:\R^2\to\Hyp^2$ be a smooth map. Expanding the derivatives
$\{\del_\alpha \bfu\}_{\alpha=0,1,2}$ (with $\del_0:=\del_t$) in the
orthonormal frame $\{\be_1,\be_2\}=\{\bfy\del_\bfx,\bfy\del_\bfy\}$
gives rise to smooth coordinate functions
$\phi^1_\alpha,\phi^2_\alpha$. In what follows, $\|\del_\alpha
\bfu\|_X$ will mean $\sum_{j=1}^2 \|\phi^j_\alpha\|_X$ for any norm
$\|\cdot\|_X$ on scalar functions. For example, the energy of $\bfu$
is
\[
 E(\bfu):= \sum_{\alpha=0}^2 \|\del_\alpha \bfu\|_2^2
\]
Similarly, suppose $\pi:\Hyp^2\to M$ is a covering map with $M$ some
hyperbolic Riemann surface with the metric that renders $\pi$ a
local isometry. In other words, $M=\Hyp^2/\Gamma$ for some discrete
subgroup $\Gamma\subset PSL(2,\R)$ which operates totally
discontinuously on~$\Hyp^2$. Now suppose $\bfu:\R^2\to M$ is a
smooth map which is constant outside of some compact set, say. It
lifts to a smooth map $\tilde\bfu:\R^2\to \Hyp^2$ uniquely, up to
composition with an element of~$\Gamma$. We now define
$\|\del_\alpha \bfu\|_X:= \|\del_\alpha \tilde\bfu\|_X$. In
particular, the energy $E(\bfu):= E(\tilde\bfu)$. Note that due to
the fact that $\Gamma$ is a group of isometries of~$\Hyp^2$, these
definitions are unambiguous. Our main result is as follows.

\begin{theorem}
\label{thm:main} The exists a function $K:(0,\infty)\to (0,\infty)$
with the following property: Let $M$ be  a hyperbolic Riemann
surface. Suppose $(\bfu_0,\bfu_1):\R^2\to M\times TM$ are smooth and
$\bfu_0=\const$, $\bfu_1=0$ outside of some compact set. Then the
wave map evolution $\bfu$ of these data as a map $\R^{1+2}\to M$
exists globally as a smooth function and, moreover, for any
$\frac{1}{p}+\frac{1}{2q}\le \frac14$ with $2\le q<\infty$,
$\gamma=1-\frac{1}{p}-\frac{2}{q}$,
\begin{equation}
 \label{eq:Sbound} \sum_{\alpha=0}^2 \|(-\Delta)^{-\frac{\gamma}{2}} \del_\alpha \bfu \|_{L^p_t L^q_x} \le C_q\,K(E)
\end{equation}
Moreover, in the case when $M\hookrightarrow \R^N$ is a compact Riemann surface, one has scattering:
\[
\max_{\alpha=0,1,2} \| \del_\alpha \bfu (t) - \del_\alpha S(t) (f,g)\|_{L^2_x} \to 0\quad \text{\ as \ }t\to\pm\infty
\]
where $S(t)(f,g)=\cos(t|\nabla|) f + \frac{\sin(t|\nabla|)}{|\nabla|}g$ and suitable $(f,g)\in  (\dot H^1\times L^2)(\R^2;\R^N)$.
Alternatively, if $M$ is non-compact, then lifting $\bfu$ to a map $\R^{1+2}\to \Hyp^2$ with derivative components $\phi_\alpha^j$
as defined above, one has
\[
 \max_{\alpha=0,1,2}\| \phi_\alpha^j(t) -  \del_\alpha S(t) (f^j,g^j)\|_{L^2_x} \to 0\quad \text{\ as \ }t\to\pm\infty
\]
where $(f^j,g^j)\in  (\dot H^1\times L^2)(\R^2;\R)$.
\end{theorem}

We emphasize that \eqref{eq:Sbound} can be strengthened considerably
in terms of the type of norm applied to the Coulomb gauged
derivative components of the wave map:
\begin{equation}
  \label{eq:Sbound2}
\sum_{\alpha=0}^2 \|\psi_\alpha  \|^2_{S} \le C\,K(E)^2
\end{equation}
The meaning $\psi_\alpha$ as well as of the $S$ norm will be
explained below.
 We now turn to describing this result and our methods in
more detail. For more background on wave maps see \cite{GrotShatah},
\cite{Tat1},  and~\cite{SStruwe}.

\subsection{Wave maps to $\Hyp^2$}

The manifold $\Hyp^2$ is the upper half-plane equipped with the
metric $ds^2=\frac{d\bfx^2+d\bfy^2}{\bfy^2}$. Expanding the
derivatives $\{\del_\alpha \bfu\}_{\alpha=0,1,2}$ (with
$\del_0:=\del_t$) of a smooth map $\bfu:\R^{1+2}\to \Hyp^2$ in the
orthonormal frame $\{\be_1,\be_2\}=\{\bfy\del_\bfx,\bfy\del_\bfy\}$
yields
\[
 \del_\alpha \bfu = (\del_\alpha\bfx,\del_\alpha \bfy) = \sum_{j=1}^2 \phi^{j}_{\alpha}\,
 \be_j
\]
whence
\begin{equation}
  \label{eq:xy_fund} \bfy = e^{\sum_{j=1,2}
  \Delta^{-1}\del_j\phi_j^2} , \quad \bfx = \sum_{j=1,2} \Delta^{-1}
  \del_j (\phi^1_j y)
\end{equation}
 Energy conservation takes the form
\begin{equation}
 \label{eq:phi_ener} \int_{\R^2} \sum_{\alpha=0}^2 \sum_{j=1}^2 |\phi^j_\alpha(t,x)|^2\, dx  = \int_{\R^2} \sum_{\alpha=0}^2
 \sum_{j=1}^2 |\phi^j_\alpha(0,x)|^2\, dx
\end{equation}
where $x=(x_1,x_2)$ and $\del_0=\del_t$.  If $\bfu(t,x)$ is a smooth
wave map, then the functions $\{\phi^j_\alpha\}$ for $0\le
\alpha\le2$ and $j=1,2$ satisfy the div--curl system
\begin{align}
\del_\beta \phi^1_\alpha -\del_\alpha \phi^1_\beta &= \phi^1_\alpha \phi^2_\beta - \phi^1_\beta \phi^2_\alpha \label{eq:compat1} \\
\del_\beta \phi^2_\alpha -\del_\alpha\phi^2_\beta &=0 \label{eq:compat2} \\
\del_\alpha \phi^{1\alpha} &= - \phi^1_\alpha \phi^{2\alpha} \label{eq:derwmp1} \\
\del_\alpha \phi^{2\alpha} &=  \phi^1_\alpha \phi^{1\alpha} \label{eq:derwmp2}
\end{align}
for all $\alpha,\beta=0,1,2$. As usual, repeated indices are being
summed over, and lowering or raising is done via the Minkowski
metric. Clearly, \eqref{eq:compat1} and~\eqref{eq:compat2} are
integrability conditions which are an expression of the curvature
of~$\Hyp^2$. On the other hand, \eqref{eq:derwmp1}
and~\eqref{eq:derwmp2} are the actual wave map system. Since the
choice of frame was arbitrary, one still has gauge freedom for the
system~\eqref{eq:compat1}--\eqref{eq:derwmp2}. We shall exclusively
rely on the Coulomb gauge which is given in terms of complex
notation by the functions
\begin{equation}
 \label{eq:psi_def} \psi_\alpha := \psi^1_\alpha + i\psi^2_\alpha = (\phi^1_\alpha + i\phi^2_\alpha) e^{-i \Delta^{-1} \sum_{j=1}^2 \del_j \phi^1_j}
\end{equation}
If $\phi^1_j$ are Schwartz functions , then $\sum_{j=1}^2 \del_j \phi^1_j$ has mean zero whence
\begin{equation} \label{eq:logpot}
 (\Delta^{-1} \sum_{j=1}^2 \del_j \phi^1_j)(z) = \frac{1}{2\pi} \int_{\R^2} \log|z-\zeta|  \sum_{j=1}^2 \del_j \phi^1_j (\zeta)\, d\zeta\wedge d\bar{\zeta}
\end{equation}
is well-defined and moreover decays like~$|z|^{-1}$ (but in general no faster).
The gauged components $\{\psi_\alpha\}_{\alpha=0,1,2}$ satisfy the new div--curl system
\begin{align}
 \del_\alpha\psi_\beta - \del_\beta\psi_\alpha &= i\psi_\beta \Delta^{-1} \sum_{j=1,2}
 \del_j(\psi^1_\alpha \psi^2_j - \psi^2_\alpha \psi^1_j) - i\psi_\alpha\Delta^{-1} \del_j(\psi^1_\beta \psi^2_j - \psi^2_\beta \psi^1_j) \label{eq:psisys1}\\
\del_\nu \psi^\nu &= i\psi^\nu \Delta^{-1} \sum_{j=1}^2 \del_j(\psi^1_\nu \psi^2_j - \psi^2_\nu \psi^1_j) \label{eq:psisys2}
\end{align}
In particular, one obtains the following system of wave equations for the $\psi_\alpha$:
\begin{equation}\label{eq:psi_wave}
\begin{aligned}
 \Box \psi_\alpha & =  i\del^\beta\big[ \psi_\alpha\Delta^{-1} \sum_{j=1,2}
 \del_j(\psi^1_\beta \psi^2_j - \psi^2_\beta \psi^1_j)\big] - i\del^\beta\big[\psi_\beta\Delta^{-1}
 \del_j(\psi^1_\alpha \psi^2_j - \psi^2_\alpha \psi^1_j)\big] \\ &\quad
 +i \del_\alpha\big[ \psi^\beta\Delta^{-1} \sum_{j=1,2} \del_j(\psi^1_\beta \psi^2_j - \psi^2_\beta \psi^1_j)   \big]
\end{aligned}
\end{equation}
Throughout this paper we shall only consider {\em admissible} wave maps~$\bfu$. These are
characterized as smooth wave maps $\bfu: I\times \R^2\to\Hyp^2$ on some time interval~$I$ so that the
derivative components~$\phi^j_\alpha$ are Schwartz functions on fixed time slices.

\noindent
By the method of {\it Hodge decompositions} from\footnote{In these
papers this decomposition is also referred to as ``dynamic
decomposition''.}~\cite{Krieger}--\cite{Krieger3}  one exhibits the
null-structure present in~\eqref{eq:psisys1}--\eqref{eq:psi_wave}.
Hodge decomposition here refers to writing
\begin{equation}\label{eq:dyn_dec}
 \psi_\beta = -R_\beta\sum_{k=1}^2 R_k \psi_k + \chi_\beta
\end{equation}
where $R_\beta:= \del_\beta|\nabla|^{-1}$ are the usual Riesz
transform. Inserting the  hyperbolic terms~$R_\beta\sum_{k=1}^2 R_k
\psi_k$ into the right-hand sides
of~\eqref{eq:psisys1}--\eqref{eq:psi_wave} leads to trilinear
nonlinearities with a null structure. As is well-known,  such null
structures are amenable to better estimates since they annihilate
``self-interactions'', or more precisely,  interactions of  waves
which propagate along the same characteristics,
cf.~\cite{KM1}--\cite{KM2}, as well as~\cite{KS1}, \cite{KS2},
\cite{FK}. Furthermore, inserting at least one ``elliptic term''
$\chi_\beta$ from~\eqref{eq:dyn_dec} leads to a higher order
nonlinearity, in fact quintic or higher which are easier to estimate
(essentially by means of Strichartz norms). To see this, note that
\begin{align*}
 \sum_{j=1}^2 \del_j \chi_j &=0\\
\del_j \chi_\beta - \del_\beta \chi_j &= \del_j \psi_\beta - \del_\beta \psi_j
\end{align*}
whence
\begin{equation}\label{eq:chi_beta}
 \chi_\beta = i\sum_{j,k=1}^2 \del_j \Delta^{-1} \big[\psi_\beta \Delta^{-1} \del_k(\psi_j^1 \psi_k^2 - \psi^1_k \psi^2_j)-
\psi_j\Delta^{-1} \del_k (\psi^1_\beta \psi^2_k - \psi^1_k \psi^2_\beta)\big]
\end{equation}
 Since we are only going to obtain apriori bounds on $\phi^j_\alpha$, it will suffice to assume throughout that the $\phi^j_\alpha$ are
Schwartz functions, whence the same holds for~$\psi_\alpha$.
In what follows, we shall never actually {\it solve the system}~\eqref{eq:psisys1}--\eqref{eq:psi_wave}.
To go further, the wave-equation~\eqref{eq:psi_wave} by {\it itself is meaningless}.
In fact, it is clear that \eqref{eq:psisys1} and~\eqref{eq:psisys2} will hold for all $t\in(-T,T)$
if and only if they hold at time~$t=0$ and~\eqref{eq:psi_wave} holds for all $t\in(-T,T)$.
This being said, we
will only use the system~\eqref{eq:psi_wave} to derive {\em apriori estimates} for $\psi_\alpha$, which will then be shown to lead to
suitable bounds on the components $\phi^j_\alpha$ of derivatives of a wave map~$\bfu$. This is done by means of Tao's device of frequency envelope,
see~\cite{T1} or~\cite{Krieger}.  This refers to a sequence~$\{c_k\}_{k\in\Z}$ of positive reals such that
\begin{equation}
 \label{eq:flat}  c_k\, 2^{-\sigma|k-\ell|} \le c_\ell \le c_k\, 2^{\sigma |k-\ell|}
\end{equation}
where $\sigma>0$ is a small number. The most relevant example is given by
\[
 c_k := \Big(\sum_{\ell\in\Z} 2^{-\sigma|k-\ell|} \|P_\ell\psi(0)\|_{2}^2\Big)^{\frac12}
\]
which controls the initial data.
 While it is of course clear that~\eqref{eq:compat1}--\eqref{eq:derwmp2}
imply the system~\eqref{eq:psisys1}--\eqref{eq:psi_wave}, the reverse implication is not such a simple matter since it involves
solving an elliptic system with large solutions. On the other hand, transferring estimates on the~$\psi_\alpha$ in~$H^s(\R^2)$ spaces
to similar bounds on the derivative components~$\phi^j_\alpha$ does not require this full implication. Indeed, assume the bound
$\|\psi\|_{L^\infty_t((-T_0,T_1); H^{\delta_1}(\R^2))}<\infty$ for some small~$\delta_1>0$ (we will obtain such bounds via frequency envelopes
 with $0<\delta_1<\sigma$). For any fixed time~$t\in (-T_0,T_1)$ one now has with $P_k$ being the usual Littlewood-Paley projections to
frequency~$2^k$,
\begin{align*}
 \| P_\ell \phi_\alpha\|_{H^{\delta_2}} & = \| P_\ell [e^{i\sum_{j=1}^2 \Delta^{-1}\del_j \phi^1_j} \psi_\alpha]\|_{H^{\delta_2}} \\
&\le \| P_\ell [P_{<\ell-10} (e^{i\sum_{j=1}^2 \Delta^{-1}\del_j \phi^1_j})  P_{[\ell-10,\ell+10]} \psi_\alpha    ]\|_{H^{\delta_2}}  \\
&\quad + \| P_\ell [P_{[\ell-10,\ell+10]} (e^{i\sum_{j=1}^2 \Delta^{-1}\del_j \phi^1_j})  P_{<\ell+15} \psi_\alpha    ]\|_{H^{\delta_2}}  \\
&\quad + \sum_{k>\ell+10} \| P_\ell [P_k (e^{i\sum_{j=1}^2 \Delta^{-1}\del_j \phi^1_j})  P_{k+O(1)} \psi_\alpha    ]\|_{H^{\delta_2}}  \\
&\les \|    P_{[\ell-10,\ell+10]} \psi_\alpha    \|_{H^{\delta_2}}  +  \| P_{[\ell-10,\ell+10]}
(e^{i\sum_{j=1}^2 \Delta^{-1}\del_j \phi^1_j}) \|_{H^{\delta_2}} \|P_{<\ell+15}\psi_\alpha\|_{\infty} \\
&\quad + \sum_{k>\ell+10} \|  P_k (e^{i\sum_{j=1}^2 \Delta^{-1}\del_j \phi^1_j})\|_{H^{\delta_2}} \|  P_{k+O(1)} \psi_\alpha \|_{\infty}
\end{align*}
Next, one has the bounds
\[
 \big\| \nabla_x\, e^{i \Delta^{-1} \sum_{j=1}^2 \del_j \phi^1_j}  \big\|_{\ener} \les \|\phi^1_j\|_{\ener},\qquad
\|P_{<\ell+15} \psi_\alpha\|_{L^\infty_x} \les 2^{(1-\delta_1)\ell} \|\psi_\alpha\|_{H^{\delta_1}}
\]
where the first one is admissible due to {\em energy conservation for the derived wave map}, see~\eqref{eq:phi_ener}.
In conclusion,
\[
 \| P_\ell \phi_\alpha\|_{H^{\delta_2}} \les \|P_{\ell+O(1)} \phi_\alpha\|_{H^{\delta_2}} + 2^{(\delta_2-\delta_1)\ell}\|\phi\|_{L^2_x} \|\psi\|_{H^{\delta_1}}
\]
Summing over $\ell\ge0$ yields
\begin{equation}\label{eq:phi_bd} \|\phi\|_{L^\infty_t((-T_0,T_1); H^{\delta_2}(\R^2))}<\infty
\end{equation}
  By the subcritical existence theory
of Klainerman and Machedon, see~\cite{KM1}--\cite{KM3} as well as~\cite{KS1}, \cite{KS2}, the solution can now be extended smoothly beyond
this time interval. More precisely, the device of frequency
envelopes allows one to place the Schwartz data in $H^s(\R^2)$ for all $s>0$ initially, and as it turns out, also for all times
provided $s>0$ is sufficiently small. The latter claim is of course the entire objective of this paper.
 We should also remark that we bring~\eqref{eq:psi_wave} into play only because it fits into the
framework of the spaces from~\cite{T1} and~\cite{Tat2}. This will
allow us to obtain the crucial energy estimate for solutions
of~\eqref{eq:psi_wave}, whereas it is not clear   how to do this
directly for the system~\eqref{eq:psisys1}, \eqref{eq:psisys2}. As
already noted in~\cite{Krieger}, the price one pays for passing
to~\eqref{eq:psi_wave} lies with the {\em initial conditions}, or
more precisely, the time derivative $\del_t \psi_\alpha(0,\cdot)$.
While $\psi_\alpha(0,\cdot)$ only involves one derivative of the
wave map~$\bfu$, this time derivative involves two. This will force
us to essentially ``randomize'' the initial time.

\subsection{The small data theory}

In this section we give a very brief introduction to the spaces
which are needed to control the $\psi$ system~\eqref{eq:psisys1},
\eqref{eq:psisys2}, and~\eqref{eq:psi_wave}. A systematic
development will be carried out in Section~\ref{sec:spaces} below,
largely following~\cite{T2} (we do need to go beyond both~\cite{T2}
and~\cite{Krieger} in some instances such as by adding the sharp
Strichartz spaces with the Klainerman-Tataru gain for small scales,
and by eventually modifying $\|\cdot\|_{S[k]}$ to the stronger
$\trip\cdot\trip_{S[k]}$ which allows for a high-high gain in the
$S\times S\to L^2_{tx}$ estimate).  First note that it is not
possible to bound the trilinear nonlinearities in this system in
Strichartz spaces due to slow dispersion in dimension two. Moreover,
it is not possible to adapt the $X^{s,b}$-space of the subcritical
theory to the scaling invariant case as this runs into logarithmic
divergences. For this reason, Tataru~\cite{Tat2} devised a class of
spaces which resolve these logarithmic divergences. His idea was to
allow characteristic frames of reference. More precisely, fix
$\omega\in S^{1}$  and define
\[
\theta_{\omega}^{\pm}:= (1,\pm\omega)/\sqrt{2},\quad t_{\omega} :=
(t,x)\cdot \theta_{\omega}^+ ,\quad x_{\omega} := (t,x)- t_{\omega}
\theta_{\omega}^+ \] which are the coordinates defined by a
generator on the light-cone. Now suppose that $\psi_i$ are free
waves such that $\psi_1$ is Fourier supported on $1\le|\xi|\le 2$,
and both $\psi_2$ and~$\psi_3$ are Fourier supported on $|\xi|\sim
2^k$ where $k$ is large and negative. Finally, we also assume that
the three wave are in ``generic position'', i.e., that their Fourier
supports make an angle of about size one. Clearly, $2^{-k} \psi_1
\psi_2\psi_3$ is then a representative model for the nonlinearities
arising in~\eqref{eq:psi_wave}. With
\[
\psi_3(t,x) = \int_{\R^2} e^{i[t|\xi|+x\cdot\xi]} f(\xi)\, d\xi
\]
we perform the {\em plane-wave} decomposition $\psi_3(t,x) =
\int\phi_\omega(\sqrt{2}t_\omega)\,d\omega$ where
\[
\phi_\omega(s) := \int e^{irs} f(r\omega)\, rdr
\]
By inspection,
\begin{equation}
 \int \|\phi_\omega\|_{L^2_{t_\omega} L^\infty_{x_\omega}}\,d\omega
 \les
 2^{\frac{k}{2}} \|\psi_3\|_{L^\infty_t L^2_x}
\label{eq:inspec} \end{equation}
 Hence,
\begin{align*}
  2^{-k} \int \| \phi_\omega\,
\psi_2\psi_3\|_{L^1_{t_\omega}L^2_{x_\omega}}\,d\omega &\les 2^{-k}
\int \| \phi_\omega \|_{L^2_{t_\omega} L^\infty_{x_\omega}}
\,d\omega \|
\psi_1\psi_2\|_{L^2_{t_\omega}L^2_{x_\omega}} \\
&\les  \| \psi_3\|_{\ener}\,  \| \psi_1\|_{\ener} \|\psi_2\|_{\ener}
\end{align*}
which is an example\footnote{Note that one does not obtain a gain in
this case. This fact will be of utmost importance in this paper,
forcing us to use a ``twisted'' wave equation resulting from these
high-low-low interactions in the linearized trilinear expressions.}
of a {\em trilinear estimate} which will be studied systematically
in Section~\ref{sec:trilin}. Here we used both~\eqref{eq:inspec} and
the standard bilinear $L^2_{tx}$ bilinear $L^2$-bound for waves with
angular separation:
\[
\| \psi_1\psi_2\|_{L^2_{t_\omega}L^2_{x_\omega}}  = \|
\psi_1\psi_2\|_{L^2_{t}L^2_{x}} \les 2^{\frac{k}{2}}
\|\psi_2\|_{\ener} \|\psi_1\|_{\ener}
\]
This suggests introducing an {\em atomic space} with atoms
$\psi_\omega$ of Fourier support $|\xi|\sim 1$ and satisfying
\[
\|\psi_\omega\|_{L^1_{t_\omega} L^2_{x_\omega}} \le 1
\]
as part of the space $N[0]$ which holds the nonlinearity (the zero
here refers to the Littlewood-Paley projection $P_0$. Below, we
refer to this space as $\NF$).  In addition, the space defined
by~\eqref{eq:inspec} is also an atomic space and should be
incorporated in the space~$S[k]$ holding the solution at
frequency~$2^k$ (we refer to this below as the $\PW$ space). By
duality to $L^1_{t_\omega} L^2_{x_\omega}$ in~$N[0]$, we then expect
to see $L^\infty_{t_\omega} L^2_{x_\omega}$ as part of~$S[0]$. The
simple observation here (originating in~\cite{Tat2}) is that one can
indeed bound the energy along a characteristic frame
$(t_\omega,x_\omega)$ of a free wave as long as its Fourier support
makes a positive angle with the direction~$\omega$. Indeed, recall
the local energy conservation identity $\del_t e -
{\mathrm{div}}(\del_t\psi \nabla\psi)=0$ for a free wave where
\[
e=\frac12(|\del_t\psi|^2+|\nabla\psi|^2)
\]
is the energy density, over a region of the form $\{-T\le t\le
T\}\cap \{t_\omega>a\}$. From the divergence theorem one obtains
that
\[
\int_{t_\omega=a} \chi_{[-T\le t\le T]} |\omega^\perp \nabla
\psi|^2\, d \calL^2 \les \|\psi\|_{\ener}^2
\]
where $\calL^2$ is the planar Lebesgue measure on $\{t_\omega=a\}$.
Sending $T\to\infty$ and letting $\rho$ denote the distance between
$\omega$ and the direction of the Fourier support of~$\psi|_{t=0}$,
one concludes that
\[
\|\psi\|_{L^\infty_{t_\omega} L^2_{x_\omega}} \les \rho^{-1}
\|\psi\|_{\ener}
\]
Hence, we should include a piece
\[
\sup_{\omega\not\in 2\kappa} d(\omega,\kappa)
\|\psi\|_{L^\infty_\omega L^2_{x_\omega}}
\]
in the norm $S[0]$ holding $P_0\psi$ provided $\psi$ is a wave
packet oriented along the cone of dimensions $1\times 2^k\times
2^{2k}$, projecting onto an angular sector  in the $\xi$-plane
associated with the cap~$\kappa\subset S^1$, where $\kappa$ is of
size~$2^k$ (this is called $\NF^*$ below).

Recall that we have made a genericity assumption which guaranteed
that the Fourier supports were well separated in the angle. In order
to relax this condition, it is essential to invoke  the usual device
of {\em nullforms} which cancel out parallel interactions. One of
the discoveries of~\cite{Krieger} is a genuinely trilinear nullform
expansion, see~\eqref{eq:nullexp} and~\eqref{eq:nullexp2}, which
exploit the relative position of all three waves simultaneously. It
seems impossible to reduce the trilinear nonlinearities
of~\eqref{eq:psi_wave} exclusively to the easier bilinear ones.

 It is shown in~\cite{Tat2} (and then also
in~\cite{T2} which develops much of the functional framework that we
use, as well as~\cite{Krieger}) that in low dimensions (especially
$n=2$ but these spaces are also needed for $n=3$), these nullframe
spaces are strong enough --- in conjunction with more traditional
scaling invariant $X^{s,b}$ spaces
--- to bound the trilinear nonlinearities,  as well as weak enough to
allow for an energy estimate to hold. This then leads modulo passing
to an appropriate gauge to the small energy theory.

The norm $\|\cdot\|_S$ in~\eqref{eq:Sbound2} is of  the form
$\|\psi\|_S:= \Big(\sum_{k\in\Z}
\|P_k\psi\|_{S[k]}^2\Big)^{\frac12}$ where $S[k]$ is built from
$\ener$, critical $X^{s,b}$, $L^4_t L^\infty_x$ Strichartz norms, as
well as the null-frame spaces which we just described.

\subsection{The Bahouri-Gerard concentrated compactness method}

We now come to the core  of the argument, namely the Bahouri-Gerard
type decomposition and the associated perturbative argument.

In~\cite{Ger} P.~G\'erard considered defocusing semilinear wave
equations in $\R^{3+1}$ of the form $\Box u + f(u)=0$ with data
given by a sequence $(\phi_n,\psi_n)$ of energy data going weakly to
zero. Denote the resulting solutions to the nonlinear problem
by~$u_n$, and the free waves with the same data by~$v_n$. G\'erard
proved that provided $f(u)$ is {\em subcritical} relative to energy
then
\[
\|u_n-v_n\|_{L^\infty(I;\calE)}\to 0\quad \text{\ as\ }n\to\infty
\]
where $\calE$ is the energy space. In contrast, for this to hold for
the energy critical problem he found via the concentrated
compactness method of P.\ L.~Lions that it is necessary and
sufficient that $\|v_n\|_{L^\infty(I;L^6(\R^3))}\to0$. In other
words, the critical problem experiences a {\em loss of compactness}.

The origin of this loss of compactness, as well as the meaning of
the $L^6$ condition were later made completely explicit by
 Bahouri-Gerard~\cite{BG}. Their result  reads as follows:
 {\em Let $\{(\phi_n,\psi_n)\}_{n=1}^\infty \subset \dot H^1\times L^2(\R^3)$
be a bounded sequence, and define $v_n$ to be a free wave with these
initial data. Then there exists a subsequence $\{v_n'\}$ of
$\{v_n\}$, a finite energy free wave $v$, as well as free waves
$V^{(j)}$ and $(\eps^{(j)}, x^{(j)})\in (\R^+,\R^3)^{\Z^+}$ for
every $j\ge1$ with the property that for all $\ell\ge1$,
\begin{equation}\label{eq:BGdecomp}
v_n'(t,x)=v(t,x) + \sum_{j=1}^\ell \frac{1}{\sqrt{\eps_n^{(j)}}}
V^{(j)} \Big(\frac{t-t_n^{(j)}}{\eps_n^{(j)}},
\frac{x-x_n^{(j)}}{\eps_n^{(j)}} \Big) + w_n^{(\ell)}(t,x)
\end{equation}
where
\[
\limsup_{n\to\infty} \|w_n^{(\ell)}\|_{L^5_t(\R,L^{10}_x(\R^3))} \to
0\qquad\text{\ \ as\ \ }\ell\to\infty
\]
and for any $j\ne k$,
\[
\frac{\eps_n^{(j)}}{\eps_n^{(k)}} +
\frac{\eps_n^{(k)}}{\eps_n^{(j)}} + \frac{|x_n^{(j)}-x_n^{(k)}|+
|t_n^{(j)}-t_n^{(k)}| }{\eps_n^{(j)}} \to \infty \qquad \text{\ \
as\ }n\to\infty
\]
Furthermore, the free energy $E_0$ satisfies the following
orthogonality property:
\[
E_0(v_n') = E_0(v) + \sum_{j=1}^\ell E_0(V^{(j)}) +
E_0(w_n^{(\ell)}) + o(1)\qquad\text{\ \ as\ \ }n\to\infty
\] }

Note that this result characterized the loss of compactness in terms
of the appearance of concentration profiles $V^{(j)}$.  Moreover,
\cite{BG} contains an analogue of this result  for so-called
Shatah-Struwe solutions of the semi-linear problem $\Box u + |u|^4
u=0$ which then leads to another proof of the main result
in~\cite{Ger}. One of the main applications of their work was to
show the existence of a function $A:[0,\infty)\to[0,\infty)$ so that
every Shatah-Struwe solution satisfies the bound
\begin{equation}\label{eq:BGenerfunc}
\|u\|_{L^5_t(\R; L^{10}_x(\R^3))}\le A(E(u))
\end{equation}
where $E(u)$ is the energy associated with the semi-linear equation.
 This is proved  by
contradiction; indeed, assuming~\eqref{eq:BGenerfunc} fails, one
then obtains sequences of bounded energy solutions with
uncontrollable Strichartz norm which is then shown to contradict the
fact the nonlinear solutions themselves converge weakly to another
solution.
 The decomposition~\eqref{eq:BGdecomp}
compensates for the aforementioned loss of compactness by reducing
it precisely to the effect of the {\em symmetries}, i.e., dilation
and scaling. This is completely analogous to the elliptic (in fact,
variational) origins of the method of concentrated compactness, see
Lions~\cite{Lions} and Struwe~\cite{Struwe}.
 See~\cite{BG} for more details
and other applications.

The importance of~\cite{BG} in the context of wave maps is made
clear by the argument of Kenig, Merle~\cite{KeM1}, \cite{KeM2}. This
method, which will be described in more detail later in this
section, represents a general method for attacking global
well-posedness problems for energy critical equations such as the
wave-map problem. Returning to the Bahouri-Gerard decomposition, we
note that any attempt at implementing this technique for wave maps
encounters numerous serious difficulties. These are of course all
rooted in the difficult nonlinear nature of the
system~\eqref{eq:compat1}--\eqref{eq:derwmp2}. Perhaps the most
salient feature of our decomposition as compared to~\cite{BG} is
that the free wave equation no longer capture the correct asymptotic
behavior for large times; rather, the atomic components $V^{(j)}$
are defined as solutions of a covariant (or ``twisted'') wave
equation of the form
\begin{equation}\label{eq:covwave}
 \Box + 2iA_\alpha \del^\alpha
\end{equation}
where the magnetic potential $A_\alpha$ arises from linearizing the wave map
equation in the Coulomb gauge. More precisely, the magnetic term here captures
the  high-low-low interactions in the trilinear nonlinearities of the wave map
system where there is no apriori smallness gain.

\noindent In keeping with the Kenig-Merle method, the Bahouri-Gerard decomposition is
used to show the following: assume that a uniform bound
of the form
\[
\|\psi\|_{S}\leq C(E)
\]
for some function $C(E)$ fails for all finite energy levels $E$.
Then there must exist a weak wave map $\bfu_{\text{critical}}: (-T_{0},
T_{1})\to S$ to a compact Riemann surface uniformized by
${\mathbb{H}}^{2}$, which enjoys certain compactness properties. In
the final part of the argument we then need to rule out the
existence of such an object. Arguing by contradiction, we now assume
there is a sequence of Schwartz class (on fixed time slices) wave
maps $\bfu^{n}: (-T_{0}^{n}, T_{1}^{n})\times\R^{2}\to{\mathbb{H}}^{2}$
with the properties that
\begin{itemize}
\item $\|\psi^{n}\|_{L_{x}^{2}}\to \Ecrit$
\item
$\lim_{n\to \infty}\|\psi^{n}\|_{S((-T_{0}^{n},
T_{1}^{n})\times\R^{2})}=\infty$
\end{itemize}
Thus all these wave maps have $t=0$ in their domain of definition.
Roughly speaking, we shall proceed along the following steps.
First, recall that the Bahouri-Gerard theorem is a genuine  phase-space result
in the sense that it identifies the main asymptotic carriers of energy {\em which are not pure radiation}, which would then sit in $w_n^{(\ell)}$.
This refers to the free waves $V^{(j)}$ above, which are ``localized'' in frequency (namely at scale $(\eps_n^{(j)})^{-1}$)
as well as in physical spaces (namely around the space-time points $(t_n^{(j)}, x_n^{(j)})$). The procedure of filtering out
the scales $\eps_n^{(j)}$ is due to Metivier-Schochet, see~\cite{MS}.

\begin{enumerate}
 \item {\it{Bahouri-Gerard I: filtering out frequency blocks.}}

If we apply the frequency localization procedure of Metivier-Schochet
 to the derivative components $\phi^{n}_{\alpha}=(\frac{\partial_{\alpha}{\bfx}^{n}}{{\bfy}^{n}},\quad
 \frac{\partial_{\alpha}{\bfy}^{n}}{{\bfy}^{n}})$ at time $t=0$, we run into the problem that the resulting frequency components
 are not necessarily related to an actual map from $\R^{2}\to  {\mathbb{H}}^{2}$. We introduce a
 procedure to obtain a frequency decomposition which is ``geometric'', i.e., the frequency localized
 pieces are themselves derivative components of maps from $\R^{2}\to{\mathbb{H}}^{2}$.
\item  {\it{Refining the considerations on frequency localization; frequency localized approximative maps}}.
In order to deal with the non-atomic (in the frequency sense)
derivative components, which may still have large energy, we need to
be able to truncate the derivative components arbitrarily in
frequency  while still retaining the geometric interpretation. Here
we shall use arguments just as in the first step to allow us to
``build up'' the components $\psi^{n}_{\alpha}$ from low frequency
ones.  In the end, we of course need to show that for some
subsequence of the $\psi^{n}_{\alpha}$, the frequency support is
essentially atomic. If this were to fail, we
 deduce an apriori bound on $\|\psi^{n}_{\alpha}\|_{S((-T_{0}^{n}, T_{1}^{n})\times\R^{2})}$.
\item {\it{Assuming the presence of a lowest energy non-atomic type component, establish an apriori estimate for its
nonlinear evolution.}} This requires deriving energy estimates for the covariant wave equation~\eqref{eq:covwave}.
\item {\it{Bahouri-Gerard II, applied to the first atomic frequency component}}. Here, assuming that  we have constructed
the first ``low frequency approximation'' in the previous step, we
need to filter out the physically localized pieces. This is where we
have to deviate from  Kenig-Merle: instead of the free wave
operator, we need to use a covariant wave operator to model the
asymptotics as $t\to\pm \infty$. Again a lot of effort needs to be
expended on showing that the components we obtain are actually the
Coulomb derivative components of Schwartz maps from
$\R^{2}\to{\mathbb{H}}^{2}$, up to arbitrarily small errors in
energy. Once we have this, we can then use the result from the
stability section in order to construct the time evolution of these
pieces and obtain their apriori dispersive behavior.
\item  {\it{Bahouri-Gerard II; completion}}. Here we repeat Steps~3 and 4 for the ensuing frequency pieces, to complete the estimate
for the $\psi^{n}_{\alpha}$. The conclusion is that upon choosing
$n$ large enough, we arrive at a contradiction, unless there is
precisely one frequency component and precisely one atomic physical
component forming that frequency component. These are the data that
then gives rise to the weak wave map with the desired compactness
properties.
\end{enumerate}

\subsection{The Kenig-Merle agument}

In \cite{KeM2}, \cite{KeM1},  Kenig and Merle developed an approach
to the global wellposendess for defocusing  energy critical
semilinear Schr\"odinger and wave equations; moreover, their
argument yields a blowup/global existence dichotomy in the focusing
case as well, provided the energy of the wave lies beneath a certain
threshold. See~\cite{CKM} for an application of these ideas to wave
maps.

Let us give a brief overview of their argument. Consider
\[
\Box u + u^5 =0
\]
in $\R^{1+3}$ with data in $\dot H^1\times L^2$. It is standard that
this equation is well-posed for small data provided we place the
solution in the energy space intersected with  suitable Strichartz
spaces. Moreover, if $I$ is the maximal interval of existence, then
necessarily $\|u\|_{L^8_t(I;L^8_x(\R^3))}=\infty$ and the energy
$E(u)$ is conserved.

Now suppose $\Ecrit$ is the minimal energy with the property that
all solutions in the above sense with $E(u)<\Ecrit$ exist globally
and satisfy $\|u\|_{L^8_t(\R;L^8_x(\R^3))}<\infty$. Then by means of
the Bahouri-Gerard decomposition, as well as the perturbation theory
for this equation one concludes that a critical solution $u_C$
exists on some interval~$I^*$ and that
$\|u\|_{L^8_t(I^*;L^8_x(\R^3))}=\infty$. Moreover, by similar
arguments one obtains the crucial property that the set
\[
K:= \{ \big(\lambda^{\frac12}(t) u(\lambda(t)
(x-y(t)),t),\lambda^{\frac32}(t) \del_t u(\lambda(t)
(x-y(t)),t)\big)\::\: t\in I\}
\]
is precompact in $\dot H^1\times L^2(\R^3)$ for a suitable path
$\lambda(t),y(t)$. To see this, one applies the Bahouri-Gerard
decomposition to a sequence $u_n$ of solutions with energy
$E(u_n)\to \Ecrit$ from above.  The logic here is that due to the
minimality assumption on~$\Ecrit$ {\em only a single limiting
profile can arise in}~\eqref{eq:BGdecomp} up to errors that go to
zero in energy as $n\to\infty$. Indeed, if this were not the case
then due to fact that the profiles diverge from each other in
physical space as $n\to\infty$ one can then apply the perturbation
theory to conclude that each of the individual nonlinear evolutions
of the limiting profiles (which exist due to the fact that their
energies are strictly below~$\Ecrit$) can be superimposed to form a
global nonlinear evolution, contradicting the choice of the sequence
$u_n$. The fact that $\ell=1$ allows one to rescale and re-translate
the unique limiting profile to a fixed position in phase space
(meaning spatial position and spatial frequency) which then gives
the desired nonlinear evolution~$u_C$. The compactness follows by
the same logic: assuming that it does not hold, one then obtains a
sequence $u_C(\cdot, t_n)$ evaluated at times $t_n\in I^*$
converging to an endpoint of~$I^*$ such that for $n\ne n'$, the
rescaled and translated versions of~$u_C(\cdot, t_n)$
and~$u_C(\cdot, t_{n'})$  remain at a minimal positive  distance
from each other in the energy norm. Again one applies Bahouri-Gerard
and finds that $\ell=1$ by the choice of~$\Ecrit$ and perturbation
theory. This gives the desired contradiction. The compactness
property is of course crucial; indeed, for illustrative purposes
suppose that $u_C$ is of the form
\[
u_C(t,x)=\lambda(t)^{\frac12}U(\lambda(t)(x-x(t)))
\]
where $\lambda(t)\to\infty$ as $t\to1$, say. Then $u_C$ blows up at
time~$t=1$ (in the sense that the energy concentrates at the tip of
a cone) and
\[
\lambda(t)^{-\frac12}u_C(\lambda(t)^{-1}x+x(t))=U(x)
\]
is compact for $0\le t<1$. Returning to the Kenig-Merle argument,
the logic is now to show that $u_C$ acts in some sense like a
blow-up solution, at least if $I^*$ is finite in one direction.

 The second half of the Kenig-Merle approach then consists of
a rigidity argument which shows that a $u_C$ with the stated
properties cannot exist. This is done mainly by means of the
conservation laws, such as the Morawetz and energy identities. More
precisely, the case where $I^*$ is finite at one end is reduced to
the self-similar blowup scenario. This, however, is excluded by
reducing to the stationary case and an elliptic analysis which
proves that the solution would have to vanish. If $I^*$ is infinite,
one basically faces the possibility of stationary solutions which
are again shown not to exist.

\medskip For the case of wave maps, we follow the same strategy.
More precisely, our adaptation of the Bahouri-Gerard decomposition
to wave maps into~$\Hyp^2$ leads to a critical wave map with the
desired compactness properties. In the course of our proof, it will
be convenient to project the wave map onto a compact Riemann
surface~$\calS$ (so that we can avail ourselves of the {\em
extrinsic formulation} of the wave map equation). However, it will
be important to work simultaneously with this object as well as the
lifted one which takes its values in~$\Hyp^2$ (since it is for the
latter that we have a meaningful wellposedness theory for maps with
energy data).

 The difference from~\cite{KeM1} lies mainly with the rigidity part.
 In fact, in our context the conservation laws are by themselves not sufficient to
yield a contradiction. This is natural, since the geometry of the
target will need to play a crucial role. As indicated above, the two
scenarios that are lead to a contradiction are the self-similar
blowup supported inside of a light-cone and the stationary weak wave
map, which is of course a weakly harmonic map (which cannot exist
since the target $\calS$ is compact with negative curvature). The
former is handled as follows: in self-similar coordinates, one
obtains a harmonic map defined on the disk with the hyperbolic
metric and with finite energy (the stationarity is derived as
in~\cite{KeM1}). Moreover, there is the added twist that one
controls the behavior of this map at the boundary in the trace sense
(in fact, one shows that this trace is constant). Therefore, one can
apply the boundary regularity version of Helein's theorem which was
obtained by Qing~\cite{Quing}. Lemaire's theorem~\cite{Lemaire} then
yields the constancy of the harmonic map, whence the contradiction
(for a version of this argument under the apriori assumption of
regularity all the way to the boundary see
Shatah-Struwe~\cite{SStruwe}).

\subsection{An overview of the paper}

The paper is essentially divided into two parts: the {\bf{modified
Bahouri-Gerard method}} is carried out in its entirety starting with
Section~\ref{sec:spaces}, and ending with Section~\ref{sec:BG}.
Indeed, all that precedes Section~\ref{sec:BG} leads to this
section, which is the core of this paper. The {\bf{Kenig-Merle
method adapted to Wave Maps}} is then performed in the much shorter
section~\ref{sec:conclusion}. We commence by describing in detail
the contents of Section~\ref{sec:spaces} to Section~\ref{sec:BG}.

\subsubsection{Preparations for the Bahouri-Gerard process.} As explained above,
we describe admissible  wave maps $u: R^{2+1}\to \Hyp^2$ mostly in terms of the associated
Coulomb derivative components $\psi_\alpha$. Our goals then are to

\begin{itemize}
\item {\it{(1): Develop a suitable functional framework}}, in
particular a space-time norm $\|\psi\|_{S(\R^{2+1})}$, together with
time-localized versions $\|\psi\|_{S([I\times\R^{2})}$ for closed
time intervals $I$, which have the property that
\[
\limsup_{I\subset \tilde{I}}\|\psi\|_{S(I\times\R^{2})}<\infty
\]
for some open interval $\tilde{I}$ implies that the underlying wave map $u$ can be extended
smoothly and admissibly beyond any endpoint of $\tilde{I}$, provided such exists.
\item
{\it{(2): Establish an a priori bound}} of the form
\[
\|\psi\|_{S(I\times\R^2)}\leq C(E)
\]
for some function $C: \R_{+}\to \R_{+}$ of the energy $E$. This
latter step will be accomplished by the Bahouri Gerard procedure,
arguing by contradiction.
\end{itemize}

We first describe (1) above in more detail: in
Section~\ref{sec:spaces}, we introduce the norms $\|\cdot\|_{S[k]}$,
$\|\cdot\|_{N[k]}$, $k\in\Z$, which are used to control the
frequency localized components of $\psi$ and the nonlinear source
terms, respectively. The norm $\|\cdot \|_{S}$ is then obtained by
square summation over all frequency blocks. The basic paradigm for
establishing estimates on $\psi$ then is to formulate a wave
equation
\[
\Box \psi= F
\]
or more accurately typically in frequency localized form
 \[
\Box P_0\psi= P_0F,
\]
and to establish bounds for $\|P_0F\|_{N[0]}$ which may then be fed
into an energy inequality, see Section~\ref{subsec:energy}, which
establishes the link between the $S$ and $N$-spaces. In order to be
able to estimate the nonlinear source terms $F$, we need to
manipulate the right-hand side of \eqref{eq:psi_wave}, making
extensive use of \eqref{eq:dyn_dec}. The precise description of the
actual nonlinear source terms that we will use for $F$ is actually
rather involved, and given in Section~\ref{sec:hodge}. In order to
estimate the collection of trilinear as well as higher order terms,
we carefully develop the necessary estimates in
Sections~\ref{sec:bilin}, \ref{sec:trilin}, as well as
\ref{sec:quintic}. We note that the estimates in \cite{Krieger},
while similar, are not quite strong enough for our purposes, since
we need to gain in the largest frequency in case of high-high
cascades. This requires us to subtly modify the spaces by comparison
to loc.\ cit. Moreover, the fact that we manage here to build in
sharp Strichartz estimates allows us to replace several arguments in
\cite{Krieger} by more natural ones, and we opted to make our
present account as self-contained as possible.

\noindent With the null-form estimates from
Sections~\ref{sec:bilin},
 \ref{sec:trilin},~\ref{sec:quintic} in hand, we establish the role
of $\|\cdot\|_S$ as a ``regularity controlling'' device in the sense
of~(1) above in Section~\ref{sec:perturb}, see
Proposition~\ref{BlowupCriterion}. The proof of this reveals a
somewhat unfortunate feature of our present setup, namely the fact
that working at the level of the differentiated wave map system
produces sometimes too many time derivatives, which forces us to use
somewhat delicate ``randomization'' of times arguments. In
particular, in the proof of all apriori estimates, we need to
distinguish between a ``small time'' case (typically called Case~1)
and a ``long time'' Case~2, by reference to a fixed frequency scale.
In the short time case, one works exclusively in terms of the
div-curl system, while in the long-time case, the wave equations
start to be essential.

\noindent Section~\ref{sec:perturb} furthermore explains the
well-posedness theory at the level of the $\psi_\alpha$, see the
most crucial Proposition~\ref{prop:ener_stable}. We do not prove
this proposition in Section~\ref{sec:perturb}, as it follows as a
byproduct of the core perturbative Proposition~\ref{PsiBootstrap} in
Section~\ref{sec:BG}. Proposition~\ref{prop:ener_stable} and the
technically difficult but fundamental Lemma~\ref{BasicStability}
allow us to define the ``Coulomb wave maps propagation'' for a tuple
$\psi_\alpha$, $\alpha=0,1,2$ which are only $L^2$ functions at time
$t=0$, provided the latter are the $L^2$-limits of the Coulomb
components of admissible maps. Indeed, this concept of propagation
is independent of the approximating sequence chosen and satisfies
the necessary continuity properties.

We also formulate the concept of a ``wave map at infinity'' at the
level of the Coulomb components, see Proposition~\ref{prop:waveops}
and the following Corollary~\ref{temporallyunbounded}. Again the
proofs of these results will follow as a byproduct of the
fundamental Proposition~\ref{PsiBootstrap} and
Proposition~\ref{BGIIHard} in the core Section~\ref{sec:BG}.

In Section~\ref{sec:bmo}, we develop some auxiliary technical
 tools from harmonic analysis which will allow us to implement the first stage of the
 Bahouri Gerard process, namely crystallizing frequency atoms from an ``essentially singular'' sequence of admissible wave maps.
These tools are derived from the imbedding $\dot
B^1_{2,\infty}(\R^2)\to\bmo$ as well as weighted (relative to $A_p$)
Coifman-Meyer commutator bounds.

As mentioned before, Section~\ref{sec:BG} is the core of the present
paper. In Section~\ref{subsec:BGstep1}, starting with an {\em
essentially singular} sequence ${\bf{u}}^n$ of admissible wave maps
with deteriorating bounds, i.e., $\|\psi^n_\alpha\|_{S}\to \infty$
as $n\to\infty$ but with the crucial criticality condition
$\lim_{n\to\infty}E({\bf{u}}^n)=\Ecrit$, we show that the derivative
components $\phi^{n}_\alpha$ may be decomposed as a sum
\[
\phi^n_\alpha=\sum_{a=1}^A\phi^{na}_\alpha+w_\alpha^{nA}
\]
where the $\phi^{na}_\alpha$ are derivative components of admissible
wave maps which have frequency supports ``drifting apart'' as
$n\to\infty$, while the error $w_\alpha^{nA}$ satisfies
\[
\limsup_{n\to\infty}\|w_\alpha^{nA}\|_{\dot{B}^{0}_{2, \infty}}<\delta
\]
provided $A\geq A_0(\delta)$ is large enough.

\noindent  In Section~\ref{subsec:BGstep2}, we then select a number
of ``principal'' frequency atoms $\phi^{na}$, $a=1,2,\ldots, A_0$,
as well as a (potentially very large) collection of ``small atoms''
$\phi^{na}$, $a=A_0+1,\ldots, A$. We order these atoms by the
frequency scale around which they are supported starting with those
of the lowest frequency. The idea now is as follows: under the
assumption that there are at least two frequency atoms, or else in
case of only one frequency atom that it has energy $<\Ecrit$, we
want to obtain a contradiction to the essential criticality of the
underlying sequence ${\bf{u}}^n$. To achieve this, we define in
Section~\ref{subsec:BGstep2}  sequence of approximating wave maps,
which are essentially obtained by carefully truncating the initial
data sequence $\phi^{na}$ in frequency space.

\noindent In Section~\ref{subsec:BGstep3}, we establish an apriori
bound for the lowest frequency approximating map which comprises all
the minimum frequency small atoms as well as the component of the
small Besov error of smallest frequency, see
Proposition~\ref{ControlNonatomicComponent1}. The proof of this
follows again by truncating the data suitably in frequency space,
and applying an inductive procedure to a sequence of approximating
wave maps. This hinges crucially on the core perturbative result
Proposition~\ref{PsiBootstrap}, which plays a fundamental role in
the paper. The  main technical difficulty encountered in the proof
of the latter comes from the issue of {\bf{
fungiblity}}\footnote{This term appears to have been coined by T.\
Tao}: given a schematically written expression
\[
\nabla_{x,t}[\partial^{\nu}\epsilon A_{\nu}]
\]
which is {\em linear} in the perturbation (so that we cannot perform
a bootstrap argument based solely on the smallness on $\epsilon$
itself),  while $A_{\nu}$ denotes some null-form depending on
apriori controlled components $\psi$. Fungibility means the property
that upon suitably truncating time into finitely many intervals
$I_j$ {\em whose number only depends on $\|\psi\|_{S}$}, one may
bound the expression by
\[
\|\nabla_{x,t}[\partial^{\nu}\epsilon A_{\nu}]\|_{N(I_j\times\R^2)}\ll \|\epsilon\|_{S}
\]
In  other words, by shrinking the time interval, we ensure that we can iterate the term away.
While this would be straightforward provided we had an estimate for $\|A_{\nu}\|_{L_t^1L_x^\infty}$
(which is possible in space dimensions $n\geq 4$), in our setting, the spaces
are much too weak and complicated. Our way out of this impasse is to build those terms for
which we have no obvious fungibility into the linear operator, and thereby form a new operator
\[
\Box_A\epsilon:=\Box\epsilon+2i\partial^{\nu}\epsilon A_{\nu}
\]
with a magnetic potential term. Fortunately, it turns out that if
$A_\nu$ is supported at much lower frequencies than $\epsilon$
(which is precisely the case where fungibility fails), one can
establish an approximate energy conservation result, which in
particular gives apriori control over a certain constituent of
$\|\cdot\|_{S}$. With this in hand, one can complete the bootstrap
argument, and obtain full control over $\|\epsilon\|_{S}$.

\noindent  Having established control over the lowest-frequency
``essentially non-atomic'' approximating wave map in
Section~\ref{subsec:BGstep3}, we face the task of ``adding the first
large atomic component'', $\phi^{n1}$. It is here that we have to
depart crucially from the original method of Bahouri-Gerard: instead
of studying the free wave evolution of the data, we extract
concentration cores by applying the ``twisted'' covariant evolution
associated with
\[
\Box_{A^n}u=0,
\]
which is essentially defined as above. The key property that makes
everything work is an almost exact energy conservation property
associated with its wave flow. This is a rather delicate point: it
relies on the mixed-Lebesgue type endpoint bilinear improved
Strichartz bound of Wolff~\cite{W}, Tao~\cite{T8}, \cite{T9},  and
Tataru~\cite{TatWolff}.

\noindent  It then requires a fair amount of work to show that the
profile decomposition at time $t=0$ in terms of covariant free waves
is ``geometric'', in the sense that the concentration profiles can
indeed by approximated by the Coulomb components of admissible
maps, up to a constant phase shift, see Proposition~\ref{ConcentrarionProfileApprox}. \\
Finally, in Proposition~\ref{BGIIHard} we show that we may evolve the data including the first
large frequency atom, provided all concentration cores have energy strictly less than $\Ecrit$.

\noindent  As most of the work has been done at this point, adding
on the remaining frequency atoms in Section~\ref{subsec:BGEnd} does
not provide any new difficulties, and can be done by the methods of
the preceding sections.

\noindent  In conjunction with the results of
Section~\ref{sec:perturb}, we can then infer that given an
essentially singular sequence of wave maps ${\bf{u}}^n$, we may
select a subsequence of them whose Coulomb components
$\psi^n_\alpha$, up to re-scalings and translations, converge to a
limiting object $\Psi^\infty_\alpha(t, x)$, which is well-defined on
some interval $I\times\R^2$ where $I$ is either a finite time
interval or (semi)-infinite,  and the limit of the Coulomb
components of admissible maps there. Moreover, most crucially for
the sequel, $\Psi^\infty_\alpha(t, x)$ satisfies a remarkable {\em
compactness property}, see Proposition~\ref{cor:compactV}. This sets
the stage for the method of Kenig-Merle, which we adopt to the
context of wave maps.

\section{The spaces $S[k]$ and $N[k]$}\label{sec:spaces}

Sections~\ref{sec:spaces}--\ref{sec:trilin} develop the functional framework
needed to prove the energy and dispersive estimates required by
the wave map system~\eqref{eq:psisys1}--\eqref{eq:psi_wave}.
 The Banach spaces which appear in this context
go back to Tataru~\cite{Tat2}, but were introduced in this form by
Tao~\cite{T1}, and developed further by Krieger~\cite{Krieger}.  We
will largely follow the latter reference although there is much
overlap with~\cite{T2}. We emphasize that this section is completely
self-contained and presents all estimates in full detail. The
spatial dimension is~two throughout.

\subsection{Preliminaries}

As usual, $P_k$ denotes a Littlewood-Paley
projection\footnote{Strictly speaking, these are not true
projections since $P_k^2\ne P_k$, but we shall  nevertheless follow
the customary abuse of language of referring to them as projections.
The same applies to smooth localizers to other regions in Fourier
space.} to frequencies of size~$2^k$. More precisely, let $m_0$ be a
nonnegative smooth, even,  bump function supported in $|\xi|<4$ and set
$m(\xi):=m_0(\xi)-m_0(2\xi)$. Then
\[
 \sum_{k\in\Z}  m(2^k \xi)=1  \quad \forall\; \xi\in\R^2\setminus\{0\}
\]
and $\widehat{P_k f}(\xi):= m(2^{-k} \xi)\hat{f}(\xi)$. The
operator $Q_j$ projects to modulation~$2^j$, i.e.,
\[
 \widehat{Q_j\phi}(\xi,\tau) := m(2^{-j}(|\tau|-|\xi|)) \wh{\phi}(\tau,\xi)
\]
with $\hat{\cdot}$ referring to the space-time Fourier transform.
Similarly,
\[
 \widehat{Q_j^{\pm}\phi}(\xi,\tau) := m\big(2^{-j}(|\tau|-|\xi|)\big) \chi_{[\pm\tau>0]} \wh{\phi}(\tau,\xi)
\]
Then the relevant $\dot X_k^{s,b,q}$ spaces here are defined as
\[
 \| \phi\|_{\dot X_k^{s,p,q}} := 2^{sk} \Big( \sum_j 2^{jq} \|P_k Q_j \phi\|_{\Ltwotx}^q \Big)^{\frac1q}
\]
If $P_k\phi=\phi$, then $\|\phi\|_{L_t^\infty L_x^2} \les \|\phi\|_{ \dot X_k^{0,\frac12,1}}$ as well
as $\|\phi\|_{L^\infty_{t,x}}\les \|\phi\|_{ \dot X_k^{1,\frac12,1}}$.

In what follows, $\calC_\ell$ is a collection of caps $\kappa\subset
S^1$ of size~$C^{-1} 2^\ell$ and finite overlap (uniformly bounded
in~$\ell$ and with $C$ some large absolute constant). There is an associated smooth partition of unity $
\sum_{\kappa\in\calC_\ell} a_\kappa(\omega) =1$ for all $\omega\in
S^1$, as well as projections $\widehat{P_\kappa f}(\xi) :=
a_\kappa\big(\wh{\xi}\big) \hat{f}(\xi)$ where $\wh{\xi}:=\frac{\xi}{|\xi|}$. By construction,
$P_{k,\kappa}:=P_k\circ P_\kappa$ is a projection to the
``rectangle'' \begin{equation}\label{eq:Rkappa_def} R_{k,\kappa}:=\{|\xi|\sim
2^k,\;\wh\xi\in\kappa\}\end{equation} in Fourier space.
\begin{figure*}[ht]
\begin{center}
\centerline{\hbox{\vbox{ \epsfxsize= 8.0 truecm \epsfysize=6.5
truecm \epsfbox{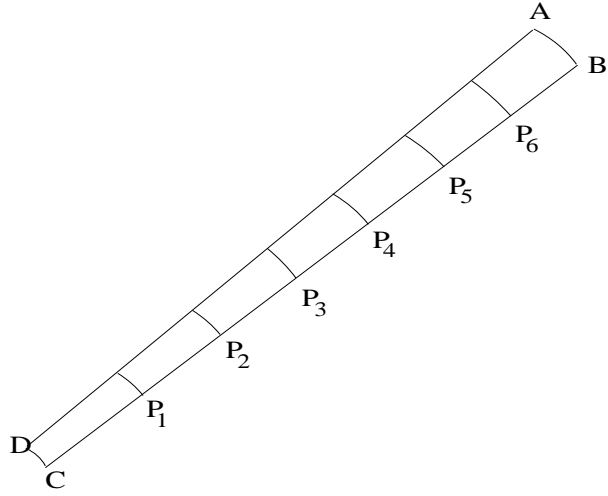}}}} \caption{Rectangles}
\end{center}
\end{figure*}
For space-time functions $F$ we shall follow the convention that
\[
 P_\kappa F =  [a_\kappa(\hat\xi) \chi_{[\tau>0]} \hat{F}(\xi,\tau)]^{\vee} +
[a_\kappa(-\hat\xi) \chi_{[\tau<0]} \hat{F}(\xi,\tau)]^{\vee}
\]
We will also encounter other rectangles $R$ which are obtained by
dividing $R_{k,\kappa}$ in the radial direction into $2^{-m}$ many
subrectangles of comparable size where $m<0$ is some integer
parameter (it will suffice for us to consider $\ell\le m\le0$ where
$\kappa\in\calC_\ell$). The collection of these rectangles will be
denoted by $\calR_{k,\kappa,m}$, and we introduce projections $P_R$
so that $\sum_{R\in\calR_{k,\kappa,m}} P_R = P_{k,\kappa}$. Figure~1
exhibits such a collection of rectangles. The sector $ABCD$ is of
length $2^k$ and width $2^{\ell+k}$, whereas the shorter segments
$\overline{AP_1}$, $\overline{P_1 P_2}$ etc.\ are of
length~$2^{k+m}$.

We shall frequently use Bernstein's inequality:  if
$\supp(\hat{\phi})\subset \Omega$, where $\Omega \subset\R^2$ is
measurable, then $\|\phi\|_p\les |\Omega|^{\frac{1}{q}-\frac{1}{p}}
\|\phi\|_q$ for any choice of~$1\le p\le q\le\infty$. We shall also
require the following variant of Bernstein's $L^2\to L^\infty$
bound, which is obtained by combining the standard form of this
bound with the $L^4_t L^\infty_x(\R^{1+2})$-Strichartz estimate for
the wave equation. This type of estimate appears in~\cite{T1},  but
the following formulation is from~\cite{Krieger}, which involves one
further localization on the Fourier side. We present the proof for
the sake of completeness.

\begin{lemma}
 \label{lem:KBern}
Let $\calD_{k,\ell}$ be a cover of $\{|\xi|\sim 2^k\}$ by disks of radius $2^{k+\ell}$. Then for all $j\le k$,
\begin{equation}
 \label{eq:JTao} \Big(\sum_{c\in \calD_{k,\ell}} \|P_c Q_j \phi\|_{L^2_t L^\infty_x}^2\Big)^{\frac12} \les
 2^{\frac{\ell}{2}} 2^{\frac{3k}{4}} 2^{\frac{j-k}{4}}\| \phi\|_{\Ltwotx}
\end{equation}
for any $\phi$ which is adapted to~$k$.
\end{lemma}
\begin{proof}
We follow the argument in~\cite{T1}, but use the small-scale
Strichartz estimate of Klainerman-Tataru at a crucial place, see
Lemma~\ref{lem:Strich} below. First,  set $j=0$, whence $k\gg1$.
Construct a Schwartz function $a(t)$ whose Fourier transform is
supported in $|\tau|\ll 1$, and which satisfies
\begin{equation}\nonumber
1=\sum_{s\in \mathbb{Z}}a^{3}(t-s)
\end{equation}
for all $t\in \R$. Then
\begin{equation}\nonumber\begin{split}
\|P_{c}Q_{0}\psi\|_{L_{t}^{2}L_{x}^{\infty}}&\leq
\|\sum_{s}a^{3}(t-s)P_{c}Q_{0}\psi\|_{L_{t}^{2}L_{x}^{\infty}}\\
&\les (\sum_{s}\|a^{2}(t-s)P_{c}Q_{
0}\psi\|_{L_{t}^{2}L_{x}^{\infty}}^{2})^{\frac{1}{2}}\lesssim
(\sum_{s}\|a(t-s)P_{c}Q_{0}\psi\|_{L_{t}^{4}L_{x}^{\infty}}^{2})^{\frac{1}{2}}\\
\end{split}\end{equation}
Now one notes that the function $a(t-s)P_{c}Q_{0}\psi$ satisfies
almost the same assumptions about modulation ($\sim 1 $) and
frequency localization as $P_{c}Q_{0}\psi$. Therefore, we can apply
Lemma~\ref{lem:Strich} to estimate
\begin{equation}\nonumber
\|a(t-s)P_{c}Q_{0}\psi\|_{L_{t}^{4}L_{x}^{\infty}}\lesssim
2^{\frac{3k}{4}}2^{\frac{l}{2}}\|P_{c}\psi\|_{\dot{X}_{k}^{0,\frac{1}{2},\infty}}
\end{equation}
Thus
\begin{equation}\nonumber\begin{split}
\|P_{c}Q_{0}\psi\|_{L_{t}^{2}L_{x}^{\infty}}&\les 2^{\frac{3k}{4}}2^{\frac{l}{2}}(\sum_{s\in\Z}\|a(t-s)Q_{0}P_{c}
\psi\|_{\dot{X}_{k}^{0,\frac{1}{2},\infty}}^{2})^{\frac{1}{2}}\\
&\les
2^{\frac{3k}{4}}2^{\frac{l}{2}}\|P_{c}Q_{0}\psi\|_{L_{t}^{2}L_{x}^{2}}\\
\end{split}\end{equation}
The lemma follows via Plancherel's theorem.
\end{proof}

The previous proof was based on the following small-scale version of
the usual $L^4_t L^\infty_x$-Strichartz estimate. It was obtained by
Klainerman and Tataru~\cite{KlainTat}.

\begin{lemma}
 \label{lem:Strich2} With $\calD_{k,\ell}$ as above,  one has
\begin{equation}
 \label{eq:KlaTat} \Big(\sum_{c\in \calD_{k,\ell}} \|P_c e^{it|\nabla|} f\|_{L^4_t L^\infty_x}^2\Big)^{\frac12} \les
 2^{\frac{\ell}{2}} 2^{\frac{3k}{4}} \| f\|_{L^2_x}
\end{equation}
for any $k$-adapted $f$.
In particular,
\begin{equation}
 \label{eq:KlaTatXsb} \Big(\sum_{c\in \calD_{k,\ell}} \|P_c \phi\|_{L^4_t L^\infty_x}^2\Big)^{\frac12} \les
 2^{\frac{\ell}{2}}  \| \phi\|_{\dot X^{\frac34,\frac12,1}_k}
\end{equation}
for any Schwartz function $\phi$ which is adapted to~$k$.
\end{lemma}

\subsection{The null-frame spaces}

In contrast to sub-critical $\dot H^s(\R^2)$ data with $s>1$,  it is
well-known that $\dot X^{s,b}$ spaces do not suffice in the critical case $s=1$. Following the aforementioned references, we
now develop Tataru's null-frame spaces which will provide sufficient control over the nonlinear interactions
in the wave-map system.
 For fixed\footnote{Henceforth, $\omega$ will always be a unit vector in the plane.}
$\omega\in S^{1}$  define
\begin{equation}\label{eq:nullframe}
\theta_{\omega}^{\pm}:= (1,\pm\omega)/\sqrt{2},\quad t_{\omega} := (t,x)\cdot \theta_{\omega}^+ ,\quad x_{\omega} := (t,x)-
t_{\omega}  \theta_{\omega}^+
\end{equation}
which are the coordinates defined by a generator on the light-cone.
Recall that a plane wave traveling in direction~$-\omega\in S^{n-1}$
is a function of the form $h(x\cdot\omega+t)$ (and $h$ sufficiently
smooth). We write a free wave $\phi$  as a superposition of such plane waves: with $\kappa\subset S^1$
and $P_{k,\kappa}$ the projection to $|\xi|\sim 2^k$ and $\wh{\xi}\in\kappa$ as defined above,
\begin{align}
  P_{k,\kappa}\phi(t,x) &= \int_{[|\xi|\sim 2^{k}]} e^{i(t|\xi|+x\cdot\xi)} \widehat{P_{k,\kappa}{f}}(\xi)\, d\xi\nn\\
  &= \int_\kappa \int_{[r\sim 2^{k}]}   e^{ir(x\cdot\omega+t)}
  \hat{f}(r\omega)\,r\, drd\omega\nn\\
  &= \int_{\kappa} \psi_{k,\omega}(t+x\cdot\omega)\, d\omega,\label{eq:plane}
\end{align}
where
\[
  \psi_{k,\omega}(s) := \int_{[r\sim 2^{k}]} e^{irs}
  \hat{f}(r\omega)r\,dr\nn
\]
The argument of $\psi_{k,\omega}$ in~\eqref{eq:plane} is
$\sqrt{2}\,t_\omega$, whence
\begin{equation}\label{eq:L2est}
\int_\kappa
\|\psi_{k,\omega}\|_{L^2_{t_\omega}L^\infty_{x_\omega}} \, d\omega
\le |\kappa|^{\frac12} 2^{\frac{k}{2}} \|P_{k,\kappa} f\|_2
\end{equation}
We now define the following pair of norms\footnote{The $\dist(\omega,\kappa)^{-1}$ factor in the $\NFA[\kappa]$-norm
arises because of a geometric property of the cone, see the proof of Lemma~\ref{lem:incl_free}.}
\begin{align}
  \|G\|_{\NFA[\kappa]} &:=
  \inf_{\omega\not\in2\kappa} \dist(\omega,\kappa)^{-1}
  \|G\|_{L^1_{t_\omega}L^2_{x_\omega}}\label{eq:NFAdef}\\
  \|\phi\|_{\PWA[\kappa]} &:= \inf_{\omega\in\kappa}
  \|\phi\|_{L^2_{t_\omega}L^\infty_{x_\omega}}\label{eq:PWAdef}
\end{align}
which are well-defined for general Schwartz functions. The notation here
derives from {\em null-frame} and {\em plane wave}, respectively.
The quantities defined in~\eqref{eq:NFAdef} and~\eqref{eq:PWAdef}
are not norms --- in fact, not even pseudo-norms --- because they violate
the triangle inequality due to the infimum. This indicates that we should be
using~\eqref{eq:NFAdef} and~\eqref{eq:PWAdef} to define {\em atomic}
Banach spaces (which is why we appended ``A'' in the norms above).
First, recall from \eqref{eq:plane} that
\[
P_{k_1,\kappa}\phi(t,x) = \int_\kappa
\psi_{k,\omega}(\sqrt{2}\,t_\omega)\,d\omega
\]
Then \eqref{eq:L2est} suggests that we define
\[
\|P_{k_1,\kappa}\phi\|_{\PW[\kappa]} := \int_\kappa
\|\psi_{k,\omega}\|_{L^2_{t_\omega}L^\infty_{x_\omega}}\,d\omega =
\int_\kappa \|\psi_{k,\omega}\|_{\PWA[\kappa]}\,d\omega
\]
In other words, $\PW[\kappa]$ is the completion of the space of all
functions $\phi$ which can be written in the form
\begin{equation}\label{eq:PWsplit}
\phi  = \sum_j \lambda_j \psi_j, \quad \sum_j |\lambda_j|<\infty,
\quad \|\psi_j\|_{\PWA[\kappa]}\le1
\end{equation} where
$\lambda_j\in\C$ and $\psi_j$ are Schwartz functions, say. The norm
of any such~$\phi$  in $\PW[\kappa]$ is then simply the infimum of
$\sum_j|\lambda_j|$ over all representations as
in~\eqref{eq:PWsplit}.
By H\"older's inequality we now obtain the simple but crucial estimate
\[
\|\phi F\|_{\NFA[\kappa]} \le \dist(\kappa,\kappa')^{-1}
\|\phi\|_{\PWA[\kappa']} \|F\|_{\Ltwotx}
\]
provided $\phi$ is a $\PWA[\kappa]$-atom. This suggests that we also
define $\NF[\kappa]$ as the atomic space obtained
from~$\NFA[\kappa]$ as usual:  the atoms of $\NF[\kappa]$ are functions $\phi$ for which there exists
$\omega\not\in 2\kappa$ such that $\|\phi\|_{L^1_{t_\omega} L^2_{x_\omega}}\le \dist(\omega,\kappa)$.
The previous estimate
then implies the bound
\begin{equation}\label{eq:mult1}
\|\phi F\|_{\NF[\kappa]} \le \dist(\kappa,\kappa')^{-1}
\|\phi\|_{\PW[\kappa']} \|F\|_{\Ltwotx}
\end{equation}
The dual space $\NF[\kappa]^*$ is characterized by the norm
\[
 \|\phi\|_{\NF[\kappa]^*} = \sup_{\omega\not\in 2\kappa} \dist(\omega,\kappa)^{-1} \|\phi\|_{L^\infty_{t_\omega} L^2_{x_\omega}}<\infty
\]
We now turn to defining the spaces which hold the wave maps.

\begin{defi}\label{def:Sk}
Let $\phi$ be a Schwarz function with
$\supp(\hat{\phi})\subset\{\xi\in\R^2\::\: |\xi|\sim 2^k\}$.
Henceforth, we shall call such a $\phi$ {\em adapted to}~$k$. Define
\begin{align}\label{eq:Skkappa}
  \|\phi\|_{S[k,\kappa]}&:=  \|\phi\|_{L_t^\infty
  L^2_x} +  |\kappa|^{-\frac12}2^{-\frac{k}{2}}\|\phi\|_{\PW[\kappa]} +
  \|\phi\|_{\NF^*[\kappa]}\\
  \|\phi\|_{S[k]} &:= \|\phi\|_{\enr} + \|Q_{\le k+2}\phi\|_{\dot
  X^{0,\frac12,\infty}} + \|Q_{\ge k}\phi\|_{\dot
  X^{-\frac12+\eps,1-\eps,2}} \label{eq:Sk1} \\
&\quad  +  \sup_{j\in\Z} \sup_{\ell\le0}\; 2^{-(\frac12-\eps)\ell}
2^{-\frac{3k}{4}} \Big(\sum_{c\in \calD_{k,\ell}} \|Q_{<j} P_c
\phi\|_{L^4_t L^\infty_x}^2\Big)^{\frac12}
\label{eq:Sk2} \\
  &\quad +  \sup_{\pm}\sup_{\ell\le-100}\;\sup_{\ell\le m\le 0}
  \Big(  \sum_{\kappa\in\caps_\ell} \sum_{R\in\calR_{k,\pm \kappa,m}}
  \|P_{R}
  Q_{\le k+2\ell}^\pm\;
  \phi\|_{S[k,\kappa]}^2 \Big)^{\frac12}\label{eq:squarefunc}
\end{align}
Here $P_\kappa$ and $P_R$ are as above, and $\eps>0$ is a small number ($\eps=\frac{1}{10}$ is sufficient).
\end{defi}

The factors $|\kappa|^{-\frac12}2^{-\frac{k}{2}}$
in~\eqref{eq:Skkappa} are from~\eqref{eq:L2est}. By inspection, the
norm of $S[k,\kappa]$ is translation invariant, and
\begin{equation}\label{eq:Skkappainfty}
\| f\phi\|_{S[k,\kappa]} \le \|f\|_{L^\infty_{tx}}
\|\phi\|_{S[k,\kappa]}
\end{equation}
One has the following scaling property:
\begin{equation}\label{eq:Sk_scale}
\| \phi\|_{S[k]} = \lambda \| \phi(\lambda\cdot)\|_{S[k+m]} ,\qquad \lambda = 2^m,\;m\in\Z
\end{equation}
It will be technically convenient to allow noninteger~$k$ in Definition~\ref{def:Sk}. The only
change required for this purpose is to allow $j ,\ell,m\in\R$  in~\eqref{eq:Sk2} and~\eqref{eq:squarefunc}.
In that case one has
\begin{equation}\label{eq:Sk_scale'}
\| \phi\|_{S[k]} = \lambda \| \phi(\lambda\cdot)\|_{S[k+\log_2\lambda]}  \qquad \forall\lambda>0
\end{equation}
Later we will need to address the question whether $\|P_k\phi\|_{S[k+h]}$ is continuous in~$h$ near $h=0$ for a fixed Schwartz function~$\phi$.
Henceforth, we shall use the operator $I:=\sum_{k\in\Z} P_k Q_{\le k}$ and $I^c:=1-I$ (we will also use $Q_{\le k+C}$ instead of $Q_{\le k}$).
Moreover, we refer to functions which belong to the range
of~$I$ as ``hyperbolic'' and to those in the range of $I^c$ as ``elliptic''. Since
\[
\|Q_{\ge k}P_k\phi\|_{\dot X_k^{0,\frac12,1}} \les \|Q_{\ge
k}P_k\phi\|_{\dot X_k^{-\frac+\eps,1-\eps,2}}
\]
one concludes that the energy norm $\ener$ in~\eqref{eq:Sk1} as well as the Strichartz norm of~\eqref{eq:Sk2} are controlled by
the final norm of~\eqref{eq:Sk1} for the case of elliptic functions (for the Strichartz norm use Lemma~\ref{lem:Strich2}).

\noindent We first verify that
temporally truncated free waves lie in these spaces (with an
imbedding constant that does not depend on the length of the
truncation interval).

\begin{lemma}
  \label{lem:incl_free} Let $\kappa\subset S^{1}$ be arbitrary.
  Then
  \begin{equation}\label{eq:Skkappaimbed}\| \phi\|_{S[k,\kappa]}\les \|P_{k,\kappa} \phi\|_{\dot X^{0,\frac12,1}}
  \end{equation} as well as
  \begin{equation}\label{eq:Skimbed}
\| \phi\|_{S[k]} \les \|P_k Q_{\le k} \phi\|_{\dot X^{0,\frac12,1}}
  \end{equation}
  In particular, if $f$ is adapted to~$k$, then
  \begin{equation}\label{eq:free-trunc}
\| \chi(t/T) e^{it\sqrt{-\Delta}} f\|_{S[k]} \le C \|f\|_{L^2}
  \end{equation}
  with a constant that depends on the Schwartz function $\chi$ but not on~$T\ge 2^{-k}$.
\end{lemma}
\begin{proof} We assume that $\phi$ is an $\dot
X^{0,\frac12,1}$-atom with $P_{0,\kappa}\phi=\phi$. Then from
Plancherel's theorem and Minkowski's and H\"older's inequalities,
\begin{align}
  \|\phi\|_{\enr} &\les \|\hat{\phi}\|_{L^2_\xi L^1_\tau} \les
  2^{\frac{j}{2}}\|\phi\|_{\Ltwo}\nn \\
  \|\phi\|_{L^2_{t_\omega} L^\infty_{x_\omega}} &\les
  (2^j|\kappa|)^{\frac12} \|\hat{\phi}\|_{\Felltwo} =
  (2^j|\kappa|)^{\frac12} \|\phi\|_{\Ltwo} \nn \\
  \|\phi\|_{\enromega} &\les \frac{2^{\frac{j}{2}}}{\dist(\kappa,\omega)}\|\phi\|_{\Felltwo} \label{eq:NFsternfree}
\end{align}
In the final estimate~\eqref{eq:NFsternfree} we used that
$\angle(\ell_\omega,T_{\omega'})\sim \angle(\omega,\omega')^2$ where
 $\ell_\omega$ is the line oriented along the generator parallel to
$(1,\omega)$ and  $T_{\omega'}$ is the tangent plane to the cone
which touches the cone along the generator $\ell_{\omega'}$. To
establish \eqref{eq:Skimbed} we begin with
  \[
\sup_{\ell\le-100}\Big( \sum_{\kappa\in\caps_\ell}
\sum_{R\in\calR_{k,\pm \kappa,\lambda}} \|P_{R}
  Q_{\le 2\ell}^\pm
  \phi\|_{ \dot
  X_0^{0,\frac12,1}  }^2 \Big)^{\frac12} \les \|\phi\|_{\dot
  X_0^{0,\frac12,1}}
  \]
which is obvious from orthogonality of the $P_{0,\pm\kappa}$. In
view of~\eqref{eq:Skkappaimbed}, this bound yields the square
function in~\eqref{eq:squarefunc}. The energy is controlled via the
imbedding $\|\phi\|_{L^\infty L^2} \les \|\phi\|_{\dot
X^{0,\frac12,1}}$, whereas the Strichartz component of~$S[k]$ is controlled
by Lemma~\ref{lem:Strich2}.

Finally, the statement
concerning the free wave reduces to the case $k=0$ for which we need
to verify the bound
\begin{align*}
&\sum_{j\in\Z} 2^{\frac{j}{2}} \| T
\hat{\chi}(T|\tau\pm|\xi\|)m(2^{-j}|\tau\pm|\xi||)
\hat{f}(\xi)\|_{L^2_\tau L^2_\xi} +\\
&  +\Big(\sum_{j\in\Z} 2^{2j} \| T
\hat{\chi}(T|\tau\pm|\xi||)m(2^{-j}|\tau\pm|\xi||)
\hat{f}(\xi)\|_{L^2_\tau L^2_\xi}^2\Big)^{\frac12}   \les \|f\|_2
\end{align*}
which are both clear provided $T\ge1$ due to the rapid decay of
$\hat{\chi}$.
\end{proof}

Naturally, $S[k]$ contains more general functions than just free
waves. One way of obtaining such functions is to take
$\phi=\Box^{-1} F$, in other words from the Duhamel formula. We will
study this in much greater generality in the context of the energy
estimate below, but for now we take $F$ to be a Schwartz function.

\begin{exse}
  \label{ex:not_free}
The bounded
function $\phi$ defined via its Fourier transform
\[
\wh{\phi}(\tau,\xi) = \chi_1(\xi) \chi_2(|\xi|-|\tau|)
(|\xi|-|\tau|)^{-1}
\]
belongs to $S[0]$ but is not a truncated free wave. Here $\chi_1$ is a smooth cut-off to
$|\xi|\sim1$, and $\chi_2(u)$ is a smooth cut-off to $|u|<1/10$. We
leave it to the reader to construct other  functions which lie in
$S[0]$ and which are not (truncated) free waves.
\end{exse}

The following basic estimates will be used repeatedly:

\begin{itemize}
  \item if $\phi$ is adapted to $k$, then
\begin{align}
 \|Q_j \phi\|_{L^2L^2} &\les \min(2^{-(j-k)(\frac12-\eps)},1) 2^{-\frac{j}{2}} \|\phi\|_{S[k]} \label{eq:Sbd1}\\
  \|Q_j \phi\|_{L^2L^\infty} &\les 2^{\frac{j-k}{4}\wedge 0}  \min(2^{-(j-k)(\frac12-\eps)},1) 2^{-\frac{j}{2}} \|\phi\|_{S[k]} \label{eq:Sbd2}
\end{align}
This follows from the $\dot X^{s,b,q}$ components of the
$S[k]$-norms, as well as the improved Bernstein's inequality of Lemma~\ref{lem:KBern}.
\item The duality between $\NF[\kappa]$ and $\NF^*[\kappa]$
implies
\begin{equation}
  |\la \phi,F \ra| \les \|\phi\|_{S[k,\kappa]} \|F\|_{\NF[\kappa]}\label{eq:SNFdual}
\end{equation}
\end{itemize}

In what follows, $ \Theta:=\sign(\tau)\wh{\xi}$, and for any
$\omega\in S^1$,  $\Pi_\omega$ denotes the orthogonal projection
onto $\NP(\omega):=\theta_\omega^\perp$ (the {\em null-plane}
of~$\omega$).

\begin{lemma}
  \label{lem:Piomega}
The projection $\Pi_\omega$ satisfies the following properties:
  \begin{itemize}
  \item
Let $\calF\subset \calC_\ell$ be a collection of disjoint caps.
Suppose that $\omega\in S^1$ satisfies
$\dist(\omega,\kappa)\in[\alpha,2\alpha]$ for any $\kappa\in\calF$
where $\alpha>2^\ell$ is arbitrary but fixed. Define\footnote{An
important detail here is that these dimensions deviate from the
usual wave-packets of dimension $1\times 2^{\ell}\times 2^{2\ell}$.}
   \begin{equation}\label{eq:Thetakappa}
   \calT_{\kappa,\alpha} := \big\{(\tau,\xi)\::\: |\xi|\sim1,\;\Theta\in \kappa,\;
  ||\xi|-|\tau||\les \alpha 2^{\ell} \big\}
  \end{equation}
  Then $\{\Pi_\omega(\calT_{\kappa,\alpha})\}_{\kappa\in\calF}\subset \NP(\omega)$ have finite
  overlap, i.e.,
  \[
\sum_{\kappa\in\calF} \chi_{\Pi_\omega(\calT_{\kappa,\alpha})} \le C
  \]
  where $C$ is some absolute constant.
\item
  Let
\[
 \calS:=\Big\{ (\pm|\xi|,\xi)\::\: \xi\in\R^2,\;  \wh{\xi}\in \pm\kappa\Big\}
\]
be a sector on the light-cone where $\kappa\subset S^1$ is any cap.
Furthermore, let $\omega\not\in 2\kappa$ and
$\wt\calS:=\Pi_\omega(\calS)$. Then on~$\wt\calS$ the Jacobian
$\frac{\del\xi}{\del\xi_\omega}$ satisfies
\begin{equation}
 \label{eq:xiJac}  \Big| \frac{\del\xi}{\del\xi_\omega} \Big|\sim d(\omega,\kappa)^{-2}
\end{equation}
The same holds on $\Pi_\omega(\calS_a)$ where
\[
\calS_a := \Big\{ (\pm|\xi|+a,\xi)\::\: \xi\in\R^2,\;|\xi|\sim 1,\;
\wh{\xi}\in \pm\kappa\Big\}
\]
provided $a$ is fixed with $|a|\le |\kappa| d(\omega,\kappa)$.
\end{itemize}
\end{lemma}
\begin{proof}
Denote
\[
\calS_{\kappa}:=\{ s(1,\omega') + \rho(1,-\omega')\::\:
\omega'\in\kappa,\; 1\le s\le 2,\; |\rho|< h \}
\]
where $h$ will be determined.  Then
\[ \{\Pi_\omega(\calS_{\kappa})\}_{\kappa\in\calF}
=\{s\vecv + \rho\vecw \::\: \omega'\in\kappa,\; 1\le s\le 2,\;
|\rho|< h \}
\]
\begin{figure*}[ht]
\label{fig:1}
\begin{center}
\centerline{\hbox{\vbox{ \epsfxsize= 11.0 truecm \epsfysize=6.5
truecm \epsfbox{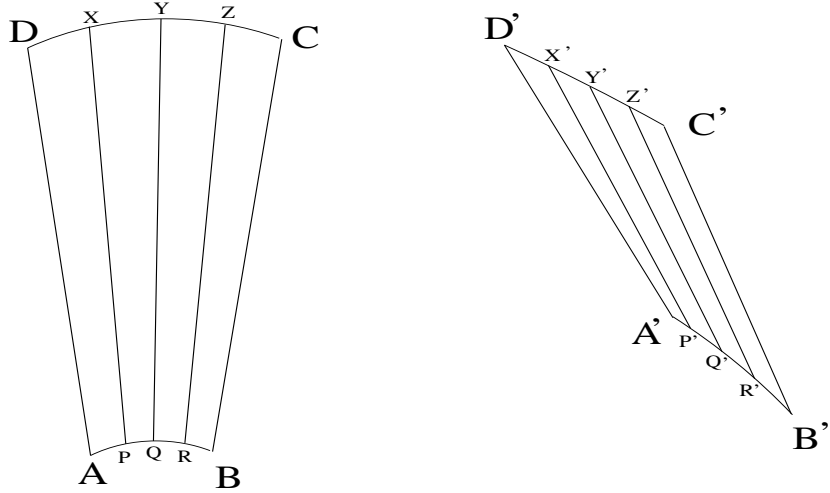}}}} \caption{The projected sectors}
\end{center}
\end{figure*}
where
\[
\vecv:= \vecv(\omega,\omega')=(1,\omega')-\lambda(1,\omega),\quad
\vecw:= (1,-\omega')-\mu (1,\omega)
\]
with $\lambda=\frac12(1+\omega\cdot\omega')$,
$\mu=\frac12(1-\omega\cdot\omega')$. Recall that
$\dist(\omega,\kappa)\sim \dist(\omega,\kappa')=:\alpha$ where
$\kappa\in\calF$ is arbitrary. Moreover, $\diam(\kappa)\sim
2^\ell=:\beta$. One checks that
\[
|\vecv|=\sqrt{2(1-\lambda^2)}\sim \sqrt{\mu} \sim \alpha
\]
Furthermore, $\del\vecv:= (0,{\omega'}^\perp)-\frac12 \omega\cdot
{\omega'}^\perp (1,\omega)$ denotes the derivative
$\del_{\theta'}\vecv$ where we have written $\omega'=e^{i\theta'}$.
Then $|\del\vecv|\sim1$ and
\[
\vecv\wedge \del\vecv = (\mu,\omega'-\omega+\mu\omega)\wedge (0,
{\omega'}^\perp) - \frac12 \omega\cdot {\omega'}^\perp
(1,\omega')\wedge (1,\omega)
\]
satisfies $|\vecv\wedge \del\vecv|\sim \alpha^2$. In conjunction
with $|\vecv|\sim\alpha$ this implies that
$|\sangle(\vecv,\del\vecv)|\sim\alpha$. Since
\[
|\vecv(\omega,\omega')-\vecv(\omega,\omega'')|\gtrsim
|\omega'-\omega''| \quad \forall\; \omega',\omega''\in\kappa'
\]
it follows that
\begin{equation}\label{eq:linedist}
\dist(\sigma(\omega,\omega'),\sigma(\omega,\omega''))\gtrsim
\alpha|\omega'-\omega''|
\end{equation}
where
\[
\sigma(\omega,\omega'):=\{s\vecv(\omega,\omega')\::\: 1\le s\le 2\}
\]
Therefore, one needs to take $h=\alpha\beta$ to insure the property
of finite overlap of the projections. This is optimal, since one can
check that $\vecv$ and $\vecw$ always satisfy
$|\cos(\sangle(\vecv,\vecw)|\le \frac12$. In Figure~2 the left-hand
side depicts four sectors as they would appear on the light-cone,
whereas the right-hand side is the projected configuration
in~$\NP(\omega)$ with $A':=\Pi_\omega(A)$ etc. Note that the
segments $\overline{A'B'}$ as well as $\overline{A'P'}$,
$\overline{P'Q'}$, $\overline{Q'R'}$, $\overline{R'B'}$ have lengths
comparable to the corresponding ones on the left, i.e.,
$\overline{AB}$ etc., whereas the lengths of $\overline{A'D'}$,
$\overline{B'C'}$ are those of $\overline{AD}$ and $\overline{BC}$
contracted by the factor $\alpha$. Finally, we have shown that
$\sangle(A'B'C')\sim \alpha$ (and similarly for the angles at the
points $P'$, $Q'$, $R'$) so that the height of the parallelogram
$A'P'X'D'$ is proportional to $\alpha$ times the length of
$\overline{A'P'}$, see~\eqref{eq:linedist}.

\noindent The second statement of the lemma follows from the
consideration of the preceding paragraph.
\end{proof}

As a consequence of Lemma~\ref{lem:Piomega}, we now show that the
square-function in \eqref{eq:squarefunc} can always be refined in
terms of the angle.

\begin{lemma}
  \label{lem:square_func} Let $\calF\subset \calC_\ell$ be a collection of disjoint caps  and let $\kappa'\in\calC_{\ell'}$ be a cap with
   $\bigcup_{\calF} \kappa\subset \kappa'$. Suppose further that for every $\kappa\in\calF$
  there is a Schwartz function
  $\phi_\kappa$ adapted to $k\in\Z$ and which is supported
  on
   \begin{equation}\nn
   \calT_{\kappa,k} := \big\{ \Theta:=\sign(\tau)\wh{\xi}\in \kappa,\;
  ||\xi|-|\tau||\les 2^{\ell+\ell'+k} \big\}
  \end{equation}
with some $k\in\Z$.
  Then
  \[
\big\| \sum_{\kappa\in\calF} \phi_{\kappa}\big\|_{S[k,\kappa']}\le C\Big(\sum_{\kappa\in \calF}  \|\phi_{\kappa}\|^2_{S[k,\kappa]}\Big)^{\frac12}
  \]
with some absolute constant $C$.
\end{lemma}
\begin{proof} First, one may take $k=0$ and $\tau>0$ (the latter by conjugation symmetry).
The $L^\infty L^2$-component of~\eqref{eq:Skkappa} satisfies the
required property due to orthogonality, whereas
the~$\PW[\kappa]$-component is reduced to Cauchy-Schwarz (via the
$|\kappa|^{-\frac12}$-factor). For the final $\NF^*[\kappa]$ (i.e.,
$L^\infty_{t_\omega} L^2_{x_\omega}$)-component one exploits
orthogonality relative to $x_\omega$ via
Lemma~\ref{lem:square_func}. Here $\omega\in S^1\setminus(2\kappa')$
is arbitrary but fixed.
\end{proof}

Later we will prove bi- and trilinear estimates involving $S$
and~$N$ space. The following bilinear bounds will be a basic
ingredient in that context.

\begin{lemma}
  \label{lem:bilin} One has the estimates
  \begin{align}
    \| \phi F\|_{\NF[\kappa]} &\les
    \frac{|\kappa'|^{\frac12} 2^{\frac{k'}{2}} }{\dist(\kappa,\kappa')}
    \|\phi\|_{S[k',\kappa']} \|F\|_{\Ltwotx} \label{eq:bilin1}\\
    \| \phi \psi\|_{\Ltwotx} &\les \frac{|\kappa|^{\frac12} 2^{\frac{k}{2}} }{ \dist(\kappa,\kappa')}
    \|\phi\|_{S[k,\kappa]} \|\psi\|_{S[k',\kappa']} \label{eq:bilin2}
  \end{align}
  For the final two bounds we require that $2\kappa\cap
  2\kappa'=\emptyset$.
\end{lemma}
\begin{proof}
  The second one follows
  from the definition of the spaces, whereas~\eqref{eq:bilin2}
  follows from~\eqref{eq:bilin1} and the duality bound~\eqref{eq:SNFdual}.
\end{proof}

Note that both of these estimates have a dispersive character, as
they involve space-time integrals. By applying ideas from the energy
estimate,  we will improve on~\eqref{eq:bilin2} in the high-high
case, see  Lemma~\ref{lem:bilin3}. Next, we define the spaces which
will hold the nonlinearities.

\begin{defi}
\label{def:Nk} $N[k]$ is generated by the following four types of atoms: with $F$ being $k$-admissible, either
\begin{itemize}
  \item $ \|F\|_{L^1_tL^2_x}\le 2^{k}$
  \item $\hat{F}$ is supported on $||\xi|-|\tau||\sim 2^j\le 2^k$ and $\|F\|_{\dot X_k^{-1,-\frac12,1}}\le1$
  \item $F=Q_{\ge k}F$,  $\|F\|_{\dot
  X_k^{-\frac12+\eps,-1-\eps,2}}\le1$ where $\eps>0$ is as in the $S[k]$ spaces
  \item $F$ is the sum of wave packets $F_\kappa$: there exists $\ell\le-100$ such that  $F=\sum_{\kappa\in\calC_\ell} F_\kappa$
  with all  $\supp(\wh{F_\kappa})$ supported on either $\tau>0$ or $\tau<0$, with $\wh{F_\kappa}$
  supported on $|\xi|\sim 2^k$, $||\xi|-|\tau||\le C^{-1} 2^{k+2\ell}$,  $\Theta:=\sign(\tau)\wh{\xi}\in \kappa$  and to that the bound
  \[
 \Big(\sum_{\kappa}  \| F_\kappa \|^2_{\NF[\kappa]}\Big)^{\frac12} \le
 2^{k}
  \]
  holds.
\end{itemize}
We refer to these types as {\em energy, $\dot X^{s,b,q}$, and
wave-packet atoms}, respectively.
\end{defi}

In what follows, we refer to functions $\phi$ adapted to some
$k\in\Z$ as ``elliptic'' iff $P_kQ_{\ge k}\phi=\phi$, whereas those
satisfying $P_kQ_{\le k}\phi=\phi$ as ``hyperbolic''. This
terminology has to do with the behavior of the wave operator~$\Box$
in these respective regimes. We now record a  fundamental duality
property of $N[k]$.

\begin{lemma}
  \label{lem:Nk} For any $\phi\in S[k]$ and $F\in N[k]$ with $F=P_{k} Q_{\le k}F$
  \begin{align}
  |\la \phi,F\ra| &\les 2^{k} \|\phi\|_{S[k]}\|F\|_{N[k]} \label{eq:SNdual}\\
\|F\|_{\dot X_k^{-1,-\frac12,\infty}} &\les \|F\|_{N[k]}\les
\|F\|_{\dot X_k^{-1,-\frac12,1}} \label{eq:Xsb_dom}
  \end{align}
\end{lemma}
\begin{proof}
The duality relation~\eqref{eq:SNdual} is proved by taking $F$ to be an atom; for the wave-packet atom use~\eqref{eq:SNFdual}.
  By definition of $N[k]$, one has $\|F\|_{N[k]}\les \|F\|_{\dot X_k^{-1,-\frac12,1}}$.
  For the left-hand bound in~\eqref{eq:Xsb_dom}
  use~\eqref{eq:Skkappaimbed} and~\eqref{eq:SNdual}.
\end{proof}

As an application of the geometric considerations of
Lemma~\ref{lem:Piomega} we now show that {\em refining} a
wave-packet atom yields another wave-packet atom.

\begin{lemma}\label{lem:Nsquare}
 Let $F=\sum_{\kappa\in\calC_\ell} F_\kappa$ be a wave-packet atom as in Definition~\ref{def:Nk}.
Then
\begin{equation}
 \label{eq:NFsquare}
\sup_{\ell'\le\ell}\; \sup_{\ell'\le j\le 0} \Big( \sum_\kappa
\sum_{\substack{\kappa'\in\calC_{\ell'}\\\kappa'\subset\kappa}}
\sum_{R\in\calR_{k,\kappa',j}} \|P_R P_{\kappa'} Q_{<\ell+\ell'+k}
F_\kappa\|^2_{\NF[\kappa']} \Big)^{\frac12} \le C\, 2^k
\end{equation}
with some absolute constant~$C$.
\end{lemma}
\begin{proof} By scaling invariance, we can set $k=0$. Moreover, fix
$\ell'\le\ell$ and $\ell'\le j\le0$.  Choose
$\omega'=\omega(\kappa')\in S^1\setminus (2\kappa')$ for
each~$\kappa'$ which attain the respective $\NF[\kappa']$ norm. Then
one has
\begin{align}
  &  \sum_\kappa
\sum_{\substack{\kappa'\in\calC_{\ell'}\\\kappa'\subset\kappa}}\sum_{R\in\calR_{0,\kappa',j}} \|P_R P_{\kappa'}
Q_{<\ell+\ell'} F_\kappa\|^2_{\NF[\kappa']} \nn \\
&\les  \sum_\kappa
\sum_{\substack{\kappa'\in\calC_{\ell'}\\\kappa'\subset\kappa}}\sum_{R\in\calR_{0,\kappa',j}}  d(\omega',\kappa')^{-2}
\|P_R P_{\kappa'} Q_{<\ell+\ell'} F_\kappa\|^2_{L^1_{t_\omega} L^2_{x_\omega}}  \nn \\
&\les  \sum_\kappa \inf_{\omega\in S^1\setminus(2\kappa)}
 d(\omega,\kappa)^{-2} \Big\| \Big( \sum_{\substack{\kappa'\in\calC_{\ell'}\\\kappa'\subset\kappa}}\sum_{R\in\calR_{0,\kappa',j}}
  \|P_R P_{\kappa'} Q_{<\ell+\ell'} F_\kappa\|^2_{
  L^2_{x_\omega}} \Big)^{\frac12}
  \Big\|_{L^1_{t_\omega}}^2 \label{eq:ell2L1}\\
&\les   \sum_{\kappa}  \inf_{\omega\in S^1\setminus(2\kappa)}
 d(\omega,\kappa)^{-2}
\|F_\kappa\|^2_{L^1_{t_\omega} L^2_{x_\omega}}
\label{eq:xomega_orth}
\end{align}
To pass to \eqref{eq:ell2L1} we used  the inclusion $\ell^2
(L^1_{t_\omega}) \supset L^1_{t_\omega}(\ell^2)$, whereas
orthogonality  implies~\eqref{eq:xomega_orth}. Indeed, first note
that
\[
\bigcup_{t_\omega\in\R} \supp \big( [P_R P_{\kappa'} Q_{<\ell+\ell'}
F_\kappa)(t_\omega,\cdot)]^{\wedge}\big) \subset
\Pi_\omega\big(\supp( \calF[P_R P_{\kappa'} Q_{<\ell+\ell'}
F_\kappa] ) \big)
\]
where the Fourier transform on the left-hand side is in~$x_\omega$
and on the right-hand side in $(t_\omega,x_\omega)$. Second, the
sets on the right-hand side enjoy a finite overlap property by
Lemma~\ref{lem:Piomega}.
\end{proof}

In what follows, we will often need to split a wave $\phi$ into $\phi^{+}+\phi^{-}$ where
\[
\phi^{+}:=
\big(\chi_{[\tau\ge0]}\hat\phi(\cdot,\tau)\big)^{\vee},\qquad
\phi^{-}:= \big(\chi_{[\tau<0]}\hat\phi(\cdot,\tau)\big)^{\vee}
\]
The question arises whether the spaces $S[k]$ and $N[k]$ are preserved under these operations.

\begin{lemma}
  \label{lem:pm} For any Schwartz function $\phi$ which is adapted to $k$,
  \[
\|\phi^{\pm}\|_{S[k]} \le C\|\phi\|_{S[k]}, \qquad \|F^{\pm}\|_{N[k]} \le C\|F\|_{N[k]}
  \]
  with some absolute constant~$C$.
\end{lemma}
\begin{proof} We set $k=0$ and assume that $\phi$ is adapted to $k=0$. Let $\chi_0$ be a bump function
on the line with $\chi_0(\tau)=1$ on $\tau\ge-C^{-1}$ and $\chi_0(\tau)=0$ if $\tau\le -2C^{-1}$ where $C>1$ is
some large constant. Then
\[
\widehat{\phi^{+}}(\tau,\xi)= \chi_0(\tau-|\xi|)\chi_{[\tau\ge0]}\wh\phi(\tau,\xi)
 + (1-\chi_0)(\tau-|\xi|)\chi_{[\tau\ge0]}\wh\phi(\tau,\xi)
\]
Denote the two functions on the right-hand side by $\phi^{(+,1)}$ and $
\phi^{(+,2)}$, respectively.
Then \begin{equation}\label{eq:+1}\phi^{(+,1)} = \phi\ast\mu\end{equation}
 where $\mu$ is a measure of bounded mass. Therefore,
\[
\|\phi^{(+,1)} \|_{S[0]} \le C\|\phi\|_{S[0]}
\]
Next,
\[
\|\phi^{(+,2)} \|_{S[0]} \le C\|\phi^{(+,2)}\|_{\dot X_0^{0,\frac12,1}}\le C\|\phi\|_{S[0]}
\]
where we used the Plancherel theorem in the final step.

For $N[k]$ it will suffice to check the case of $L^1_t L^2_x$-atoms. For these, we write
\[
F^{+}=F^{(+,1)}+F^{(+,2)}
\]
as above. The first term here is fine from~\eqref{eq:+1}, whereas the second is placed in~$\Ltwotx$ and bounded by means
of~\eqref{eq:Xsb_dom}.
\end{proof}

Another piece of terminology used by Tao is the following:
\begin{defi}
We shall say
that a family $\{m_\alpha\}_{\alpha}$ is {\em disposable}, if $T_\alpha f:= (m_\alpha \hat{f})^{\vee} = f\ast \mu_\alpha$ where
$\mu_\alpha$ are measures with uniformly bounded mass:
\[
\sup_\alpha \|\mu_\alpha\|\le C<\infty
\]
with some universal constant $C$.
\end{defi}

Clearly, disposable multipliers give rise to bounded operators on
any translation invariant Banach space. Thus, if $X$ is a Banach
space of functions on $\R^{n+1}$ with the property that for all
$f\in X$ one has
\[
\|f(\cdot-y)\|_{X}=\|f\|_X \qquad \forall\; y\in \R^{n+1}
\]
then $\sup_{\alpha}\|T_\alpha f\|_X\le C\|f\|_X$.
The following observation will be a useful device for removing frequency cut-offs.

\begin{lemma}
  \label{lem:dispose} The families  \[\big\{ P_{k,\kappa}  \big\}_{k,\kappa},\quad \big\{P_k Q_j\big\}_{j\ge k},\quad
  \big\{ P_kQ_{<j}\big\}_{j\ge k},\quad
  \big\{P_k Q_{<k-C}^{\pm}\big\}_{k}\]
  are disposable. In the first family $\kappa$ is any cap, whereas in the last family $C>0$ has to be chosen such that
  the support of the multiplier associated with $P_k Q_{<k-C}$ does not intersect $\tau=0$.
    In addition,
  \[
\big\{ P_{k,\kappa} Q_{<k+2\ell}\big\}
  \] is disposable
  where $k\in\Z$ and  $\kappa$ is any cap with
  $\diam(\kappa)\sim 2^\ell$ with $\ell\le-100$ arbitrary.
\end{lemma}
\begin{proof}
  Without loss of generality one may take $k=0$. Then these
  statements reduce to simple exercises in harmonic analysis.
\end{proof}

The following fact will serve as a substitute for the previous
problem in a non-disposable context.

\begin{lemma}
  \label{lem:QLp}
   $Q_j$, $Q_{<j}$ are bounded on $L^pL^2$ for every $1\le p\le\infty$ with a constant independent
  of~$j\in\Z$.
\end{lemma}
\begin{proof}
 The inverse Fourier transform of $Q_{<j}$ with respect to time alone is
\begin{align*}
 & \int e^{it\tau} m_0(2^{-j} (\tau-|\xi|)) \wh{F}(\tau,\xi)  \, d\tau \\
&=  2^{j} \int \wh{m_0}(2^{j} (t-s) ) e^{i|\xi|(t-s)} {F}(s,\hat{\xi})  \, ds
\end{align*}
where ${F}(s,\hat{\xi})$ in the second line denotes the Fourier transform with respect to
the second variable. Consequently,
\[
 \|Q_{<j} F\|_{L^1_t L^2_x} \le \|\wh{m_0}\|_1 \|F\|_{L^1_t L^2_x}
\]
as claimed.
\end{proof}

The previous result, combined with Lemma~\ref{lem:square_func},
implies the following square-function bound.

\begin{cor}
  \label{cor:S_cut} For all $j,k\in\Z$ and all $k$--adapted Schwartz functions $\phi$ one has
  $\|Q_{<j} \phi\|_{S[k]} \le C \|\phi\|_{S[k]}$ with some absolute constant~$C$.
\end{cor}
\begin{proof} We may again take $k=0$.
   The  $L^\infty_t L^2_x$-component of the $S[0]$-norm is covered by
   Lemma~\ref{lem:QLp}. The $\dot X^{s,b,q}$-components are obvious ,
the Strichartz norms as well by construction, and
   the square-function is a consequence of
   Lemma~\ref{lem:square_func}.
\end{proof}

We remark that the analogous statement for~$N[k]$ holds as well, see Corollary~\ref{cor:N_cut} below.
Next, for the sake of completeness we state the full range of Strichartz estimates that follow from~\eqref{eq:Sk2}.

\begin{lemma}
 \label{lem:Strich}  For any $4\le p\le \infty$ and $2\le q\le\infty$ which satisfy $\frac{1}{p}+\frac{1}{2q}\le\frac14$,
\begin{equation}
 \label{eq:Strich} \Big(\sum_{c\in \calD_{k,\ell}} \|  P_c
\phi\|_{L^p_t L^q_x}^2\Big)^{\frac12} \le C 2^{\ell(1-\frac{2}{p}-\frac{2}{q}-\frac{4\eps}{p})} 2^{k(1-\frac{1}{p}-\frac{2}{q})}
\|\psi\|_{S[k]}
\end{equation}
for any $k\in\Z$, $\ell\le0$, and with an absolute constant~$C$.
\end{lemma}
\begin{proof} Assume first that $\frac{1}{p}+\frac{1}{2q}=\frac14$.
 By interpolation, and with $\theta=\frac{4}{p}$,
\begin{align}
 \Big(\sum_{c\in \calD_{0,\ell}} \|  P_c
\phi\|_{L^4_t L^\infty_x}^2\Big)^{\frac12} &\le \Big(\sum_{c\in \calD_{0,\ell}} \|P_c \psi\|_{L^4_t L^\infty_x}^{2\theta}   \|  P_c
\phi\|_{L^\infty_t L^2_x}^{2(1-\theta)} \Big)^{\frac12} \nn \\
&\le  \Big(\sum_{c\in \calD_{0,\ell}} \|P_c \psi\|_{L^4_t L^\infty_x}^{2 }   \Big)^{\frac{\theta}{2}}
\Big(\sum_{c\in \calD_{0,\ell}} \|P_c \psi\|_{L^\infty_t L^2_x}^{2 }   \Big)^{\frac{1-\theta}{2}} \label{eq:zweifach} \\
&\les   2^{\theta(\frac12-\eps)\ell}  \|\psi\|_{S[0]}  \nn
\end{align}
To pass from~\eqref{eq:zweifach} to the last line, one uses~\eqref{eq:Sk2} as well as the energy component of~\eqref{eq:squarefunc}.
For larger $q$, one gains a factor $2^{2\ell(\frac12-\frac{2}{p}-\frac1q)}$ by Bernstein's inequality, and rescaling to frequency $2^k$ yields a factor
of~$2^{k(1-\frac{1}{p}-\frac{2}{q})}$ as claimed.
\end{proof}

Finally, we conclude this section with the following useful fact.

\begin{lemma}
 \label{lem:enersquaresum} Let $\phi$ be adapted to~$0$. Then for any $m_0\le -10$,
\[
 \Big(\sum_{\kappa\in\calC_{m_0} }\|P_{0,\kappa} \phi\|_{\ener}^2\Big)^{\frac12} \les |m_0| \|\phi\|_{S[0]}
\]
\end{lemma}
\begin{proof}
 First,
\begin{align*}
 \Big(\sum_{\kappa\in\calC_{m_0}} \|P_{0,\kappa} Q_{\le 2m_0} \phi\|_{\ener}^2\Big)^{\frac12} &\les
\Big(\sum_{\kappa\in\calC_{m_0}} \|P_{0,\kappa} Q_{\le 2m_0} \phi\|_{S[0,\kappa]}^2\Big)^{\frac12}\les  \|\phi\|_{S[0]}
\end{align*}
by~\eqref{eq:squarefunc}.  Second,
\begin{align*}
 \sum_{2m_0\le \ell\le 0} \Big(\sum_{\kappa\in\calC_{m_0}} \|P_{0,\kappa} Q_{\ell} \phi\|_{\ener}^2\Big)^{\frac12} &\les
\sum_{2m_0\le \ell\le 0} \Big(\sum_{\kappa\in\calC_{m_0}} \|P_{0,\kappa} Q_{\ell} \phi\|_{\dot X_0^{0,\frac12,\infty}}^2\Big)^{\frac12} \\
&\les  |m_0| \|\phi\|_{S[0]}
\end{align*}
And third,
\[
 \| P_{0,\kappa} Q_{\ge 0} \phi\|_{\ener} \les \| P_{0,\kappa} Q_{\ge 0} \phi\|_{\dot X_0^{0,\frac12,1}} \les
 \| P_{0,\kappa} Q_{\ge 0} \phi\|_{\dot X_0^{0,1-\eps,2}}
\]
whence
\begin{align*}
 \Big(\sum_{\kappa\in\calC_{m_0}} \|P_{0,\kappa} Q_{\ge 0} \phi\|_{\ener}^2\Big)^{\frac12} &\les
\Big(\sum_{\kappa\in\calC_{m_0}} \|P_{0,\kappa} Q_{\ge 0} \phi\|_{ \dot X_0^{0,1-\eps,2}}^2\Big)^{\frac12}  \\
&\les  \sum_{\ell\ge0}  \Big(\sum_{\kappa\in\calC_{m_0}} \|P_{0,\kappa} Q_{\ell} \phi\|_{ \dot X_0^{0,1-\eps,\infty}}^2\Big)^{\frac12} \\
&\les \sum_{\ell\ge0} \|P_{0,\kappa} Q_{\ell} \phi\|_{ \dot
X_0^{0,1-\eps,\infty}}\le \|P_{0,\kappa} Q_{\ell} \phi\|_{ \dot
X_0^{0,1-\eps,2}} \les \|\phi\|_{S[0]}
\end{align*}
as claimed.
\end{proof}

The central problems concerning
the $S[k]$ and $N[k]$ spaces are how to obtain an energy estimate and how to control the
trilinear nonlinearities appearing in the gauged wave-map system.   We begin with the energy
estimate, and then develop bilinear bounds which are preliminary to the central trilinear
bounds.

\subsection{The energy estimate}

\label{subsec:energy}

The purpose of this section is to prove the  energy estimate in the
context of the $S[k]$ and $N[k]$ spaces, see
Proposition~\ref{prop:energy} below. First, we require some
technical lemmas. The first two of these lemmas will arise in the
Duhamel integral.

\begin{lemma}
  \label{lem:Ncut} For any $F$ which is $k$-adapted and satisfies $F=Q_{\le k+C}F$,
  \begin{equation}
\label{eq:rplus_remov} \|\chi_{\R^+}F\|_{N[k]}\les \|F\|_{N[k]}
\end{equation}
where $\chi_{\R^+}$ acts only in time.
\end{lemma}
\begin{proof} We may assume that $k=0$.
  This is clear if $F$ is an energy atom. Next, we consider the $\dot
X^{0,-\frac12,1}$-atoms. First, let $\hat{F}$ be supported on
$|\xi|~\sim1, ||\xi|-|\tau||\sim 2^j$ with $\|F\|_{\Ltwotx}\les
2^{j/2}$. Then
\begin{align*}
\|\chi_{\R^+}F\|_{N[0]}  &\les   \| P_{<j}(\chi_{\R^+}) F\|_{N[0]} + \|P_{\ge j} (\chi_{\R^+}) F\|_{N[0]} \\
&\les  2^{-j/2}  \| P_{<j}(\chi_{\R^+}) F\|_{\Ltwotx} + \| P_{\ge j} (\chi_{\R^+}) F\|_{L^1 L^2}  \\
&\les 2^{-j/2}  \| P_{<j}(\chi_{\R^+})\|_{L_t^\infty} \|F\|_{\Ltwotx} + \|P_{\ge j} (\chi_{\R^+})\|_{L_t^2} \|F\|_{L^2 L^2}\\
&\les 2^{-j/2}   \|F\|_{\Ltwotx} \les 1
\end{align*}
Now let $F$ be a wave-packet atom, i.e., for some $\ell\le-100$,
\[
F=\sum_{\kappa\in\calC_\ell} F_\kappa,\quad \supp(\wh{F_\kappa})
\subset \big\{\tau>0,\; |\xi|\sim1, \; ||\xi|-\tau|\sim 2^{2\ell},\;
\Theta\in\kappa \big \}
\]
and $\sum_{\kappa} \|F_\kappa\|_{\NF[\kappa]}^2\le1$. We write, with
$j=2\ell$,
\[
\chi_{\R^+} = P_{<j} (\chi_{\R^+}) +  P_{\ge j}( \chi_{\R^+})
\]
as before. Then $P_{<j} (\chi_{\R^+})$ does not significantly change
the support properties of~$F_\kappa$. Moreover, since $\|P_{<j}
(\chi_{\R^+})\|_\infty\les1$, we see that $P_{<j} (\chi_{\R^+}) F$
is essentially a wave-packet atom. On the other hand, since
$\|F\|_{\Ltwotx}\les 2^{j/2}$ from~\eqref{eq:Xsb_dom} we conclude
that
\begin{equation}
  \label{eq:PjFener}\| P_{\ge j}  (\chi_{\R^+}) F\|_{L^1 L^2} \les  2^{-j/2}
\|F\|_{\Ltwotx}\les 1
\end{equation}
which proves~\eqref{eq:rplus_remov}.
\end{proof}

It is important to note that the previous lemma {\em fails} for
functions in $N[0]$ which are ``elliptic'' since the $\dot
X_k^{-\frac12+\eps,-1-\eps,2}$-norm is finite on functions which are
too singular.  But in the elliptic regime, there will be no need for
the Duhamel formula and thus for Lemma~\ref{lem:Ncut}.

\noindent The Duhamel formula (in other words, $\Box^{-1}$) introduces a
Hilbert transform  in the normal direction to the light-cone. The
following lemma is of this type.

\begin{lemma}
  \label{lem:etaT} Let $\eta$ be a smooth function on~$\R$ such that $0\le\eta\le1$, $\eta(u)=1$ on $-1\le u\le 1$, $\supp(\eta)\subset
[-2,2]$, and $\eta'(u)\ge0$ on $u\le0$,  and $\eta'(u)\le0$ on
$u\ge0$. Define $\eta_T^+(t):= \chi_{[0,\infty)} \eta(t/T)$ for each
$T\ge1$. Then, with $\chi:=\chi_{[0,\infty)}\eta'$,
\begin{equation}\label{eq:whF}
\wh{\eta_T^+}(\tau) = -\frac{1}{i\tau} (\wh{\chi}(T\tau)+1)
\end{equation}
In particular, $\wh{\eta_T^+}(\tau)=
a\wh{\eta_{\frac{T}{a}}^+}(a\tau)$ for all $0<a<1$ and
\[
\big| \wh{\eta_T^+}(\tau) \big|\les |\tau|^{-1},\quad  \big|
\frac{d}{d\tau}{\wh{\eta_T^+}}(\tau)\big|\les |\tau|^{-2}
\]
Moreover, let $\mu=\mu(\tau)$ be a smooth function on $[-1,1]$ with
$\mu(0)=1$ and $\mu\ge 1$ on~$[-1,1]$. Then
\[
\sup_{|\tau|\le1}\big| \wh{\eta_T^+}(\tau) -
\wh{\eta_T^+}(\mu(\tau)\tau)\big|\le C \|\mu'\|_\infty
\]
with an absolute constant $C$. Finally, if $T'\in [T/2,2T]$, then
\[
\big| \wh{\eta_T^+}(\tau)-\wh{\eta_{T'}^+}(\tau) \big|\le C
T\min\big(1,(T|\tau|)^{-100}\big)
\]
with a constant $C$ that only depends on~$\chi$.
\end{lemma}
\begin{proof} Integrating by parts in
\[
\wh{\eta_T^+}(\tau) = \int e^{-i\tau u} \eta_T^+(u)\, du
\]
yields \eqref{eq:whF}. In particular,
\[
\big| \wh{\eta_T^+}(\tau) \big|\les |\tau|^{-1}
\min(T|\tau|,(T|\tau|)^{-100})
\]
and similarly for the derivatives. Next, write
\begin{align*}
\wh{\eta_T^+}(\tau)-\wh{\eta_T^+}(\mu(\tau)\tau) = -\frac{1}{i\tau}
(\wh{\chi}(T\tau)+1) + \frac{1}{i\mu(\tau)\tau}
(\wh{\chi}(T\mu(\tau)\tau)+1)
\end{align*}
In view of our assumptions on~$\mu$,
\[
\big| \tau^{-1} (1-\mu(\tau)^{-1})\big|\les \|\mu'\|_\infty
\]
and similarly for the terms involving $\wh{\chi}(T\tau)$. The final
statement is an immediate consequence of~\eqref{eq:whF}.
\end{proof}

The following representation of waves $\Box^{-1}F$ with $F$ a
null-frame atom will be useful in several instances. Hence, we state
it as a separate fact.

\begin{lemma}
  \label{lem:Null_rep}
Assume that $F\in N[0]$ is a wave-packet atom, i.e.,
$F=F^+=\sum_{\kappa\in\calC_\ell} F_{\kappa}$ with
\begin{equation}
 \nn
\sum_{\kappa\in\calC_\ell} \|F_\kappa\|_{\NF[\kappa]}^2\le1
\end{equation}
for some $\ell\le -100$, see Definition~\ref{def:Nk}. Then
\[
\phi(t):= \Box^{-1} F(t) =  \int_0^t
\frac{\sin((t-s)|\nabla|}{|\nabla|} F(s)\, ds
\]
admits a decomposition of the form
\begin{equation}\label{eq:phi_rep}
\phi = \Box^{-1} F_1 + \sum_{\kappa\in\calC_\ell} \int_{\R}
\big(\Psi^1_{\kappa,a} +  B_{\kappa,a}\, \Psi^2_{\kappa,a}\big)\, da
\end{equation}
where $\|F_1\|_{L^1_t L^2_x}\les 1$ and
\[
\sup_{\kappa,a} \|B_{\kappa,a}\|_{L_{t,x}^\infty}\le C,\quad
\sum_{j=1}^2 \int\|\Psi^j_{\kappa,a}\|_{\dot X_0^{0,\frac12,1}}\, da
\le C \|F_\kappa\|_{\NF[\kappa]}
\]
with an absolute constant $C$ whence
\[
\sup_{j=1,2} \sum_{\kappa} \Big(\int\|\Psi^j_{\kappa,a}\|_{\dot
X_0^{0,\frac12,1}}\, da\Big)^2 \les 1
\]
Finally\footnote{Here $cE$  denotes the  dilation of the convex set
$E$ about its center of mass by the constant~$c$.}, for $j=1,2$
\begin{equation}\label{eq:PsiFsupp}
\supp\big(\wh{\Psi^j_{\kappa,a}} \big) \subset
C\supp\big(\wh{F_\kappa}\big), \quad  \supp\big(\wh{B_{\kappa,a}\, \Psi^2_{\kappa,a}} \big) \subset
C\supp\big(\wh{F_\kappa}\big)
\end{equation}
for all $a$ and $\kappa$ and some absolute constant~$C$.
\end{lemma}
\begin{proof}
  As in the proof of Lemma~\ref{lem:Ncut},
we first write \[\chi_{\R^+} = P_{\ge 2\ell}  (\chi_{\R^+}) + P_{<
2\ell} (\chi_{\R^+}) =: \chi_1 + \chi_2
\]
Then $F_1:= \chi_1 F$ satisfies $\|F_1\|_{L^1_t L^2_x}\les 1$,
see~\eqref{eq:PjFener}. On the other hand, $F_2:= \chi_2 F$ is again
a wave-packet atom at essentially the same scale as $F$, i.e.,
$F_2=\sum_{\kappa\in\calC_\ell} \wt F_{\kappa}$ with
\begin{equation}
 \nn
\sum_{\kappa\in\calC_\ell} \|\wt F_\kappa\|_{\NF[\kappa]}^2\le1
\end{equation}
  Define
$\Phi:= \Box^{-1} F_2$. Then
$\Phi=\sum_{\kappa}\lim_{T\to\infty}\Phi_{T,\kappa}$ with
\[
\Phi_{T,\kappa}(t) :=  \int_{-\infty}^\infty
\frac{\sin((t-s)|\nabla|)}{|\nabla|} \eta_T^+(t-s) \wt F_\kappa(s)\,
ds
\]
It suffices to prove that \[ \Phi_{T,\kappa}=\int \big(
\Psi_{T,\kappa,a}^1 + B_{T,\kappa,a} \Psi^2_{T,\kappa,a}\big)\, da\]
where
\[\sup_{T\ge1} \sup_a \|B_{T,\kappa,a}\|_\infty\les1\] and
\begin{equation}
  \label{eq:tildeFkappa}
  \sup_{j=1,2} \sup_{T\ge1} \int \| \Psi^j_{T,\kappa,a}\|_{\dot X_0^{0,\frac12,1}}\, da \les
  \|\tilde F_\kappa\|_{\NF[\kappa]}
\end{equation}
both uniformly in~$\kappa$.

\noindent Fix $\kappa\in\calC_\ell$ and $\omega=\omega(\kappa)\in
S^1\setminus(2\kappa)$ so that
\[
d(\omega,\kappa)^{-1} \|
  \wt F_\kappa\|_{L^1_{t_\omega} L^2_{x_\omega}} \le 2\|\wt F_\kappa\|_{\NF[\kappa]}
\]
As usual, we foliate relative to~$t_\omega$. More precisely, define
\[
f_a(x_\omega) = \wt F_\kappa((t,x)(a,x_\omega))
\]
where $t_\omega=a$ means that
\[
(t,x)(a,x_\omega) := a \theta_\omega^+ + x_\omega
\]
By Lemma~\ref{lem:Piomega}
\begin{equation}\label{eq:f_omega}
\supp(\hat{f_a})\subset \Pi_\omega\big(\big\{(\tau,\xi)\::\:
|\xi|\sim1,\; \wh{\xi}\in \kappa,\; |\tau-|\xi||\les
2^{2\ell}\big\}\big) =: R_{\kappa,\omega}
\end{equation}
Let  $(t_\omega,x_\omega)$ denote the null-frame coordinates. Then
\begin{equation}\label{eq:tauminxi}
\tau-|\xi| = \frac{2\xi_\omega^1}{\tau+|\xi|}
(\tau_\omega-h(\xi_\omega))
\end{equation}
where  $\xi_\omega^1:=\xi_\omega\cdot \theta_\omega^{-}$ and, with
$|\xi_\omega|^2=(\xi_\omega^1)^2 + (\xi_\omega^2)^2$, one has
$h(\xi_\omega) := \frac{(\xi_\omega^2)^2}{2\xi_\omega^1}$. Moreover,
 $|\xi^1_\omega|\sim d(\omega,\kappa)^2$ and  $|\xi^2_\omega|\les d(\omega,\kappa)$
by elementary geometry (cf.~Lemma~\ref{lem:Piomega}). We define
  \[
P_{\kappa,\omega} f := \calF^{-1} [
\chi_{R_{\kappa,\omega}}(\xi_\omega) f(\tau_\omega,\xi_\omega)]
  \]
  where $\chi_{R_{\kappa,\omega}}$ is a smooth cut-off adapted to the
  rectangle $\chi_{R_{\kappa,\omega}}$ in the $\xi_\omega$-plane.
Furthermore, we set
\[
Q_{\le j,\omega}^+ f := \calF^{-1} \big[ m_0 \big (2^{-j-C}
d^{2}(\kappa,\omega) ( \tau_\omega-h(\xi_\omega))   \big)
 f(\tau_\omega,\xi_\omega) \big]
\]
By construction, $P_{\kappa,\omega} Q_{\le 2\ell,\omega}^+$ is
essentially the same as $P_{0,\kappa} Q_{\le 2\ell}$, see
Lemma~\ref{lem:Piomega}. In fact, one has
\[
\wt F_\kappa = P_{\kappa,\omega} Q_{\le 2\ell,\omega}^+ \wt F_\kappa
\]
and $P_{\kappa,\omega} Q_{\le 2\ell,\omega}^+$ is disposable.
Clearly,
\[
\Phi_{T,\kappa} = \int \Phi_{T,\kappa,a}\, da
\]
where
\[
\Phi_{T,\kappa,a}(t) :=  P_{\kappa,\omega} Q_{\le 2\ell,\omega}^+
\int_{-\infty}^\infty \frac{\sin((t-s)|\nabla|)}{|\nabla|}
\eta_T^+(t-s) \delta(s_\omega-a) f_a \, ds
\]
Then $\Phi_{T,\kappa,a}= P_{0,C\kappa}
  Q_{\le 2\ell+C}^+ \Phi_{T,\kappa,a}$ and
\begin{equation}\label{eq:Sf2}
\begin{aligned}
 & \Phi_{T,\kappa,a} = P_{\kappa,\omega} Q_{\le 2\ell,\omega}^+
  \;\calF^{-1}[  (\wh{\eta_T^+}(\tau-|\xi|) - \wh{\eta_T^+}(|\xi|+\tau))
  e^{-i\tau_\omega a} \hat{f_a}(\xi_\omega)]
  \end{aligned}
\end{equation}
We claim that the contribution of $|\wh{\eta_T^+}(|\xi|+\tau)|\les
1$ to~\eqref{eq:Sf2} can be added to~$\Psi^1_{T,\kappa,a}$. In fact,
\begin{equation}\label{eq:O1calc}\begin{aligned}
  &\big\|
  Q_{\le 2\ell+C}^+ \;\calF^{-1}[  O(1)
  e^{-i\tau_\omega a} \chi_{R_{\kappa,\omega}}(\xi_\omega)\hat{f_a}(\xi_\omega)]
  \big\|_{\dot X_0^{0,\frac12,1}} \\
  &\les  2^{\ell}\big\| m_0(2^{-2\ell-C}(\tau-|\xi|)) \chi_{R_{\kappa,\omega}}(\xi_\omega)
  \hat{f_a}(\xi_\omega)
  \big\|_{L^2_\tau L^2_\xi}\\
  &\les 2^{2\ell} d(\omega,\kappa)^{-1} \| f_a\|_{L^2}
\end{aligned}
\end{equation}
which is better than needed. To pass to the final estimate here we
used Lemma~\ref{lem:Piomega}, especially~\eqref{eq:xiJac}; the
latter estimate can be applied for fixed~$\tau$, since then
$\xi_\omega=\xi_\omega(\xi,\tau)$. Next, we split the contribution
of $\wh{\eta_T^+}(|\xi|-\tau)$ to~\eqref{eq:Sf2} into several
pieces.
 Since $\tau_\omega-h(\xi_\omega)=0$ implies that
\[
2|\xi|=\tau+|\xi|=h(\xi_\omega)/\sqrt{2} + \xi_\omega\cdot e_1 +
\sqrt{h(\xi_\omega)^2+|\xi_\omega|^2} =: g(\xi_\omega)
\]
where $e_1=(1,0,0)$, one has by Lemma~\ref{lem:etaT}
\begin{equation}\label{eq:whdiff}
\wh{\eta_T^+}(\tau-|\xi|) - \wh{\eta_T^+}\big(
\frac{2\xi_\omega^1}{g(\xi_\omega)} (\tau_\omega-h(\xi_\omega)\big)
= O\big(\frac{1}{\xi_\omega^1}\big)=
O\big(d(\omega,\kappa)^{-2}\big)
\end{equation}
In view of~\eqref{eq:O1calc} (which {\em gains} a factor of
$2^{2\ell}\les d(\omega,\kappa)^{2}$), the contribution
of~\eqref{eq:whdiff} to~\eqref{eq:Sf2} can again be added
to~$\Psi^1_{T,\kappa,a}$. Set
$b=b(\xi_\omega)=\frac{2\xi_\omega^1}{g(\xi_\omega)}$. Furthermore,
set $b_0:= b(\xi_\omega^{(0)})$ where $\xi_\omega^{(0)}\in
R_{\kappa,\omega}$ is fixed, cf.~\eqref{eq:f_omega}.
 In view of
Lemma~\ref{lem:etaT},
\begin{align}
& \wh{\eta_T^+}\big( \frac{2\xi_\omega^1}{g(\xi_\omega)}
(\tau_\omega-h(\xi_\omega)\big) = b^{-1} \wh{\eta_{bT}^+}
 (\tau_\omega-h(\xi_\omega)) \nn \\
 & \quad = b^{-1} \wh{\eta_{b_0T}^+}
 (\tau_\omega-h(\xi_\omega)) + O\big[T\min(1,(T|\tau-|\xi||)^{-100})  \big]
 \label{eq:whdiff2}
\end{align}
where we used that $b\sim b_0$ on~$R_{\kappa,\omega}$. The
computation from~\eqref{eq:O1calc} above now shows that the
$O(\cdot)$ term in~\eqref{eq:whdiff2} can be added
to~$\Psi^1_{T,\kappa,a}$. It therefore remains to analyze the
contribution of the first term in~\eqref{eq:whdiff2}
to~\eqref{eq:Sf2}. Define
\[
B_{T,a,\kappa}(t,x):= \int \eta_{b_0T}^+(t_\omega-s_\omega)
e^{-ia(t_\omega-s_\omega)}  \lambda \wh{m_0} (\lambda s_ \omega)\, d
s_ \omega
\]
where $2^{j+C} d^{-2}(\kappa,\omega)=:\lambda$ (recall that $m_0$ is
even). On the one hand, $\|B_{T,a,\kappa}\|\le \|\wh{m_0}\|_1$ and
on the other hand,
\begin{align*}
\nn
 & P_{\kappa,\omega} Q_{\le 2\ell,\omega}^+  \;\calF^{-1}[ b^{-1} \wh{\eta_{b_0T}^+} (
 \tau_\omega-h(\xi_\omega) )
  e^{-i\tau_\omega a} \hat{f}(\xi_\omega)]  \\
& = B_{T,a,\kappa}
   \;\calF^{-1}\Big[ \delta (
 \tau_\omega-h(\xi_\omega) )\,\chi_{R_{\kappa,\omega}}(\xi_\omega) \frac{g(\xi_\omega)}{2\xi_\omega^1}
  e^{-ih(\xi_\omega) a} \hat{f_a}(\xi_\omega)\Big] \\
&=: B_{T,a,\kappa} \Psi^2_{T,\kappa,a}
\end{align*}
By inspection,  the Fourier support of $\Psi^2_{T,\kappa,a}$ as well
as that of $B_{T,a,\kappa} \Psi^2_{T,\kappa,a}$  are no larger than
that of the original wave-packet $F_\kappa$ (up to a dilation by a
constant). Finally, by a calculation similar to~\eqref{eq:O1calc},
\begin{align*}
&\|\Psi^2_{T,\kappa,a}\|_{\dot X_0^{0,\frac12,1}}\les \big\|
\calF^{-1}\big[ \delta (
 \tau_\omega-h(\xi_\omega) )\, \frac{g(\xi_\omega)}{2\xi_\omega^1}
  e^{-ih(\xi_\omega) a}\chi_{R_{\kappa,\omega}}(\xi_\omega) \hat{f_a}(\xi_\omega)\big]
  \big\|_{\dot X_0^{0,\frac12,1}} \nn \\
&\les d(\omega,\kappa)^2 \limsup_{M\to\infty} \big\|\calF^{-1}[ M
\eta(M(
 \tau-|\xi| ))\, \frac{g(\xi_\omega)}{2\xi_\omega^1}
  e^{-ih(\xi_\omega) a} \chi_{R_{\kappa,\omega}}(\xi_\omega) \hat{f_a}(\xi_\omega)]
  \big\|_{\dot X_0^{0,\frac12,1}} \\
  & \les  d(\omega,\kappa)^{-1} \| f_a\|_{L^2}
\end{align*}
This concludes the proof of the lemma.
\end{proof}

In passing, we now prove the analogue of Lemma~\ref{lem:QLp} for null-frame coordinates,
which then gives Corollary~\ref{cor:S_cut} for the~$N[k]$ spaces.

\begin{cor}
 \label{cor:N_cut}  For all $F\in N[k]$ and all $j\in\Z$ one has
$\|Q_{\le j}F\|_{N[k]}\le C\|F\|_{N[k]}$ with some absolute constant~$C$.
\end{cor}
\begin{proof}
 This is clear if $F$ is either an energy or a $\dot X^{s,b}$-atom. Therefore, suppose that
$F=\sum_{\kappa} F_\kappa$ is a wave-packet atom with $k=0$. It suffices to prove that
\[
 \|Q_{\le j}F_\kappa\|_{\NF[\kappa]} \le C \|F_\kappa\|_{\NF[\kappa]}
\]
This in turn follows from
\begin{equation}
 \label{eq:Qjframe}
\|Q_{\le j}F_\kappa\|_{L^1_{t_\omega} L^2_{x_\omega}} \le C \|F_\kappa\|_{L^1_{t_\omega} L^2_{x_\omega}}
\end{equation}
which holds uniformly in $\omega\in S^1\setminus(2\kappa)$. Fix such
an~$\omega$ and apply Plancherel's theorem in~$x_\omega$.
By~\eqref{eq:tauminxi},
\begin{align*}
 \calF_2\,  Q_{<j}F_\kappa (t_\omega,\xi_\omega) =\int m_0\Big( 2^{-j} \frac{2\xi_\omega^1}{\tau+|\xi|}
(\tau_\omega-h(\xi_\omega)) \Big)
e^{i\tau_\omega(t_\omega-s_\omega)} \, d\tau_\omega\;
 \calF_2 F_\kappa(s_\omega,\xi_\omega)\, ds_\omega
\end{align*}
where for our purposes here $\calF_2$  refers to a partial Fourier
transform relative to the second variable~$x_\omega$. In view of
$|\xi^1_\omega|\sim d(\omega,\kappa)^2$,
\begin{align*}
 \big|\calF_2\, Q_{<j} F_\kappa (t_\omega,\xi_\omega)\big| &\le \Big|\int m_0\Big( 2^{-j} \frac{2\xi_\omega^1}{\tau+|\xi|}
\tau_\omega \Big) e^{i\tau_\omega(t_\omega-s_\omega)} \,
d\tau_\omega\Big|
 \big|\calF_2 F_\kappa(s_\omega,\xi_\omega)\big|\, ds_\omega \\
 &\les 2^j d(\omega,\kappa)^{-2} \int \la 2^j d(\omega,\kappa)^{-2}
 (t_\omega-s_\omega)\ra \big|\calF_2 F_\kappa(s_\omega,\xi_\omega)\big|\, ds_\omega
\end{align*}
Performing an $L^2_{\xi_\omega}$ estimate followed by an
$L^1_{t_\omega}$ bound yields~\eqref{eq:Qjframe}.
\end{proof}

Finally, there is the following simple fact that will play a role in
the proof of the Strichartz component of~$\|\cdot\|_{S[k]}$.

\begin{lemma}
 \label{lem:trivsum} Let $a_{km}\ge0$ for all $1\le m\le M$ and $1\le k\le K$. Suppose
$\sum_{k=1}^K a_{km}\le \sigma$ for all $k$ where $\sigma\ge0$ is
arbitrary. Then
\begin{equation}\label{eq:Sdef}
 \sum_{k=1}^K \Big( \sum_{m=1}^M a_{km}^2\Big)^{\frac12} \le \sigma M^{\frac{1+\theta}{2}} K^{\frac{1-\theta}{2}}
\end{equation}
for all $0\le\theta\le 1$.
\end{lemma}
\begin{proof}
 Denote the sum in \eqref{eq:Sdef} by~$S$. On the one hand,
\[
 S\le \sum_{k=1}^K \sum_{m=1}^M a_{km} \le \sigma\, M
\]
 On the other hand,
\[
 S\le \sqrt{K} \Big(  \sum_{k=1}^K \sum_{m=1}^M a_{km}^2\Big)^{\frac12}\le \sigma \sqrt{KM}
\]
and the lemma is proved.
\end{proof}

Now we can state the main energy bound.
We begin with the easier elliptic regime.

\begin{lemma}
 \label{lem:ener_ell} Let $F$ be a space-time Schwartz function which is adapted to~$k\in\Z$. Assume
furthermore that $F=I^c F$ and set $\phi:=\Box^{-1} F$, which is defined via division by~$\tau^2-|\xi|^2$ on the Fourier side. Then
\[
 \|\phi\|_{S[k]}\les \|F\|_{N[k]}
\]
with an absolute implicit constant.
\end{lemma}
\begin{proof}
We may again assume that $k=0$. We then need to prove that
\begin{equation}
 \label{eq:phiFell}
\|\phi\|_{\dot X_k^{0,1-\eps,2}} \les \min\big(
\|F\|_{L^1_tL^2_x},\|F\|_{\dot X_0^{0,-1-\eps,2}} \big)
\end{equation}
since, as we observed after Definition~\ref{def:Sk}, the norm on the
left-hand side dominates the other norms which make
up~$\|\cdot\|_{S[k]}$.  If we select $\|F\|_{\dot
X_0^{0,-1-\eps,2}}$ on the right-hand side of~\eqref{eq:phiFell},
then this inequality is obvious. On the other hand, if we select
$\|F\|_{L^1_tL^2_x}$, then one concludes via Bernstein's inequality
in time.
\end{proof}

Next, we deal with the hyperbolic regime.

\begin{prop}
\label{prop:energy} Let $k\in\Z$ and suppose $\phi_0, \phi_1$ are Schwartz functions in~$\R^2$ which
are adapted to~$k$. Further, suppose $F$ is a space-time Schwartz function which is adapted to~$k$, and
which is moreover hyperbolic, i.e., $F=IF$. Then the unique smooth solution of
\[
 \Box\phi=F,\quad (\phi(0),\partial_t \phi(0))=(\phi_0,\phi_1)
\]
satisfies
\begin{equation}\label{eq:energy}
\|\phi\|_{S[k]} \les \|(\phi_0,\phi_1)\|_{L^2 \times  \dot H^{-1}}
+  \|F\|_{N[k]}
\end{equation}
with an absolute implicit constant.
\end{prop}
\begin{proof} By scaling we may assume that $k=0$.
We first assume that $F=0$. Then
\[
\wh{\phi(t)}(\xi) = \cos(t|\xi|)\wh{\phi_0}(\xi) + \frac{\sin(t|\xi|)}{|\xi|} \wh{\phi_1}(\xi)
\]
Consequently,  \eqref{eq:energy} follows from \eqref{eq:free-trunc} upon sending~$T\to\infty$.

Next, we assume that $\phi_0=\phi_1=0$. By the Duhamel formula,
\[
\wh{\phi(t)}=\int_0^t   \frac{\sin((t-s)|\xi|)}{|\xi|} \wh{F(s)}(\xi)\, ds
\]
In other words, we need to show that
\[
\Big\|  \int_{-\infty}^t   \frac{\sin((t-s)|\nabla|)}{|\nabla|}\chi_{[0,\infty)}(s) F(s)\, ds \Big\|_{S[0]} \les
\| F\|_{N[0]}
\]
In view of Lemma~\ref{lem:Ncut}, we may remove the indicator
function $\chi_{[0,\infty)}(t)=\chi_{\R^+}(t)$ on the left-hand
side. This is where we use that $F=IF$, but after this point we may
no longer assume that $F=IF$ since $\chi_{\R^+}F$ loses this property.

\noindent
The goal is now to prove uniformly in $T\ge1$
\begin{equation}\label{eq:ener_goal}
\Big\|  \int_{-\infty}^\infty    \sin((t-s)|\nabla|) \eta_T^+(t-s) F(s)\, ds \Big\|_{S[0]} \les
\| F\|_{N[0]}
\end{equation}
where $\eta_T^+(u):= \eta(u/T)\chi_{\R^+}(u)$ is a bump function as
specified in Lemma~\ref{lem:etaT}.  Denote
\begin{equation}
 \phi(t) =   \int_{-\infty}^\infty  \sin((t-s)|\nabla|) \eta_T^+(t-s) F(s)\, ds
\label{eq:phiDuhamel}
\end{equation}
Then the space-time Fourier transform of $\phi$ equals (up to a multiplicative constant)
\begin{equation}\label{eq:phiFourier}
   \wh{\phi}(\tau,\xi) = \big(\wh{\eta_T^+}(\tau-|\xi|)-\wh{\eta_T^+}(\tau+|\xi|)\big) \wh{F}(\tau,\xi)
\end{equation}
whence,  by Lemma~\ref{lem:etaT},
\begin{equation}\label{eq:etaTest}
 |\wh{\phi}(\tau,\xi)|\les \Big(\big||\tau|-|\xi|\big|^{-1}\chi_{[|\tau|<10]}+\tau^{-2}\chi_{[|\tau|\ge10]} \Big) |\wh{F}(\tau,\xi)|
\end{equation}
and thus also
\begin{equation}\label{eq:phi_Xsb}
\|Q_{\le 0} \phi\|_{\dot X_0^{0,\frac12,\infty}}+ \|Q_{> 0}
\phi\|_{\dot X_0^{0,1-\eps,2}} \les \|F\|_{N[0]}
\end{equation}
from \eqref{eq:Xsb_dom} and Lemma~\ref{lem:ener_ell}.   By
Lemma~\ref{lem:incl_free} it suffices to assume that $F$ is either
an energy or a wave-packet atom. Moreover, in each of these cases
the $\dot X^{0,\frac12,\infty}$ and~$\dot X^{0,1-\eps,2}$-components
of the $S[k]$ norm of $\phi$ can be ignored due
to~\eqref{eq:phi_Xsb}. Moreover, since the $\dot
X^{0,1-\eps,2}$-norm controls the entire $S[0]$-norm in the elliptic
regime, it suffices to consider only~$Q_{\le 0}\phi$.

In case $F$ is an energy atom, i.e., $\|F\|_{L^1 L^2}\le1$ standard
$X^{s,b}$ and Strichartz norms for the wave equation bound the norms
in~\eqref{eq:Sk1} and~\eqref{eq:Sk2}, see Lemma~\ref{lem:Strich2}.
We are therefore reduced to bounding~\eqref{eq:squarefunc}, for
which it suffices to verify that
\[
\sup_{\ell\le-100}\sup_{T\ge1}\sup_{\kappa\in\calC_\ell} \big\|P_{0,\kappa}Q^+_{<2\ell} \sin(t|\nabla|)\eta_T^+(t) f\big \|_{S[0,\kappa]}\les \|f\|_{L^2_x}
\]
for any $f$ which is $0$-adapted (the case of $Q^-$ being analogous). We can ignore the further localization to
the rectangle~$R$ due to orthogonality, cf.~\eqref{eq:squarefunc}.
The Fourier transform of the function inside the norms on the left-hand side is
\[
 \chi_\kappa(\hat\xi) m_0(2^{-2\ell}(|\xi|-\tau)) (\wh{\eta_T^+}(\tau-|\xi|) - \wh{\eta_T^+}(|\xi|+\tau)) \hat{f}(\xi)
\]
where $\chi_\kappa$ is a cut-off adapted to the cap~$\kappa$. The contribution by $\wh{\eta_T^+}(|\xi|+\tau)$ is controlled
by~\eqref{eq:Skkappaimbed}. As for $\wh{\eta_T^+}(|\xi|-\tau)$, one needs to show that
\[
 \big\| [\eta_T^+\ast 2^{2\ell}\wh{m_0}(2^{2\ell}\cdot)] P_{0,\kappa} e^{it|\nabla|} f\big \|_{S[0,\kappa]}\les \|f\|_{L^2_x}
\]
However, since the term in brackets is a bounded function uniformly in~$\ell$, one can again apply Lemma~\ref{lem:incl_free}.

Now assume that $F$ is a wave-packet atom, i.e., $F=\sum_{\kappa\in\calC_\ell} F_{\kappa}$ with
\begin{equation}
 \label{eq:Fnull}
\sum_{\kappa\in\calC_\ell} \|F_\kappa\|_{\NF[\kappa]}^2\le1
\end{equation}
where the $F_\kappa$ have the wave-packet form as specified in
Definition~\ref{def:Nk}.  We need to show that
\begin{equation}\label{eq:FnullS}
\sup_{\pm}\sup_{\ell'\le-100}\;\sup_{\ell'\le m\le 0}
  \Big(  \sum_{\kappa'\in\caps_{\ell'}} \sum_{R\in\calR_{0,\pm \kappa',m}}
  \|P_{R}
  Q_{\le 2\ell'-C}^\pm\;
  \phi\|_{S[0,\kappa']}^2 \Big)^{\frac12} \les 1
\end{equation}
We first consider the case $\ell'\le\ell$. Lemma~\ref{lem:Nsquare}
implies that it suffices to assume that $\ell'=\ell$ and to show
that, uniformly in $\kappa\in\calC_{\ell}$,
\[
\|
  \phi_\kappa\|_{S[0,\kappa]}  \le C \|F_\kappa\|_{\NF[\kappa]}
\]
with an absolute constant~$C$ where
\[
\phi_\kappa:=\int_{-\infty}^\infty    \sin((t-s)|\nabla|)
\eta_T^+(t-s) F_\kappa(s)\, ds
\]
However, this follows immediately from Lemma~\ref{lem:Null_rep}
applied to~$\phi_\kappa$, the stability
property~\eqref{eq:Skkappainfty}, and the
imbedding~\eqref{eq:Skkappaimbed}; note that the term $\Box^{-1}F_1$
in~\eqref{eq:phi_rep} can be ignored as it was dealt with in the
beginning of this proof. Finally, the case $\ell'\ge\ell$ is reduced
the to $\ell'=\ell$ by means of Lemma~\ref{lem:square_func} (note
that the Fourier-support of $\phi_\kappa$ equals that
of~$F_\kappa$).

\noindent It remains to control the Strichartz norms~\eqref{eq:Sk2}. Due to Corollary~\ref{cor:N_cut},
we may ignore the projection~$Q_{<j}$.
We split the argument into two parts: First, we will prove the
estimate
\begin{equation}
 \label{eq:Duheasy}
\Big(\sum_{c\in \calD_{0,\ell}} \Big\| P_c \int_{-\infty}^\infty  e^{\pm i(t-s)|\nabla|} F(s)\, ds
 \Big\|^2_{L^4_t L^\infty_x}\Big)^{\frac12} \les 2^{\lhalb} \|F\|_{N[0]}
\end{equation}
for any $F$ as in~\eqref{eq:Fnull}, cf.~Lemma~\ref{lem:KBern}.
Second, we take the $\eta_T^+$ cut-off as in~\eqref{eq:phiDuhamel}
into account which then yields the full result. This second step is done by an adaptation of the Christ-Kiselev argument
and will result in the loss of a power $2^{\ell\delta}$ where $\delta>0$ can be made arbitrarily small.
Lemma~\ref{lem:Strich2} reduces the proof of~\eqref{eq:Duheasy} to
the bound
\[
 \Big\|  \int_{-\infty}^\infty  e^{\mp is|\nabla|} F(s)\, ds \Big\|_{L^2_x} \les  \|F\|_{N[0]}
\]
By orthogonality, it suffices to show that uniformly in $\kappa$
\[
 \Big\|  \int_{-\infty}^\infty  e^{\mp is|\nabla|} F_\kappa(s)\, ds \Big\|_{L^2_x} \les
\inf_{\omega\not\in 2\kappa} d(\omega,\kappa)^{-1} \|F_\kappa\|_{L^1_{t_{\omega}}L^2_{x_{\omega}}}
\]
with $F_\kappa$ as in~\eqref{eq:Fnull}. By Plancherel, this is the same as
\[
 \big\|  \wh F_\kappa(\pm |\xi|,\xi) \big\|_{L^2_\xi} \les
 d(\omega,\kappa)^{-1} \|F_\kappa\|_{L^1_{t_{\omega}}L^2_{x_{\omega}}}
\]
where we choose an arbitrary $\omega\not\in 2\kappa$. As above, we
may set $F_\kappa=\delta(  t_{\omega} -t_{\omega}^{(0)}
)f_\kappa(x_{\omega})$ where $t_{\omega}^{(0)}\in\R$ is an arbitrary
number and $f_\kappa\in L^2_{x_{\omega}}$ is an arbitrary function
whose Fourier support is contained in the projection of the Fourier
support of $F_\kappa$ onto the $\xi_\omega$-plane. This reduces us
further to the bound
\begin{equation}\label{eq:L2NFdom}
 \| \hat f_\kappa(\xi_\omega) \|_{L^2_{\xi}} \les d(\omega,\kappa)^{-1} \|f_\kappa\|_{L^2_{\xi_\omega}}
\end{equation}
where on the left-hand side we regard $\xi_\omega$ as a function
of~$\xi$. By Lemma~\ref{lem:Piomega} the Jacobian obeys
$\Big|\frac{\del\xi}{\del \xi_\omega}\Big|\sim
d(\omega,\kappa)^{-2}$ which implies~\eqref{eq:L2NFdom}. This
concludes the first step, i.e., the proof of~\eqref{eq:Duheasy}.
Note  that our proof of~\eqref{eq:Duheasy} applies to any $F$ which
can be written in the form $F=\sum_\kappa F_\kappa$ provided
$F_\kappa$ satisfy
\[ \supp({\wh F_\kappa})\subset \{\xi\in\R^2\::\: |\xi|\sim 1, \; \wh{\xi} \in\kappa\}
\]
In other words, one does not need any condition on the modulations
of~$F_\kappa$. This fact will be most important for the remainder of
the proof (since we will need to multiply $F$ by cutoff functions in
time). Our next goal is to establish the estimate, with $\delta>0$
arbitrarily small,
\begin{equation}
 \label{eq:DuhStrich}
 \Big(\sum_{c\in\calD_{0,\ell}}\Big\| P_c \int_{-\infty}^\infty  e^{\pm i(t-s)|\nabla|}
  \eta_T^+(t-s) F(s)\, ds \Big\|^2_{L^4_t L^\infty_x}\Big)^2 \les 2^{(\frac12-\delta)\ell} \| F\|_{N[0]}
\end{equation}
for any $F$ as in~\eqref{eq:Fnull} but without any restriction on the modulations of each $F_\kappa$. For the remainder
of the proof we will {\em fix} such a Schwartz function~$F$.  Moreover, $\|\cdot\|^2$ without any subscripts will mean the sum in~\eqref{eq:Fnull}.
As mentioned before, we prove~\eqref{eq:DuhStrich} by an adaptation of the Christ-Kiselev lemma. The latter
does not apply directly since the null-frame norm in~\eqref{eq:Fnull} is not of pure Lebesgue type.
We make the following preliminary observation:  the map $\mu(E):= \| \chi_E  F\|^2$ is a $\sigma$-subadditive
set function on the Borel sets of~$\R$; here $\chi_E=\chi_E(t)$ acts only in the time variable.
To prove this, let $\{E_j\}\subset\R$ be an arbitrary collection of disjoint Borel sets.  Then
\begin{align*}
 \sum_j \mu(E_j) & = \sum_j \| \chi_{E_j}  F\|^2 = \sum_\kappa \sum_j  \inf_{\omega\notin2\kappa} d(\omega,\kappa)^{-2}
 \| \chi_{E_j}  F_\kappa\|^2_{L^1_{t_\omega} L^2_{x_\omega}} \\
& \le   \sum_\kappa   \inf_{\omega\notin2\kappa} d(\omega,\kappa)^{-2}  \sum_j \| \chi_{E_j}  F_\kappa\|^2_{L^1_{t_\omega} L^2_{x_\omega}} \\
&\le   \sum_\kappa   \inf_{\omega\notin2\kappa} d(\omega,\kappa)^{-2} \Big(
 \sum_j \| \chi_{E_j}  F_\kappa\|_{L^1_{t_\omega} L^2_{x_\omega}} \Big)^2 \\
& =  \sum_\kappa   \inf_{\omega\notin2\kappa} d(\omega,\kappa)^{-2}  \|   F_\kappa\|_{L^1_{t_\omega} L^2_{x_\omega}}^2 \\
& = \sum_\kappa  \| F_\kappa\|_{\NF[\kappa]}^2 \le 1
\end{align*}
as claimed.  In view of this property it suffices to prove~\eqref{eq:DuhStrich}  for $F$ which are supported
on intervals\footnote{Strictly speaking, one would need to choose something like $10T$ here to accommodate the
support of $\eta_T^+$, but we ignore this issue here.}  of size~$T$ in time and we may also replace $\eta_T^+$ by the indicator $\chi_{[s<t]}$.
We now perform a Whitney decomposition  of the triangle
\[
 \Delta_T:=\{ (t,s)\::\: 0\le s\le t\le T \}
\]
by means of squares (we have shifted the support of $F$ to be contained in~$[0,T]$). This yields finitely many disjoint squares of
the form $$\calQ:=\Big\{ I_{m,n}\times J_{m,n}\}_{n\ge0,\;1\le m\le M_n}$$ with intervals $I_{m,n}, J_{m,n}$  such that  $M_n\le 2^n$ and
\begin{align*}
  \Delta &= \bigcup_{n\ge0} \bigcup_{1\le m\le 2^n}  I_{m,n}\times J_{m,n} \\
 |I_{m,n}|&=|J_{m,n}|=T2^{-n} \quad\forall\; 1\le m\le M_n,\; n\ge0 \\
 \dist(  I_{m,n}\times J_{m,n},\{s=t\})&\in (T2^{-n}/10,10 T2^{-n}) \quad\forall\; 1\le m\le M_n,\; n\ge0
\end{align*}
We call any two intervals $I,J$ of length $T2^{-n}$ {\em related} provided $I\times J\in\calQ$.
Note that any $I$ can be related to at most $20$ of the $J$ intervals.
To each $n\ge0$ we now also associate $2^n$ pairwise disjoint intervals $\{\wt J_{m,n}\}_{1\le m\le 2^n}$ which partition
$[0,T]$ and
with the property that
\[
 \mu(\wt J_{m,n}) =  \mu(\wt J_{m',n})  \quad \forall\; 1\le m,m'\le 2^n
\]
The subadditivity of~$\mu$ implies that $\mu(\wt J_{m,n}) \le
2^{-n}$. Finally, we introduce an auxiliary function $\Phi$ which is
piece-wise linear, strictly increasing on $[0,T]$ and which has the
property that $\Phi( J_{m,n})= \wt J_{m,n}$. In view of all these
properties
\begin{align*}
 & \sum_{c\in\calD_{0,\ell}}\Big\| P_c \int_{-\infty}^t  e^{\pm i(t-s)|\nabla|}  F(s)\, ds \Big\|^2_{L^4_t L^\infty_x} \\
& \les \sum_{c\in\calD_{0,\ell}}\Big(  \sum_{n=0}^\infty \Big\|
\sum_{m=1}^{M_n} \chi_{\Phi(I_{m,n})} (t) P_c \int_{-\infty}^t
e^{\pm i(t-s)|\nabla|} \chi_{\Phi({ J}_{m,n})} (s) F(s)\, ds
\Big\|_{L^4_t L^\infty_x} \Big)^2
\end{align*}
Applying Cauchy-Schwarz to the sum over $n$ allows one to bound this
further as
\begin{align}
 & \les \sum_{c\in\calD_{0,\ell}}
\sum_{n=0}^\infty (1+n)^2  \Big\| \sum_{m=1}^{M_n}
\chi_{\Phi(I_{m,n})} (t) P_c \int_{-\infty}^t  e^{\pm
i(t-s)|\nabla|}
\chi_{\Phi({ J}_{m,n})} (s) F(s)\, ds \Big\|^2_{L^4_t L^\infty_x}  \nn \\
&  \les  \sum_{n=0}^\infty  (1+n)^2  \sum_{c\in\calD_{0,\ell}}
\Big(\sum_{m=1}^{2^n} \Big\|  P_c \int_{-\infty}^t  e^{\pm
i(t-s)|\nabla|} \chi_{ {\wt J}_{m,n}} (s)  F(s)\, ds \Big\|_{L^4_t
L^\infty_x}^4 \Big)^{\frac12}  \label{eq:L4exp}
\end{align}
Label the disks $c\in\calD_{0,\ell}$ by $\{c_k\}_{k=1}^K$, $K\sim 2^{-2\ell}$, and denote for fixed $n$,
\begin{align*}
 a_{km,n} &:= \Big\|  P_{c_k} \int_{-\infty}^t  e^{\pm i(t-s)|\nabla|}
\chi_{ {\wt J}_{m,n}} (s)  F(s)\, ds \Big\|_{L^4_t L^\infty_x}^2
\\& = \Big\|  P_{c_k} \int_{-\infty}^\infty  e^{\pm i(t-s)|\nabla|} \chi_{
{\wt J}_{m,n}} (s)  F(s)\, ds \Big\|_{L^4_t L^\infty_x}^2
\end{align*}
The previous bound now takes the form
\[
 \sum_{c\in\calD_{0,\ell}}\Big\| P_c \int_{-\infty}^t  e^{\pm i(t-s)|\nabla|}  F(s)\, ds \Big\|^2_{L^4_t L^\infty_x}
\les \sum_{n=0}^\infty  (1+n)^2  \sum_{k=1}^K \Big(\sum_{m=1}^{2^n} a_{km,n}^2 \Big)^{\frac12}
\]
In view of~\eqref{eq:Duheasy} (and the remark at the end of its
proof concerning time cutoffs)
\begin{align*}
\sum_{k=1}^K a_{km,n} &  \les  2^{\ell} \big\|  \chi_{{\wt J}_{m,n}}
F \big\|^2
 = 2^{\ell}  \mu(\wt J_{m,n})
\le 2^{\ell}   2^{-n}
\end{align*}
By Lemma~\ref{lem:trivsum} with $\sigma=2^{\ell-n}$, $M=2^{n}$,
$K=2^{-2\ell}$,
\[
 \sum_{k=1}^K \Big(\sum_{m=1}^{2^n} a_{km,n}^2 \Big)^{\frac12} \les 2^{(1-2\delta)\ell} 2^{-\delta n}
\]
for any $0\le\delta\le 1$.  In view of~\eqref{eq:L4exp}, one
obtains~\eqref{eq:DuhStrich}.
\end{proof}

As a simple corollary, we now obtain the following continuity
result. Recall that the norm of~$S[k]$ can also be defined for
non-integer~$k$, cf.~\eqref{eq:Sk_scale'}. The continuity in~$k$ is
not obvious due to the various Fourier multipliers in~\eqref{eq:Sk2}
and~\eqref{eq:squarefunc} over infinitely many scales.

\begin{cor}
 \label{cor:contSk} Let $\phi$ be a   Schwartz function in~$\R^{1+2}$ which is adapted to~$k\in\R$. Then
\[
 \lim_{h\to0}\| \phi\|_{S[k+h]}  = \|\phi\|_{S[k]}
\]
\end{cor}
\begin{proof}
 By \eqref{eq:Sk_scale'},
\[
 \lambda^{-1} \| \phi(\lambda^{-1}\cdot)\|_{S[k]}= \|\phi\|_{S[k+\log_2\lambda]}
\]
It therefore suffices to note that by the energy estimate
\begin{align*}
 \big| \lambda^{-1}\| \phi(\lambda^{-1}\cdot)\|_{S[k]} - \| \phi \|_{S[k]} \big| &\le \|\lambda^{-1}\phi(\lambda^{-1}\cdot)- \phi\|_{S[k]} \\
&\le \|(\lambda^{-1}\phi(\lambda^{-1}\cdot)- \phi)[0]\|_{L^2\times \dot H^{-1}} + \|\Box(\lambda^{-1}\phi(\lambda^{-1}\cdot)- \phi)\|_{N[k]}  \\
&\les \|(\lambda^{-1}\phi(\lambda^{-1}\cdot)- \phi)[0]\|_{L^2\times \dot H^{-1}} + \|\Box(\lambda^{-1}\phi(\lambda^{-1}\cdot)- \phi)\|_{L^1_t \dot H^{-1}}
\to 0
\end{align*}
as $\lambda\to1$.
\end{proof}

\subsection{A stronger $S[k]$-norm, and time localizations}

The energy estimate of Proposition~\ref{prop:energy} and Lemma~\ref{lem:ener_ell} can be summarized as the statement that
$\|\phi\|_{S[k]} \les  \trip \phi\trip_{S[k]}$ where
\begin{equation}\label{eq:trip_def}
 \trip\phi\trip_{S[k]} :=   \|\phi\|_{\ener}
+  \|\Box \phi\|_{N[k]}
\end{equation}
for any space-time Schwartz function $\phi$ which is adapted to~$k\in\Z$.
To see this, one estimates
\begin{align*}
 \|\phi\|_{S[k]} &\les \|I\phi\|_{S[k]} + \|I^c \phi\|_{S[k]} \\
&\les \|((I\phi)(0),(\del_t I\phi)(0))\|_{L^2\times \dot H^1} + \|I\Box\phi\|_{N[k]} + \|I^c \Box\phi\|_{N[k]} \\
&\les \|\phi\|_{\ener} + \|\Box\phi\|_{N[k]} = \trip\phi\trip_{S[k]}
\end{align*}
To remove $I$ from the right-hand side here one uses Corollary~\ref{cor:N_cut}.

We shall henceforth use this stronger norm and the resulting smaller $S[k]$-space.
We introduce this norm because it leads to an improvement over the
bilinear bound~\eqref{eq:bilin2} in the case of high-high
interactions, see Lemma~\ref{lem:bilin3} below. This improvement
reflects a smoothing effect of convolutions of measures supported on
the light cone. It thus {\em cannot} be obtained using the $S[k,\kappa]$
norms alone, since~\eqref{eq:bilin2} is based on H\"older's
inequality
\[ L^2_{t_\omega}L^\infty_{x_\omega}\cdot  L^\infty_{t_\omega}
L^2_{x_\omega}\hookrightarrow \Ltwotx
\]
which does not allow for any gain in regularity. It will be essential to
note that Corollary~\ref{cor:S_cut} still applies to the stronger
norm~$\trip\cdot\trip$:

\begin{lemma}
  \label{lem:tripQj} For all $\phi$ which are adapted to $k\in\Z$ and all $j\in\Z$ one has
$\trip Q_{\le j}\phi\trip_{S[k]}\le C\trip \phi\trip_{S[k]}$ with
some absolute constant~$C$.
\end{lemma}
\begin{proof} This follows immediately from Lemma~\ref{lem:QLp} and
Corollary~\ref{cor:N_cut}.
\end{proof}

Another property which the stronger norm inherits is that it is finite on free wave, cf.~Lemma~\ref{lem:incl_free}. More
precisely,  for any $\phi$ which is adapted to~$k$ and satisfies $\phi=Q_{\le k}\phi$,
\begin{align*}
 \trip\phi\trip_{S[k]} &=  \|\phi\|_{\ener}
+  \|\Box \phi\|_{N[k]} \\
&\le  \|\phi \|_{\ener}
+  \|\Box \phi\|_{\dot X_k^{0,-\frac12,1}}
\les \|\phi\|_{\dot X_k^{0,\frac12,1}}
\end{align*}
As in \cite{Krieger}, one needs to allow for time-localized versions of $S[k]$, both relative to the
original $\|\cdot\|_{S[k]}$, as well as the stronger $\trip\cdot\trip$-norm. This has to do with the fact
that the we need to derive apriori bounds in these spaces for Schwartz functions~$\psi_\alpha$ which
satisfy~\eqref{eq:psisys1}--\eqref{eq:psi_wave} on some time interval~$[-T,T]$. Since the norms of the $S[k]$
and~$N[k]$ spaces are defined in phase space, one cannot simply define these norms by time truncations.
Rather, one proceeds as in~\cite{T1} and~\cite{Krieger} by means of Schwartz extensions: with $\psi$ and~$\tilde\psi$ both
Schwarz functions, and $T\ge0$,
\begin{equation}\label{eq:Skloc}
\begin{aligned}
 \|\psi\|_{S[k]([-T,T]\times\R^2)}  &:= \inf_{\tilde\psi|_{[-T,T]} = \psi|_{[-T,T]}} \|P_k \tilde \psi\|_{S[k]} \\
\trip \psi\trip_{S[k]([-T,T]\times\R^2)}  &:= \inf_{\tilde\psi|_{[-T,T]} = \psi|_{[-T,T]}} \trip P_k \tilde \psi\trip_{S[k]}
\end{aligned}
\end{equation}
It is easy to see that the triangle inequality holds for these expressions and that they are actually norms.
 Moreover, it is clear that these norms are nondecreasing in~$T$.
Following~\cite{Krieger}, we now verify that these norms are
continuous in~$T$.

\begin{lemma}
 \label{lem:Tcont}  Let $\psi$ be the restriction of some Schwartz function~$\psi_0$ in~$\R^{1+2}$ to the time
interval~$[-T_0,T_0]$ where $T_0>0$. Then
\[
 \|\psi\|_{S[k]([-T,T]\times\R^2)} \text{\ \ and\ \ } \trip \psi\trip_{S[k]([-T,T]\times\R^2)}
\]
are nondecreasing and continuous in~$0\le T<T_0$.
\end{lemma}
\begin{proof}
The definition of $S[k]$ with respect to either norm can be extended to non-integer~$k$. Given $T>0$, let $|\eps|$ be very small
and set $\lambda:=\frac{T+\eps}{T}$. Then
\[
 \|P_k  \psi\|_{S[k]([-T-\eps,T+\eps]\times\R^2)} =  \| P_{k+\mu} \psi_\lambda\|_{S[k+\mu]([-T,T]\times\R^2)}
\]
where $\mu:=\log_2\lambda$ and $\psi_\lambda(t,x):= \lambda \psi(\lambda t,\lambda x)$, and similarly for~$\trip\cdot\trip$.
Clearly,  for $\eps>0$,
\begin{align*}
&    \big| \|P_k\psi\|_{ S[k]([-T-\eps,T+\eps]\times\R^2)} - \|P_k\psi\|_{ S[k]([-T,T]\times\R^2)} \big| \\
& = \big|  \|P_{k+\mu} \psi_\lambda \|_{ S[k+\mu]([-T,T]\times\R^2)}  - \|P_k\psi\|_{ S[k]([-T,T]\times\R^2)} \big|     \\
&\le  \big| \|P_{k+\mu} \psi  \|_{ S[k+\mu]([-T,T]\times\R^2)} -
\|P_k\psi\|_{ S[k]([-T,T]\times\R^2)} \big| +  \|P_{k+\mu}(\psi-
\psi_\lambda ) \|_{ S[k+\mu]([-T,T]\times\R^2)}
\end{align*}
By the energy estimate,
\begin{align*}
  \|P_{k+\mu}(\psi- \psi_\lambda ) \|_{ S[k+\mu]([-T,T]\times\R^2)} &\les \|P_{k+\mu}(\psi- \psi_\lambda ) \|_{ S[k+\mu]} \\
&\les \|(\psi-\psi_\lambda)[0]\|_{L^2\times\dot H^{-1}} + \| \Box P_{k+\mu}(\psi- \psi_\lambda ) \|_{ N[k+\mu]}\\
& \les \|(\psi-\psi_\lambda)[0]\|_{L^2\times\dot H^{-1}} + \| \Box P_{k+\mu}(\psi- \psi_\lambda ) \|_{L^1_t \dot H^{-1}(\R^{1+2}) }\to0
\end{align*}
as $\lambda\to1$. By Corollary~\ref{cor:contSk},
\[
 \lim_{\lambda\to1}\|P_{k+\mu} \psi_\lambda \|_{ S[k+\mu]([-T,T]\times\R^2)} = \|P_k\psi\|_{ S[k]([-T,T]\times\R^2)}
\]
which implies that
\[
 \lim_{\eps\to0+} \|P_k\psi\|_{ S[k]([-T-\eps,T+\eps]\times\R^2)} = \|P_k\psi\|_{ S[k]([-T,T]\times\R^2)}
\]
as claimed. The case of $T=0$ follows directly from the energy estimate. The case of~$\trip\cdot\trip$ is essentially
the same.
\end{proof}

We define localized $N[k]$-norms similarly, i.e.,
\[
 \|\psi\|_{N[k]([-T,T]\times\R^2)}  := \inf_{\tilde\psi|_{[-T,T]} = \psi|_{[-T,T]}} \|P_k \tilde \psi\|_{N[k]} \\
\]
for Schwartz functions. In particular, one has a localized version
of~\eqref{eq:trip_def} \begin{equation}\nn
 \trip\phi\trip_{S[k]([-T,T]\times\R^2)} :=   \| \phi \|_{L^2(I;L^2(\R^2))}
+  \|\Box \phi\|_{N[k]([-T,T]\times\R^2)}
\end{equation}
Furthermore, later we will also need localized norms on asymmetric
time intervals $[-T',T]$ for which the results here of course
continue to hold.

\subsection{Solving the inhomogeneous wave equation in the Coulomb gauge}
\label{subsec:waveeq}

Consider the wave equation \eqref{eq:psi_wave}, i.e., $\Box\psi_\alpha=F_\alpha$.
Here $F_\alpha$ is a nonlinear expression in~$\psi$, but we will not pay attention to
this now. In the sequel, we shall require apriori bounds on $\psi_\alpha$ in the~$S[k]$-space.
To do so, we reduce matters to the energy estimates  of Section~\ref{subsec:energy} as follows:
writing  (suppressing $\alpha$ for simplicity)
\[
 \Box \psi  = IF  + I^c F
\]
one concludes (with both $\psi$ and $F$ global space-time Schwartz functions adapted to frequency~$1$),
\begin{equation}
\label{eq:wave_solution}
 \psi (t) = S(t-t_0)(I\psi)[t_0] + \int_{t_0}^t U(t-s) IF(s)\, ds + \Box^{-1} I^c F
\end{equation}
where the final term is obtained by division by the symbol
of~$\Box$, and the first two terms represent the free wave and the
Duhamel integral, respectively. Note that the first term here
implicitly depends on all of $\psi$, not just $\psi[t_0]$, and so in
order to actually obtain a bound on $\|\psi\|_S$, one needs to
implement a bootstrap argument. Specifically, assume that we apriori
have a bound on
\[
 \|\psi|_{[-T_0, T_0]}\|_S
\]
for some $T_{0}>0$. Also, assume that we define $I=\sum_{k\in\Z}P_k Q_{<k+C}$ where $2^{C}\gg T_{0}^{-1}$.
Then, using the energy estimate from Section~\ref{subsec:energy},  we claim that
\begin{equation}
  \label{eq:IIc}
\|\psi\|_S \les T_0^{-1}\|\psi|_{[-T_0, T_0]}\|_S + \|F\|_N
\end{equation}
where the implied constant is absolute (the $T_0^{-1}$ here comes from the time-derivative in the initial data).
Indeed, this follows from
\[
 (I\psi)[t_{0}]=(I(\chi_{[-T_{0}, T_{0}]}\psi)) [t_{0}]+ (I([1-\chi_{[-T_{0}, T_{0}]}]\psi) ) [t_{0}]
\]
and
\[
\min_{t_{0}\in [-T_0, T_0]}\|\big((I([1-\chi_{[-T_{0}, T_{0}]}]\psi)\big)) [t_{0}]\|_{L_{x}^{2}\times \dot{H}^{-1}}\ll \|\psi\|_{S}
\]
The above energy inequality then follows immediately.

\noindent It is apparent that in order to use this energy
inequality, one needs to establish an apriori bound for $\psi$ on a
small time interval $[-T_0, T_0]$. In fact, in later applications we
will always split the estimates for $P_k \psi$ into the small-time
case  $|t-t_0|\le \eps_1 2^{-k}$ and the large  time case
$|t-t_0|\ge \eps_1 2^{-k}$ (with a small $\eps_1$ that is determined
by the specific context - this then requires the constant~$C$ in the
definition of~$I$ to be large). In the small time case, the
necessary apriori bound is derived from the div-curl
system~\eqref{eq:psisys1}, \eqref{eq:psisys2} for the gauged
components. This information is then fed into the large-time case as
described above.

\section{Hodge decomposition and null-structures}
\label{sec:hodge}

Here we introduce the actual system of wave equations for which our
$S$ and $N$-spaces allow us to deduce apriori estimates. From the
discussion at the very beginning, we recall that the Coulomb
components $\psi_{\alpha}$ satisfy the system \eqref{eq:psi_wave},
which has the schematic form
\begin{equation}\label{eq:psiwave2}
\Box \psi_\alpha = i\partial^{\beta}[\psi_\alpha
A_{\beta}]-i\partial^{\beta}[\psi_\beta
A_{\alpha}]+i\partial_{\alpha}[\psi^{\beta} A_{\beta}]
\end{equation}
where $A_{\beta}$ denotes the Coulomb gauge potential
\[
A_{\beta}=\sum_{j=1,2}\triangle^{-1}\partial_{j}[\psi^{1}_{\beta}\psi^{2}_{j}-\psi^{2}_{\beta}\psi^{1}_{j}]
\]
This system in an of itself does not appear to lend itself to good
estimates, and to overcome this we have to use a key additional
feature, namely the fact that {\em the flow of \eqref{eq:psi_wave}
preserves the div-curl system} \eqref{eq:psisys1},
\eqref{eq:psisys2} in the obvious sense: if the $\psi_{\alpha}$ at time $t=0$ are the
Coulomb derivative components of an actual map, whence \eqref{eq:psisys1},
\eqref{eq:psisys2} holds at time $t=0$, then the corresponding
solution of~\eqref{eq:psi_wave} satisfies this system on its entire
time interval of existence.
 The div-curl system allows us to decompose the
components $\psi_\alpha$ as the sum of a gradient term and an error
term solving an elliptic equation, see \eqref{eq:dyn_dec}. Thus we
have schematic identities of the form
\[
\psi_\alpha=R_{\alpha}\psi+\chi_\alpha
\]
Substituting the gradient terms introduces the desired
null-structure. The present section serves to make this
decomposition of the nonlinear source terms precise. We now describe
this procedure for each of the three terms on the right-hand side
of~\eqref{eq:psiwave2}. First, define
$\del_j^{-1}:=\Delta^{-1}\del_j$ and
\begin{align*}
 \calN_{\beta j}(\psi,\psi) &= R_\beta \psi^1 R_j \psi^2 - R_j \psi^1 R_\beta
\psi^2 \\
 \calN_{\beta j}(\psi,\chi) &= R_\beta \psi^1 \chi_j^2 - R_j \psi^1
\chi^2_\beta \\
\calN_{\beta j}(\chi,\psi) &= \chi^1_\beta R_j\psi^2 - \chi^1_j
R_\beta \psi^2 \\
\calN_{\beta j}(\chi,\chi) &=  \chi^1_\beta \chi_j^2 - \chi^1_j
\chi_\beta^2
\end{align*}
Then,  adopting the Einstein summation convention,
\begin{align}
i\partial^{\beta}[\psi_\alpha A_{\beta}] &=
i\partial^{\beta}[\psi_\alpha \, I^c \del_j^{-1}\calN_{\beta
j}(\psi,\psi)] + i\partial^{\beta}[\psi_\alpha \, I
\del_j^{-1}\calN_{\beta j}(\psi,\psi)]  \label{eq:I}\\ & +
i\partial^{\beta}[\psi_\alpha \,  \del_j^{-1}\calN_{\beta
j}(\psi,\chi)] + i\partial^{\beta}[\psi_\alpha \,
\del_j^{-1}\calN_{\beta j}(\chi,\psi)]+
i\partial^{\beta}[\psi_\alpha \, \del_j^{-1}\calN_{\beta
j}(\chi,\chi)] \nn
\end{align}
The two main terms here are the trilinear ones in~$\psi$. We
introduced the modulation cutoff~$I$ in front of~$\calN_{\beta j}$
since the two resulting expressions are estimated differently: for
the first, one uses a trilinear null-form structure,
see~\eqref{eq:nullexp} below, whereas for the second the bilinear
null-form~$\calN_{\beta j}$ suffices. Note that the other three
terms involving~$\chi$ are quintilinear and septilinear in~$\psi$,
respectively, due to~\eqref{eq:chi_beta}. These are discussed in
greater detail below, under the heading ``higher order errors''.

\noindent  Next,
\begin{align*}
   -i\partial^{\beta}[\psi_\beta
A_{\alpha}] &= -i\partial^{\beta}[\psi_\beta \,\del_j^{-1}
\calN_{\alpha j}(\psi,\psi)] -i\partial^{\beta}[\psi_\beta
\,\del_j^{-1} \calN_{\alpha j}(\psi,\chi)]
-i\partial^{\beta}[\psi_\beta \,\del_j^{-1} \calN_{\alpha
j}(\chi,\psi)] \\
&-i\partial^{\beta}[\psi_\beta \,\del_j^{-1} \calN_{\alpha
j}(\chi,\chi)]
\end{align*}
The $\chi$-terms need to be decomposed further, whereas the main
term here is again the trilinear one in~$\psi$, which we now rewrite
as follows:
\begin{align}
-i\partial^{\beta}[\psi_\beta \,\del_j^{-1} \calN_{\alpha
j}(\psi,\psi)] &= -i\partial^{\beta}[\psi_\beta \,\del_j^{-1} I^c
\calN_{\alpha j}(\psi,\psi)]  -i\partial^{\beta}[\psi_\beta
\,\del_j^{-1} I \calN_{\alpha j}(\psi,\psi)] \label{eq:II}
\end{align}
The first term on the right-hand side will be estimated as is,
whereas the second term now needs to be rewritten according to the
Littlewood-Paley trichotomy, in order to make it amenable to our
estimates:
\begin{align}
&-i\partial^{\beta}[\psi_\beta \,\del_j^{-1} I \calN_{\alpha
j}(\psi,\psi)] = \nn \\
&= -i\sum_{k} P_k\partial^{\beta}\psi_\beta \,\del_j^{-1} I
P_{<k-5}\calN_{\alpha j}(\psi,\psi)-i\sum_{k} P_k\psi_\beta\,
\,\del_j^{-1} I P_{<k-5}\partial^{\beta}\calN_{\alpha j}(\psi,\psi) \label{eq:IIHL} \\
& -i \sum_{k}\partial^{\beta}P_k[P_{>k}R_\beta\psi\, \del_j^{-1} I
P_{>k}\calN_{\alpha j}(\psi,\psi)]  -i
\sum_{k}\partial^{\beta}P_k[P_{>k}\chi_\beta \,\del_j^{-1} I
P_{>k}\calN_{\alpha j}(\psi,\psi)]\label{eq:IIHH} \\
& -i\sum_{k}
\partial^{\beta} [P_{<k-5}R_\beta\psi\, \del_j^{-1} I
P_{k}\calN_{\alpha j}(\psi,\psi)] -i\sum_{k}
\partial^{\beta} [P_{<k-5}\chi_\beta \,\del_j^{-1} I
P_{k}\calN_{\alpha j}(\psi,\psi) ] \label{eq:IILH}
\end{align}
The terms involving $\chi$ are expanded further as explained below.
For the first term on the right-hand side of~\eqref{eq:IIHL} one
replaces~$\del^\beta \psi_\beta$ by the right-hand side
of~\eqref{eq:psisys2} which leads to a quintilinear term. The second
term can be estimated since the $\del^\beta$-term falls on the small
frequencies.

\noindent  Finally, the third term in~\eqref{eq:psiwave2} is treated
as follows:
\begin{align}
i\partial_{\alpha}[\psi^{\beta} A_{\beta}] &=
i\partial_{\alpha}[\psi^{\beta} I^c A_{\beta}] +
i\partial_{\alpha}[\psi^{\beta} I A_{\beta}] \label{eq:III}
\end{align}
The first term on the right-hand side of~\eqref{eq:III} is estimated
as is; in fact, it is essential that one does not perform the Hodge
decomposition in the first slot since otherwise $\beta=0$ would
create problems if~$\psi$ has large modulation. For the second term, one needs to distinguish
frequency interactions as before:
\begin{align}
i\partial_{\alpha}[\psi^{\beta} I A_{\beta}] &= i\sum_{k}
 P_k\partial_{\alpha}\psi^\beta \,\del_j^{-1} I
P_{<k-5}\calN_{\beta j}(\psi,\psi) +i\sum_{k}  P_k\psi^\beta\,
\,\del_j^{-1} I P_{<k-5}\partial_{\alpha}\calN_{\beta j}(\psi,\psi)  \label{eq:IIIHL} \\
& +i \sum_{k}\partial_{\alpha}P_k[P_{>k}R^\beta\psi\, \del_j^{-1} I
P_{>k}\calN_{\beta j}(\psi,\psi)]  +i
\sum_{k}\partial_{\alpha}P_k[P_{>k}\chi^\beta \,\del_j^{-1} I
P_{>k}\calN_{\beta j}(\psi,\psi)]\label{eq:IIIHH} \\
& +i\sum_{k}
\partial_\alpha [ P_{<k-5}R^\beta\psi\, \del_j^{-1} I
P_{k}\calN_{\beta j}(\psi,\psi)] + i\sum_{k}
\partial_\alpha  [P_{<k-5}\chi^\beta \,\del_j^{-1} I
P_{k}\calN_{\beta j}(\psi,\psi)]  \label{eq:IIILH}
\end{align}
The $\chi$-terms need to be expanded further, see below,
whereas the $\psi$-terms in~\eqref{eq:IIIHH} and~\eqref{eq:IIILH}
are estimated as they are. The second term on right-hand side of~\eqref{eq:IIIHL}
is expanded by means of the Hodge decomposition:
\begin{align}
 i\sum_{k} P_k\psi^\beta\,
\,\del_j^{-1} I P_{<k-5}\partial_{\alpha}\calN_{\beta j}(\psi,\psi) &= i\sum_{k} P_k R^\beta \psi\,
\,\del_j^{-1} I P_{<k-5}\partial_{\alpha}\calN_{\beta j}(\psi,\psi) \label{eq:IIIHLpsi} \\ &+ i\sum_{k} P_k\chi^\beta\,
\,\del_j^{-1} I P_{<k-5}\partial_{\alpha}\calN_{\beta j}(\psi,\psi) \label{eq:IIIHLchi}
\end{align}
The trilinear estimates of Section~\ref{sec:trilin} cover \eqref{eq:IIIHLpsi}, and~\eqref{eq:IIIHLchi} is handled
below]. Finally, the first term on the right-hand side of~\eqref{eq:IIIHL} is rewritten by means of~\eqref{eq:psisys1}:
\begin{equation}
 \label{eq:IIIHLswitch}
i\sum_{k}
P_k\partial_{\alpha}\psi^\beta \,\del_j^{-1} I
P_{<k-5}\calN_{\beta j}(\psi,\psi) = i\sum_{k}
P_k\partial^{\beta}\psi_\alpha  \,\del_j^{-1} I
P_{<k-5}\calN_{\beta j}(\psi,\psi) + \mathtt{quintilinear\  terms}
\end{equation}
where the quintilinear terms arise by using the curl identity for
$\partial_{\alpha}\psi^\beta-\partial^{\beta}\psi_\alpha$ into this expression. Note that we have
switched the derivatives $\del_\alpha$ and $\del_\beta$.

\noindent  We still have to explain how to deal with the higher
order terms involving at least one factor of~$\chi$.
\\

{\bf{Higher order errors.}}
\\
Note that these arise in two ways: first, we generate errors by replacing the Gauge potential $A_{\beta}$ in
\[
i\partial^\beta[\psi_\alpha A_{\beta}]
\]
by a $\calN_{\beta j}(\psi,\psi)$ null-form, and similarly for the remaining types of terms
\[
 i\partial^\beta[\psi_\beta A_{\alpha}],\quad i\partial_\alpha[\psi^\beta A_{\beta}]
\]
We shall call the higher order terms generated by this process (and later further Hodge decompositions applied to them) {\em of the first type or kind}.
\\
Second, we generate errors of the schematic form
\[
\chi \nabla^{-1}IQ_{\beta j}(\psi, \psi),
\]
and we call these together with all the terms generated by them upon applying further Hodge decompositions
{\em of the second type or kind}.  For simplicity, we omit frequency localizations in the ensuing discussion.
\\
Considering the errors of the first kind, these are of the schematic form
\[
\nabla_{x,t}[\psi\nabla^{-1}[\chi\psi]],\quad
\nabla_{x,t}[\psi\nabla^{-1}[\chi\chi]],
\]
where we recall from the very beginning, section 1, that
\[
\chi=\nabla^{-1}[\psi\nabla^{-1}(\psi^{2})],
\]
whence the above terms may be thought of a quintilinear and septilinear. Now as they are written,
we cannot yet quite estimate these expressions, and we need to introduce more null-structure, by expanding the $\nabla^{-1}(\psi^2)$ in
\[
\chi=\nabla^{-1}[\psi\nabla^{-1}(\psi^{2})],
\]
into a $Q_{\nu j}$-null-form as well as even higher order error
terms.  To keep track of things we associate an expansion graph,
i.e., a simple binary tree with the expressions generated: represent
the original terms
\[
\nabla_{x,t}[\psi A_{\beta}]
\]
by a simple node, and whenever we replace one of the factors in the (schematically written)
\[
A_\beta=\nabla^{-1}(\psi^{2})
\]
by the corresponding $\chi$, we draw a downward edge pointing left or right corresponding to which factor we replace.
\begin{figure*}[ht]
\begin{center}
\centerline{\hbox{\vbox{ \epsfxsize= 8.0 truecm \epsfysize=4.5
truecm \epsfbox{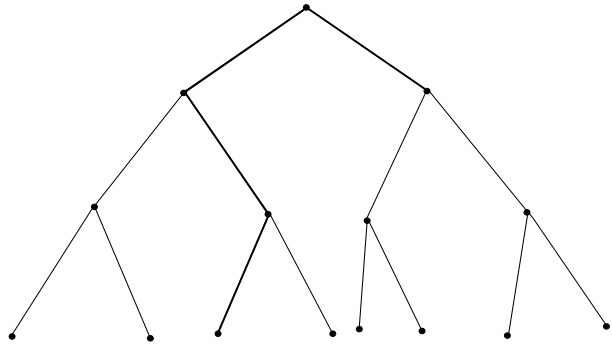}}}} \caption{An example of an expansion
graph}
\end{center}
\end{figure*}
We  can now exactly specify the full expansion of the higher order errors of first type:
\\
{\bf{Precise description of expansion for errors of first type:}}
\\
{\em keep applying Hodge decompositions to the inner $\nabla^{-1}(\psi^{2})$ in all factors
\[
\chi=\nabla^{-1}[\psi\nabla^{-1}(\psi^{2})],
\]
generated until the associated expansion graph has a directed subgraph of length four. Then the process stops.}
Note that formally, the terms thereby generated are up to the 11th degree in $\psi$.
\\

Next, we apply a similar process to the errors of the second type. We represent the first such error, schematically given by
\[
\nabla_{x,t}[\chi\nabla^{-1}IQ_{\nu j}(\psi, \psi)]
\]
by a simple node, and whenever we apply a Hodge decomposition to one of the factors  of $\nabla^{-1}(\psi^2)$ in
\[
\chi=\nabla^{-1}(\psi\nabla^{-1}(\psi^2))
\]
we draw a downward edge pointing left or right, thereby generating an associated expansion graph. Then we have
\\
{\bf{Precise description of expansion for errors of second type:}}
\\
{\em Keep applying Hodge decompositions as above until the associated expansion graph has a directed subgraph of
length three. Then the process stops.} Again we generated a list of errors of degree of multilinearity up to order 11 in $\psi$.
\\
\\

To summarize this discussion, we have now recast our system of equations in the form
\[
 \Box \psi_\alpha=\sum_{i=1}^{5}F_{\alpha}^{2i+1}
\]
where the superscript indicates the degree of multilinearity of the corresponding term in $\psi$, and the
{\em {leading cubic terms}} $F^{3}_{\alpha}$ can be expressed as
\begin{equation}\label{eq:Fthree}
\begin{aligned}
F_\alpha^{3}=&i\partial^{\beta}[\psi_\alpha \, I^c \del_j^{-1}\calN_{\beta
j}(\psi,\psi)] + i\partial^{\beta}[\psi_\alpha \, I
\del_j^{-1}\calN_{\beta j}(\psi,\psi)]-i\partial^{\beta}[\psi_\beta \,\del_j^{-1} I^c
\calN_{\alpha j}(\psi,\psi)] \\
& -i\sum_{k} P_k R_\beta \psi\,
\,\del_j^{-1} I P_{<k-5}\partial^{\beta}\calN_{\alpha j}(\psi,\psi)-i \sum_{k}\partial^{\beta}P_k[P_{>k}R_\beta\psi\, \del_j^{-1} I
P_{>k}\calN_{\alpha j}(\psi,\psi)]\\
&- i\sum_{k}
\partial^{\beta} [P_{<k-5}R_\beta\psi\, \del_j^{-1} I
P_{k}\calN_{\alpha j}(\psi,\psi)]
+i\partial_{\alpha}[\psi^{\beta} I^c \del_j^{-1}\calN_{\beta
j}(\psi,\psi)]\\
&+i\sum_{k}
 P_k\partial^{\beta}\psi_\alpha \,\del_j^{-1} I
P_{<k-5}\calN_{\beta j}(\psi,\psi) +i\sum_{k}  P_kR^\beta \psi\,
\,\del_j^{-1} I P_{<k-5}\partial_{\alpha}\calN_{\beta j}(\psi,\psi)\\
&+i \sum_{k}\partial_{\alpha}P_k[P_{>k}R^\beta\psi\, \del_j^{-1} I
P_{>k}\calN_{\beta j}(\psi,\psi)]+i\sum_{k}
\partial_\alpha [ P_{<k-5}R^\beta\psi\, \del_j^{-1} I
P_{k}\calN_{\beta j}(\psi,\psi)]
\end{aligned}
\end{equation}
Here it is very important to note that the second as well as the ninth term on the right contribute a {\em{magnetic potential interaction term}} of the form
\begin{equation}\label{eq:magnetic}
2i\sum_{k}
 P_k\partial^{\beta}\psi_\alpha \,\del_j^{-1} I
P_{<k-5}\calN_{\beta j}(\psi,\psi),
\end{equation}
the idea being that we interpret the low-frequency term $\del_j^{-1}
I P_{<k-5}\calN_{\beta j}(\psi,\psi)$ as a magnetic gauge potential.
The main issue here is that these high-low interactions cannot be
made small in general which creates problems for a bootstrap
argument. Hence, in order to prove the core perturbative results in
Section~\ref{sec:BG} we shall have to move these interaction terms
to the left-hand side, i.e., build them into the linear operator.
For later reference, we shall denote by $F_{\alpha}^{3k}$, $k= 1, 2,
3$, those trilinear terms contributed by the first, second or third
term in \eqref{eq:psiwave2}; thus for example, we write
\begin{equation}\begin{aligned}\label{eq:F32}
 F_{\alpha}^{32} &=i\partial^{\beta}[\psi_\beta \,\del_j^{-1} I^c
\calN_{\alpha j}(\psi,\psi)]
-i \sum_{k}\partial^{\beta}P_k[P_{>k}R_\beta\psi\, \del_j^{-1} I
P_{>k}\calN_{\alpha j}(\psi,\psi)]\\
&\quad -i\sum_{k} P_k R_\beta \psi\,
\,\del_j^{-1} I P_{<k-5}\partial^{\beta}\calN_{\alpha j}(\psi,\psi) - i\sum_{k}
\partial^{\beta} [P_{<k-5}R_\beta\psi\, \del_j^{-1} I
P_{k}\calN_{\alpha j}(\psi,\psi)]
\end{aligned}
\end{equation}
Furthermore, we denote by
\[
 F_{\alpha}^{3k}(\psi_{1}, \psi_{2}, \psi_{3})
\]
the corresponding multilinear expressions. We also introduce frequency localized versions
\[
 F_{\alpha}^{3k}(\psi_{1}; P_{<\ell};\psi_{2}, \psi_{3})
\]
in which one includes a cutoff $P_{<\ell}$ in front of all instances
of $\calN_{\alpha j}(\psi_2,\psi_3)$, and similarly for other
multipliers $P_{\leq \ell}$ etc.

\section{Bilinear estimates involving  $S$ and $N$ spaces}
\label{sec:bilin}

In this section we develop some of the required bilinear bounds.
First, we present some bounds from $S\times S$~into~$L^2_{tx}$, in
particular one which involves a gain in the high-high case and which
does  not appear in~\cite{T1} or~\cite{Krieger}, see
Lemma~\ref{lem:bilinbasic} below. This result allows for better
control on products $\phi_1\,\phi_2$ of $S$-waves and will be most
useful in the trilinear case. In addition, as in the aforementioned
references we consider the case of $\phi_1\in S$ and $\phi_2\in N$.
This section concludes with bilinear estimates for null-forms.

\subsection{Basic $L^2$-bounds}
To begin with, we present the following geometric lemma for cones,
see~\cite{T2} for a similar result. It will be used repeatedly.

\begin{lemma}
  \label{lem:cone} Suppose $\phi_1, \phi_2$ are such that
  \[
\supp(\wh{\phi_j})\subset \big\{(\xi,\tau)\:|\: |\xi|\sim 2^{k_j},\;
\big||\xi|-|\tau|\big|\sim 2^{\ell_j}\big\}
  \]
for $j=1,2$. Let $\ell_0,k_0\in\Z$ and assume that there exists
$j_0\in\{0,1,2\}$ so that
\begin{equation}\label{eq:imbalance}
\ell_{j_0}>\ell_j+C \qquad\forall\; j\in\{0,1,2\}\setminus \{j_0\}
\end{equation}
Then there is the following dichotomy:

(A) If $k_0= k_{\max}+O(1)$, then
\begin{align}
P_{k_0} Q_{\ell_0} \big( \phi_1 \phi_2  \big) &= P_{k_0} Q_{\ell_0}
\Big( \sum_{\kappa_1,\kappa_2} P_{k_1,\kappa_1} \phi_1 \cdot
P_{k_2,\kappa_2} \phi_2 \Big)
\end{align}
where $\kappa_1,\kappa_2$ are caps of size $C^{-1} r$ and separation
$\dist(\kappa_1,\kappa_2)\sim r$ with
\[
r:=  2^{(\ell_{\max}-k_{\min})/2}
\]
In particular, $\ell_{\max}\le k_{\min}+O(1)$.

\medskip

(B) If $k_0<k_{\max}-C$, then
\begin{align}
P_{k_0} Q_{\ell_0} \big( \phi_1 \phi_2  \big) &= \sum_{\eps=\pm}
P_{k_0} Q_{\ell_0} \Big( \sum_{\kappa} P_{k_1,\kappa}
\phi_1^{(\eps)}
\cdot P_{k_2,-\kappa} \phi_2^{(\eps)} \Big) \label{eq:plussum'}\\
&+\sum_{\eps=\pm} P_{k_0} Q_{\ell_0} \Big( \sum_{\kappa_1,\kappa_2}
P_{k_1,\kappa_1} \phi_1^{(\eps)} \cdot P_{k_2,-\kappa_2}
\phi_2^{(-\eps)} \Big) \label{eq:minsum'}
\end{align}
the sum in~\eqref{eq:minsum'} runs over caps of size $C^{-1} r$ with
\[
r:= 2^{k_0-k_{\max}} 2^{(\ell_{\max}-k_{\min})/2}
\]
and with separation $\dist(\kappa_1,\kappa_2)\sim r$, whereas the
sum in~\eqref{eq:plussum'} runs over caps of size~$r'$ where
$2^{k_0-k_{\max}}\le r'\le 1$ is arbitrary but fixed. The
sum~\eqref{eq:plussum'} is empty if $\ell_{\max}<k_{\max}-C$ and
\eqref{eq:minsum'} is nonzero only if $\ell_{\max}\le
k_{\min}+O(1)$. Finally, if \eqref{eq:imbalance} fails, then the
same representations hold provided $r\les1$ and one replaces
$\dist(\kappa_1,\kappa_2)\sim r$ with $\dist(\kappa_1,\kappa_2)\les
r$.
\end{lemma}
\begin{proof}
We consider first the $(++)$ and $(--)$ cases, i.e., when
$\tau_1,\tau_2$ have the same sign. Then
\begin{equation}\label{eq:xi1etc}
|\xi_1|+|\xi_2|-|\xi_1+\xi_2|\sim 2^{\ell_{\max}}
\end{equation}
whence
\[
(|\xi_1|+|\xi_2|)^2 - |\xi_1+\xi_2|^2 \sim 2^{\ell_{\max}+k_{\max}}
\]
and thus
\[
\sangle(\xi_1,\xi_2)\sim 2^{(\ell_{\max}+k_{\max}-k_1-k_2)/2}
\]
Now assume further that $k_0=k_{\max}+O(1)$. Then it follows that
\[
\sangle(\xi_1,\xi_2)\sim 2^{(\ell_{\max}-k_{\min})/2}
\]
If on the other hand $k_0<k_{\max}-C$, then
$k_1=k_2+O(1)=k_{\max}+O(1)$ and from~\eqref{eq:xi1etc},
$\ell_{\max}=k_{\max}+O(1)$. Furthermore, $\xi_2=-\xi_1+O(2^{k_0})$
implies that
\[
|\sangle(\xi_1,-\xi_2)|\sim \frac{|\xi_1\wedge
\xi_2|}{|\xi_1||\xi_2|} = O(2^{k_0-k_{\max}})
\]
Next, consider the $(+-)$ or $(-+)$ cases.
\begin{figure*}[ht]
\begin{center}
\centerline{\hbox{\vbox{ \epsfxsize= 7.0 truecm \epsfysize=6.5
truecm \epsfbox{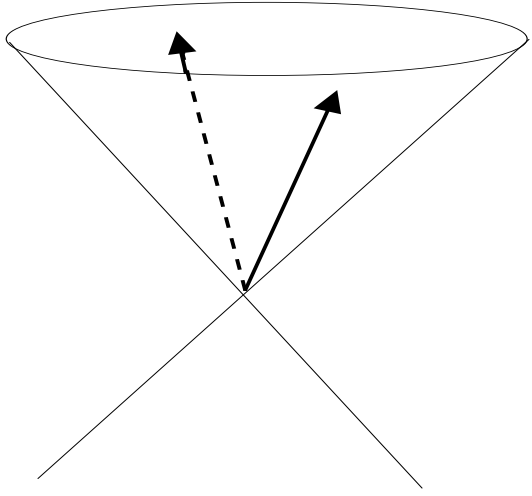}}}} \caption{Opposing $(++)$ waves}
\end{center}
\end{figure*}
Then
\begin{equation}\label{eq:xi1etc'}
|\xi_1+\xi_2|-||\xi_1|-|\xi_2||\sim 2^{\ell_{\max}}
\end{equation}
which implies that
\[
|\xi_1+\xi_2|^2-||\xi_1|-|\xi_2||^2\sim 2^{\ell_{\max}}
(|\xi_1+\xi_2|+||\xi_1|-|\xi_2||)
\] or equivalently,
\begin{equation}
  \label{eq:ang2}
2^{k_1+k_2} \sangle^2(\xi_1,-\xi_2) \sim 2^{\ell_{\max}+k_0}
  \end{equation}
If $k_0=k_{\max}+O(1)$, then
\[
\sangle(\xi_1,-\xi_2)\sim 2^{(\ell_{\max}-k_{\min})/2}
\]
If, on the other hand, $k_0\le k_{\max}-C$, then
\[
\sangle(\xi_1,-\xi_2)\sim 2^{k_0-k_{\max}}
2^{(\ell_{\max}-k_{\min})/2}
\]
and we are done. While it is clear that $\ell_{\max}\le k_{\min}+O(1)$ if $k_0=k_{\max}+O(1)$,
some proof is needed in case $k_0<k_{\max}-C$. Thus, suppose $|\xi_1|\ge|\xi_2|$ whence
\[
|\xi_1+\xi_2|-|\xi_1|+|\xi_2|\sim 2^{\ell_{\max}}
\]
which implies that
\[
2^{k_0+k_1} \sangle^2(\xi_1+\xi_2,-\xi_2)\sim 2^{\ell_{\max}+k_{\max}}
\]
since $2^{k_0+k_1}\sim 2^{k_{\min}+k_{\max}}$, the claim follows.

Finally, if \eqref{eq:imbalance} fails, then \eqref{eq:xi1etc} turns into
\[
|\xi_1|+|\xi_2|-|\xi_1+\xi_2|\les 2^{\ell_{\max}}
\]
which then leads to the claimed loss of separation between the sectors. However, their maximal distances are
controlled by the same quantities as before.
\end{proof}

The special appearance of \eqref{eq:plussum'} derives from the contributions of waves
which lie on opposing sides of the light-cone. In fact,
Figure~3 shows two vectors on the same half (i.e., $\tau>0$) but opposing sides of the light
cone. They add up to produce a wave of small frequency but large modulation, as described by~\eqref{eq:plussum'}.
This is the mechanism by which nonlinearities can turn free waves into ``elliptic objects''. This phrase
refers to functions whose Fourier support has large separation from the characteristic variety of~$\Box$.
Also, following Tao we refer to \eqref{eq:imbalance} as the {\em modulation imbalanced case}, whereas its opposite is the
{\em modulation balanced case}.

\begin{exse}
 \label{ex:cone} Lemma~\ref{lem:cone} is optimal in the following sense:
\begin{itemize}
 \item Given $\ell_0\le k_0\le -10$  there exist $\xi_1,\xi_2\in\R^n$ with $1\le |\xi_1|,|\xi_2|\le2$,
$\sangle(\xi_1,\xi_2)\sim 2^{(\ell_0+k_0)/2}$ and such that
\begin{align*}
 |\xi_1+\xi_2|-||\xi_1|-|\xi_2|| &\sim 2^{\ell_0} \\
|\xi_1+\xi_2| &\sim 2^{k_0}
\end{align*}
\item Given $\ell_0\le k_1\le -10$ there exist $\xi_1, \xi_2\in\R^n$ with $2^{k_1-1}\le |\xi_1|\le2^{k_1}$, $1\le |\xi_1|,|\xi_2|\le2$
and $\sangle(\xi_1,\xi_2)\sim 2^{(\ell_0-k_1)/2}$ and
so that \[|\xi_1|+|\xi_2|-|\xi_1+\xi_2|\sim 2^{\ell_0}.\]
\end{itemize}
\end{exse}

Our  immediate goal now is the proof of Lemma~\ref{lem:bilin3}.
It is important to note
that the improvement of~$2^{\frac{k}{2}}$ over~\eqref{eq:bilin2}
which is obtained in Lemma~\ref{lem:bilin3} coincides with the gain
for the case of free waves.
In order to accomplish this, we require three preparatory lemmas, all of which are well-known. The first is
Mockenhaupt's ``square function estimate'' (more precisely, its geometric content), see~\cite{Gerd}, \cite{MSS}.  Recall that
$\Theta=\sign(\tau)\hat{\xi}$.

\begin{lemma}
  \label{lem:gerd}
Let $\kappa,\tilde\kappa\in\calC_{\ell}$ with
$\dist(\kappa,\tilde\kappa)\sim|\kappa|\ll 1$ and suppose that
$\calF_i\subset \calC_{\ell_i}$ for $i=1,2$ are partitions of
$\kappa$ and $\tilde\kappa$, respectively, by
 pairwise disjoint caps.
Further, let $r\in(0,1)$, $\mu\in(1,2)$,  and  define for any cap $\kappa'\subset S^1$
\begin{align*}
 \calT_{\kappa',\mu,r} := \{ (\tau,\xi)\::\: ||\xi|-\mu|\le r,\; \Theta
 \in\kappa', \; ||\tau|-|\xi||\le |\kappa'|^2\}
\end{align*}
Set $M_i:=\# \calF_i$.  Then
\begin{equation}
  \label{eq:gerd}
\sup_{\mu_1,\mu_2\sim1} \Big\| \sum_{\kappa_1\in\calF_1}
  \sum_{\kappa_2\in\calF_2} \chi_{\calT_{\kappa_1,\mu_1,r} +
  \calT_{\kappa_2,\mu_2,r}} \Big\|_{L^\infty(\R^3)} \le C\max(1,r(M_1+M_2))
\end{equation}
where $C$ is some absolute constant.
\end{lemma}
\begin{proof}
Fix $r\in(0,1)$, and $\mu_1,\mu_2\sim1$.
 Applying a Lorentz transform, one may assume that $\ell=-10$, say.
 Also, suppose without loss of generality that $\ell_1\le \ell_2$ whence $M_1\ge M_2$.
We first consider the case where $rM_1\ge1$.
Fix $(\tau,\xi)\in\R^3$ such that\footnote{The
 $\frac12$-
 factor is a convenient modification that can be made due to
 scaling.}
\[
\sum_{\kappa_1\in\calF_1}
  \sum_{\kappa_2\in\calF_2} \chi_{\frac12(\calT_{\kappa_1,\mu_1,r} +
  \calT_{\kappa_2,\mu_2,r})} (\tau,\xi)\ge1
\]
Suppose $\calT_{\kappa_1,\mu_1,r}$ with $\kappa_1\in\calF_1$
contributes to the sum on the left-hand side. Define a mirror-image
$\calT_{\kappa_1,\mu_1,r}^*$ of $\calT_{\kappa_1,\mu_1,r}$ by
reflecting  $\calT_{\kappa_1,\mu_1,r}$ about the point~$(\tau,\xi)$.
Due to $\ell=-10$ and the dimensions of the tubes $\calT$, the
mirror images of all
$\{\calT_{\kappa_1,\mu_1,r}\}_{\kappa_1\in\calF_1}$ have uniformly
bounded overlap.  The same applies with the role of $\calF_1$ and
$\calF_2$ reversed. In conclusion, each  $\calT_{\kappa_1,\mu_1,r}$
can pair up with at most $O(1)$-many $\calT_{\kappa_2,\mu_2,r}$ so
as to give a contribution to~\eqref{eq:gerd}, whence the bound
of~$M_1$ for~\eqref{eq:gerd}. To obtain the factor~$r$ improvement,
we further note that due to fixed $\mu_1$ and~$\mu_2$, only those
contributions to~\eqref{eq:gerd} need to be counted which derive
from pairs  $(\calT_{\kappa_1,\mu_1,r}, \calT_{\kappa_2,\mu_2,r})$
which lie in fixed cylinders $||\xi_i|-\mu_i|<r$, $i=1,2$. In terms
of equations, we are given $(\sigma,\zeta)\in\R^3$ and we need to
consider the sets of $(\tau_i,r_i\omega_i)$, $i=1,2$ with
$\omega_i\in S^1$ satisfying the transversality condition
$\sangle(\omega_1,\omega_2)\in[\frac{1}{100},\frac{1}{50}]$, say,
and such that
\begin{align*}
  \tau_1+\tau_2&= \sigma,\quad r_1\omega_1+r_2\omega_2= \zeta \\
 |r_1-\mu_1| &< r, \quad |r_2-\mu_2|<r\\
 ||\tau_1|-r_1| &< 2^{2\ell_1},\quad ||\tau_2|-r_2| < 2^{2\ell_2}
\end{align*}
It follows from the second, third, and fourth conditions that
\[
\mu_1\omega_1+\mu_2\omega_2= \zeta +O(r)
\]
and  since the circular arcs containing $\omega_1$ and $\omega_2$
are transverse to each other, they must be of lengths~$\les r$.
Consequently, we can only count tubes which correspond to an
$r\times r$ disk on the light-cone and of those there are at most
$rM_1$-many. In case $rM_1\le1$, then the number of the allowed
pairs is~$\les1$ in light of this construction and we are done.
\end{proof}

Next, we present a standard bilinear $L^2$ bound for free waves.

\begin{lemma}
 \label{lem:free_bilin}
Let $\kappa,\tilde\kappa\in\calC_{\ell}$ with $\dist(\kappa,\tilde\kappa)\sim|\kappa|:=\beta$
and suppose $\kappa_1\subset\kappa$, $\kappa_2\subset\tilde\kappa$ are arbitrary caps.
Let $r\in(0,1)$ and $\mu_1,\mu_2\sim1$. Then
\begin{equation}
 \label{eq:wave_est}
 \|e^{it|\nabla|}f_1\,  e^{\pm it|\nabla|}f_2\|_{\Ltwotx}  \les  \beta^{-1}
\sqrt{\min\big(r\beta, |\kappa_1|, |\kappa_2|  \big)}\; \|f_1\|_2\|f_2\|_2
\end{equation}
provided
\begin{align*}
 \supp(\wh f_1) &\subset \{ \xi\in\R^2\::\: \hat{\xi}\in\kappa_1, \, ||\xi|-\mu_1|\les r\}\\
  \supp(\wh f_2) &\subset \{ \xi\in\R^2\::\: \hat{\xi}\in \pm\kappa_2, \, ||\xi|-\mu_2|\les r\}\\
\end{align*}
and the sign in the last sign is chosen to be the same as in~\eqref{eq:wave_est}.
\end{lemma}
\begin{proof}
 The proof reduces to the following well-known property of convolutions:
suppose
\begin{align*}
 \Gamma_1 &:=  \{ (|\xi|,\xi)\in\R^3\::\: \hat{\xi}\in\kappa_1, \, ||\xi|-\mu_1|\les r\} \\
 \Gamma_2 &:=  \{ (\pm |\xi|,\xi)\in\R^3\::\: \hat{\xi}\in\pm\kappa_2, \, ||\xi|-\mu_2|\les r\}
\end{align*}
Note that $\angle(\xi,\pm \eta)\gtrsim \beta$ for any
$(|\xi|,\xi)\in\Gamma_1$ and $(\pm|\eta|,\eta)\in \Gamma_2$. Then
\begin{equation}\label{eq:mesconv}
 \| f\sigma_{\Gamma_1}\ast g\sigma_{\Gamma_2} \|_{L^2(\R^{3})}
 \les \beta^{-1}
\sqrt{\min\big(r\beta, |\kappa_1|, |\kappa_2|  \big)} \;
\|f\|_{L^2(d\sigma_{\Gamma_1})}
 \|g\|_{L^2(d\sigma_{\Gamma_2})}
\end{equation}
where $\sigma_{\Gamma_1}$ and $\sigma_{\Gamma_2}$ are the lifts of
the measure in~$\R^2$ to the sectors $\Gamma_1,\Gamma_2$ on the
light-cones.
To prove~\eqref{eq:mesconv}, interpolate between $L^1$ and
$L^\infty$. On $L^1$ we have the standard fact that
$\|\mu\ast\nu\|\le \|\mu\|\|\nu\|$ for measures and their total
variation norms. This fact does not use the angular separation of
the supports nor their sizes. On $L^\infty$, however, this separation and size are crucial and
yield
\[\| f\sigma_{\Gamma_1}\ast g\sigma_{\Gamma_2}
\|_{L^\infty(\R^{3})} \les \beta^{-1} \min(r,|\kappa_1|\beta^{-1})
\|f\|_{L^\infty(d\sigma_{\Gamma_1})}\|g\|_{L^\infty(d\sigma_{\Gamma_2})}
\]
assuming as we may that $|\kappa_1|\le |\kappa_2|$.
To obtain this bound, consider $\delta$-neighborhoods of $\Gamma_1$
and $\Gamma_2$, respectively. In other words, replace $d\sigma_1$ by
\[
d\wt\sigma^{(\delta)}_j:=\delta^{-1}\chi_{[\dist((\xi,\tau),\Gamma_j)<\delta]}d\xi
d\tau
\]
for small $\delta>0$ and observe that
\begin{equation}\label{eq:Linftymeas} \limsup_{\delta\to0+}\|
d\wt\sigma^{(\delta)}_1\ast
d\wt\sigma^{(\delta)}_2\|_{L^\infty_{\xi,\tau}} \les
\beta^{-1}\min(r,|\kappa_1|\beta^{-1})
\end{equation}
by elementary geometry. To pass
from~\eqref{eq:mesconv}  to estimates for the
wave equation  use Plancherel's theorem.
\end{proof}

We can now state the aforementioned improved bilinear $L^2$ bound.
The norm $\trip\cdot\trip$ is the one from~\eqref{eq:trip_def}.

\begin{lemma}
 \label{lem:bilin3}  Let
$\phi_i$  be adapted to~$k_i$ for $i=1,2$.  Assume further that we
are in the high-high case $k_1=k_2+O(1)$ and that $\phi_i=Q_{\le j+k-2k_1-C}\phi_i$ for $i=1,2$. Then
\begin{equation}
 \label{eq:bilin3}
\|P_k Q_j (\phi_1\,\phi_2)\|_{\Ltwotx} \les 2^{\frac{k_1}{2}} 2^{\frac{k-j}{4}} \trip\phi_1\trip_{S[k_1]} \trip\phi_2\trip_{S[k_2]}
\end{equation}
for any $j\le k\le k_1+O(1)$. Moreover, in the same range of $j$,
\begin{equation}
 \label{eq:bilin3null}
\|P_k Q_j (R_\alpha \phi_1 R_\beta \phi_2 - R_\beta \phi_1 R_\alpha
\phi_2 )\|_{\Ltwotx} \les 2^{\frac{k_1}{2}} 2^{\frac{3k+j}{4}}
\trip\phi_1\trip_{S[k_1]} \trip\phi_2\trip_{S[k_2]}
\end{equation}
for any $\alpha,\beta=0,1,2$.
\end{lemma}
\begin{proof}
We assume that $k_1=k_2+O(1)=0$. At first, we also assume that $k\le -C$ so as to exclude the
opposing $(++)$ and $(--)$ waves in Lemma~\ref{lem:cone}.
We need to prove that
\begin{align}
 \|P_k (\phi_1\phi_2)\|_{\Ltwotx} &\les 2^{\frac{k-j}{4}} \big( \|(f_1,g_1)\|_{L^2 \times  \dot H^{-1}} + \|F_1\|_{N[k_1]}\big)
\cdot \label{eq:F1F2} \big( \|(f_2,g_2)\|_{L^2 \times  \dot H^{-1}}
+ \|F_2\|_{N[k_2]}\big)
\end{align} for any $k_i$-adapted Schwartz functions $f_i,g_i,F_i$, $i=1,2$ and
\begin{equation}\label{eq:phiirep}
 \phi_i(t) = \cos(t|\nabla|) f_i + \frac{\sin(t|\nabla|)}{|\nabla|} g_i + \int_0^t  \frac{\sin((t-s)|\nabla|)}{|\nabla|} F_i(s)\, ds
\end{equation}
We reduce this to three cases:
\begin{align}
 \|P_k Q_j (\phi_1\phi_2)\|_{\Ltwotx} &\les 2^{\frac{k-j}{4}}  \|(f_1,g_1)\|_{L^2 \times  \dot H^{-1}}
  \|(f_2,g_2)\|_{L^2 \times  \dot H^{-1}} \label{eq:freefree} \\
\|P_k Q_j (\phi_1\phi_2)\|_{\Ltwotx} &\les 2^{\frac{k-j}{4}}  \|(f_1,g_1)\|_{L^2 \times  \dot H^{-1}}  \|F_2 \|_{N[k_2]} \label{eq:freeinh} \\
\|P_k Q_j (\phi_1\phi_2)\|_{\Ltwotx} &\les 2^{\frac{k-j}{4}}  \|F_1\|_{N[k_1]} \|F_2\|_{N[k_2]} \label{eq:inhinh}
\end{align}
where the absence of terms on the right-hand side implies that the corresponding functions are zero (thus, $F_1=F_2=0$ in~\eqref{eq:freefree}
etc.)  We begin with~\eqref{eq:freefree} which follows easily from Lemma~\ref{lem:gerd}. To see this,
we decompose $\phi_i$ into caps of size $\ell=(j+k)/2$ as in Lemma~\ref{lem:cone}. Adopting the convention that
$\kappa_1\sim\kappa_2$ means that $\dist(\kappa_1,\kappa_2)\sim 2^\ell$, and setting $g_1=g_2=0$ for simplicity,
one has\footnote{Recall our convention about $P_{k_i,\kappa}$ which
takes the sign of~$\tau$ into account.}
\begin{align}
 &\|P_k Q_j (\phi_1\phi_2)\|_{\Ltwotx} \les \sum_{\kappa_1\sim\kappa_2\in\calC_{\ell}}
\|P_k Q_j (P_{k_1,\kappa_1} \phi_1 \, P_{k_2,\kappa_2} \phi_2)\|_{\Ltwotx} \nn \\
&\les \sum_{c\in\calD_{0,k}}\sum_{\kappa_1\sim\kappa_2\in\calC_{\ell}}    \|P_{k_1,\kappa_1}P_{c} \phi_1 \, P_{k_2,\kappa_2}P_{-c} \phi_2\|_{\Ltwotx} \nn \\
&\les \sum_{c\in\calD_{0,k}}\sum_{\kappa_1\sim\kappa_2\in\calC_{\ell}}  2^{\frac{k-\ell}{2}}
  \|P_{k_1,\kappa_1}P_{c} f_1\|_2 \| P_{k_2,\kappa_2}P_{-c} f_2\|_{\Ltwotx} \label{eq:scheiss23} \\
&\les 2^{\frac{k-\ell}{2}}   \|f_1\|_2 \|f_2\|_2 \nn
\end{align}
as needed. The estimate in~\eqref{eq:scheiss23} follows from~\eqref{eq:gerd} since $k\ge\ell$.

To prove~\eqref{eq:freeinh} and~\eqref{eq:inhinh} it will suffice as usual to assume that  $F_i$ are $N[k_i]$-atoms for~$i=1,2$.
In fact, if $F_2$ in~\eqref{eq:freeinh} is either and energy or an $\dot X^{s,b}$-atom, then one again reduces matters to
the free case. Consequently, we may restrict ourselves to~\eqref{eq:inhinh} when both $F_1$ and~$F_2$ are null-frame atoms.
Using Lemma~\ref{lem:Nsquare} to refine these null-frame atoms one can thus assume that
\begin{equation}\label{eq:refine}
 F_1= \sum_{\kappa'\in\calC_{\ell'}} F_{\kappa'},\quad  F_1= \sum_{\kappa''\in\calC_{\ell''}} \wt F_{\kappa''}
\end{equation}
where $\ell',\ell''\le \ell$. Again by Lemma~\ref{lem:Nsquare}, we can further assume that there exists a fixed $c\in\calD_{0,k}$
so that $P_c F_1=F_1$ and $P_{-c} F_2=F_2$.  Applying the same decomposition as in~\eqref{eq:scheiss23}, fix $\kappa_1\sim \kappa_2$.
In view of Lemma~\ref{lem:Null_rep},
\begin{align}
 P_{k_1,\kappa_1} \phi_1 &=   \Box^{-1} G_{\kappa_1} + \sum_{\substack{\kappa'\in\calC_{\ell'}\\ \kappa'\subset \kappa_1}}
 \int_{\R}
\big(\Psi^1_{\kappa',a} +  B_{\kappa',a}\, \Psi^2_{\kappa',a}\big)\, da \label{eq:scheiss27} \\
P_{k_2,\kappa_2} \phi_2 &=   \Box^{-1} \wt G_{\kappa_2} + \sum_{\substack{\kappa''\in\calC_{\ell''}\\ \kappa''\subset \kappa_2}}
 \int_{\R}
\big(\wt\Psi^1_{\kappa'',a} +  \wt B_{\kappa'',a}\, \wt\Psi^2_{\kappa'',a}\big)\, da  \label{eq:scheiss28}
\end{align}
where the functions on the right-hand side satisfy the bounds
specified in that lemma. Moreover, the Fourier supports of the
functions appearing inside the integral in \eqref{eq:scheiss27}
and~\eqref{eq:scheiss28} satisfy~\eqref{eq:PsiFsupp}, and they also
retain the $P_{c}$ and $P_{-c}$ localization property, respectively,
due to the fact that $k\ge\ell$.  We can ignore the terms involving
$G_{\kappa_1}$ and $\wt G_{\kappa_2}$ as they are reducible to free
waves. For simplicity, we also set $\Psi^1_{\kappa',a} =
\wt\Psi^1_{\kappa'',a}=0$. By Plancherel's theorem and
Lemma~\ref{lem:gerd},
\begin{align*}
 &\| P_{k_1,\kappa_1} \phi_1 P_{k_2,\kappa_2} \phi_2 \|_{\Ltwotx} \\
&\les  \sqrt{1+2^k(M_1+M_2)} \Bigg(\sum_{ \substack{\kappa'\in\calC_{\ell'}\\ \kappa'\subset \kappa_1}}
\sum_{\substack{\kappa''\in\calC_{\ell''}\\ \kappa''\subset \kappa_2}}
\|P_{k_1,\kappa'} \phi_1 P_{k_1,\kappa''} \phi_2  \|_{\Ltwotx}^2 \Bigg)^{\frac12}
\end{align*}
where  $M_1=2^{\ell-\ell'}, M_2=2^{\ell-\ell''}$.  On the other hand, Lemma~\ref{lem:free_bilin} implies that
\begin{align*}
 & \|P_{k_1,\kappa'} \phi_1 P_{k_1,\kappa''} \phi_2  \|_{\Ltwotx}  \\
&\les \int_{\R^2}  \| B_{\kappa',a}\, \Psi^2_{\kappa',a}
\wt B_{\kappa'',b}\, \wt\Psi^2_{\kappa'',b} \|_{\Ltwotx} \, dadb\\
&\les  \int_{\R^2}  \|   \Psi^2_{\kappa',a}
 \, \wt\Psi^2_{\kappa'',b} \|_{\Ltwotx}  \, dadb \\
&\les   2^{-\ell} \sqrt{\min(2^{k+\ell}, 2^{\ell'}, 2^{\ell''})}   \int_{\R^2}
 \| \Psi^2_{\kappa',a} \|_{\dot X_0^{0,\frac12,1}}
\, \| \wt\Psi^2_{\kappa'',b} \|_{\dot X_0^{0,\frac12,1}}\, dadb \\
&\les   2^{-\ell} \sqrt{\min(2^{k+\ell}, 2^{\ell'}, 2^{\ell''})}
 \| F_{\kappa'} \|_{\NF[\kappa']}
\, \| \wt F_{\kappa''} \|_{\NF[\kappa''] }
\end{align*}
One checks that
\[
 \sqrt{1+2^k(M_1+M_2)}\;  2^{-\ell} \sqrt{\min(2^{k+\ell}, 2^{\ell'}, 2^{\ell''})} \les 2^{\frac{k-\ell}{2}}
\]
whence
\[
 \| P_{k_1,\kappa_1} \phi_1 P_{k_2,\kappa_2} \phi_2 \|_{\Ltwotx} \les 2^{\frac{k-\ell}{2}}
\Bigg(\sum_{ \substack{\kappa'\in\calC_{\ell'}\\ \kappa'\subset \kappa_1}} \sum_{\substack{\kappa''\in\calC_{\ell''}\\
\kappa''\subset \kappa_2}}  \| F_{\kappa'} \|_{\NF[\kappa']}^2
\, \| \wt F_{\kappa''} \|_{\NF[\kappa''] }^2 \Bigg)^{\frac12}
\]
In conclusion,
\begin{align*}
 \|P_k Q_j (\phi_1\phi_2)\|_{\Ltwotx} &\les \sum_{\kappa_1\sim\kappa_2\in\calC_{\ell}}
  \| P_{k_1,\kappa_1} \phi_1 P_{k_2,\kappa_2} \phi_2 \|_{\Ltwotx} \\
&\les  2^{\frac{k-\ell}{2}} \sum_{\kappa_1\sim\kappa_2\in\calC_{\ell}}
\Bigg(\sum_{ \substack{\kappa'\in\calC_{\ell'}\\ \kappa'\subset \kappa_1}} \sum_{\substack{\kappa''\in\calC_{\ell''}\\
 \kappa''\subset \kappa_2}}  \| F_{\kappa'} \|_{\NF[\kappa']}^2
\, \| \wt F_{\kappa''} \|_{\NF[\kappa''] }^2 \Bigg)^{\frac12} \\
&\les 2^{\frac{k-\ell}{2}} \Big(\sum_{\kappa'\in\calC_{\ell'}} \|F_{\kappa'}\|_{\NF[\kappa']}^2\Big)^{\frac12}
\Big(\sum_{\kappa''\in\calC_{\ell''}} \|\wt F_{\kappa''}\|_{\NF[\kappa'']}^2\Big)^{\frac12}
\end{align*}
as desired.  This concludes the proof of~\eqref{eq:inhinh} for the case of null-frame atoms $F_1, F_2$.
As indicated, the other cases are easier since they can be reduced to free waves.

\noindent Finally, if $k=O(1)$, then the proof is easier. In fact, it follows
via a cap-decomposition from the basic bilinear
bound~\eqref{eq:bilin2}. We leave those details to the reader.

\noindent The second bound~\eqref{eq:bilin3null} follows by the same argument.
The only difference from~\eqref{eq:bilin3} lies with an additional
gain of~$2^\ell$ which is precisely the size of the angle in the
above decompositions into caps.
\end{proof}

Later, we shall require the following
  technical variant of the previous bound.

\begin{cor}
  \label{cor:squarebilin3}
   Under the assumptions of
  Lemma~\ref{lem:bilin3}, for any $j\le k\le k_1+O(1)$ and any $m_0\le -10$,
  \begin{align}
 \sum_{\substack{\kappa_1,\kappa_2\in\calC_{m_0}\\ \dist(\kappa_1,\kappa_2)\le 2^{m_0}}}
   \|P_k Q_j (P_{k_1,\kappa_1} \phi_1\, P_{k_2,\kappa_2}
\phi_2)\|_{\Ltwotx} &\les 2^{\frac{k_1}{2}} 2^{\frac{k-j}{4}}
\trip\phi_1\trip_{S[k_1]} \trip\phi_2\trip_{S[k_2]} \label{eq:sum12}  \\
\Big(\sum_{\kappa\in\calC_{m_0}} \|P_k Q_j (P_{k_1,\kappa}\phi_1\,
 \phi_2)\|_{\Ltwotx}^2 \Big)^{\frac12} &\les |m_0|\,
2^{\frac{k_1}{2}} 2^{\frac{k-j}{4}} \trip\phi_1\trip_{S[k_1]}
\trip\phi_2\trip_{S[k_2]} \label{eq:sum13}
\end{align}
Moreover, analogous bounds hold for the null form in~\eqref{eq:bilin3null} with an extra gain of~$2^{\frac{j+k}{2}}$.
Finally, the left-hand side in~\eqref{eq:sum12} vanishes unless
$j+k\le
  2m_0\le -100$.
\end{cor}
\begin{proof} The final statement here is  due to Lemma~\ref{lem:cone}.
Note that one cannot simply square sum the bounds of
Lemma~\ref{lem:bilin3} applied to $P_{k_1,\kappa_1}\phi_1$
and~$P_{k_2,\kappa_2}\phi_2$ due to the fact that~$\sum_{\kappa}
\trip P_{k,\kappa} \phi\trip_{S[k]}^2$ (or $\sum_{\kappa} \|
P_{k,\kappa} \phi\|_{S[k]}^2$ for that matter) cannot be controlled.
However, since we may assume that $\frac{j+k}{2}\le m_0$, the angular decomposition
induced by the frequency and modulation cutoffs $P_k Q_j$ is {\em
finer} than the one superimposed by $\kappa_1$ and~$\kappa_2$.
Inspection of the proof now reveals that either by orthogonality or
by organizing the finer caps into subsets of the
$\kappa_1,\kappa_2\in\calC_{m_0}$, and applying the Cauchy-Schwarz inequality
yields the stated bound.
For~\eqref{eq:sum13} one needs to  distinguish two cases: either $m_0\ge \frac{j+k}{2}$ or
not. In the former case, the decomposition into caps in~$\calC_{m_0}$ is coarser than
the one coming from Lemma~\ref{lem:cone} and one
can again argue by means of Cauchy-Schwarz as before. In the latter case, however,
we split the modulation of the first input as follows:
 $$Q_{<j+k-C}= Q_{<2m_0-C} + Q_{2m_0-C\le \cdot <j+k-C}$$
The contribution of $Q_{<2m_0-C}\phi_1 $ is handled exactly as in the Lemma~\ref{lem:bilin3}
since one may always refine the null-frame representation,
cf.~\eqref{eq:refine}.
On the other hand, $Q_{2m_0-C\le \cdot <j+k-C}\phi_1$ is controlled by means of Lemma~\ref{lem:incl_free}. More precisely,
for any $2m_0-C\le \ell <j+k-C$ one has
$
 Q_\ell \phi_1 = Q_\ell\Box^{-1} F_1,
$
see \eqref{eq:phiirep}.  Since \eqref{eq:Xsb_dom} implies that
\[
 \|Q_\ell \phi_1\|_{\dot X_0^{0,\frac12,\infty}}  = \| Q_\ell\Box^{-1} F_1\|_{\dot X_0^{0,\frac12,\infty}} \les
\| Q_\ell  F_1\|_{\dot X_0^{0,-\frac12,\infty}} \les \|F_1\|_{N[0]}
\]
one can reduce the contribution of $Q_\ell \phi_1$ to the case of free waves as in the proof
of Lemma~\ref{lem:bilin3}. Summing over all~$\ell$ in this range loses a factor of at most~$|m_0|$,
as claimed. Finally, the claim concerning the null-forms is immediate.
\end{proof}

\noindent
Removing the modulation restrictions on the inputs in Lemma~\ref{lem:bilin3} results in the following estimates.

\begin{lemma}
 \label{lem:bilinbasic} If $\phi_1$ and $\phi_2$ are adapted to~$k_1$ and~$k_2$, respectively,
then for $j\le k\le k_1+O(1)=k_2+O(1)$,
\begin{equation}\label{eq:phi12hh}
 \| P_k Q_j (\phi_1\phi_2)\|_{\Ltwotx} \les   2^{\frac{k-j}{4}} 2^{\frac{k_1}{2}} \trip\phi_1\trip_{S[k_1]} \trip\phi_2\trip_{S[k_2]}
\end{equation}
whereas for $j\le k_2\le k=k_1+O(1)$,
\begin{equation}\label{eq:phi12hl}
 \| P_k Q_j (\phi_1\phi_2)\|_{\Ltwotx} \les   2^{\frac{3k_2}{4}} 2^{-\frac{j}{4}} \|\phi_1\|_{S[k_1]} \|\phi_2\|_{S[k_2]}
\end{equation}
\end{lemma}
\begin{proof}
Consider the high-high case $j\le k\le k_1+O(1)=k_2+O(1)=0$. On the
one hand, there is  the bound
\begin{equation}\label{eq:tripuse}
 \|P_k Q_j (Q_{\le j+k-C}\phi_1 Q_{\le j+k-C}\phi_2)\|_{\Ltwotx} \les 2^{\frac{k-j}{4}} \trip
 \phi_1\trip_{S[k_1]}\,
 \trip\phi_2\trip_{S[k_2]}
\end{equation}
which is given by Lemma~\ref{lem:bilin3}. On the other hand, by the
improved Bernstein bound of Lemma~\ref{lem:KBern},
\begin{align}
 \|P_k Q_j (Q_{>j+k-C}\phi_1\cdot \phi_2)\|_{\Ltwotx} & \les
 2^{\frac{j-k}{4}\wedge0}2^k \|Q_{>j+k-C}\phi_1\cdot \phi_2\|_{L_t^2
 L^1_x} \nn
 \\
 &\les 2^{\frac{j-k}{4}}2^k \|Q_{>j+k-C}\phi_1\|_{\Ltwotx} \| \phi_2\|_{L_t^\infty
 L^2_x} \nn \\
&\les  2^{\frac{k-j}{4}} \|\phi_1\|_{S[k_1]} \|\phi_2\|_{S[k_2]}
\label{eq:highhigh2}
\end{align}
In the high-low case $j\le k_2\le k=k_1+O(1)=0$ consider the
following three subcases. First,
\begin{align*}
 \| P_k Q_j (Q_{<j-C} \phi_1 Q_{<j-C} \phi_2)\|_{\Ltwotx} &\les  2^{-\frac{j-k_2}{4}} 2^{\frac{k_2}{2}} \|\phi_1\|_{S[k_1]} \|\phi_2\|_{S[k_2]}
\end{align*}
by a decomposition into caps of size~$2^{\frac{j-k_2}{2}}$ and the
$L^2$-bilinear bound~\eqref{eq:bilin2}. Next, by the improved
Bernstein estimate Lemma~\ref{lem:KBern},
\begin{align*}
 \| P_k Q_j ( \phi_1 Q_{\ge j-C} \phi_2)\|_{\Ltwotx} &\les  \|\phi_1\|_{L^\infty L^2} \|Q_{\ge j-C}\phi_2\|_{L^2_t L^\infty_x}\\
&\les 2^{\frac{j-k_2}{4}} 2^{k_2} 2^{-\frac{j}{2}}
\|\phi_1\|_{S[k_1]} \|\phi_2\|_{S[k_2]}
\end{align*}
And third,
\begin{align*}
 \| P_k Q_j (Q_{\ge j-C} \phi_1 Q_{<j-C} \phi_2)\|_{\Ltwotx} &\les  \sum_{m= j+O(1)}
\| P_k Q_j (Q_{m} \phi_1 Q_{<j-C} \phi_2)\|_{\Ltwotx}\\
&\les \sum_{m= j+O(1)} \sum_{\kappa_1,\kappa_2} \| P_{k_1,\kappa_1}
Q_{m} \phi_1 P_{k_2,\kappa_2} Q_{<j-C} \phi_2\|_{\Ltwotx}\\
 & \les
\sum_{m= j+O(1)} \sum_{\kappa_1,\kappa_2}
\| P_{k_1,\kappa_1} Q_{m} \phi_1\|_{\Ltwotx}\| P_{k_2,\kappa_2} Q_{<j-C} \phi_2\|_{L^\infty_{t,x}}\\
&\les  \sum_{m= j+O(1)} 2^{-\frac{m}{2}} 2^{k_2} 2^{\frac{m-k_2}{4}}
\|\phi_1\|_{S[k_1]} \|\phi_2\|_{S[k_2]} \\&\les  2^{\frac{3k_2}{4}}
2^{-\frac{j}{4}} \|\phi_1\|_{S[k_1]} \|\phi_2\|_{S[k_2]}
\end{align*}
as claimed. The inner sums run over
$\kappa_1,\kappa_2\in\calC_{\frac{m-k_2}{2}}$ with
$\dist(\kappa_1,\kappa_2)\les 2^{\frac{m-k_2}{2} }$.
\end{proof}

Later we shall also need the following technical variants, both of which are in the same
spirit as Corollary~\ref{cor:squarebilin3}.

\begin{cor}
 \label{cor:bilinbasicsquare}
Let $\phi$ be adapted to $k_1$ and suppose for every $\kappa\in\calC_{m_0}$ with $m_0\le -100$
there is a Schwarz function~$\psi_\kappa$ which is adapted to~$k_2$. Then, provided $j\le k_2\le k=k_1+O(1)$,
\begin{equation}\label{eq:phi12hl'}
\sum_{\kappa\in\calC_{m_0}}  \| P_k Q_j (P_{k_1,\kappa} \phi \: \psi_\kappa)\|_{\Ltwotx} \les
 |m_0| 2^{\frac{3k_2}{4}} 2^{-\frac{j}{4}} \|\phi\|_{S[k_1]} \Big(\sum_{\kappa\in\calC_{m_0}} \|\psi_\kappa\|^2_{S[k_2]} \Big)^{\frac12}
\end{equation}
\end{cor}
\begin{proof}
 One uses the argument for the high-low case of Lemma~\ref{lem:bilinbasic}. In particular, $k=k_1+O(1)=0$. First,
with $m=\frac{j-k_2}{2}$,
\begin{align*}
\sum_{\kappa\in\calC_{m_0}} P_k Q_j (Q_{<j-C} P_{k_1,\kappa}\phi \: Q_{<j-C} \psi_\kappa) &= \sum_{\kappa\in\calC_{m_0}}\sum_{\kappa_1,\kappa_2\in\calC_{m}}
P_k Q_j (Q_{<j-C} P_{k_1,\kappa_1}P_{k_1,\kappa}\phi\: Q_{<j-C} P_{k_2,\kappa_2} \psi_\kappa)
\end{align*}
If $m\le m_0$, then by the
$L^2$-bilinear bound~\eqref{eq:bilin2}
\begin{align}
 & \sum_{\kappa\in\calC_{m_0}} \| P_k Q_j (Q_{<j-C} P_{k_1,\kappa} \phi\: Q_{<j-C} \psi_\kappa)\|_{\Ltwotx} \nn \\
&\les
\sum_{\kappa\in\calC_{m_0}}\sum_{\substack{\kappa_1,\kappa_2\in\calC_{m}\\ \kappa_1\subset\kappa}}  \|
P_k Q_j (Q_{<j-C} P_{k_1,\kappa_1} \phi\: Q_{<j-C} P_{k_2,\kappa_2} \psi_\kappa)\|_{\Ltwotx} \nn \\
&\les \sum_{\kappa\in\calC_{m_0}}\sum_{\substack{\kappa_1,\kappa_2\in\calC_{m}\\ \kappa_1\subset\kappa}}  2^{-\frac{j-k_2}{4}} 2^{\frac{k_2}{2}}
\| Q_{<j-C} P_{k_1,\kappa_1} \phi\|_{S[k_1,\kappa_1]}  \| Q_{<j-C} P_{k_2,\kappa_2} \psi_\kappa\|_{S[k_2,\kappa_2]} \nn  \\
&\les 2^{-\frac{j-k_2}{4}} 2^{\frac{k_2}{2}} \|\phi\|_{S[k_1]} \Big(\sum_{\kappa\in\calC_{m_0}}\|\psi_\kappa\|_{S[k_2]}^2\Big)^{\frac12}
\nn
\end{align}
where we applied Cauchy-Schwarz twice to pass to the last line.
If, on the other hand, $m>m_0$, then we first consider smaller modulations of~$\phi$. In fact, dropping the $Q_{<j-C}$ on~$\phi$
as we may one has
\begin{align*}
 & \sum_{\kappa\in\calC_{m_0}} \| P_k Q_j (Q_{<2m_0-C} P_{k_1,\kappa} \phi\: Q_{<j-C} \psi_\kappa)\|_{\Ltwotx} \\&\les
\sum_{\substack{\kappa_1,\kappa_2\in\calC_{m}\\ \kappa\subset\kappa_1}} \sum_{\kappa\in\calC_{m_0}} \|
P_k Q_j (Q_{<2m_0-C} P_{k_1,\kappa_1} \phi\: Q_{<j-C} P_{k_2,\kappa_2} \psi_\kappa)\|_{\Ltwotx} \\
&\les \sum_{\substack{\kappa_1,\kappa_2\in\calC_{m}\\ \kappa\subset\kappa_1}} \sum_{\kappa\in\calC_{m_0}}  2^{-\frac{j-k_2}{4}} 2^{\frac{k_2}{2}}
\| Q_{<2m_0-C} P_{k_1,\kappa} \phi\|_{S[k_1,\kappa]}  \| Q_{<j-C} P_{k_2,\kappa_2} \psi_\kappa\|_{S[k_2,\kappa_2]} \\
&\les 2^{-\frac{j-k_2}{4}} 2^{\frac{k_2}{2}} \|\phi\|_{S[k_1]} \Big(\sum_{\kappa\in\calC_{m_0}}\|\psi_\kappa\|_{S[k_2]}^2\Big)^{\frac12}
\end{align*}
where we again applied Cauchy-Schwarz twice to pass to the last line.
Finally, we need to account for~$Q_{2m_0-C\le \cdot<j-C}\phi$.  Fix $\ell$ with $2m_0-C\le \ell <j-C$ and
repeat the previous estimate. This yields
\[
\sum_{\kappa\in\calC_{m_0}} \| P_k Q_j (Q_{\ell} P_{k_1,\kappa} \phi\: Q_{<j-C} \psi_\kappa)\|_{\Ltwotx} \les
2^{-\frac{j-k_2}{4}} 2^{\frac{k_2}{2}} \|\phi\|_{S[k_1]} \Big(\sum_{\kappa\in\calC_{m_0}}\|\psi_\kappa\|_{S[k_2]}^2\Big)^{\frac12}
\]
which, upon summing in~$\ell$ yields the same bound with the loss of a factor of~$(j-2m_0)_+$. Replacing this by the larger~$|m_0|$ then
implies the bound of the corollary.
Next, by the improved
Bernstein estimate of Lemma~\ref{lem:KBern},  and Lemma~\ref{lem:enersquaresum},
\begin{align*}
\sum_{\kappa\in\calC_{m_0}}  \| P_k Q_j ( P_{k_1,\kappa}\phi\: Q_{\ge j-C} \psi_\kappa)\|_{\Ltwotx} &\les
\sum_{\kappa\in\calC_{m_0}}  \|P_{k_1,\kappa} \phi\|_{L^\infty L^2} \|Q_{\ge j-C}\psi_\kappa\|_{L^2_t L^\infty_x}\\
&\les |m_0|\, 2^{\frac{j-k_2}{4}} 2^{k_2} 2^{-\frac{j}{2}}
\|\phi\|_{S[k_1]} \Big(\sum_{\kappa\in\calC_{m_0}}\|\psi_\kappa\|_{S[k_2]}^2\Big)^{\frac12}
\end{align*}
And third,
\begin{align*}
&\sum_{\kappa\in\calC_{m_0}}  \| P_k Q_j (Q_{\ge j-C} P_{k_1,\kappa}\phi\: Q_{<j-C} \psi_\kappa )\|_{\Ltwotx} \\
&\les  \sum_{m= j+O(1)} \sum_{\kappa\in\calC_{m_0}}
\| P_k Q_j (Q_{m} P_{k_1,\kappa} \phi\: Q_{<j-C} \psi_\kappa)\|_{\Ltwotx}\\
&\les \sum_{m= j+O(1)} \sum_{\kappa\in\calC_{m_0}}  \sum_{\kappa_1\sim \kappa_2}
\| P_{k_1,\kappa_1} Q_{m} P_{k_1,\kappa} \phi\: P_{k_2,\kappa_2} Q_{<j-C} \psi_\kappa\|_{\Ltwotx}\\
&\les \sum_{m= j+O(1)} \sum_{\kappa\in\calC_{m_0}} \sum_{\kappa_1\sim \kappa_2}
\| P_{k_1,\kappa_1} Q_{m} P_{k_1,\kappa} \phi\|_{\Ltwotx} \;
2^{k_2} 2^{\frac{m-k_2}{4}} \| P_{k_2,\kappa_2} Q_{<j-C} \psi_\kappa\, \|_{\ener}\\
&\les 2^{k_2} 2^{\frac{j-k_2}{4}} \sum_{m= j+O(1)}  \Big(\sum_{\kappa_1,\kappa}
\| P_{k_1,\kappa_1} Q_{m} P_{k_1,\kappa} \phi\|_{\Ltwotx}^2\Big)^{\frac12} \;
 \Big(\sum_{\kappa_2,\kappa'} \| P_{k_2,\kappa_2} Q_{<j-C} \psi_{\kappa'}\, \|^2_{\ener}\Big)^{\frac12}\\
&\les    2^{\frac{3k_2}{4}} 2^{-\frac{j}{4}} \|\phi\|_{S[k_1]}
 \Big(\sum_{\kappa\in\calC_{m_0}}\|\psi_\kappa\|_{S[k_2]}^2\Big)^{\frac12}
\end{align*}
as claimed. The inner sums run over
$\kappa_1,\kappa_2\in\calC_{\frac{m-k_2}{2}}$ and $\kappa_1\sim\kappa_2$ denotes
$\dist(\kappa_1,\kappa_2)\les 2^{\frac{m-k_2}{2} }$.
\end{proof}

We shall also require the following estimates which gain something in terms of the small angle.

\begin{cor}
  \label{cor:capsbilin}
Given $\delta>0$ small and $L\gg1$, there exists $m_0(\delta,L)\ll-1$ with
the following property:  let $k,k_1,k_2\in\Z$ so that $\max_{i=1,2} |k-k_i|\le L$. For any
 $\phi_1$ and $\phi_2$ which are adapted to $k_1,k_2$, respectively, and $j\le k+C$,
\begin{equation}
\sum_{\substack{\kappa_1,\kappa_2\in\calC_{m_0}\\
\dist(\kappa_1,\kappa_2)\le 2^{m_0}}}  \| P_k Q_j
(P_{k_1,\kappa_1}\phi_1\: P_{k_2,\kappa_2}\phi_2)\|_{\Ltwotx} \le
\delta 2^{\frac{k-j}{3}} 2^{\frac{k_1}{2}} \trip\phi_1\trip_{S[k_1]}
\trip\phi_2\trip_{S[k_2]} \label{eq:L2bdhh_delta}
\end{equation}
In the high-low case $k=k_1+O(1)$, $k_2\le k_1-C$,
\begin{equation}
\sum_{\substack{\kappa_1,\kappa_2\in\calC_{m_0}\\
\dist(\kappa_1,\kappa_2)\le 2^{m_0}}}  \| P_k Q_j
(P_{k_1,\kappa_1}\phi_1\: P_{k_2,\kappa_2}\phi_2)\|_{\Ltwotx} \le
\delta 2^{\frac{k_2-j}{3}} 2^{\frac{k_2}{2}} \|\phi_1\|_{S[k_1]}
\|\phi_2\|_{S[k_2]} \label{eq:L2bdhl_delta}
\end{equation}
as well as
\begin{equation}
\Big(\sum_{\kappa_2\in\calC_{m_0}}  \| P_k Q_j (\phi_1\:
P_{k_2,\kappa}\phi_2)\|^2_{\Ltwotx}\Big)^{\frac12} \le \delta
2^{\frac{k_2-j}{3}} 2^{\frac{k_2}{2}} \|\phi_1\|_{S[k_1]}
\|\phi_2\|_{S[k_2]} \label{eq:L2bdhl_delta'}
\end{equation}
\end{cor}
\begin{proof} Let $k_1=0$ whence $|k|\le L$ and $|k_2|\le 2L$.
Implicit constants will be allowed to depend on~$L$.
By Corollary~\ref{cor:squarebilin3}
and~\eqref{eq:tripuse},
\begin{align*}
& \sum_{\substack{\kappa_1,\kappa_2\in\calC_{m_0}\\
\dist(\kappa_1,\kappa_2)\le 2^{m_0}}}   \|P_k Q_j (Q_{\le j+k-C}P_{k_1,\kappa_1} \phi_1\,
 Q_{\le j+k-C}P_{k_2,\kappa_2} \phi_2)\|_{\Ltwotx}\\
 & \les 2^{\frac{k+j}{4}} 2^{-\frac{j}{2}} \trip
 \phi_1\trip_{S[k_1]}\,
 \trip\phi_2\trip_{S[k_2]}
\le \delta 2^{-\frac{j}{4-}} \trip
 \phi_1\trip_{S[k_1]}\,
 \trip\phi_2\trip_{S[k_2]}
\end{align*}
which is sufficient.
Note that we used interpolation and $2^{\frac{k+j}{4}}\le 2^{\frac{m_0}{2}}$ which gives
the desired gain of~$\delta$ provided $m_0$ is small enough relative to~$\delta$ and~$L$.
For the remaining cases we use a variant of~\eqref{eq:highhigh2}:  with $2>r>1$, $\theta=\frac{2}{r}-1$,
and $\frac{1}{p}=\frac{1}{r}-\frac12$,
\begin{align}
 &\|P_k Q_j (Q_{>j+k-C}P_{k_1,\kappa_1}\phi_1\cdot P_{k_2,\kappa_2}\phi_2)\|_{\Ltwotx} \\ & \les
 2^{\frac{j}{4}\theta}   \|Q_{>j+k-C}\phi_1\cdot P_{k_2,\kappa_2}\phi_2\|_{L_t^2
 L^r_x} \nn
 \\
 &\les 2^{\frac{j}{4}\theta}  \|Q_{>j+k-C}P_{k_1,\kappa_1}\phi_1\|_{\Ltwotx} \| P_{k_2,\kappa_2}\phi_2\|_{L_t^\infty
 L^p_x} \nn \\
&\les  2^{\frac{j}{4}(\theta-2)} 2^{m_0(\frac12-\frac1p)} \|P_{k_1,\kappa_1}\phi_1\|_{\dot X_0^{0,\frac12,\infty}}
 \|P_{k_2,\kappa_2}\phi_2\|_{\ener}
\label{eq:highhigh2'}
\end{align}
Taking $\theta$ close to~$1$, one can make this $\le \delta 2^{-\frac{j}{4-}}$ as desired. This bound can be
 summed over $\kappa_1,\kappa_2$ by Cauchy-Schwarz and the definition of the~$S[k]$-norm; see also Lemma~\ref{lem:enersquaresum}.

In the high-low case $j\le k_2\le k=k_1+O(1)=0$ we proceed as follows.  First,
\begin{align*}
 & \| P_k Q_j (Q_{<j-C} P_{k_1,\kappa_1} \phi_1 Q_{<j-C} P_{k_2,\kappa_2} \phi_2)\|_{\Ltwotx} \\
 &\les  2^{-\frac{j-k_2}{2}} \min\big(2^{\frac{m_0}{2}}, 2^{\frac{j-k_2}{4}}\big) 2^{\frac{k_2}{2}}
 \|P_{k_1,\kappa_1}Q_{<j-C}  \phi_1\|_{S[k_1]} \|P_{k_2,\kappa_2}Q_{<j-C} \phi_2\|_{S[k_2]} \\
&\les \delta 2^{-\frac{j-k_2}{3}}2^{\frac{k_2}{2}} \|P_{k_1,\kappa_1} Q_{<j-C} \phi_1\|_{S[k_1]} \|P_{k_2,\kappa_2}Q_{<j-C} \phi_2\|_{S[k_2]}
\end{align*}
by a decomposition into caps of size~$2^{\frac{j-k_2}{2}}$ and the
$L^2$-bilinear bound~\eqref{eq:bilin2}.  The summation over $\kappa_1$ and~$\kappa_2$ can be carried out since it leads to
the square function~\eqref{eq:squarefunc}.
Next, by the improved Bernstein estimate, see Lemma~\ref{lem:KBern},
\begin{align*}
& \| P_k Q_j ( P_{k_1,\kappa_1}\phi_1 Q_{\ge j-C} P_{k_2,\kappa_2}\phi_2)\|_{\Ltwotx} \\
&\les  \|P_{k_1,\kappa_1}\phi_1\|_{L^\infty L^2} \|Q_{\ge j-C}P_{k_2,\kappa_2}\phi_2\|_{L^2_t L^\infty_x}\\
&\les \min(2^{\frac{j-k_2}{4}},2^{\frac{m_0}{2}}) 2^{k_2} 2^{-\frac{j}{2}}
\|P_{k_1,\kappa_1}\phi_1\|_{\ener} \|P_{k_2,\kappa_2}\phi_2\|_{\dot X_{k_2}^{0,\frac12,\infty}} \\
&\le   \delta 2^{-\frac{j-k_2}{3}}2^{\frac{k_2}{2}}  \|P_{k_1,\kappa_1}\phi_1\|_{\ener} \|P_{k_2,\kappa_2}\phi_2\|_{\dot X_{k_2}^{0,\frac12,\infty}}
\end{align*}
and summation over $\kappa_1,\kappa_2$ is again admissible.
Finally,
\begin{align}
& \| P_k Q_j (Q_{\ge j-C} P_{k_1,\kappa_1} \phi_1 Q_{<j-C}
P_{k_2,\kappa_2} \phi_2)\|_{\Ltwotx} \nn\\
&\les  \sum_{m= j+O(1)}
\| P_k Q_j (Q_{m} P_{k_1,\kappa_1} \phi_1 Q_{<j-C} P_{k_2,\kappa_2} \phi_2)\|_{\Ltwotx}\nn \\
&\les \sum_{m= j+O(1)} \sum_{\kappa_1',\kappa_2'}
\| P_{k_1,\kappa_1} Q_{m}P_{k_1,\kappa_1'} \phi_1 P_{k_2,\kappa_2} Q_{<j-C} P_{k_2,\kappa_2'}\phi_2\|_{\Ltwotx}\nn \\
&\les \sum_{m= j+O(1)} \sum_{\kappa_1',\kappa_2'}
\| Q_{m} P_{k_1,\kappa_1} P_{k_1,\kappa_1'}\phi_1\|_{\Ltwotx}\|  Q_{<j-C} P_{k_2,\kappa_2}P_{k_2,\kappa_2'} \phi_2\|_{L^\infty_{t,x}}\nn \\
&\les  \sum_{m= j+O(1)} 2^{-\frac{m}{2}} 2^{k_2} \min(2^{\frac{m-k_2}{4}},
2^{\frac{m_0}{2}}) \|P_{k_1,\kappa_1} \phi_1\|_{S[k_1]} \|P_{k_2,\kappa_2} \phi_2\|_{S[k_2]}\nn \\
&\les  \delta 2^{\frac{k_2}{2}} 2^{-\frac{j-k_2}{3}}
\|Q_{j+O(1)}P_{k_1,\kappa_1} \phi_1\|_{\dot
X_{k_1}^{0,\frac12,\infty}} \|P_{k_2,\kappa_2} \phi_2\|_{\ener}
\label{eq:sch123}
\end{align}
as claimed. The inner sums run over
$\kappa_1',\kappa_2'\in\calC_{\frac{j-k_2}{2}}$ with
$\dist(\kappa_1',\kappa_2')\les 2^{\frac{j-k_2}{2} }$. The bound
in~\eqref{eq:sch123} can be summed over the caps~$\kappa_1,\kappa_2$
by definition of the~$S[k]$ norm.

Finally, \eqref{eq:L2bdhl_delta'} follows from the preceding since
the gain of~$\delta$ was obtained only from the low-frequency
function~$\phi_2$. We can therefore square-sum the final estimate to
obtain the desired conclusion.
\end{proof}

\subsection{An algebra estimate for $S[k]$}

The following bilinear bound expresses something close to an algebra
property of the $S[k]$ spaces. It is obtained by removing the
restriction on the modulation of the output in
Lemma~\ref{lem:bilinbasic}.

\begin{lemma}
  \label{lem:Sk_prod2}  For any $j,k\in\Z$,
  \begin{equation}
    \label{eq:Sk_prod2}
    \| P_kQ_j (\phi\psi) \|_{\dot X^{0,\frac12,\infty}} \les 2^{k_1\wedge k_2} 2^{\frac{k-k_1\vee k_2}{4}}2^{\frac{j-k_1\wedge
    k_2}{4}\wedge0}\;
       2^{(k_1\vee  k_2-j)(\frac12-\eps)\wedge 0}\; \trip\phi\trip_{S[k_1]} \trip\psi\trip_{S[k_2]}
  \end{equation}
provided $\phi,\psi$ are Schwartz functions which are adapted to~$k_1$ and $k_2$, respectively.
\end{lemma}
\begin{proof}
 We commence with the high-high case $k_1=k_2+O(1)=0$ and $k\le O(1)$.
We need to prove that
\[
\| P_kQ_j (\phi\psi) \|_{\dot X^{0,\frac12,\infty}} \les
 2^{\frac{k}{4}} \min\big( 2^{\frac{j}{4}},
       2^{-j(\frac12-\eps)}\big)\; \trip\phi\trip_{S[k_1]} \trip\psi\trip_{S[k_2]}
\]
 To begin with,  one has
 \begin{align*}
   2^{j/2}  \| P_k Q_j  (Q_{>j-C} \phi\cdot \psi)\|_{\Ltwotx}
   &\les 2^{k} 2^{j/2} 2^{\frac{j-k}{4}\wedge0} \| P_k Q_j  (Q_{>j-C} \phi\cdot \psi)\|_{L_t^2L^1_x}\\
   &\les 2^{k} 2^{j/2} 2^{\frac{j-k}{4}\wedge0} \| Q_{>j-C}  \phi\|_{L^2 L^2} \| \psi\|_{L^\infty L^2}\\
    &\les 2^{k}  2^{\frac{j-k}{4}\wedge0}  \min(1,2^{-(\frac12-\eps)j}) \|\phi\|_{S[k_1]} \|\psi\|_{S[k_2]}\\
 \end{align*}
which is admissible. So it suffices to estimate $P_k Q_j (Q_{\le
j-C} \phi\cdot Q_{\le j-C}\psi)$. As usual, we perform a wave-packet
decomposition by means of Lemma~\ref{lem:cone}. Note
that~\eqref{eq:imbalance} holds here. We begin
with~\eqref{eq:plussum'} where we choose $r':= 2^{k}$. Thus, $k<-C$
and $j=O(1)$, and in view of \eqref{eq:bilin2}
\begin{align*}
    \| P_k Q_j (Q_{\le j-C} \phi^{+}\cdot Q_{\le j-C}\psi^{+}) \|_{\Ltwotx}
&\les  \sum_{\kappa\in \calC_{k}} \| P_\kappa Q_{\le j-C} \phi^{+}
\cdot P_{-\kappa} Q_{\le j-C}\psi^{+}\|_{\Ltwotx}\\
&\les   \sum_{\kappa\in \calC_{k}} |\kappa|^{\frac12} \| P_\kappa Q_{\le j-C} \phi^{+}\|_{S[k_1,\kappa]}
\| P_{-\kappa} Q_{\le j-C}\psi^{+}\|_{S[k_2,-\kappa]}\\
&\les 2^{k/2} \|\phi\|_{S[k_1]} \|\psi\|_{S[k_2]}
\end{align*}
where we invoked Lemma~\ref{lem:pm} in the final step. The same
estimate applies to~$\phi^{-}$ and~$\psi^{-}$. It therefore suffices
to assume that $j\le k+O(1)$; but then Lemma~\ref{lem:bilinbasic}
applies.

\noindent
Next, we consider the low-high case $k=k_2+O(1)=0$, $k_1<-C$. We need to prove that
\[
2^{\frac{j}{2}} \|P_0 Q_j (\phi\psi)\|_{\Ltwotx} \les 2^{k_1}
2^{\frac{j-k_1}{4}\wedge 0} \min(1, 2^{-j(\frac12-\eps)})
\|\phi\|_{S[k_1]} \|\psi\|_{S[0]}
\]
In view of Lemma~\ref{lem:bilinbasic} we can assume that $j\ge k_1$.
From the $\dot X^{s,b,q}$ components of the $S[k]$ norm,
\begin{align*}
& 2^{j/2} \| P_0 Q_j (Q_{\ge j-C}\phi\cdot \psi)\|_{\Ltwotx}
 \les 2^{j/2} \| Q_{\ge j-C}\phi\|_{L^2 L^\infty} \| \psi\|_{L^\infty L^2} \\
 &\les  2^{k_1}2^{\frac{j-k_1}{4}\wedge0} \min(1,2^{-(\frac12-\eps)j}) \|\phi\|_{S[k_1]} \|\psi\|_{S[k_2]}\\
&\les 2^{k_1} 2^{\frac{j-k_1}{4}\wedge0}
\min(1,2^{-(\frac12-\eps)j}) \|\phi\|_{S[k_1]} \|\psi\|_{S[k_2]}
\end{align*}
Finally, it remains to bound
\begin{equation}\nn
2^{j/2} \| P_0 Q_j (Q_{< j-C}\phi\cdot Q_{\ge j-C } \psi)\|_{\Ltwotx}
\end{equation}
which will be done using the usual angular decomposition.
In fact, from Lemma~\ref{lem:cone}, and provided $j\le C$, and with $\ell=\frac{m-k_1}{2}\wedge 0$,
\begin{align}
 2^{j/2} \| P_0 Q_j (Q_{< j-C}\phi\cdot Q_{\ge j-C } \psi)\|_{\Ltwotx}
&\les \sum_{m\ge j-C} 2^{j/2} \sum_{\kappa\in\calC_{\ell}} \| P_0 Q_j (P_{k_1,\kappa}
Q_{< j-C}\phi\cdot P_{k_2,\kappa'} Q_{m } \psi)\|_{\Ltwotx} \label{eq:m_sum2}\\
&\les \sum_{m\ge j-C} 2^{j/2} \sum_{\kappa\in\calC_{\ell}} \| P_{k_1,\kappa} Q_{< j-C}\phi\|_{L^\infty L^\infty}
\| P_{k_2,\kappa'} Q_{m } \psi \|_{\Ltwotx} \nn \\
&\les  \sum_{m\ge j-C} 2^{j/2} \sum_{\kappa\in\calC_{\ell}} 2^{k_1 }
2^{\frac{m-k_1}{4}\wedge0 }  \| P_{k_1,\kappa} \phi\|_{L^\infty L^2}
\| P_{k_2,\kappa'} Q_{m } \psi\|_{\Ltwotx}  \nn \\
&\les 2^{k_1} 2^{\frac{j-k_1}{4}\wedge 0}
 \|\phi\|_{S[k_1]} \|\psi\|_{S[k_2]}\nn
\end{align}
where we used Corollary~\ref{cor:S_cut} in the final inequality. If
$j\ge C$, then only $m=j+O(1)$ contributes to the sum
in~\eqref{eq:m_sum2}. The $\dot X^{0,1-\eps,2}$ component of the
$S[k]$-norm then leads to a gain of $\min(1,2^{-(\frac12-\eps)j})$
and we are done.
\end{proof}

\begin{cor} Under the same conditions as in the previous lemma and provided
 $k_1\ll k_2$ one has
\begin{equation}\label{eq:joch2}
 \| P_k (\phi Q_{<a}\psi) \|_{\dot X^{0,\frac12,1}} \les 2^{k_1} \big(1+ (k_2\wedge a -k_1)_+\big)
     \|\phi\|_{S[k_1]} \|\psi\|_{S[k_2]}
\end{equation}
where $k=O(1)+k_2$.
\end{cor}
\begin{proof}
 Summing \eqref{eq:Sk_prod2} over $j$ yields \eqref{eq:joch2} with $a\ge k_2$. It thus suffices to consider $a\le k_2$.
If $a\le k_1$  we use $Q_{<a} = Q_{<k_1}Q_{<a}$ to reduce matters to $a=k_1$ (see Corollary~\ref{cor:S_cut}). If $a=k_1$, then
\begin{align*}
&\sum_j 2^{j/2} \| P_k Q_j (\phi Q_{<a}\psi) \|_{\Ltwotx} \\&\le \sum_{j\le a+10} 2^{j/2} \| P_k Q_j (Q_{<a} \phi Q_{<a}\psi) \|_{\Ltwotx} \\
&\qquad + \sum_{j\ge a+10} 2^{j/2} \| P_k Q_j (Q_{j+O(1)} \phi Q_{<a}\psi) \|_{\Ltwotx} \\
&\les 2^{k_1}  \|\phi\|_{S[k_1]} \|\psi\|_{S[k_2]}  + \sum_{j\ge a+10} 2^{j/2} \|Q_{j+O(1)} \phi\|_{L^2 L^\infty} \| Q_{<a}\psi \|_{L^\infty L^2}\\
&\les 2^{k_1}  \|\phi\|_{S[k_1]} \|\psi\|_{S[k_2]}  + \sum_{j\ge a+10} 2^{j/2} 2^{(\frac32-\eps)k_1} 2^{-j(1-\eps)}
 \|\phi\|_{S[k_1]} \| Q_{<a}\psi \|_{S[k_2]} \\
&\les 2^{k_1}  \|\phi\|_{S[k_1]} \|\psi\|_{S[k_2]}
\end{align*}
as desired. The sum over $j\le a+10$ was estimated via
Lemma~\ref{lem:Sk_prod2}. If $k_1<a\le k_2$, one proceeds similarly.
\end{proof}

\subsection{Bilinear estimates involving both $S[k_1]$ and~$N[k_2]$ waves.}

The following lemma is a crucial tool. In essence, it expresses the property $N\times S \hookrightarrow N$.

\begin{lemma}
  \label{lem:core} For $\phi$ and $F$ which are $k_1$ and $k_2$-adapted, respectively, one has
  \begin{equation}\label{eq:core}
\|P_k (\phi F)\|_{N[k]} \les 2^{k_1\wedge k_2}\, 2^{\frac{j-k\wedge k_1\wedge k_2}{4}\wedge 0}  \|\phi\|_{S[k_1]} \|F\|_{N[k_2]}
  \end{equation}
  provided $P_{k_2} Q_j F=F$ and under the following condition
\begin{equation}\label{eq:weird_cond}
 P_k Q_{\le j-C} (Q_{<j-C} \phi\cdot F) = P_k Q_{\le j-C} (Q_{<j+k-k_1-C} \phi\cdot F)
\end{equation}
in the case $k_1=k_2+O(1)\ge k+O(1)\ge j$.  If \eqref{eq:weird_cond}
fails, then one loses a factor of $1+(k_1-k)_+$ on the right-hand
side of~\eqref{eq:core}; alternatively, one has the following weaker
version of~\eqref{eq:core}
  \begin{equation}\label{eq:weaker_core}
\|P_k (\phi F)\|_{N[k]} \les 2^{k_1\wedge k_2}\, 2^{\frac{j-k\wedge k_1\wedge k_2}{4}\wedge 0}  \|\phi\|_{\dot X_{k_1}^{0,\frac12,1}} \|F\|_{N[k_2]}
  \end{equation}
\end{lemma}
\begin{proof} We remark beforehand that this proof will only use the $\dot X_{k_1}^{0,\frac12,\infty}$-norm for the elliptic regime
$\phi= P_k Q_{\ge k}\phi$ of the $S[k]$ norm. In particular, the imbedding
$\|\phi\|_{S[k_1]}\les \|\phi\|_{\dot X_{k_1}^{0,\frac12,\infty}}$ holds without any restrictions on the modulation, cf.~\eqref{eq:Skimbed}.
   We start with the high-high case $k_1=k_2+O(1)=0$.  Throughout this proof, we shall freely use
 Lemma~\ref{lem:QLp} in order to remove~$Q_{<j-C}$ from various estimates.
First, we consider the case where  $\phi=Q_{\ge j-C} \phi$. If $j\ge k$,  then by Bernstein's and H\"older's inequalities
   \begin{align*}
    \|P_k(\phi F)\|_{N[k]} &\les  2^{-k}\|P_k(\phi F)\|_{L^1L^2}  \les  \|\phi F\|_{L^1 L^1} \\
     &\les  \|\phi\|_{\Ltwotx} \|F\|_{\Ltwotx} \les \|\phi\|_{S[k_1]} 2^{-\frac{j}{2}} \|F\|_{\Ltwotx}\\
     &\les \|\phi\|_{S[k_1]} \|F\|_{N[k_2]}
   \end{align*}
which is admissible. If on the other hand $j\le k$, then we again have to consider several subcases.
If $\phi=Q_{\ge k}\phi$, then
\begin{align*}
 2^{-k} \| Q_{\ge k} \phi\cdot F\|_{L^1_t L^2_x} &\les \| Q_{\ge k}\phi\cdot F\|_{L^1_t L^1_x}
\les \|Q_{\ge k}\phi\|_{\Ltwotx} \| F\|_{\Ltwotx}  \\
&\les 2^{\frac{j-k}{2}} \|\phi\|_{S[k_1]} \|F\|_{N[k_2]}
\end{align*}
which is admissible. Hence it suffice to assume that $\phi=Q_{ j-C\le\cdot\le k} \phi$.
Furthermore, we can assume that the output is at modulation $\le j$. In fact,
by the improved Bernstein's inequality,
\begin{align*}
     \| P_kQ_{> j} (\phi F)\|_{N[k]} & \les 2^{-k}\sum_{\ell>j} 2^{-\frac{\ell}{2}} \|P_k Q_\ell(\phi F)\|_{\Ltwotx} \\
     &\les \sum_{\ell>j} 2^{-\frac{\ell}{2}} 2^{\frac{\ell-k}{4}\wedge 0} \|\phi F\|_{L^2 L^1} \\
  &\les  \sum_{\ell>j} 2^{\frac{j-\ell}{2}} 2^{\frac{\ell-k}{4}\wedge 0}     \|\phi\|_{L_t^\infty L^2_x} 2^{-\frac{j}{2}}\|F\|_{\Ltwotx} \\
     &\les  2^{\frac{j-k}{4}}  \|\phi\|_{S[k_1]} \|F\|_{N[k_2]}
   \end{align*}
as desired. Now consider the output of modulation at most $j$. We also first restrict ourselves to the
contributions by~$Q_{j-C<\cdot<j+C} \phi$.
Thus,  by Lemma~\ref{lem:cone} and Lemma~\ref{lem:QLp},
\begin{align*}
 &\|P_k Q_{\le j} ( Q_{j-C<\cdot<j+C} \phi \cdot F)\|_{N[k]} \\&\les 2^{-k}\sum_{\ell= j+O(1)}  \| P_k Q_{\le j} ( Q_{\ell} \phi \cdot F)\|_{L^1_t L^2_x} \\
&\les   2^{-k}\sum_{\ell= j+O(1)}  \sum_{D\in\calD_{k}}\sum_{\kappa\sim\kappa'\in \calC_{(\ell+k)/2}}
\| P_k Q_{\le j} ( P_D P_{\kappa} Q_{\ell}\, \phi \cdot P_{\kappa'} P_{-D} F)\|_{L^1_t L^2_x} \\
&\les  2^{-k}\sum_{\ell= j+O(1)}  \sum_{D\in\calD_{k}} \sum_{\kappa\sim\kappa'\in \calC_{(\ell+k)/2}}
\| P_D P_{\kappa} Q_{\ell}\, \phi \cdot P_{\kappa'} P_{-D} F\|_{L^1_t L^2_x}
\end{align*}
where $\kappa\sim\kappa'$ means $\dist(\kappa,\kappa')\les \diam(\kappa)$. Moreover, $\calD_k$ is a cover of $\{|\xi|\sim1\}$ by disks
of diameter $2^k$ and
with  overlap uniformly bounded in~$k$; the associated projections are $P_\kappa$.
Hence, one can further estimate (recall $j\le k$)
\begin{align*}
 &\|P_k Q_{\le j} ( \phi \cdot F)\|_{N[k]} \\
&\les  2^{-k} \sum_{\ell= j+O(1)}  \sum_{D\in\calD_{k}} \sum_{\kappa\sim\kappa'\in \calC_{(\ell+k)/2}}
\| P_D P_{\kappa} Q_{\ell}\, \phi\|_{L^2_t L_x^\infty} \|P_{\kappa'} P_{-D} F\|_{L^2_t L^2_x} \\
&\les  2^{-k}\sum_{\ell= j+O(1)}  \sum_{D\in\calD_{k}} \sum_{\kappa\sim\kappa'\in \calC_{(\ell+k)/2}}   2^{\frac{\ell+3k}{4}}
\| P_D P_{\kappa} Q_{\ell}\, \phi\|_{L^2_t L_x^2} \|P_{\kappa'} P_{-D} F\|_{L^2_t L^2_x}\\
&\les 2^{(j-k)/4} \|\phi\|_{S[k_1]} \|F\|_{N[k_2]}
\end{align*}
Next, we consider the output of modulation at most $j$ and~$\phi=Q_{j+C\le\cdot\le k} \phi$.
Then we are in the ``imbalanced case'' of Lemma~\ref{lem:cone} whence
\[
 P_k Q_{\le j} (\phi F) = \sum_{k\ge \ell\ge j+C} \sum_{\kappa,\kappa',\kappa''}   P_{k,\kappa} Q_{\le j} (P_{k_1,\kappa'} \phi \cdot P_{k_2,\kappa''} F)
\]
where $\kappa\in\calC_{\frac{\ell-k}{2}}$, $\kappa',\kappa''\in \calC_{(\ell+k)/2}$ and $\dist(\kappa,\kappa')\sim 2^{\frac{\ell-k}{2}}$,
$\dist(\kappa',\kappa'')\sim 2^{\frac{\ell+k}{2}}$. Using~\eqref{eq:bilin1} one obtains
\begin{align*}
 &\|P_k Q_{\le j} (\phi F)\|_{N[k]}  \\
&\le  2^{-k} \sum_{k\ge\ell\ge j+C} \sum_{\kappa,\kappa',\kappa''}
\|  P_{k,\kappa} Q_{\le j} (P_{k_1,\kappa'}Q_\ell \phi \cdot P_{k_2,\kappa''} F)\|_{\NF[\kappa]} \\
&\les  \sum_{k\ge\ell\ge j+C} 2^{-k} \frac{2^{\frac{\ell+k}{4}}}{2^{\frac{\ell-k}{2}\wedge 0}}
\sum_{\kappa',\kappa''} \| P_{k_1,\kappa'}  Q_\ell\phi\|_{S[k_1,\kappa]}  \|P_{k_2,\kappa''} F\|_{\Ltwotx} \\
&\les  \sum_{k\ge\ell\ge j+C} 2^{\frac{j-k}{4}} 2^{\frac{j-\ell}{4}}\|Q_\ell \phi\|_{S[k_1]} 2^{-\frac{j}{2}}\|F\|_{\Ltwotx}
\les 2^{\frac{j-k}{4}} \|\phi\|_{S[k_1]} \|F\|_{N[k_2]}
\end{align*}
as desired.  To pass to the final inequality one uses~\eqref{eq:Skimbed} as well as
$\|Q_\ell \phi\|_{S[k_1]}\sim \|Q_\ell \phi\|_{\dot X_{k_1}^{0,\frac12,\infty}}$.

\noindent
Now assume $Q_{<j-C} \phi=\phi$. We first dispose of outputs of modulation exceeding $j-C$. If $j\ge k$, then
\begin{align*}
  &\|P_k Q_{>j-C}(\phi F)\|_{N[k]} \les  2^{-\frac{j}{2}} \|P_k Q_{>j-C} (\phi F)\|_{\Ltwotx} \les
  2^{-\frac{j}{2}} \|\phi\|_{L^\infty L^2} \|F\|_{L^2L^\infty}\\
  &\les  2^{-\frac{j}{2}} \|\phi\|_{S[k_1]}  \|F\|_{\Ltwotx} \les  \|\phi\|_{S[k_1]}  \|F\|_{N[k_2]}
\end{align*}
which is admissible. On the other hand, if $j\le k$, then
\begin{align*}
 & \|P_k Q_{>j-C}(\phi F)\|_{N[k]}
\les  2^{-k} \sum_{\ell>j-C} 2^{-\frac{\ell}{2}} \|P_k Q_{\ell} (\phi F)\|_{\Ltwotx} \\
& \les
   \sum_{k\ge \ell>j-C} 2^{-\frac{\ell}{2}} 2^{\frac{\ell-k}{4}} \|P_k Q_{\ell} (\phi F)\|_{L^2_t L^1_x}  + 2^{-\frac{k}{2}}
   \|P_kQ_{\ge k} (\phi F)\|_{\Ltwotx}
\les  2^{\frac{j-k}{4}} \| \phi\|_{S[k_1]} \|F\|_{N[k_2]}
\end{align*}
as desired.
It therefore remains to consider
\[
P_k Q_{\le j-C} (Q_{<j-C} \phi\cdot F) = \sum_{\pm} P_k Q_{\le j-C} (Q_{<j-C} \phi^{\pm}\cdot F^{\pm})
\]
where all four possibilities $(++),(+-),(-+),(--)$ are allowed on the right-hand side.
We first dispose of the contributions ``opposing waves'' as described by~\eqref{eq:plussum'}. This occurs only if
$k<-C$ and~$j=O(1)$, in fact,
\[
P_k Q_{\le j-C} (Q_{<j-C} \phi^{+}\cdot F^{+}) = P_k Q_{-C<\cdot<C} (Q_{<j-C} \phi^{+}\cdot F^{+})
\]
whence
\begin{align*}
  &\|P_k Q_{\le j-C} (Q_{<j-C} \phi^{+}\cdot F^{+}) \|_{N[k]} \\&\les  \|\phi^{+} F^{+}\|_{L^2 L^1}
  \les  \|\phi^{+}\|_{L^\infty L^2} \| F^{+}\|_{\Ltwotx} \\&\les  \|\phi\|_{S[k_1]} \| F\|_{N[k_2]}
\end{align*}
which is admissible. Therefore, we can now ignore the contribution of~\eqref{eq:plussum'}.
Let us now also assume without loss of generality that $\phi=\phi^{+}$, see Lemma~\ref{lem:pm}.
Using duality and Lemma~\ref{lem:cone},  one obtains  in view of~\eqref{eq:weird_cond}
\begin{equation}\label{eq:wavepack}
P_k Q_{\le j-C} (Q_{<j-C} \phi\cdot F) =
\sum_{\pm}\sum_{\kappa,\kappa',\kappa'' }P_{k,\pm\kappa} Q_{\le j-C}^{\pm} (P_{k_1,\kappa'}Q_{<j+k-C} \phi\cdot P_{k_2,\kappa''} F)
\end{equation}
with caps $\kappa\in\calC_{\ell}$, $\kappa',\kappa''\in\calC_{\ell}$  satisfying $\dist(\kappa',\kappa'')\sim 2^{m}$, $\dist(\kappa,\kappa')\sim 2^\ell$ where
 $\ell=(j-k)/2$, $m=(j+k)/2$. Note that  Lemma~\ref{lem:cone} also implies that $j\le k+O(1)$. Since
\begin{equation}\label{eq:wavepack_disp}
P_{k,\pm\kappa} Q_{\le j-C}^{\pm} = P_{k,\pm\kappa} Q_{\le k+2\ell-C}^{\pm}
\end{equation}
the right-hand side of \eqref{eq:wavepack} represents a wave-packet decomposition in the sense of Definition~\ref{def:Nk}.
Moreover, the operators in~\eqref{eq:wavepack_disp} are disposable in the sense of Lemma~\ref{lem:dispose}. Therefore,
\begin{align*}
 \|P_k Q_{\le j-C} (Q_{<j-C} \phi\cdot F)\|_{N[k]}
&\les 2^{-k}\max_{\pm} \sum_{\kappa',\kappa''} \big\|P_{k_1,\kappa'}Q_{<j+k-C} \phi\cdot F\|_{\NF[\kappa]}
\end{align*}
We could discard $\kappa$ here since the choice of $\kappa'$ leaves only a finite number of choice of~$\kappa$.
Invoking~\eqref{eq:bilin1}, this can be further estimated by
\begin{align*}
& \les 2^{-k}
\max_{\pm}\sum_{\kappa',\kappa''}\frac{2^{\frac{j+k}{4}}}{2^{\frac{j-k}{4}}}  \|P_{k_1,\kappa'} Q_{<j+k-C} \phi \|_{S[k_1,\kappa']}
\| P_{k_2,\kappa''} F\|_{\Ltwotx} \\
 &\les 2^{(j-k)/4} \|\phi\|_{S[k_1]} \|F\|_{N[k_2]}
\end{align*}
To pass to the final inequality here it was essential that \eqref{eq:weird_cond} reduced the modulation of~$\phi$ from $<j-C$
to~$<j+k-C$. Indeed, if~\eqref{eq:weird_cond} fails, then we need to write
$Q_{<j-C}\phi=Q_{<j+k-C}\phi+ Q_{j+k-C\le \cdot <j-C}\phi$. For the first summand here one applies the argument we just gave,
whereas for the second summand the best one can do is to invoke~\eqref{eq:Skimbed}
which results in the loss of of factor of~$k$ a claimed.
This concludes the high-high case.

Let us now consider the low-high case $k_1<-C$, $k_2=k=O(1)$. Since \eqref{eq:Sbd2} implies that
\begin{align*}
\|Q_{\ge j-C}\phi\cdot F\|_{L^1L^2}  &\les \|Q_{\ge j-C}\phi\|_{L^2L^\infty} \|F\|_{\Ltwotx} \\
&\les 2^{-\frac{j}{2}} 2^{\frac{j-k_1}{4}\wedge0} 2^{k_1} \|\phi\|_{S[k_1]} \|F\|_{\Ltwotx} \\
&\les 2^{\frac{j-k_1}{4}\wedge0}  2^{k_1} \|\phi\|_{S[k_1]} \|F\|_{N[k_2]}
\end{align*}
it will suffice to bound $\|Q_{\le j-C}\phi\cdot F\|_{N[k]}$. Moreover, if $j\ge k_1+C$, then the modulation of
$Q_{\le j-C}\phi\cdot F$ is on the order of~$j$ whence
\begin{align*}
  \| Q_{\le j-C}\phi\cdot F\|_{N[k]} &\les 2^{-\frac{j}{2}}\| Q_{\le j-C}\phi\cdot F\|_{\Ltwotx} \\
  &\les 2^{-\frac{j}{2}}\| Q_{\le j-C}\phi\|_{L^\infty L^\infty} \| F\|_{\Ltwotx} \\
  &\les 2^{-\frac{j}{2}}2^{k_1}  \|\phi\|_{S[k_1]} \| F\|_{\Ltwotx} \les 2^{k_1} \|\phi\|_{S[k_1]} \| F\|_{N[k_2]}
\end{align*}
as desired. We may therefore assume that $j\le k_1+C$. We first consider the case where the output
has modulation~$\ge j-C$.
More precisely, let
 $j-C\le m\le k_1+C$, $\ell=(m-k_1)/2$, as well as without loss of generality $F=F^{+}$. Then
 by the balanced modulation case of Lemma~\ref{lem:cone},
\[
Q_{m}(Q_{\le j-C}\phi\cdot F)  = \sum_{\kappa,\kappa'} P_{k,\kappa} Q_{m}^{+} (P_{k_1,\kappa'} Q_{\le j-C}\phi\cdot F)
\]
where $\kappa,\kappa'$ are caps of size $C^{-1}2^{\ell}$ and with
$\dist(\kappa,\kappa')\les 2^\ell$. Therefore,
\begin{align*}
 \|Q_{m}(Q_{\le j-C}\phi\cdot F) \|_{N[k]}
  &\les 2^{-m/2} \| \sum_{\kappa,\kappa'} P_{k,\kappa} Q_{m}^{+} (P_{k_1,\kappa'} Q_{\le j-C}\phi\cdot F) \|_{\Ltwotx} \\
&\les 2^{-m/2}
 \Big( \sum_{\kappa,\kappa'} \|P_{k_1,\kappa'} Q_{\le j-C}\phi\|^2_{L^\infty L^\infty}
\| F\|^2_{L^2 L^2}\Big)^{\frac12} \\
&\les 2^{-m/2}
 \Big( \sum_{\kappa,\kappa'} 2^{2k_1} 2^{\ell} \|P_{k_1,\kappa} Q_{\le j-C}\phi\|^2_{L^\infty L^2}\;
\| F\|^2_{L^2 L^2}\Big)^{\frac12} \\
&\les 2^{-m/2} 2^{-(k_1-m)/4} 2^{k_1} \|Q_{\le j-C}\phi\|_{L^\infty
L^2} \|F\|_{\Ltwotx}\\& \les 2^{-m/2} 2^{-(k_1-m)/4} 2^{k_1}
\|\phi\|_{S[k_1]} \|F\|_{\Ltwotx}
\end{align*}
where we used Lemma~\ref{lem:QLp}  in the final step.
Summing over $m\ge j-C$ implies that
\begin{align*}
\|Q_{>j-C}(Q_{\le j-C}\phi\cdot F) \|_{N[k]} &\les 2^{k_1} 2^{-j/2} 2^{-(k_1-j)/4} \|\phi\|_{S[k_1]} \|F\|_{\Ltwotx} \\
&\les  2^{k_1} 2^{-(k_1-j)/4} \|\phi\|_{S[k_1]} \|F\|_{N[k_2]}
\end{align*}
which is admissible.

It therefore remains to bound $\|Q_{<j-C}(Q_{\le j-C}\phi\cdot F)
\|_{N[k]}$ for which we shall again apply a wave-packet
decomposition as in Lemma~\ref{lem:cone}. Since $j\le k_1+C$ and
$k_1\ll -1$, we can assume that $j\le -C$ in applying
Lemma~\ref{lem:cone} (which allows us to ignore the opposing $(++)$
or $(--)$ contributions in~\eqref{eq:plussum'}). Without loss of
generality, we assume further that $\phi=\phi^{+}$ (see
Lemma~\ref{lem:pm}).  Then with caps $\kappa,\kappa'$ of size
$C^{-1} 2^m$ and separation $\sim 2^m$ where $m:=(j-k_1)/2$,
\begin{equation}\begin{aligned}
  \|Q_{<j-C}(Q_{\le j-C}\phi\cdot F) \|_{N[k]}
&\les  \Big\| \sum_{\kappa\sim \kappa'} P_{k,\kappa} Q_{<j-C}
(P_{k_1,\kappa'}Q_{\le j-C}
 \phi\cdot  F) \Big\|_{N[k]}  \\
 &\les  \Big(\sum_{\kappa\sim\kappa'} \Big\| P_{k,\kappa} Q_{<j-k_1}
(P_{k_1,\kappa'}Q_{\le j-C}
 \phi\cdot  F) \Big\|_{N[k]}^2\Big)^{\frac12}
\end{aligned}\label{eq:low-highNF}
\end{equation}
where we used Corollary~\ref{cor:N_cut} to dispose of~$Q_{\le j-C}$.
In view of~\eqref{eq:bilin1} this is further bounded by
\begin{align*}
 &\les 2^{\frac{k_1}{2}} 2^{-\frac{j-k_1}{4}}    \Big( \sum_{\kappa'} \|P_{k_1,\kappa'}Q_{\le j-C}
 \phi\|^2_{S[k_1,\kappa']} \|F\|^2_{\Ltwotx}\Big)^{\frac12}\\
 &\les 2^{(j-k_1)/4}  2^{k_1} \|\phi\|_{S[k_1]} \|F\|_{N[k_2]}
\end{align*}
as desired.

It remains to consider the high-low case $k=k_1+O(1)=0$, $k_2<-C$. First,
\begin{align*}
  \| Q_{>j+k_2} \phi \cdot F\|_{N[k]} &\les \|  Q_{>j+k_2} \phi \cdot F\|_{L^1 L^2} \\
  &\les \| Q_{>j+k_2} \phi \|_{L^2 L^2} \| F\|_{L^2 L^\infty} \\
  &\les 2^{-(j+k_2)/2} \|\phi\|_{S[k_1]} 2^{k_2} 2^{\frac{j-k_2}{4}\wedge0}  \|F\|_{\Ltwotx} \\
  &\les 2^{k_2} 2^{\frac{j-k_2}{4}\wedge0}   \|\phi\|_{S[k_1]} 2^{-\frac{j}{2}} 2^{-k_2} \|F\|_{\Ltwotx}
\end{align*}
which is acceptable with a factor of $2^{\frac{k_2}{2}}$ to spare. The reason for using $Q_{>j+k_2}$ rather
than $Q_{>j}$ will become clear momentarily.  Next,
\begin{align}
 \| Q_{\le j+k_2} \phi \cdot F\|_{N[k]} &\les  \| Q_{\ge j+k_2-C} [ Q_{\le j+k_2} \phi \cdot F]\|_{N[k]} \label{eq:QQ1}\\
 &\quad +  \| Q_{< j+k_2-C} [ Q_{\le j+k_2} \phi \cdot F]\|_{N[k]} \label{eq:QQ2}
\end{align}
As usual,  \eqref{eq:QQ1}  is controlled in the $\dot X^{-1,-\frac12,1}$  norm whence
\begin{align*}
  \eqref{eq:QQ1} &\les 2^{-(j+k_2)/2} \|  Q_{\le j+k_2} \phi \cdot F \|_{\Ltwotx} \\
  &\les 2^{-(j+k_2)/2} \|  Q_{\le j+k_2} \phi \|_{L^\infty L^2}  \| F \|_{L^2 L^\infty} \\
  &\les 2^{-(j+k_2)/2} \|  Q_{\le j+k_2} \phi \|_{L^\infty L^2} 2^{k_2} 2^{\frac{j-k_2}{4}\wedge0}  \| F \|_{\Ltwotx}\\
  &\les 2^{k_2} 2^{\frac{j-k_2}{4}\wedge0} \|\phi\|_{S[k_1]} 2^{-\frac{j}{2}} 2^{-k_2} \|F\|_{\Ltwotx}
\end{align*}
which is again acceptable. Finally, we perform a wave-packet decomposition on~\eqref{eq:QQ2} via Lemma~\ref{lem:cone}
in the imbalanced case and duality.
Thus,   one has
\[
Q_{< j+k_2-C} ( Q_{\le j+k_2} \phi \cdot F )  =  \sum_{\kappa,\kappa'} P_{k,\kappa} Q_{< j+k_2-C}^{+} (P_{k_1,\kappa'}
 Q_{\le j+k_2} \phi \cdot F )
\]
where the sum runs over pairs of caps $\kappa,\kappa'$ of size~$C^{-1}2^\ell$
with $\ell:=(j+k_2)/2$ and $\dist(\kappa,\kappa')\sim 2^\ell$. Moreover, $j\le k_2+O(1)$ since the only other
possibility  $j=O(1)$ allowed by~\eqref{eq:plussum'}  contributes
a vanishing term (as does $Q_{< j+k_2-C}^{-}$).
Therefore, with $\kappa'\sim\kappa$ denoting the admissible pairs,
\begin{align*}
  \label{eq:QQ2} &\les  \Big(\sum_{\kappa} \big\|
\sum_{\kappa'\sim \kappa} P_{k,\kappa} Q_{< j+k_2-C}^{+} (P_{k_1,\kappa'}
 Q_{\le j+k_2} \phi \cdot F ) \big\|_{N[\kappa]}^2\Big)^{\frac12} \\
 &\les  2^{-\frac{\ell}{2}} 2^{j/2} 2^{k_2} \|\phi\|_{S[k_1]} \|F\|_{N[k_2]} \\
 &\les 2^{k_2} 2^{(j-k_2)/4} \|\phi\|_{S[k_1]}\|F\|_{N[k_2]}
\end{align*}
as desired.
\end{proof}

There is the following general estimate that does not require~\eqref{eq:weird_cond} since
we restrict ourselves to $k\ge k_1+O(1)$.

\begin{cor}
  \label{cor:core_max} For $\phi$ and $F$ which are $k_1$ and $k_2$-adapted, respectively, one has
  \begin{equation}\label{eq:core_max}
\|P_{k} (\phi F)\|_{N[k]} \les 2^{k_1\wedge k_2}\,
2^{\frac{j-k\wedge k_1\wedge k_2}{4}\wedge 0}  \|\phi\|_{S[k_1]}
\|F\|_{N[k_2]}
  \end{equation}
  provided $P_{k_2} Q_j F=F$ and $k=k_1\vee k_2+O(1)$.
\end{cor}
\begin{proof}
 This is an immediate consequence of Lemma~\ref{lem:core}.
\end{proof}

Another important technical variant of Lemma~\ref{lem:core} has to
do with an additional angular localization of the inputs. This will
be important later in the trilinear section. Its statement is
somewhat technically cumbersome, but this is precisely the form in
which we shall use it later.

\begin{cor}
  \label{cor:angcore} Let $\phi$ be $k_1$-adapted, and
assume that   for some $m_0\le -100$, for every $
\kappa\in\calC_{m_0}$ there is a Schwarz function~$F_{\kappa}$ which is
adapted to~$k_2$ and so that $P_{k_2} Q_j F_\kappa=F_\kappa$. Then
  \begin{align}\label{eq:coresquare}
\sum_{\kappa\in\calC_{m_0}} \|P_k (P_{k_1,\kappa} \phi\, F_{\kappa}
)\|_{N[k]} &\les |m_0| \, 2^{k_1}\, 2^{\frac{j-k_1}{4}\wedge 0}
\|\phi\|_{S[k_1]} \Big(\sum_{\kappa\in\calC_{m_0}}  \|
 F_\kappa \|_{N[k_2]}^2\Big)^{\frac12}
  \end{align}
   provided we are in the low-high
  case $k=k_2+O(1)\ge k_1$. The sum here runs over caps with
  $\dist(\kappa_1,\kappa_2)\les 2^{m_0}$.
\end{cor}
\begin{proof}
  For this, one simply repeats the proof of the low-high case of
  Lemma~\ref{lem:core} with one additional twist: since
  $\sum_{\kappa} \|P_{k_1,\kappa}\phi\|_{S[k]}^2$ cannot be
  controlled by~$\|\phi\|_{S[k_1]}$, one has to check carefully that
  the square summation --- which~\eqref{eq:coresquare} leads to after Chauchy-Schwarz ---  is compatible with
  the estimates we are making (the norm for~$F$ is always~$\Ltwotx$). This is the case if we place
  $P_{k_1,\kappa}\phi$ in~$\ener$ or an~$\dot X^{s,b}$-norm. In the latter case one does not
incur any loss due to orthogonality, whereas in the former case there is a loss of~$|m_0|$, see Lemma~\ref{lem:enersquaresum}.
  The only place where one cannot use either of these norms is~\eqref{eq:low-highNF}.
  Indeed, if~$k_1+2m_0\le j -C$, then the caps of sizes~$2^{m_0}$
  are  smaller than those of size~$2^\ell=2^{\frac{j-k_1}{2}}$ in the
  wave-packet decomposition of~\eqref{eq:low-highNF}. In this case,
  however, one considers a wave-packet decomposition induced by the
  projections~$P_{k_1,\kappa}Q_{<k_1+2m_0}$
  with~$\kappa\in\calC_{m_0}$ which leads to the desired bound; the
  remaining projection~$P_{k_1,\kappa}Q_{k_1+2m_0\le \cdot\le j-C}$
  is then controlled by means of Lemma~\ref{lem:square_func} leading
  to a loos of~$|m_0|$ as claimed. If, on the other hand, $k_1+2m_0> j -C$,
  then this issue does not arise at all and the
  estimate~\eqref{eq:low-highNF} is performed essentially as in
  Lemma~\ref{lem:core} --- the only difference being that the caps
  in the wave-packet decomposition are grouped together inside the larger $\calC_{m_0}$--caps.
\end{proof}

\subsection{Nullform bounds in the high-high case}

Henceforth,  $\|\cdot\|_{S[k]}$  will mean the stronger
norm~$\trip\cdot\trip_{S[k]}$.
 The following
definition introduces the basic nullforms as well as the method of
``pulling out a derivative''.

\begin{defi}\label{def:Nab}
 The nullforms $\calN_{\alpha\beta}$ for $0\le\alpha,\beta\le 2$, $\alpha\ne\beta$, are defined as
\[
 \calN_{\alpha\beta}(\phi,\psi):= R_\alpha \phi R_\beta \psi - R_\beta \phi R_\alpha \psi
\]
whereas
\[
 \calN_0(\phi,\psi):=  R_\alpha \phi R^\alpha \psi
\]
By ``pulling out a derivative from'' from $\calN_{\alpha\beta}$ we
mean writing
\[
 \calN_{\alpha\beta}  (\phi,\psi) = \del_{\alpha} (|\nabla|^{-1} \phi R_\beta\psi) -  \del_{\beta} (|\nabla|^{-1} \phi R_\alpha\psi)
\]
or the analogous expression with $\phi$ and $\psi$ interchanged.
\end{defi}

Recall the $L^2$-bound~\eqref{eq:bilin3null} of
Lemma~\ref{lem:bilin3} for $\calN_{\alpha\beta}$-nullforms. We
separate the nullform bounds according to high-high vs.\ high-low
and low-high interactions. The high-high case is slightly more
involved due to the possibility of opposing $(++)$ or $(--)$ waves
with comparable frequencies and very small modulations which produce
a wave of small frequency but very large modulation.

\begin{lemma}
  \label{lem:Nablowmod}
For any $\ell\le k+O(1)$, and $\phi_j$ adapted to $k_j$ with
$k_1=k_2+O(1)$,
\begin{equation}\label{eq:Nablowell}
\|P_kQ_{\ell} \calN_{\alpha\beta} (\phi_1,\phi_2)\|_{\Ltwotx} \les
2^{\frac{\ell-k}{4+}} 2^{\frac{k}{2}}
2^{\frac{k-k_1}{2}} \|\phi_1\|_{S[k_1]} \|\phi_2\|_{S[k_2]}
\end{equation}
In particular,
\begin{equation}\label{eq:Nablow}
\|P_kQ_{\le k+C} \calN_{\alpha\beta} (\phi_1,\phi_2)\|_{\Ltwotx} \les
2^{k-\frac{k_1}{2}} \|\phi_1\|_{S[k_1]} \|\phi_2\|_{S[k_2]}
\end{equation}
Finally, for any $m_0\le -10$,
\begin{equation}
\Big(\sum_{ \kappa \in\calC_{m_0 } } \|P_kQ_{\le k+C}
\calN_{\alpha\beta} ( P_{k_1,\kappa}\phi_1,
\phi_2)\|_{\Ltwotx}^2\Big)^{\frac12} \les |m_0| \, 2^{\frac{k_1}{2}}
\|\phi_1\|_{S[k_1]} \|\phi_2\|_{S[k_2]} \label{eq:angphi2}
\end{equation}
\end{lemma}
\begin{proof}
   We can take $k_1=k_2+O(1)=0$. First, by~\eqref{eq:bilin3null},
   \[
\|P_k Q_\ell \calN_{\alpha\beta} (Q_{\le k+\ell-C}\phi_1,Q_{\le
k+\ell-C}\phi_2)\|_{\Ltwotx} \les  2^{\frac{\ell+k}{4}} 2^{\frac{k}{2}}
\|\phi_1\|_{S[k_1]} \|\phi_2\|_{S[k_2]}
   \]
Second, by an angular decomposition into caps of
size~$2^{\frac{\ell+k}{2}}$,
\begin{align}
& \sum_{\ell+k-C\le m\le \ell} \| P_k Q_\ell \calN_{\alpha\beta}
(Q_{m}\phi_1,Q_{\le m}\phi_2)\|_{\Ltwotx} \nn  \\& \les   \sum
_{\ell+k-C\le m\le \ell} 2^{\frac{\ell+k}{2}} 2^{\frac{\ell-k}{4}}
2^k \|Q_{m}\phi_1\|_{\Ltwotx} \|Q_{\le
m}\phi_2\|_{L^\infty_t L^2_x} \label{eq:sch79}  \\
&\les  2^{\frac{\ell+3k}{4}}  \|\phi_1\|_{S[k_1]}
\|\phi_2\|_{S[k_2]} \nn
\end{align}
To pass to \eqref{eq:sch79} one uses the improved Bernstein
inequality, which yields  a factor of $2^k 2^{\frac{\ell-k}{4}}$,
whereas the~$2^{\frac{\ell+k}{2}}$ corresponds to the angular gain
from the nullform (note that the error coming from the modulation is
at most~$2^m\le 2^\ell$ which is less than this gain).  And third,
by the improved Bernstein inequality and a decomposition into caps
of size~$2^{\frac{m+k}{2}}$,
\begin{align*}
 \sum_{\ell\le m\le C} \| P_k Q_\ell \calN_{\alpha\beta} (Q_{m}\phi_1,Q_{\le m}\phi_2)\|_{\Ltwotx}
&\les   \sum _{\ell\le m\le C}  2^{\frac{\ell-k}{4}} 2^k \big( 2^{\frac{m+k}{2}}+2^m\big) \|Q_{m}\phi_1\|_{\Ltwotx} \|Q_{\le
m}\phi_2\|_{L^\infty_t L^2_x}  \\
&\les  \sum _{\ell\le m\le C} 2^{\frac{\ell-k}{4}} 2^k  \big( 2^{\frac{m+k}{2}}+2^m\big) 2^{-\frac{m}{2}}   \|\phi_1\|_{S[k_1]} \|\phi_2\|_{S[k_2]}\\
&\les  2^{\frac{\ell-k}{4+}} 2^k  \|\phi_1\|_{S[k_1]} \|\phi_2\|_{S[k_2]}
\end{align*}
The factor $2^{\frac{m+k}{2}}+2^m$ here is made up out of the angular gain $2^{\frac{m+k}{2}}$ and the
loss of~$2^m$ in modulation (in case $\beta=0$).
And finally, due to $\eps<\frac12$,
\begin{align*}
  \| P_k Q_\ell \calN_{\alpha\beta} (Q_{\ge C}\phi_1,\phi_2)\|_{\Ltwotx} &\les   2^k 2^{\frac{\ell}{2}}
  \|  \calN_{\alpha\beta} (Q_{\ge C}\phi_1,\phi_2)\|_{L^1_t L^1_x}  \\
&\les \sum_{m\ge C} 2^k 2^{\frac{\ell}{2}} 2^m  \|  Q_{m}\phi_1\|_{\Ltwotx} \| \tilde Q_m \phi_2\|_{\Ltwotx} \\
&\les  2^k 2^{\frac{\ell}{2}}  \sum_{m\ge C} 2^{m} 2^{-2m(1-\eps)} \|\phi_1\|_{S[k_1]} \|\phi_2\|_{S[k_2]}\\
&\les 2^{\frac{\ell-k}{4}} 2^k  \|\phi_1\|_{S[k_1]} \|\phi_2\|_{S[k_2]}
\end{align*}
as desired.

Next, we consider~\eqref{eq:angphi2}. Here one essentially repeats the  proof of~\eqref{eq:Nablow} verbatim.
The only difference being that instead of Lemma~\ref{lem:bilin3} one uses Corollary~\ref{cor:squarebilin3}, in fact
the null-form version of~\eqref{eq:sum13}. Note that this loses a factor of~$|m_0|$. To sum over the caps one also needs
to invoke Lemma~\ref{lem:enersquaresum} in case of a~$\ener$-norm,  which incurs the same loss.
\end{proof}

We shall also require the following technical variant of the estimate of Lemma~\ref{lem:Nablowmod}.
It obtains an improvement for the case of angular alignment in the Fourier supports of the inputs.

\begin{lemma}
  \label{lem:Nablowmod'} Let $\delta>0$ be small and $L>1$ be large. Then there exists $m_0=m_0(\delta,L)<0$ large and negative  such that
for any  $\phi_j$ adapted to $k_j$ for $j=1,2$,
\begin{align}
\sum_{\substack{\kappa_1,\kappa_2\in\calC_{m_0 }\\
\dist(\kappa_1,\kappa_2)\le 2^{m_0}}} \|P_kQ_{\le k+C}
\calN_{\alpha\beta} (P_{k_1,\kappa_1} \phi_1,
P_{k_2,\kappa_2}\phi_2)\|_{\Ltwotx} &\le \delta \, 2^{\frac{k_1}{2}}
\|\phi_1\|_{S[k_1]} \|\phi_2\|_{S[k_2]} \label{eq:angphi}
\end{align}
provided $\max_{j=1,2}|k-k_j|\le L$. The constant~$C$
  is an absolute constant which does not
depend on~$L$ or~$\delta$.
\end{lemma}
\begin{proof}
Set $k=0$. We first note that summing~\eqref{eq:Nablowell} over~$\ell\le -B$
already yields an improvement over~\eqref{eq:Nablow} provided~$B$ is large enough (in relation to~$\delta$ and~$L$).
Hence it suffices to consider the contribution of $P_0Q_\ell\calN_{\alpha\beta} (P_{k_1,\kappa_1} \phi_1, P_{k_2,\kappa_2}\phi_2)$
with $-B\le\ell\le O(1)$ fixed.
First, if we choose $m_0 $ to be a sufficiently large negative integer, then
\[
\sum_{\substack{\kappa_1,\kappa_2\in\calC_{m_0 }\\ \dist(\kappa_1,\kappa_2)\le 2^{m_0}}}
 P_0 Q_\ell \calN_{\alpha\beta} (Q_{\le \ell-C} P_{k_1,\kappa_1}  \phi_1,\: Q_{\le
\ell-C} P_{k_2,\kappa_2} \phi_2) = 0
\]
by Lemma~\ref{lem:cone}.
Second, by an angular decomposition into caps of
size~$2^{\frac{\ell}{2}}$,
\begin{align*}
 & \sum_{\substack{\kappa_1,\kappa_2\in\calC_{m_0}\\ \dist(\kappa_1,\kappa_2)\le 2^{m_0}}}  \sum_{\ell-C\le m\le C} \| P_0 Q_\ell \calN_{\alpha\beta}
(Q_{m}P_{k_1,\kappa_1} \phi_1,Q_{\le m}P_{k_2,\kappa_2} \phi_2)\|_{\Ltwotx}  \\
&\le C(L,\delta)  \sum_{\substack{\kappa_1,\kappa_2\in\calC_{m_0}\\ \dist(\kappa_1,\kappa_2)\le 2^{m_0}}}  \sum _{\ell-C\le
m\le C}
\|Q_{m} P_{k_1,\kappa_1} \phi_1\|_{\Ltwotx} \|Q_{\le
m} P_{k_2,\kappa_2} \phi_2\|_{L^\infty_t L^\infty_x}  \\
&\le C(L,\delta)\,|m_0|\,  2^{\frac{m_0}{2}}    \|\phi_1\|_{S[k_1]} \|\phi_2\|_{S[k_2]} \le \delta \, \|\phi_1\|_{S[k_1]} \|\phi_2\|_{S[k_2]}
\end{align*}
To pass to the last line we applied Cauchy-Schwarz to the sum over the caps as well as Lemma~\ref{lem:enersquaresum}.
The case dealing with~$Q_{\le
m} P_{k_1,\kappa_1} \phi_1$ and~$Q_{m} P_{k_2,\kappa_2} \phi_2$ is analogous.
And finally, due to $\eps<\frac12$,
\begin{align*}
 & \sum_{\substack{\kappa_1,\kappa_2\in\calC_{m_0}\\ \dist(\kappa_1,\kappa_2)\le 2^{m_0}}}
 \| P_0 Q_\ell \calN_{\alpha\beta} (Q_{\ge C} P_{k_1,\kappa_1} \phi_1, P_{k_2,\kappa_2} \phi_2)\|_{\Ltwotx} \\
& \les \sum_{\substack{\kappa_1,\kappa_2\in\calC_{m_0}\\ \dist(\kappa_1,\kappa_2)\le 2^{m_0}}}
 \| P_0 Q_\ell \calN_{\alpha\beta} (Q_{\ge C} P_{k_1,\kappa_1} \phi_1, P_{k_2,\kappa_2} \phi_2)\|_{L^1_t L^2_x} \\
&\le C(L,\delta) \sum_{m\ge C}  2^m  \sum_{\substack{\kappa_1,\kappa_2\in\calC_{m_0}\\ \dist(\kappa_1,\kappa_2)\le 2^{m_0}}}
\|  Q_{m}P_{k_1,\kappa_1} \phi_1\|_{\Ltwotx} \| \tilde Q_m P_{k_2,\kappa_2} \phi_2\|_{L^2_t L^\infty_x} \\
&\le C(L,\delta)  2^{\frac{m_0}{2}}   \sum_{m\ge C} 2^{m} 2^{-2m(1-\eps )}    \|  \phi_1\|_{S[k_1]} \|\phi_2\|_{S[k_2]}
\le \delta  \|\phi_1\|_{S[k_1]} \|\phi_2\|_{S[k_2]}
\end{align*}
as desired.
\end{proof}

In case the output has ``elliptic'' rather than hyperbolic character, there
is the following bound.

\begin{lemma}
 \label{lem:Nabhighmod}
For any  $\phi_j$ adapted to $k_j$ with
$k_1=k_2+O(1)$,
\[
 \sum_{\ell\ge k+C} 2^{-\eps\ell}  \|P_kQ_{\ell} \calN_{\alpha\beta} (\phi_1,\phi_2)\|_{\Ltwotx} \les
2^{\frac{k}{2}} 2^{-\eps k_1} \la k_1-k\ra^2 \|\phi_1\|_{S[k_1]} \|\phi_2\|_{S[k_2]}
\]
Furthermore,
\begin{equation}\label{eq:sch113}
 \sum_{\ell\ge k+C}   \|P_kQ_{\ell} \calN_{\alpha\beta} (Q_{\le k_1+C} \phi_1, Q_{\le k_2+C} \phi_2)\|_{\Ltwotx} \les
2^{\frac{k}{2}}  \la k_1-k\ra^2 \|\phi_1\|_{S[k_1]} \|\phi_2\|_{S[k_2]}
\end{equation}
\end{lemma}
\begin{proof} We set $k_1=k_2+O(1)=0$. One has the decomposition
 \begin{align}
  \|P_kQ_{\ell} \calN_{\alpha\beta} (\phi_1,\phi_2)\|_{\Ltwotx}
&\les \|P_k Q_{\ell} \calN_{\alpha\beta} (Q_{\ge \ell-C} \phi_1, Q_{\le k_1+C} \phi_2)\|_{\Ltwotx} \label{eq:largemod1}\\
& \quad +\|P_k Q_{\ell} \calN_{\alpha\beta} (Q_{\ge \ell-C} \phi_1, Q_{> k_1+C}\phi_2)\|_{\Ltwotx}  \label{eq:largemod2} \\
&\quad + \|P_k Q_{\ell} \calN_{\alpha\beta} (Q_{<\ell-C}\phi_1, Q_{\ge\ell-C}\phi_2)\|_{\Ltwotx} \label{eq:largemod3} \\
&\quad + \|P_k Q_{\ell} \calN_{\alpha\beta} (Q_{<\ell-C}\phi_1, Q_{<\ell-C}\phi_2)\|_{\Ltwotx} \label{eq:largemod4}
 \end{align}
We begin with the estimate
\[
\sum_{\substack{\ell\ge k+C\\ \ell\ne k_1+O(1)}}  2^{-\eps\ell}
\|P_kQ_{\ell} \calN_{\alpha\beta} (\phi_1,\phi_2)\|_{\Ltwotx} \les
2^{k}    \|\phi_1\|_{S[k_1]} \|\phi_2\|_{S[k_2]}
\]
which is stronger than what we claim -- this is due to the fact that
the the case of opposing~$(++)$ and~$(--)$ waves is excluded in this
sum.  We first consider the case $\ell\le k_1- C'$ where $C'$ is
large but still smaller than the constant~$C$
in~\eqref{eq:largemod1}--\eqref{eq:largemod4}. Then the term
in~\eqref{eq:largemod4} vanishes. On the one hand,
\begin{align*}
  \eqref{eq:largemod1} &\les 2^k \|P_k Q_{\ell} [\del_\beta( Q_{\ge \ell-C} |\nabla|^{-1} \phi_1\cdot Q_{\le k_2+C}\del_\alpha|\nabla|^{-1} \phi_2)
\\ &\qquad - \del_\alpha( Q_{\ge \ell-C} |\nabla|^{-1} \phi_1 \cdot Q_{\le k_2+C}\del_\beta|\nabla|^{-1} \phi_2) ] \|_{\Lzweins} \\
&\les  2^{k+\ell} \| Q_{\ge \ell-C} \phi_1\|_{\Ltwotx} \|\phi_2\|_{S[k_2]}
\end{align*}
Here we used that
\[
 \|Q_{\le k_2+C}\del_\beta|\nabla|^{-1} \phi_2\|_{\ener} \les \|\phi_2\|_{S[k_2]}
\]
Furthermore,
\begin{align*}
 \sum_{k+C\le \ell\le k_1-C}\sum_{m\ge \ell-C}  2^{(1-\eps)\ell}
2^{k-k_1}  \| Q_{m} \phi_1\|_{\Ltwotx} \|\phi_2\|_{S[k_2]} &\les
2^{k} \|\phi_1\|_{\dot X^{0,1-\eps,2}_{k_1}} \|\phi_2\|_{S[k_2]}\\
&\les 2^{k} \|\phi_1\|_{S[k_1]} \|\phi_2\|_{S[k_2]}
\end{align*}
as desired. The term \eqref{eq:largemod3} satisfies the same bound, more precisely, it can be
reduced to~\eqref{eq:largemod1}, \eqref{eq:largemod2}.
Next, note that due to $\ell\le k_1-C'$ it suffices to
consider~$\phi_1=Q_{\ge k_1+C}\phi_1$ in~\eqref{eq:largemod2}.
Consequently,
\begin{align*}
\sum_{k+C\le \ell\le k_1-C} 2^{-\eps\ell} \eqref{eq:largemod2} &\les
\sum_{k+C\le \ell\le k_1-C}2^{-\eps\ell}
\sum_{m\ge k_1+C} \|P_k Q_{\ell} \calN_{\alpha\beta} (Q_{m} \phi_1, \tilde Q_m\phi_2)\|_{\Ltwotx}\\
&\les \sum_{k+C\le \ell\le k_1-C}2^{-\eps\ell} \sum_{m\ge k_1+C}  2^k 2^{\frac{\ell}{2}} \|\calN_{\alpha\beta} (Q_{m} \phi_1, \tilde Q_m\phi_2)\|_{\Leins} \\
&\les \sum_{k+C\le \ell\le k_1-C} 2^{-\eps\ell} \sum_{m\ge k_1+C}  2^k 2^{\frac{\ell}{2}} 2^{m} 2^{-2m(1-\eps)}  \|\phi_1\|_{S[k_1]} \|\phi_2\|_{S[k_2]}\\
&\les 2^{k}  \|\phi_1\|_{S[k_1]} \|\phi_2\|_{S[k_2]}
\end{align*}
 where we used that $\eps<\frac14$ in the final step.
Second, suppose that $\ell\ge k_1+C'$. Then
\begin{align}
& \sum_{\ell\ge k_1+C'}  2^{-\eps\ell}  \|P_kQ_{\ell}
\calN_{\alpha\beta} (\phi_1,\phi_2)\|_{\Ltwotx} \nn \\ &\les
\sum_{\ell\ge k_1+C'}  2^{-\eps\ell}  \|P_kQ_{\ell}
\calN_{\alpha\beta} (\tilde
Q_\ell\phi_1,Q_{\le\ell-5}\phi_2)\|_{\Ltwotx}
\label{eq:high1} \\
&\quad  +\sum_{\ell\ge k_1+C'}  2^{-\eps\ell}  \|P_kQ_{\ell} \calN_{\alpha\beta} (Q_{\le\ell-5}\phi_1,\tilde Q_\ell\phi_2)\|_{\Ltwotx} \label{eq:high2} \\
&\quad  +\sum_{\ell\ge k_1+C'}  2^{-\eps\ell} \sum_{m\ge\ell-5}
\|P_kQ_{\ell} \calN_{\alpha\beta} (Q_{m}\phi_1,\tilde
Q_m\phi_2)\|_{\Ltwotx}  \label{eq:high3}
\end{align}
which are in turn estimated as follows:
\begin{align*}
 \eqref{eq:high1} &\les  \sum_{\ell\ge k_1+C'}  2^{-\eps\ell} 2^{k}  \|P_kQ_{\ell} \calN_{\alpha\beta} (\tilde Q_\ell\phi_1,
 Q_{\le\ell-5}\phi_2)\|_{L^2_t L^1_x} \\
&\les \sum_{\ell\ge k_1+C'}  2^{-\eps\ell} 2^{k} 2^{\ell}  \|\tilde Q_\ell\phi_1\|_{\Ltwotx} \|\phi_2\|_{\ener}\\
&\les 2^{k} \|\phi_1\|_{S[k_1]} \|\phi_2\|_{S[k_2]}
\end{align*}
and similarly for~\eqref{eq:high2}, whereas \eqref{eq:high3} is bounded by
\begin{align*}
 &\les  \sum_{\ell\ge k_1+C'}  2^{-\eps\ell} \sum_{m\ge\ell-5} 2^{\frac{\ell}{2}} 2^k  \| \calN_{\alpha\beta} (Q_{m}\phi_1,\tilde Q_m\phi_2)\|_{\Leins} \\
&\les  \sum_{\ell\ge k_1+C'}  2^{-\eps\ell} \sum_{m\ge\ell-5}
2^{\frac{\ell}{2}} 2^k 2^{m} 2^{-2m(1-\eps)} \|\phi_1\|_{S[k_1]}
\|\phi_2\|_{S[k_2]} \\
&\les  2^{k}  \|\phi_1\|_{S[k_1]}
\|\phi_2\|_{S[k_2]}
\end{align*}
as desired.

\noindent It remains to consider the case $|\ell-k_1|=|\ell|\le C'$,
which gives us the weaker bound stated in the lemma. We use the
decomposition \eqref{eq:largemod1}--\eqref{eq:largemod4}. The terms
\eqref{eq:largemod1}--\eqref{eq:largemod3} give a bound of $2^k$ as
before. The main difference lies with~\eqref{eq:largemod4} which is
nonzero only due to the contribution to opposing $(++)$ or $(--)$
waves, see Lemma~\ref{lem:cone}. In fact, one has
with~$\ell=k_1+O(1)=O(1)$,
\begin{align}
 & \|P_k Q_{O(1)} \calN_{\alpha\beta} (Q_{<-C}\phi_1, Q_{<-C}\phi_2)\|_{\Ltwotx} \nn \\
&\les \sum_{\pm} \sum_{\kappa\in\calC_{k}}  \|P_k Q_{O(1)} \calN_{\alpha\beta} (Q_{<-C}P_{\kappa}
\phi_1^{\pm}, Q_{<-C}P_{-\kappa}\phi_2^{\pm})\|_{\Ltwotx} \nn \\
&\les \sum_{\pm} \sum_{\kappa\in\calC_{k}}  2^{\frac{k}{2}} \|Q_{<-C}P_{\kappa}
\phi_1^{\pm}\|_{S[k_1,\kappa]} \|Q_{<-C}P_{-\kappa}\phi_2^{\pm}\|_{S[k_2,-\kappa]} \label{eq:oppose} \\
&\les 2^{\frac{k}{2}}  \sum_{\pm} \Big(\sum_{\kappa\in\calC_{k}}
\|Q_{<-C}P_{\kappa} \phi_1^{\pm}\|^2_{S[k_1,\kappa]}\Big)^{\frac12}
\Big(\sum_{\kappa\in\calC_{k}}
\|Q_{<-C}P_{-\kappa}\phi_2^{\pm}\|_{S[k_2,-\kappa]}^2\Big)^{\frac12} \nn \\
&\les 2^{\frac{k}{2}}  k^2   \|\phi_1\|_{S[k_1]}
\|\phi_2\|_{S[k_2]}\nn
\end{align}
To pass to the last line,  we wrote
\begin{align}
 \sum_{\kappa\in\calC_{k}} \|Q_{<-C}P_{\kappa}
 \phi_1^{\pm}\|_{S[k_1,\kappa]}^2
 &\les \sum_{\kappa\in\calC_{k}} \|Q_{<2k}P_{\kappa} \phi_1^{\pm}\|_{S[k_1,\kappa]}^2
+\sum_{\kappa\in\calC_{k}} \Big(\sum_{2k\le j\le -C} \|Q_{j}P_{\kappa} \phi_1^{\pm}\|_{S[k_1,\kappa]}\Big)^{2} \nn \\
&\les \| \phi_1^{\pm}\|_{S[k_1]}^2
+|k| \sum_{2k\le j\le -C} \sum_{\kappa\in\calC_{k}} \|Q_j P_{\kappa} \phi_1^{\pm}\|^2_{\dot X^{0,\frac12,\infty}_{k_1}}\label{eq:opposeloss}\\
&\les  |k|^2 \| \phi_1^{\pm}\|_{S[k_1]}^2 \nn
\end{align}
and the result follows.

\noindent The second statement~\eqref{eq:sch113} follows by essentially the same proof.
\end{proof}

\begin{remark}\label{rem:log_loss}
It is important to note that the logarithmic loss of $\la k_1-k\ra^2$ in~\eqref{eq:sch113} only
results from the case of opposing waves in the high-high case. Later we will use~\eqref{eq:sch113}
without this loss in those cases where these interactions are excluded.
\end{remark}

Later, we shall also require the following  technical
refinement of Lemma~\ref{lem:Nabhighmod} dealing with a further angular restriction of the first input.

\begin{cor}
  \label{cor:Nabhighmod} Under the assumptions of Lemma~\ref{lem:Nabhighmod} and for any $m_0\le -10$,
\begin{align*}
 \Big(\sum_{\kappa\in\calC_{m_0} } \Big(\sum_{\ell\ge k+C} 2^{-\eps\ell}
\|P_kQ_{\ell} \calN_{\alpha\beta}
(P_{k_1,\kappa}\phi_1,\phi_2)\|_{\Ltwotx}\Big)^2 \Big)^{\frac12}
 & \les |m_0|\, 2^{\frac{k}{2}} 2^{-\eps k_1} \la k_1-k\ra^2
\|\phi_1\|_{S[k_1]} \|\phi_2\|_{S[k_2]}
\end{align*}
with an absolute implicit constant.
\end{cor}
\begin{proof}
This can be seen by reviewing
the proof\footnote{It is important to observe that one
{\em cannot} square sum the bound of Lemma~\ref{lem:Nabhighmod}
directly due to the fact that $\sum_{\kappa\in\calC_{m_0}}
\|P_{k_1,\kappa} \phi_1\|_{S[k]}^2$ cannot be controlled.} of Lemma~\ref{lem:Nabhighmod}. Specifically, up until~\eqref{eq:oppose},
one places $P_{k_1,\kappa} \phi_1$ either in the $\dot X^{s,b}$ or~$\ener$ norms. The norms
are amenable to square summation, in the latter case at the expense of a factor~$|m_0|$, see Lemma~\ref{lem:enersquaresum}.
However, as far as~\eqref{eq:oppose} is concerned, we distinguish two cases: $k\le m_0$ and $k>m_0$. In the former case,
the caps in~$\calC_{k}$ are smaller than those in~$\calC_{m_0}$ and~\eqref{eq:oppose} applies directly (one organizes the caps
in~$\calC_{k}$ into subsets of the larger~$\calC_{m_0}$--caps). In the latter case, however, the
$\calC_{m_0}$--caps are smaller which forces us to write
\[
 Q_{<-C}\phi_1= Q_{<2m_0} \phi_1 + Q_{2m_0<\cdot<-C}\phi_1
\]
The former is subsumed in a square-function bound as
in~\eqref{eq:oppose}, whereas the latter leads to a loss of~$|m_0|$
as in~\eqref{eq:opposeloss} and the corollary is proved.
\end{proof}

Next, we obtain  an improvement in case of angular alignment of the inputs. This is analogous the case of low modulations,
see Lemma~\ref{lem:Nablowmod'}.

\begin{lemma}
 \label{lem:Nabhighmod'}
Let $\delta>0$ be small and $L>1$ be large. Then there exists $m_0=m_0(\delta,L)<0$ large and negative  such that
for any  $\phi_j$ adapted to $k_j$ for $j=1,2$,
\begin{equation}\label{eq:Nabhighmoddelta}
 \sum_{\substack{\kappa_1,\kappa_2\in\calC_{m_0 }\\ \dist(\kappa_1,\kappa_2)\le 2^{m_0}}} \sum_{\ell\ge k+C} 2^{-\eps\ell}
 \|P_kQ_{\ell} \calN_{\alpha\beta} (P_{k_1,\kappa_1} \phi_1, P_{k_2,\kappa_2}  \phi_2)\|_{\Ltwotx} \le \delta \,
  2^{(\frac12-\eps) k_1}   \|\phi_1\|_{S[k_1]} \|\phi_2\|_{S[k_2]}
\end{equation}
provided $\max_{j=1,2}|k-k_j|\le L$. The constant~$C$
in~\eqref{eq:Nabhighmoddelta} is an absolute constant which does not
depend on~$L$ or~$\delta$.
\end{lemma}
\begin{proof} The proof consists of checking that one can glean a
gain from angular alignment by following the proof of
Lemma~\ref{lem:Nabhighmod}. In effect, this will always be done by
means of  Bernstein's  inequality. The only case where this is not
possible is~\eqref{eq:oppose}, but that case is excluded by the
angular alignment assumption.

We set $k_1=0$ whence $|k|\le L$ and~$|k_2|\le 2L$. Implicit
constants here will be allowed to depend on~$L$, but not the
constants $C$ appearing in modulation cutoffs. As before, one has
the decomposition
 \begin{align}
  \|P_kQ_{\ell} \calN_{\alpha\beta} (P_{k_1,\kappa_1}\phi_1, P_{k_2,\kappa_2} \phi_2)\|_{\Ltwotx}
&\les \|P_k Q_{\ell} \calN_{\alpha\beta} (Q_{\ge \ell-C}P_{k_1,\kappa_1} \phi_1, Q_{\le k_1+C} P_{k_2,\kappa_2} \phi_2)\|_{\Ltwotx} \label{eq:largemod1'}\\
& \quad +\|P_k Q_{\ell} \calN_{\alpha\beta} (Q_{\ge \ell-C} P_{k_1,\kappa_1}\phi_1, Q_{> k_1+C} P_{k_2,\kappa_2} \phi_2)\|_{\Ltwotx}  \label{eq:largemod2'} \\
&\quad + \|P_k Q_{\ell} \calN_{\alpha\beta} (Q_{<\ell-C}P_{k_1,\kappa_1}\phi_1, Q_{\ge\ell-C} P_{k_2,\kappa_2} \phi_2)\|_{\Ltwotx} \label{eq:largemod3'} \\
&\quad + \|P_k Q_{\ell} \calN_{\alpha\beta}
(Q_{<\ell-C}P_{k_1,\kappa_1} \phi_1, Q_{<\ell-C}P_{k_2,\kappa_2}
\phi_2)\|_{\Ltwotx} \label{eq:largemod4'}
 \end{align}
We first consider the case $\ell\le k_1- C'$ where $C'$ is large but
still smaller than the constant~$C$
in~\eqref{eq:largemod1'}--\eqref{eq:largemod4'}. Then the term
in~\eqref{eq:largemod4'} vanishes by Lemma~\ref{lem:cone}. By
Bernstein's inequality,
\begin{align*}
 \sum_{k+C\le \ell\le k_1-C} 2^{-\eps\ell} \eqref{eq:largemod1'} &\les  \sum_{k+C\le \ell\le
 k_1-C}
  \|P_k Q_{\ell} [\del_\beta( Q_{\ge \ell-C} |\nabla|^{-1} P_{k_1,\kappa_1} \phi_1\cdot Q_{\le k_2+C}\del_\alpha|\nabla|^{-1}
  P_{k_2,\kappa_2} \phi_2)
\\ &\qquad - \del_\alpha( Q_{\ge \ell-C} |\nabla|^{-1} P_{k_1,\kappa_1} \phi_1 \cdot Q_{\le k_2+C}\del_\beta|\nabla|^{-1}
P_{k_2,\kappa_2} \phi_2) ] \|_{\Ltwotx} \\
&\les  \sum_{k+C\le \ell\le k_1-C}   \| Q_{\ge \ell-C}
P_{k_1,\kappa_1}  \phi_1\|_{\Ltwotx} \| P_{k_2,\kappa_2}
\phi_2\|_{\Linf}\\
&\les 2^{\frac{m_0}{2}} \sum_{k+C\le \ell\le k_1-C}   \| Q_{\ge
\ell-C} P_{k_1,\kappa_1}  \phi_1\|_{\Ltwotx} \|
P_{k_2,\kappa_2} \phi_2\|_{\ener} \\
&\le \delta \|P_{k_1,\kappa_1} \phi_1\|_{\dot
X_0^{0,\frac12,\infty}} \|P_{k_2,\kappa_2} \phi_2\|_{\ener}
\end{align*}
Summing over $\kappa_1,\kappa_2$ now yields the desired bound by
Cauchy-Schwarz (see also Lemma~\ref{lem:enersquaresum}).   The term \eqref{eq:largemod3'} satisfies the same
bound. Next, note that due to $\ell\le k_1-C'$ it suffices to
consider~$\phi_1=Q_{\ge k_1+C}\phi_1$ in~\eqref{eq:largemod2'}.
Consequently,
\begin{align*}
\sum_{k+C\le \ell\le k_1-C} 2^{-\eps\ell} \eqref{eq:largemod2'}
&\les \sum_{k+C\le \ell\le k_1-C}  \sum_{m\ge k_1+C}
\|P_k Q_{\ell} \calN_{\alpha\beta} (Q_{m} P_{k_1,\kappa_1} \phi_1, \tilde Q_m P_{k_2,\kappa_2}\phi_2)\|_{\Ltwotx}\\
&\les \sum_{k+C\le \ell\le k_1-C}  \sum_{m\ge k_1+C}
  \|\calN_{\alpha\beta} (Q_{m}P_{k_1,\kappa_1}
\phi_1, \tilde Q_m
P_{k_2,\kappa_2} \phi_2)\|_{L^1_t L^2_x} \\
&\les \sum_{k+C\le \ell\le k_1-C}   \sum_{m\ge k_1+C}
  2^{m} 2^{-2m(1-\eps)}  \|P_{k_1,\kappa_1} Q_m
\phi_1\|_{\Ltwotx}
 \|P_{k_2,\kappa_2}\tilde Q_m \phi_2\|_{L^2_t L^\infty_x}\\
&\les \delta   \|P_{k_1,\kappa_1}\phi_1\|_{\dot X_0^{0,1-\eps,2}}
\|P_{k_2,\kappa_2}\phi_2\|_{\dot X_0^{0,1-\eps,2}}
\end{align*}
Summing over $\kappa_1,\kappa_2$ again leads to the desired bound.
Second, suppose that $\ell\ge k_1+C'$. Then
\begin{align}
& \sum_{\ell\ge k_1+C'}  2^{-\eps\ell}  \|P_kQ_{\ell}
\calN_{\alpha\beta}
(P_{k_1,\kappa_1}\phi_1,P_{k_2,\kappa_2}\phi_2)\|_{\Ltwotx} \nn
\\ &\les \sum_{\ell\ge k_1+C'}  2^{-\eps\ell}      \|P_kQ_{\ell} \calN_{\alpha\beta}
(\tilde Q_\ell P_{k_1,\kappa_1}
\phi_1,Q_{\le\ell-5}P_{k_2,\kappa_2}\phi_2)\|_{\Ltwotx}
\label{eq:high1'} \\
&\quad  +\sum_{\ell\ge k_1+C'} 2^{-\eps\ell}      \|P_kQ_{\ell} \calN_{\alpha\beta}
(Q_{\le\ell-5}P_{k_1,\kappa_1}\phi_1,\tilde Q_\ell P_{k_2,\kappa_2}\phi_2)\|_{\Ltwotx} \label{eq:high2'} \\
&\quad  +\sum_{\ell\ge k_1+C'}  2^{-\eps\ell}    \sum_{m\ge\ell-5}
\|P_kQ_{\ell} \calN_{\alpha\beta}
(Q_{m}P_{k_1,\kappa_1}\phi_1,\tilde
Q_mP_{k_2,\kappa_2}\phi_2)\|_{\Ltwotx} \label{eq:high3'}
\end{align}
which are in turn estimated as follows:
\begin{align*}
 \sum_{\substack{\kappa_1,\kappa_2\in\calC_{m_0 }\\ \dist(\kappa_1,\kappa_2)\le 2^{m_0}}} \eqref{eq:high1'}
 &\les  \sum_{\substack{\kappa_1,\kappa_2\in\calC_{m_0 }\\ \dist(\kappa_1,\kappa_2)\le 2^{m_0}}}
 \sum_{\ell\ge k_1+C'}  2^{-\eps\ell}    \|P_kQ_{\ell} \calN_{\alpha\beta} (\tilde Q_\ell P_{k_1,\kappa_1} \phi_1,
 Q_{\le\ell-5}P_{k_2,\kappa_2}\phi_2)\|_{L^2_t L^2_x} \\
&\les \sum_{\substack{\kappa_1,\kappa_2\in\calC_{m_0 }\\
\dist(\kappa_1,\kappa_2)\le 2^{m_0}}}
\sum_{\ell\ge k_1+C'}  2^{(1-\eps)\ell}     \|\tilde Q_\ell P_{k_1,\kappa_1} \phi_1\|_{\Ltwotx} \| P_{k_2,\kappa_2}\phi_2\|_{\ener}\\
&\le \delta  \|\phi_1\|_{S[k_1]} \|\phi_2\|_{S[k_2]}
\end{align*}
and similarly for~\eqref{eq:high2'}, whereas
\begin{align*}
\sum_{\substack{\kappa_1,\kappa_2\in\calC_{m_0 }\\
\dist(\kappa_1,\kappa_2)\le 2^{m_0}}} \eqref{eq:high3'}  &\les
 \sum_{\substack{\kappa_1,\kappa_2\in\calC_{m_0 }\\
\dist(\kappa_1,\kappa_2)\le 2^{m_0}}}\!\!\!\! \sum_{\ell\ge k_1+C'}
2^{-\eps\ell} \sum_{m\ge\ell-5} 2^{\frac{\ell}{2}}
 \| \calN_{\alpha\beta} (Q_{m}P_{k_1,\kappa_1}\phi_1,\tilde Q_m P_{k_2,\kappa_2}\phi_2)\|_{L^1_t L^2_x} \\
&\les 2^{\frac{m_0}{2}} \!\!\!\!\!\!\!\!\!\sum_{\substack{\kappa_1,\kappa_2\in\calC_{m_0 }\\
\dist(\kappa_1,\kappa_2)\le 2^{m_0}}} \!\!\!\sum_{m+5\ge \ell\ge
k_1+C'} 2^{(\frac12-\eps)\ell}      2^{m} 2^{-2m(1-\eps)} \|
P_{k_1,\kappa_1}\phi_1\|_{S[k_1]}
\|P_{k_2,\kappa_2} \phi_2\|_{S[k_2]} \\
&\le \delta  \| \phi_1\|_{S[k_1]} \| \phi_2\|_{S[k_2]}
\end{align*}
as desired.

\noindent It the remaining case $|\ell-k_1|=|\ell|\le C'$ we use the
decomposition \eqref{eq:largemod1'}--\eqref{eq:largemod4'}. The
terms \eqref{eq:largemod1'}--\eqref{eq:largemod3'} give a bound of
$\delta$ as before. The main difference lies
with~\eqref{eq:largemod4'} which is nonzero only due to the
contribution to opposing $(++)$ or $(--)$ waves, see
Lemma~\ref{lem:cone}. However, this case is excluded due to the
angular alignment assumption.
\end{proof}

\subsection{Nullform bounds in the low-high and high-low cases}

We now derive analogues of the previous two lemmas in the high-low case,
with the low-high case being completely analogous.

\begin{lemma}
  \label{lem:Nablowmod2}
For any
$\phi_j$ adapted to $k_j$ with
$k_2\le k_1+O(1)=k$ one has
\[
\|P_kQ_{\le k+C} \calN_{\alpha\beta} (\phi_1,\phi_2)\|_{\Ltwotx} \les
2^{(\frac{1}{2}-\eps)k_2} 2^{\eps k_1}  \|\phi_1\|_{S[k_1]} \|\phi_2\|_{S[k_2]}
\]
\end{lemma}
\begin{proof} We may take $k=k_1+O(1)=0$ and $k_2\le -C$.
  Assume first that $Q_{<k_2}\phi_i=\phi_i$ for $i=1,2$.  Then the modulation
of the output does not exceed $2^{k_2}$, and we are
reduced to bounding the following three expressions:
 \begin{align}
 & \sum_{j\le k_2+O(1)}\!\!\! \|P_0 Q_j ( R_\alpha Q_{<j-C} \phi_1 R_\beta Q_{<j-C}
  \phi_2 - R_\beta Q_{<j-C} \phi_1R_\alpha Q_{<j-C}
  \phi_2) \|_{\Ltwotx} \label{eq:3sum1} \\
  &+ \sum_{j\le k_2+O(1)} \|P_0 Q_{<j+C}  (R_\alpha Q_{j} \phi_1 R_\beta Q_{\le j}
  \phi_2 - R_\beta Q_{j} \phi_1R_\alpha Q_{\le j}
  \phi_2) \|_{\Ltwotx} \label{eq:3sum2} \\
&+ \sum_{j\le k_2+O(1)} \|P_0 Q_{<j+C}  (R_\alpha Q_{<j} \phi_1 R_\beta Q_{ j}
  \phi_2 - R_\beta Q_{<j} \phi_1R_\alpha Q_{j}
  \phi_2) \|_{\Ltwotx} \label{eq:3sum3}
\end{align}
Each of the summands here is bounded by $2^{(j+k_2)/4}$.  For the first, one decomposes into caps of size~$2^{\frac{j-k_2}{2}}$:
\begin{align*}
 \eqref{eq:3sum1} &\les \sum_{j\le k_2+O(1)}\;\;  \sum_{\kappa\sim\kappa'\in\calC_{\frac{j-k_2}{2}}}
 \|P_{0} Q_j ( R_\alpha Q_{<j-C} P_\kappa \phi_1 R_\beta Q_{<j-C} P_{\kappa'}
  \phi_2 \\
&\qquad - R_\beta Q_{<j-C} P_{\kappa} \phi_1R_\alpha Q_{<j-C}  P_{\kappa'}
  \phi_2) \|_{\Ltwotx} \\
&\les  \sum_{j\le k_2+O(1)} 2^{\frac{j-k_2}{4}} 2^{\frac{k_2}{2}}   \| Q_{<j-C} P_\kappa \phi_1\|_{S[k_1,\kappa]} \| Q_{<j-C} P_{\kappa'}
  \phi_2 \|_{S[k_2,\kappa']} \\
&\les 2^{\frac{k_2}{2}} \|\phi_1\|_{S[k_1]} \|\phi_2\|_{S[k_2]}
\end{align*}
where we applied \eqref{eq:bilin1} in the last step. Note that the nullform gains a factor of the angle
in this bound.  As for~\eqref{eq:3sum2}, one performs a similar cap decomposition but without the separation
between the caps:
\begin{align*}
 \eqref{eq:3sum2} &\les \sum_{j\le k_2+O(1)}\;\;  \sum_{\kappa,\kappa'\in\calC_{\frac{j-k_2}{2}}}
 \|P_{0} Q_{<j+C} ( R_\alpha Q_{j} P_\kappa \phi_1 R_\beta Q_{\le j} P_{\kappa'}
  \phi_2 \\
&\qquad - R_\beta Q_{j} P_{\kappa} \phi_1R_\alpha Q_{\le j}  P_{\kappa'}
  \phi_2) \|_{\Ltwotx} \\
&\les  \sum_{j\le k_2+O(1)} 2^{\frac{j-k_2}{2}}    \| Q_{j} P_\kappa \phi_1\|_{\Ltwotx} \| Q_{\le j} P_{\kappa'}
  \phi_2 \|_{L^\infty_t L^\infty_x} \\
&\les  \sum_{j\le k_2+O(1)} 2^{\frac{j-k_2}{2}} 2^{-\frac{j}{2}}
2^{\frac{j-k_2}{4}} 2^{k_2}    \| Q_{j} P_\kappa
\phi_1\|_{S[k_1,\kappa]} \|  P_{\kappa'}
  \phi_2 \|_{L^\infty_t L^2_x} \\
&\les \sum_{j\le k_2+O(1)}  2^{\frac{j+k_2}{4}} \|\phi_1\|_{S[k_1]} \|\phi_2\|_{S[k_2]}  \les 2^{\frac{k_2}{2}} \|\phi_1\|_{S[k_1]} \|\phi_2\|_{S[k_2]}
\end{align*}
Finally,
\begin{align*}
 \eqref{eq:3sum3} &\les \sum_{j\le k_2+O(1)}\;\;  \sum_{\kappa,\kappa'\in\calC_{\frac{j-k_2}{2}}}
 \|P_{0} Q_{<j+C} ( R_\alpha Q_{<j} P_\kappa \phi_1 R_\beta Q_{j} P_{\kappa'}
  \phi_2 \\
&\qquad - R_\beta Q_{<j} P_{\kappa} \phi_1R_\alpha Q_{j}  P_{\kappa'}
  \phi_2) \|_{\Ltwotx} \\
&\les  \sum_{j\le k_2+O(1)} 2^{\frac{j-k_2}{2}}    \| Q_{<j} P_\kappa \phi_1\|_{L^\infty_t L^2_x} \| Q_{j} P_{\kappa'}
  \phi_2 \|_{L^2_t L^\infty_x} \\
&\les  \sum_{j\le k_2+O(1)} 2^{\frac{j-k_2}{2}}  2^{\frac{j-k_2}{4}} 2^{k_2}    \| Q_{<j} P_\kappa \phi_1\|_{S[k_1,\kappa]} \|Q_j  P_{\kappa'}
  \phi_2 \|_{\Ltwotx} \\
&\les \sum_{j\le k_2+O(1)}  2^{\frac{j+k_2}{4}} \|\phi_1\|_{S[k_1]} \|\phi_2\|_{S[k_2]}  \les 2^{\frac{k_2}{2}} \|\phi_1\|_{S[k_1]} \|\phi_2\|_{S[k_2]}
\end{align*}
If $Q_{k_2\le \cdot\le C} \phi_2=\phi_2$, then we may take $\phi_1=Q_{\le C} \phi_1$ whence
\begin{align*}
  \|P_0Q_{\le O(1)} \calN_{\alpha\beta} (\phi_1,\phi_2)\|_{\Ltwotx}
&\les \|\phi_1\|_{\ener} \|R_0 \phi_2\|_{L^2_t L^\infty_x} \\
&\les \|\phi_1\|_{\ener}  \sum_{k_2\le j\le C} 2^j \|Q_j
\phi_2\|_{\Ltwotx}  \\&  \les 2^{(\frac12-\eps)k_2}
\|\phi_1\|_{S[k_1]} \|\phi_2\|_{S[k_2]}
\end{align*}
On the other hand, if $\phi_2=Q_{\ge C} \phi_2$, the necessarily also $\phi_1=Q_{\ge C} \phi_1$
so that
\begin{align*}
  \|P_0Q_{\le O(1)} \calN_{\alpha\beta} (\phi_1,\phi_2)\|_{\Ltwotx}
&\les \|P_0Q_{\le O(1)} \calN_{\alpha\beta} (\phi_1,\phi_2)\|_{L^1_t L^1_x} \\
&\les \sum_{m\ge C} \|Q_m \phi_1\|_{\Ltwotx} 2^m \|\tilde Q_m \phi_2\|_{\Ltwotx} \\
&\les \|\phi_1\|_{S[k_1]}  \sum_{m\ge C} 2^m 2^{-2m(1-\eps)} 2^{k_2(\frac12-\eps)} \| \phi_2\|_{S[k_2]}   \\
&\les 2^{(\frac12-\eps)k_2} \|\phi_1\|_{S[k_1]} \|\phi_2\|_{S[k_2]}
\end{align*}
and the lemma is proved.
\end{proof}

Next, we deal with the case of outputs with large modulation.

\begin{lemma}
  \label{lem:Nabhighmod2}
For any
$\phi_j$ adapted to $k_j$ with
$k_2\le k_1+O(1)=k$ one has
\[
\sum_{\ell\ge k+C } 2^{-\eps\ell} \|P_kQ_{\ell} \calN_{\alpha\beta} (\phi_1,\phi_2)\|_{\Ltwotx} \les
  2^{(\frac{1}{2}-\eps)k_2}     \|\phi_1\|_{S[k_1]} \|\phi_2\|_{S[k_2]}
\]
\end{lemma}
\begin{proof}
 Set $k=k_1+O(1)=0$ and $k_2\le -C$. Then
\begin{align}
 &\sum_{\ell\ge C } 2^{-\eps\ell} \|P_0 Q_{\ell} \calN_{\alpha\beta} (\phi_1,\phi_2)\|_{\Ltwotx} \nn \\
&\les   \sum_{\ell\ge C } 2^{-\eps\ell} \|P_0 Q_{\ell} \calN_{\alpha\beta} (\tilde Q_\ell \phi_1, Q_{<\ell-C}\phi_2)\|_{\Ltwotx}
\label{eq:highmodsum1}\\
&\quad +  \sum_{\ell\ge C } 2^{-\eps\ell} \|P_0 Q_{\ell} \calN_{\alpha\beta} (Q_{\ge \ell-C} \phi_1, Q_{\ge\ell-C}\phi_2)\|_{\Ltwotx}
\label{eq:highmodsum2}
\end{align}
First, taking $\alpha=0$ and $ \beta=1$,
\begin{align*}
 \eqref{eq:highmodsum1} &\les \sum_{\ell\ge C } 2^{-\eps\ell} \|P_0 Q_{\ell} (R_0 \tilde Q_\ell \phi_1 R_1 Q_{<\ell-C}\phi_2
 -  R_1 \tilde Q_\ell \phi_1 R_0 Q_{<\ell-C}\phi_2   )\|_{\Ltwotx}  \\
&\les  \sum_{\ell\ge C } 2^{-\eps\ell} (2^\ell \|Q_\ell\phi_1\|_{\Ltwotx} \|\phi_2\|_{L^\infty_t L^\infty_x} + \|\tilde Q_\ell\phi_1\|_{\Ltwotx}
\|R_0 Q_{<\ell-C} \phi_2\|_{L^\infty_t L^\infty_x} )\\
&\les \|\phi_1\|_{S[k_1]} \|\phi_2\|_{S[k_2]} 2^{k_2} + \sum_{\ell\ge C }\| \tilde Q_\ell\phi_1\|_{\Ltwotx}
 2^{(\frac12+\eps)\ell} 2^{(\frac12-\eps)k_2} \|\phi_2\|_{S[k_2]} \\
&\les 2^{(\frac12-\eps)k_2} \|\phi_1\|_{S[k_1]}\|\phi_2\|_{S[k_2]}
\end{align*}
To pass to the second to last line we used the estimate
\begin{align*}
 \|R_0 Q_{<\ell-C} \phi_2\|_{L^\infty_t L^\infty_x} &\les \|Q_{\le k_2} \phi_2\|_{\ener} + \sum_{k_2<j<\ell-C}
2^{\frac{3j}{2}} \|Q_j\phi_2\|_{\Ltwotx} \\
&\les 2^{k_2}\|Q_{\le k_2} \phi_2\|_{\ener} + \sum_{k_2<j<\ell-C}
2^{\frac{3j}{2}} 2^{-(1-\eps)j} 2^{(\frac12-\eps)k_2} \|\phi_2\|_{S[k_2]} \\
&\les 2^{k_2} \|Q_{\le k_2} \phi_2\|_{\ener} + \sum_{k_2<j<\ell-C} 2^{(\frac12+\eps) j} 2^{(\frac12-\eps)k_2}
\|\phi_2\|_{S[k_2]} \\
&\les 2^{k_2} \|Q_{\le k_2} \phi_2\|_{\ener} +  2^{(\frac12+\eps)\ell} 2^{(\frac12-\eps)k_2}
\|\phi_2\|_{S[k_2]}
\end{align*}
On \eqref{eq:highmodsum2} one has the bound (again for $\alpha=0$ and $ \beta=1$)
\begin{align*}
 \eqref{eq:highmodsum2} &\les
\sum_{m\ge\ell\ge C } 2^{-\eps\ell} \|P_0 Q_{\ell} (R_0 Q_{m} \phi_1 R_1 \tilde Q_m\phi_2
 -  R_1 Q_{m} \phi_1 R_0 \tilde Q_m\phi_2   )\|_{\Ltwotx}  \\
&\les  \sum_{m\ge\ell\ge C } 2^{(\frac12-\eps)\ell} 2^m \|Q_m\phi_1\|_{\Ltwotx} \|\tilde Q_m \phi_2\|_{\Ltwotx} \\
& \les \sum_{m\ge C } 2^{-(\frac12-\eps)m}   \|\phi_1\|_{S[k_1]}\; 2^{(\frac12-\eps)k_2} \|\phi_2\|_{S[k_2]} \\
&\les 2^{(\frac12-\eps)k_2} \|\phi_1\|_{S[k_1]}\|\phi_2\|_{S[k_2]}
\end{align*}
as claimed.
\end{proof}

\section{Trilinear estimates}\label{sec:trilin}

The purpose of this section is to derive the estimates on the
trilinear nonlinearities which govern the wave map system. In
addition to the bilinear estimates of the previous sections, we
will also heavily use the Strichartz component of the~$S[k]$-norm,
see~\eqref{eq:Sk2}.  As already in~\cite{Krieger}, we will partially
rely on  Tao's trilinear estimate from~\cite{T1} which states that
(relative to our norms in the $S[k]$-spaces)
\begin{equation}
  \label{eq:Taotrilin} \|P_0[ \psi_1 \, R^\beta \psi_2 \,  R_\beta
  \psi_3 ]
  \|_{N[0]} \les 2^{\sigma_1(k_2\wedge k_3-k_1)\wedge  0} 2^{k_1}
  \prod_{i=1}^3 \|\psi_i\|_{S[k_i]}
\end{equation}
for some $\sigma_1>0$. To obtain~\eqref{eq:Taotrilin}
from~\cite{T1}, one observes that $\|\nabla \psi\|_{S[k]}$
(strictly) dominates the $S[k]$-norms of~\cite{T1}, whereas $\|P_0
\cdot\|_{N[0]}$ is dominated by the respective $N[0]$-norm used
in~\cite{T1}. Because of this property, the trilinear bound
from~\cite{T1} can be adapted to this setting provided the correct
scaling is taken into account. Moreover, throughout this section we
define, with $k_{\max}:=k_1\vee k_2\vee k_3$, $k_{\min}:=k_1\wedge
k_2\wedge k_3$, and $k_{\med}$ the median of $k_1,k_2,k_3$,
\[
w(k_1,k_2,k_3) := \left\{ \begin{array}{ll}
2^{-\sigma_0 k_{\max}}\, 2^{\sigma_0 k_{\min}\wedge 0} & \text{if\ \ } k_{\max}\ge C \\
2^{-\sigma_0(k_{\med}-k_{\min}) } & \text{if\ \ } k_1= k_{\max}=O(1)
\\
2^{\sigma_0(k_1+k_2\wedge k_3)} & \text{if\ \ } k_1< k_{\max}=O(1)
\end{array} \right.\]
where $\sigma_0>0$ is some fixed small constant.

\noindent We split our argument into two cases, depending on whether
all inputs are ``hyperbolic'' or not. This distinction is based on
modulation vs.~frequency.

\subsection{Reduction to the hyperbolic case}
The following lemma deals with the case where at least one of the
inputs or the interior null-form have ``elliptic'' character. Recall that
 $I:= \sum_{k\in\Z} P_k Q_{\le k+C}$ and $I^c:=1-I$
(here $C$ is an absolute constant, $C=10$ will suffice). Throughout
this section, we will write $\tilde P_k$ to denote a
projection~$\sum_{k'=k+O(1)} P_{k'}$, and similarly with~$\tilde
Q_k$.

\begin{lemma}
 \label{lem:hyp_redux}
Let $\psi_i$ be Schwarz functions adapted to  $k_i$ for $i=0,1,2$.
Then for any $\alpha=0,1,2$, and $j=1,2$,
\[
\| P_0 \del^\beta A_0[A_1 R_\alpha \psi_1 \Delta^{-1}\del_j \tilde
A_1 \calN_{\beta j}( A_2\psi_2,A_3\psi_3) ] \|_{N[0]} \les
w(k_1,k_2,k_3) \prod_{i=1}^3 \|\psi_i\|_{S[k_i]}
\]
where  $A_i$ and $\tilde A_1$ are either $I$ or $I^c$, with at least
one being~$I^c$. Moreover, we impose the condition that $A_1=\tilde
A_1=I^c$ implies~$\alpha\ne0$.
\end{lemma}
\begin{proof}
{\em Case 1:} $\mathit{0\le k_1\le k_2+O(1)=k_3+O(1) }$.
 We begin with $A_0=I^c$ and $A_1=I$. Then we can drop $I R_\alpha$
 from~$\psi_1$ and estimate\footnote{It is convenient to prove the
 somewhat stronger bound with $A_0=Q_{\ge0}$ here.}
\begin{align}
  \| P_0 Q_{\ge 0} \del^\beta [\psi_1 \Delta^{-1}\del_j  \calN_{\beta j}( \psi_2,\psi_3)
] \|_{N[0]} &\les  \| P_0 Q_{\ge 0} \del^\beta [\psi_1
\Delta^{-1}\del_j Q_{\le k_1+C}  \calN_{\beta j}( \psi_2,\psi_3)
] \|_{N[0]} \label{eq:A}\\
& + \| P_0 Q_{\ge 0} \del^\beta [\psi_1 \Delta^{-1}\del_j Q_{\ge k_1+C} \calN_{\beta j}( \psi_2,\psi_3)
] \|_{N[0]}  \label{eq:B}
\end{align}
By Lemma~\ref{lem:Nablowmod}, placing \eqref{eq:A} into $\dot
X_0^{0,-1-\eps,2}$ implies
\begin{align*}
 \| \psi_1 \Delta^{-1}\del_j Q_{\le k_1+C}  \calN_{\beta j}(
\psi_2,\psi_3) ] \|_{L^2_t L^1_x}  &\les \|\psi_1\|_{L^\infty_t
L^2_x}  2^{-\frac{k_2}{2}}  \|\psi_2\|_{S[k_2]} \|\psi_3\|_{S[k_3]}
\\
&\les 2^{-\frac{k_2}{2}} \|\psi_1\|_{S[k_1]} \|\psi_2\|_{S[k_2]}
\|\psi_3\|_{S[k_3]}
\end{align*}
whereas
\begin{align}
 \eqref{eq:B} &\les \sum_{m\ge k_1+C}  \| P_0 Q_{m} \del^\beta [Q_{\le m-C} \psi_1 \Delta^{-1}\del_j \tilde Q_m \calN_{\beta j}( \psi_2,\psi_3)
] \|_{N[0]} \label{eq:B1} \\
&\quad + \sum_{m\ge k_1+C}   \| P_0 Q_{\ge 0} \del^\beta [Q_{\ge m-C} \psi_1 \Delta^{-1}\del_j Q_{m} \calN_{\beta j}( \psi_2,\psi_3)
] \|_{N[0]} \label{eq:B2}
\end{align}
Lemma~\ref{lem:Nabhighmod} yields the following bound on~\eqref{eq:B1}:
\begin{align*}
 & \sum_{m\ge k_1+C}  \| P_0 Q_{m} \del^\beta [Q_{\le m-C} \psi_1 \Delta^{-1}\del_j \tilde Q_m \calN_{\beta j}(   \psi_2,\psi_3)
] \|_{N[0]} \nn \\
&\les \sum_{m\ge k_1+C}  \| P_0 Q_{m} \del^\beta [Q_{\le m-C} \psi_1
\Delta^{-1}\del_j \tilde Q_m \calN_{\beta j}(   \psi_2,\psi_3)
] \|_{\dot X^{0,-1-\eps,2}_0}  \nn    \\
&\les  \sum_{m\ge k_1+C} 2^{-m\eps} \| \psi_1\|_{\ener} 2^{-k_1} \| \tilde P_{k_1} \tilde Q_m \calN_{\beta j}( \psi_2,\psi_3) \|_{L^2_t L^2_x}\nn \\
&\les  2^{-\eps k_2}2^{-\frac{k_1}{2}} \la k_2-k_1\ra^2
\prod_{i=1}^3 \|\psi_i\|_{S[k_i]}
\end{align*}
The bound on~\eqref{eq:B2} proceeds similarly:
\begin{align*}
 \eqref{eq:B2} &\les  \sum_{m\ge k_1+C}  \| P_0 Q_{\ge0 } \del^\beta [Q_{> m-C} \psi_1 \Delta^{-1}\del_j \tilde Q_m \calN_{\beta j}(\psi_2,\psi_3)
] \|_{N[0]} \nn \\
&  \les \sum_{m\ge k_1+C} \sum_{0\le\ell\le m+C}  \| P_0 Q_{\ell} \del^\beta [Q_{> m-C} \psi_1 \Delta^{-1}\del_j \tilde Q_m \calN_{\beta j}(\psi_2,\psi_3)
] \|_{N[0]} \nn \\
& \quad + \sum_{m\ge k_1+C} \sum_{\ell\ge m+C}  \| P_0 Q_{\ell} \del^\beta [\tilde Q_\ell \psi_1 \Delta^{-1}\del_j \tilde Q_m \calN_{\beta j}(\psi_2,\psi_3)
] \|_{N[0]} \nn  \\
&  \les \sum_{m\ge k_1+C} \sum_{0\le\ell\le m+C} 2^{(\frac12-\eps)\ell} \| Q_{> m-C} \psi_1\|_{\Ltwotx}
\| \tilde P_{k_1}\Delta^{-1}\del_j \tilde Q_m \calN_{\beta j}(\psi_2,\psi_3)
] \|_{\Ltwotx} \nn \\
& \quad + \sum_{m\ge k_1+C} \sum_{\ell\ge m+C}  2^{-\eps\ell} \| \tilde Q_\ell \psi_1\|_{\Ltwotx}
\|  \tilde P_{k_1}\Delta^{-1}\del_j \tilde Q_m \calN_{\beta j}(\psi_2,\psi_3)
] \|_{\ener}
\end{align*}
In the second to last line we applied Bernstein's inequality in the
time variable to switch from $L^2_t$ to~$L^1_t$.  We now replace the
$L^\infty_t$ on the right-hand side of the last line by an $L^2_t$
at the expense of a factor of~$2^{\frac{m}{2}}$. Together with
Lemma~\ref{lem:Nabhighmod}  this yields
\begin{align*}
 \eqref{eq:B2} &\les  \sum_{m\ge k_1+C} \sum_{0\le\ell\le m+C} 2^{-k_1+(\frac12-\eps)\ell}
 \| Q_{> m-C} \psi_1\|_{\Ltwotx} \| \tilde P_{k_1} \tilde Q_m \calN_{\beta j}(\psi_2,\psi_3)
] \|_{\Ltwotx} \nn \\
& \quad + \sum_{m\ge k_1+C} \sum_{\ell\ge m+C}  2^{-k_1-\eps\ell}
2^{\frac{m}{2}} \| \tilde Q_\ell \psi_1\|_{\Ltwotx} \|\tilde P_{k_1}
\tilde Q_m \calN_{\beta j}(\psi_2,\psi_3)
] \|_{\Ltwotx}  \\
&\les  \sum_{m\ge k_1+C}  2^{-(\frac12+\eps)k_1} 2^{-\frac12m}
\| \psi_1\|_{S[k_1]}
\| \tilde P_{k_1} \tilde Q_m \calN_{\beta j}(\psi_2,\psi_3)
] \|_{\Ltwotx} \nn \\
& \quad + \sum_{m\ge k_1+C} \sum_{\ell\ge m+C}  2^{-k_1-\eps\ell} 2^{\frac{m}{2}} 2^{-(1-\eps)\ell} 2^{(\frac12-\eps)k_1} \| \psi_1\|_{S[k_1]}
 \| \tilde P_{k_1} \tilde Q_m \calN_{\beta j}(\psi_2,\psi_3)
] \|_{\Ltwotx}  \\
&\les 2^{-\frac{k_1}{2}} 2^{-\eps k_2} \la k_2-k_1\ra^2
\prod_{i=1}^3 \|\psi_i\|_{S[k_i]}
\end{align*}
Next, we consider the case where both $A_0=I^c$ and $A_1=I^c$. If
$\alpha\ne0$, then $A_1=I$ and one can drop~$R_\alpha$ altogether so
that  the previous analysis applies. Otherwise, if~$\alpha=0$, then
by assumption $\tilde A_1=I$ and
\begin{align}
   & \| P_0 Q_{\ge 0} \del^\beta [Q_{\ge k_1+C} R_\alpha \psi_1 \Delta^{-1}\del_j
\tilde P_{k_1} Q_{\le k_1+C}  \calN_{\beta j}( A_2\psi_2,A_3\psi_3) ]
\|_{N[0]} \nn \\
&\le \sum_{m\ge k_1+10C} \| P_0 Q_{m} \del^\beta [\tilde Q_m
R_\alpha \psi_1 \Delta^{-1}\del_j   \tilde P_{k_1} Q_{\le k_1+C}
\calN_{\beta j}( A_2\psi_2,A_3\psi_3) ]
\|_{N[0]} \label{eq:sch6}\\
& +  \| P_0 Q_{0\le\cdot \le k_1+10C} \del^\beta [Q_{0\le \cdot \le
k_1+10C} R_\alpha \psi_1 \Delta^{-1}\del_j   \tilde P_{k_1}
Q_{\le k_1+C} \calN_{\beta j}( A_2\psi_2,A_3\psi_3) ] \|_{N[0]}
\label{eq:sch7}
\end{align}
By Lemma~\ref{lem:Nablowmod}, \eqref{eq:sch7} is bounded by
\begin{align*}
 & \| P_0 Q_{0\le\cdot\le k_1+10C} \del^\beta [Q_{0\le \cdot\le k_1+10C}
R_\alpha \psi_1 \Delta^{-1}\del_j   \tilde P_{k_1} Q_{\le k_1+C}
\calN_{\beta j}( A_2\psi_2,A_3\psi_3) ] \|_{\dot X_0^{0,-1-\eps,2}} \\
&\les \| Q_{0\le \cdot\le k_1+10C} R_\alpha \psi_1\;
\Delta^{-1}\del_j   \tilde P_{k_1} Q_{\le k_1+C} \calN_{\beta j}(
A_2\psi_2,A_3\psi_3)
\|_{L^2_t L^1_x} \\
&\les  \| Q_{0\le \cdot\le k_1+10C} R_\alpha
\psi_1\|_{\ener}\,
 \| \Delta^{-1}\del_j   \tilde P_{k_1} Q_{\le
k_1+C} \calN_{\beta j}( A_2\psi_2,A_3\psi_3) \|_{L^2_t L^2_x} \\
&\les 2^{-\frac{k_2}{2}} \prod_{i=1}^3 \|\psi_i\|_{S[k_i]}
\end{align*}
 On the other hand,~\eqref{eq:sch6} is
estimated as follows:
\begin{align*}
 & \sum_{m\ge k_1+10C} \| P_0 Q_{m} \del^\beta [\tilde Q_m
R_\alpha \psi_1 \Delta^{-1}\del_j   \tilde P_{k_1} Q_{\le k_1+C}
\calN_{\beta j}( A_2\psi_2,A_3\psi_3) ] \|_{\dot
X_0^{0,-1-\eps,2}}\\
& \les \sum_{m\ge k_1} 2^{-m\eps}\| \tilde Q_m R_\alpha \psi_1
\|_{\Ltwotx} \| \Delta^{-1}\del_j   \tilde P_{k_1} Q_{\le k_1+C}
\calN_{\beta j}( A_2\psi_2,A_3\psi_3) ] \|_{L^\infty_t L^2_x}\\
& \les  2^{-(\frac12+\eps)k_1}  \| \psi_1 \|_{S[k_1]}\, 2^{-\frac{k_1}{2}} \|
 \tilde P_{k_1} Q_{\le k_1+C} \calN_{\beta j}(
A_2\psi_2,A_3\psi_3) ] \|_{L^2_t L^2_x} \\
&\les 2^{-\eps k_1 -\frac{k_2}{2}} \prod_{i=1}^3 \|\psi_i\|_{S[k_i]}
\end{align*}
where we applied Bernstein's inequality relative to~$t$ as well as
Lemma~\ref{lem:Nablowmod} to pass to the last line.

\noindent Now suppose $A_0=I$ (in fact, $A_0=Q_{\le0}$), but at
least one of $A_1$ or~$\tilde A_1$ equals $I^c$. But then the
modulations of $\psi_1$ and $\calN_{\beta j}$ essentially agree,
whence $\alpha\ne0$ and
\begin{align*}
 & \sum_{m\ge k_1+C} \| P_0 Q_{\le 0} \del^\beta [ Q_m R_\alpha \psi_1  \Delta^{-1}\del_j \tilde Q_m  \calN_{\beta j}( \psi_2,\psi_3)
] \|_{N[0]} \\
& \les \sum_{m\ge k_1+C} \| P_0 Q_{\le 0} \del^\beta [ Q_m R_\alpha \psi_1
\Delta^{-1}\del_j \tilde Q_m  \calN_{\beta j}( \psi_2,\psi_3)
] \|_{L^1_{tx}} \\
& \les  \sum_{m\ge k_1+C} \|  Q_m \psi_1 \|_{\Ltwotx} 2^{-k_1} \|
\tilde P_{k_1}  \tilde Q_m  \calN_{\beta j}( \psi_2,\psi_3)
 \|_{\Ltwotx} \\
&\les \sum_{m\ge k_1+C} 2^{(\frac12-\eps)k_1} 2^{-m(1-2\eps)} \|
\psi_1 \|_{S[k_1]} \, 2^{-k_1} 2^{-m\eps}\| \tilde P_{k_1}  \tilde Q_m
\calN_{\beta j}( \psi_2,\psi_3)
 \|_{\Ltwotx} \\
 &\les 2^{-\frac{k_1}{2}-\eps k_2}  \la k_2-k_1\ra^2
\prod_{i=1}^3 \|\psi_i\|_{S[k_i]}
\end{align*}
The final estimate here uses Lemma~\ref{lem:Nabhighmod2}.  The last
case which we need to consider is $A_0=A_1=\tilde A_1 =I$ and either
one of $A_2, A_3$ equal to~$I^c$. But then necessarily $A_2=A_3=I^c$
whence
\begin{align*}
 & \| P_0 \del^\beta I[IR_\alpha\psi_1 \Delta^{-1}\del_j I \calN_{\beta j}( Q_{\ge k_2+C} \psi_2, Q_{\ge k_2+C}\psi_3)
] \|_{N[0]}  \\
 & \les\| IR_\alpha\psi_1 \Delta^{-1}\del_j I \calN_{\beta j}( Q_{\ge k_2+C} \psi_2, Q_{\ge k_2+C}\psi_3)
 \|_{L^1_t L^1_x}  \\
&\les \|\psi_1\|_{\ener} \sum_{m\ge k_2+C} \|\tilde P_{k_1}Q_{\le
k_1+C} \calN_{\beta j}( Q_{m} \psi_2, \tilde Q_m\psi_3) ] \|_{L^1_t
L^1_x} \\
&\les \|\psi_1\|_{\ener} \sum_{m\ge k_2+C} 2^{m-k_2} 2^{-2m(1-\eps)}
2^{(1-2\eps)k_2} \|\psi_2\|_{S[k_2]} \|\psi_3\|_{S[k_3]} \les
2^{-k_2} \prod_{i=1}^3 \|\psi_i\|_{S[k_i]}
\end{align*}
which concludes Case~1.

\medskip
\noindent {\em Case 2:} $\mathit{0\le k_1= k_2+O(1), k_3\le k_2-C }$. We
again begin with $A_0=I^c$, $A_1=I$ and the
representation~\eqref{eq:A} and~\eqref{eq:B} (dropping $IR_\alpha$ from $\psi_1$ as before). By
Lemma~\ref{lem:Nablowmod2}, \eqref{eq:A} is bounded by
\begin{align*}
& \| \psi_1 \Delta^{-1}\del_j Q_{\le k_1+C}  \calN_{\beta j}(
\psi_2,\psi_3) ] \|_{L^2_t L^1_x}  \les \|\psi_1\|_{L^\infty_t
L^2_x}  2^{(\frac12-\eps)k_3} 2^{-(1-\eps)k_1}   \|\psi_2\|_{S[k_2]}
\|\psi_3\|_{S[k_3]}
\end{align*}
whereas
\begin{align}
 \eqref{eq:B} &\les \sum_{m\ge k_1+C}  \| P_0 Q_{m} \del^\beta [Q_{\le m-C} \psi_1 \Delta^{-1}\del_j \tilde Q_m \calN_{\beta j}( \psi_2,\psi_3)
] \|_{N[0]} \label{eq:B12} \\
&\quad + \sum_{m\ge k_1+C}   \| P_0 Q_{\ge 0} \del^\beta [Q_{\ge m-C} \psi_1 \Delta^{-1}\del_j Q_{m} \calN_{\beta j}( \psi_2,\psi_3)
] \|_{N[0]} \label{eq:B22}
\end{align}
Lemma~\ref{lem:Nabhighmod2} yields the following bound on~\eqref{eq:B12}:
\begin{align*}
 & \sum_{k_1+C\le m}  \| P_0 Q_{m} \del^\beta [Q_{\le m-C} \psi_1 \Delta^{-1}\del_j \tilde Q_m \calN_{\beta j}(   \psi_2,\psi_3)
] \|_{N[0]} \nn \\
&\les \sum_{k_1+C\le m}  \| P_0 Q_{m} \del^\beta [Q_{\le m-C} \psi_1 \Delta^{-1}\del_j \tilde Q_m \calN_{\beta j}(   \psi_2,\psi_3)
] \|_{\dot X^{0,-1-\eps,2}_0}  \nn    \\
&\les  \sum_{k_1+C\le m} 2^{-m\eps} \| \psi_1\|_{\ener} 2^{-k_1} \| \tilde P_{k_1} \tilde Q_m \calN_{\beta j}( \psi_2,\psi_3) \|_{L^2_t L^2_x}\nn \\
&\les  2^{(\frac12-\eps) k_3}2^{-k_1}   \prod_{i=1}^3 \|\psi\|_{S[k_i]}
\end{align*}
The bound on~\eqref{eq:B22} proceeds similarly:
\begin{align}
 \eqref{eq:B22} &\les  \sum_{m\ge k_1+C}  \| P_0 Q_{\ge0 } \del^\beta [Q_{> m-C} \psi_1 \Delta^{-1}\del_j \tilde Q_m \calN_{\beta j}(\psi_2,\psi_3)
] \|_{N[0]} \nn \\
&  \les \sum_{m\ge k_1+C} \sum_{0\le\ell\le m+C}  \| P_0 Q_{\ell} \del^\beta [Q_{> m-C} \psi_1 \Delta^{-1}\del_j \tilde Q_m \calN_{\beta j}(\psi_2,\psi_3)
] \|_{\dot X_0^{0,-1-\eps,2}} \nn \\
& \quad + \sum_{m\ge k_1+C} \sum_{\ell\ge m+C}  \| P_0 Q_{\ell} \del^\beta [\tilde Q_\ell \psi_1 \Delta^{-1}\del_j \tilde Q_m \calN_{\beta j}(\psi_2,\psi_3)
] \|_{\dot X_0^{0,-1-\eps,2}} \nn  \\
&  \les \sum_{m\ge k_1+C} \sum_{0\le\ell\le m+C}
2^{(\frac12-\eps)\ell} \| Q_{> m-C} \psi_1\|_{\Ltwotx} \|
\tilde P_{k_1}\Delta^{-1}\del_j \tilde Q_m \calN_{\beta
j}(\psi_2,\psi_3)
] \|_{\Ltwotx} \label{eq:sch3} \\
& \quad + \sum_{m\ge k_1+C} \sum_{\ell\ge m+C}  2^{-\eps\ell} \| \tilde Q_\ell \psi_1\|_{\Ltwotx}
\|  \tilde P_{k_1}\Delta^{-1}\del_j \tilde Q_m \calN_{\beta j}(\psi_2,\psi_3)
] \|_{\ener}
\end{align}
To pass to~\eqref{eq:sch3} we  used Bernstein's inequality to switch
from~$L^2_t$ to~$L^1_t$, which costs~$2^{\frac{\ell}{2}}$. We now
replace the $L^\infty_t$ on the right-hand side of the last line by
an $L^2_t$ at the expense of a factor of~$2^{\frac{m}{2}}$. In view
of Lemma~\ref{lem:Nabhighmod2} one concludes that
\begin{align*}
 \eqref{eq:B22} &\les  \sum_{m\ge k_1+C} \sum_{0\le\ell\le m+C} 2^{-k_1+(\frac12-\eps)\ell}  \| Q_{> m-C} \psi_1\|_{\Ltwotx}
  \| \tilde P_{k_1} \tilde Q_m \calN_{\beta j}(\psi_2,\psi_3)
] \|_{\Ltwotx} \nn \\
& \quad + \sum_{m\ge k_1+C} \sum_{\ell\ge m+C}  2^{-k_1-\eps\ell}
2^{\frac{m}{2}} \| \tilde Q_\ell \psi_1\|_{\Ltwotx} \|\tilde P_{k_1}
\tilde Q_m \calN_{\beta j}(\psi_2,\psi_3)
] \|_{\Ltwotx}  \\
&\les  \sum_{m\ge k_1+C}  2^{-(\frac12+\eps)k_1 } 2^{-\frac12 m}
\| \psi_1\|_{S[k_1]}
\| \tilde P_{k_1} \tilde Q_m \calN_{\beta j}(\psi_2,\psi_3)
] \|_{\Ltwotx} \nn \\
& \quad + \sum_{m\ge k_1+C} \sum_{\ell\ge m+C}  2^{-k_1-\eps\ell} 2^{\frac{m}{2}} 2^{-(1-\eps)\ell} 2^{(\frac12-\eps)k_1} \| \psi_1\|_{S[k_1]}
 \| \tilde P_{k_1} \tilde Q_m \calN_{\beta j}(\psi_2,\psi_3)
] \|_{\Ltwotx}  \\
&\les 2^{-k_1} 2^{(\frac12-\eps)k_3}
\prod_{i=1}^3 \|\psi_i\|_{S[k_i]}
\end{align*}
Next, we consider the case where both $A_0=I^c$ and $A_1=I^c$. If
$\alpha\ne0$, then $\tilde A_1=I$ and one can drop~$R_\alpha$ altogether so
that  the previous analysis applies. Otherwise, if~$\alpha=0$, then
by assumption $\tilde A_1=I$ and as in Case~1 one
obtains~\eqref{eq:sch6} and~\eqref{eq:sch7}. By
Lemma~\ref{lem:Nablowmod2}, \eqref{eq:sch7} is bounded by
\begin{align*}
 & \| P_0 Q_{0\le\cdot\le k_1+10C} \del^\beta [Q_{0\le \cdot\le k_1+10C}
R_\alpha \psi_1 \Delta^{-1}\del_j   \tilde P_{k_1} Q_{\le k_1+C}
\calN_{\beta j}( A_2\psi_2,A_3\psi_3) ] \|_{\dot X_0^{0,-1-\eps,2}} \\
&\les \| Q_{0\le \cdot\le k_1+10C} R_\alpha \psi_1\;
\Delta^{-1}\del_j   \tilde P_{k_1} Q_{\le k_1+C} \calN_{\beta j}(
A_2\psi_2,A_3\psi_3)
\|_{L^2_t L^1_x} \\
&\les 2^{\frac{k_1}{2}} \| Q_{0\le \cdot\le k_1+10C} R_\alpha
\psi_1\|_{\Ltwotx}\,
 \| \Delta^{-1}\del_j   \tilde P_{k_1} Q_{\le
k_1+C} \calN_{\beta j}( A_2\psi_2,A_3\psi_3) \|_{L^2_t L^2_x} \\
&\les 2^{-(1-\eps)k_1 } 2^{(\frac12-\eps)k_3} \prod_{i=1}^3
\|\psi_i\|_{S[k_i]}
\end{align*}
On the other hand,~\eqref{eq:sch6} is estimated as follows:
\begin{align*}
 & \sum_{m\ge k_1+10C} \| P_0 Q_{m} \del^\beta [\tilde Q_m
R_\alpha \psi_1 \Delta^{-1}\del_j   \tilde P_{k_1} Q_{\le k_1+C}
\calN_{\beta j}( A_2\psi_2,A_3\psi_3) ] \|_{\dot
X_0^{0,-1-\eps,2}}\\
& \les \sum_{m\ge k_1} 2^{-m\eps}\| \tilde Q_m \nabla_{t,x}|\nabla|^{-1} \psi_1
\|_{\Ltwotx}\; 2^{-k_1} \|    \tilde P_{k_1} Q_{\le k_1+C}
\calN_{\beta j}( A_2\psi_2,A_3\psi_3) ] \|_{L^\infty_t L^2_x}\\
& \les  2^{-(\frac12+\eps)k_1}  \| \psi_1 \|_{S[k_1]}\, 2^{-\frac{k_1}{2}} \|
 \tilde P_{k_1} Q_{\le k_1+C} \calN_{\beta j}(
A_2\psi_2,A_3\psi_3) ] \|_{L^2_t L^2_x} \\
&\les 2^{-k_1 +(\frac12-\eps)k_3} \prod_{i=1}^3 \|\psi_i\|_{S[k_i]}
\end{align*}
where we applied Bernstein's inequality relative to~$t$ as well as
Lemma~\ref{lem:Nablowmod2}.

\noindent We now turn to the case where $A_0=I$, but at least one of
$A_1$ or~$\tilde A_1$ equals $I^c$. But then the modulations of
$\psi_1$ and $\calN_{\beta j}$ essentially agree whence
$\alpha\ne0$. Bounding $N[0]$ by~$L^1_t L^2_x$ and invoking
Lemma~\ref{lem:Nabhighmod2} yields
\begin{align*}
 & \sum_{m\ge k_1+C} \| P_0 Q_{\le 0} \del^\beta [ Q_m \psi_1  \Delta^{-1}\del_j \tilde Q_m  \calN_{\beta j}( \psi_2,\psi_3)
] \|_{N[0]} \\
&\les  \sum_{m\ge k_1+C} \|  Q_m \psi_1 \|_{\Ltwotx}\, 2^{-k_1} \|
\tilde P_{k_1}  \tilde Q_m  \calN_{\beta j}( \psi_2,\psi_3)
] \|_{\Ltwotx} \\
&\les 2^{-(\frac32-\eps)k_1+(\frac12-\eps) k_3}  \prod_{i=1}^3 \|\psi_i\|_{S[k_i]}
\les 2^{-k_1} 2^{(\frac12-\eps)k_3}
\prod_{i=1}^3 \|\psi_i\|_{S[k_i]}
\end{align*}
The last case which we need to consider is $A_0=A_1=\tilde A_1 =I$
and either one of $A_2, A_3$ equal to~$I^c$. We begin with
$A_2=I^c$. But then necessarily $A_2=A_3=I^c$ whence
\begin{align*}
 & \| P_0 \del^\beta I[I\psi_1 \Delta^{-1}\del_j I \calN_{\beta j}( Q_{\ge k_2+C} \psi_2, Q_{\ge k_2+C}\psi_3)
] \|_{N[0]}  \\
 & \les\| I\psi_1 \Delta^{-1}\del_j I \calN_{\beta j}( Q_{\ge k_2+C} \psi_2, Q_{\ge k_2+C}\psi_3)
 \|_{L^1_t L^1_x}  \\
&\les \|\psi_1\|_{\ener} 2^{-k_1} \sum_{m\ge k_2+C} \|\tilde P_{k_1}Q_{\le
k_1+C} \calN_{\beta j}( Q_{m} \psi_2, \tilde Q_m\psi_3) ] \|_{L^1_t
L^2_x} \\
&\les \|\psi_1\|_{\ener} \sum_{m\ge k_2+C} 2^{m-k_2} 2^{-2(1-\eps)m} 2^{(\frac12-\eps)k_2} 2^{(\frac12-\eps)k_3}   \|\psi_2\|_{S[k_2]}
\|\psi_3\|_{S[k_3]}\\
&\les  2^{-(\frac32-\eps)k_1+(\frac12-\eps) k_3}  \prod_{i=1}^3 \|\psi_i\|_{S[k_i]}
\les 2^{-k_1} 2^{(\frac12-\eps)k_3}
\prod_{i=1}^3 \|\psi_i\|_{S[k_i]}
\end{align*}
It remains to consider the case $A_2=I$ and $A_3=I^c$. We begin by
reducing the modulation of the entire output. Indeed, by
Lemma~\ref{lem:Nablowmod2},
\begin{align*}
 & \| P_0 \del^\beta Q_{(1-3\eps)k_3\le\cdot\le C} [Q_{\le k_1+C} \psi_1 \Delta^{-1}\del_j Q_{\le k_1+C}  \calN_{\beta j}( I \psi_2,I^c \psi_3)
] \|_{N[0]} \nn \\
&\les 2^{-\frac12(1-3\eps)k_3} \| \psi_1\|_{L^\infty_t L^2_x} \,
2^{-k_1}  \| I\calN_{\beta j}( I \psi_2,I^c \psi_3) \|_{\Ltwotx}\\
&\les 2^{-k_1}  2^{-\frac12(1-3\eps)k_3} \| \psi_1\|_{L^\infty_t
L^2_x} \, 2^{(\frac12-\eps)k_3} 2^{\eps k_2} \|\psi_2\|_{S[k_2]}
\|\psi_3\|_{S[k_3]}\\
&\les 2^{-(1-\eps)k_1} 2^{\frac{\eps}{2}k_3} \prod_{i=1}^3
\|\psi_i\|_{S[k_i]}
\end{align*}
Next, we reduce the modulation of $\psi_1$:
\begin{align*}
 & \| P_0 \del^\beta Q_{\le(1-3\eps)k_3} [Q_{\ge (1-3\eps)k_3-k_1} \psi_1 \Delta^{-1}\del_j Q_{\le k_1+C}  \calN_{\beta j}( I \psi_2,I^c \psi_3)
] \|_{N[0]} \nn \\
 &\les \| P_0 \del^\beta Q_{\le(1-3\eps)k_3} [Q_{\ge (1-3\eps)k_3-k_1} \psi_1 \Delta^{-1}\del_j Q_{\le k_1+C}  \calN_{\beta j}( I \psi_2,I^c \psi_3)
] \|_{L^1_{t}L^1_x} \nn \\
&\les  \|Q_{\ge (1-3\eps)k_3-k_1} \psi_1\|_{L^2_t L^2_x} \,
2^{-k_1}  \| I\calN_{\beta j}( I \psi_2,I^c \psi_3) \|_{\Ltwotx}\\
&\les 2^{-\frac{k_1}{2}}  2^{-\frac12(1-3\eps)k_3} \|
\psi_1\|_{S[k_1]} \, 2^{(\frac12-\eps)k_3} 2^{\eps k_2}
\|\psi_2\|_{S[k_2]}
\|\psi_3\|_{S[k_3]}\\
&\les 2^{-(\frac12-\eps)k_1} 2^{\frac{\eps}{2}k_3} \prod_{i=1}^3
\|\psi_i\|_{S[k_i]}
\end{align*}
Finally, we reduce the modulation of the interior null-form using
Lemma~\ref{lem:core}:
\begin{align*}
& \| P_0 \del^\beta Q_{\le (1-3\eps)k_3} [Q_{\le (1-3\eps)k_3-k_1}
\psi_1 \Delta^{-1}\del_j Q_{k_3\le \cdot\le k_1+C} \calN_{\beta j}(
I \psi_2,
I^c\psi_3) ] \|_{N[0]} \\
&\les  \|\psi_1\|_{S[k_1]} \sum_{k_3\le\ell\le k_1+C }
2^{-\frac{\ell}{4}} 2^{-k_1} \la k_1\ra \| P_{k_1}
 Q_{\ell} \calN_{\beta j}( I \psi_2,
I^c\psi_3)  \|_{\Ltwotx}\\
&\les 2^{-(1-2\eps)k_1} 2^{(\frac14-\eps)k_3}  \prod_{i=1}^3
\|\psi_i\|_{S[k_i]}
\end{align*}
which is again admissible.  After these preparations, we are faced
with the following decomposition:
\begin{align*}
 &  P_0 \del^\beta Q_{\le (1-3\eps)k_3} [Q_{\le (1-3\eps)k_3-k_1}
 \psi_1 \Delta^{-1}\del_j Q_{\le k_3}  \calN_{\beta j}( I \psi_2, I^c\psi_3)
]  \nn \\
&=  P_0 \del^\beta Q_{\le (1-3\eps)k_3} [Q_{\le (1-3\eps)k_3-k_1 }
\psi_1 \Delta^{-1}\del_j Q_{\le k_3} \calN_{\beta j}( Q_{k_3\le
\cdot\le k_2+C} \psi_2, Q_{k_3+C\le\cdot\le k_2+C}\psi_3) ]\\
& = \sum_{\kappa,\kappa'\in\calC_{\ell}} P_{0,\kappa} \del^\beta
Q_{\le (1-3\eps)k_3} [P_{k_1,\kappa'} Q_{\le (1-3\eps)k_3-k_1 }
\psi_1 \Delta^{-1}\del_j Q_{\le k_3} \calN_{\beta j}( Q_{k_3\le
\cdot\le k_2+C} \psi_2, Q_{k_3+C\le\cdot\le k_2+C}\psi_3) ]
\end{align*}
where $\ell=\frac12 (1-3\eps)k_3$ and $\dist(\kappa,\kappa')\les
2^\ell$. Placing the entire expression in~$L^1_t L^2_x$ and using
Bernstein's inequality results in the following estimate: with $J:=Q_{k_3\le \cdot\le k_2+C} $,
\begin{align*}
 &  \|P_0 \del^\beta Q_{\le (1-3\eps)k_3} [Q_{\le (1-3\eps)k_3-k_1}
 \psi_1 \Delta^{-1}\del_j Q_{\le k_3}  \calN_{\beta j}( I \psi_2, I^c\psi_3)
] \|_{L^1_t L^2_x}  \nn \\
& \le \Big\|  \Big( \sum_{\kappa,\kappa'\in\calC_{\ell}}
\|P_{0,\kappa}  [P_{k_1,\kappa'} Q_{\le (1-3\eps)k_3-k_1 } \psi_1
\Delta^{-1}\del_j Q_{\le k_3} \calN_{\beta j}( J \psi_2, J\psi_3) ]\|_{L^2_x}^2
\Big)^{\frac12}
\Big\|_{L^1_t} \\
& \le 2^{\frac{\ell}{2}} \Big\|  \Big(
\sum_{\kappa,\kappa'\in\calC_{\ell}} \|P_{0,\kappa} [P_{k_1,\kappa'}
Q_{\le (1-3\eps)k_3-k_1 } \psi_1 \Delta^{-1}\del_j Q_{\le k_3}
\calN_{\beta j}( J \psi_2,
J\psi_3) ]\|_{L^1_x}^2 \Big)^{\frac12}
\Big\|_{L^1_t}\\
& \le 2^{\frac{\ell}{2}} \Big\|  \Big( \sum_{\kappa'\in\calC_{\ell}}
\|P_{k_1,\kappa'} Q_{\le (1-3\eps)k_3-k_1 } \psi_1\|_{L^2_x}^2 \|
\Delta^{-1}\del_j Q_{\le k_3} \calN_{\beta j}( J \psi_2, J\psi_3) ]\|_{L^2_x}^2
\Big)^{\frac12}
\Big\|_{L^1_t}\\
& \le 2^{\frac{\ell}{2}}  \| Q_{\le (1-3\eps)k_3-k_1 }
\psi_1\|_{L^\infty_t L^2_x} \| \Delta^{-1}\del_j Q_{\le k_3}
\calN_{\beta j}( J \psi_2,
J\psi_3) ]\|_{L^1_t L^2_x}\\
&\les 2^{\frac14(1-3\eps)k_3} \|\psi_1\|_{S[k_1]} \; 2^{-k_1}
\|\nabla_{t,x}|\nabla|^{-1}  J
\psi_2\|_{\Ltwotx} \| \nabla_{t,x}|\nabla|^{-1} J\psi_3\|_{L^2_t L^\infty_x} \\
&\les 2^{\frac14(1-3\eps)k_3} \|\psi_1\|_{S[k_1]} \; 2^{-k_1}
2^{-\frac{k_3}{2}} \| \psi_2\|_{S[k_2]} \, 2^{(\frac12-\eps)k_3}
2^{\eps k_2}  \| \psi_3\|_{S[k_3]}
\end{align*}
which is again admissible for small $\eps>0$.

\medskip
\noindent {\em Case 3:} $\mathit{0\le k_1= k_3+O(1), k_2\le k_3-C .}$ This
case is symmetric to the previous one.

\medskip
\noindent {\em Case 4:} $\mathit{O(1)\le k_2= k_3+O(1), k_1\le -C }$. This case proceeds
similarly to Case~1. We again begin with $A_0=I^c$ and $A_1=I$. Then we can drop $I R_\alpha$
 from~$\psi_1$ and estimate
\begin{align}
  \| P_0 Q_{\ge 0} \del^\beta [\psi_1 \Delta^{-1}\del_j  \calN_{\beta j}( \psi_2,\psi_3)
] \|_{N[0]} &\les  \| P_0 Q_{\ge 0} \del^\beta [\psi_1
\Delta^{-1}\del_j \wt P_0 Q_{<C}  \calN_{\beta j}( \psi_2,\psi_3)
] \|_{N[0]} \label{eq:A'}\\
& + \| P_0 Q_{\ge 0} \del^\beta [\psi_1 \Delta^{-1}\del_j \wt P_0 Q_{\ge C} \calN_{\beta j}( \psi_2,\psi_3)
] \|_{N[0]}  \label{eq:B'}
\end{align}
where we write $\wt P_0=P_{[-C,C]}$ for simplicity. By
Lemma~\ref{lem:Nablowmod}, placing \eqref{eq:A'} into $\dot
X_0^{0,-1-\eps,2}$ implies
\begin{align*}
 \| \psi_1 \Delta^{-1}\del_j \wt P_0 Q_{< C}   \calN_{\beta j}(
\psi_2,\psi_3) ] \|_{L^2_t L^2_x}  &\les \|\psi_1\|_{L^\infty_t
L^\infty_x} \; 2^{-\frac{k_2}{2}}  \|\psi_2\|_{S[k_2]} \|\psi_3\|_{S[k_3]}
\\
&\les 2^{k_1} 2^{-\frac{k_2}{2}} \prod_{i=1}^3 \|\psi_i\|_{S[k_i]}
\end{align*}
whereas
\begin{align}
 \eqref{eq:B'} &\les \sum_{m\ge C}  \| P_0 Q_{m} \del^\beta [Q_{\le m-C} \psi_1 \Delta^{-1}\del_j \wt P_0 \tilde Q_m \calN_{\beta j}( \psi_2,\psi_3)
] \|_{N[0]} \label{eq:B'1} \\
&\quad + \sum_{m\ge C}   \| P_0 Q_{\ge 0} \del^\beta [Q_{\ge m-C} \psi_1 \Delta^{-1}\del_j \wt P_0 Q_{m} \calN_{\beta j}( \psi_2,\psi_3)
] \|_{N[0]} \label{eq:B'2}
\end{align}
Lemma~\ref{lem:Nabhighmod} yields the following bound on~\eqref{eq:B'1}:
\begin{align*}
 & \sum_{m\ge C}  \| P_0 Q_{m} \del^\beta [Q_{\le m-C} \psi_1 \Delta^{-1}\del_j \wt P_0 \tilde Q_m \calN_{\beta j}(   \psi_2,\psi_3)
] \|_{N[0]} \nn \\
&\les \sum_{m\ge C}  \| P_0 Q_{m} \del^\beta [Q_{\le m-C} \psi_1
\Delta^{-1}\del_j \tilde Q_m \calN_{\beta j}(   \psi_2,\psi_3)
] \|_{\dot X^{0,-1-\eps,2}_0}  \nn    \\
&\les  \sum_{m\ge  C} 2^{-m\eps} \| \psi_1\|_{\Linf}  \|   \wt P_0 \tilde Q_m \calN_{\beta j}( \psi_2,\psi_3) \|_{L^2_t L^2_x}\nn \\
&\les 2^{k_1} 2^{-\eps k_2} \la k_2-k_1\ra^2
\prod_{i=1}^3 \|\psi_i\|_{S[k_i]}
\end{align*}
which is admissible.
The bound on~\eqref{eq:B'2} proceeds similarly:
\begin{align*}
 \eqref{eq:B'2} &\les  \sum_{m\ge C}  \| P_0 Q_{\ge0 } \del^\beta [Q_{> m-C} \psi_1 \Delta^{-1}\del_j \tilde Q_m \calN_{\beta j}(\psi_2,\psi_3)
] \|_{N[0]} \nn \\
&  \les \sum_{m\ge  C} \sum_{0\le\ell\le m+C}  \| P_0 Q_{\ell} \del^\beta [Q_{> m-C} \psi_1 \Delta^{-1}\del_j \tilde Q_m \calN_{\beta j}(\psi_2,\psi_3)
] \|_{\dot X^{0,-1-\eps,2}_0 } \nn \\
& \quad + \sum_{m\ge  C} \sum_{\ell\ge m+C}  \| P_0 Q_{\ell} \del^\beta [\tilde Q_\ell \psi_1 \Delta^{-1}\del_j \tilde Q_m \calN_{\beta j}(\psi_2,\psi_3)
] \|_{ \dot X^{0,-1-\eps,2}_0} \nn  \\
&  \les \sum_{m\ge  C} \sum_{0\le\ell\le m+C} 2^{(\frac12-\eps)\ell} \| Q_{> m-C} \psi_1\|_{L^2_t L^\infty_x} \| \wt P_0    \tilde Q_m \calN_{\beta j}(\psi_2,\psi_3)
] \|_{\Ltwotx} \nn \\
& \quad + \sum_{m\ge  C} \sum_{\ell\ge m+C}  2^{-\eps\ell} \| \tilde Q_\ell \psi_1\|_{L^2_t L^\infty_x} \| \wt P_0    \tilde Q_m \calN_{\beta j}(\psi_2,\psi_3)
] \|_{\ener}
\end{align*}
In the second to last line we applied Bernstein's inequality in the
time variable to switch from $L^2_t$ to~$L^1_t$.  We now replace the
$L^\infty_t$ on the right-hand side of the last line by an $L^2_t$
at the expense of a factor of~$2^{\frac{m}{2}}$. Together with
Lemma~\ref{lem:Nabhighmod}  this yields
\begin{align*}
 \eqref{eq:B2} &\les  \sum_{m\ge  C} \sum_{0\le\ell\le m+C} 2^{ (\frac12-\eps)\ell}  2^{k_1} \| Q_{> m-C} \psi_1\|_{\Ltwotx} \| \wt P_0 \tilde Q_m \calN_{\beta j}(\psi_2,\psi_3)
] \|_{\Ltwotx} \nn \\
& \quad + \sum_{m\ge  C} \sum_{\ell\ge m+C}  2^{ -\eps\ell} 2^{k_1}
2^{\frac{m}{2}} \| \tilde Q_\ell \psi_1\|_{\Ltwotx} \| \wt P_0
\tilde Q_m \calN_{\beta j}(\psi_2,\psi_3)
] \|_{\Ltwotx}  \\
&\les  2^{k_1} \sum_{m\ge C}    2^{-\frac12m}
\| \psi_1\|_{S[k_1]}
\| \wt P_0 \tilde Q_m \calN_{\beta j}(\psi_2,\psi_3)
] \|_{\Ltwotx} \nn \\
& \quad +2^{k_1} \sum_{m\ge  C} \sum_{\ell\ge m+C}  2^{ -\eps\ell} 2^{\frac{m}{2}} 2^{-(1-\eps)\ell}   \| \psi_1\|_{S[k_1]}
 \|  \wt P_0 \tilde Q_m \calN_{\beta j}(\psi_2,\psi_3)
] \|_{\Ltwotx}  \\
&\les  2^{(\frac32-\eps)k_1}  2^{-\eps k_2} \la k_2 \ra^2
\prod_{i=1}^3 \|\psi_i\|_{S[k_i]}
\end{align*}
which is admissible. Next, we consider the case where both $A_0=I^c$
and $A_1=I^c$. If $\alpha\ne0$, then $\tilde A_1=I$ and one can
drop~$R_\alpha$ altogether so that  the previous analysis applies.
Otherwise, if~$\alpha=0$, then by assumption $\tilde A_1=I$ and
\begin{align}
   & \| P_0 Q_{\ge 0} \del^\beta [Q_{\ge C} R_\alpha \psi_1 \Delta^{-1}\del_j
 \wt P_0 Q_{\le C}  \calN_{\beta j}( A_2\psi_2,A_3\psi_3) ]
\|_{N[0]} \nn \\
&\le \sum_{m\ge 10C} \| P_0 Q_{m} \del^\beta [\tilde Q_m
R_\alpha \psi_1 \Delta^{-1}\del_j  \wt P_0  Q_{\le C}
\calN_{\beta j}( A_2\psi_2,A_3\psi_3) ]
\|_{N[0]} \label{eq:sch6'}\\
& +  \| P_0 Q_{0\le\cdot \le  10C} \del^\beta [Q_{0\le \cdot \le
 10C} R_\alpha \psi_1 \Delta^{-1}\del_j \wt P_0
Q_{\le  C} \calN_{\beta j}( A_2\psi_2,A_3\psi_3) ] \|_{N[0]}
\label{eq:sch7'}
\end{align}
By Lemma~\ref{lem:Nablowmod}, \eqref{eq:sch7'} is bounded by
\begin{align*}
 & \| P_0 Q_{0\le\cdot\le  10C} \del^\beta [Q_{0\le \cdot\le  10C}
R_\alpha \psi_1 \Delta^{-1}\del_j  \wt P_0  Q_{\le C}
\calN_{\beta j}( A_2\psi_2,A_3\psi_3) ] \|_{\dot X_0^{0,-1-\eps,2}} \\
&\les \| Q_{0\le \cdot\le  10C} R_\alpha \psi_1\;
\Delta^{-1}\del_j   \wt P_0  Q_{\le C} \calN_{\beta j}(
A_2\psi_2,A_3\psi_3)
\|_{L^2_t L^2_x} \\
&\les  \| Q_{0\le \cdot\le k_1+10C} R_\alpha
\psi_1\|_{L^\infty_t L^\infty_x}\,
 \| \Delta^{-1}\del_j   \wt P_0  Q_{\le
C} \calN_{\beta j}( A_2\psi_2,A_3\psi_3) \|_{L^2_t L^2_x} \\
&\les 2^{k_1}  2^{-\frac{k_2}{2}} \prod_{i=1}^3 \|\psi_i\|_{S[k_i]}
\end{align*}
 On the other hand,~\eqref{eq:sch6'} is
estimated as follows:
\begin{align*}
 & \sum_{m\ge  10C} \| P_0 Q_{m} \del^\beta [\tilde Q_m
R_\alpha \psi_1 \Delta^{-1}\del_j  \wt P_0 Q_{\le C}
\calN_{\beta j}( A_2\psi_2,A_3\psi_3) ] \|_{\dot
X_0^{0,-1-\eps,2}}\\
& \les \sum_{m\ge 0} 2^{-m\eps}\| \tilde Q_m \nabla_{t,x} |\nabla|^{-1} \psi_1
\|_{L^2_t L^\infty_x } \| \wt P_0 Q_{\le k_1+C}
\calN_{\beta j}( A_2\psi_2,A_3\psi_3) ] \|_{L^\infty_t L^2_x}\\
& \les  2^{(\frac12-\eps)k_1}  \| \psi_1 \|_{S[k_1]} \|
\wt  P_{0} Q_{\le C} \calN_{\beta j}(
A_2\psi_2,A_3\psi_3) ] \|_{L^2_t L^2_x} \\
&\les 2^{(\frac12-\eps)k_1 -\frac{k_2}{2}} \prod_{i=1}^3 \|\psi_i\|_{S[k_i]}
\end{align*}
where we applied Bernstein's inequality relative to~$t$ as well as
Lemma~\ref{lem:Nablowmod}.
\noindent Now suppose $A_0=I$ (in fact, $A_0=Q_{\le0}$), but at
least one of $A_1$ or~$\tilde A_1$ equals $I^c$. If $\tilde A_1=I^c$, then the
modulations of $\psi_1$ and $\calN_{\beta j}$ essentially agree,
whence $\alpha\ne0$ and
\begin{align*}
 & \sum_{m\ge C} \| P_0 Q_{\le 0} \del^\beta [ Q_m R_\alpha \psi_1  \Delta^{-1}\del_j \tilde Q_m  \calN_{\beta j}( \psi_2,\psi_3)
] \|_{N[0]} \\
& \les \sum_{m\ge k_1+C} \| P_0 Q_{\le 0} \del^\beta [ Q_m R_\alpha \psi_1
\Delta^{-1}\del_j \tilde Q_m  \calN_{\beta j}( \psi_2,\psi_3)
] \|_{L^1_{t}L^2_x} \\
& \les  \sum_{m\ge C} \|  Q_m \psi_1 \|_{L^2_t L^\infty_x}   \|
\wt P_0  \tilde Q_m  \calN_{\beta j}( \psi_2,\psi_3)
 \|_{\Ltwotx} \\
&\les \sum_{m\ge C} 2^{(\frac32-\eps)k_1} 2^{-m(1-2\eps)} \|
\psi_1 \|_{S[k_1]} \,  2^{-m\eps}\| \wt P_0  \tilde Q_m
\calN_{\beta j}( \psi_2,\psi_3)
 \|_{\Ltwotx} \\
 &\les 2^{ (\frac32-\eps)k_1-\eps k_2}  \la k_2 \ra^2
\prod_{i=1}^3 \|\psi_i\|_{S[k_i]}
\end{align*}
The final estimate here uses Lemma~\ref{lem:Nabhighmod2}.
Now suppose that $\tilde A_1=I$ and $A_1=I^c$. Then
\begin{align*}
 & \| P_0 Q_{\le 0} \del^\beta [I^c R_\alpha\psi_1 \Delta^{-1}\del_j  I \calN_{\beta j}( \psi_2,\psi_3)
] \|_{N[0]} \\
&\les  \|I^c R_\alpha\psi_1 \Delta^{-1}\del_j  \, \wt P_0 I \calN_{\beta j}( \psi_2,\psi_3)
 \|_{L^1_t L^2_x} \\
&\les  \|I^c \nabla_{t,x}|\nabla|^{-1} \psi_1\|_{L^2_t L^\infty_x} \| \wt P_0 Q_{\le C}  \calN_{\beta j}( \psi_2,\psi_3)
 \|_{L^2_t L^2_x} \\
 &\les 2^{(\frac12-\eps)k_1} \|\psi_1\|_{\Ltwotx} 2^{-\frac{k_2}{2}} \|\psi_2\|_{S[k_2]} \|\psi_3\|_{S[k_3]}
\end{align*}
The last
case which we need to consider is $A_0=A_1=\tilde A_1 =I$ and either
one of $A_2, A_3$ equal to~$I^c$. But then necessarily $A_2=A_3=I^c$
whence
\begin{align*}
 & \| P_0 \del^\beta I[I\psi_1 \Delta^{-1}\del_j I \calN_{\beta j}( Q_{\ge k_2+C} \psi_2, Q_{\ge k_2+C}\psi_3)
] \|_{N[0]}  \\
 & \les\| I\psi_1 \Delta^{-1}\del_j I \calN_{\beta j}( Q_{\ge k_2+C} \psi_2, Q_{\ge k_2+C}\psi_3)
 \|_{L^1_t L^1_x}  \\
&\les \|\psi_1\|_{\Linf} \sum_{m\ge k_2+C} \|\wt P_0 Q_{\le C}
\calN_{\beta j}( Q_{m} \psi_2, \tilde Q_m\psi_3)  \|_{L^1_t
L^1_x} \\
&\les 2^{k_1} \|\psi_1\|_{\ener} \sum_{m\ge k_2+C} 2^{m-k_2} 2^{-2m(1-\eps)}
2^{(1-2\eps)k_2} \|\psi_2\|_{S[k_2]} \|\psi_3\|_{S[k_3]} \\
& \les
2^{k_1-k_2} \prod_{i=1}^3 \|\psi_i\|_{S[k_i]}
\end{align*}
which concludes Case~4.

\medskip
\noindent {\em Case 5:} $\mathit{O(1)=k_1,\; k_2= k_3+O(1)}$.
We begin with $A_0=I^c$ and $\tilde A_1=I$ (in fact, $A_0=Q_{\ge0}$ suffices here as usual). Moreover, we will drop $R_\alpha$ from~$\psi_1$ which amounts
to excluding the case $A_1=I^c$ and~$\alpha=0$ but nothing else.
Then, from Lemma~\ref{lem:Nablowmod},
\begin{align*}
 & \| P_0 I^c \del^\beta [ \psi_1 \Delta^{-1}\del_j I \calN_{\beta j}( \psi_2, \psi_3)
] \|_{N[0]} \\
&\les \sum_{k\le k_2\wedge0+O(1)}  \|  \psi_1 \Delta^{-1}\del_j I P_k \calN_{\beta j}(  \psi_2, \psi_3)
] \|_{\Ltwotx} \\
&\les  \sum_{k\le k_2\wedge0+O(1)}   \|\psi_1\|_{\ener} \, 2^{-k} \| I P_k \calN_{\beta j}(  \psi_2, \psi_3) \|_{L^2_t L^\infty_x} \\
&\les  \sum_{k\le k_2\wedge0+O(1)}   \|\psi_1\|_{\ener}  2^{k-\frac{k_2}{2}} \|\psi_2\|_{S[k_2]} \|\psi_3\|_{S[k_3]}
 \les 2^{-\frac{|k_2|}{2}} \prod_{i=1}^3 \|\psi_i\|_{S[k_i]}
\end{align*}
which is better than needed.  Now suppose $\alpha=0$ and~$A_0=A_1=I^c$, which implies that $\tilde A_1=I$. Then
\begin{align*}
 & \| P_0 I^c \del^\beta [ I^c R_0 \psi_1 \Delta^{-1}\del_j I \calN_{\beta j}( \psi_2, \psi_3)
] \|_{N[0]} \\
& \les \sum_{k\le k_2\wedge0+O(1)} \sum_{m\ge0} 2^{-\eps m}  \| P_0 Q_m [ \tilde Q_m R_0 \psi_1\, \Delta^{-1}\del_j P_k I \calN_{\beta j}( \psi_2, \psi_3)
] \|_{\Ltwotx} \\
&\les \sum_{k\le k_2\wedge0+O(1)}  \sum_{m\ge0} 2^{(1-\eps)m}   \|  \tilde Q_m  \psi_1\|_{\Ltwotx}  2^{-k} \|  P_k I \calN_{\beta j}(  \psi_2, \psi_3)
] \|_{\Linf} \\
&\les  \sum_{k\le k_2\wedge 0 + O(1)}    \|\psi_1\|_{S[k_1]} \, 2^{\frac{k}{2}} \| I P_k \calN_{\beta j}(  \psi_2, \psi_3) \|_{\Ltwotx} \\
&\les  \sum_{k\le k_2\wedge 0 + O(1)}  2^{\frac{k}{2}}  \|\psi_1\|_{S[k_1]}  2^{k-\frac{k_2}{2}} \|\psi_2\|_{S[k_2]} \|\psi_3\|_{S[k_3]}  \les 2^{-\frac{|k_2|}{2}}  \prod_{i=1}^3 \|\psi_i\|_{S[k_i]}
\end{align*}
Next, consider the case $A_0=I^c$,  and $\tilde A_1=I^c$. Since $I^c R_0\psi_1 $ is now excluded, we may drop $A_1 R_\alpha$ altogether. Then
\begin{align}
 & \| P_0 I^c \del^\beta [ \psi_1 \Delta^{-1}\del_j I^c \calN_{\beta j}( \psi_2, \psi_3)
] \|_{N[0]} \nn \\
&\les \sum_{k\le k_2\wedge 0 + O(1)}   \|P_0 Q_{\ge0}  [   \psi_1 \Delta^{-1}\del_j Q_{k\le \cdot\le C} P_k \calN_{\beta j}(  \psi_2, \psi_3)
] \|_{\Ltwotx} \label{eq:sch13}\\
&\quad+ \sum_{k\le k_2\wedge 0 + O(1)}  \sum_{m\ge C}  2^{-\eps m} \|P_0 Q_m [    Q_{\le m-C}  \psi_1 \Delta^{-1}\del_j \tilde Q_m P_k \calN_{\beta j}(  \psi_2, \psi_3)
] \|_{\Ltwotx} \label{eq:sch14} \\
&\quad+ \sum_{k\le k_2\wedge 0 + O(1)}  \sum_{m\ge C}  \|P_0 Q_{\ge0}  [    Q_{> m-C}  \psi_1 \Delta^{-1}\del_j \tilde Q_m P_k \calN_{\beta j}(  \psi_2, \psi_3)
] \|_{\Ltwotx}  \label{eq:sch15}
\end{align}
First,  by Lemma~\ref{lem:Nabhighmod},
\begin{align*}
 \eqref{eq:sch13} &\les \sum_{k\le k_2\wedge 0 + O(1)}   \|  \psi_1\|_{\ener}\,  2^{-k} \|  Q_{k\le \cdot\le C} P_k \calN_{\beta j}(  \psi_2, \psi_3)
] \|_{L^2_t L^\infty_x} \\
 &\les \sum_{k\le k_2\wedge 0 + O(1)}   \|  \psi_1\|_{\ener}\,   \|  Q_{k\le \cdot\le C} P_k \calN_{\beta j}(  \psi_2, \psi_3)
] \|_{\Ltwotx} \\
&\les \sum_{k\le k_2\wedge 0 + O(1)}   \|  \psi_1\|_{\ener}\, 2^{\frac{k}{2}} 2^{-\eps k_2} \la k_2-k\ra^2  \|\psi_2\|_{S[k_2]} \|\psi_3\|_{S[k_3]} \les 2^{-\eps|k_2|}\la k_2\ra    \prod_{i=1}^3 \|\psi_i\|_{S[k_i]}
\end{align*}
Second, again  by Lemma~\ref{lem:Nabhighmod},
\begin{align*}
 \eqref{eq:sch14} &\les \sum_{k\le k_2\wedge 0 + O(1)}  \sum_{m\ge C}  2^{-\eps m} \|P_0 Q_m [    Q_{\le m-C}  \psi_1 \Delta^{-1}\del_j \tilde Q_m P_k \calN_{\beta j}(  \psi_2, \psi_3)
] \|_{\Ltwotx} \\
&\les \sum_{k\le k_2\wedge 0 + O(1)}  \sum_{m\ge C}  2^{-\eps m} \|\psi_1\|_{\ener}  \| \tilde Q_m P_k \calN_{\beta j}(  \psi_2, \psi_3)
] \|_{\Ltwotx} \\
&\les \sum_{k\le k_2\wedge 0 + O(1)}     2^{\frac{k}{2}} 2^{-\eps k_2} \la k_2-k\ra^2 \prod_{i=1}^3 \|\psi_i\|_{S[k_i]}
\les 2^{-\eps|k_2|}\la k_2\ra   \prod_{i=1}^3 \|\psi_i\|_{S[k_i]}
\end{align*}
and third,
\begin{align*}
\eqref{eq:sch15} &\les \sum_{k\le k_2\wedge 0 + O(1)}  \sum_{m\ge C}  \|   Q_{> m-C}  \psi_1\|_{\Ltwotx} 2^{-k} \| \tilde Q_m P_k \calN_{\beta j}(  \psi_2, \psi_3)
] \|_{\Linf} \\
&\les \sum_{k\le k_2\wedge 0 + O(1)}  \sum_{m\ge C}  2^{-(1-\eps)m} \|   \psi_1\|_{S[k_1]} \, 2^{\frac{m}{2}} \|   \tilde Q_m P_k \calN_{\beta j}(  \psi_2, \psi_3)
] \|_{\Ltwotx} \\
&\les \sum_{k\le k_2\wedge 0 + O(1)}  \sum_{m\ge C}  2^{-(\frac12-2\eps)m} \|   \psi_1\|_{S[k_1]} \,  2^{-\eps m} \|   \tilde Q_m P_k \calN_{\beta j}(  \psi_2, \psi_3)
] \|_{\Ltwotx} \\
&\les 2^{-\eps|k_2|}\la k_2\ra   \prod_{i=1}^3 \|\psi_i\|_{S[k_i]}
\end{align*}
where one argues as in the previous two cases to pass to the last line.

\noindent Thus, $A_0=Q_{\le0}$ for the remainder of Case~5. If $A_1=I^c$, then necessarily~$\tilde A_1=I^c$ which implies $\alpha\ne0$.
Therefore,
\begin{equation}\label{eq:wejustdid}
\begin{aligned}
 & \| P_0 Q_{\le 0} \del^\beta [ I^c \psi_1 \Delta^{-1}\del_j I^c \calN_{\beta j}( \psi_2, \psi_3)
] \|_{N[0]} \\
&\les \sum_{k\le k_2\wedge 0 + O(1)}  \sum_{m\ge C}  \| Q_m \psi_1 \Delta^{-1}\del_j  P_k \tilde Q_m \calN_{\beta j}(  \psi_2, \psi_3)
] \|_{\enerN} \\
&\les  \sum_{k\le k_2\wedge 0 + O(1)}   \sum_{m\ge C}  \|Q_m \psi_1\|_{\Ltwotx} \,  \| P_k \tilde Q_m \calN_{\beta j}(  \psi_2, \psi_3) \|_{\Ltwotx}\\
&\les  \sum_{k\le k_2\wedge 0 + O(1)}   \sum_{m\ge C} 2^{-(1-2\eps)m}  \| \psi_1\|_{S[k_1]} \, 2^{-\eps m} \| P_k \tilde Q_m \calN_{\beta j}(  \psi_2, \psi_3) \|_{\Ltwotx} \\
&\les \sum_{k\le k_2\wedge 0 + O(1)}     \| \psi_1\|_{S[k_1]} \; 2^{\frac{k}{2}} 2^{-\eps k_2} \la k_2-k\ra^2 \|\psi_2\|_{S[k_2]} \|\psi_3\|_{S[k_3]}
\\ &\les 2^{-\eps|k_2|}\la k_2\ra   \prod_{i=1}^3 \|\psi_i\|_{S[k_i]}
\end{aligned}
\end{equation}
which is again admissible. So we may assume also that $A_1=I$ which means that we can drop~$R_\alpha$ from~$\psi_1$. Furthermore,
in view of~\eqref{eq:wejustdid} it suffices to set~$\tilde A_1=I$.
Finally, suppose at least one choice of $j=2,3$ satisfies~$A_j=I^c$. Then necessarily, $A_2=A_3=I^c$ and
\begin{align}
 & \| P_0 I \del^\beta [ I \psi_1 \Delta^{-1}\del_j I \calN_{\beta j}( I^c \psi_2, I^c \psi_3)
] \|_{N[0]}\nn  \\
&\les \sum_{k\le k_2\wedge 0 + O(1)}   \| I \psi_1 \Delta^{-1}\del_j I P_k \calN_{\beta j}(I^c  \psi_2, I^c\psi_3)
] \|_{\enerN} \nn \\
&\les  \sum_{k\le k_2\wedge 0 + O(1)}    \|\psi_1\|_{\ener} \, 2^{-k} \| I P_k \calN_{\beta j}( I^c \psi_2, I^c\psi_3) \|_{L^1_t L^\infty_x} \label{eq:5oppose}\\
&\les  \sum_{k\le k_2\wedge 0 + O(1)}    \|\psi_1\|_{\ener}  \,2^k   \| I P_k \calN_{\beta j}( I^c \psi_2, I^c \psi_3) \|_{\Leins} \nn \\
&\les 2^{k_2\wedge 0} \|\psi_1\|_{\ener}  \, \sum_{m\ge k_2+C}  \|  \calN_{\beta j}( Q_m  \psi_2, \tilde Q_m \psi_3) \|_{\Leins} \nn \\
&\les 2^{k_2\wedge 0} \|\psi_1\|_{\ener}  \, \sum_{m\ge k_2+C} 2^{m-k_2} 2^{-2(1-\eps)m} 2^{(1-2\eps)k_2} \|  \psi_2\|_{S[k_2]} \| \psi_3\|_{S[k_3]} \nn  \\
&\les 2^{-k_2\vee 0}  \prod_{i=1}^3 \|\psi_i\|_{S[k_i]}
\end{align}
as claimed.

\medskip
\noindent {\em Case 6:} $\mathit{O(1)= k_1\ge k_2+O(1)\ge k_3+C}$.
This case proceeds similarly to Case~5. We begin with $A_0=I^c$ and
$\tilde A_1=I$ (in fact, $A_0=Q_{\ge0}$ suffices here as usual).
Moreover, we will drop $R_\alpha$ from~$\psi_1$ which amounts to
excluding the case $A_1=I^c$ and~$\alpha=0$ but nothing else. Then,
from Lemma~\ref{lem:Nablowmod2},
\begin{align*}
&  \| P_0 I^c \del^\beta [ \psi_1 \Delta^{-1}\del_j I \calN_{\beta j}( \psi_2, \psi_3)
] \|_{N[0]}
\les   \|  \psi_1 \Delta^{-1}\del_j I  \calN_{\beta j}(  \psi_2, \psi_3)
] \|_{\Ltwotx} \\
&\les      \|\psi_1\|_{\ener} \, 2^{-k_2} \| I \tilde P_{k_2} \calN_{\beta j}(  \psi_2, \psi_3) \|_{L^2_t L^\infty_x}
\les     2^{(\frac12-\eps)k_3+\eps k_2} \prod_{i=1}^3 \|\psi_i\|_{S[k_i]}
\end{align*}
which is better than needed.  Now suppose $\alpha=0$ and~$A_0=A_1=I^c$, which implies that $\tilde A_1=I$. Then
\begin{align*}
 & \| P_0 I^c \del^\beta [ I^c R_0 \psi_1 \Delta^{-1}\del_j I \calN_{\beta j}( \psi_2, \psi_3)
] \|_{N[0]} \\
& \les   \sum_{m\ge0} 2^{-\eps m}  \| P_0 Q_m [ \tilde Q_m R_0 \psi_1\, \Delta^{-1}\del_j \tilde P_{k_2} I \calN_{\beta j}( \psi_2, \psi_3)
] \|_{\Ltwotx} \\
&\les    \sum_{m\ge0} 2^{(1-\eps)m}   \|  \tilde Q_m  \psi_1\|_{\Ltwotx}  2^{-k_2} \|  \tilde P_{k_2} I \calN_{\beta j}(  \psi_2, \psi_3)
] \|_{\Linf} \\
&\les   2^{\frac{k_2}{2}} 2^{(\frac12-\eps)k_3+\eps k_2} \prod_{i=1}^3 \|\psi_i\|_{S[k_i]}
\end{align*}
Next, consider the case $A_0=I^c$,  and $\tilde A_1=I^c$. As before, we can drop $A_1 R_\alpha$ in this case. Then
\begin{align}
 & \| P_0 I^c \del^\beta [ \psi_1 \Delta^{-1}\del_j I^c \calN_{\beta j}( \psi_2, \psi_3)
] \|_{N[0]} \nn \\
&\les  \|P_0 Q_{\ge0}  [   \psi_1 \Delta^{-1}\del_j Q_{k\le \cdot\le C} \tilde P_{k_2} \calN_{\beta j}(  \psi_2, \psi_3)
] \|_{\Ltwotx} \label{eq:sch131}\\
&\quad+   \sum_{m\ge C}  2^{-\eps m} \|P_0 Q_m [    Q_{\le m-C}  \psi_1 \Delta^{-1}\del_j \tilde Q_m
\tilde P_{k_2} \calN_{\beta j}(  \psi_2, \psi_3)
] \|_{\Ltwotx} \label{eq:sch141} \\
&\quad+  \sum_{m\ge C}  \|P_0 Q_{\ge0}  [    Q_{> m-C}  \psi_1 \Delta^{-1}\del_j \tilde Q_m \tilde P_{k_2} \calN_{\beta j}(  \psi_2, \psi_3)
] \|_{\Ltwotx}  \label{eq:sch151}
\end{align}
First,  by Lemma~\ref{lem:Nabhighmod2},
\begin{align*}
 \eqref{eq:sch131} &\les    \|  \psi_1\|_{\ener}\,  2^{-k_2} \|  Q_{k_2-O(1)\le \cdot\le C} \tilde P_{k_2} \calN_{\beta j}(  \psi_2, \psi_3)
] \|_{L^2_t L^\infty_x} \\
 &\les   \|  \psi_1\|_{\ener}\,   \|  Q_{k_2-O(1)\le \cdot\le C} \tilde P_{k_2} \calN_{\beta j}(  \psi_2, \psi_3)
] \|_{\Ltwotx} \\
&\les   2^{(\frac12-\eps)k_3 }  \prod_{i=1}^3 \|\psi_i\|_{S[k_i]}
\end{align*}
Second, again  by Lemma~\ref{lem:Nabhighmod2},
\begin{align*}
 \eqref{eq:sch141} &\les   \sum_{m\ge C}  2^{-\eps m} \|P_0 Q_m [    Q_{\le m-C}  \psi_1 \Delta^{-1}\del_j \tilde Q_m
\tilde P_{k_2} \calN_{\beta j}(  \psi_2, \psi_3)
] \|_{\Ltwotx} \\
&\les  \sum_{m\ge C}  2^{-\eps m} \|\psi_1\|_{\ener}  \| \tilde Q_m \tilde P_{k_2} \calN_{\beta j}(  \psi_2, \psi_3)
] \|_{\Ltwotx} \\
&\les  2^{(\frac12-\eps)k_3 }  \prod_{i=1}^3 \|\psi_i\|_{S[k_i]}
\end{align*}
and third,
\begin{align*}
\eqref{eq:sch151} &\les   \sum_{m\ge C}  \|   Q_{> m-C}  \psi_1\|_{\Ltwotx} 2^{-k_2} \| \tilde Q_m \tilde P_{k_2} \calN_{\beta j}(  \psi_2, \psi_3)
] \|_{\Linf} \\
&\les  \sum_{m\ge C}  2^{-(1-\eps)m} \|   \psi_1\|_{S[k_1]} \, 2^{\frac{m}{2}} \|   \tilde Q_m \tilde P_{k_2} \calN_{\beta j}(  \psi_2, \psi_3)
] \|_{\Ltwotx} \\
&\les  \sum_{m\ge C}  2^{-(\frac12-2\eps)m} \|   \psi_1\|_{S[k_1]} \,  2^{-\eps m} \|   \tilde Q_m \tilde P_{k_2} \calN_{\beta j}(  \psi_2, \psi_3)
] \|_{\Ltwotx} \\
&\les  2^{(\frac12-\eps)k_3 }  \prod_{i=1}^3 \|\psi_i\|_{S[k_i]}
\end{align*}
where one argues as in the previous two cases to pass to the last line.

\noindent Thus, $A_0=Q_{\le0}$ for the remainder of Case~5. If $A_1=I^c$, then necessarily~$\tilde A_1=I^c$ which implies $\alpha\ne0$.
Therefore,
\begin{equation}\label{eq:wejustdid2}
\begin{aligned}
 & \| P_0 Q_{\le 0} \del^\beta [ I^c \psi_1 \Delta^{-1}\del_j I^c \calN_{\beta j}( \psi_2, \psi_3)
] \|_{N[0]} \\
&\les  \sum_{m\ge C}  \| Q_m \psi_1 \Delta^{-1}\del_j  \tilde P_{k_2} \tilde Q_m \calN_{\beta j}(  \psi_2, \psi_3)
] \|_{\enerN} \\
&\les     \sum_{m\ge C}  \|Q_m \psi_1\|_{\Ltwotx} \,  \| \tilde P_{k_2} \tilde Q_m \calN_{\beta j}(  \psi_2, \psi_3) \|_{\Ltwotx}\\
&\les     \sum_{m\ge C} 2^{-(1-2\eps)m}  \| \psi_1\|_{S[k_1]} \, 2^{-\eps m} \| \tilde P_{k_2}
 \tilde Q_m \calN_{\beta j}(  \psi_2, \psi_3) \|_{\Ltwotx} \\
&\les     2^{(\frac12-\eps)k_3}  \prod_{i=1}^3 \|\psi_i\|_{S[k_i]}
\end{aligned}
\end{equation}
which is again admissible. So we may assume also that $A_1=I$ which
means that we can drop~$R_\alpha$ from~$\psi_1$. Furthermore, in
view of~\eqref{eq:wejustdid2} it suffices to set~$\tilde A_1=I$.
This leaves the cases $A_2=I^c$ or $A_3=I^c$ to be considered. In
the former case, necessarily $A_2=A_3=I^c$ and
\begin{align*}
 & \| P_0 I \del^\beta [ I \psi_1 \Delta^{-1}\del_j I \calN_{\beta j}( I^c \psi_2, I^c \psi_3)
] \|_{N[0]}\les   \| I \psi_1 \Delta^{-1}\del_j I \tilde P_{k_2}
\calN_{\beta j}(I^c  \psi_2, I^c\psi_3)
] \|_{\enerN} \\
&\les      \|\psi_1\|_{\ener} \, 2^{-k_2} \| I \tilde P_{k_2}
\calN_{\beta j}( I^c \psi_2, I^c\psi_3) \|_{L^1_t L^\infty_x}
\les      \|\psi_1\|_{\ener}  \,   \| I  \tilde P_{k_2} \calN_{\beta j}( I^c \psi_2, I^c \psi_3) \|_{\enerN} \\
&\les   \|\psi_1\|_{\ener}  \, \sum_{m\ge k_2 }  \|  \calN_{\beta j}( Q_m  \psi_2, \tilde Q_m \psi_3) \|_{\enerN} \\
&\les   \|\psi_1\|_{\ener}  \, \sum_{m\ge k_2 }  \big( \|  \nabla_{t,x}|\nabla|^{-1} Q_m  \psi_2\|_{\Ltwotx}
\| \tilde Q_m \psi_3 \|_{L^2_t L^\infty_x} + \|  Q_m  \psi_2\|_{\Ltwotx}
\|  \nabla_{t,x}|\nabla|^{-1} \tilde Q_m \psi_3 \|_{L^2_t L^\infty_x} \big) \\
&\les   \|\psi_1\|_{\ener}  \, \sum_{m\ge k_2 }
 2^{-(1-2\eps)m} 2^{(\frac12-\eps)k_2}  2^{(\frac12-\eps)k_3}
\|  \psi_2\|_{S[k_2]} \| \psi_3\|_{S[k_3]}   \les
2^{(\frac12-\eps)(k_3-k_2)}  \prod_{i=1}^3 \|\psi_i\|_{S[k_i]}
\end{align*}
which is acceptable. The one remaining case is $A_0=A_1=\tilde
A_1=A_2=I$ and~$A_3=I^c$. Of course one may also assume that
$\psi_3=Q_{\le k_2+C}\psi_3$.  Then we write
\begin{align}
   P_0 I \del^\beta [ I \psi_1 \Delta^{-1}\del_j   I \calN_{\beta j}( I \psi_2, I^c \psi_3)
]  &=  P_0 I [ \del^\beta  I \psi_1 \, \tilde
P_{k_2}\Delta^{-1}\del_j I \calN_{\beta j}( I \psi_2, I^c \psi_3)
]  \label{eq:sch47} \\
& + P_0 I [  I \psi_1\,  \Delta^{-1}\del_j \del^\beta \tilde P_{k_2}
I \calN_{\beta j}( I \psi_2, I^c \psi_3)] \label{eq:sch48}
\end{align}
The term on the right-hand side of \eqref{eq:sch47} is  difficult.
More specifically, the methods that we have employed up to this
point do not seem to yield the necessary bound. However,  Tao's
trilinear estimate~\eqref{eq:Taotrilin} implies that
\begin{equation}\label{eq:Taoappl1}
\|  \del^\beta P_0   \psi_1 \, R_\beta \psi_3 \,  \psi_2 \|_{N[0]}
\les 2^{\sigma(k_3-k_2)} 2^{k_2} \prod_{i=1}^3 \|\psi_i\|_{S[k_i]}
\end{equation}
for some constant $\sigma>0$ as well as
\begin{equation}\label{eq:Taoappl2}
\|  \del^\beta P_0   \psi_1 \, R_\beta \psi_2 \,  \psi_3 \|_{N[0]}
\les 2^{k_3 } \prod_{i=1}^3 \|\psi_i\|_{S[k_i]}
\end{equation}
Since  $2^{k_2} \tilde P_{k_2}\Delta^{-1}\del_j I$ can be replaced
by the convolution by a measure and all norms involved are
translation invariant, these estimates imply~\eqref{eq:sch47}.

The analysis of~\eqref{eq:sch48} is easier and  similar to the
considerations  at the end of Case~2. More precisely, we first
reduce the modulation of the entire output by means of
Lemma~\ref{lem:Nablowmod2}:
\begin{align*}
 & \| P_0  Q_{(1-3\eps)k_3\le\cdot\le C} [I \psi_1 \Delta^{-1}\del_j\del^\beta \tilde P_{k_2} I  \calN_{\beta j}( I \psi_2,I^c \psi_3)
] \|_{N[0]} \nn \\
&\les 2^{-\frac12(1-3\eps)k_3} \| \psi_1\|_{L^\infty_t L^2_x} \,
  \| I\calN_{\beta j}( I \psi_2,I^c \psi_3) \|_{L^2_t L^\infty_x}\\
&\les   2^{-\frac12(1-3\eps)k_3} \| \psi_1\|_{S[k_1]} \,
2^{(\frac12-\eps)k_3} 2^{(1+\eps) k_2} \|\psi_2\|_{S[k_2]}
\|\psi_3\|_{S[k_3]}\\
&\les 2^{(1+\eps) k_2} 2^{\frac{\eps}{2}k_3} \prod_{i=1}^3
\|\psi_i\|_{S[k_i]}
\end{align*}
Next, we reduce the modulation of $\psi_1$:
\begin{align*}
 & \| P_0  Q_{\le(1-3\eps)k_3} [Q_{\ge (1-3\eps)k_3 } \psi_1 \Delta^{-1}\del_j\del^\beta I \tilde P_{k_2}  \calN_{\beta j}( I \psi_2,I^c \psi_3)
] \|_{N[0]} \nn \\
 &\les \| P_0 \del^\beta Q_{\le(1-3\eps)k_3} [Q_{\ge (1-3\eps)k_3 } \psi_1 \Delta^{-1}\del_j\del^\beta I \tilde P_{k_2} \calN_{\beta j}( I \psi_2,I^c \psi_3)
] \|_{\enerN} \nn \\
&\les  2^{k_2} \|Q_{\ge (1-3\eps)k_3 } \psi_1\|_{L^2_t L^2_x} \,
 \| \tilde P_{k_2} I \calN_{\beta j}( I \psi_2,I^c \psi_3) \|_{\Ltwotx}\\
&\les    2^{-\frac12(1-3\eps)k_3} \| \psi_1\|_{S[k_1]} \,
2^{(\frac12-\eps)k_3} 2^{(1+\eps) k_2} \|\psi_2\|_{S[k_2]}
\|\psi_3\|_{S[k_3]}\\
&\les   2^{\frac{\eps}{2}k_3+(1+\eps) k_2} \prod_{i=1}^3
\|\psi_i\|_{S[k_i]}
\end{align*}
Finally, we reduce the modulation of the interior null-form using
Corollary~\ref{cor:core_max}:
\begin{align*}
& \| P_0   Q_{\le (1-3\eps)k_3} [Q_{\le (1-3\eps)k_3 } \psi_1
\Delta^{-1}\del_j\del^\beta\tilde P_{k_2} Q_{k_3\le \cdot\le k_2+C}
\calN_{\beta j}( I \psi_2,
I^c\psi_3) ] \|_{N[0]} \\
&\les  \|\psi_1\|_{S[k_1]} \sum_{k_3\le\ell\le  k_2+C }
2^{\frac{\ell-k_2}{4}}  2^{-\frac{\ell}{2}}  \| \tilde P_{k_2}
 Q_{\ell} \calN_{\beta j}( I \psi_2,
I^c\psi_3)  \|_{\Ltwotx}\\
&\les  2^{(\frac14-\eps)(k_3-k_2)}  \prod_{i=1}^3
\|\psi_i\|_{S[k_i]}
\end{align*}
which is again admissible.  After these preparations, we are faced
with the following decomposition:
\begin{align*}
 &  P_0  Q_{\le (1-3\eps)k_3} [Q_{\le (1-3\eps)k_3 }
 \psi_1 \Delta^{-1}\del_j\del^\beta Q_{\le k_3}  \calN_{\beta j}( I \psi_2, I^c\psi_3)
]  \nn \\
&=  P_0   Q_{\le (1-3\eps)k_3} [Q_{\le (1-3\eps)k_3  } \psi_1
\Delta^{-1}\del_j  Q_{\le k_3} \calN_{\beta j}( Q_{k_3\le
\cdot\le k_2+C} \psi_2, Q_{k_3+C\le\cdot\le k_2+C}\psi_3) ]\\
& = \sum_{\kappa,\kappa'\in\calC_{\ell}} P_{0,\kappa}  Q_{\le
(1-3\eps)k_3} [P_{k_1,\kappa'} Q_{\le (1-3\eps)k_3  } \psi_1
\Delta^{-1}\del_j\del^\beta Q_{\le k_3} \calN_{\beta j}( Q_{k_3\le
\cdot\le k_2+C} \psi_2, Q_{k_3+C\le\cdot\le k_2+C}\psi_3) ]
\end{align*}
where $\ell=\frac12[k_2+ (1-3\eps)k_3]$ and
$\dist(\kappa,\kappa')\les 2^\ell$. Placing the entire expression
in~$L^1_t L^2_x$ and using Bernstein's inequality results in the
following estimate: with $J:= Q_{k_3\le
\cdot\le k_2+C}$,
\begin{align*}
 &  \|P_0  Q_{\le (1-3\eps)k_3} [Q_{\le (1-3\eps)k_3-k_1}
 \psi_1 \Delta^{-1}\del_j\del^\beta Q_{\le k_3}  \calN_{\beta j}( I \psi_2, I^c\psi_3)
] \|_{L^1_t L^2_x}  \nn \\
& \le \Big\|  \Big( \sum_{\kappa,\kappa'\in\calC_{\ell}}
\|P_{0,\kappa}  [P_{k_1,\kappa'} Q_{\le (1-3\eps)k_3-k_1 } \psi_1
\Delta^{-1}\del_j\del^\beta Q_{\le k_3} \calN_{\beta j}( J \psi_2, J\psi_3)
]\|_{L^2_x}^2 \Big)^{\frac12}
\Big\|_{L^1_t} \\
& \le 2^{\frac{\ell}{2}} \Big\|  \Big(
\sum_{\kappa,\kappa'\in\calC_{\ell}} \|P_{0,\kappa} [P_{k_1,\kappa'}
Q_{\le (1-3\eps)k_3-k_1 } \psi_1 \Delta^{-1}\del_j\del^\beta Q_{\le
k_3} \calN_{\beta j}( J \psi_2,
J\psi_3) ]\|_{L^1_x}^2 \Big)^{\frac12}
\Big\|_{L^1_t}\\
& \le 2^{\frac{\ell}{2}}   \Big\|  \Big(
\sum_{\kappa'\in\calC_{\ell}} \|P_{k_1,\kappa'} Q_{\le
(1-3\eps)k_3-k_1 } \psi_1\|_{L^2_x}^2 \| \Delta^{-1}\del_j\del^\beta
Q_{\le k_3} \calN_{\beta j}( J \psi_2,
J\psi_3) ]\|_{L^2_x}^2 \Big)^{\frac12}
\Big\|_{L^1_t}\\
& \le  2^{\frac{\ell}{2}}  \| Q_{\le (1-3\eps)k_3-k_1 }
\psi_1\|_{L^\infty_t L^2_x} \| \Delta^{-1}\del_j\del^\beta Q_{\le
k_3} \calN_{\beta j}( J \psi_2,
J\psi_3) ]\|_{L^1_t L^2_x}\\
&\les  2^{\frac{\ell}{2}} \|\psi_1\|_{S[k_1]} \;
\|\nabla_{t,x}|\nabla|^{-1}  J
\psi_2\|_{\Ltwotx} \| \nabla_{t,x}|\nabla|^{-1} J\psi_3\|_{L^2_t L^\infty_x} \\
&\les 2^{\frac14(1-3\eps)k_3+\frac{k_2}{4}} \|\psi_1\|_{S[k_1]} \;
2^{-\frac{k_3}{2}} \| \psi_2\|_{S[k_2]} \, 2^{(\frac12-\eps)k_3}
2^{\eps k_2}  \| \psi_3\|_{S[k_3]}
\end{align*}
which is again admissible for small $\eps>0$.

\medskip
\noindent {\em Case 7:} $\mathit{ k_1=O(1)\ge k_3+O(1)\ge k_2+C}$.
This case is symmetric to the previous one.

\medskip
\noindent {\em Case 8:} $\mathit{k_2=O(1), \max(k_1, k_3)\le -C}$.
 We  begin with
$A_0=Q_{\ge0}$ and $\tilde A_1=I$, and  we  drop $R_\alpha$
from~$\psi_1$  excluding the case $A_1=I^c$ and~$\alpha=0$ but
nothing else. Then, from Lemma~\ref{lem:Nablowmod2},
\begin{align*}
&  \| P_0 I^c \del^\beta [ \psi_1 \Delta^{-1}\del_j I \calN_{\beta
j}( \psi_2, \psi_3) ] \|_{N[0]} \les   \|  \psi_1 \Delta^{-1}\del_j
I  \calN_{\beta j}(  \psi_2, \psi_3)
] \|_{\Ltwotx} \\
&\les      \|\psi_1\|_{\Linf} \,  \| I \tilde P_{0} \calN_{\beta j}(
\psi_2, \psi_3) \|_{L^2_t L^2_x} \les 2^{k_1} 2^{(\frac12-\eps)k_3}
\prod_{i=1}^3 \|\psi_i\|_{S[k_i]}
\end{align*}
which is better than needed.  Now suppose $\alpha=0$
and~$A_0=A_1=I^c$, which implies that $\tilde A_1=I$. Then by
Lemma~\ref{lem:Nablowmod2}
\begin{align*}
 & \| P_0 I^c \del^\beta [ I^c R_0 \psi_1 \Delta^{-1}\del_j I \calN_{\beta j}( \psi_2, \psi_3)
] \|_{N[0]} \\
& \les   \sum_{m\ge0} 2^{-\eps m}  \| P_0 Q_m [ \tilde Q_m R_0
\psi_1\, \Delta^{-1}\del_j \tilde P_{0} I \calN_{\beta j}( \psi_2,
\psi_3)
] \|_{\Ltwotx} \\
& \les   \sum_{m\ge0} 2^{-\eps m}  \|  \tilde Q_m
\nabla_{t,x}|\nabla|^{-1} \psi_1\|_{L^2_t L^\infty_x} \|
\Delta^{-1}\del_j \tilde P_{0} I \calN_{\beta j}( \psi_2, \psi_3)
] \|_{\ener} \\
&\les    \sum_{m\ge0} 2^{(1-\eps)m}    \|  \tilde Q_m
\psi_1\|_{\Ltwotx}    \|  \tilde P_{0} I \calN_{\beta j}( \psi_2,
\psi_3)
] \|_{\Ltwotx} \\
&\les    2^{(\frac12-\eps)(k_1+k_3) } \prod_{i=1}^3
\|\psi_i\|_{S[k_i]}
\end{align*}
Next, consider the case $A_0=I^c$,  and $\tilde A_1=I^c$. As before,
we can drop $A_1 R_\alpha$ in this case. Then
\begin{align}
 & \| P_0 I^c \del^\beta [ \psi_1 \Delta^{-1}\del_j I^c \calN_{\beta j}( \psi_2, \psi_3)
] \|_{N[0]} \nn \\
&\les     \sum_{m\ge C}  2^{-\eps m} \|P_0 Q_m [    Q_{\le m-C}
\psi_1 \Delta^{-1}\del_j \tilde Q_m \tilde P_{0} \calN_{\beta j}(
\psi_2, \psi_3)
] \|_{\Ltwotx} \label{eq:sch142} \\
&\quad+  \sum_{m\ge C}  \|P_0 Q_{\ge0}  [    Q_{> m-C}  \psi_1
\Delta^{-1}\del_j \tilde Q_m \tilde P_{0} \calN_{\beta j}( \psi_2,
\psi_3) ] \|_{\Ltwotx}  \label{eq:sch152}
\end{align}
First,   by Lemma~\ref{lem:Nabhighmod2},
\begin{align*}
 \eqref{eq:sch142} &\les   \sum_{m\ge C}  2^{-\eps m} \|P_0 Q_m [    Q_{\le m-C}  \psi_1 \Delta^{-1}\del_j \tilde Q_m
\tilde P_{0} \calN_{\beta j}(  \psi_2, \psi_3)
] \|_{\Ltwotx} \\
&\les  \sum_{m\ge C}  2^{-\eps m} 2^{k_1} \|\psi_1\|_{\ener}  \|
\tilde Q_m \tilde P_{0} \calN_{\beta j}(  \psi_2, \psi_3)
] \|_{\Ltwotx} \\
&\les  2^{k_1+(\frac12-\eps)k_3 }  \prod_{i=1}^3 \|\psi_i\|_{S[k_i]}
\end{align*}
and second,
\begin{align*}
\eqref{eq:sch152} &\les   \sum_{m\ge C} 2^{k_1} \|   Q_{> m-C}
\psi_1\|_{\Ltwotx}   \| \tilde Q_m \tilde P_{0} \calN_{\beta j}(
\psi_2, \psi_3)
] \|_{\ener} \\
&\les  \sum_{m\ge C} 2^{(\frac32-\eps)k_1} 2^{-(1-\eps)m} \|
\psi_1\|_{S[k_1]} \, 2^{\frac{m}{2}} \|   \tilde Q_m \tilde P_{0}
\calN_{\beta j}( \psi_2, \psi_3)
] \|_{\Ltwotx} \\
&\les  2^{k_1} \sum_{m\ge C}  2^{-(\frac12-2\eps)m} \|
\psi_1\|_{S[k_1]} \,  2^{-\eps m} \|   \tilde Q_m \tilde P_{0}
\calN_{\beta j}( \psi_2, \psi_3)
] \|_{\Ltwotx} \\
&\les  2^{k_1+(\frac12-\eps)k_3 }  \prod_{i=1}^3 \|\psi_i\|_{S[k_i]}
\end{align*}
where one argues as in the previous two cases to pass to the last
line.

\noindent Thus, $A_0=Q_{\le0}$ for the remainder of Case~8. If
$\tilde A_1=I^c$, then necessarily~$A_1=I^c$ which implies
$\alpha\ne0$. Therefore,
\begin{equation}\nn
\begin{aligned}
 & \| P_0 Q_{\le 0} \del^\beta [ I^c \psi_1 \Delta^{-1}\del_j I^c \calN_{\beta j}( \psi_2, \psi_3)
] \|_{N[0]} \\
&\les  \sum_{m\ge C}  \| \tilde Q_m \psi_1 \Delta^{-1}\del_j  \tilde
P_{0}  Q_m \calN_{\beta j}(  \psi_2, \psi_3)
] \|_{\enerN} \\
&\les     \sum_{m\ge C}  2^{k_1} \|\tilde Q_m \psi_1\|_{\Ltwotx} \,  \| \tilde P_{0}  Q_m \calN_{\beta j}(  \psi_2, \psi_3) \|_{\Ltwotx}\\
&\les     2^{k_1}  \sum_{m\ge C} 2^{-(1-2\eps)m}  \|
\psi_1\|_{S[k_1]} \, 2^{-\eps m} \| \tilde P_{0}
 \tilde Q_m \calN_{\beta j}(  \psi_2, \psi_3) \|_{\Ltwotx} \\
&\les   2^{k_1}  2^{(\frac12-\eps)k_3}  \prod_{i=1}^3
\|\psi_i\|_{S[k_i]}
\end{aligned}
\end{equation}
which is again admissible. So we may assume also that $\tilde
A_1=I$. Now suppose that $A_1=I^c$. Then we can take
$A_1=Q_{k_1\le\cdot\le C}$ whence
\begin{align*}
 & \| P_0 Q_{\le 0} \del^\beta [ Q_{k_1\le\cdot\le C} R_\alpha \psi_1 \Delta^{-1}\del_j I \calN_{\beta j}( \psi_2, \psi_3)
] \|_{N[0]} \\
&\les    \| Q_{k_1\le\cdot\le C}  R_\alpha \psi_1 \Delta^{-1}\del_j
\tilde P_{0}  I \calN_{\beta j}(  \psi_2, \psi_3)
] \|_{\enerN} \\
&\les        \|Q_{k_1\le\cdot\le C} \nabla_{t,x} |\nabla|^{-1} \psi_1\|_{L^2_t L^\infty_x} \,
\| \tilde P_{0}  I  \calN_{\beta j}(  \psi_2, \psi_3) \|_{\Ltwotx}\\
&\les        2^{(\frac12-\eps)k_1}  \| \psi_1\|_{S[k_1]} \, 2^{(\frac12-\eps)k_3} \| \psi_2\|_{S[k_2]} \| \psi_3\|_{S[k_3]} \\
&\les      2^{(\frac12-\eps)(k_1+k_3)}  \prod_{i=1}^3
\|\psi_i\|_{S[k_i]}
\end{align*}
So we may assume for the remainder of this case that $A_1=I$ which
means that we can drop~$R_\alpha$ from~$\psi_1$. This leaves the
cases $A_2=I^c$ or $A_3=I^c$ to be considered. In the former case,
necessarily $A_2=A_3=I^c$ and
\begin{align*}
 & \| P_0 I \del^\beta [ I \psi_1 \Delta^{-1}\del_j I \calN_{\beta j}( I^c \psi_2, I^c \psi_3)
] \|_{N[0]}  \les   \| I \psi_1 \Delta^{-1}\del_j I \tilde P_{0}
\calN_{\beta j}(I^c  \psi_2, I^c\psi_3)
] \|_{\enerN} \\
&\les   2^{k_1}    \|\psi_1\|_{\ener} \,   \| I \tilde P_{0}
\calN_{\beta j}( I^c \psi_2, I^c\psi_3) \|_{\enerN}
\les     2^{k_1}   \|\psi_1\|_{\ener}  \,   \| I \tilde P_{0}  \calN_{\beta j}( I^c \psi_2, I^c \psi_3) \|_{\enerN} \\
&\les   2^{k_1}  \|\psi_1\|_{\ener}  \, \sum_{m\ge 0 }  \|  \calN_{\beta j}( Q_m  \psi_2, \tilde Q_m \psi_3) \|_{\enerN} \\
&\les  2^{k_1}   \|\psi_1\|_{\ener}  \, \sum_{m\ge 0 }  \big( \|
\nabla_{t,x}|\nabla|^{-1} Q_m  \psi_2\|_{\Ltwotx} \| \tilde Q_m
\psi_3 \|_{L^2_t L^\infty_x} + \|  Q_m  \psi_2\|_{\Ltwotx}
\|  \nabla_{t,x}|\nabla|^{-1} \tilde Q_m \psi_3 \|_{L^2_t L^\infty_x} \big) \\
&\les  2^{k_1}  \|\psi_1\|_{\ener}  \, \sum_{m\ge 0 }
 2^{-(1-2\eps)m}    2^{(\frac12-\eps)k_3}
\|  \psi_2\|_{S[k_2]} \| \psi_3\|_{S[k_3]}   \les
2^{k_1+(\frac12-\eps)k_3}  \prod_{i=1}^3 \|\psi_i\|_{S[k_i]}
\end{align*}
which is acceptable. The one remaining case is $A_0=A_1=\tilde
A_1=A_2=I$ and~$A_3=I^c$. Of course one may also assume that
$\psi_3=Q_{\le  C}\psi_3$. The analysis in this case is similar to
the considerations at the end of Case~2. More precisely, we first
reduce the modulation of the entire output by means of
Lemma~\ref{lem:Nablowmod2}:
\begin{align*}
 & \| P_0 \del^\beta Q_{(1-3\eps)k_3\le\cdot\le C} [I \psi_1 \Delta^{-1}\del_j \tilde P_{0} I  \calN_{\beta j}( I \psi_2,I^c \psi_3)
] \|_{N[0]} \nn \\
&\les 2^{-\frac12(1-3\eps)k_3} \| \psi_1\|_{L^\infty_t L^\infty_x}
\,
  \| I\calN_{\beta j}( I \psi_2,I^c \psi_3) \|_{L^2_t L^2_x}\\
&\les  2^{k_1} 2^{-\frac12(1-3\eps)k_3} \| \psi_1\|_{S[k_1]} \,
2^{(\frac12-\eps)k_3}   \|\psi_2\|_{S[k_2]}
\|\psi_3\|_{S[k_3]}\\
&\les 2^{  k_1} 2^{\frac{\eps}{2}k_3} \prod_{i=1}^3
\|\psi_i\|_{S[k_i]}
\end{align*}
Next, we reduce the modulation of $\psi_1$:
\begin{align*}
 & \| P_0  Q_{\le(1-3\eps)k_3}\del^\beta [Q_{\ge (1-3\eps)k_3  } \psi_1 \Delta^{-1}\del_j I \tilde P_{0}  \calN_{\beta j}( I \psi_2,I^c \psi_3)
] \|_{N[0]} \nn \\
 &\les \| P_0 \del^\beta Q_{\le(1-3\eps)k_3} [Q_{\ge (1-3\eps)k_3 } \psi_1 \Delta^{-1}\del_j  I \tilde P_{0} \calN_{\beta j}( I \psi_2,I^c \psi_3)
] \|_{\enerN} \nn \\
&\les  2^{k_1} \|Q_{\ge (1-3\eps)k_3 } \psi_1\|_{L^2_t L^2_x} \,
 \| \tilde P_{0} I \calN_{\beta j}( I \psi_2,I^c \psi_3) \|_{\Ltwotx}\\
&\les    2^{k_1-\frac12(1-3\eps)k_3} \| \psi_1\|_{S[k_1]} \,
2^{(\frac12-\eps)k_3}   \|\psi_2\|_{S[k_2]}
\|\psi_3\|_{S[k_3]}\\
&\les   2^{\frac{\eps}{2}k_3+k_1} \prod_{i=1}^3 \|\psi_i\|_{S[k_i]}
\end{align*}
Finally, we reduce the modulation of the interior null-form using
Corollary~\ref{cor:core_max}:
\begin{align*}
& \| P_0   Q_{\le (1-3\eps)k_3} [Q_{\le (1-3\eps)k_3 } \psi_1
\Delta^{-1}\del_j \tilde P_{0} Q_{k_3\le \cdot\le  C} \calN_{\beta
j}( I \psi_2,
I^c\psi_3) ] \|_{N[0]} \\
&\les  2^{k_1} \|\psi_1\|_{S[k_1]} \sum_{k_3\le\ell\le   C }
2^{\frac{\ell-k_1}{4}}  2^{-\frac{\ell}{2}}  \| \tilde P_{0}
 Q_{\ell} \calN_{\beta j}( I \psi_2,
I^c\psi_3)  \|_{\Ltwotx}\\
&\les  2^{\frac{3k_1}{4}+(\frac14-\eps)k_3}  \prod_{i=1}^3
\|\psi_i\|_{S[k_i]}
\end{align*}
which is again admissible.  After these preparations, we are faced
with the following decomposition:
\begin{align*}
 &  P_0 \del^\beta Q_{\le (1-3\eps)k_3} [Q_{\le (1-3\eps)k_3 }
 \psi_1 \Delta^{-1}\del_j  Q_{\le k_3}  \calN_{\beta j}( I \psi_2, I^c\psi_3)
]  \nn \\
&=  P_0   Q_{\le (1-3\eps)k_3}\del^\beta [Q_{\le (1-3\eps)k_3  }
\psi_1 \Delta^{-1}\del_j  Q_{\le k_3} \calN_{\beta j}( Q_{k_3\le
\cdot\le k_2+C} \psi_2, Q_{k_3+C\le\cdot\le k_2+C}\psi_3) ]\\
& = \sum_{\kappa,\kappa'\in\calC_{\ell}} P_{0,\kappa}  Q_{\le
(1-3\eps)k_3}\del^\beta [P_{k_1,\kappa'} Q_{\le (1-3\eps)k_3  }
\psi_1 \Delta^{-1}\del_j  Q_{\le k_3} \calN_{\beta j}( Q_{k_3\le
\cdot\le k_2+C} \psi_2, Q_{k_3+C\le\cdot\le k_2+C}\psi_3) ]
\end{align*}
where $\ell=\frac12[(1-3\eps)k_3-k_1]\wedge 0$ and
$\dist(\kappa,\kappa')\les 2^\ell$. Placing the entire expression
in~$L^1_t L^2_x$ and using Bernstein's inequality results in the
following estimate:
\begin{align*}
 &  \|P_0  Q_{\le (1-3\eps)k_3} [Q_{\le (1-3\eps)k_3 }
 \psi_1 \Delta^{-1}\del_j  Q_{\le k_3}  \calN_{\beta j}( I \psi_2, I^c\psi_3)
] \|_{L^1_t L^2_x}  \nn \\
& \le \Big\|  \Big( \sum_{\kappa,\kappa'\in\calC_{\ell}}
\|P_{0,\kappa}  [P_{k_1,\kappa'} Q_{\le (1-3\eps)k_3  } \psi_1
\Delta^{-1}\del_j  Q_{\le k_3} \calN_{\beta j}( Q_{k_3\le \cdot\le
 C} \psi_2, Q_{k_3+C\le\cdot\le  C}\psi_3) ]\|_{L^2_x}^2
\Big)^{\frac12}
\Big\|_{L^1_t} \\
& \le   \Big\|  \Big( \sum_{\kappa,\kappa'\in\calC_{\ell}} \|
P_{k_1,\kappa'} Q_{\le (1-3\eps)k_3  } \psi_1 \Delta^{-1}\del_j
Q_{\le k_3} \calN_{\beta j}( Q_{k_3\le \cdot\le C} \psi_2,
Q_{k_3+C\le\cdot\le
 C}\psi_3) \|_{L^2_x}^2 \Big)^{\frac12}
\Big\|_{L^1_t}\\
& \le     \Big\|  \Big( \sum_{\kappa'\in\calC_{\ell}}
\|P_{k_1,\kappa'} Q_{\le (1-3\eps)k_3 } \psi_1\|_{L^\infty_x}^2 \|
\Delta^{-1}\del_j  Q_{\le k_3} \calN_{\beta j}( Q_{k_3\le \cdot\le
C} \psi_2, Q_{k_3+C\le\cdot\le C}\psi_3) ]\|_{L^2_x}^2
\Big)^{\frac12}
\Big\|_{L^1_t}\\
& \le  2^{\frac{\ell}{2}}  2^{k_1} \| Q_{\le (1-3\eps)k_3  }
\psi_1\|_{L^\infty_t L^2_x} \| \Delta^{-1}\del_j  Q_{\le k_3}
\calN_{\beta j}( Q_{k_3\le \cdot\le C} \psi_2,
Q_{k_3+C\le\cdot\le  C}\psi_3) ]\|_{L^1_t L^2_x}\\
&\les  2^{k_1+\frac{\ell}{2}} \|\psi_1\|_{S[k_1]} \;
\|\nabla_{t,x}|\nabla|^{-1}  Q_{k_3\le \cdot\le  C}
\psi_2\|_{\Ltwotx} \| \nabla_{t,x}|\nabla|^{-1} Q_{k_3+C\le\cdot\le
 C}\psi_3\|_{L^2_t L^\infty_x} \\
&\les 2^{\frac{3k_1}{4}+\frac14(1-3\eps)k_3} \|\psi_1\|_{S[k_1]} \;
2^{-\frac{k_3}{2}} \| \psi_2\|_{S[k_2]} \, 2^{(\frac12-\eps)k_3}
   \| \psi_3\|_{S[k_3]}
\end{align*}
which is again admissible for small $\eps>0$.

\medskip
\noindent {\em Case 9:} $\mathit{k_3=O(1), \max(k_1, k_2)\le -C}$.
Symmetric to Case~8.
\end{proof}

It is important to realize that Lemma~\ref{lem:hyp_redux} yields the
following statement, which is really a corollary of its proof rather
than its lemma.

\begin{cor}
 \label{cor:hyp_redux}
Let $\psi_i$ be Schwarz functions adapted to  $k_i$ for $i=0,1,2$.
Then for any $\alpha,\beta=0,1,2$, and $j=1,2$,
\[
\| P_0 \nabla_{t,x} A_0[A_1 R_\alpha  \psi_1 \Delta^{-1}\del_j
\tilde A_1 \calN_{\beta j}( A_2\psi_2,A_3\psi_3) ] \|_{N[0]} \les
w(k_1,k_2,k_3) \prod_{i=1}^3 \|\psi_i\|_{S[k_i]}
\]
where  $A_i$ and $\tilde A_1$ are either $I$ or $I^c$, with at least
one being~$I^c$. Moreover, we impose the following restrictions:
\begin{itemize}
   \item if $A_1=\tilde A_1=I^c$ then $\alpha=0$ is excluded
   \item  if $k_1=O(1)> k_2\ge  k_3+C$, then  $A_0=A_1=\tilde
   A_1=A_2=I$, $A_3=I^c$ is excluded
 \item  if $k_1=O(1)> k_3\ge  k_2+C$, then  $A_0=A_1=\tilde
   A_1=A_3=I$, $A_2=I^c$ is excluded
\end{itemize}
In particular,
 \begin{equation}
\| P_0 \nabla_{t,x} [  \psi_1 \Delta^{-1}\del_j
I^c \calN_{\beta j}( \psi_2,\psi_3) ] \|_{N[0]} \les
w(k_1,k_2,k_3) \prod_{i=1}^3 \|\psi_i\|_{S[k_i]}
\label{eq:JoachimIc}
\end{equation}
\end{cor}
\begin{proof} Note that the first exclusion in our list is precisely
the exclusion in Lemma~\ref{lem:hyp_redux}.  The only real
difference between this statement and that of
Lemma~\ref{lem:hyp_redux} lies with the fact that we no longer
require the outer most derivative to be~$\del^\beta$. But this
mattered only in one case, namely when we applied Tao's
bound~\eqref{eq:Taotrilin} in Cases~6 and~7 above. Moreover,
inspection of the argument in those cases reveals that the
$\del^\beta \phi \del_\beta \psi$ null-form was needed only in those
instances which are excluded as the second and third conditions of
our above list (in fact, the modulations were narrowed down much
more before any need for~\eqref{eq:Taotrilin} arose).
The final statement is an immediate consequence of the first one, since we removed
$R_\alpha$ altogether (which eliminates the first exclusion) and since the
other two exclusions do not arise due to~$\tilde A_1=I^c$. Therefore, one simply
sums over all choices of~$A_0,A_1,A_2$ and~$A_3$.
\end{proof}

In fact, the proof of Lemma~\ref{lem:hyp_redux} makes no use of the
fact that $\Delta^{-1}\del_j$ contains the same index as the
null-form~$\calN_{\beta j}$.  But the strengthening resulting from
replacing $\Delta^{-1}\del_j$ by $|\nabla|^{-1}$, say, is of no
benefit to us so we do not carry it out. The following variant of
Lemma~\ref{lem:hyp_redux} covers the other two types of trilinear
nonlinearities arising in the Coulomb gauged wave-map system.

\begin{lemma}
 \label{lem:hyp_redux2}
Let $\psi_i$ be Schwarz functions adapted to  $k_i$ for $i=0,1,2$.
Then for any $\alpha=0,1,2$,  $j=1,2$,
\begin{align}
\| P_0 \del^\beta A_0[A_1 R_\beta \psi_1 \Delta^{-1}\del_j I
\calN_{\alpha j}( A_2\psi_2,A_3\psi_3) ] \|_{N[0]} &\les
w(k_1,k_2,k_3) \prod_{i=1}^3 \|\psi_i\|_{S[k_i]}
\label{eq:trilin_art2} \\
\| P_0 \del^\alpha A_0[A_1 R^\beta \psi_1 \Delta^{-1}\del_j I
\calN_{\beta j}( A_2\psi_2,A_3\psi_3) ] \|_{N[0]} &\les
w(k_1,k_2,k_3) \prod_{i=1}^3 \|\psi_i\|_{S[k_i]}
\label{eq:trilin_art3}
\end{align}
where  $A_i$   are either $I$ or $I^c$, with at least
one being~$I^c$.
\end{lemma}
\begin{proof}
Both these bounds follow from Corollary~\ref{cor:hyp_redux} provided
we are not in those cases described as Items~2 and~3 in the list of
exclusions (observe that the first exclusion does not arise due to
our limitation to~$\tilde A_1=I$). So let us consider the second exclusion
$k_1=O(1)> k_2\ge  k_3+C$ and $A_0=A_1=\tilde
   A_1=A_2=I$, $A_3=I^c$ (the third one being symmetric to this
   case). Then~\eqref{eq:trilin_art3} is an immediate consequence
   of~\eqref{eq:Taotrilin}, see~\eqref{eq:Taoappl1} and~\eqref{eq:Taoappl2} above.
   As for~\eqref{eq:trilin_art2},   observe that due to the analysis
   of~\eqref{eq:sch48} we may assume that the outer~$\del^\beta$
   derivative hits~$\psi_1$. Hence, it suffices to bound
   \[
\|P_0 I[I  \del^\beta  R_\beta \psi_1 \Delta^{-1}\del_j
I\calN_{\alpha j}( I\psi_2, I^c\psi_3) ] \|_{N[0]}
   \]
   However, due to the property that $\| \Box I P_0 \phi\|_{\Ltwotx}
   \les \|\phi\|_{S[0]}$ and $\del^\beta\del_\beta=\Box$,
   this is easy:
   \begin{align*}
     & \|P_0 I[Q_{\le C}  \del^\beta  R_\beta \psi_1
\Delta^{-1}\del_j I\calN_{\alpha j}( I\psi_2, I^c\psi_3) ] \|_{N[0]}
\\
& \les \|P_0 I[Q_{\le C}  \del^\beta  R_\beta \psi_1
\Delta^{-1}\del_j I\calN_{\alpha j}( I\psi_2, I^c\psi_3) ]
\|_{\enerN}\\
&\les \| Q_{\le C}  \del^\beta  R_\beta \psi_1 \|_{\Ltwotx} \|
\Delta^{-1}\del_j \tilde P_{k_2} I\calN_{\alpha j}(
I\psi_2, I^c\psi_3) ] \|_{L^2_t L^\infty_x}\\
&\les  \|   \psi_1 \|_{S[k_1]} \, 2^{\frac12(1-3\eps)k_3} 2^{\eps
k_2} \|\psi_2\|_{S[k_2]} \|\psi_3\|_{S[k_3]}
   \end{align*}
   as desired.
\end{proof}

The following technical corollary will be important later.

\begin{cor}
 \label{cor:epstrilin}  For some absolute constant $\sigma_0>0$,  and arbitrary
Schwartz functions $\psi_i$,
\begin{equation}\label{eq:smallIctrilin}
 \sum_{j=1}^2 \| P_0 \nabla_{t,x} [   \psi_1 \Delta^{-1}\del_j
I^c \calN_{\beta j}(  \psi_2, \psi_3) ] \|_{N[0]} \les  K^2
   \sup_{k\in\Z}\max_{i=1,2,3} 2^{-\sigma_0|k|} \|P_k \psi_i\|_{S[k]}
\end{equation}
provided  $\max_{i=1,2,3} \sum_{k\in\Z} \|P_k \psi_i\|_{S[k]}^2\le
K^2$ and with an absolute implicit constant. Moreover, given any
$\delta>0$ there exists a constant~$L=L(\delta)\gg1$ such that
\[
 {\sum}'_{k_1,k_2,k_3} \sum_{j=1}^2 \| P_0 \nabla_{t,x} [ P_{k_1}   \psi_1 \Delta^{-1}\del_j
I^c  \calN_{\beta j}( P_{k_2} \psi_2, P_{k_3} \psi_3) ] \|_{N[0]}
\le \delta\,  K^2
   \sup_{k\in\Z}\max_{i=1,2,3} 2^{-\sigma_0 |k|} \|P_k \psi_i\|_{S[k]}
\]
where the sum ${\sum}'_{k_1,k_2,k_3}$ extends over all $k_1,k_2,k_3$ {\em outside of the range}
\begin{equation}
 \label{eq:outside}
 |k_1|\le L,\quad k_2,k_3\le L,\quad |k_2-k_3|\le L
\end{equation}
Finally, if ${\sum}''_{k_1,k_2,k_3}$ denotes the sum {\em over} this
range, then
\begin{align*}
& {\sum}''_{k_1,k_2,k_3} \sum_{k\le k_2-L'} \sum_{j=1}^2 \| P_0
\nabla_{t,x} [ P_{k_1}   \psi_1 \Delta^{-1}\del_j I^c P_k
\calN_{\beta j}( P_{k_2} \psi_2, P_{k_3} \psi_3) ] \|_{N[0]} \\ &\le
\delta\,  K^2
   \sup_{k\in\Z}\max_{i=1,2,3} 2^{-\sigma_0 |k|} \|P_k \psi_i\|_{S[k]}
\end{align*}
where $L'=L'(L,\delta)$  is a large constant.
\end{cor}
\begin{proof}
Write $\psi_i=\sum_{k_i\in\Z} P_{k_i}\psi_i$ for $1\le i\le 3$.
 In view of the definition of the weights~$w(k_1,k_2,k_3)$, summing~\eqref{eq:JoachimIc} over all choices of $k_1,k_2,k_3$
yields~\eqref{eq:smallIctrilin}. The second statement follows
immediately from the fact that the weights $w(k_1,k_2,k_3)$ gain
some smallness outside of the range~\eqref{eq:outside}
(namely~$2^{-\delta L}$). For the third statement one needs to
observe that in Case~5 --- which is the one specified
by~\eqref{eq:outside} but of course with a range specified by the
constant~$L$ --- an extra gain can be obtained by restricting $k$ to
sufficiently small values compared to~$k_2,k_3$.
\end{proof}

\subsection{Trilinear estimates for hyperbolic $S$-waves}
\label{subsec:hyp_trilin}

The following lemma finally proves the trilinear estimates in the
``hyperbolic'' case. The argument will rely on the following trilinear
null-form expansion from~\cite{Krieger}:
\begin{equation}\label{eq:nullexp}
\begin{aligned}
  2\del^\beta \psi_1 \;\Delta^{-1}\del_j \calN_{\beta j} (\psi_2,\psi_3 )
& = (\Box\psi_1) |\nabla|^{-1}\psi_2 |\nabla|^{-1} \psi_3 - \Box (\psi_1 |\nabla|^{-1} \psi_2) |\nabla|^{-1}\psi_3  \\
& + \psi_1 \Box(|\nabla|^{-1} \psi_2) |\nabla|^{-1}\psi_3  + \Box(\psi_1 \Delta^{-1} \del_j (R_j \psi_2 |\nabla|^{-1}\psi_3)) \\
& - (\Box \psi_1) \Delta^{-1} \del_j (R_j\psi_2 |\nabla|^{-1}\psi_3)
- \psi_1 \Box\Delta^{-1}\del_j (R_j \psi_2 |\nabla|^{-1}\psi_3)
\end{aligned}
\end{equation}
as well as its ``dual'' form
\begin{equation}\label{eq:nullexp2}
\begin{aligned}
  2\del^\beta [\psi_1 \;\Delta^{-1}\del_j \calN_{\beta j} (\psi_2,\psi_3 ) ]
& = -\Box (\psi_1 |\nabla|^{-1}\psi_2 |\nabla|^{-1} \psi_3) + \Box (\psi_1  |\nabla|^{-1} \psi_3) |\nabla|^{-1}\psi_2  \\
& -  \psi_1 \Box(|\nabla|^{-1} \psi_2) |\nabla|^{-1}\psi_3  - (\Box \psi_1) \Delta^{-1} \del_j (R_j \psi_2 |\nabla|^{-1}\psi_3) \\
& + \Box (\psi_1  \Delta^{-1} \del_j (R_j\psi_2 |\nabla|^{-1}\psi_3))
+ \psi_1 \Box\Delta^{-1}\del_j (R_j \psi_2 |\nabla|^{-1}\psi_3)
\end{aligned}
\end{equation}
Strictly speaking, we shall want to apply these identities to the
trilinear expression
\[
\del^\beta [\psi_1 \;\Delta^{-1}\del_jIP_k \calN_{\beta j}
(\psi_2,\psi_3 ) ]
\]
for some~$P_k$. In the case of~\eqref{eq:nullexp2} the operator
$IP_k$ can be inserted in front of any product involving $\psi_2$
and~$\psi_3$ which is the case for all but the second term on the
right-hand side of~\eqref{eq:nullexp2}, i.e., $\Box (\psi_1
|\nabla|^{-1} \psi_3) |\nabla|^{-1}\psi_2$ (and similarly
for~\eqref{eq:nullexp}). Since $IP_k$ is disposable, it takes the
form of convolution with a measure~$\nu_k$ with
mass~$\|\nu_k\|\les1$. Thus, the second term needs to be replaced by
the convolution
\begin{equation}\label{eq:nuk_conv}
\int \Box (\psi_1 |\nabla|^{-1} \psi_3(\cdot-y))
|\nabla|^{-1}\psi_2(\cdot-y) \nu_k(dy)
\end{equation}
The logic will be that any estimate that we make on~$\Box (\psi_1
|\nabla|^{-1} \psi_3) |\nabla|^{-1}\psi_2$ in the context of the
$S[k]$ and~$N[k]$ spaces will equally well apply to this convolution
since all norms are translation invariant. We shall use this
observation repeatedly in what follows without any further comment.
Finally, the weights $ w(k_1,k_2,k_3)$ are those specified at the
beginning of this section.

\begin{lemma}
 \label{lem:tri_hyp} Let $\psi_j$ be adapted to $k_j$, for
 $j=1,2,3$.
Then
\begin{align}
\big\| \sum_{j=1}^2 P_0 I \del^\beta [I R_\alpha \psi_1 \Delta^{-1}\del_j I \calN_{\beta j}(I
\psi_2, I\psi_3)] \big\|_{N[0]} &  \les w(k_1,k_2,k_3) \prod_{i=1}^3
\|\psi_i\|_{S[k_i]} \label{eq:keytri1} \\
\big\| \sum_{j=1}^2 P_0 I \del_\alpha [I R^\beta \psi_1 \Delta^{-1}\del_j I \calN_{\beta j}(I
\psi_2, I\psi_3)] \big\|_{N[0]} &\les  w(k_1,k_2,k_3) \prod_{i=1}^3
\|\psi_i\|_{S[k_i]} \label{eq:keytri2} \\
 \big\| \sum_{j=1}^2 P_0 I \del^\beta [I R_\beta \psi_1 \Delta^{-1}\del_j I \calN_{\alpha j}(I
\psi_2, I\psi_3)] \big\|_{N[0]} &\les w(k_1,k_2,k_3) \prod_{i=1}^3 \label{eq:keytri3}
\|\psi_i\|_{S[k_i]}
\end{align}
for any $\alpha=0,1,2$.
\end{lemma}
\begin{proof}
We begin with~\eqref{eq:keytri1}.
Due to the $I$ in front of~$\psi_1$ we shall drop the~$R_\alpha$ operator.
Also, it will be understood in this proof that $\psi_i=Q_{\le
k_i+C}\psi_i$ for $1\le i\le 3$ and we will often drop
the~$I$-operator in front of the input functions.

\noindent
{\em Case 1:}  $\mathit{0\le k_1\le k_2+O(1)=k_3+O(1)}$. By
Lemma~\ref{lem:Nablowmod},
\begin{equation}
\begin{aligned}
 \| P_0 I \del^\beta [Q_{\ge0}  \psi_1 \Delta^{-1}\del_j I \calN_{\beta j}(I
\psi_2, I\psi_3)] \|_{N[0]} &\les \| P_0 I \del^\beta [Q_{\ge0}
\psi_1 \Delta^{-1}\del_j I \calN_{\beta j}(I \psi_2, I\psi_3)]
\|_{\enerN}\\
&\les \| Q_{\ge0} \psi_1\|_{\Ltwotx} \, 2^{-k_1} \|\tilde P_{k_1}  I
\calN_{\beta j}(I \psi_2, I\psi_3)] \|_{\Ltwotx}\\
&\les 2^{-\frac{k_2}{2}} \prod_{i=1}^3 \|\psi_i\|_{S[k_i]}
\end{aligned}\label{eq:case1_red}
\end{equation}
So it suffices to consider
\begin{align}
 P_0 I \del^\beta [Q_{<0}  \psi_1 \Delta^{-1}\del_j I \calN_{\beta j}(I
\psi_2, I\psi_3)] & =   P_0 Q_{\le   C} \del^\beta
[Q_{ <0} \psi_1
\Delta^{-1}\del_j Q_{\le C}  \calN_{\beta j}(I \psi_2, I\psi_3)] \label{eq:c11}
\end{align}
One can also limit the modulations of $\psi_2,\psi_3$ further. Indeed, by \eqref{eq:weaker_core} of Lemma~\ref{lem:core} and Corollary~\ref{cor:core_max},
\begin{equation}\begin{aligned}
 &\| P_0 Q_{\le   C} \del^\beta [Q_{<0}  \psi_1 \Delta^{-1}\del_j I \tilde P_{k_1}  \calN_{\beta j}( Q_{\ge \eps k_2} I
\psi_2, I\psi_3)] \|_{N[0]} \\
&\les 2^{-k_1}   \|
\psi_1 \nabla_{x,t}|\nabla|^{-1} I \psi_3\|_{\dot X^{0,\frac12,1}_{k_3}} \| Q_{\ge\eps k_2}\nabla_{x,t}|\nabla|^{-1} I \psi_2\|_{\dot X_{k_2}^{0,-\frac12,1}} \\
&\les 2^{-\frac12 \eps k_2} \la k_2-k_1\ra \prod_{i=1}^3 \|\psi_i\|_{S[k_i]}
\end{aligned}\label{eq:blabla}
\end{equation}
which is admissible. Note that we replaced~$\Delta^{-1}\del_j\tilde
P_{k_1}$ by $2^{-k_1}$ as explained in the paragraph preceding this
lemma.  Thus,  assume that $\psi_1=Q_{\le C}\psi_1$, $\psi_j=Q_{\le
\eps k_j}\psi_j$ for $j=2,3$, apply the
identity~\eqref{eq:nullexp2}, and estimate the six terms on the
right-hand side of~\eqref{eq:nullexp2} in the order in which they
appear. First, by the Strichartz component~\eqref{eq:Sk2},
\begin{align*}
 \|P_0 I \Box (\psi_1 |\nabla|^{-1}\psi_2 |\nabla|^{-1} \psi_3)\|_{N[0]} &\les \|P_0 I \Box (\psi_1 |\nabla|^{-1}\psi_2 |\nabla|^{-1} \psi_3)\|_{\dot X_0^{0,-\frac12,1}} \\
&\les \|\psi_1 |\nabla|^{-1}\psi_2 |\nabla|^{-1} \psi_3\|_{\Ltwotx} \\
&\les \|\psi\|_{\ener} 2^{-k_2} \|\psi_2\|_{L^4_t L^\infty_x} 2^{-k_3} \|\psi_3\|_{L^4_t L^\infty_x}  \\
&\les 2^{-\frac{k_2}{2}} \prod_{i=1}^3 \|\psi_i\|_{S[k_i]}
\end{align*}
Second, by \eqref{eq:core} of Lemma~\ref{lem:core} and Lemma~\ref{lem:Sk_prod2},
\begin{align*}
 \|P_0 I [\Box (\psi_1  |\nabla|^{-1} \psi_3) |\nabla|^{-1}\psi_2] \|_{N[0]} &\les
 \la k_3\ra \| \tilde P_{k_3} Q_{\le\eps k_3}  \Box (\psi_1  |\nabla|^{-1} \psi_3) \|_{\dot X_{k_3}^{0,-\frac12,1}} \| |\nabla|^{-1}\psi_2\|_{S[k_2]} \\
&\les  \la k_3\ra \| \tilde P_{k_3} Q_{\le\eps k_3}  (\psi_1  |\nabla|^{-1} \psi_3) \|_{\dot X_{k_3}^{0,\frac12,1}} \|  \psi_2\|_{S[k_2]}\\
&\les 2^{k_1-k_3} 2^{\frac14(\eps k_3 -k_1)}    \la k_3\ra \prod_{i=1}^3 \|\psi_i\|_{S[k_i]}
\end{align*}
Third, by \eqref{eq:weaker_core}  and Lemma~\ref{lem:Sk_prod2},
\begin{align*}
 & \| P_0 I [ \psi_1 \Box(|\nabla|^{-1} \psi_2) |\nabla|^{-1}\psi_3]\|_{N[0]}  \\
 & \les
   \|\tilde P_{k_3} Q_{\le\eps k_3} ( \psi_1 |\nabla|^{-1}\psi_3) \|_{\dot X_{k_3}^{0,\frac12,1}} \sum_{j\le \eps k_2} 2^{\frac14 j\wedge 0}
    \| \Box Q_j(|\nabla|^{-1} \psi_2) \|_{\dot X_{k_2}^{0,-\frac12,\infty}} \\
&\les 2^{k_1-k_3} 2^{\frac14(3\eps k_3 -k_1)}    \la k_3\ra^2
\prod_{i=1}^3 \|\psi_i\|_{S[k_i]}
\end{align*}
Fourth, again by \eqref{eq:weaker_core}  and Lemma~\ref{lem:Sk_prod2},
\begin{align*}
& \| P_0 I [ (\Box \psi_1) \Delta^{-1} \del_j (R_j \psi_2 |\nabla|^{-1}\psi_3)] \|_{N[0]} \\ &\les 2^{-k_1} \sum_{\ell\le C} 2^{\frac{\ell}{4}}
\| \Box Q_\ell \psi_1\|_{\dot X_{k_1}^{0,-\frac12,\infty}} \|\tilde P_{k_1}Q_{\le C} [ R_j \psi_2 |\nabla|^{-1}\psi_3] \|_{\dot X_{k_3}^{0,\frac12,1}} \\
&\les  2^{\frac{k_1-k_2}{4}} 2^{-\frac{k_2}{4}}  \prod_{i=1}^3 \|\psi_i\|_{S[k_i]}
\end{align*}
Fifth, with $\ell=k_1-k_2$,
\begin{align*}
& \| P_0 I  \Box [\psi_1  \Delta^{-1} \del_j (R_j\psi_2 |\nabla|^{-1}\psi_3)]  \|_{N[0]} \\
&\les
\| \psi_1  \Delta^{-1} \del_j\tilde P_{k_1}  (R_j\psi_2 |\nabla|^{-1}\psi_3) \|_{\Ltwotx} \\
&\les \|\psi_1\|_{\ener} 2^{-k_1} \sum_{c\in\calD_{k_2,\ell}} \| P_c R_j\psi_2 |\nabla|^{-1}P_{-c} \psi_3 \|_{L^2_t L^\infty_x} \\
&\les \|\psi_1\|_{\ener} 2^{-k_1-k_3} \Big(\sum_{c\in\calD_{k_2,\ell}} \| P_c R_j\psi_2\|_{L^4_t L^\infty_x}^2\Big)^{\frac12}  \Big(\sum_{c\in\calD_{k_2,\ell}}
\| P_{-c} \psi_3 \|_{L^4_t L^\infty_x}^2 \Big)^{\frac12} \\
&\les \|\psi_1\|_{\ener} 2^{-k_1-k_2} 2^{(1-2\eps)\ell} 2^{\frac{3k_2}{2}}  \|\psi_2\|_{S[k_2]} \|\psi_3\|_{S[k_3]}
\end{align*}
which is admissible for small $\eps>0$.
The sixth and final term is estimated by means of~\eqref{eq:core} and Lemma~\ref{lem:Sk_prod2}:
\begin{align*}
 \|P_0 I [\psi_1 \Box\Delta^{-1}\del_j (R_j \psi_2 |\nabla|^{-1}\psi_3)] \|_{N[0]}
&\les  \la k_1\ra \|\psi_1\|_{S[k_1]}  \| \tilde P_{k_1} Q_{\le C} \Box\Delta^{-1}\del_j (R_j \psi_2 |\nabla|^{-1}\psi_3) \|_{\dot X_{k_1}^{0,-\frac12,1} }   \\
&  \les  \la k_1\ra \|\psi_1\|_{S[k_1]}  \| \tilde P_{k_1} Q_{\le C}   (R_j \psi_2 |\nabla|^{-1}\psi_3) \|_{\dot X_{k_1}^{0,\frac12,1} }  \\
&  \les   \la k_1\ra \|\psi_1\|_{S[k_1]} \,   2^{\frac{5(k_1-k_2)}{4}} 2^{-\frac{k_2}{4}}   \|\psi_2\|_{S[k_2]} \|\psi_3\|_{S[k_3]}
\end{align*}
which concludes Case~1.

\medskip
\noindent {\em Case 2:} $\mathit{0\le k_1= k_3+O(1), k_2\le k_3-C}$.
By
Lemma~\ref{lem:Nablowmod2},
\begin{align*}
 \| P_0 I \del^\beta [Q_{\ge0}  \psi_1 \Delta^{-1}\del_j I \calN_{\beta j}(I
\psi_2, I\psi_3)] \|_{N[0]} &\les \| P_0 I \del^\beta [Q_{\ge0}
\psi_1 \Delta^{-1}\del_j I \calN_{\beta j}(I \psi_2, I\psi_3)]
\|_{\enerN}\\
&\les \| Q_{\ge0} \psi_1\|_{\Ltwotx} \, 2^{-k_1} \|\tilde P_{k_1}  I
\calN_{\beta j}(I \psi_2, I\psi_3)] \|_{\Ltwotx}\\
&\les 2^{(\frac12-\eps)k_2} 2^{-(1-\eps)k_1}  \prod_{i=1}^3 \|\psi_i\|_{S[k_i]}
\end{align*}
So it suffices to consider
\begin{align}
 P_0 I \del^\beta [Q_{<0}  \psi_1 \Delta^{-1}\del_j I \calN_{\beta j}(I
\psi_2, I\psi_3)] & =   P_0 Q_{\le   C} \del^\beta
[Q_{ <0} \psi_1
\Delta^{-1}\del_j Q_{\le C}  \calN_{\beta j}(I \psi_2, I\psi_3)] \label{eq:c21}
\end{align}
One can also limit the modulations of $\psi_3$ further. Indeed, by \eqref{eq:weaker_core} of Lemma~\ref{lem:core} and Corollary~\ref{cor:core_max},
\begin{equation}\begin{aligned}
 &\| P_0 Q_{\le   C} \del^\beta [Q_{<0}  \psi_1 \Delta^{-1}\del_j I \tilde P_{k_1}  \calN_{\beta j}(  I
\psi_2, Q_{\ge \eps k_3}I\psi_3)] \|_{N[0]} \\
&\les 2^{-k_1}   \| \tilde P_{k_3}
\psi_1 \nabla_{x,t}|\nabla|^{-1} I \psi_2\|_{\dot X^{0,\frac12,1}_{k_3}} \| Q_{\ge\eps k_3}\nabla_{x,t}|\nabla|^{-1} I \psi_3\|_{\dot X_{k_3}^{0,-\frac12,1}} \\
&\les 2^{k_2-k_1}
\la k_1-k_2\ra 2^{-\frac12 \eps k_1}  \prod_{i=1}^3 \|\psi_i\|_{S[k_i]}
\end{aligned}\label{eq:blabla2}
\end{equation}
which is admissible. As explained in Case~1, we replaced~$\Delta^{-1}\del_j\tilde P_{k_1}$ by $2^{-k_1}$.
If $0\le k_2$, then we can similarly reduce the modulation of the small frequency term, cf.~\eqref{eq:blabla}:
\begin{equation}\begin{aligned}
 &\| P_0 Q_{\le   C} \del^\beta [Q_{<0}  \psi_1 \Delta^{-1}\del_j I \tilde P_{k_1}  \calN_{\beta j}( Q_{\ge \eps k_2} I
\psi_2, I\psi_3)] \|_{N[0]} \\
&\les 2^{-k_1}   \| \tilde P_{k_2}[
\psi_1 \nabla_{x,t}|\nabla|^{-1} I \psi_3]\|_{\dot X^{0,\frac12,1}_{k_2}} \| Q_{\ge\eps k_2}\nabla_{x,t}|\nabla|^{-1} I \psi_2\|_{\dot X_{k_2}^{0,-\frac12,1}} \\
&\les  2^{\frac{k_2-k_1}{4}}   2^{-\frac12 \eps k_2}   \prod_{i=1}^3 \|\psi_i\|_{S[k_i]}
\end{aligned}\nn
\end{equation}
As a final preparation, we limit the modulation of the output in case $k_2\le0$. In fact,
by
Lemma~\ref{lem:Nablowmod2},
\begin{align*}
& \| P_0 Q_{(1-3\eps)k_2\le \cdot\le C}  \del^\beta [   \psi_1 \Delta^{-1}\del_j I \calN_{\beta j}(I
\psi_2, I\psi_3)] \|_{N[0]} \\&\les \| P_0 Q_{(1-3\eps)k_2\le \cdot\le C}  \del^\beta [
\psi_1 \Delta^{-1}\del_j I \calN_{\beta j}(I \psi_2, I\psi_3)]
\|_{\dot X_{k_2}^{0,-\frac12,1}}\\
&\les 2^{-\frac12(1-3\eps)k_2} \|  \psi_1\|_{\ener} \, 2^{-k_1} \|\tilde P_{k_1}  I
\calN_{\beta j}(I \psi_2, I\psi_3)] \|_{\Ltwotx}\\
&\les 2^{\frac12 \eps k_2} 2^{-(1-\eps)k_1}  \prod_{i=1}^3 \|\psi_i\|_{S[k_i]}
\end{align*}
Thus, for the remainder of this case we   assume that $\psi_1=Q_{\le C}\psi_1$, $\psi_2=Q_{\le \eps k_2\wedge k_2}\psi_2$,
and $\psi_3=Q_{\le \eps k_3}\psi_3$. Moreover, the output is restricted by~$Q_{\le (1-3\eps)k_2\wedge C}$.
We now stimate the six terms on the right-hand side of~\eqref{eq:nullexp2}. First, by the Strichartz component~\eqref{eq:Sk2},
\begin{align*}
 \|P_0 Q_{\le (1-3\eps)k_2\wedge C} \Box (\psi_1 |\nabla|^{-1}\psi_2 |\nabla|^{-1} \psi_3)\|_{N[0]} &\les \|P_0 Q_{\le (1-3\eps)k_2\wedge C} \Box (\psi_1 |\nabla|^{-1}\psi_2 |\nabla|^{-1} \psi_3)\|_{\dot X_0^{0,-\frac12,1}} \\
&\les 2^{\frac12(1-3\eps)k_2\wedge 0} \|\psi_1 |\nabla|^{-1}\psi_2 |\nabla|^{-1} \psi_3\|_{\Ltwotx} \\
&\les 2^{\frac12(1-3\eps)k_2\wedge 0} \|\psi\|_{\ener} 2^{-k_2} \|\psi_2\|_{L^4_t L^\infty_x} 2^{-k_3} \|\psi_3\|_{L^4_t L^\infty_x}  \\
&\les 2^{\frac12(1-3\eps)k_2\wedge 0} 2^{-\frac{k_1}{4}}2^{-\frac{k_2}{4}} \prod_{i=1}^3 \|\psi_i\|_{S[k_i]}
\end{align*}
which is admissible.
Second, by \eqref{eq:core} of Lemma~\ref{lem:core} and Lemma~\ref{lem:Sk_prod2},
\begin{align*}
 &\|P_0 Q_{\le (1-3\eps)k_2\wedge C} [\Box (\psi_1  |\nabla|^{-1} \psi_3) |\nabla|^{-1}\psi_2] \|_{N[0]} \\&\les
 2^{k_2\wedge 0} \la k_1\ra \| \tilde P_{k_2\vee0} Q_{\le (1-3\eps)k_2\wedge \eps k_2}  \Box (\psi_1  |\nabla|^{-1} \psi_3) \|_{\dot X_{k_2\vee0}^{0,-\frac12,1}} \| |\nabla|^{-1}\psi_2\|_{S[k_2]} \\
&\les  2^{-k_2\vee0} \la k_1\ra \| \tilde P_{k_2\vee 0} Q_{\le (1-3\eps)k_2\wedge \eps k_2}  (\psi_1  |\nabla|^{-1} \psi_3) \|_{\dot X_{k_2\vee0}^{0,\frac12,1}} \|  \psi_2\|_{S[k_2]}\\
&\les   2^{-k_2\vee0} \la k_1\ra 2^{\frac{k_2\vee0-k_1}{4}} 2^{\frac14 [(1-3\eps)k_2\wedge \eps k_2 -k_1]}  \prod_{i=1}^3 \|\psi_i\|_{S[k_i]}
\end{align*}
which is again admissible.
Third, by \eqref{eq:weaker_core}  and Lemma~\ref{lem:Sk_prod2},
\begin{align*}
 & \| P_0 Q_{\le (1-3\eps)k_2\wedge C}  [ \psi_1 \Box(|\nabla|^{-1} \psi_2) |\nabla|^{-1}\psi_3]\|_{N[0]} \\
 & \les
   \|\tilde P_{k_2\vee0} Q_{\le (1-3\eps)k_2\wedge \eps k_2}
    ( \psi_1 |\nabla|^{-1}\psi_3) \|_{\dot X_{k_2\vee0}^{0,\frac12,1}}\;\la k_2\vee0\ra  \|
    \Box(|\nabla|^{-1}Q_{\le \eps k_2\wedge k_2} \psi_2) \|_{\dot X_{k_2}^{0,-\frac12,\infty}} \\
&\les \la k_2\vee0\ra  2^{\frac{k_2\vee0-k_1}{4}} 2^{\frac14
[(1-3\eps)k_2\wedge \eps k_2 -k_1]} 2^{\frac{\eps k_2\wedge k_2}{2}
}    \prod_{i=1}^3 \|\psi_i\|_{S[k_i]}
\end{align*}
Fourth, again by \eqref{eq:weaker_core}  and Lemma~\ref{lem:Sk_prod2},
\begin{align*}
& \| P_0 I [ (\Box \psi_1) \Delta^{-1} \del_j (R_j \psi_2 |\nabla|^{-1}\psi_3)] \|_{N[0]} \\ &\les 2^{-k_1} \sum_{\ell\le C} 2^{\frac{\ell}{4}}
\| \Box Q_\ell \psi_1\|_{\dot X_{k_1}^{0,-\frac12,\infty}} \|\tilde P_{k_1}Q_{\le C} [ R_j \psi_2 |\nabla|^{-1}\psi_3] \|_{\dot X_{k_3}^{0,\frac12,1}} \\
&\les  2^{k_2-k_1} 2^{-\frac{k_1}{4}}  \prod_{i=1}^3 \|\psi_i\|_{S[k_i]}
\end{align*}
Fifth,
\begin{align*}
 & \| P_0 I  \Box [\psi_1  \Delta^{-1} \del_j (R_j\psi_2 |\nabla|^{-1}\psi_3)]  \|_{N[0]}
\les
\| \psi_1  \Delta^{-1} \del_j\tilde P_{k_1}  (R_j\psi_2 |\nabla|^{-1}\psi_3) \|_{\Ltwotx} \\
&\les \|\psi_1\|_{\ener}\, 2^{-k_1}   \|   R_j\psi_2 |\nabla|^{-1}  \psi_3 \|_{L^2_t L^\infty_x}
\les 2^{-2k_1} \|\psi_1\|_{\ener}    \|   R_j\psi_2\|_{L^4_t L^\infty_x}
\|  \psi_3 \|_{L^4_t L^\infty_x}  \\
&\les  2^{\frac{3k_2}{4}}    2^{-\frac{5k_1}{4}} \prod_{i=1}^3 \|\psi_i\|_{S[k_i]}
\end{align*}
The sixth and final term is estimated by means of~\eqref{eq:core} and Lemma~\ref{lem:Sk_prod2}:
\begin{align*}
 &\|P_0 I [\psi_1 \Box\Delta^{-1}\del_j (R_j \psi_2 |\nabla|^{-1}\psi_3)] \|_{N[0]}  \\
&\les  \la k_1\ra \|\psi_1\|_{S[k_1]}  \| \tilde P_{k_1} Q_{\le \eps k_3} \Box\Delta^{-1}\del_j (R_j \psi_2 |\nabla|^{-1}\psi_3) \|_{\dot X_{k_1}^{0,-\frac12,1} }   \\
&  \les  \la k_1\ra \|\psi_1\|_{S[k_1]}  \| \tilde P_{k_1} Q_{\le \eps k_3}   (R_j \psi_2 |\nabla|^{-1}\psi_3) \|_{\dot X_{k_1}^{0,\frac12,1} }  \\
&  \les   \la k_1\ra  2^{k_2-k_1} 2^{-\frac14(1-\eps)k_1}     \prod_{i=1}^3 \|\psi_i\|_{S[k_i]}
\end{align*}
which concludes Case~2.

\medskip
\noindent {\em Case 3:} $\mathit{0\le k_1= k_2+O(1), k_3\le k_2-C}$.  This
is symmetric to the preceding.

\medskip
\noindent {\em Case 4:} $\mathit{O(1)\le k_2= k_3+O(1), k_1\le -C }$.  This case proceeds similarly to Case~1.
Following~\eqref{eq:blabla}, we begin by limiting the modulations of $\psi_2,\psi_3$ to $2^{\eps k_2}$. Indeed, by \eqref{eq:weaker_core} of
Lemma~\ref{lem:core} and Corollary~\ref{cor:core_max},
\begin{equation}\begin{aligned}
 &\| P_0 Q_{\le   C} \del^\beta [I  \psi_1 \Delta^{-1}\del_j I \tilde P_{0}  \calN_{\beta j}( Q_{\ge \eps k_2} I
\psi_2, I\psi_3)] \|_{N[0]} \\
&\les  2^{k_1-k_2}      \|
\psi_1 \nabla_{x,t}|\nabla|^{-1} I \psi_3\|_{\dot X^{0,\frac12,1}_{k_3}} \| Q_{\ge\eps k_2}\nabla_{x,t}|\nabla|^{-1} I \psi_2\|_{\dot X_{k_2}^{0,-\frac12,1}} \\
&\les 2^{2k_1-k_2} 2^{-\frac12 \eps k_2} \la k_2-k_1\ra \prod_{i=1}^3 \|\psi_i\|_{S[k_i]}
\end{aligned}\nn
\end{equation}
which is admissible.  Next, we limit the modulation of the output:
by
Lemma~\ref{lem:Nablowmod},
\begin{align*}
 \| P_0 Q_{k_1 \le \cdot\le C}  \del^\beta [   \psi_1 \Delta^{-1}\del_j I \calN_{\beta j}(I
\psi_2, I\psi_3)] \|_{N[0]}
&\les \| P_0 Q_{k_1\le \cdot\le C}  \del^\beta [
\psi_1 \Delta^{-1}\del_j I \calN_{\beta j}(I \psi_2, I\psi_3)]
\|_{\dot X_{k_2}^{0,-\frac12,1}}\\
\les 2^{-\frac{k_1}{2}}  \|  \psi_1\|_{\Linf} \, \|\tilde P_{0}  I
\calN_{\beta j}(I \psi_2, I\psi_3)] \|_{\Ltwotx} &\les 2^{\frac{k_1-k_2}{2}} \prod_{i=1}^3 \|\psi_i\|_{S[k_i]}
\end{align*}
We now again estimate the six terms on the right-hand side of~\eqref{eq:nullexp2}. First, by the Strichartz component~\eqref{eq:Sk2},
\begin{align*}
 \|P_0 Q_{\le k_1} \Box (\psi_1 |\nabla|^{-1}\psi_2 |\nabla|^{-1} \psi_3)\|_{N[0]} &\les \|P_0 Q_{\le k_1}
 \Box (\psi_1 |\nabla|^{-1}\psi_2 |\nabla|^{-1} \psi_3)\|_{\dot X_0^{0,-\frac12,1}} \\
&\les 2^{\frac{k_1}{2}} \|\psi_1 |\nabla|^{-1}\psi_2 |\nabla|^{-1} \psi_3\|_{\Ltwotx} \\
&\les 2^{\frac{k_1}{2}} \|\psi\|_{\ener} 2^{-k_2} \|\psi_2\|_{L^4_t L^\infty_x} 2^{-k_3} \|\psi_3\|_{L^4_t L^\infty_x}  \\
&\les 2^{\frac{k_1-k_2}{2}} \prod_{i=1}^3 \|\psi_i\|_{S[k_i]}
\end{align*}
Second, by \eqref{eq:core} of Lemma~\ref{lem:core} and Lemma~\ref{lem:Sk_prod2},
\begin{align*}
& \|P_0 I [\Box (\psi_1  |\nabla|^{-1} \psi_3) |\nabla|^{-1}\psi_2] \|_{N[0]} \\&\les
 \la k_3\ra \| \tilde P_{k_3} Q_{\le\eps k_3}  \Box (\psi_1  |\nabla|^{-1} \psi_3) \|_{\dot X_{k_3}^{0,-\frac12,1}} \| |\nabla|^{-1}\psi_2\|_{S[k_2]} \\
&\les  \la k_3\ra \| \tilde P_{k_3} Q_{\le\eps k_3}  (\psi_1  |\nabla|^{-1} \psi_3) \|_{\dot X_{k_3}^{0,\frac12,1}} \|  \psi_2\|_{S[k_2]}\\
&\les 2^{k_1-k_3} 2^{\frac14(\eps k_3 -k_1)}    \la k_3-k_1\ra \prod_{i=1}^3 \|\psi_i\|_{S[k_i]}
\end{align*}
Third, by \eqref{eq:weaker_core}  and Lemma~\ref{lem:Sk_prod2},
\begin{align*}
 &\| P_0 I [ \psi_1 \Box(|\nabla|^{-1} \psi_2) |\nabla|^{-1}\psi_3]\|_{N[0]} \\ &
 \les \la k_2\ra
   \|\tilde P_{k_3} Q_{\le\eps k_3} ( \psi_1 |\nabla|^{-1}\psi_3) \|_{\dot X_{k_3}^{0,\frac12,1}}  \| \Box(|\nabla|^{-1} \psi_2)
   \|_{\dot X_{k_2}^{0,-\frac12,\infty}} \\
&\les 2^{k_1-k_3} 2^{\frac14(3\eps k_3 -k_1)} \la k_2\ra ^2 \la
k_1\ra \prod_{i=1}^3 \|\psi_i\|_{S[k_i]}
\end{align*}
Fourth, again by \eqref{eq:weaker_core}  and Lemma~\ref{lem:Sk_prod2},
\begin{align*}
& \| P_0 I [ (\Box \psi_1) \Delta^{-1} \del_j (R_j \psi_2 |\nabla|^{-1}\psi_3)] \|_{N[0]} \\ &\les  \sum_{\ell\le k_1+C} 2^{\frac{\ell-k_1}{4}}
\| \Box Q_\ell \psi_1\|_{\dot X_{k_1}^{0,-\frac12,\infty}} \|\tilde P_{0}Q_{\le C} [ R_j \psi_2 |\nabla|^{-1}\psi_3] \|_{\dot X_{k_3}^{0,\frac12,1}} \\ &
\les  2^{k_1-\frac{k_2}{4}}   \prod_{i=1}^3 \|\psi_i\|_{S[k_i]}
\end{align*}
Fifth, with $\ell=-k_2$,
\begin{align*}
& \| P_0 Q_{\le k_1}    \Box [\psi_1  \Delta^{-1} \del_j (R_j\psi_2 |\nabla|^{-1}\psi_3)]  \|_{N[0]} \\
&\les  2^{k_1}
\| \psi_1  \Delta^{-1} \del_j\tilde P_{k_1}  (R_j\psi_2 |\nabla|^{-1}\psi_3) \|_{\Ltwotx} \\
&\les 2^{k_1} \|\psi_1\|_{\ener}   \sum_{c\in\calD_{k_2,\ell}} \| P_c R_j\psi_2 |\nabla|^{-1}P_{-c} \psi_3 \|_{L^2_t L^\infty_x} \\
&\les 2^{k_1} \|\psi_1\|_{\ener} 2^{-k_3} \Big(\sum_{c\in\calD_{k_2,\ell}} \| P_c R_j\psi_2\|_{L^4_t L^\infty_x}^2\Big)^{\frac12}
\Big(\sum_{c\in\calD_{k_2,\ell}}
\| P_{-c} \psi_3 \|_{L^4_t L^\infty_x}^2 \Big)^{\frac12} \\
&\les  2^{k_1-k_2} 2^{(1-2\eps)\ell} 2^{\frac{3k_2}{2}}   \prod_{i=1}^3 \|\psi_i\|_{S[k_i]}
\end{align*}
which is admissible for small $\eps>0$.
The sixth and final term is estimated by means of~\eqref{eq:core} and Lemma~\ref{lem:Sk_prod2}:
\begin{align*}
& \|P_0 I [\psi_1 \Box\Delta^{-1}\del_j (R_j \psi_2 |\nabla|^{-1}\psi_3)] \|_{N[0]} \\
&\les  2^{k_1} \|\psi_1\|_{S[k_1]}  \| \tilde P_{0} Q_{\le k_1} \Box\Delta^{-1}\del_j (R_j \psi_2 |\nabla|^{-1}\psi_3) \|_{\dot X_{k_1}^{0,-\frac12,1} }   \\
&  \les 2^{2k_1}  \|\psi_1\|_{S[k_1]}  \| \tilde P_{0}     (R_j \psi_2 |\nabla|^{-1}\psi_3) \|_{\dot X_{0}^{0,\frac12,1} }  \\
&  \les   2^{2k_1} 2^{-\frac{k_2}{4}}     \prod_{i=1}^3 \|\psi_i\|_{S[k_i]}
\end{align*}
which concludes Case~4.

\medskip
\noindent {\em Case 5:} $\mathit{O(1)= k_1,\; k_2=k_3+O(1) }$.  We
start with the decomposition
\begin{equation}\label{eq:case5begin}
  P_0 \del^\beta [\psi_1 \Delta^{-1}\del_j I\calN_{\beta j}( \psi_2, \psi_3)]  = \sum_{k\le k_2\wedge 0 + O(1)}
    P_0 \del^\beta [\psi_1 \Delta^{-1}\del_j P_k I\calN_{\beta j}( \psi_2, \psi_3)]
\end{equation}
We first limit the modulation of~$\psi_1$:
\begin{align}
 & \sum_{k\le k_2\wedge 0 + O(1)}  \|P_0 \del^\beta Q_{>k} I[ Q_{>k+C} I\psi_1 \;
 \Delta^{-1}\del_j P_k I(R_\beta \psi_2 R_j\psi_3 - R_j \psi_2 R_\beta \psi_3)]\|_{N[0]}  \nn  \\
&\les   \sum_{k\le k_2\wedge 0 + O(1)}  \|P_0 \del^\beta Q_{>k} I [ Q_{>k+C}I
\psi_1 \;[\Delta^{-1}\del_{j\beta}^2 P_k I(|\nabla|^{-1} \psi_2
R_j\psi_3)
 - P_kI(|\nabla|^{-1} \psi_2 R_\beta \psi_3)]]\|_{\dot X^{0,-\frac12,1}} \nn \\
&\les  \sum_{k\le k_2\wedge 0 + O(1)}   2^{-\frac{k}{2}}
\|Q_{>k+C}\psi_1\|_{\Ltwotx} \| \Delta^{-1}\del_{j\beta}^2 P_k
I(|\nabla|^{-1} \psi_2 R_j\psi_3)
 - P_k I(|\nabla|^{-1} \psi_2 R_\beta \psi_3)]\|_{L^\infty_t L^\infty_x}  \nn  \\
&\les  \sum_{k\le k_2\wedge 0 + O(1)}   2^{k}
\| \psi_1\|_{S[k_1]} \| \Delta^{-1}\del_{j\beta}^2 P_k
I(|\nabla|^{-1} \psi_2 R_j\psi_3)
 - P_k I(|\nabla|^{-1} \psi_2 R_\beta \psi_3)]\|_{L^\infty_t L^1_x}  \nn \\
&\les  \sum_{k\le k_2\wedge 0 + O(1)}   2^{k-k_2} \|\psi_1\|_{S[k_1]}  \|
\psi_2\|_{L^\infty_t L^2_x} \| \psi_3\|_{L^\infty_t L^2_x}  \les 2^{-k_2\vee0}
\prod_{i=1}^3 \|\psi_i\|_{S[k_i]} \label{eq:psi1_modred}
\end{align}
Hence, if the inner output has frequency $\sim 2^k$ then  we may
assume that $\psi_1$ has modulation~$\les 2^k$. As usual, we
apply~\eqref{eq:nullexp2}. First, by the Strichartz
component~\eqref{eq:Sk2},
\begin{align*}
& \sum_{k\le k_2\wedge 0 + C}   \|P_0 I \Box (Q_{\le k}\psi_1 P_k
I[|\nabla|^{-1}\psi_2 |\nabla|^{-1} \psi_3])\|_{N[0]} \\
& \les \sum_{k\le k_2\wedge 0 + C}   \|P_0 Q_{\le k+C} \Box (Q_{\le k} \psi_1 P_k I[ |\nabla|^{-1}\psi_2 |\nabla|^{-1} \psi_3])\|_{\dot X_0^{0,-\frac12,1}} \\ &
\les \sum_{k\le k_2\wedge 0 + C}   2^{\frac{k}{2}}  \|Q_{\le k} \psi_1 P_k I[ |\nabla|^{-1}\psi_2 |\nabla|^{-1} \psi_3] \|_{\Ltwotx}
\les \sum_{k\le k_2\wedge 0 + C}   2^{\frac{k-k_2}{2}} \prod_{i=1}^3
\|\psi_i\|_{S[k_i]} \\
& \les 2^{-\frac12 k_2\vee0}  \prod_{i=1}^3 \|\psi_i\|_{S[k_i]}
\end{align*}
For the second term, we can assume that $\psi_1=Q_{\le
k_2\wedge 0+C}\psi_1$, see above.  Then, by~\eqref{eq:core} of
Lemma~\ref{lem:core} and Lemma~\ref{lem:Sk_prod2},
\begin{align*}
 & \|P_0 I [\Box ( \psi_1  |\nabla|^{-1} \psi_3) |\nabla|^{-1}\psi_2]
 \|_{N[0]} \\ &\les 2^{k_2\wedge0} \sum_{j\le k_2\wedge0+C} 2^{\frac{j-k_2\wedge 0}{4}}
  \| \tilde P_{k_2\vee 0}   \Box Q_j (\psi_1  |\nabla|^{-1} \psi_3) \|_{\dot X_{ k_2\vee 0}^{0,-\frac12,\infty}} \| |\nabla|^{-1}\psi_2\|_{S[k_2]} \\
&\les   2^{-k_2\vee0}    \| \tilde P_{0} Q_{\le k_2\wedge0+C}  (\psi_1  |\nabla|^{-1}
\psi_3) \|_{\dot X_{k_2\vee 0}^{0,\frac12,1}} \|  \psi_2\|_{S[k_2]} \\ &\les  2^{-2k_2\vee0}
\prod_{i=1}^3 \|\psi_i\|_{S[k_i]}
\end{align*}
Third, by \eqref{eq:weaker_core}  and Lemma~\ref{lem:Sk_prod2},
\begin{align*}
 & \| P_0 I [Q_{\le k_2\wedge0+C}  \psi_1 \Box(|\nabla|^{-1} \psi_2) |\nabla|^{-1}\psi_3]\|_{N[0]}  \\
 & \les \sum_{j\le k_2\wedge0+C} 2^{\frac{j-k_2\wedge 0}{4}}
   \|\tilde P_{k_2\vee 0} Q_{\le k_2\wedge0+C } ( \psi_1 |\nabla|^{-1}\psi_3) \|_{\dot X_{k_2\vee 0}^{0,\frac12,1}}
   \| \Box Q_j(|\nabla|^{-1} \psi_2) \|_{\dot X_{k_2}^{0,-\frac12,\infty}}\\ &
\les 2^{-k_2\vee0}  \prod_{i=1}^3 \|\psi_i\|_{S[k_i]}
\end{align*}
Fourth, again by \eqref{eq:weaker_core}  and
Lemma~\ref{lem:Sk_prod2},
\begin{align*}
& \sum_{k\le k_2\wedge 0 + C}   \| P_0 I [ (\Box Q_{\le k+C} \psi_1) \Delta^{-1}
\del_j P_k I(R_j \psi_2 |\nabla|^{-1}\psi_3)] \|_{N[0]} \\ &\les
\sum_{k\le k_2\wedge 0 + C}   \sum_{\ell\le C} 2^{\frac{\ell}{4}}
\| \Box Q_\ell \psi_1\|_{\dot X_{k_1}^{0,-\frac12,\infty}} \|\tilde P_{k}Q_{\le k+C} [ R_j \psi_2 |\nabla|^{-1}\psi_3] \|_{\dot X_{k}^{0,\frac12,1}} \\
&\les  \sum_{k\le k_2\wedge 0 + C}    \| \psi_1\|_{\dot
X_{k_1}^{0,\frac12,\infty}}   2^{\frac{k-k_2}{2}}
\|\psi_2\|_{S[k_2]} \|\psi_3\|_{S[k_3]}  2^{-\frac12 k_2\vee0} \les \prod_{i=1}^3
\|\psi_i\|_{S[k_i]}
\end{align*}
Fifth, with $\ell=k-k_2$,
\begin{align*}
& \sum_{k\le k_2\wedge 0 + C}   \| P_0 I  \Box [Q_{\le k+C} \psi_1  \Delta^{-1} \del_j P_k I (R_j\psi_2 |\nabla|^{-1}\psi_3)]  \|_{N[0]} \\
&\les \sum_{k\le k_2\wedge 0 + C}   2^{\frac{k}{2}}
\| Q_{\le k+C} \psi_1  \Delta^{-1} \del_j P_k I  (R_j\psi_2 |\nabla|^{-1}\psi_3) \|_{\Ltwotx} \\
&\les \|\psi_1\|_{\ener} \sum_{k\le k_2\wedge 0 + C}     2^{-\frac{k}{2}} \sum_{c\in\calD_{k_2,\ell}} \| P_c R_j\psi_2 |\nabla|^{-1}P_{-c} \psi_3 \|_{L^2_t L^\infty_x} \\
&\les \|\psi_1\|_{\ener} \sum_{k\le k_2\wedge 0 + C}    2^{-\frac{k}{2}-k_3}
\Big(\sum_{c\in\calD_{k_2,\ell}} \| P_c R_j\psi_2\|_{L^4_t
L^\infty_x}^2\Big)^{\frac12} \Big(\sum_{c\in\calD_{k_2,\ell}}
\| P_{-c} \psi_3 \|_{L^4_t L^\infty_x}^2 \Big)^{\frac12} \\
&\les \|\psi_1\|_{\ener} \sum_{k\le k_2\wedge 0 + C}
2^{(\frac12-2\eps)(k-k_2)}    \|\psi_2\|_{S[k_2]}
\|\psi_3\|_{S[k_3]} \les 2^{-(\frac12-2\eps)k_2\vee0}  \prod_{i=1}^3 \|\psi_i\|_{S[k_i]}
\end{align*}
which is admissible for small $\eps>0$. The sixth and final term is
estimated by means of~\eqref{eq:core} and Lemma~\ref{lem:Sk_prod2}:
\begin{align*}
&\sum_{k\le k_2\wedge 0 + C}    \|P_0 I [Q_{\le k+C}\psi_1 \Box\Delta^{-1}\del_j
P_k I (R_j \psi_2 |\nabla|^{-1}\psi_3)] \|_{N[0]} \\
&\les    \|\psi_1\|_{S[k_1]} \sum_{k\le k_2\wedge 0 + C}  2^k  \|   P_{k} Q_{\le
k+C}
\Box\Delta^{-1}\del_j (R_j \psi_2 |\nabla|^{-1}\psi_3) \|_{\dot X_{k}^{0,-\frac12,1} }   \\
&  \les     \|\psi_1\|_{S[k_1]}  \sum_{k\le k_2\wedge 0 + C}  2^k \| P_{k} Q_{\le k+ C}   (R_j \psi_2 |\nabla|^{-1}\psi_3) \|_{\dot X_{k}^{0,\frac12,1} }  \\
&  \les    \|\psi_1\|_{S[k_1]} \sum_{k\le k_2\wedge 0 + C}  2^k \,
 2^{\frac{k-k_2}{2}}   \|\psi_2\|_{S[k_2]}
\|\psi_3\|_{S[k_3]} \les  2^{-\frac12 k_2\vee0} \prod_{i=1}^3 \|\psi_i\|_{S[k_i]}
\end{align*}
which concludes Case~5.

\medskip
\noindent {\em Case 6:} $\mathit{O(1)= k_1\ge k_2+O(1)\ge k_3+C}$. Since Lemma~\ref{lem:Nablowmod2} implies that
\begin{align*}
 \|P_0 \del^\beta [Q_{>k_2}\psi_1 \Delta^{-1}\del_j I\calN_{\beta j}( \psi_2, \psi_3)] \|_{\enerN} &\les
\| Q_{>k_2}\psi_1 \|_{\Ltwotx} 2^{-k_2} \|I\calN_{\beta j}( \psi_2, \psi_3)\|_{L^2_t L^\infty_x} \\
&\les 2^{(\frac12-\eps)(k_3-k_2)} \prod_{i=1}^3 \|\psi_i\|_{S[k_i]}
\end{align*}
we may assume that $\psi_1=Q_{\le k_2} \psi_1$.  Next, we reduce matters to~\eqref{eq:Taotrilin}. More
precisely,
\begin{align}
   P_0 I \del^\beta [  \psi_1 \Delta^{-1}\del_j   I \calN_{\beta j}(  \psi_2,  \psi_3)
]  &=  P_0 I [ \del^\beta   \psi_1 \, \tilde
P_{k_2}\Delta^{-1}\del_j I\tilde P_{k_2} \calN_{\beta j}(  \psi_2,  \psi_3)
]  \label{eq:sch471} \\
& + P_0 I [  I \psi_1\,  \Delta^{-1}\del_j \del^\beta \tilde P_{k_2}
I \calN_{\beta j}(  \psi_2,  \psi_3)] \label{eq:sch481}
\end{align}
The term in~\eqref{eq:sch471} satisfies the bounds~\eqref{eq:Taoappl1} and~\eqref{eq:Taoappl2},
whereas~\eqref{eq:sch481} is expanded further:
\begin{align}
 P_0 I [  I \psi_1\,  \Delta^{-1}\del_j \del^\beta \tilde P_{k_2}
I \calN_{\beta j}(  \psi_2,  \psi_3)] &=  P_0 I [  I \psi_1\,  \Delta^{-1}\del_j \tilde P_{k_2}
I (\Box |\nabla|^{-1}\psi_2 R_j\psi_3  - R_j \psi_2 \Box|\nabla|^{-1}\psi_3  \label{eq:zweiterm}  \\
& \qquad\qquad   + R_\beta \psi_2 \del^\beta R_j\psi_3  -
\del^\beta R_j\psi_2 R_\beta\psi_3    ) ] \label{eq:zweiterm'}
\end{align}
The two terms in~\eqref{eq:zweiterm'} are again controlled by~\eqref{eq:Taotrilin}. Consider the first term
on the right-hand side of~\eqref{eq:zweiterm}. Replacing~$\Delta^{-1}\del_j \tilde P_{k_2}$ by~$2^{-k_2}$ as usual, one
obtains from Lemmas~\ref{lem:core} and~\ref{lem:Sk_prod2},
\begin{align*}
 \| \psi_1\, \Box |\nabla|^{-1}\psi_2 \, R_j\psi_3 \|_{N[0]} &\les  2^{k_2}  \sum_{j\le k_2+C} 2^{\frac{j-k_2}{4}}  \| \Box Q_j |\nabla|^{-1}\psi_2\|_{\dot X_{k_2}^{0,-\frac12,\infty}}  \|\psi_1\,  R_j\psi_3 \|_{\dot X_{0}^{0,\frac12,1}} \\
&\les  2^{k_2}    \|  \psi_2\|_{\dot X_{k_2}^{0,\frac12,\infty}}  2^{k_3} \la k_3\ra \|\psi_1\|_{S[k_1]} \|\psi_3 \|_{S[k_3]}
\end{align*}
which is more than enough. The second term in~\eqref{eq:zweiterm} is estimated similarly:
\begin{align*}
 \| \psi_1\, \Box |\nabla|^{-1}\psi_3 \, R_j\psi_2 \|_{N[0]} &\les  2^{k_3}  \sum_{j\le k_3+C} 2^{\frac{j-k_3}{4}}  \| \Box Q_j |\nabla|^{-1}\psi_3\|_{\dot X_{k_3}^{0,-\frac12,\infty}}  \|\psi_1\,  R_j\psi_2 \|_{\dot X_{0}^{0,\frac12,1}} \\
&\les  2^{k_3}    \|  \psi_2\|_{\dot X_{k_2}^{0,\frac12,\infty}}  2^{k_2} \la k_2\ra \|\psi_1\|_{S[k_1]} \|\psi_2 \|_{S[k_2]}
\end{align*}
which concludes Case~6.

\medskip
\noindent {\em Case 7:} $\mathit{ k_1=O(1)\ge k_3+O(1)\ge k_2+C}$.
This case is symmetric to the previous one.

\medskip
\noindent {\em Case 8:} $\mathit{k_3=O(1), \max(k_1, k_2)\le -C}$.
By Lemma~\ref{lem:Nablowmod2},
\begin{align*}
 &\| P_0 I \del^\beta [Q_{\ge k_1+ (1-3\eps)k_2}  \psi_1 \Delta^{-1}\del_j I \calN_{\beta j}(I
\psi_2, I\psi_3)] \|_{N[0]} \\&\les \| P_0 I \del^\beta [Q_{\ge k_1+(1-3\eps)k_2}
\psi_1 \Delta^{-1}\del_j I \calN_{\beta j}(I \psi_2, I\psi_3)]
\|_{\enerN}\\
&\les 2^{k_1} \| Q_{\ge k_1+(1-3\eps)k_2} \psi_1\|_{\Ltwotx} \,   \|\tilde P_{0}  I
\calN_{\beta j}(I \psi_2, I\psi_3)] \|_{\Ltwotx}\\
&\les 2^{\frac12\eps k_2} 2^{\frac{k_1}{2}}   \prod_{i=1}^3 \|\psi_i\|_{S[k_i]}
\end{align*}
A similar calculation shows that one can place $Q_{\le k_1+(1-3\eps)k_2}$ in front of the entire output.
So it suffices to consider
\begin{align*}
 & P_0 Q_{\le k_1+(1-3\eps)k_2} \del^\beta [Q_{< k_1+ (1-3\eps)k_2}  \psi_1 \Delta^{-1}\del_j I \calN_{\beta j}(I
\psi_2, I\psi_3)] \\& =   P_0 Q_{\le k_1+(1-3\eps)k_2} \del^\beta
[Q_{ < k_1+ (1-3\eps)k_2 } \psi_1
\Delta^{-1}\del_j Q_{\le k_1+ C}  \tilde P_0 \calN_{\beta j}(I \psi_2, I\psi_3)] \nn
\end{align*}
We now stimate the six terms on the right-hand side of~\eqref{eq:nullexp2}. First, by the Strichartz component~\eqref{eq:Sk2},
\begin{align*}
& \|P_0 Q_{\le k_1+(1-3\eps)k_2} \Box (\psi_1 |\nabla|^{-1}\psi_2 |\nabla|^{-1} \psi_3)\|_{N[0]} \\&\les \|P_0 Q_{\le k_1+(1-3\eps)k_2} \Box (\psi_1 |\nabla|^{-1}\psi_2 |\nabla|^{-1} \psi_3)\|_{\dot X_0^{0,-\frac12,1}} \\
&\les 2^{\frac12[(1-3\eps)k_2+k_1]} \|\psi_1 |\nabla|^{-1}\psi_2 |\nabla|^{-1} \psi_3\|_{\Ltwotx} \\
&\les 2^{\frac12[(1-3\eps)k_2+3k_1]} \|\psi_1\|_{\ener}  \|\psi_2\|_{L^4_t L^\infty_x} 2^{-\frac{k_2}{4}}  \|\psi_3\|_{L^4_t L^\infty_x}  \\
&\les 2^{\frac12[(1-3\eps)k_2+3k_1]}  2^{-\frac{k_2}{4}} \prod_{i=1}^3 \|\psi_i\|_{S[k_i]}
\end{align*}
which is sufficient.
Second, by \eqref{eq:core} of Lemma~\ref{lem:core} and Lemma~\ref{lem:Sk_prod2},
\begin{align*}
 &\|P_0 Q_{\le k_1+(1-3\eps)k_2} [\Box ( \psi_1  |\nabla|^{-1} \psi_3) |\nabla|^{-1}\psi_2] \|_{N[0]} \\&\les
  2^{k_2}  \| \tilde P_{0} Q_{\le (1-3\eps)k_2}  \Box (\psi_1  |\nabla|^{-1} \psi_3) \|_{\dot X_{0}^{0,-\frac12,1}} \| |\nabla|^{-1}\psi_2\|_{S[k_2]} \\
&\les   \| \tilde P_{ 0} Q_{\le (1-3\eps)k_2 }  (\psi_1  |\nabla|^{-1} \psi_3) \|_{\dot X_{ 0}^{0,\frac12,1}} \|  \psi_2\|_{S[k_2]}\\
&\les    2^{k_1} 2^{\frac{(1-3\eps)k_2-k_1 }{4}\wedge 0} \la k_1\ra   \prod_{i=1}^3 \|\psi_i\|_{S[k_i]}
\end{align*}
which is admissible.
Third, by \eqref{eq:weaker_core}  and Lemma~\ref{lem:Sk_prod2},
\begin{align*}
 & \| P_0 Q_{\le k_1+ (1-3\eps)k_2 }  [ \psi_1 \Box(|\nabla|^{-1} I\psi_2) |\nabla|^{-1}\psi_3]\|_{N[0]} \\
 & \les
   \|\tilde P_{0} Q_{\le (1-3\eps)k_2 }
    ( \psi_1 |\nabla|^{-1}\psi_3) \|_{\dot X_{0}^{0,\frac12,1}}\;  \|
    \Box(|\nabla|^{-1}I \psi_2) \|_{\dot X_{k_2}^{0,-\frac12,\infty}} \\
&\les  2^{k_1} 2^{\frac{(1-3\eps)k_2-k_1 }{4}\wedge 0} \la k_1\ra     \prod_{i=1}^3 \|\psi_i\|_{S[k_i]}
\end{align*}
Fourth, again by \eqref{eq:weaker_core}  and Lemma~\ref{lem:Sk_prod2},
\begin{align*}
& \| P_0 Q_{\le k_1+ (1-3\eps)k_2 }  [ (\Box Q_{\le k_1+ (1-3\eps)k_2 } \psi_1) \Delta^{-1} \del_j (R_j \psi_2 |\nabla|^{-1}\psi_3)] \|_{N[0]} \\ &\les 2^{k_1} \sum_{\ell\le k_1+ (1-3\eps)k_2} 2^{\frac{\ell-k_1}{4}}
\| \Box Q_\ell \psi_1\|_{\dot X_{k_1}^{0,-\frac12,\infty}} \|\tilde P_{0}Q_{\le k_1+C} [ R_j \psi_2 |\nabla|^{-1}\psi_3] \|_{\dot X_{k_3}^{0,\frac12,1}} \\
&\les  2^{2k_1} 2^{\frac14(1-3\eps)k_2} 2^{k_2} \la k_2\ra  \prod_{i=1}^3 \|\psi_i\|_{S[k_i]}
\end{align*}
Fifth,
\begin{align*}
 & \| P_0 Q_{\le k_1+ (1-3\eps)k_2 }  \Box [\psi_1  \Delta^{-1} \del_j (R_j\psi_2 |\nabla|^{-1}\psi_3)]  \|_{N[0]} \\
& \les 2^{\frac12[(1-3\eps)k_2+k_1]} \| \psi_1  \Delta^{-1} \del_j\tilde P_{k_1}  (R_j\psi_2 |\nabla|^{-1}\psi_3) \|_{\Ltwotx} \\
&\les 2^{\frac12[(1-3\eps)k_2+k_1]} \|\psi_1\|_{\ener}\,     \|   R_j\psi_2 |\nabla|^{-1}  \psi_3 \|_{L^2_t L^\infty_x}\\
&\les 2^{\frac12[(1-3\eps)k_2+k_1]}  \|\psi_1\|_{\ener}    \|   R_j\psi_2\|_{L^4_t L^\infty_x}
\|  \psi_3 \|_{L^4_t L^\infty_x}  \\
&\les 2^{\frac12[(1-3\eps)k_2+k_1]}     2^{\frac{3k_2}{4}}  \prod_{i=1}^3 \|\psi_i\|_{S[k_i]}
\end{align*}
The sixth and final term is estimated by means of~\eqref{eq:core} and Lemma~\ref{lem:Sk_prod2}:
\begin{align*}
 &\|P_0 Q_{\le k_1+ (1-3\eps)k_2 } [Q_{\le k_1+ (1-3\eps)k_2 } \psi_1 \Box\Delta^{-1}\del_j (R_j \psi_2 |\nabla|^{-1}\psi_3)] \|_{N[0]}  \\
&\les  2^{k_1}   \|\psi_1\|_{S[k_1]}  \| \tilde P_{0} Q_{\le k_1+C} \Box\Delta^{-1}\del_j (R_j \psi_2 |\nabla|^{-1}\psi_3) \|_{\dot X_{k_1}^{0,-\frac12,1} }   \\
&  \les  2^{k_1}  \|\psi_1\|_{S[k_1]}  \| \tilde P_{k_3} Q_{\le k_1}   (R_j \psi_2 |\nabla|^{-1}\psi_3) \|_{\dot X_{k_3}^{0,\frac12,1} }  \\
&  \les     2^{k_1+k_2} \la k_2\ra     \prod_{i=1}^3 \|\psi_i\|_{S[k_i]}
\end{align*}
which concludes Case~8.

\medskip
\noindent {\em Case 9:} $\mathit{k_2=O(1), \max(k_1, k_3)\le -C}$.
Symmetric to Case~8.

\medskip Hence we are done with \eqref{eq:keytri1}. Next, we turn to~\eqref{eq:keytri2} which is similar; basically,
one uses~\eqref{eq:nullexp} instead of~\eqref{eq:nullexp2}.  First,
one observes that any reductions in modulation which preceded
application of~\eqref{eq:nullexp2} to~\eqref{eq:keytri1} can equally
well be carried out for~\eqref{eq:keytri2} since these bounds only
use Lemmas~\ref{lem:Nablowmod} and~\ref{lem:Nablowmod2}.  Second,
observe that the last four terms of~\eqref{eq:nullexp} reappear as
the last four terms of~\eqref{eq:nullexp2} up to the order and the
choice of signs, both of which are irrelevant. Consequently, one
only needs to verify that the first two terms of~\eqref{eq:nullexp}
satisfy the desired bounds.

\medskip
\noindent {\em Case 1:}  $\mathit{0\le k_1\le k_2+O(1)=k_3+O(1)}$.
In this case the second terms in~\eqref{eq:nullexp}
and~\eqref{eq:nullexp2} satisfy the same bounds, whence it will
suffice to bound the first term in~\eqref{eq:nullexp}. However, by
\eqref{eq:weaker_core}  and Lemma~\ref{lem:Sk_prod2},
\begin{align*}
& \| P_0 I [ (\Box \psi_1) |\nabla|^{-1}  \psi_2
|\nabla|^{-1}\psi_3] \|_{N[0]} \\ &\les   \sum_{\ell\le C}
2^{\frac{\ell}{4}}
\| \Box Q_\ell \psi_1\|_{\dot X_{k_1}^{0,-\frac12,\infty}} \|\tilde P_{k_1}Q_{\le C} [ |\nabla|^{-1}
\psi_2 |\nabla|^{-1}\psi_3] \|_{\dot X_{k_3}^{0,\frac12,1}} \\
&\les  2^{\frac{k_1-k_2}{4}} 2^{-\frac{k_2}{4}}  \prod_{i=1}^3
\|\psi_i\|_{S[k_i]}
\end{align*}
which is admissible.

\medskip
\noindent {\em Case 2:} $\mathit{0\le k_1= k_3+O(1), k_2\le k_3-C}$.
Using the arguments from Case~2 above, we may assume that
$\psi_1=Q_{\le (1-3\eps)k_2\wedge0-k_1}\psi_1$. In addition, it was
shown there that it suffices to assume that  $\psi_2=Q_{\le \eps
k_2\wedge k_2}\psi_2$, $\psi_3=Q_{\le \eps k_2}\psi_3$. First,
\begin{align*}
& \|P_0 I  (Q_{\le (1-3\eps)k_2 \wedge0-k_1 } \Box\psi_1 \;Q_{\le C}[ |\nabla|^{-1}\psi_2 |\nabla|^{-1} \psi_3])\|_{N[0]} \\
 &\les 2^{\frac{(1-3\eps)k_2 \wedge 0 -k_1}{4}}   \|\Box Q_{\le (1-3\eps)k_2\wedge0-k_1} \psi_1\|_{\dot X_{k_1}^{0,-\frac12,\infty}}
   \| Q_{\le C}[ |\nabla|^{-1}\psi_2 |\nabla|^{-1} \psi_3]\|_{\dot X_{k_1}^{0,\frac12,1}}  \\
&\les  2^{\frac{(1-3\eps)k_2 \wedge 0 -k_1}{4}} \la k_1-k_2\ra
\prod_{i=1}^3 \|\psi_i\|_{S[k_i]}
\end{align*}
which is admissible. One may also restrict the modulation of the entire output by~$Q_{\le (1-3\eps)k_2\wedge0-k_1}$. Therefore,
applying Lemma~\ref{lem:core} and Lemma~\ref{lem:Sk_prod2} to the
second expression in~\eqref{eq:nullexp} yields
\begin{align*}
 &\|P_0 Q_{\le (1-3\eps)k_2\wedge0-k_1}  [\Box Q_{\le (1-3\eps)k_2\wedge \eps k_2 }  (\psi_1
 |\nabla|^{-1} \psi_2) |\nabla|^{-1}\psi_3] \|_{N[0]} \\&\les     2^{\frac{(1-3\eps)k_2\wedge0-k_1}{4}}
   \| \tilde P_{k_1} Q_{\le (1-3\eps)k_2\wedge \eps k_2}  \Box (\psi_1  |\nabla|^{-1} \psi_2) \|_{\dot X_{k_1}^{0,-\frac12,1}}
\| |\nabla|^{-1}\psi_3\|_{S[k_3]} \\
&\les   2^{\frac{(1-3\eps)k_2\wedge0-k_1}{4}}
   \| Q_{\le (1-3\eps)k_2\wedge \eps k_2} (\psi_1  |\nabla|^{-1} \psi_2) \|_{\dot X_{k_1}^{0,\frac12,1}}
\|  \psi_3\|_{S[k_3]}    \\
&\les  2^{\frac{(1-3\eps)k_2\wedge0-k_1}{4}} \la k_1-k_2\ra    \prod_{i=1}^3 \|\psi_i\|_{S[k_i]}
\end{align*}
which is admissible.

\medskip
\noindent {\em Case 3:} $\mathit{0\le k_1= k_2+O(1), k_3\le k_2-C}$.  This
is symmetric to Case~2.

\medskip
\noindent {\em Case 4:} $\mathit{O(1)\le k_2= k_3+O(1), k_1\le -C
}$. This is similar to Case~1. Indeed,  the second terms
in~\eqref{eq:nullexp} and~\eqref{eq:nullexp2} satisfy the same
bounds, whence it will suffice to bound the first term
in~\eqref{eq:nullexp}. However, by \eqref{eq:weaker_core}  and
Lemma~\ref{lem:Sk_prod2},
\begin{align*}
& \| P_0 I [ (\Box I \psi_1) |\nabla|^{-1}  \psi_2
|\nabla|^{-1}\psi_3] \|_{N[0]} \\ &\les   \sum_{\ell\le k_1+ C}
2^{\frac{\ell-k_1}{4}}
\| \Box Q_\ell \psi_1\|_{\dot X_{k_1}^{0,-\frac12,\infty}} \|\tilde P_{0}Q_{\le C} [
|\nabla|^{-1} \psi_2 |\nabla|^{-1}\psi_3] \|_{\dot X_{k_2}^{0,\frac12,1}} \\
&\les   2^{k_1-k_2}    \prod_{i=1}^3 \|\psi_i\|_{S[k_i]}
\end{align*}
which is admissible.

\medskip
\noindent {\em Case 5:} $\mathit{O(1)= k_1,\; k_2=k_3+O(1) }$.
Here again it suffices to only consider the first term
in~\eqref{eq:nullexp}. Moreover, \eqref{eq:case5begin}
and~\eqref{eq:psi1_modred} apply whence that first term is bounded
by the Strichartz component~\eqref{eq:Sk2}:
\begin{align*}
& \sum_{k\le k_2\wedge 0 + C}   \|P_0 I (Q_{\le k}\Box \psi_1 P_k
I[|\nabla|^{-1}\psi_2 |\nabla|^{-1} \psi_3])\|_{N[0]} \\
& \les \sum_{k\le k_2\wedge 0 + C}   \|P_0 Q_{\le k+C}  (Q_{\le k}\Box \psi_1
P_k I[ |\nabla|^{-1}\psi_2 |\nabla|^{-1} \psi_3])\|_{\enerN} \\ &
\les \sum_{k\le k_2\wedge 0 + C}     \|Q_{\le k}\Box \psi_1\|_{\Ltwotx} \| P_k
I[ |\nabla|^{-1}\psi_2 |\nabla|^{-1} \psi_3] \|_{L^2_t L^\infty_x} \\
&
\les \sum_{k\le k_2\wedge 0 + C}   2^{k-\frac{k_2}{2}} \prod_{i=1}^3
\|\psi_i\|_{S[k_i]}
 \les 2^{-\frac{k_2\vee0}{2}} \prod_{i=1}^3 \|\psi_i\|_{S[k_i]}
\end{align*}

\medskip
\noindent {\em Case 6:} $\mathit{O(1)= k_1\ge k_2+O(1)\ge k_3+C}$.
Here one basically starts from~\eqref{eq:sch471}, which can be
handled via~\eqref{eq:Taotrilin}.

\medskip
\noindent {\em Case 7:} $\mathit{ k_1=O(1)\ge k_3+O(1)\ge k_2+C}$.
This case is symmetric to the previous one.

\medskip
\noindent {\em Case 8:} $\mathit{k_3=O(1), \max(k_1, k_2)\le -C}$.
As in Case~8 above, one first shows that one can place $Q_{\le
k_1+(1-3\eps)k_2}$ in front of the entire output, as well as in
front of~$\psi_1$. So it suffices to consider
\begin{align*}
 & P_0 Q_{\le k_1+(1-3\eps)k_2} [Q_{< k_1+ (1-3\eps)k_2}  \del^\beta \psi_1 \Delta^{-1}\del_j I \calN_{\beta j}(I
\psi_2, I\psi_3)] \\& =   P_0 Q_{\le k_1+(1-3\eps)k_2} [Q_{ < k_1+
(1-3\eps)k_2 } \del^\beta \psi_1 \Delta^{-1}\del_j Q_{\le k_1+ C}
\tilde P_0 \calN_{\beta j}(I \psi_2, I\psi_3)] \nn
\end{align*}
We now stimate the first two terms on the right-hand side
of~\eqref{eq:nullexp}. First, by the Strichartz
component~\eqref{eq:Sk2},
\begin{align*}
& \|P_0 I( Q_{\le k_1+(1-3\eps)k_2} \Box \psi_1 |\nabla|^{-1}\psi_2 |\nabla|^{-1} \psi_3)\|_{N[0]} \\
&\les \|P_0 I ( Q_{\le k_1+(1-3\eps)k_2} \Box \psi_1 |\nabla|^{-1}\psi_2 |\nabla|^{-1} \psi_3)\|_{\enerN} \\
&\les  2^{ (1-3\eps)k_2+k_1}  \|\psi_1 \|_{\ener} \| |\nabla|^{-1}\psi_2 |\nabla|^{-1} \psi_3\|_{L^2_t L^\infty_x} \\
&\les 2^{(1-3\eps)k_2+k_1} \|\psi_1\|_{\ener}  \|\psi_2\|_{L^4_t L^\infty_x} 2^{-\frac{k_2}{4}}  \|\psi_3\|_{L^4_t L^\infty_x}  \\
&\les 2^{(1-3\eps)k_2+k_1}  2^{-\frac{k_2}{4}} \prod_{i=1}^3
\|\psi_i\|_{S[k_i]}
\end{align*}
which is sufficient. Second, by \eqref{eq:core} of
Lemma~\ref{lem:core} and Lemma~\ref{lem:Sk_prod2}, and assuming
first that $k_1=k_2+O(1)$,
\begin{align*}
 &\|P_0 Q_{\le k_1+(1-3\eps)k_2} [\Box  ( Q_{\le k_1+(1-3\eps)k_2} \psi_1  |\nabla|^{-1} \psi_2) |\nabla|^{-1}\psi_3] \|_{N[0]} \\
 &\les \sum_{k\le  k_1 +C} \|P_0 [\Box  Q_{\le k_1+C} P_k ( Q_{\le k_1+(1-3\eps)k_2} \psi_1  |\nabla|^{-1} \psi_2) |\nabla|^{-1}\psi_3] \|_{N[0]} \\
 &\les    \sum_{k\le  k_1+C}  \|   P_{k} Q_{\le k_1+C}  \Box (\psi_1  |\nabla|^{-1} \psi_2) \|_{\dot X_{k}^{0,-\frac12,1}}
 \|\, |\nabla|^{-1}\psi_3\|_{S[k_3]} \\
&\les  \sum_{k\le  k_1+C}  2^k    \|  P_k  Q_{\le  k_1+C }  (\psi_1
|\nabla|^{-1} \psi_2) \|_{\dot X_{ k}^{0,\frac12,1}} \|
\psi_3\|_{S[k_3]}\\
& \les \sum_{k\le  k_1+C}  2^k 2^{\frac{k-k_1}{4}} \prod_{i=1}^3
\|\psi_i\|_{S[k_i]}
\end{align*}
which is admissible. If  $k_2<k_1-C$, then by the same lemmas,
\begin{align*}
 &\|P_0 Q_{\le k_1+(1-3\eps)k_2} [\Box  ( Q_{\le k_1+(1-3\eps)k_2} \psi_1  |\nabla|^{-1} \psi_2) |\nabla|^{-1}\psi_3] \|_{N[0]} \\
 &\les   \|P_0 [\Box  Q_{\le (1-3\eps)k_2+C} \tilde P_{k_1} ( Q_{\le k_1+(1-3\eps)k_2} \psi_1  |\nabla|^{-1} \psi_2) |\nabla|^{-1}\psi_3] \|_{N[0]} \\
 &\les       \| \tilde P_{k_1} Q_{\le (1-3\eps)k_2+C}  \Box (\psi_1  |\nabla|^{-1} \psi_2) \|_{\dot X_{k_1}^{0,-\frac12,1}}
 \|\, |\nabla|^{-1}\psi_3\|_{S[k_3]} \\
&\les    2^{k_1}     \|  \tilde P_{k_1}  Q_{\le  (1-3\eps)k_2+C }  (\psi_1  |\nabla|^{-1} \psi_2) \|_{\dot X_{ k_1}^{0,\frac12,1}} \|  \psi_3\|_{S[k_3]}\\
&\les   2^{k_1} 2^{\frac{ (1-3\eps)k_2-k_1}{4}}    \prod_{i=1}^3
\|\psi_i\|_{S[k_i]}
\end{align*}
which is again admissible. Finally, if  $k_1<k_2-C$, then arguing
analogously yields
\begin{align*}
 &\|P_0 Q_{\le k_1+(1-3\eps)k_2} [\Box  ( Q_{\le k_1+(1-3\eps)k_2} \psi_1  |\nabla|^{-1} \psi_2) |\nabla|^{-1}\psi_3] \|_{N[0]} \\
 &\les   \|P_0 [\Box  Q_{\le  k_2} \tilde P_{k_2} ( Q_{\le k_1+(1-3\eps)k_2} \psi_1  |\nabla|^{-1} \psi_2) |\nabla|^{-1}\psi_3] \|_{N[0]} \\
 &\les       \| \tilde P_{k_2} Q_{\le k_2}  \Box (\psi_1  |\nabla|^{-1} \psi_2) \|_{\dot X_{k_2}^{0,-\frac12,1}}
 \|\, |\nabla|^{-1}\psi_3\|_{S[k_3]} \\
&\les    2^{k_2}     \|  \tilde P_{k_2}  Q_{\le   k_2}  (\psi_1  |\nabla|^{-1} \psi_2) \|_{\dot X_{ k_2}^{0,\frac12,1}} \|  \psi_3\|_{S[k_3]}\\
&\les   2^{k_1}    \prod_{i=1}^3 \|\psi_i\|_{S[k_i]}
\end{align*}
which concludes this case.

\medskip
\noindent {\em Case 9:} $\mathit{k_2=O(1), \max(k_1, k_3)\le -C}$.
Symmetric to Case~8.
This concludes the analysis of~\eqref{eq:keytri2}.

\medskip
Neither of the identities \eqref{eq:nullexp} or \eqref{eq:nullexp2} applies to~\eqref{eq:keytri3}. Hence,
\eqref{eq:keytri3} requires somewhat different arguments.

\noindent
{\em Case 1:}  $\mathit{0\le k_1\le k_2+O(1)=k_3+O(1)}$.
As in~\eqref{eq:case1_red} one sees that it suffices to consider $\psi_1=Q_{\le 0}\psi_1$.
Then $\calN_{\alpha j}=Q_{\le C}\calN_{\alpha j}$ and we split
\begin{align}
 &  P_0 I \del^\beta [I R_\beta \psi_1 \Delta^{-1}\del_j I \calN_{\alpha j}(I
\psi_2, I\psi_3)] \nn \\
 & = \sum_{\ell\le C} P_0 Q_{\le \ell-C}  \del^\beta [ R_\beta Q_{\le \ell-C} \psi_1 \Delta^{-1}\del_j Q_\ell \calN_{\alpha j}(I
\psi_2, I\psi_3)] \label{eq:Q01} \\
& + \sum_{\ell-C\le\ell_1\le  C} P_0 Q_{\ell_1}  \del^\beta [ R_\beta Q_{\le \ell_1} \psi_1 \Delta^{-1}\del_j Q_\ell  \calN_{\alpha j}(I
\psi_2, I\psi_3)] \label{eq:Q02}\\
&  + \sum_{\ell-C\le\ell_2\le  C} P_0 Q_{< \ell_2}  \del^\beta [ R_\beta Q_{\ell_2} \psi_1 \Delta^{-1}\del_j Q_\ell  \calN_{\alpha j}(I
\psi_2, I\psi_3)] \label{eq:Q03}
\end{align}
Decomposing~\eqref{eq:Q01} via Lemma~\ref{lem:cone} into caps of size~$2^{\frac{\ell}{2}}$ yields
\begin{align*}
 & P_0 Q_{\le \ell-C}  \del^\beta [ R_\beta Q_{\le \ell-C} \psi_1 \Delta^{-1}\del_j Q_\ell \calN_{\alpha j}(I
\psi_2, I\psi_3)] \\
&= \sum_{\kappa\sim \kappa'\in\calC_{\frac{\ell}{2}}} P_{0,\kappa} Q_{\le \ell-C}  \del^\beta
[ R_\beta Q_{\le \ell-C} P_{k_1,\kappa'}  \psi_1 \Delta^{-1}\del_j Q_\ell \calN_{\alpha j}(I
\psi_2, I\psi_3)]
\end{align*}
where $\kappa\sim\kappa'$ denotes that these caps have distance about~$2^{\frac{\ell}{2}}$. Hence we
gain a factor of~$2^\ell$ from the nullform involving $\del^\beta$ and~$R_\beta$. From~\eqref{eq:bilin1} one now obtains
\begin{align*}
\|\eqref{eq:Q01}\|_{N[0]} &\les \sum_{\ell\le C}\Big( \sum_{\kappa\sim\kappa'\in\calC_{\frac{\ell}{2}}}
\|P_{0,\kappa} Q_{\le \ell-C}  \del^\beta [ R_\beta Q_{\le \ell-C} P_{k_1,\kappa'}  \psi_1 \Delta^{-1}\del_j Q_\ell \calN_{\alpha j}(I
\psi_2, I\psi_3)]\|_{\NF[\kappa]}^2\Big)^{\frac12}  \\
&\les  \sum_{\ell\le C} 2^\ell 2^{-\frac{\ell}{4}}  2^{\frac{k_1}{2}} \Big( \sum_{\kappa\in\calC_{\frac{\ell}{2}}}
\|P_{k_1,\kappa} Q_{\le \ell-C} \psi_1\|_{S[\kappa]}^2  2^{-2k_1} \|Q_\ell \calN_{\alpha j} (\psi_2,\psi_3)\|_{\Ltwotx}^2\Big)^{\frac12} \\
&\les \sum_{\ell\le C} 2^{\frac{3\ell}{4}} 2^{\frac{k_1}{2}} \|\psi_1\|_{S[k_1]}  2^{-\frac{k_1}{4+}} 2^{-\frac{k_2}{2}}
\| \psi_2\|_{S[k_2]} \|\psi_3\|_{S[k_3]} \\
&\les 2^{\frac{k_1}{4-}-\frac{k_2}{2}} \prod_{i=1}^3 \|\psi_i\|_{S[k_i]}
\end{align*}
Here we also used Lemma~\ref{lem:square_func} as well as Lemma~\ref{lem:Nablowmod}.
The expressions in~\eqref{eq:Q02} are decomposed into caps of size~$2^{\frac{\ell_1}{2}}$ but without separation.
Therefore,  with a gain of~$2^{\ell_1}$ from the outer null-form,
\begin{align}
\|\eqref{eq:Q02}\|_{N[0]} &\les \sum_{\ell-C\le \ell_1\le C} \Big( \sum_{\kappa,\kappa'\in\calC_{\frac{\ell_1}{2}}}
 \|P_{0,\kappa} Q_{\ell_1}  \del^\beta [ R_\beta Q_{\le \ell_1} P_{k_1,\kappa'}  \psi_1 \Delta^{-1}\del_j Q_\ell \calN_{\alpha j}(I
\psi_2, I\psi_3)]\|_{\dot X_0^{0,-\frac12,1}}^2\Big)^{\frac12}  \nn \\
&\les  \sum_{\ell-C\le \ell_1\le C} 2^{-\frac{\ell_1}{2}} \Big( \sum_{\kappa,\kappa'\in\calC_{\frac{\ell_1}{2}}}
\|P_{0,\kappa} Q_{\ell_1}  \del^\beta [ R_\beta Q_{\le \ell_1} P_{k_1,\kappa'}  \psi_1 \Delta^{-1}\del_j Q_\ell \calN_{\alpha j}(I
\psi_2, I\psi_3)]\|_{\Ltwotx}^2\Big)^{\frac12}  \nn \\
&\les  \sum_{\ell-C\le \ell_1\le C}   \Big( \sum_{\kappa,\kappa'\in\calC_{\frac{\ell_1}{2}}}
\|P_{0,\kappa} Q_{\ell_1}  \del^\beta [ R_\beta Q_{\le \ell_1} P_{k_1,\kappa'}  \psi_1 \Delta^{-1}\del_j Q_\ell \calN_{\alpha j}(I
\psi_2, I\psi_3)]\|_{L^2_t L^1_x}^2\Big)^{\frac12} \label{eq:Bern23}
\end{align}
To pass to~\eqref{eq:Bern23} one invokes the improved Bernstein estimate of Lemma~\ref{lem:KBern}.
 Hence, this can be further bounded
by
\begin{align}
&\les  \sum_{\ell-C\le \ell_1\le C}  2^{\ell_1}  \Big( \sum_{\kappa'\in\calC_{\frac{\ell_1}{2}}}
\| Q_{\le \ell_1} P_{k_1,\kappa'}  \psi_1\|_{\ener}^2 \, 2^{-2k_1} \|  Q_\ell \calN_{\alpha j}(I
\psi_2, I\psi_3)]\|_{L^2_t L^2_x}^2\Big)^{\frac12}  \nn \\
&\les  \sum_{\ell-C\le \ell_1\le C}  2^{\ell_1}  \Big( \sum_{\kappa'\in\calC_{\frac{\ell_1}{2}}}
\| Q_{\le \ell_1} P_{k_1,\kappa'}  \psi_1\|_{S[k_1,\kappa']}^2\Big)^{\frac12}  \, 2^{-k_1}
2^{\frac{\ell-k_1}{4+}} 2^{k_1-\frac{k_2}{2}} \|  Q_\ell \calN_{\alpha j}(I
\psi_2, I\psi_3)]\|_{L^2_t L^2_x}  \nn \\
&\les 2^{-\frac{k_1}{4-}-\frac{k_2}{2}} \prod_{i=1}^3 \|\psi_i\|_{S[k_i]}  \nn
\end{align}
For~\eqref{eq:Q03} one proceeds similarly, performing a cap
decomposition and placing the entire expression in~$\enerN$. We skip
the details.

\medskip
\noindent {\em Case 2:} $\mathit{0\le k_1= k_3+O(1), k_2\le k_3-C}$.  This is essentially the same
as the preceding with Lemma~\ref{lem:Nablowmod2} replacing Lemma~\ref{lem:Nablowmod}.

\medskip
\noindent {\em Case 3:} $\mathit{0\le k_1= k_2+O(1), k_3\le k_2-C}$.  This
is symmetric to the preceding.

\medskip
\noindent {\em Case 4:} $\mathit{O(1)\le k_2= k_3+O(1), k_1\le -C }$.  This is very similar to Case~1.
First, one checks that the entire output can be restricted by~$Q_{\le k_1}$.
This implies that  due to the $I$-operator in front of~$\psi_1$, the decomposition
\eqref{eq:Q01}--\eqref{eq:Q03} continues to hold but with $\ell\le k_1+C$:
\begin{align}
 &  P_0 I \del^\beta [I R_\beta \psi_1 \Delta^{-1}\del_j I \calN_{\alpha j}(I
\psi_2, I\psi_3)] \nn \\
 & = \sum_{\ell\le k_1+C} P_0 Q_{\le \ell-C}  \del^\beta [ R_\beta Q_{\le \ell-C} \psi_1 \Delta^{-1}\del_j Q_\ell \calN_{\alpha j}(I
\psi_2, I\psi_3)] \label{eq:Q01'} \\
& + \sum_{\ell-C\le\ell_1\le k_1+ C} P_0 Q_{\ell_1}  \del^\beta [ R_\beta Q_{\le \ell_1} \psi_1 \Delta^{-1}\del_j Q_\ell  \calN_{\alpha j}(I
\psi_2, I\psi_3)] \label{eq:Q02'}\\
&  + \sum_{\ell-C\le\ell_2\le k_1+ C} P_0 Q_{< \ell_2}  \del^\beta [ R_\beta Q_{\ell_2} \psi_1 \Delta^{-1}\del_j Q_\ell  \calN_{\alpha j}(I
\psi_2, I\psi_3)] \label{eq:Q03'}
\end{align}
One can again decompose~\eqref{eq:Q01'}
into caps, but of size~$2^{\frac{\ell-k_1}{2}}$. Therefore,
\begin{align*}
\|\eqref{eq:Q01'}\|_{N[0]} &\les \sum_{\ell\le k_1+C}\Big( \sum_{\kappa\sim\kappa'\in\calC_{\frac{\ell-k_1}{2}}}
 \|P_{0,\kappa} Q_{\le \ell-C}  \del^\beta [ R_\beta Q_{\le \ell-C} P_{k_1,\kappa'}  \psi_1 \Delta^{-1}\del_j Q_\ell \tilde P_0 \calN_{\alpha j}(I
\psi_2, I\psi_3)]\|_{\NF[\kappa]}^2\Big)^{\frac12}  \\
&\les  \sum_{\ell\le  k_1+ C} 2^{\frac{3(\ell-k_1)}{4}}  2^{\frac{k_1}{2}} \Big( \sum_{\kappa\in\calC_{\frac{\ell}{2}}}
\|P_{k_1,\kappa} Q_{\le \ell-C} \psi_1\|_{S[\kappa]}^2    \|Q_\ell \tilde P_0 \calN_{\alpha j} (\psi_2,\psi_3)\|_{\Ltwotx}^2\Big)^{\frac12} \\
&\les \sum_{\ell\le  k_1+ C} 2^{\frac{3(\ell-k_1)}{4}} 2^{\frac{k_1}{2}} \|\psi_1\|_{S[k_1]}    2^{-\frac{k_2}{2}}
\| \psi_2\|_{S[k_2]} \|\psi_3\|_{S[k_3]} \\
&\les 2^{\frac{k_1}{2}-\frac{k_2}{2}} \prod_{i=1}^3 \|\psi_i\|_{S[k_i]}
\end{align*}
which is admissible.  Furthermore, $\|\eqref{eq:Q02'}\|_{N[0]}$ is bounded by
\begin{align}
 &\les \sum_{\ell-C\le \ell_1\le k_1+ C} \Big( \sum_{\kappa,\kappa'\in\calC_{\frac{\ell_1-k_1}{2}}}
   \|P_{0,\kappa} Q_{\ell_1}  \del^\beta [ R_\beta Q_{\le \ell_1} P_{k_1,\kappa'}  \psi_1 \Delta^{-1}\del_j Q_\ell \calN_{\alpha j}(I
\psi_2, I\psi_3)]\|_{\dot X_0^{0,-\frac12,1}}^2\Big)^{\frac12}  \nn \\
&\les  \sum_{\ell-C\le \ell_1\le k_1+ C} 2^{-\frac{\ell_1 }{2}} \Big( \sum_{\kappa,\kappa'\in\calC_{\frac{\ell_1-k_1}{2}}}
 \|P_{0,\kappa} Q_{\ell_1}  \del^\beta [ R_\beta Q_{\le \ell_1} P_{k_1,\kappa'}  \psi_1 \Delta^{-1}\del_j Q_\ell \calN_{\alpha j}(I
\psi_2, I\psi_3)]\|_{\Ltwotx}^2\Big)^{\frac12}  \nn \\
&\les  \sum_{\ell-C\le \ell_1\le k_1+ C}  2^{\frac{\ell_1}{2}-k_1}   \Big( \sum_{\kappa'\in\calC_{\frac{\ell_1-k_1}{2}}}
\| Q_{\le \ell_1} P_{k_1,\kappa'}  \psi_1\|_{\Linf}^2 \,  \|  Q_\ell \tilde P_0 \calN_{\alpha j}(I
\psi_2, I\psi_3)]\|_{L^2_t L^2_x}^2\Big)^{\frac12}  \nn \\
&\les  \sum_{\ell-C\le \ell_1\le k_1+ C}   2^{\frac{\ell_1}{2} }  2^{\frac{\ell_1-k_1}{4}}
 \Big( \sum_{\kappa'\in\calC_{\frac{\ell_1-k_1}{2}}}  \| Q_{\le \ell_1} P_{k_1,\kappa'}
   \psi_1\|_{S[k_1,\kappa']}^2\Big)^{\frac12}  \,   2^{\frac{\ell}{4+}-\frac{k_2}{2}} \|  Q_\ell\tilde P_0 \calN_{\alpha j}(I
\psi_2, I\psi_3)]\|_{L^2_t L^2_x}  \nn \\
&\les 2^{\frac{3k_1}{4+}-\frac{k_2}{2}} \prod_{i=1}^3 \|\psi_i\|_{S[k_i]}  \nn
\end{align}
Finally, \eqref{eq:Q03'} is similar to the previous estimate and we skip it.

\medskip
\noindent {\em Case 5:} $\mathit{O(1)= k_1,\; k_2=k_3+O(1) }$.
We apply~\eqref{eq:case5begin} and reduce the modulation of~$\psi_1$ via~\eqref{eq:psi1_modred}
to~$\psi_1=Q_{\le k}\psi_1$.
Furthermore,
\begin{align}
 P_0\del^\beta I [R_\beta Q_{\le k} \psi_1 \, \Delta^{-1}\del_j P_k I\calN_{\alpha j}(\psi_2,\psi_3)]
&=  P_0 I [\Box|\nabla|^{-1}  Q_{\le k} \psi_1 \, \Delta^{-1}\del_j P_k I\calN_{\alpha j}(\psi_2,\psi_3)] \label{eq:Case5term1} \\
& + P_0 I [R_\beta Q_{\le k} \psi_1 \, \Delta^{-1}\del_j\del^\beta P_k I \calN_{\alpha j}(\psi_2,\psi_3)] \label{eq:Case5term2}
\end{align}
Lemmas~\ref{lem:core} and~\ref{lem:Nablowmod} imply the following bound on~\eqref{eq:Case5term1}:
\begin{align*}
 &\sum_{k\le k_2\wedge 0 + C}   \| P_0 I [\Box|\nabla|^{-1}  Q_{\le k} \psi_1 \, \Delta^{-1}\del_j P_k I\calN_{\alpha j}(\psi_2,\psi_3)] \|_{N[0]} \\
&\les \sum_{k\le k_2\wedge 0 + C}   \sum_{m\le k} 2^{\frac{m-k}{4}} \|\Box|\nabla|^{-1}  Q_{m} \psi_1\|_{\dot X_0^{0,-\frac12,1}}
\|  \Delta^{-1}\del_j P_k I\calN_{\alpha j}(\psi_2,\psi_3)] \|_{\dot X_k^{0,\frac12,1}} \\
&\les \sum_{k\le k_2\wedge 0 + C}    \, 2^{\frac{k-k_2}{2}} \prod_{i=1}^3 \|\psi_i\|_{S[k_i]} \les
2^{-\frac12 k_2\vee0} \prod_{i=1}^3 \|\psi_i\|_{S[k_i]}
\end{align*}
which is admissible. The second term \eqref{eq:Case5term2} needs to be expanded as follows:
\begin{align}
   2 R_\beta Q_{\le k} \psi_1 \, \Delta^{-1}\del_j\del^\beta P_k I \calN_{\alpha j}(\psi_2,\psi_3)
&= \Box [Q_{\le k} |\nabla|^{-1}\psi_1 \, \Delta^{-1}\del_j  P_k I \calN_{\alpha j}(\psi_2,\psi_3) ] \label{eq:sch97}\\
& - \Box  Q_{\le k} |\nabla|^{-1}\psi_1 \, \Delta^{-1}\del_j  P_k I \calN_{\alpha j}(\psi_2,\psi_3) \label{eq:sch98} \\
& -  Q_{\le k} |\nabla|^{-1} \psi_1 \, \Box  \Delta^{-1}\del_j  P_k I \calN_{\alpha j}(\psi_2,\psi_3) \label{eq:sch99}
\end{align}
We just dealt with the term~\eqref{eq:sch98}. Since the modulation
of the entire output is~$\les 2^k$, one concludes that
\begin{align*}
\eqref{eq:sch97} & \les \sum_{k\le k_2\wedge 0 + C}   \| \Box [Q_{\le k} |\nabla|^{-1}\psi_1 \, \Delta^{-1}\del_j
P_k I \calN_{\alpha j}(\psi_2,\psi_3) ]\|_{\dot X_0^{0,-\frac12,1}} \\
&\les \sum_{k\le k_2\wedge 0 + C}   2^{\frac{k}{2}} \|\psi_1\|_{\ener}  \|\Delta^{-1}\del_j  P_k I \calN_{\alpha j}(\psi_2,\psi_3)\|_{\Ltwotx}\\
&\les \sum_{k\le k_2\wedge 0 + C}   2^{\frac{k-k_2}{2}} \prod_{i=1}^3 \|\psi_i\|_{S[k_i]}  \les 2^{-\frac12 k_2\vee0}  \prod_{i=1}^3 \|\psi_i\|_{S[k_i]}
\end{align*}
as well as, from Lemma~\ref{lem:core},
\begin{align*}
\eqref{eq:sch99} & \les \sum_{k\le k_2\wedge 0 + C}   \|  Q_{\le k} |\nabla|^{-1}\psi_1 \, \Box \Delta^{-1}\del_j
 P_k I \calN_{\alpha j}(\psi_2,\psi_3) ]\|_{N[0]} \\
&\les \sum_{k\le k_2\wedge 0 + C}    \|\psi_1\|_{\ener}  \|\Box \Delta^{-1}\del_j  P_k I \calN_{\alpha j}(\psi_2,\psi_3)\|_{\dot X_0^{0,-\frac12,1} }\\
&\les \sum_{k\le k_2\wedge 0 + C}   2^{\frac{3k-k_2}{2}} \prod_{i=1}^3 \|\psi_i\|_{S[k_i]}  \les 2^{-\frac12 k_2\vee0}  \prod_{i=1}^3 \|\psi_i\|_{S[k_i]}
\end{align*}
which is sufficient.

\medskip
\noindent {\em Case 6:} $\mathit{O(1)= k_1\ge k_2+O(1)\ge k_3+C}$.
As before, one reduces the modulation of~$\psi_1$ to~$\psi_1=Q_{\le k_2}\psi_1$.
Furthermore,
\begin{align}
 P_0\del^\beta I [R_\beta Q_{\le k_2} \psi_1 \, \Delta^{-1}\del_j \tilde P_{k_2} I\calN_{\alpha j}(\psi_2,\psi_3)]
&=  P_0 I [\Box|\nabla|^{-1}  Q_{\le k_2} \psi_1 \, \Delta^{-1}\del_j \tilde P_{k_2} I\calN_{\alpha j}(\psi_2,\psi_3)] \label{eq:Case6term1} \\
& + P_0 I [R_\beta Q_{\le k_2} \psi_1 \, \Delta^{-1}\del_j\del^\beta \tilde P_{k_2} I \calN_{\alpha j}(\psi_2,\psi_3)] \label{eq:Case6term2}
\end{align}
Lemmas~\ref{lem:core} and~\ref{lem:Nablowmod2} imply the following bound on~\eqref{eq:Case6term1}:
\begin{align*}
 &  \| P_0 I [\Box|\nabla|^{-1}  Q_{\le k_2} \psi_1 \, \Delta^{-1}\del_j \tilde P_{k_2} I\calN_{\alpha j}(\psi_2,\psi_3)] \|_{N[0]} \\
&\les   \sum_{m\le k_2} 2^{\frac{m-k_2}{4}} \|\Box|\nabla|^{-1}  Q_{m} \psi_1\|_{\dot X_0^{0,-\frac12,1}}
\|  \Delta^{-1}\del_j \tilde P_{k_2} I\calN_{\alpha j}(\psi_2,\psi_3)] \|_{\dot X_{k_2}^{0,\frac12,1}} \\
&\les   2^{(\frac12-\eps)(k_3-k_2)}     \prod_{i=1}^3 \|\psi_i\|_{S[k_i]}
\end{align*}
which is admissible. The second term \eqref{eq:Case6term2} needs to be expanded as follows:
\begin{align}
   2 R_\beta Q_{\le k_2} \psi_1 \, \Delta^{-1}\del_j\del^\beta \tilde P_{k_2} I \calN_{\alpha j}(\psi_2,\psi_3)
&= \Box [Q_{\le k_2} |\nabla|^{-1}\psi_1 \, \Delta^{-1}\del_j \tilde P_{k_2} I \calN_{\alpha j}(\psi_2,\psi_3) ] \label{eq:sch97'}\\
& - \Box  Q_{\le k_2} |\nabla|^{-1}\psi_1 \, \Delta^{-1}\del_j  \tilde P_{k_2} I \calN_{\alpha j}(\psi_2,\psi_3) \label{eq:sch98'} \\
& -  Q_{\le k_2} |\nabla|^{-1} \psi_1 \, \Box  \Delta^{-1}\del_j  \tilde P_{k_2} I \calN_{\alpha j}(\psi_2,\psi_3) \label{eq:sch99'}
\end{align}
We just dealt with the term~\eqref{eq:sch98'}. Since the modulation
of the entire output is~$\les 2^{k_2}$, one concludes that
\begin{align*}
\eqref{eq:sch97'} & \les  \| \Box [Q_{\le k_2} |\nabla|^{-1}\psi_1 \, \Delta^{-1}\del_j
 \tilde P_{k_2} I \calN_{\alpha j}(\psi_2,\psi_3) ]\|_{\dot X_0^{0,-\frac12,1}} \\
&\les   2^{\frac{k_2}{2}} \|\psi_1\|_{\ener}  \|\Delta^{-1}\del_j  \tilde P_{k_2} I \calN_{\alpha j}(\psi_2,\psi_3)\|_{\Ltwotx}\\
&\les   2^{(\frac12-\eps)(k_3-k_2)} \prod_{i=1}^3 \|\psi_i\|_{S[k_i]}  \les \prod_{i=1}^3 \|\psi_i\|_{S[k_i]}
\end{align*}
as well as, from Lemma~\ref{lem:core},
\begin{align*}
\eqref{eq:sch99'} & \les   \|  Q_{\le k} |\nabla|^{-1}\psi_1 \, \Box \Delta^{-1}\del_j  \tilde P_{k_2} I \calN_{\alpha j}(\psi_2,\psi_3) ]\|_{N[0]} \\
&\les    \|\psi_1\|_{S[k_1]}  \|\Box \Delta^{-1}\del_j  \tilde P_{k_2} I \calN_{\alpha j}(\psi_2,\psi_3)\|_{\dot X_0^{0,-\frac12,1} }\\
&\les   2^{(\frac12-\eps)(k_3-k_2)} 2^{k_2}  \prod_{i=1}^3 \|\psi_i\|_{S[k_i]}
\end{align*}
which concludes Case~6.

\medskip
\noindent {\em Case 7:} $\mathit{ k_1=O(1)\ge k_3+O(1)\ge k_2+C}$.
This case is symmetric to the previous one.

\medskip
\noindent {\em Case 8:} $\mathit{k_3=O(1), \max(k_1, k_2)\le -C}$.
The modulation of the output can be reduced to $Q_{\le  k_1}$:
\begin{align*}
 & \| P_0\del^\beta I Q_{\ge   k_1}  [R_\beta Q_{\le k_1} \psi_1 \, \Delta^{-1}\del_j \tilde P_{0} I\calN_{\alpha j}(\psi_2,\psi_3)]   \|_{N[0]}  \\
&\les   2^{-\frac{ k_1}{2}}  \|\psi_1\|_{\Linf} \| \Delta^{-1}\del_j \tilde P_{0} I\calN_{\alpha j}(\psi_2,\psi_3)] \|_{\Ltwotx} \\
&\les 2^{\frac{k_1}{2}}  2^{(\frac12-\eps)k_3} \prod_{i=1}^3 \|\psi_i\|_{S[k_i]}
\end{align*}
Similarly, the input $\psi_1$ can be reduced to $Q_{\le k_1}\psi_1$.
As in Case~6,
\begin{align}
& P_0\del^\beta Q_{\le   k_1} [R_\beta Q_{\le   k_1} \psi_1 \, \Delta^{-1}\del_j \tilde P_{0} Q_{\le k_1}  \calN_{\alpha j}(\psi_2,\psi_3)]\nn \\
&=  P_0  Q_{\le   k_1} [\Box|\nabla|^{-1}  Q_{\le   k_1} \psi_1 \, \Delta^{-1}\del_j \tilde P_{0} Q_{\le   k_1} \calN_{\alpha j}(\psi_2,\psi_3)] \label{eq:Case8term1} \\
& + P_0 Q_{\le   k_1} [R_\beta Q_{\le   k_1} \psi_1 \, \Delta^{-1}\del_j\del^\beta \tilde P_{0}  Q_{\le   k_1} \calN_{\alpha j}(\psi_2,\psi_3)] \label{eq:Case8term2}
\end{align}
Lemmas~\ref{lem:core} and~\ref{lem:Nablowmod2} imply the following bound on~\eqref{eq:Case8term1}:
\begin{align*}
 &  \| P_0 Q_{\le  k_1}[\Box|\nabla|^{-1}  Q_{\le  k_1} \psi_1 \, \Delta^{-1}\del_j \tilde P_{0} I\calN_{\alpha j}(\psi_2,\psi_3)] \|_{N[0]} \\
&\les   \sum_{m\le  k_1} 2^{\frac{m-k_1}{4}} \|\Box|\nabla|^{-1}  Q_{m} \psi_1\|_{\dot X_{k_1}^{0,-\frac12,1}} \|    \tilde P_{0} Q_{\le   k_1} \calN_{\alpha j}(\psi_2,\psi_3)] \|_{\dot X_0^{0,\frac12,1}} \\
&\les  2^{\frac{k_1}{2}}  2^{(\frac12-\eps)k_2}      \prod_{i=1}^3 \|\psi_i\|_{S[k_i]}
\end{align*}
which is admissible. The second term \eqref{eq:Case8term2} needs to be expanded as follows:
\begin{align}
   2 R_\beta Q_{\le k_1} \psi_1 \, \Delta^{-1}\del_j\del^\beta \tilde P_{0} I \calN_{\alpha j}(\psi_2,\psi_3)
&= \Box [Q_{\le k_1} |\nabla|^{-1}\psi_1 \, \Delta^{-1}\del_j \tilde P_{0} I \calN_{\alpha j}(\psi_2,\psi_3) ] \label{eq:sch97''}\\
& - \Box  Q_{\le k_1} |\nabla|^{-1}\psi_1 \, \Delta^{-1}\del_j  \tilde P_{0} I \calN_{\alpha j}(\psi_2,\psi_3) \label{eq:sch98''} \\
& -  Q_{\le k_1} |\nabla|^{-1} \psi_1 \, \Box  \Delta^{-1}\del_j  \tilde P_{0} I \calN_{\alpha j}(\psi_2,\psi_3) \label{eq:sch99''}
\end{align}
We just dealt with the term~\eqref{eq:sch98''}. Next,
\begin{align*}
\eqref{eq:sch97''} & \les  \| \Box Q_{\le k_1} [Q_{\le k_1} |\nabla|^{-1}\psi_1 \, \Delta^{-1}\del_j  \tilde P_{0} I \calN_{\alpha j}(\psi_2,\psi_3) ]\|_{\dot X_0^{0,-\frac12,1}} \\
&\les   2^{\frac{k_1}{2}}   \||\nabla|^{-1} \psi_1\|_{\Linf}  \|\Delta^{-1}\del_j  \tilde P_{0} I \calN_{\alpha j}(\psi_2,\psi_3)\|_{\Ltwotx}\\
&\les  2^{\frac{k_1}{2}}  2^{(\frac12-\eps)k_3} \prod_{i=1}^3 \|\psi_i\|_{S[k_i]}
\end{align*}
as well as, from Lemma~\ref{lem:core},
\begin{align*}
\eqref{eq:sch99''} & \les   \|  Q_{\le k_1} |\nabla|^{-1}\psi_1 \, \Box \Delta^{-1}\del_j  \tilde P_{0} Q_{\le k_1} \calN_{\alpha j}(\psi_2,\psi_3) ]\|_{N[0]} \\
&\les      \|\psi_1\|_{S[k_1]}  \|\Box \Delta^{-1}\del_j  \tilde P_{0} Q_{\le k_1} \calN_{\alpha j}(\psi_2,\psi_3)\|_{\dot X_0^{0,-\frac12,1} }\\
&\les   2^{(\frac12-\eps)k_3} 2^{\frac{k_1}{2}}  \prod_{i=1}^3 \|\psi_i\|_{S[k_i]}
\end{align*}
which concludes Case~8.

\medskip
\noindent {\em Case 9:} $\mathit{k_2=O(1), \max(k_1, k_3)\le -C}$.
Symmetric to Case~8.
\end{proof}

\begin{remark}
 \label{rem:trilder} It follows from the high-low-low interaction case of the proof of  Lemma~\ref{lem:tri_hyp}
that for some $\sigma>0$,
\begin{equation}
 \label{eq:dergain}
\big\|  P_0 I [P_{k_1} I \psi_1\, \del^\beta \del_j \Delta^{-1} P_{k} I\calN_{\beta j}(P_{k_2}
 \psi_2, P_{k_3}\psi_3)\big \|_{N[0]} \les 2^{\sigma k}\, w(k_1,k_2,k_3)\prod_{i=1}^3 \|\psi_i\|_{S[k_i]}
\end{equation}
provided $k_1=O(1)$, $k\le k_2=k_3+O(1)\le O(1)$.
\end{remark}

\subsection{Improved trilinear estimates with angular alignment}\label{subsec:improvetrilin}

We conclude this section on trilinear bounds with a technical result
which we shall require in several instances, such as the blow-up
criterion of the following section.  By
Corollary~\ref{cor:epstrilin}, one gains extra smallness outside of
the parameter range~\eqref{eq:outside}; note that the latter
describes precisely Case~5 in the proof of
Lemmas~\ref{lem:hyp_redux} and~\ref{lem:tri_hyp} which is the
high-low-low case of interactions. In fact, the exact same gain as
in that corollary can also be obtained for the trilinear expressions
of Lemma~\ref{lem:tri_hyp}.

\begin{cor}
  \label{cor:epstrilinhyp} The nonlinearities of
  Lemma~\ref{lem:tri_hyp} satisfy the estimates of
  Corollary~\ref{cor:epstrilin}. I.e., given $\delta>0$ there exist
  $L,L'$ large so that the $\delta$--gains in the sum over
  $\sum_{k_1,k_2,k_3}'$ as well as  $\sum_{k_1,k_2,k_3}''$ with $k\le
  k_2-L'$, are obtained for the three types of trilinear null-forms in
  Lemma~\ref{lem:tri_hyp}.
\end{cor}
\begin{proof}
  As in the case of Corollary~\ref{cor:epstrilin}, this follows from
  the form of the weights $w(k_1,k_2,k_3)$ as well as from the fact
  that an extra gain in Case~5 of Lemma~\ref{lem:tri_hyp} was obtained when $k<k_2-L'$.
\end{proof}

However, one cannot gain smallness in the high-low-low case without
further assumptions. In this section we shall prove that {\em
angular alignment} between the Fourier support of at least two of
the inputs implies smallness in this case.

We start with the contributions by~$I^c \calN_{\beta j}$.  In
Corollary~\ref{cor:epstrilin} we isolated one case where smallness
cannot be obtained without any further assumptions. It was given by
the sum ${\sum}''_{k_1,k_2,k_3\in\Z}$ over the
range~\eqref{eq:outside} together with $k_2-L'\le k\le k_2+O(1)$.
Recall that $L$ and~$L'$ are very large depending on~$\delta$.
Throughout this section, $\psi_i$ will be Schwartz functions
satisfying
\[
\max_{i=1,2,3} \sum_{k\in\Z} \|P_k \psi_i\|_{S[k]}^2\le K^2
\]
for some constant~$K$. We shall use ${\sum}''_{k_1,k_2,k_3\in\Z}$
repeatedly in the sense that it was defined earlier.

\begin{lemma}
  \label{lem:trilindelta}
Given any $\delta>0$ there exists~$m_0(\delta)$ large and negative
such that
\begin{align*}
 & \sum_{\substack{\kappa_2,\kappa_3\in\calC_{m_0 }\\
\dist(\kappa_2,\kappa_3)\le 2^{m_0}}} {\sum}''_{k_1,k_2,k_3\in\Z}
\sum_{k=k_2-L'}^{k_2+O(1)} \, \sum_{j=1}^2 \| P_0 \nabla_{t,x}  [
P_{k_1}  \psi_1 \Delta^{-1}\del_j  I^c P_k \calN_{\beta
j}( P_{k_2,\kappa_2}  \psi_2, P_{k_3,\kappa_3} \psi_3) ] \|_{N[0]} \\
&
 \le \delta\,  K^2
   \sup_{k\in\Z}\max_{i=1,2,3} 2^{-\sigma_0 |k|} \|P_k \psi_i\|_{S[k]}
\end{align*}
as well as
\begin{align*}
& \sum_{\substack{\kappa_1,\kappa_2\in\calC_{m_0 }\\
\dist(\kappa_1,\kappa_2)\le 2^{m_0}}} {\sum}''_{k_1,k_2,k_3\in\Z}
\sum_{k=k_2-L'}^{k_2+O(1)}  \sum_{j=1}^2 \| P_0 \nabla_{t,x}   [
P_{k_1,\kappa_1}   \psi_1 \Delta^{-1}\del_j I^c  P_k \calN_{\beta
j}(  P_{k_2,\kappa_2} \psi_2, P_{k_3} \psi_3) ] \|_{N[0]} \\
& + \sum_{\substack{\kappa_1,\kappa_3\in\calC_{m_0 }\\
\dist(\kappa_1,\kappa_3)\le 2^{m_0}}} {\sum}''_{k_1,k_2,k_3\in\Z}
\sum_{k=k_2-L'}^{k_2+O(1)} \, \sum_{j=1}^2 \| P_0 \nabla_{t,x} [
P_{k_1,\kappa_1}   \psi_1 \Delta^{-1}\del_j I^c P_k \calN_{\beta
j}( P_{k_2} \psi_2, P_{k_3,\kappa_3} \psi_3) ] \|_{N[0]}\\
&
 \le \delta\,  K^2
   \sup_{k\in\Z}\max_{i=1,2,3} 2^{-\sigma_0 |k|} \|P_k \psi_i\|_{S[k]}
\end{align*}
\end{lemma}
\begin{proof} The proof simply consists in verifying that the
argument in  Case~5 of Lemma~\ref{lem:hyp_redux} allows for this
extra gain. We first consider angular alignment between $\psi_1$
and~$\psi_2$. In this case, we will need to repeat the argument of
Case~5, obtaining the gain from Bernstein's inequality. First,
restrict the output by $Q_{\ge0}$ and assume that
$\psi_1=P_{k_1,\kappa_1}\psi_1$ and~$\psi_2=P_{k_2,\kappa_2}\psi_2$
with fixed caps~$\kappa_1,\kappa_2$.
In the end, one verifies that it is possible to sum  over these caps.
Then
\begin{align}
 & \| P_0 Q_{\ge0}  \del^\beta [ \psi_1 \Delta^{-1}\del_j I^c \calN_{\beta j}( \psi_2, \psi_3)
] \|_{N[0]} \nn \\
&\les    \sum_{k=k_2-L'}^{k_2+O(1)}    \|P_0 Q_{\ge0}  [   \psi_1
\Delta^{-1}\del_j Q_{k\le \cdot\le C} P_k \calN_{\beta j}(  \psi_2,
\psi_3)
] \|_{\Ltwotx} \label{eq:sch13'}\\
&\quad+ \sum_{k=k_2-L'}^{k_2+O(1)}  \sum_{m\ge C}  2^{-\eps m} \|P_0
Q_m [    Q_{\le m-C}  \psi_1 \Delta^{-1}\del_j \tilde Q_m P_k
\calN_{\beta j}(  \psi_2, \psi_3)
] \|_{\Ltwotx} \label{eq:sch14'} \\
&\quad+ \sum_{k=k_2-L'}^{k_2+O(1)}  \sum_{m\ge C}  \|P_0 Q_{\ge0}  [
Q_{> m-C}  \psi_1 \Delta^{-1}\del_j \tilde Q_m P_k \calN_{\beta j}(
\psi_2, \psi_3) ] \|_{\Ltwotx}  \label{eq:sch15'}
\end{align}
First,  by Lemma~\ref{lem:Nabhighmod}, and with~$M$ large but finite
and $\frac{1}{p}+\frac{1}{M}=\frac12$,
\begin{align*}
 \eqref{eq:sch13'}
 &\les \sum_{k=k_2-L'}^{k_2+O(1)}   \|  \psi_1\|_{L^\infty_t L^p_x}\,   \|  Q_{k\le \cdot\le C} P_k \calN_{\beta j}(  \psi_2, \psi_3)
] \|_{L^2_t L^M_x} \\
&\les \sum_{k=k_2-L'}^{k_2+O(1)} 2^{m_0(\frac12-\frac{1}{p})} \|
\psi_1\|_{\ener}\, 2^{\frac{k}{2}} 2^{-\eps k_2}
2^{|k_2|\frac{2}{M}}
 \|\psi_2\|_{S[k_2]} \|\psi_3\|_{S[k_3]} \le \delta  \prod_{i=1}^3 \|\psi_i\|_{S[k_i]}
\end{align*}
Since $p>2$ one can take $m_0$ large and negative to obtain the
final estimate here.
 Second, again  by Lemma~\ref{lem:Nabhighmod},
\begin{align*}
 \eqref{eq:sch14'} &\les \sum_{k=k_2-L'}^{k_2+O(1)} \sum_{m\ge C}  2^{-\eps m} \|P_0 Q_m [    Q_{\le m-C}
  \psi_1 \Delta^{-1}\del_j \tilde Q_m P_k \calN_{\beta j}(  \psi_2, \psi_3)
] \|_{\Ltwotx} \\
&\les \sum_{k=k_2-L'}^{k_2+O(1)}  \sum_{m\ge C}  2^{-\eps m}
\|\psi_1\|_{L^\infty_t L^p_x} \; 2^{-k} \| \tilde Q_m P_k
\calN_{\beta j}( \psi_2, \psi_3)
] \|_{L^2_t L^M_x} \\
&\les \sum_{k=k_2-L'}^{k_2+O(1)}  2^{m_0(\frac12-\frac{1}{p})}
2^{\frac{k}{2}} 2^{-\eps k_2} 2^{|k_2|\frac{2}{M}} \prod_{i=1}^3
\|\psi_i\|_{S[k_i]} \le \delta \prod_{i=1}^3 \|\psi_i\|_{S[k_i]}
\end{align*}
and third,
\begin{align*}
\eqref{eq:sch15'} &\les \sum_{k=k_2-L'}^{k_2+O(1)}  \sum_{m\ge C} \|
Q_{> m-C}  \psi_1\|_{L^2_t L^p_x} 2^{-k} \| \tilde Q_m P_k
\calN_{\beta j}( \psi_2, \psi_3)
] \|_{L^\infty_t L^M_x} \\
&\les 2^{|k_2|\frac{2}{M}} \sum_{k=k_2-L'}^{k_2+O(1)}
2^{m_0(\frac12-\frac{1}{p})} \sum_{m\ge C}  2^{-(1-\eps)m} \|
\psi_1\|_{S[k_1]} \, 2^{\frac{m}{2}} \|   \tilde Q_m P_k
\calN_{\beta j}(  \psi_2, \psi_3)
] \|_{L^2_t L^M_x} \\
&\les 2^{m_0(\frac12-\frac{1}{p})} 2^{|k_2|\frac{2}{M}}
\sum_{k=k_2-L'}^{k_2+O(1)} \sum_{m\ge C} 2^{-(\frac12-2\eps)m} \|
\psi_1\|_{S[k_1]} \, 2^{-\eps m} \| \tilde Q_m P_k \calN_{\beta j}(
\psi_2, \psi_3)
] \|_{\Ltwotx} \\
&\le \delta   \prod_{i=1}^3 \|\psi_i\|_{S[k_i]}
\end{align*}
where one argues as in the previous two cases to pass to the last
line. Next, suppose the output is limited by~$Q_{\le0}$. Then
\begin{equation}\label{eq:wejustdid'}
\begin{aligned}
 & \| P_0 Q_{\le 0} \del^\beta [  I^c \psi_1 \Delta^{-1}\del_j I^c \calN_{\beta j}( \psi_2, \psi_3)
] \|_{N[0]} \les \sum_{k=k_2-L'}^{k_2+O(1)}  \sum_{m\ge C}  \| Q_m
\psi_1 \Delta^{-1}\del_j  P_k \tilde Q_m \calN_{\beta j}(  \psi_2,
\psi_3)
] \|_{\enerN} \\
&\les  \sum_{k=k_2-L'}^{k_2+O(1)}  \sum_{m\ge C}  \|Q_m \psi_1\|_{L^2_t L^p_x} \, 2^{-k} \| P_k \tilde Q_m \calN_{\beta j}(  \psi_2, \psi_3) \|_{L^2_t L^M_x}\\
&\les  2^{m_0(\frac12-\frac{1}{p})} \sum_{k=k_2-L'}^{k_2+O(1)}
\sum_{m\ge C} 2^{-(1-2\eps)m}  \| \psi_1\|_{S[k_1]} \, 2^{-\eps m}
2^{|k_2|\frac{2}{M}}
 \| P_k \tilde Q_m \calN_{\beta j}(  \psi_2, \psi_3) \|_{\Ltwotx} \\
&\les 2^{m_0(\frac12-\frac{1}{p})} \sum_{k=k_2-L'}^{k_2+O(1)}     \|
\psi_1\|_{S[k_1]} \; 2^{\frac{k}{2}} 2^{-\eps k_2}
2^{|k_2|\frac{2}{M}} \|\psi_2\|_{S[k_2]} \|\psi_3\|_{S[k_3]} \le
\delta  \prod_{i=1}^3 \|\psi_i\|_{S[k_i]}
\end{aligned}
\end{equation}
which is again admissible. To conclude the case of angular alignment
between $\psi_1,\psi_2$, we sum over~$\kappa_1,\kappa_2$ using
Cauchy-Schwarz, Lemma~\ref{lem:enersquaresum}, and Corollary~\ref{cor:Nabhighmod}.

Finally, consider the case where $\psi_2$ and $\psi_3$ are aligned
on the Fourier side. Using Lemmas~\ref{lem:Nablowmod'}
   and~\ref{lem:Nabhighmod'} instead of Lemmas~\ref{lem:Nablowmod}
   and~\ref{lem:Nabhighmod}, respectively,  one immediately verifies that the
   desired gain can indeed be obtained. The only exception here is
   the estimate~\eqref{eq:5oppose}. But this case is excluded here
   as it involves~$I\calN_{\beta j}$ and not~$I^c\calN_{\beta j}$.
\end{proof}

Next, we need to obtain an analogous statement in the hyperbolic
regime of the inner nullform. As in Corollary~\ref{cor:epstrilin},
Lemma~\ref{lem:tri_hyp} implies the following result.

\begin{cor}
  \label{cor:epstrilin'} Let $\delta>0$ be small.
Then
\begin{align*}
&  {\sum}''_{k_1,k_2,k_3\in\Z} \sum_{k\le k_2-L'}  \, \big\|
\sum_{j=1}^2 P_0  \del^\beta [ R_\alpha P_{k_1} \psi_1
\Delta^{-1}\del_j I P_k  \calN_{\beta j}(
P_{k_2 } \psi_2, P_{k_3 } \psi_3)] \big\|_{N[0]} \\
& +  {\sum}''_{k_1,k_2,k_3\in\Z} \sum_{k\le k_2-L'}  \, \big\|
\sum_{j=1}^2 P_0  \del_\alpha [ R^\beta P_{k_1} \psi_1
\Delta^{-1}\del_j I P_k \calN_{\beta j}(P_{k_2 } \psi_2,
P_{k_3 }\psi_3)] \big\|_{N[0]} \\
& +  {\sum}''_{k_1,k_2,k_3\in\Z} \sum_{k\le k_2-L'}  \, \big\|
\sum_{j=1}^2 P_0  \del_\alpha [ R^\beta P_{k_1} \psi_1
\Delta^{-1}\del_j I P_k \calN_{\beta j}( P_{k_2 } \psi_2,
P_{k_3 }\psi_3)] \big\|_{N[0]} \\
&\le  \delta\,  K^2
   \sup_{k\in\Z}\max_{i=1,2,3} 2^{-\sigma_0 |k|} \|P_k \psi_i\|_{S[k]}
\end{align*}
where $L'=L'(L,\delta)$  is a large constant.
\end{cor}

Next, we need to obtain an improvement in the
range~\eqref{eq:outside} under the additional assumption of angular
alignment.

\begin{lemma}
  \label{lem:hyp_trilindelta} For any $\delta>0$ there exists
  $m_0(\delta)$, a large negative constant, such that
\begin{align*}
& \sum_{\substack{\kappa_2,\kappa_3\in\calC_{m_0 }\\
\dist(\kappa_2,\kappa_3)\le 2^{m_0}}} {\sum}''_{k_1,k_2,k_3\in\Z}
\sum_{k=k_2-L'}^{k_2+O(1)} \,  \big\| \sum_{j=1}^2 P_0  \del^\beta [
R_\alpha P_{k_1} \psi_1 \Delta^{-1}\del_j I P_k \calN_{\beta j}(
P_{k_2,\kappa_2} \psi_2, P_{k_3,\kappa_3} \psi_3)] \big\|_{N[0]} \\
&\le  \delta\,  K^2
   \sup_{k\in\Z}\max_{i=1,2,3} 2^{-\sigma_0 |k|} \|P_k \psi_i\|_{S[k]}
\end{align*}
as well as
\begin{align*}
& \sum_{\substack{\kappa_2,\kappa_3\in\calC_{m_0 }\\
\dist(\kappa_2,\kappa_3)\le 2^{m_0}}} {\sum}''_{k_1,k_2,k_3\in\Z}
\sum_{k=k_2-L'}^{k_2+O(1)} \, \big\| \sum_{j=1}^2 P_0  \del_\alpha [
R^\beta P_{k_1} \psi_1 \Delta^{-1}\del_j I P_k \calN_{\beta
j}(P_{k_2,\kappa_2} \psi_2,
P_{k_3,\kappa_3}\psi_3)] \big\|_{N[0]} \\
& +\sum_{\substack{\kappa_2,\kappa_3\in\calC_{m_0 }\\
\dist(\kappa_2,\kappa_3)\le 2^{m_0}}} {\sum}''_{k_1,k_2,k_3\in\Z}
\sum_{k=k_2-L'}^{k_2+O(1)} \, \big\| \sum_{j=1}^2 P_0  \del_\beta [
R^\beta P_{k_1} \psi_1 \Delta^{-1}\del_j I P_k \calN_{\beta j}(
P_{k_2,\kappa_2} \psi_2,
P_{k_3,\kappa_3}\psi_3)] \big\|_{N[0]} \\
&\le  \delta\,  K^2
   \sup_{k\in\Z}\max_{i=1,2,3} 2^{-\sigma_0 |k|} \|P_k \psi_i\|_{S[k]}
\end{align*}
for any $\alpha=0,1,2$. An analogous statement holds in case
$\psi_1,\psi_2$ or $\psi_1,\psi_3$ are similarly aligned.
\end{lemma}
\begin{proof} We begin with the first trilinear form, and also assume alignment between~$\psi_2$ and~$\psi_3$.
We first reduce ourselves to the purely hyperbolic case, i.e., when
all inputs are restricted by the operator~$I$, as well as the entire
output. Without further mention, implicit constants are allowed to
depend on~$L,L'$. In particular, we assume that $k,k_1,k_2,k_3$ are
fixed in the range we are summing over. In the notation of
Lemma~\ref{lem:hyp_redux}, if $A_0=I^c$, then~$A_1=I^c$ and by
Lemma~\ref{lem:Nablowmod'},
\begin{align*}
 &\sum_{\substack{\kappa_2,\kappa_3\in\calC_{m_0 }\\
\dist(\kappa_2,\kappa_3)\le 2^{m_0}}}  \| P_0 I^c \del^\beta [ I^c
\nabla_{t,x} \psi_1 \Delta^{-1}\del_j I \calN_{\beta j}(
P_{k_2,\kappa_2}\psi_2, P_{k_3,\kappa_3} \psi_3)
] \|_{N[0]} \\
& \les \sum_{\substack{\kappa_2,\kappa_3\in\calC_{m_0 }\\
\dist(\kappa_2,\kappa_3)\le 2^{m_0}}} \sum_{m\ge0} 2^{-\eps m}  \|
P_0 Q_m [ \tilde Q_m \nabla_{t,x} \psi_1\, \Delta^{-1}\del_j P_k I
\calN_{\beta j}( P_{k_2,\kappa_2} \psi_2, P_{k_3,\kappa_3}\psi_3)
] \|_{\Ltwotx} \\
&\les  \sum_{\substack{\kappa_2,\kappa_3\in\calC_{m_0 }\\
\dist(\kappa_2,\kappa_3)\le 2^{m_0}}}   \sum_{m\ge0} 2^{(1-\eps)m}
\| \tilde Q_m \psi_1\|_{\Ltwotx}  2^{-k} \|  P_k I \calN_{\beta
j}(P_{k_2,\kappa_2} \psi_2, P_{k_3,\kappa_3}\psi_3)
] \|_{\Linf} \\
&\les    \sum_{\substack{\kappa_2,\kappa_3\in\calC_{m_0 }\\
\dist(\kappa_2,\kappa_3)\le 2^{m_0}}}    \|\psi_1\|_{S[k_1]} \, 2^{\frac{k}{2}}
\| I P_k \calN_{\beta j}( P_{k_2,\kappa_2} \psi_2, P_{k_3,\kappa_3}
 \psi_3) \|_{\Ltwotx} \\
&\les   \delta  2^{\frac{k}{2}} \|\psi_1\|_{S[k_1]}
2^{\frac{k_2}{2}} \|P_{k_2} \psi_2\|_{S[k_2]} \| P_{k_3}
\psi_3\|_{S[k_3]} \le \delta
 \|\psi_1\|_{S[k_1]} \|P_{k_2} \psi_2\|_{S[k_2]} \| P_{k_3}
\psi_3\|_{S[k_3]}
\end{align*}
Summing over $k_1=O(1)$, $k_2=k_3+O(1)$ yields the desired gain.
Hence, we can assume that $A_0=I$ as well as $A_1=I$. If $A_2=I^c$,
then also~$A_3=I^c$ and
\begin{align*}
 & \| P_0 I \del^\beta [ I \psi_1 \Delta^{-1}\del_j I \calN_{\beta j}( I^c P_{k_2,\kappa_2} \psi_2, I^c P_{k_3,\kappa_3} \psi_3)
] \|_{N[0]}\nn  \\
&\les  \sum_{m\ge k_2+C}   \| I \psi_1 \Delta^{-1}\del_j I P_k
\calN_{\beta j}(I^c P_{k_2,\kappa_2} \psi_2, I^c P_{k_3,\kappa_3}
\psi_3)
] \|_{\enerN} \nn \\
&\les      \|\psi_1\|_{L^\infty_t L^2_x} \, 2^{-k} \| I P_k \calN_{\beta j}( I^c P_{k_2,\kappa_2}\psi_2, I^c P_{k_3,\kappa_3}\psi_3) \|_{L^1_t L^\infty_x}  \\
&\les  \|\psi_1\|_{\ener}  \, \sum_{m\ge k_2+C} 2^{-k_2} \|
\calN_{\beta j}( Q_m P_{k_2,\kappa_2} \psi_2, \tilde Q_m
P_{k_3,\kappa_3} \psi_3) \|_{L^1_t L^\infty_x}
\end{align*}
Splitting the modulations of the last two inputs dyadically yields
\begin{align*}
&\les \|\psi_1\|_{\ener}  \, \sum_{m\ge k_2+C} 2^{m-k_2}
  \| P_{k_2,\kappa_2} Q_m
\psi_2\|_{L^2_t L^\infty_x}
\| P_{k_3,\kappa_3}\tilde Q_m \psi_3\|_{L^2_t L^\infty_x} \nn  \\
&\les  2^{m_0+k_2}  \|\psi_1\|_{\ener}  \, \sum_{m\ge k_2+C}
2^{m-k_2} 2^{-2(1-\eps)m} 2^{(1-2\eps)k_2}   \| P_{k_2,\kappa_2}Q_m
\psi_2\|_{\dot X_{k_2}^{-\frac12+\eps,1-\eps,\infty}}
\| P_{k_3,\kappa_3} \tilde Q_m\psi_3\|_{\dot X_{k_3}^{-\frac12+\eps,1-\eps,\infty}} \nn  \\
&\les  2^{m_0}   \|\psi_1\|_{S[k_1]}
\|P_{k_2,\kappa_2}\psi_2\|_{\dot X_{k_2}^{-\frac12+\eps,1-\eps,2} }
\|P_{k_3,\kappa_3}\psi_3\|_{\dot X_{k_3}^{-\frac12+\eps,1-\eps,2} }
\end{align*}
Summing over the caps $\kappa_2,\kappa_3$ and $k_1=O(1)$,
$k_2=k_3+O(1)$ yields the desired gain.

\noindent We may therefore assume that $A_0=A_1=A_2=A_3=I$, which reduces us to
the trilinear nullform expansion~\eqref{eq:nullexp2} restricted to
Case~5 of
Lemma~\ref{lem:tri_hyp}. Beginning with the first of the trilinear nonlinearities
and for the case of aligned $\psi_2,\psi_3$, we now modify the analysis of Case~5 from
that lemma. For ease of notation
we will fix caps~$\kappa_2,\kappa_3$ and drop the projections $P_{k_i,\kappa_i}$. In the end,
an application of the Cauchy-Schwarz inequality will allow for summation over the caps.
We first limit the modulation of~$\psi_1$:
\begin{equation} \begin{aligned}
 &    \|P_0 \del^\beta Q_{>k} I[ Q_{>k+C} I\psi_1 \;
 \Delta^{-1}\del_j P_k I(R_\beta \psi_2 R_j\psi_3 - R_j \psi_2 R_\beta \psi_3)]\|_{N[0]}  \nn  \\
&\les      \|P_0 \del^\beta Q_{>k} I [ Q_{>k+C}I
\psi_1 \;[\Delta^{-1}\del_{j\beta}^2 P_k I(|\nabla|^{-1} \psi_2
R_j\psi_3)
 - P_kI(|\nabla|^{-1} \psi_2 R_\beta \psi_3)]]\|_{\dot X^{0,-\frac12,1}} \nn \\
&\les     2^{-\frac{k}{2}} \|Q_{>k+C}\psi_1\|_{\Ltwotx} \|
\Delta^{-1}\del_{j\beta}^2 P_k I(|\nabla|^{-1} \psi_2 R_j\psi_3)
 - P_k I(|\nabla|^{-1} \psi_2 R_\beta \psi_3)\|_{L^\infty_t L^\infty_x}  \nn  \\
&\les      2^{-k}
\| \psi_1\|_{S[k_1]}  2^{k_2} 2^{m_0} \|\psi_1\|_{S[k_1]}  \|
\psi_2\|_{L^\infty_t L^2_x} \| \psi_3\|_{L^\infty_t L^2_x}  \le \delta
\prod_{i=1}^3 \|\psi_i\|_{S[k_i]}
\end{aligned}\label{eq:psi1limit}
\end{equation}
where the gain is a result of Bernstein's inequality.
Summation over $\kappa_2,\kappa_3$ is admissible here in view of Lemma~\ref{lem:enersquaresum}.
Hence, if the inner output has frequency $\sim 2^k$ then  we may
assume that $\psi_1$ has modulation~$\les 2^k$. Next, we
apply~\eqref{eq:nullexp2} and bound the six terms on the right-hand side of that identity one by one. Previously, we estimated the
first term by means of the Strichartz
component~\eqref{eq:Sk2}. However, this does not seem to yield the angular improvement so we use a different argument:
\begin{equation}
 \label{eq:firsta}
\begin{aligned}
&    \|P_0 I \Box (Q_{\le k}\psi_1 P_k
I[|\nabla|^{-1}\psi_2 |\nabla|^{-1} \psi_3])\|_{N[0]} \\
& \les \sum_{a\le k+C}    \|P_0 Q_{a}  (Q_{\le k} \psi_1 P_k I[ |\nabla|^{-1}\psi_2 |\nabla|^{-1} \psi_3])\|_{\dot X_0^{0,\frac12,1}} \\
&\les  \sum_{a\le j\le k+C}  2^{\frac{a}{2}}  \|Q_{\le k} \psi_1\|_{\ener} \, 2^{k}
  \| P_k Q_j [ |\nabla|^{-1}\psi_2 |\nabla|^{-1} \psi_3] \|_{\Ltwotx} \\
 & + \sum_{j\le a\le  k+C}  2^{k} 2^{\frac{a-k}{4}}  \|Q_{\le k} \psi_1\|_{S[k_1]}
  \| P_k Q_j [ |\nabla|^{-1}\psi_2 |\nabla|^{-1} \psi_3] \|_{ \dot X_k^{0,\frac12,1}}
\end{aligned}
\end{equation}
Lemma~\ref{lem:Sk_prod2} was used to pass to the last line. By Corollary~\ref{cor:capsbilin}
one can continue as follows:
\begin{equation}\begin{aligned}
 &\les \sum_{ j\le k+C}  2^{\frac{j}{2}}  \|  \psi_1\|_{\ener} \, 2^{k} \delta 2^{-\frac{j-k_2}{3}} 2^{-\frac{3k_2}{2}}
  \|  \psi_2 \|_{S[k_2]} \| \psi_3 \|_{S[k_3]} \\
 & + \sum_{j\le  k+C}  2^{k}  \|  \psi_1\|_{S[k_1]}  \: \delta 2^{-\frac{j-k_2}{3}} 2^{-\frac{3k_2}{2}} 2^{\frac{j}{2}}
  \| \psi_2\|_{S[k_2]} \|\psi_3\|_{S[k_3]}
\le \delta \prod_{i=1}^3 \|\psi_i\|_{S[k_i]}
\end{aligned}
\label{eq:firstb}
\end{equation}
Moreover, Corollary~\ref{cor:capsbilin} shows that this bound allows
for summation over the caps.

For the second term, we can assume that $\psi_1=Q_{\le k_2
+C}\psi_1$, see above.  Then, by Corollary~\ref{cor:angcore} as well
as Corollary~\ref{cor:capsbilin}, and some large constant~$M$,
\begin{align*}
 & \sum_{\substack{\kappa_2,\kappa_3\in\calC_{m_0 }\\
\dist(\kappa_2,\kappa_3)\le 2^{m_0}}}  \|P_0 I [\Box ( \psi_1
|\nabla|^{-1} P_{k_3,\kappa_3}\psi_3) |\nabla|^{-1}P_{k_2,\kappa_2}
\psi_2]
 \|_{N[0]} \\
 &\les 2^{k_2} |m_0| \sum_{j\le k_2+C} 2^{\frac{j-k_2}{4}} \Big(\sum_{\kappa_3\in\calC_{m_0
 }}
  \| \tilde P_{0}   \Box Q_j (\psi_1  |\nabla|^{-1} P_{k_3,\kappa_3} \psi_3) \|^2_{\dot X_{0}^{0,-\frac12,\infty}}\Big)^{\frac12} \| |\nabla|^{-1}\psi_2\|_{S[k_2]} \\
&\les  2^{\frac{m_0}{M}} |m_0| \sum_{j\le k_2+C} 2^{\frac{j-k_2}{4}}
2^{\frac{k_2-j}{3}} 2^{\frac{k_2}{2}} 2^{\frac{j}{2}}
2^{-k_2}\prod_{i=1}^3 \|\psi_i\|_{S[k_i]} \les \delta \prod_{i=1}^3
\|\psi_i\|_{S[k_i]}
\end{align*}
Third, by Lemma~\ref{lem:core} and~\eqref{eq:L2bdhl_delta'} of
Corollary~\ref{cor:capsbilin},
\begin{align*}
 & \sum_{\substack{\kappa_2,\kappa_3\in\calC_{m_0 }\\
\dist(\kappa_2,\kappa_3)\le 2^{m_0}}} \| P_0 I [Q_{\le k_2+C} \psi_1
\Box(|\nabla|^{-1} P_{k_2,\kappa_2}\psi_2)
|\nabla|^{-1}P_{k_3,\kappa_3}
\psi_3]\|_{N[0]}  \\
 & \les \sum_{j\le k_2+C} 2^{\frac{j-k_2}{4}} \sum_{\substack{\kappa_2,\kappa_3\in\calC_{m_0 }\\
\dist(\kappa_2,\kappa_3)\le 2^{m_0}}}
   \|\tilde P_{0} Q_{\le k_2+C} ( \psi_1 |\nabla|^{-1}P_{k_3,\kappa_3}\psi_3) \|_{\dot X_{0}^{0,\frac12,1}}
   \| \Box Q_j(|\nabla|^{-1} P_{k_2,\kappa_2}\psi_2) \|_{\dot X_{k_2}^{0,-\frac12,\infty}}\\
& \les \Big( \sum_{\kappa_3\in\calC_{m_0}}   \|\tilde P_{0} Q_{\le
k_2+C} ( \psi_1 |\nabla|^{-1}P_{k_3,\kappa_3}\psi_3) \|_{\dot X_{0}^{0,\frac12,1}}^2    \Big)^{\frac12} \|\psi_2\|_{S[k_2]} \\
&\les \sum_{\ell\le k_2+C} \Big( \sum_{\kappa_3\in\calC_{m_0}}
\|\tilde P_{0} Q_{\ell} ( \psi_1
|\nabla|^{-1}P_{k_3,\kappa_3}\psi_3) \|_{\dot X_{0}^{0,\frac12,1}}^2    \Big)^{\frac12} \|\psi_2\|_{S[k_2]} \\
&\les \sum_{\ell\le k_2+C} \delta 2^{\frac{\ell}{2}}
2^{\frac{k_2-\ell}{3}} 2^{\frac{k_3}{2}} 2^{-k_3}  \prod_{i=1}^3
\|\psi_i\|_{S[k_i]} \les \delta \prod_{i=1}^3 \|\psi_i\|_{S[k_i]}
\end{align*}
The summation over the caps was carried out explicitly for the
second and third terms since it requires some care.  Fourth, by
\eqref{eq:weaker_core} and Corollary~\ref{cor:capsbilin},
\begin{align*}
&     \| P_0 I [ (\Box Q_{\le k+C} \psi_1) \Delta^{-1}
\del_j P_k I(R_j \psi_2 |\nabla|^{-1}\psi_3)] \|_{N[0]} \\ &\les
   \sum_{\ell\le k+C} 2^{\frac{\ell}{4}}
\| \Box Q_\ell \psi_1\|_{\dot X_{k_1}^{0,-\frac12,\infty}} \|\tilde P_{k}Q_{\le k+C} [ R_j \psi_2 |\nabla|^{-1}\psi_3] \|_{\dot X_{k}^{0,\frac12,1}} \\
&\les  \delta  \sum_{\ell\le k+C} \sum_{m\le k+C}  2^{\frac{\ell}{4}}  \| \psi_1\|_{\dot
X_{k_1}^{0,\frac12,\infty}}  \: 2^{\frac{k-j}{3}} 2^{\frac{k_2}{2}} 2^{\frac{m}{2}}
\|\psi_2\|_{S[k_2]} \||\nabla|^{-1} \psi_3\|_{S[k_3]}   \les \delta 2^{\frac{k_2}{4}} \prod_{i=1}^3
\|\psi_i\|_{S[k_i]}
\end{align*}
Since $k=k_1+O(1)=k_2+O(1)$, the fifth term
\begin{align*}
&     \| P_0 I  \Box [Q_{\le k+C} \psi_1  \Delta^{-1} \del_j P_k I (R_j\psi_2 |\nabla|^{-1}\psi_3)]  \|_{N[0]}
\end{align*}
is bounded exactly like the first, see~\eqref{eq:firsta}, \eqref{eq:firstb}.
The sixth and final term is
estimated by means of~\eqref{eq:core} and Corollary~\ref{cor:capsbilin}:
\begin{align*}
&     \|P_0 I [Q_{\le k+C}\psi_1 \Box\Delta^{-1}\del_j
P_k I (R_j \psi_2 |\nabla|^{-1}\psi_3)] \|_{N[0]} \\
&\les    \|\psi_1\|_{S[k_1]}  \:  2^k  \|   P_{k} Q_{\le
k+C}
\Box\Delta^{-1}\del_j (R_j \psi_2 |\nabla|^{-1}\psi_3) \|_{\dot X_{k}^{0,-\frac12,1} }   \\
&  \les     \|\psi_1\|_{S[k_1]}  \sum_{m\le k + C}  2^k \| P_{k} Q_{m}   (R_j \psi_2 |\nabla|^{-1}\psi_3) \|_{\dot X_{k}^{0,\frac12,1} }  \\
&  \les  \delta   \|\psi_1\|_{S[k_1]} \sum_{m\le k + C}  2^k \,
 2^{\frac{k-m}{3}} 2^{\frac{m}{2}}  2^{-\frac{k_2}{2}}   \|\psi_2\|_{S[k_2]}
\|\psi_3\|_{S[k_3]} \les \delta 2^{k_2}  \prod_{i=1}^3 \|\psi_i\|_{S[k_i]}
\end{align*}
as claimed.

We now repeat this analysis for the case of alignment between $\psi_1$ and~$\psi_3$ (the remaining case being symmetric).
We again begin with the reduction of various modulations.
Using the notation of Lemma~\ref{lem:hyp_redux}, if $A_0=I^c$, then~$A_1=I^c$. By~\eqref{eq:angphi2} of
Lemma~\ref{lem:Nablowmod} and with $\frac12=\frac1p+\frac1q$ where $q<\infty$ is very large,
\begin{align*}
 &\sum_{\substack{\kappa_1,\kappa_3\in\calC_{m_0 }\\
\dist(\kappa_1,\kappa_3)\le 2^{m_0}}}  \| P_0 I^c \del^\beta [ I^c
\nabla_{t,x} P_{k_1,\kappa_1} \psi_1 \Delta^{-1}\del_j I \calN_{\beta j}(
\psi_2, P_{k_3,\kappa_3} \psi_3)
] \|_{N[0]} \\
& \les \sum_{\substack{\kappa_1,\kappa_3\in\calC_{m_0 }\\
\dist(\kappa_1,\kappa_3)\le 2^{m_0}}} \sum_{m\ge0} 2^{-\eps m}  \|
P_0 Q_m [ \tilde Q_m \nabla_{t,x} P_{k_1,\kappa_1} \psi_1\, \Delta^{-1}\del_j P_k I
\calN_{\beta j}(   \psi_2, P_{k_3,\kappa_3}\psi_3)
] \|_{\Ltwotx} \\
&\les  \sum_{m\ge0} 2^{(1-\eps)m}    \sum_{\substack{\kappa_1,\kappa_3\in\calC_{m_0 }\\
\dist(\kappa_1,\kappa_3)\le 2^{m_0}}}
\| \tilde Q_m P_{k_1,\kappa_1} \psi_1\|_{L^2_t L^p_x}\;  2^{-k} \|  P_k I \calN_{\beta
j}(  \psi_2, P_{k_3,\kappa_3}\psi_3)
] \|_{L^\infty_t L^q_x} \\
&\les  2^{m_0(\frac12-\frac1p)}  \sum_{m\ge0} 2^{(1-\eps)m} \Big(  \sum_{ \kappa_1 \in\calC_{m_0 } }
 \|P_{k_1,\kappa_1} \tilde Q_m \psi_1\|_{\Ltwotx}^2\Big)^{\frac12} \: 2^{(\frac12-\frac{2}{q})k} \Big(  \sum_{ \kappa_3 \in\calC_{m_0 } }
\| I P_k \calN_{\beta j}(   \psi_2, P_{k_3,\kappa_3}
 \psi_3) \|_{\Ltwotx}^2\Big)^{\frac12}  \\
& \les  |m_0|\, 2^{m_0(\frac12-\frac1p)}  2^{(1-\frac{2}{q})k}
\|\psi_1\|_{\dot X_{0}^{0,1-\eps,2}} \|P_{k_2} \psi_2\|_{S[k_2]} \|
P_{k_3} \psi_3\|_{S[k_3]} \les    \delta \prod_{i=1}^3
 \|\psi_i\|_{S[k_i]}
\end{align*}
Hence, we can assume that $A_0=I$ as well as $A_1=I$. If $A_2=I^c$,
then also~$A_3=I^c$ and
\begin{align*}
 & \| P_0 I \del^\beta [ I P_{k_1,\kappa_1} \psi_1 \Delta^{-1}\del_j I \calN_{\beta j}( I^c  \psi_2, I^c P_{k_3,\kappa_3} \psi_3)
] \|_{N[0]}\nn  \\
&\les  \sum_{m\ge k_2+C}   \| I P_{k_1,\kappa_1} \psi_1 \Delta^{-1}\del_j I P_k
\calN_{\beta j}(I^c   \psi_2, I^c P_{k_3,\kappa_3}
\psi_3)
] \|_{\enerN} \nn \\
&\les      \|P_{k_1,\kappa_1} \psi_1\|_{L^\infty_t L^2_x} \, 2^{-k} \| I P_k \calN_{\beta j}( I^c  \psi_2, I^c P_{k_3,\kappa_3}\psi_3) \|_{L^1_t L^\infty_x}  \\
&\les  \|P_{k_1,\kappa_1} \psi_1\|_{\ener}  \, \sum_{m\ge k_2+C} 2^{-k_2} \|
\calN_{\beta j}( Q_m   \psi_2, \tilde Q_m
P_{k_3,\kappa_3} \psi_3) \|_{L^1_t L^\infty_x}\\
&\les \|P_{k_1,\kappa_1} \psi_1\|_{\ener}  \, \sum_{m\ge k_2+C} 2^{m-k_2}
  \|   Q_m
\psi_2\|_{L^2_t L^\infty_x}
\| P_{k_3,\kappa_3}\tilde Q_m \psi_3\|_{L^2_t L^\infty_x} \nn  \\
&\les  2^{\frac{m_0}{2}+k_2}  \|P_{k_1,\kappa_1} \psi_1\|_{\ener}  \, \sum_{m\ge k_2+C}
2^{m-k_2} 2^{-2(1-\eps)m} 2^{(1-2\eps)k_2}   \|  Q_m
\psi_2\|_{\dot X_{k_2}^{-\frac12+\eps,1-\eps,\infty}}
\| P_{k_3,\kappa_3} \tilde Q_m\psi_3\|_{\dot X_{k_3}^{-\frac12+\eps,1-\eps,\infty}} \nn  \\
&\les  2^{\frac{m_0}{2}}    \|P_{k_1,\kappa_1}\psi_1\|_{\ener} \|
\psi_2\|_{\dot X_{k_2}^{-\frac12+\eps,1-\eps,2} }
\|P_{k_3,\kappa_3}\psi_3\|_{\dot X_{k_3}^{-\frac12+\eps,1-\eps,2} }
\end{align*}
Summing over the caps $\kappa_1,\kappa_3$ and $k_1=O(1)$,
$k_2=k_3+O(1)$ yields the desired gain. For $\psi_1$ one uses Lemma~\ref{lem:enersquaresum}.
As before, this reduces us  to
the trilinear nullform expansion~\eqref{eq:nullexp2}. By the estimate~\eqref{eq:psi1limit}, it
suffices to consider~$P_{\le k+C}\psi_1$  if the inner output has frequency $\sim 2^k$.
Beginning with the first term on the right-hand side of~\eqref{eq:nullexp2},  one has
\begin{equation}
 \label{eq:firsta'}
\begin{aligned}
&  \sum_{\substack{\kappa_1,\kappa_3\in\calC_{m_0 }\\
\dist(\kappa_1,\kappa_3)\le 2^{m_0}}}  \|P_0 I \Box (Q_{\le k}P_{k_1,\kappa_1}\psi_1\: P_k
I[|\nabla|^{-1}\psi_2 |\nabla|^{-1}P_{k_3,\kappa_3} \psi_3])\|_{N[0]} \\
& \les  \sum_{a\le k+C}   \sum_{\substack{\kappa_1,\kappa_3\in\calC_{m_0 }\\
\dist(\kappa_1,\kappa_3)\le 2^{m_0}}}   \|P_0 Q_{a}  (Q_{\le k}P_{k_1,\kappa_1} \psi_1 P_k I[ |\nabla|^{-1}\psi_2 |\nabla|^{-1} P_{k_3,\kappa_3}\psi_3])\|_{\dot X_0^{0,\frac12,1}} \\
&\les  \sum_{a\le j\le k+C}  2^{\frac{a}{2}} \, 2^{k} \Big(\sum_{\kappa_1}  \|Q_{\le k} P_{k_1,\kappa_1}\psi_1\|_{\ener}^2\Big)^{\frac12}
  \Big(\sum_{\kappa_1} \| P_k Q_j [ |\nabla|^{-1}\psi_2 |\nabla|^{-1} P_{k_3,\kappa_3}\psi_3] \|_{\Ltwotx}^2\Big)^{\frac12} \\
 & + \sum_{j\le a\le  k+C}  2^{\frac{3k}{4}}  2^{\frac{a}{4}}  \|  \psi_1\|_{S[k_1]}
  \Big(\sum_{\kappa_3} \| P_k Q_j [ |\nabla|^{-1}\psi_2 |\nabla|^{-1} P_{k_3,\kappa_3}\psi_3] \|_{\dot X_k^{0,\frac12,\infty}}^2\Big)^{\frac12}
\end{aligned}
\end{equation}
Corollary~\ref{cor:bilinbasicsquare} was used to pass to the last line. By Lemma~\ref{lem:enersquaresum} and Corollary~\ref{cor:capsbilin}
one can continue as follows:
\begin{equation}\begin{aligned}
 &\les \delta \sum_{ j\le k+C}  2^{\frac{j}{2}}  \|  \psi_1\|_{\ener} \, 2^{k}   2^{-\frac{j-k_2}{3}} 2^{-\frac{3k_2}{2}}
  \|  \psi_2 \|_{S[k_2]} \| \psi_3 \|_{S[k_3]} \\
 & + \delta \sum_{j\le  k+C}  2^{k}  \|  \psi_1\|_{S[k_1]}  \:   2^{-\frac{j-k_2}{3}} 2^{-\frac{3k_2}{2}} 2^{\frac{j}{2}}
  \| \psi_2\|_{S[k_2]} \|\psi_3\|_{S[k_3]}
\le \delta \prod_{i=1}^3 \|\psi_i\|_{S[k_i]}
\end{aligned}
\label{eq:firstb'}
\end{equation}
Moreover, Corollary~\ref{cor:capsbilin} shows that this bound allows
for summation over the caps.
For the second term, we can assume that $\psi_1=Q_{\le k_2
+C}\psi_1$, see above.  Then, by Lemma~\ref{lem:core} as well
as Corollary~\ref{cor:capsbilin},
\begin{align*}
 & \sum_{\substack{\kappa_1,\kappa_3\in\calC_{m_0 }\\
\dist(\kappa_1,\kappa_3)\le 2^{m_0}}}  \|P_0 I [\Box ( P_{k_1,\kappa_1} \psi_1
|\nabla|^{-1} P_{k_3,\kappa_3}\psi_3) |\nabla|^{-1}
\psi_2]
 \|_{N[0]} \\
 &\les 2^{k_2}  \sum_{j\le k_2+C} 2^{\frac{j-k_2}{4}}   \sum_{\substack{\kappa_1,\kappa_3\in\calC_{m_0 }\\
\dist(\kappa_1,\kappa_3)\le 2^{m_0}}}
  \| \tilde P_{0}   \Box Q_j (P_{k_1,\kappa_1} \psi_1\,   |\nabla|^{-1} P_{k_3,\kappa_3} \psi_3) \|_{\dot X_{0}^{0,-\frac12,\infty}} \: \| |\nabla|^{-1}\psi_2\|_{S[k_2]} \\
&\les
\delta \sum_{j\le k_2+C} 2^{\frac{j-k_2}{4}}
2^{\frac{k_2-j}{3}} 2^{\frac{k_2}{2}} 2^{\frac{j}{2}}
2^{-k_2}\prod_{i=1}^3 \|\psi_i\|_{S[k_i]} \les \delta \prod_{i=1}^3
\|\psi_i\|_{S[k_i]}
\end{align*}
Third, by Lemma~\ref{lem:core} and~\eqref{eq:L2bdhl_delta'} of
Corollary~\ref{cor:capsbilin},
\begin{align*}
 & \sum_{\substack{\kappa_1,\kappa_3\in\calC_{m_0 }\\
\dist(\kappa_1,\kappa_3)\le 2^{m_0}}} \| P_0 I [Q_{\le k_2+C} P_{k_1,\kappa_1} \psi_1
\Box(|\nabla|^{-1} \psi_2)
|\nabla|^{-1}P_{k_3,\kappa_3}
\psi_3]\|_{N[0]}  \\
 & \les \sum_{j\le k_2+C} 2^{\frac{j-k_2}{4}} \sum_{\substack{\kappa_1,\kappa_3\in\calC_{m_0 }\\
\dist(\kappa_1,\kappa_3)\le 2^{m_0}}}
   \|\tilde P_{0} Q_{\le k_2+C} (P_{k_1,\kappa_1} \psi_1 |\nabla|^{-1}P_{k_3,\kappa_3}\psi_3) \|_{\dot X_{0}^{0,\frac12,1}}
   \| \Box Q_j(|\nabla|^{-1}   \psi_2) \|_{\dot X_{k_2}^{0,-\frac12,\infty}}\\
& \les   \sum_{\substack{\kappa_1,\kappa_3\in\calC_{m_0 }\\
\dist(\kappa_1,\kappa_3)\le 2^{m_0}}}   \|\tilde P_{0} Q_{\le
k_2+C} ( P_{k_1,\kappa_1} \psi_1 |\nabla|^{-1}P_{k_3,\kappa_3}\psi_3) \|_{\dot X_{0}^{0,\frac12,1}} \|\psi_2\|_{S[k_2]} \\
&\les \sum_{\ell\le k_2+C}   \sum_{\substack{\kappa_1,\kappa_3\in\calC_{m_0 }\\
\dist(\kappa_1,\kappa_3)\le 2^{m_0}}}
\|\tilde P_{0} Q_{\ell} ( P_{k_1,\kappa_1} \psi_1
|\nabla|^{-1}P_{k_3,\kappa_3}\psi_3) \|_{\dot X_{0}^{0,\frac12,1}}     \|\psi_2\|_{S[k_2]} \\
&\les \sum_{\ell\le k_2+C} \delta 2^{\frac{\ell}{2}}
2^{\frac{k_2-\ell}{3}} 2^{\frac{k_3}{2}} 2^{-k_3}  \prod_{i=1}^3
\|\psi_i\|_{S[k_i]} \les \delta \prod_{i=1}^3 \|\psi_i\|_{S[k_i]}
\end{align*}
Fourth, by Lemma~\ref{lem:core}, Cauchy-Schwarz applied to  the cap-sum, and Corollary~\ref{cor:capsbilin},
\begin{align*}
&   \sum_{\substack{\kappa_1,\kappa_3\in\calC_{m_0 }\\
\dist(\kappa_1,\kappa_3)\le 2^{m_0}}}  \| P_0 I [ (\Box Q_{\le k+C}P_{k_1,\kappa_1} \psi_1) \Delta^{-1}
\del_j P_k I(R_j \psi_2\, |\nabla|^{-1} P_{k_3,\kappa_3}\psi_3)] \|_{N[0]} \\
&\les |m_0|\!\!\!\
   \sum_{\ell\le k+C} 2^{\frac{\ell}{4}} \| \Box Q_\ell \psi_1\|_{\dot X_{k_1}^{0,-\frac12,\infty}} \Big( \sum_{ \kappa_3\in\calC_{m_0 }}
 \|\tilde P_{k}Q_{\le k+C} [ R_j \psi_2 |\nabla|^{-1}P_{k_3,\kappa_3}\psi_3] \|_{\dot X_{k}^{0,\frac12,1}}^2\Big)^{\frac12}  \\
&\les  \delta  \sum_{\ell\le k+C} \sum_{m\le k+C}  2^{\frac{\ell}{4}}  \| \psi_1\|_{\dot
X_{k_1}^{0,\frac12,\infty}}  \: 2^{\frac{k-m}{3}} 2^{\frac{k_2}{2}} 2^{\frac{m}{2}}
\|\psi_2\|_{S[k_2]} \||\nabla|^{-1} \psi_3\|_{S[k_3]}   \les \delta 2^{\frac{k_2}{4}} \prod_{i=1}^3
\|\psi_i\|_{S[k_i]}
\end{align*}
Since $k=k_1+O(1)=k_2+O(1)$, the fifth term
\begin{align*}
&     \| P_0 I  \Box [Q_{\le k+C} \psi_1  \Delta^{-1} \del_j P_k I (R_j\psi_2 |\nabla|^{-1}\psi_3)]  \|_{N[0]}
\end{align*}
is bounded exactly like the first, see~\eqref{eq:firsta}, \eqref{eq:firstb}.
The sixth and final term is
estimated by means of Corollary~\ref{cor:angcore} and Corollary~\ref{cor:capsbilin}:
\begin{align*}
&  \sum_{\substack{\kappa_1,\kappa_3\in\calC_{m_0 }\\
\dist(\kappa_1,\kappa_3)\le 2^{m_0}}}    \|P_0 I [Q_{\le k+C}P_{k_1,\kappa_1}  \psi_1 \Box\Delta^{-1}\del_j
P_k I (R_j \psi_2 |\nabla|^{-1}P_{k_3,\kappa_3} \psi_3)] \|_{N[0]} \\
&\les  |m_0|   \|\psi_1\|_{S[k_1]}  \:  2^k \sum_{m\le k + C}  \Big(\sum_{\kappa_3\in\calC_{m_0}} \|   P_{k} Q_{m}
\Box\Delta^{-1}\del_j (R_j \psi_2 |\nabla|^{-1}P_{k_3,\kappa_3} \psi_3) \|_{\dot X_{k}^{0,-\frac12,1} }^2 \Big)^{\frac12}  \\
&  \les  \delta   \|\psi_1\|_{S[k_1]} \sum_{m\le k + C}  2^k \,
 2^{\frac{k-m}{3}} 2^{\frac{m}{2}}  2^{-\frac{k_2}{2}}   \|\psi_2\|_{S[k_2]}
\|\psi_3\|_{S[k_3]} \les \delta 2^{k_2}  \prod_{i=1}^3 \|\psi_i\|_{S[k_i]}
\end{align*}
as claimed.
The other two types of trilinear null-forms are similar and left to the reader.
\end{proof}

\begin{remark}\label{rem:slightimprov}
The proof of the preceding estimates actually leads to a slightly better result: letting $P_0F(P_{k_1}\psi_1, P_{k_2}\psi_2, P_{k_3}\psi_3)$ be a frequency localized trilinear null-form as above, then given any $\delta>0$, there exists some $l_0\leq -100$ such that we can write
\[
P_0F(P_{k_1}\psi_1, P_{k_2}\psi_2, P_{k_3}\psi_3)=F_1+F_2
\]
where $F_1$ is a sum of energy, $\dot{X}^{s, b, q}$, as well as wave-packet atoms of scale $l\geq l_0$ (where scale refers to the size $2^l$ of the caps $\kappa$ used), with the bound
\[
\|F_1\|_{N[0]}\les w(k_1, k_2, k_3)\prod_{j=1}^{3}\|P_{k_j}\psi_j\|_{S[k_j]}
\]
and universal implied constant (independent of $\delta$), while we also have
\[
\|F_1\|_{N[0]}\les \delta w(k_1, k_2, k_3)\prod_{j=1}^{3}\|P_{k_j}\psi_j\|_{S[k_j]}
\]
The reason for this is that whenever a wave-packet atom of extremely fine scale is being used to estimate some constituent of $P_0F$, one gains a small exponential power in that scale.
\end{remark}

\section{Quintilinear and higher nonlinearities}\label{sec:quintic}

Here we detail the estimates needed in order to control the higher
order error terms generated by the process described in
Section~\ref{sec:hodge}.  This section is quite technical but the
main point here is that the higher order terms, while still somewhat
complicated, are much easier to estimate than the trilinear
null-forms, and only require a very mild null-structure. We start
with the lowest order errors, of quintilinear type. These are either
of first or second type, see the discussion in
Section~\ref{sec:hodge}. We commence with those of the first type,
which can be schematically written as
\[
\nabla_{x,t}[\psi\nabla^{-1}(R_{\nu}\psi \nabla^{-1}(\psi\nabla^{-1}Q_{\mu j}(\psi,\psi)))],
\]
where not both $\nu, \mu$ are simultaneously zero. Assume that
$\nu=0, \mu\neq 0$, the remaining cases being treated analogously.
The following lemma is then representative for the higher order
errors, for a universal $\delta>0$.

\begin{lemma}\label{quintilinear1} We have the estimates
\begin{equation}\nonumber\begin{split}
&\|\nabla_{x,t}[P_{0}\psi_{0}\nabla^{-1}P_{r_{1}}(R_{0}P_{k_{1}}\psi_{1} \nabla^{-1}P_{r_{2}}(P_{k_{2}}\psi_{2}\nabla^{-1}P_{r_{3}}Q_{jk}(P_{k_{3}}\psi_{3},P_{k_{4}}\psi_{4})))]\|_{N[0]}\\
&\lesssim 2^{\delta[\min_{j\neq 0}\{r_{j}, k_{j}\}-\max_{j\neq
0}\{r_{j},
k_{j}\}]}\prod_{i=0}^{4}\|P_{k_{i}}\psi_{i}\|_{S[k_{i}]},\quad
r_{1}<-10
\end{split}\end{equation}
\begin{equation}\nonumber\begin{split}
&\|\nabla_{x,t}[P_{k_{0}}\psi_{0}\nabla^{-1}P_{r_1}(R_{0}P_{k_{1}}\psi_{1} \nabla^{-1}P_{r_{2}}(P_{k_{2}}\psi_{2}\nabla^{-1}P_{r_{3}}Q_{jk}(P_{k_{3}}\psi_{3},P_{k_{4}}\psi_{4})))]\|_{N[0]}\\
& \lesssim 2^{\delta k_{0}}2^{\delta[\min\{r_{j},
k_{j}\}-\max\{r_{j},
k_{j}\}]}\prod_{i=0}^{4}\|P_{k_{i}}\psi_{i}\|_{S[k_{i}]},\quad
r_{1}\in [-10, 10]
\end{split}\end{equation}
\begin{equation}\nonumber\begin{split}
&\|\nabla_{x,t}P_{0}[P_{k_{0}}\psi_{0}\nabla^{-1}P_{r_{1}}(R_{0}P_{k_{1}}\psi_{1} \nabla^{-1}P_{r_{2}}(P_{k_{2}}\psi_{2}\nabla^{-1}P_{r_{3}}Q_{jk}(P_{k_{3}}\psi_{3},P_{k_{4}}\psi_{4})))]\|_{N[0]}\\
& \lesssim 2^{-\delta k_{0}}2^{\delta[\min\{r_{j},
k_{j}\}-\max\{r_{j},
k_{j}\}]}\prod_{i=0}^{4}\|P_{k_{i}}\psi_{i}\|_{S[k_{i}]},\quad
r_{1}>10
\end{split}\end{equation}
All implied constants are universal.
\end{lemma}
\begin{proof}
All three inequalities are proved similarly, and we treat here the
high-low case in detail, i.e., the first of them. We first deal with
the elliptic cases:
\\

{\it{(i): Output in elliptic regime.}} This is the expression (we have included the gratuitous cutoff $P_{[-5,5]}$ in light of $r_{1}<-10$)
\begin{equation}\nonumber\begin{split}
&\nabla_{x,t}P_{[-5,5]}Q_{>10}[P_{0}\psi_{0}\nabla^{-1}P_{r_{1}}(R_{0}P_{k_{1}}\psi_{1} \nabla^{-1}P_{r_{2}}(P_{k_{2}}\psi_{2}\nabla^{-1}P_{r_{3}}Q_{jk}(P_{k_{3}}\psi_{3},P_{k_{4}}\psi_{4})))]\\
&=\sum_{l>10}\nabla_{x,t}P_{[-5,5]}Q_{l}[P_{0}\psi_{0}\nabla^{-1}P_{r_{1}}(R_{0}P_{k_{1}}\psi_{1} \nabla^{-1}P_{r_{2}}(P_{k_{2}}\psi_{2}\nabla^{-1}P_{r_{3}}Q_{jk}(P_{k_{3}}\psi_{3},P_{k_{4}}\psi_{4})))]\\
\end{split}\end{equation}
Now distinguish between further cases:
\\

{\it{(i1): $\max\{k_{1},\ldots,k_{4}\}\ll l$,
$R_{0}P_{k_{1}}\psi_{1}=R_{0}P_{k_{1}}Q_{<l-100}\psi_{1}$.}} In this
case at least one other factor $P_{k_{j}}\psi_{j}$ has modulation at
least $2^{l-10}$. For argument' s sake, let this be
$P_{k_{2}}\psi_{2}=P_{k_{2}}Q_{>l-10}\psi_{2}$ (the other cases
being similar), so we now reduce to estimating
\[
\sum_{l>10}\nabla_{x,t}P_{[-5,5]}Q_{l}[P_{0}\psi_{0}\nabla^{-1}P_{r_{1}}(R_{0}P_{k_{1}}Q_{<l-100}\psi_{1} \nabla^{-1}P_{r_{2}}(P_{k_{2}}Q_{>l-10}\psi_{2}\nabla^{-1}P_{r_{3}}Q_{jk}(P_{k_{3}}\psi_{3},P_{k_{4}}\psi_{4})))],
\]
where we also make the further assumptions of case (i1). Freezing $l$ for now, we estimate this expression as follows: first, note that we get
\begin{equation}\nonumber\begin{split}
&\|\nabla^{-1}P_{r_{2}}(P_{k_{2}}Q_{>l-10}\psi_{2}\nabla^{-1}P_{r_{3}}Q_{jk}(P_{k_{3}}\psi_{3},P_{k_{4}}\psi_{4}))\|_{L_{t}^{2}\dot{H}_{x}^{\frac{1}{2}}}\\&\lesssim 2^{(1-\epsilon)(k_{2}-l)}2^{[\min\{r_{2,3}, k_{2,3,4}\}-\max\{r_{2,3},k_{2,3,4}\}]}\prod_{j=2}^{4}\|P_{k_{j}}\psi_{j}\|_{S[k_{j}]}
\end{split}\end{equation}
This follows by straightforward usage of Bernstein's inequality and the definition of $S[k]$, as well as exploiting the null-structure of $Q_{jk}$. Furthermore, we have
\[
\|R_{0}P_{k_{1}}Q_{<l-100}\psi_{1}\|_{L_{t,x}^{2}}\lesssim 2^{\epsilon(l-k_{1})}2^{-\frac{k_{1}}{2}}\|P_{k_{1}}\psi\|_{S[k_{1}]},
\]
where $\epsilon>0$ is as in the definition of $S[k]$, which implies that
\[
\|R_{0}P_{k_{1}}Q_{<l-100}\psi_{1}\|_{L_{t}^{\infty}L_{x}^{2}}\lesssim 2^{\epsilon(l-k_{1})}2^{\frac{l-k_{1}}{2}}\|P_{k_{1}}\psi\|_{S[k_{1}]}
\]
From here we get
\begin{equation}\nonumber\begin{split}
&\|\nabla_{x,t}P_{[-5,5]}Q_{l}[P_{0}\psi_{0}\nabla^{-1}P_{r_{1}}(R_{0}P_{k_{1}}Q_{<l-100}\psi_{1}\\&\hspace{4cm}\times\nabla^{-1}P_{r_{2}}(P_{k_{2}}Q_{>l-10}\psi_{2}\nabla^{-1}P_{r_{3}}Q_{jk}(P_{k_{3}}\psi_{3},P_{k_{4}}\psi_{4})))]\|_{\dot{X}_{0}^{-\frac{1}{2}+\epsilon,-1-\epsilon,1}}\\
&\lesssim 2^{-\epsilon l}\|P_{0}\psi_{0}\|_{L_{t}^{\infty}L_{x}^{2}}\\&\times\|\nabla^{-1}P_{r_{1}}(R_{0}P_{k_{1}}Q_{<l-100}\psi_{1} \nabla^{-1}P_{r_{2}}(P_{k_{2}}Q_{>l-10}\psi_{2}\nabla^{-1}P_{r_{3}}Q_{jk}(P_{k_{3}}\psi_{3},P_{k_{4}}\psi_{4})))\|_{L_{t}^{2}L_{x}^{\infty}}\\
&\lesssim 2^{\frac{r_{1}}{2}}2^{\min\{k_{1}-\min\{r_{1,2}\},0\}}2^{\frac{\min\{\min\{r_{1,2}\}-k_{1},0\}}{2}}\|R_{0}P_{k_{1}}Q_{<l-100}\psi_{1}\|_{L_{t}^{\infty}L_{x}^{2}}\\
&\hspace{5cm}\times\|\nabla^{-1}P_{r_{2}}(P_{k_{2}}Q_{>l-10}\psi_{2}\nabla^{-1}P_{r_{3}}Q_{jk}(P_{k_{3}}\psi_{3},P_{k_{4}}\psi_{4}))\|_{L_{t}^{2}\dot{H}_{x}^{\frac{1}{2}}}
\end{split}\end{equation}
Substituting the bounds from before, this is bounded by
\begin{equation}\nonumber\begin{split}
\lesssim &2^{\frac{r_{1}}{2}}2^{\min\{k_{1}-\min\{r_{1,2}\},0\}}2^{\frac{\min\{\min\{r_{1,2}\}-k_{1},0\}}{2}}\\&\times 2^{\epsilon(l-k_{1})}2^{\frac{l-k_{1}}{2}}\|P_{k_{1}}\psi\|_{S[k_{1}]} 2^{(1-\epsilon)(k_{2}-l)}2^{[\min\{r_{2,3}, k_{2,3,4}\}-\max\{r_{2,3},k_{2,3,4}\}]}\prod_{j=2}^{4}\|P_{k_{j}}\psi_{j}\|_{S[k_{j}]}
\end{split}\end{equation}
This is equivalent to an estimate of the form claimed in the lemma, with an extra gain $2^{-\epsilon l}$ which allows us to sum over $l>10$.
\\
{\it{(i2): $\max\{k_{1},\ldots,k_{4}\}\ll l$,
$R_{0}P_{k_{1}}\psi_{1}=R_{0}P_{k_{1}}Q_{[l-100.l+100]}\psi_{1}$.}}
The estimate here is similar except that summation over $l$ is made
possible since we have
\[
P_{k_{1}}Q_{>k_{1}}\psi_{1}\in \dot{X}_{k_{1}}^{-\frac{1}{2}+\epsilon, 1-\epsilon, 1}
\]
{\it{(i3): $\max\{k_{1},\ldots,k_{4}\}\ll l$,
$R_{0}P_{k_{1}}\psi_{1}=R_{0}P_{k_{1}}Q_{\gg l+100]}\psi_{1}$.}}
This is again similar. Fixing the modulation of
$R_{0}P_{k_{1}}Q_{\gg l+100]}\psi_{1}$ to size $2^{l_{1}}$,
$l_{1}>l+100$, there is at least one other input which has
modulation at least comparable to $2^{l_{1}}$. Then one proceeds as
in case (i1).
\\

{\it{(i4):   $\max\{k_{1},\ldots,k_{4}\}>l+O(1)$,  $R_{0}P_{k_{1}}\psi_{1}=R_{0}P_{k_{1}}Q_{<\max\{k_{1,2,3,4}\}}\psi_{1}$}}.
\\
Here we obtain a gain in $\min\{r_{1,2,3}, k_{1,2,3,4}\}-\max\{r_{1,2,3}, k_{1,2,3,4}\}$, which suffices to offset the loss due to the possibly large modulation of $R_{0}P_{k_{1}}Q_{<\max\{k_{1,2,3,4}\}}\psi_{1}$.  Specifically, write
\begin{equation}\nonumber\begin{split}
&\nabla_{x,t}Q_{l}[P_{0}\psi_{0}\nabla^{-1}P_{r_{1}}(R_{0}Q_{<\max\{k_{1,2,3,4}\}}P_{k_{1}}\psi_{1} \\&\hspace{3cm}\times\nabla^{-1}P_{r_{2}}(P_{k_{2}}\psi_{2}\nabla^{-1}P_{r_{3}}Q_{jk}(P_{k_{3}}\psi_{3},P_{k_{4}}\psi_{4})))]\\
&=\nabla_{x,t}Q_{l}[P_{0}\psi_{0}\nabla^{-1}P_{r_{1}}(R_{0}Q_{<k_{1}}P_{k_{1}}\psi_{1}  \\&\hspace{3cm}\times\nabla^{-1}P_{r_{2}}(P_{k_{2}}\psi_{2}\nabla^{-1}P_{r_{3}}Q_{jk}(P_{k_{3}}\psi_{3},P_{k_{4}}\psi_{4})))]\\
&+\nabla_{x,t}Q_{l}[P_{0}\psi_{0}\nabla^{-1}P_{r_{1}}(R_{0}Q_{[k_{1},\max\{k_{1,2,3,4}\}]}P_{k_{1}}\psi_{1} \\&\hspace{3cm}\times \nabla^{-1}P_{r_{2}}(P_{k_{2}}\psi_{2}\nabla^{-1}P_{r_{3}}Q_{jk}(P_{k_{3}}\psi_{3},P_{k_{4}}\psi_{4})))]\\
\end{split}\end{equation}
Here we use the inequalities
\begin{align*}
&\|\nabla^{-1}P_{r_{2}}(P_{k_{2}}\psi_{2}\nabla^{-1}P_{r_{3}}Q_{jk}(P_{k_{3}}\psi_{3},P_{k_{4}}\psi_{4}))\|_{L_{t}^{2}\dot{H}_{x}^{\frac{1}{2}}}\\&\lesssim2^{\frac{1}{2}[\min\{r_{2,3}, k_{1,2,3,4}\}-\max\{r_{2,3}, k_{1,2,3,4}\}]} \prod_{i=1}^{4}\|P_{k_{i}}\psi_{i}\|_{S[k_{i}]},
\end{align*}
\begin{align*}
&\|\nabla^{-1}P_{r_{2}}(P_{k_{2}}\psi_{2}\nabla^{-1}P_{r_{3}}Q_{jk}(P_{k_{3}}\psi_{3},P_{k_{4}}\psi_{4}))\|_{L_{t}^{\infty}L_{x}^{2}}\\&\lesssim2^{\min\{r_{2,3}, k_{1,2,3,4}\}-\max\{r_{2,3}, k_{1,2,3,4}\}} \prod_{i=1}^{4}\|P_{k_{i}}\psi_{i}\|_{S[k_{i}]}.
\end{align*}
Then we can estimate
\begin{align*}
&\|\nabla_{x,t}Q_{l}[P_{0}\psi_{0}\nabla^{-1}P_{r_{1}}(R_{0}Q_{[k_{1},\max\{k_{1,2,3,4}\}]}P_{k_{1}}\psi_{1} \\&\hspace{3cm}\times\nabla^{-1}P_{r_{2}}(P_{k_{2}}\psi_{2}\nabla^{-1}P_{r_{3}}Q_{jk}(P_{k_{3}}\psi_{3},P_{k_{4}}\psi_{4})))]\|_{\dot{X}_{0}^{-\frac{1}{2}+\epsilon, -1-\epsilon,1}}\\
&\lesssim  2^{-\epsilon l}\|P_{0}\psi_{0}\|_{L_{t}^{\infty}L_{x}^{2}}\\&\times \|\nabla^{-1}P_{r_{1}}(R_{0}Q_{[k_{1},\max\{k_{1,2,3,4}\}]}P_{k_{1}}\psi_{1} \nabla^{-1}P_{r_{2}}(P_{k_{2}}\psi_{2}\nabla^{-1}P_{r_{3}}Q_{jk}(P_{k_{3}}\psi_{3},P_{k_{4}}\psi_{4})))\|_{L_{t}^{2}L_{x}^{\infty}}
\end{align*}
To conclude the contribution of this term, one then checks, using  standard Littlewood-Paley trichotomy, that
\begin{align*}
&\|\nabla^{-1}P_{r_{1}}(R_{0}Q_{[k_{1},\max\{k_{1,2,3,4}\}]}P_{k_{1}}\psi_{1} \nabla^{-1}P_{r_{2}}(P_{k_{2}}\psi_{2}\nabla^{-1}P_{r_{3}}Q_{jk}(P_{k_{3}}\psi_{3},P_{k_{4}}\psi_{4})))\|_{L_{t}^{2}L_{x}^{\infty}}\\&\lesssim 2^{\min\{r_{1,2},k_{1}\}}2^{\frac{\min\{r_{1,2},k_{1}\}-\max\{r_{1,2}, k_{1}\}}{2}}\|R_{0}Q_{[k_{1},\max\{k_{1,2,3,4}\}]}P_{k_{1}}\psi_{1}\|_{L_{t}^{2}\dot{H}_{x}^{\frac{1}{2}}}\\&\hspace{5cm}\times\|\nabla^{-1}P_{r_{2}}(P_{k_{2}}\psi_{2}\nabla^{-1}P_{r_{3}}Q_{jk}(P_{k_{3}}\psi_{3},P_{k_{4}}\psi_{4}))\|_{L_{t}^{\infty}L_{x}^{2}}
\end{align*}
Combining with the bound from above, and furthermore assuming the $\epsilon$ in the definition of $\|.\|_{S}$ to be small enough, we conclude that for suitable $\delta>0$ we have
\begin{align*}
&\|\nabla_{x,t}Q_{l}[P_{0}\psi_{0}\nabla^{-1}P_{r_{1}}(R_{0}Q_{[k_{1},\max\{k_{1,2,3,4}\}]}P_{k_{1}}\psi_{1} \\&\hspace{3cm}\times\nabla^{-1}P_{r_{2}}(P_{k_{2}}\psi_{2}\nabla^{-1}P_{r_{3}}Q_{jk}(P_{k_{3}}\psi_{3},P_{k_{4}}\psi_{4})))]\|_{\dot{X}_{0}^{-\frac{1}{2}+\epsilon, -1-\epsilon,1}}\\
&\lesssim
 2^{-\epsilon l}2^{\delta[\min_{j\neq 0}\{r_{j}, k_{j}\}-\max_{j\neq 0}\{r_{j}, k_{j}\}]}\prod_{i=0}^{4}\|P_{k_{i}}\psi_{i}\|_{S[k_{i}]}
\end{align*}
and summing over $10<l<\max\{k_{1,2,3,4}\}$ yields the desired bound.
\\
For the term
\begin{align*}
&\nabla_{x,t}Q_{l}[P_{0}\psi_{0}\nabla^{-1}P_{r_{1}}(R_{0}Q_{<k_{1}}P_{k_{1}}\psi_{1}\\&\hspace{3cm}\times\nabla^{-1}P_{r_{2}}(P_{k_{2}}\psi_{2}\nabla^{-1}P_{r_{3}}Q_{jk}(P_{k_{3}}\psi_{3},P_{k_{4}}\psi_{4})))]
\end{align*}
from further above,  estimate
\begin{align*}
&\|\nabla_{x,t}Q_{l}[P_{0}\psi_{0}\nabla^{-1}P_{r_{1}}(R_{0}Q_{<k_{1}}P_{k_{1}}\psi_{1} \\&\hspace{3cm}\times\nabla^{-1}P_{r_{2}}(P_{k_{2}}\psi_{2}\nabla^{-1}P_{r_{3}}Q_{jk}(P_{k_{3}}\psi_{3},P_{k_{4}}\psi_{4})))]\|_{\dot{X}_{0}^{-\frac{1}{2}+\epsilon, -1-\epsilon,1}}\\
&\lesssim 2^{-\epsilon l}\|P_{0}\psi_{0}\|_{L_{t}^{\infty}L_{x}^{2}}\\&\times\|\nabla^{-1}P_{r_{1}}(R_{0}Q_{<k_{1}}P_{k_{1}}\psi_{1} \nabla^{-1}P_{r_{2}}(P_{k_{2}}\psi_{2}\nabla^{-1}P_{r_{3}}Q_{jk}(P_{k_{3}}\psi_{3},P_{k_{4}}\psi_{4})))\|_{L_{t}^{2}L_{x}^{\infty}},
\end{align*}
and we have
\begin{align*}
&\|\nabla^{-1}P_{r_{1}}(R_{0}Q_{<k_{1}}P_{k_{1}}\psi_{1} \nabla^{-1}P_{r_{2}}(P_{k_{2}}\psi_{2}\nabla^{-1}P_{r_{3}}Q_{jk}(P_{k_{3}}\psi_{3},P_{k_{4}}\psi_{4})))\|_{L_{t}^{2}L_{x}^{\infty}}\\
&\lesssim  2^{\frac{\min\{r_{1,2},k_{1}\}}{2}}2^{\frac{\min\{r_{1,2},k_{1}\}-\max\{r_{1,2}, k_{1}\}}{2}}\|P_{k_{1}}\psi_{1}\|_{L_{t}^{\infty}L_{x}^{2}}\\&\hspace{2cm}\times\| \nabla^{-1}P_{r_{2}}(P_{k_{2}}\psi_{2}\nabla^{-1}P_{r_{3}}Q_{jk}(P_{k_{3}}\psi_{3},P_{k_{4}}\psi_{4}))\|_{L_{t}^{2}\dot{H}^{\frac{1}{2}}},
\end{align*}
which in conjunction with the bound from above
\begin{align*}
&\|\nabla^{-1}P_{r_{2}}(P_{k_{2}}\psi_{2}\nabla^{-1}P_{r_{3}}Q_{jk}(P_{k_{3}}\psi_{3},P_{k_{4}}\psi_{4}))\|_{L_{t}^{2}\dot{H}_{x}^{\frac{1}{2}}}\\&\lesssim2^{\frac{1}{2}[\min\{r_{2,3}, k_{1,2,3,4}\}-\max\{r_{2,3}, k_{1,2,3,4}\}]} \prod_{i=1}^{4}\|P_{k_{i}}\psi_{i}\|_{S[k_{i}]},
\end{align*}
implies that
\begin{align*}
&\|\nabla_{x,t}Q_{l}[P_{0}\psi_{0}\nabla^{-1}P_{r_{1}}(R_{0}Q_{<k_{1}}P_{k_{1}}\psi_{1} \\&\hspace{3cm}\nabla^{-1}P_{r_{2}}(P_{k_{2}}\psi_{2}\nabla^{-1}P_{r_{3}}Q_{jk}(P_{k_{3}}\psi_{3},P_{k_{4}}\psi_{4})))]\|_{\dot{X}_{0}^{-\frac{1}{2}+\epsilon, -1-\epsilon,1}}\\
&\lesssim 2^{-\epsilon l}2^{\frac{\min\{r_{1,2},k_{1}\}}{2}}2^{\frac{\min\{r_{1,2},k_{1}\}-\max\{r_{1,2}, k_{1}\}}{2}}
\\&\hspace{2cm}\times2^{\frac{1}{2}[\min\{r_{2,3}, k_{1,2,3,4}\}-\max\{r_{2,3}, k_{1,2,3,4}\}]} \prod_{i=0}^{4}\|P_{k_{i}}\psi_{i}\|_{S[k_{i}]}
\end{align*}
Summing over $l>10$ yields a bound as claimed in the lemma with $\delta=\frac{1}{2}$.
\\

{\it{(i5):  $\max\{k_{1},\ldots,k_{4}\}>l+O(1)$,  $R_{0}P_{k_{1}}\psi_{1}=R_{0}P_{k_{1}}Q_{\gg\max\{k_{1,2,3,4}\}}\psi_{1}$}}.
\\

Freeze the modulation of $R_{0}P_{k_{1}}\psi_{1}$ to dyadic value $2^{l_{1}}\gg 2^{\max\{k_{1,2,3,4}\}}$. Here there must be at least one other input with modulation at least comparable to $2^{l_{1}}$. Let this input be $P_{k_{2}}\psi_{2}$ for definitiveness' sake, the other cases being treated similarly. Thus consider the term
\[
\nabla_{x,t}Q_{l}[P_{0}\psi_{0}\nabla^{-1}P_{r_{1}}(R_{0}Q_{l_{1}}P_{k_{1}}\psi_{1} \nabla^{-1}P_{r_{2}}(P_{k_{2}}Q_{\geq l_{1}+O(1)}\psi_{2}\nabla^{-1}P_{r_{3}}Q_{jk}(P_{k_{3}}\psi_{3},P_{k_{4}}\psi_{4})))]
\]
Assuming a high-low frequency cascade $r_{1}\ll k_{1}\ll k_{2}$, we can estimate this by (using Bernstein's inequality)
\begin{align*}
&\|\nabla_{x,t}Q_{l}[P_{0}\psi_{0}\nabla^{-1}P_{r_{1}}(R_{0}Q_{l_{1}}P_{k_{1}}\psi_{1}\\&\hspace{3cm}\times\nabla^{-1}P_{r_{2}}(P_{k_{2}}Q_{>l_{1}+O(1)}\psi_{2}\nabla^{-1}P_{r_{3}}Q_{jk}(P_{k_{3}}\psi_{3},P_{k_{4}}\psi_{4})))]\|_{\dot{X}_{0}^{-\frac{1}{2}+\epsilon, -1-\epsilon, 1}}\\
&\lesssim 2^{l(\frac{1}{2}-\epsilon)}\|P_{0}\psi_{0}\|_{L_{t}^{\infty}L_{x}^{2}}\\&\times\|\nabla^{-1}P_{r_{1}}(R_{0}Q_{l_{1}}P_{k_{1}}\psi_{1} \nabla^{-1}P_{r_{2}}(P_{k_{2}}Q_{>l_{1}+O(1)}\psi_{2}\nabla^{-1}P_{r_{3}}Q_{jk}(P_{k_{3}}\psi_{3},P_{k_{4}}\psi_{4})))\|_{L_{t}^{1}L_{x}^{\infty}}\\
&\lesssim  2^{l(\frac{1}{2}-\epsilon)+r_{1}}2^{\min\{r_{3}, k_{3,4}\}-\max\{r_{3}, k_{3,4}\}}\\&\hspace{2cm}\times\|P_{0}\psi_{0}\|_{L_{t}^{\infty}L_{x}^{2}}\|R_{0}Q_{l_{1}}P_{k_{1}}\psi_{1}\|_{L_{t,x}^{2}}\|P_{k_{2}}Q_{>l_{1}+O(1)}\psi_{2}\|_{L_{t,x}^{2}}\prod_{j=3,4}\|P_{k_{j}}\psi_{j}\|_{S[k_{j}]}\\
&\lesssim 2^{l(\frac{1}{2}-\epsilon)}2^{\min\{r_{3}, k_{3,4}\}-\max\{r_{3}, k_{3,4}\}}2^{\epsilon(l_{1}-k_{1})}2^{r_{1}-\frac{k_{1}+k_{2}}{2}}2^{k_{1}-k_{2}}2^{(1-\epsilon)(k_{2}-l_{1})}\prod_{j=0}^{4}\|P_{k_{j}}\psi_{j}\|_{S[k_{j}]}
\end{align*}
Summing over $l_{1}\gg\max\{k_{1,2,3,4}\}>l+O(1)$, one obtains  a bound of the form claimed in the lemma with $\delta=\frac{1}{2}$ in the particular case at hand. The remaining frequency interactions, while keeping our assumptions on the modulations, are treated similarly.
This concludes the {\it{elliptic case (i)}}.
\\

{\it{(ii): Output in hyperbolic regime}}.  Now we consider the expression
\begin{align*}
&\nabla_{x,t}P_{[-5,5]}Q_{<10}[P_{0}\psi_{0}\nabla^{-1}P_{r_{1}}(R_{0}P_{k_{1}}\psi_{1}\\&\hspace{3cm} \times\nabla^{-1}P_{r_{2}}(P_{k_{2}}\psi_{2}\nabla^{-1}P_{r_{3}}Q_{jk}(P_{k_{3}}\psi_{3},P_{k_{4}}\psi_{4})))]
\end{align*}
We decompose this into
\begin{align}
&\nabla_{x,t}P_{[-5,5]}Q_{<10}[P_{0}\psi_{0}\nabla^{-1}P_{r_{1}}(R_{0}P_{k_{1}}\psi_{1}\nonumber \\&\hspace{3cm}\times\nabla^{-1}P_{r_{2}}(P_{k_{2}}\psi_{2}\nabla^{-1}P_{r_{3}}Q_{jk}(P_{k_{3}}\psi_{3},P_{k_{4}}\psi_{4})))]\nonumber\\
&=\nabla_{x,t}P_{[-5,5]}Q_{<10}[P_{0}\psi_{0}\nabla^{-1}P_{r_{1}}(R_{0}P_{k_{1}}Q_{<k_{1}}\psi_{1}\label{eq:decomp201}\\&\hspace{3cm}\times\nabla^{-1}P_{r_{2}}(P_{k_{2}}\psi_{2}\nabla^{-1}P_{r_{3}}Q_{jk}(P_{k_{3}}\psi_{3},P_{k_{4}}\psi_{4})))]\nonumber\\
&+\nabla_{x,t}P_{[-5,5]}Q_{<10}[P_{0}\psi_{0}\nabla^{-1}P_{r_{1}}(R_{0}P_{k_{1}}Q_{[k_{1},\max\{k_{1,2,3,4}\}+O(1)]}\psi_{1}\label{eq:decomp202} \\&\hspace{3cm}\times \nabla^{-1}P_{r_{2}}(P_{k_{2}}\psi_{2}\nabla^{-1}P_{r_{3}}Q_{jk}(P_{k_{3}}\psi_{3},P_{k_{4}}\psi_{4})))]\nonumber\\
&+\nabla_{x,t}P_{[-5,5]}Q_{<10}[P_{0}\psi_{0}\nabla^{-1}P_{r_{1}}(R_{0}P_{k_{1}}Q_{\gg\max\{k_{1,2,3,4}\}}\psi_{1}\label{eq:decomp203}\\&\hspace{3cm}\times\nabla^{-1}P_{r_{2}}(P_{k_{2}}\psi_{2}\nabla^{-1}P_{r_{3}}Q_{jk}(P_{k_{3}}\psi_{3},P_{k_{4}}\psi_{4})))]\nonumber
\end{align}
To estimate the first expression \eqref{eq:decomp201} on the right, we exploit the fact that we control sharp Strichartz norms, in addition to the basic null-form bilinear estimate controlling $Q_{jk}(P_{k_{3}}\psi_{3},P_{k_{4}}\psi_{4})$.
The key is the fact that we have the almost sharp Klainerman-Tataru built into $S$. To see this, consider the most difficult case, a high-low frequency cascade corresponding to $r_{1}\ll k_{1}\ll k_{2}$. We estimate the expression by starting from the inside:
\begin{align*}
&\|\nabla^{-1}P_{r_{2}}(P_{k_{2}}\psi_{2}\nabla^{-1}P_{r_{3}}Q_{jk}(P_{k_{3}}\psi_{3},P_{k_{4}}\psi_{4}))\|_{L_{t}^{\frac{4}{3}}L_{x}^{2}}\\
&=2^{-r_{2}}\sum_{c_{1,2}\in R_{k_{2}, r_{2}}, \text{dist}(c_{1}, -c_{2})\lesssim r_{2}}\|\nabla^{-1}P_{r_{2}}(P_{k_{2}, c_{1}}\psi_{2}\nabla^{-1}P_{r_{3}, c_{2}}Q_{jk}(P_{k_{3}}\psi_{3},P_{k_{4}}\psi_{4}))\|_{L_{t}^{\frac{4}{3}}L_{x}^{2}}\\
&\lesssim 2^{-r_{2}}(\sum_{c_{1}\in  R_{k_{2}, r_{2}}}\|P_{k_{2}, c_{1}}\psi_{2}\|_{L_{t}^{4}L_{x}^{\infty}}^{2})^{\frac{1}{2}}\|\nabla^{-1}P_{r_{3}, c_{2}}Q_{jk}(P_{k_{3}}\psi_{3},P_{k_{4}}\psi_{4})\|_{L_{t,x}^{2}}
\end{align*}
Here we have used Cauchy-Schwarz and Plancherel's theorem. Then using the definition of $\|.\|_{S}$, we can bound this by
\begin{align*}
& 2^{-r_2}(\sum_{c_{1}\in  R_{k_{2}, r_{2}}}\|P_{k_{2}, c_{1}}\psi_{2}\|_{L_{t}^{4}L_{x}^{\infty}}^{2})^{\frac{1}{2}}\|\nabla^{-1}P_{r_{3}, c_{2}}Q_{jk}(P_{k_{3}}\psi_{3},P_{k_{4}}\psi_{4})\|_{L_{t,x}^{2}}\\
&\lesssim
2^{-r_{2}}2^{\frac{k_{2}}{4}}2^{\frac{r_{2}-k_{2}}{2+}}2^{\frac{\min\{r_{3},
k_{3,4}\}-\max\{r_{3},
k_{3,4}\}}{2}}\prod_{j=2}^{4}\|P_{k_{j}}\psi_{2}\|_{S[k_{j}]}
\end{align*}
Turning to the full expression further above, we then get for the contribution of this term to the hyperbolic part of the output
\begin{align*}
&\|\nabla_{x,t}Q_{<10}[P_{0}\psi_{0}\nabla^{-1}P_{r_{1}}(R_{0}Q_{<k_{1}}P_{k_{1}}\psi_{1} \\&\hspace{3cm}\times \nabla^{-1}P_{r_{2}}(P_{k_{2}}\psi_{2}\nabla^{-1}P_{r_{3}}Q_{jk}(P_{k_{3}}\psi_{3},P_{k_{4}}\psi_{4})))]\|_{L_{t}^{1}\dot{H}_{x}^{-1}}\\
&\lesssim \|P_{0}\psi_{0}\|_{L_{t}^{\infty}L_{x}^{2}}\|\nabla^{-1}P_{r_{1}}(R_{0}Q_{<k_{1}}P_{k_{1}}\psi_{1} \nabla^{-1}P_{r_{2}}(P_{k_{2}}\psi_{2}\nabla^{-1}P_{r_{3}}Q_{jk}(P_{k_{3}}\psi_{3},P_{k_{4}}\psi_{4})))\|_{L_{t}^{1}L_{x}^{\infty}}\\
&\lesssim \|P_{0}\psi_{0}\|_{L_{t}^{\infty}}\\&\times\|\nabla^{-1}P_{r_{1}}(\sum_{c_{1,2}\in R_{k_{1}, r_{1}},\,\text{dist}(c_{1},-c_{2})\lesssim 2^{r_{1}}}R_{0}Q_{<k_{1}}P_{k_{1}, c_{1}}\psi_{1} \nabla^{-1}P_{r_{2}, c_{2}}(P_{k_{2}}\psi_{2}\nabla^{-1}P_{r_{3}}Q_{jk}(P_{k_{3}}\psi_{3},P_{k_{4}}\psi_{4})))\|_{L_{t}^{1}L_{x}^{\infty}}
\end{align*}
We intend to substitute the intermediate bound from above for
\[
\|\nabla^{-1}P_{r_{2}}(P_{k_{2}}\psi_{2}\nabla^{-1}P_{r_{3}}Q_{jk}(P_{k_{3}}\psi_{3},P_{k_{4}}\psi_{4}))\|_{L_{t}^{\frac{4}{3}}L_{x}^{2}},
\]
where we can exploit that, by Minkowski's and Plancherel's inequality, we have
\[
(\sum_{c_{1}\in  R_{r_{2}, c_{2}}}\| P_{r_{2}, c_{2}}F\|_{L_{t}^{\frac{4}{3}}L_{x}^{2}}^{2})^{\frac{1}{2}}\lesssim \|P_{r_{2}}F\|_{L_{t}^{\frac{4}{3}}L_{x}^{2}}.
\]
Thus we can estimate, using Cauchy-Schwarz, Bernstein's inequality
and the preceding observation
\begin{align*}
&\|\nabla^{-1}P_{r_{1}}(\sum_{c_{1,2}\in R_{k_{1}, r_{1}},\,\text{dist}(c_{1},-c_{2})\lesssim 2^{r_{1}}}R_{0}Q_{<k_{1}}P_{k_{1}, c_{1}}\psi_{1} \\&\hspace{4cm}\times\nabla^{-1}P_{r_{2}, c_{2}}(P_{k_{2}}\psi_{2}\nabla^{-1}P_{r_{3}}Q_{jk}(P_{k_{3}}\psi_{3},P_{k_{4}}\psi_{4})))\|_{L_{t}^{1}L_{x}^{\infty}}\\
&\lesssim (\sum_{c_{1}\in  R_{k_{1}, r_{1}}}\|P_{k_{1}, c_{1}}\psi_{1}\|_{L_{t}^{4}L_{x}^{\infty}}^{2})^{\frac{1}{2}} \|\nabla^{-1}P_{r_{2}, c_{2}}(P_{k_{2}}\psi_{2}\nabla^{-1}P_{r_{3}}Q_{jk}(P_{k_{3}}\psi_{3},P_{k_{4}}\psi_{4})))\|_{L_{t}^{\frac{4}{3}}L_{x}^{2}}\\
&\lesssim
2^{\frac{r_{1}-k_{1}}{2+}}2^{\frac{3}{4}k_{1}}2^{-r_{2}}2^{\frac{k_{2}}{4}}2^{\frac{r_{2}-k_{2}}{2+}}2^{\frac{\min\{r_{3},
k_{3,4}\}-\max\{r_{3},
k_{3,4}\}}{2}}\prod_{j=1}^{4}\|P_{k_{j}}\psi_{2}\|_{S[k_{j}]}
\end{align*}
But by our assumption $r_{1}\ll k_{1}\ll k_{2}$ we have $r_{2}=k_{1}+O(1)$, whence we can replace the above bound by
\begin{align*}
&\|\nabla_{x,t}Q_{<10}[P_{0}\psi_{0}\nabla^{-1}P_{r_{1}}(R_{0}Q_{<k_{1}}P_{k_{1}}\psi_{1} \\&\hspace{3cm}\times \nabla^{-1}P_{r_{2}}(P_{k_{2}}\psi_{2}\nabla^{-1}P_{r_{3}}Q_{jk}(P_{k_{3}}\psi_{3},P_{k_{4}}\psi_{4})))]\|_{L_{t}^{1}\dot{H}_{x}^{-1}}\\
&\lesssim
2^{\frac{r_{1}-k_{1}}{2+}}2^{\frac{r_{2}-k_{2}}{4+}}2^{\frac{\min\{r_{3},
k_{3,4}\}-\max\{r_{3},
k_{3,4}\}}{2}}\prod_{j=0}^{4}\|P_{k_{j}}\psi_{2}\|_{S[k_{j}]},
\end{align*}
and this is again enough to yield the statement of the lemma (here with $\delta=\frac{1}{4+}$). The remaining frequency interactions can be handled similarly.
\\

Next, consider the second term \eqref{eq:decomp202} above, i.e.,
\begin{align*}
&\nabla_{x,t}P_{[-5,5]}Q_{<10}[P_{0}\psi_{0}\nabla^{-1}P_{r_{1}}(R_{0}P_{k_{1}}Q_{[k_{1},\max\{k_{1,2,3,4}\}+O(1)]}\psi_{1}\\&\hspace{6cm}\times\nabla^{-1}P_{r_{2}}(P_{k_{2}}\psi_{2}\nabla^{-1}P_{r_{3}}Q_{jk}(P_{k_{3}}\psi_{3},P_{k_{4}}\psi_{4})))]
\end{align*}
This is much simpler: we get
\begin{align*}
&\|\nabla_{x,t}P_{[-5,5]}Q_{<10}[P_{0}\psi_{0}\nabla^{-1}P_{r_{1}}(R_{0}P_{k_{1}}Q_{[k_{1},\max\{k_{1,2,3,4}\}+O(1)]}\psi_{1} \\&\hspace{6cm}\times\nabla^{-1}P_{r_{2}}(P_{k_{2}}\psi_{2}\nabla^{-1}P_{r_{3}}Q_{jk}(P_{k_{3}}\psi_{3},P_{k_{4}}\psi_{4})))]\|_{L_{t}^{1}\dot{H}^{-1}}\\
&\lesssim \|P_{0}\psi_{0}\|_{L_{t}^{\infty}L_{x}^{2}}\nabla^{-1}P_{r_{1}}(R_{0}P_{k_{1}}Q_{[k_{1},\max\{k_{1,2,3,4}\}+O(1)]}\psi_{1}\\&\hspace{5cm}\times\nabla^{-1}P_{r_{2}}(P_{k_{2}}\psi_{2}\nabla^{-1}P_{r_{3}}Q_{jk}(P_{k_{3}}\psi_{3},P_{k_{4}}\psi_{4})))\|_{L_{t}^{1}L_{x}^{\infty}}
\end{align*}
For definitiveness' sake, we again assume that $r_{1}\ll k_{1}\ll k_{2}$, the remaining cases being similar. Then we get
\begin{align*}
&\nabla^{-1}P_{r_{1}}(R_{0}P_{k_{1}}Q_{[k_{1},\max\{k_{1,2,3,4}\}+O(1)]}\psi_{1} \nabla^{-1}P_{r_{2}}(P_{k_{2}}\psi_{2}\nabla^{-1}P_{r_{3}}Q_{jk}(P_{k_{3}}\psi_{3},P_{k_{4}}\psi_{4})))\|_{L_{t}^{1}L_{x}^{\infty}}\\
&\lesssim 2^{r_{1}-k_{1}}\|R_{0}P_{k_{1}}Q_{[k_{1},\max\{k_{1,2,3,4}\}+O(1)]}\psi_{1}\|_{L_{t,x}^{2}}\|\nabla^{-1}P_{r_{2}}(P_{k_{2}}\psi_{2}\nabla^{-1}P_{r_{3}}Q_{jk}(P_{k_{3}}\psi_{3},P_{k_{4}}\psi_{4})\|_{L_{t}^{2}\dot{H}^{\frac{1}{2}}}\\
&\lesssim
2^{r_{1}-k_{1}}2^{\epsilon(\max\{k_{1,2,3,4}\}-k_{1})}2^{\frac{r_{2}-\max\{k_{3,4}\}}{2}}\prod_{j=0}^{4}\|P_{k_{j}}\psi_{j}\|_{S[k_{j}]}
\end{align*}
This corresponds to a bound as in the lemma with
$\delta=\frac{1}{2}-\epsilon$, where we recall $\epsilon$ is as in
the definition of $\|\cdot \|_{S[k]}$. The remaining frequency
interactions for this term are treated similarly.
\\

Finally, consider the last term above
\begin{align*}
&\nabla_{x,t}P_{[-5,5]}Q_{<10}[P_{0}\psi_{0}\nabla^{-1}P_{r_{1}}(R_{0}P_{k_{1}}Q_{\gg\max\{k_{1,2,3,4}\}}\psi_{1}
 \\&\hspace{5cm}\times\nabla^{-1}P_{r_{2}}(P_{k_{2}}\psi_{2}\nabla^{-1}P_{r_{3}}Q_{jk}(P_{k_{3}}\psi_{3},P_{k_{4}}\psi_{4})))]
\end{align*}
Here we again need to compensate for the losses coming from estimating $R_{0}P_{k_{1}}Q_{\gg\max\{k_{1,2,3,4}\}}\psi_{1} $.
Freeze its modulation to dyadic size $2^{l}$. Then either at least one other input has at least comparable modulation,
or else the output has modulation $\sim 2^{l}$ (in which case necessarily $l<O(1)$. In the latter case, one then estimates
(where $l\gg \max\{k_{1,2,3,4}\}$ and we assume all other inputs to be at much lower modulation)
\begin{align*}
&\|\nabla_{x,t}P_{[-5,5]}Q_{<10}[P_{0}\psi_{0}\nabla^{-1}P_{r_{1}}(R_{0}P_{k_{1}}Q_{l}\psi_{1} \\&\hspace{5cm}\times\nabla^{-1}P_{r_{2}}(P_{k_{2}}\psi_{2}\nabla^{-1}P_{r_{3}}Q_{jk}(P_{k_{3}}\psi_{3},P_{k_{4}}\psi_{4})))]\|_{N[0]}\\
&=\|\nabla_{x,t}P_{[-5,5]}Q_{l}[P_{0}\psi_{0}\nabla^{-1}P_{r_{1}}(R_{0}P_{k_{1}}Q_{l}\psi_{1} \\&\hspace{5cm}\times \nabla^{-1}P_{r_{2}}(P_{k_{2}}\psi_{2}\nabla^{-1}P_{r_{3}}Q_{jk}(P_{k_{3}}\psi_{3},P_{k_{4}}\psi_{4})))]\|_{N[0]}\\
&\leq \|\nabla_{x,t}P_{[-5,5]}Q_{l}[P_{0}\psi_{0}\nabla^{-1}P_{r_{1}}(R_{0}P_{k_{1}}Q_{l}\psi_{1}  \\&\hspace{5cm}\times\nabla^{-1}P_{r_{2}}(P_{k_{2}}\psi_{2}\nabla^{-1}P_{r_{3}}Q_{jk}(P_{k_{3}}\psi_{3},P_{k_{4}}\psi_{4})))]\|_{\dot{X}_{0}^{-1,-\frac{1}{2},1}}\\
&\lesssim 2^{-\frac{l}{2}}\|P_{0}\psi_{0}\|_{L_{t}^{\infty}L_{x}^{2}}\|\nabla^{-1}P_{r_{1}}(R_{0}P_{k_{1}}Q_{l}\psi_{1}  \\&\hspace{5cm}\times\nabla^{-1}P_{r_{2}}(P_{k_{2}}\psi_{2}\nabla^{-1}P_{r_{3}}Q_{jk}(P_{k_{3}}\psi_{3},P_{k_{4}}\psi_{4})))\|_{L_{t}^{2}L_{x}^{\infty}}
\end{align*}
Here the second factor above is estimated by
\begin{align*}
&\|\nabla^{-1}P_{r_{1}}(R_{0}P_{k_{1}}Q_{l}\psi_{1} \nabla^{-1}P_{r_{2}}(P_{k_{2}}\psi_{2}\nabla^{-1}P_{r_{3}}Q_{jk}(P_{k_{3}}\psi_{3},P_{k_{4}}\psi_{4})))\|_{L_{t}^{2}L_{x}^{\infty}}\\
&\lesssim 2^{\frac{\min\{r_{1,2}, k_{1}\}}{2}}2^{\frac{\min\{r_{1}, r_{2}, k_{1}\}-\max\{r_{1}, r_{2}, k_{1}\}}{2}}2^{\frac{\min\{r_{2,3}, k_{2,3,4}\}-\max\{r_{2,3}, k_{2,3,4}\}}{2}}2^{\epsilon(l-k_{1})}\prod_{j=1}^{4}\|P_{k_{j}}\psi_{j}\|_{S[k_{j}]}
\end{align*}
Inserting this bound into the last inequality but one and summing over $l\gg\max\{k_{1,2,3,4}\}$ results in a bound as
in the lemma with $\delta=\frac{1}{2}-\epsilon$.
\\
The case when at least one further input has at least modulation at least comparable to $2^{l}$ is similar, one places the output into $L_{t}^{1}\dot{H}^{-1}$.
\\
This completes the proof of the first inequality of the lemma. The remaining ones are treated by an identical procedure.
\end{proof}

In a similar vein, one has estimates controlling the second kind of quintilinear term. We state the

\begin{lemma}\label{quintilinear2} For the second type of quintilinear null-form, we have the following estimates for suitable $\delta>0$:
\begin{align*}
&\|\nabla_{x,t}[\big(P_{0}[\nabla^{-1}(P_{k_{1}}\psi_{1} P_{s_{1}}\nabla^{-1}Q_{\nu j}(P_{k_{2}}\psi_{2}, P_{k_{3}}\psi_{3}))] P_{r_{1}}\nabla^{-1}IQ_{\mu j}(P_{k_{4}}\psi_{4}, P_{k_{5}}\psi_{5})\big)\|_{N[0]}
\\&\lesssim 2^{\delta[\min\{0, k_{1,2,3}, s_{1}\}-\max\{0, k_{1,2,3}, s_{1}\}]}2^{\delta[\min\{r_{1}, k_{4,5}\}-\max\{r_{1}, k_{4,5}\}]}\prod_{j=1}^{5}\|P_{k_{j}}\psi_{j}\|_{S[k_{j}]},\quad r_{1}<-10
\end{align*}
\begin{align*}
&\|\nabla_{x,t}P_{0}[\big(P_{s_{1}}[\nabla^{-1}(P_{k_{1}}\psi_{1} P_{s_{2}}\nabla^{-1}Q_{\nu j}(P_{k_{2}}\psi_{2}, P_{k_{3}}\psi_{3}))] P_{r_{1}}\nabla^{-1}IQ_{\mu j}(P_{k_{4}}\psi_{4}, P_{k_{5}}\psi_{5})\big)\|_{N[0]}
\\&\lesssim 2^{\delta s_{1}}2^{\delta[\min\{s_{1,2}, k_{1,2,3}\}-\max\{s{1,2}, k_{1,2,3}\}]}2^{\delta[\min\{r_{1}, k_{4,5}\}-\max\{r_{1}, k_{4,5}\}]}\prod_{j=1}^{5}\|P_{k_{j}}\psi_{j}\|_{S[k_{j}]},\quad r_{1}\in [-10, 10]
\end{align*}
\begin{align*}
&\|\nabla_{x,t}P_{0}[\big(P_{s_{1}}[\nabla^{-1}(P_{k_{1}}\psi_{1} P_{s_{2}}\nabla^{-1}Q_{\nu j}(P_{k_{2}}\psi_{2}, P_{k_{3}}\psi_{3}))] P_{r_{1}}\nabla^{-1}IQ_{\mu j}(P_{k_{4}}\psi_{4}, P_{k_{5}}\psi_{5})\big)\|_{N[0]}
\\&\lesssim 2^{-\delta s_{1}}2^{\delta[\min\{s_{1,2}, k_{1,2,3}\}-\max\{s_{1,2}, k_{1,2,3}\}]}2^{\delta[\min\{r_{1}, k_{4,5}\}-\max\{r_{1}, k_{4,5}\}]}\prod_{j=1}^{5}\|P_{k_{j}}\psi_{j}\|_{S[k_{j}]},\quad r_{1}>10
\end{align*}
\end{lemma}

\begin{proof} We verify this again for the first inequality above, the other ones following a similar pattern.
As usual, we distinguish between elliptic and hyperbolic output components:
\\

{\it{(i): Output in elliptic regime.}} This is the expression
\begin{align*}
&\nabla_{x,t}P_{[-5,5]}Q_{>10}[\big(P_{0}[\nabla^{-1}(P_{k_{1}}\psi_{1} P_{s_{1}}\nabla^{-1}Q_{\nu j}(P_{k_{2}}\psi_{2}, P_{k_{3}}\psi_{3}))] \\&\hspace{5cm}\times P_{r_{1}}\nabla^{-1}IQ_{\mu j}(P_{k_{4}}\psi_{4}, P_{k_{5}}\psi_{5})\big)
\end{align*}
As usual the only slight complication arises due to the fact that w may have $\nu=0$. Freeze the modulation of the output
 to dyadic size $2^{l}$, $l>10$. Then one re-iterates the same steps as in the preceding proof:
\\

{\it{(i1): $\max\{k_{1,2,3}\}\ll l$, time derivative falls on term
with modulation $<2^{l-100}$}}. In this case at least one additional
input (which is not hit by a time derivative) has modulation
$>2^{l-10}$. For example, assume this is
$P_{k_{1}}\psi_{1}=P_{k_{1}}Q_{>l-10}\psi_{1}$, the other cases
being treated similarly. Then assuming a high-low scenario, say,
i.e., $k_{1}\gg 1$, we have (using Bernstein's inequality)
\begin{align*}
&\|P_{0}[\nabla^{-1}(P_{k_{1}}Q_{>l-10}\psi_{1} P_{s_{1}}\nabla^{-1}Q_{\nu j}(P_{k_{2}}\psi_{2}, P_{k_{3}}\psi_{3}))]\|_{L_{t}^{1}L_{x}^{2}}\\&
\lesssim \|P_{k_{1}}Q_{>l-10}\psi_{1}\|_{L_{t,x}^{2}}\|P_{s_{1}}\nabla^{-1}Q_{\nu j}(P_{k_{2}}\psi_{2}, P_{k_{3}}\psi_{3})\|_{L_{t,x}^{2}}\\
&\lesssim 2^{-\frac{3k_{1}}{2}}2^{(1-\epsilon)(k_{1}-l)}2^{\epsilon(l-k_{2})}\prod_{j=1}^{3}\|P_{k_{j}}\psi_{j}\|_{S[k_{j}]}
\end{align*}
Substituting this into the full expression, we obtain for the output the bound
\begin{align*}
&\|\nabla_{x,t}P_{[-5,5]}Q_{l}[\big(P_{0}[\nabla^{-1}(P_{k_{1}}Q_{>l-10}\psi_{1} P_{s_{1}}\nabla^{-1}Q_{\nu j}(P_{k_{2}}\psi_{2}, P_{k_{3}}\psi_{3}))] \\&\hspace{7cm}\times P_{r_{1}}\nabla^{-1}IQ_{\mu j}(P_{k_{4}}\psi_{4}, P_{k_{5}}\psi_{5})\big)\|_{\dot{X}_{0}^{-\frac{1}{2}+\epsilon, -1-\epsilon, 1}}\\
&\lesssim 2^{-\epsilon l}\|P_{0}[\nabla^{-1}(P_{k_{1}}Q_{>l-10}\psi_{1} P_{s_{1}}\nabla^{-1}Q_{\nu j}(P_{k_{2}}\psi_{2}, P_{k_{3}}\psi_{3}))]\|_{L_{t,x}^{2}}\\&\hspace{7cm}\times\|P_{r_{1}}\nabla^{-1}IQ_{\mu j}(P_{k_{4}}\psi_{4}, P_{k_{5}}\psi_{5})\|_{L_{t,x}^{\infty}}\\
&\lesssim 2^{(\frac{1}{2}-\epsilon)l}\|P_{0}[\nabla^{-1}(P_{k_{1}}Q_{>l-10}\psi_{1} P_{s_{1}}\nabla^{-1}Q_{\nu j}(P_{k_{2}}\psi_{2}, P_{k_{3}}\psi_{3}))]\|_{L_{t}^{1}L_{x}^{2}}\\&\hspace{7cm}\times\|P_{r_{1}}\nabla^{-1}IQ_{\mu j}(P_{k_{4}}\psi_{4}, P_{k_{5}}\psi_{5})\|_{L_{t,x}^{\infty}}\\
&\lesssim  2^{(\frac{1}{2}-\epsilon)l}2^{-\frac{3k_{1}}{2}}2^{(1-\epsilon)(k_{1}-l)}2^{\epsilon(l-k_{2})}2^{\frac{\min\{r_{1}, k_{4,5}\}-\max\{r_{1}, k_{4,5}\}}{2}}\prod_{j=1}^{5}\|P_{k_{j}}\psi_{j}\|_{S[k_{j}]}
\end{align*}
Summing over $l\gg \max\{k_{1,2,3}\}$, the desired inequality of the first type of the lemma follows in this case. The remaining frequency interactions within
\[
P_{0}[\nabla^{-1}(P_{k_{1}}Q_{>l-10}\psi_{1} P_{s_{1}}\nabla^{-1}Q_{\nu j}(P_{k_{2}}\psi_{2}, P_{k_{3}}\psi_{3}))]
\]
are handled similarly.
\\

{\it{(i2):  $\max\{k_{1,2,3}\}\ll l$, time derivative falls on term
with modulation $\sim 2^{l}$}}.  In this case, we place the time
derivative term into $L_{t,x}^{2}$, and are guaranteed gains in the
maximal occurring frequency: for example, consider the term (arising
upon unraveling the inner $Q_{\nu j}$ null-structure with $\nu=0$)
\[
P_{0}[\nabla^{-1}(P_{k_{1}}\psi_{1} P_{s_{1}}\nabla^{-1}(P_{k_{2}}Q_{l+O(1)}R_{0}\psi_{2} P_{k_{3}}\psi_{3}))]
\]
In the high-high case $k_{2}\gg s_{1}$, one can then estimate
\begin{align*}
&\|P_{0}[\nabla^{-1}(P_{k_{1}}\psi_{1}
P_{s_{1}}\nabla^{-1}(P_{k_{2}}Q_{l+O(1)}R_{0}\psi_{2}
P_{k_{3}}\psi_{3}))]\|_{L_{t,x}^{2}}\\&\lesssim 2^{\frac{\min\{0,
k_{1}, s_{1}\}-\max\{0, k_{1}, s_{1}\}}{2}}
 2^{\frac{\min\{s_{1}, k_{2}, k_{3}\}-\max\{s_{1}, k_{2}, k_{3}\}}{2}}\\&\hspace{2cm}\times\|P_{k_{1}}\psi_{1}\|_{L_{t}^{\infty}L_{x}^{2}}\|P_{k_{2}}Q_{l+O(1)}R_{0}\psi_{2}\|_{L_{t}^{2}\dot{H}^{\frac{1}{2}}}\|P_{k_{3}}\psi_{3}\|_{L_{t}^{\infty}L_{x}^{2}}\\
\end{align*}
From here one estimates the full expression by
\begin{align*}
&\|\nabla_{x,t}P_{[-5,5]}Q_{l}[\big(P_{0}[\nabla^{-1}(P_{k_{1}}\psi_{1} P_{s_{1}}\nabla^{-1}(P_{k_{2}}Q_{l+O(1)}R_{0}\psi_{2} P_{k_{3}}\psi_{3}))] \\&\hspace{6cm}\times P_{r_{1}}\nabla^{-1}IQ_{\mu j}(P_{k_{4}}\psi_{4}, P_{k_{5}}\psi_{5})\big)\|_{\dot{X}_{0}^{-\frac{1}{2}+\epsilon, -1-\epsilon, 1}}\\
&\lesssim 2^{-\epsilon l}\|P_{0}[\nabla^{-1}(P_{k_{1}}\psi_{1} P_{s_{1}}\nabla^{-1}(P_{k_{2}}Q_{l+O(1)}R_{0}\psi_{2} P_{k_{3}}\psi_{3}))]\|_{L_{t,x}^{2}}\\&\hspace{6cm}\times\|P_{r_{1}}\nabla^{-1}IQ_{\mu j}(P_{k_{4}}\psi_{4}, P_{k_{5}}\psi_{5})\|_{L_{t,x}^{\infty}}\\
&\lesssim 2^{-\epsilon l}2^{\epsilon (l-k_{2})}2^{\frac{\min\{0, k_{1}, s_{1}\}-\max\{0, k_{1}, s_{1}\}}{2}} 2^{\frac{\min\{s_{1}, k_{2}, k_{3}\}-\max\{s_{1}, k_{2}, k_{3}\}}{2}}\\&\hspace{1.5cm}\times\|P_{k_{1}}\psi_{1}\|_{L_{t}^{\infty}L_{x}^{2}}\|P_{k_{2}}Q_{l+O(1)}\psi_{2}\|_{\dot{X}_{k_{2}}^{-\frac{1}{2}+\epsilon, 1-\epsilon, 1}}\|P_{k_{3}}\psi_{3}\|_{L_{t}^{\infty}L_{x}^{2}}\\&\hspace{3cm}\times 2^{\min\{r_{1}, k_{4,5}\}-\max\{r_{1}, k_{4,5}\}}\prod_{j=4,5}\|P_{k_{j}}\psi_{j}\|_{S[k_{j}]}
\end{align*}
One may sum here to obtain a bound of the type as in the first inequality of the lemma, with $\delta=\frac{1}{2}$. \\
The remaining frequency interactions within
\[
P_{0}[\nabla^{-1}(P_{k_{1}}\psi_{1} P_{s_{1}}\nabla^{-1}(P_{k_{2}}Q_{l+O(1)}R_{0}\psi_{2} P_{k_{3}}\psi_{3}))]
\]
are again handled similarly.
\\

{\it{(i3):  $\max\{k_{1,2,3}\}\ll l$, time derivative falls on term
with modulation $\gg 2^{l}$}}. In this case at least one additional
term has at least comparable modulation, and one argues as in case
(i1).
\\

{\it{(i4): $\max\{k_{1,2,3}\}>l+O(1)$, time derivative falls on term with modulation $<\max\{k_{1,2,3}\}+O(1)$.}}
Here the losses coming from the time derivative are easily counteracted by the gains in the large frequencies:  first, one reduces the inputs $P_{k_{2,3}}\psi_{2,3}$ to the elliptic regimes. To do so, note that we have
\begin{align*}
&\|P_{0}[\nabla^{-1}(P_{k_{1}}\psi_{1} P_{s_{1}}\nabla^{-1}(P_{k_{2}}Q_{[k_{2}, \max\{k_{1,2,3}\}}R_{0}\psi_{2} P_{k_{3}}\psi_{3}))]\|_{L_{t,x}^{2}}\\&\lesssim 2^{\frac{\min\{0, s_{1}, k_{1,2,3}\}-\max\{0, s_{1}, k_{1,2,3}\}}{2}}2^{\epsilon(\max\{k_{1,2,3}\}-k_{2})}\prod_{j=1}^{3}\|P_{k_{j}}\psi_{j}\|_{S[k_{j}]},
\end{align*}
and inserting this into the full expression is easily seen to yield the desired inequality. Hence we have reduced this case to the expression
\begin{align*}
&\nabla_{x,t}P_{[-5,5]}Q_{l}[\big(P_{0}[\nabla^{-1}(P_{k_{1}}\psi_{1} P_{s_{1}}\nabla^{-1}Q_{\nu j}(P_{k_{2}}Q_{<k_{2}}\psi_{2} P_{k_{3}}Q_{<k_{3}}\psi_{3}))] \\&\hspace{7cm}\times P_{r_{1}}\nabla^{-1}IQ_{\mu j}(P_{k_{4}}\psi_{4}, P_{k_{5}}\psi_{5})\big)
\end{align*}
Of course in the present case at least one of $k_{2,3}>l+O(1)$. Assume that we have a high-high-low type situation in
\[
P_{s_{1}}\nabla^{-1}Q_{\nu j}(P_{k_{2}}Q_{<k_{2}}\psi_{2} P_{k_{3}}Q_{<k_{3}}\psi_{3}),
\]
i.e., $s_{1}\ll k_{2}$, this being the most delicate case. We
distinguish between two cases:
\\

{\it{(a): modulation of $P_{s_{1}}\ldots$ is less than $2^{l+O(1)}$.}} In this case, we may "pull out" a (time)- derivative from the $Q_{\nu j}$-null-form, using the simple identity
\[
R_{\nu}\psi^{1}R_{j}\psi^{2}-R_{j}\psi^{1}R_{\nu}\psi^{2}=\partial_{\nu}[\nabla^{-1}\psi^{1}R_{j}\psi_{2}]
-\partial_{j}[\nabla^{-1}\psi^{1}R_{\nu}\psi_{2}]
\]
Hence in this case we can estimate
\begin{align*}
&\|\nabla_{x,t}P_{[-5,5]}Q_{l}[\big(P_{0}[\nabla^{-1}(P_{k_{1}}\psi_{1} P_{s_{1}}Q_{<l+O(1)}\nabla^{-1}Q_{\nu j}(P_{k_{2}}Q_{<k_{2}}\psi_{2} P_{k_{3}}Q_{<k_{3}}\psi_{3}))] \\&\hspace{8cm}\times P_{r_{1}}\nabla^{-1}IQ_{\mu j}(P_{k_{4}}\psi_{4}, P_{k_{5}}\psi_{5})\big)\|_{\dot{X}_{0}^{-\frac{1}{2}+\epsilon, -1-\epsilon,1}}\\
&\lesssim 2^{-\epsilon l}\|P_{0}[\nabla^{-1}(P_{k_{1}}\psi_{1} P_{s_{1}}Q_{<l+O(1)}\nabla^{-1}Q_{\nu j}(P_{k_{2}}Q_{<k_{2}}\psi_{2} P_{k_{3}}Q_{<k_{3}}\psi_{3}))]\|_{L_{t}^{\infty}L_{x}^{2}}\\&\hspace{8cm}\times\|P_{r_{1}}\nabla^{-1}IQ_{\mu j}(P_{k_{4}}\psi_{4}, P_{k_{5}}\psi_{5})\|_{L_{t}^{2}L_{x}^{\infty}}\\
&\lesssim  2^{-\epsilon l}2^{l-k_{2}}2^{\frac{r_{1}}{2}}2^{\frac{\min\{r_{1}, k_{3,4}\}-\max\{r_{1}, k_{3,4}\}}{2}}2^{\min\{s_{1}, k_{1},0\}}\prod_{j=1}^{5}\|P_{k_{j}}\psi_{j}\|_{S[k_{j}]}
\end{align*}
One can sum over $l<\max\{k_{1,2,3}\}$ to get the desired first inequality of the lemma in the case at hand.
\\

{\it{(b): modulation of $P_{s_{1}}\ldots$ is $\gg 2^{l}$.}} In this case the modulation of the first input $P_{k_{1}}\psi_{1}$ needs to be comparable to that of
\[
P_{s_{1}}\nabla^{-1}Q_{\nu j}(P_{k_{2}}Q_{<k_{2}}\psi_{2} P_{k_{3}}Q_{<k_{3}}\psi_{3}))]
\]
Hence we can write this contribution as
\begin{align*}
&\sum_{l_1\gg l}\nabla_{x,t}P_{[-5,5]}Q_{l}[\big(P_{0}[\nabla^{-1}(P_{k_{1}}Q_{l_1+O(1)}\psi_{1} P_{s_{1}}Q_{l_1}\nabla^{-1}Q_{\nu j}(P_{k_{2}}Q_{<k_{2}}\psi_{2} P_{k_{3}}Q_{<k_{3}}\psi_{3}))]\\
&\hspace{8cm} \times P_{r_{1}}\nabla^{-1}IQ_{\mu j}(P_{k_{4}}\psi_{4}, P_{k_{5}}\psi_{5})\big)
\end{align*}
To estimate it, we use
\[
\|P_{s_{1}}Q_{l_1}\nabla^{-1}Q_{\nu j}(P_{k_{2}}Q_{<k_{2}}\psi_{2} P_{k_{3}}Q_{<k_{3}}\psi_{3})\|_{L_{t,x}^2}
\lesssim 2^{l_1-k_2}2^{-\frac{s_1}{2}}\prod_{i=2,3}\|P_{k_{i}}Q_{<k_{i}}\psi_{2}\|_{S[k_i]}
\]
We the insert this bound into the full expression. In case that $k_{1}>O(1)$, we can estimate
\begin{align*}
\|&\nabla_{x,t}P_{[-5,5]}Q_{l}[\big(P_{0}[\nabla^{-1}(P_{k_{1}}Q_{l_1+O(1)}\psi_{1} P_{s_{1}}Q_{l_1}\nabla^{-1}Q_{\nu j}(P_{k_{2}}Q_{<k_{2}}\psi_{2} P_{k_{3}}Q_{<k_{3}}\psi_{3}))]\\
&\hspace{6cm} \times P_{r_{1}}\nabla^{-1}IQ_{\mu j}(P_{k_{4}}\psi_{4}, P_{k_{5}}\psi_{5})\big)]
]\|_{\dot{X}_{0}^{-\frac{1}{2}+\epsilon, -1-\epsilon,1}}\\
&\lesssim 2^{-\epsilon l}2^{\frac{l}{2}}\|P_{k_{1}}Q_{l_1+O(1)}\psi_{1}\|_{L_{t,x}^{2}}\|P_{s_{1}}Q_{l_1}\nabla^{-1}Q_{\nu j}(P_{k_{2}}Q_{<k_{2}}\psi_{2} P_{k_{3}}Q_{<k_{3}}\psi_{3})\|_{L_{t,x}^{2}}
\\
&\hspace{6cm} \times \|P_{r_{1}}\nabla^{-1}IQ_{\mu j}(P_{k_{4}}\psi_{4}, P_{k_{5}}\psi_{5})\|_{L_{t,x}^{\infty}}\\
\end{align*}
In case that $l_1<k_2+O(1)$, we can bound this by
\begin{align*}
&2^{-\epsilon l}2^{\frac{l}{2}}\|P_{k_{1}}Q_{l_1+O(1)}\psi_{1}\|_{L_{t,x}^{2}}\|P_{s_{1}}Q_{l_1}\nabla^{-1}Q_{\nu j}(P_{k_{2}}Q_{<k_{2}}\psi_{2} P_{k_{3}}Q_{<k_{3}}\psi_{3})\|_{L_{t,x}^{2}}
\\
&\hspace{6cm} \times \|P_{r_{1}}\nabla^{-1}IQ_{\mu j}(P_{k_{4}}\psi_{4}, P_{k_{5}}\psi_{5})\|_{L_{t,x}^{\infty}}\\
&\lesssim 2^{-\epsilon l}2^{\frac{l}{2}}2^{-\frac{l_1}{2}}\|P_{k_1}\psi_1\|_{S[k_1]}2^{l_1-k_2}2^{-\frac{s_1}{2}}\prod_{i=2,3}\|P_{k_{i}}Q_{<k_{i}}\psi_{2}\|_{S[k_i]}\\
&\hspace{6cm} \times \|P_{r_{1}}\nabla^{-1}IQ_{\mu j}(P_{k_{4}}\psi_{4}, P_{k_{5}}\psi_{5})\|_{L_{t,x}^{\infty}}\\
&\lesssim 2^{-\epsilon l}2^{\frac{l}{2}}2^{-\frac{l_1}{2}}2^{l_1-k_2}2^{-\frac{s_1}{2}}2^{r_1}2^{\min\{k_3, k_4\}-\max\{k_3, k_4\}}\prod_{i=1}^{5}\|P_{k_i}\psi_i\|_{S[k_i]}
\end{align*}
Summing over $k_1+O(1)>l_1\gg l$ and then over $l>O(1)$ results in a bound as in the first inequality of the lemma with $\delta=\frac{1}{2}+\epsilon$.
\\
Next, still in the case $k_{1}>O(1)$, if $l_1<k_2+O(1)$, one proceeds as before but uses
\[
\|P_{k_{1}}Q_{l_1+O(1)}\psi_{1}\|_{L_{t,x}^{2}}\lesssim 2^{-\frac{k_1}{2}}2^{(1-\epsilon)(k_1-l_1)}\|P_{k_1}\psi\|_{S[k_1]}
\]
One obtains a final bound with the same $\delta=\frac{1}{2}+\epsilon$ as in the preceding case.
\\

In the case $k_1<O(1)$, one simply places $P_{k_{1}}Q_{l_1+O(1)}\psi_{1}$ into $L_t^2L_x^\infty$, thereby gaining an additional factor $2^{k_1}$. We omit the details.
\\
This concludes case (i), when the output is in the elliptic regime.
\\

{\it{(ii): Output in hyperbolic regime}}.
This is the expression
\[
\nabla_{x,t}P_{[-5,5]}Q_{<10}[\big(P_{0}[\nabla^{-1}(P_{k_{1}}\psi_{1} P_{s_{1}}\nabla^{-1}Q_{\nu j}(P_{k_{2}}\psi_{2}, P_{k_{3}}\psi_{3}))]P_{r_{1}}\nabla^{-1}IQ_{\mu j}(P_{k_{4}}\psi_{4}, P_{k_{5}}\psi_{5})\big)
\]
To treat it, we decompose
\begin{align*}
P_{0}[\nabla^{-1}(P_{k_{1}}\psi_{1} P_{s_{1}}\nabla^{-1}Q_{\nu j}(P_{k_{2}}\psi_{2}, P_{k_{3}}\psi_{3}))]
=&P_{0}[\nabla^{-1}(P_{k_{1}}\psi_{1} P_{s_{1}}\nabla^{-1}I^cQ_{\nu j}(P_{k_{2}}\psi_{2}, P_{k_{3}}\psi_{3}))]\\&
+P_{0}[\nabla^{-1}(P_{k_{1}}\psi_{1} P_{s_{1}}\nabla^{-1}IQ_{\nu j}(P_{k_{2}}\psi_{2}, P_{k_{3}}\psi_{3}))]
\end{align*}

{\it{(iia): contribution of the elliptic type term.}} This is the expression
\[
\nabla_{x,t}P_{[-5,5]}Q_{<10}[\big(P_{0}[\nabla^{-1}(P_{k_{1}}\psi_{1} P_{s_{1}}\nabla^{-1}I^cQ_{\nu j}(P_{k_{2}}\psi_{2}, P_{k_{3}}\psi_{3}))]P_{r_{1}}\nabla^{-1}IQ_{\mu j}(P_{k_{4}}\psi_{4}, P_{k_{5}}\psi_{5})\big)
\]
We shall treat the case $s_1\ll -10$, i.e., the case of a high-low
interaction within
\[
P_{0}[\nabla^{-1}(P_{k_{1}}\psi_{1} P_{s_{1}}\nabla^{-1}Q_{\nu j}(P_{k_{2}}\psi_{2}, P_{k_{3}}\psi_{3}))]
\]
The remaining cases are again more of the same. Now freeze the modulation of the expression
\[
P_{s_{1}}\nabla^{-1}I^cQ_{\nu j}(P_{k_{2}}\psi_{2}, P_{k_{3}}\psi_{3}))
\]
to size $2^{l}$, $l\gg s_1$. Then decompose the corresponding full
expression into the following:
\begin{align}
&\nabla_{x,t}P_{[-5,5]}Q_{<10}[\big(P_{0}[\nabla^{-1}(P_{k_{1}}\psi_{1} P_{s_{1}}Q_l\nabla^{-1}I^cQ_{\nu j}(P_{k_{2}}\psi_{2}, P_{k_{3}}\psi_{3}))]P_{r_{1}}\nabla^{-1}IQ_{\mu j}(P_{k_{4}}\psi_{4}, P_{k_{5}}\psi_{5})\big)\nonumber\\
&=\nabla_{x,t}P_{[-5,5]}Q_{<10}[\big(P_{0}[\nabla^{-1}(P_{k_{1}}Q_{>l-10}\psi_{1} P_{s_{1}}Q_l\nabla^{-1}I^cQ_{\nu j}(P_{k_{2}}\psi_{2}, P_{k_{3}}\psi_{3}))]\nonumber\\&\hspace{8cm}\times P_{r_{1}}\nabla^{-1}IQ_{\mu j}(P_{k_{4}}\psi_{4}, P_{k_{5}}\psi_{5})\big)\label{eq:decomp101}\\
&+\nabla_{x,t}P_{[-5,5]}Q_{<10}[\big(P_{0}[\nabla^{-1}(P_{k_{1}}Q_{<l-10}\psi_{1} P_{s_{1}}Q_l\nabla^{-1}I^cQ_{\nu j}(P_{k_{2}}\psi_{2}, P_{k_{3}}\psi_{3}))]\nonumber\\&\hspace{8cm}\times P_{r_{1}}\nabla^{-1}IQ_{\mu j}(P_{k_{4}}\psi_{4}, P_{k_{5}}\psi_{5})\big)\label{eq:decomp111}
\end{align}
The first term \eqref{eq:decomp101} on the right can then be estimated by
\begin{align*}
&\|\nabla_{x,t}P_{[-5,5]}Q_{<10}[\big(P_{0}[\nabla^{-1}(P_{k_{1}}Q_{>l-10}\psi_{1} P_{s_{1}}Q_l\nabla^{-1}I^cQ_{\nu j}(P_{k_{2}}\psi_{2}, P_{k_{3}}\psi_{3}))]\\&\hspace{8cm}\times P_{r_{1}}\nabla^{-1}IQ_{\mu j}(P_{k_{4}}\psi_{4}, P_{k_{5}}\psi_{5})\big)\|_{L_{t}^{1}\dot{H}^{-1}}\\
&\lesssim \|P_{k_{1}}Q_{>l-10}\psi_{1}\|_{L_{t,x}^{2}}\|P_{s_{1}}Q_l\nabla^{-1}I^cQ_{\nu j}(P_{k_{2}}\psi_{2}, P_{k_{3}}\psi_{3}))\|_{L_{t}^2L_x^\infty}\\&\hspace{8cm}\times \|P_{r_{1}}\nabla^{-1}IQ_{\mu j}(P_{k_{4}}\psi_{4}, P_{k_{5}}\psi_{5})\|_{L_{t,x}^{\infty}}
\end{align*}
Then from Lemma 4.18 and Bernstein's inequality we infer that
provided $k_2\gg s_1$, we have
\[
2^{-\frac{s_1}{2}}\|P_{s_{1}}Q_l\nabla^{-1}I^cQ_{\nu j}(P_{k_{2}}\psi_{2}, P_{k_{3}}\psi_{3}))\|_{L_{t}^2L_x^\infty}
\lesssim 2^{\epsilon l} 2^{-\epsilon \max\{s_1, k_{2,3}\}}\max\{k_2-s_1, 1\}^{2}\prod_{j=2,3}\|P_{k_j}\psi_j\|_{S[k_j]}
\]
Inserting this into the preceding bound we infer that
\begin{align*}
&\|P_{k_{1}}Q_{>l-10}\psi_{1}\|_{L_{t,x}^{2}}\|P_{s_{1}}Q_l\nabla^{-1}I^cQ_{\nu j}(P_{k_{2}}\psi_{2}, P_{k_{3}}\psi_{3}))\|_{L_{t}^2L_x^\infty}\\&\hspace{6cm}\times \|P_{r_{1}}\nabla^{-1}IQ_{\mu j}(P_{k_{4}}\psi_{4}, P_{k_{5}}\psi_{5})\|_{L_{t,x}^{\infty}}\\
&\lesssim 2^{\frac{s_1-l}{2}}2^{\epsilon l} 2^{-\epsilon \max\{s_1, k_{2,3}\}}\max\{k_2-s_1, 1\}^{2}2^{\min\{r_1, k_{4,5}\}-\max\{r_1, k_{4,5}\}}\prod_{j=1}^{5}\|P_{k_{j}}\psi_j\|_{S[k_j]}
\end{align*}
Summing over $l>s_1$ yields the bound of the first inequality of the lemma with $\delta=\epsilon-$. On the other hand, when $k_2=s_1+O(1)$, say, one can use Lemma 4.23 instead, which then gives the desired inequality with $\delta=\frac{1}{2}-\epsilon$.
\\
Next consider \eqref{eq:decomp111}. Here we distinguish between the cases $l<r_1+O(1)$ and $l\gg r_1$. In the former case,  as before assuming $s_1<-10$, we get
\begin{align*}
&\|\nabla_{x,t}P_{[-5,5]}Q_{<10}[\big(P_{0}[\nabla^{-1}(P_{k_{1}}Q_{<l-10}\psi_{1} P_{s_{1}}Q_l\nabla^{-1}I^cQ_{\nu j}(P_{k_{2}}\psi_{2}, P_{k_{3}}\psi_{3}))]\nonumber\\&\hspace{8cm}\times P_{r_{1}}\nabla^{-1}IQ_{\mu j}(P_{k_{4}}\psi_{4}, P_{k_{5}}\psi_{5})\big)\|_{L_t^1\dot{H}^1}\\
&\lesssim \|P_{k_{1}}Q_{<l-10}\psi_{1}\|_{L_t^\infty L_x^2}\|P_{s_{1}}Q_l\nabla^{-1}I^cQ_{\nu j}(P_{k_{2}}\psi_{2}, P_{k_{3}}\psi_{3})\|_{L_t^2L_x^\infty}\\&\hspace{8cm}\times\|P_{r_{1}}\nabla^{-1}IQ_{\mu j}(P_{k_{4}}\psi_{4}, P_{k_{5}}\psi_{5})\|_{L_t^2L_x^\infty}
\end{align*}
Using Lemma 4.18-4.23 again, we obtain the bound
\[
\lesssim 2^{\frac{s_1+r_1}{2}}2^{\epsilon(l-k_2)}|s_1-k_2|^{2}2^{\frac{\min\{r_1, k_{4,5}\}-\max\{r_1, k_{4,5}\}}{2}}\prod_{j=1}^5\|P_{k_j}\psi_j\|_{S[k_j]}
\]
One may sum here over $s_1<l<r_1+O(1)$ to get the desired first inequality of the lemma with $\delta=\epsilon-$. \\
Next, consider the case $l\gg r_1$. But in this case we can write
\begin{align*}
&\nabla_{x,t}P_{[-5,5]}Q_{<10}[\big(P_{0}[\nabla^{-1}(P_{k_{1}}Q_{<l-10}\psi_{1} P_{s_{1}}Q_l\nabla^{-1}I^cQ_{\nu j}(P_{k_{2}}\psi_{2}, P_{k_{3}}\psi_{3}))]\nonumber\\&\hspace{8cm}\times P_{r_{1}}\nabla^{-1}IQ_{\mu j}(P_{k_{4}}\psi_{4}, P_{k_{5}}\psi_{5})\big)\\
&=\nabla_{x,t}P_{[-5,5]}Q_{[l-10, 10]}[\big(P_{0}[\nabla^{-1}(P_{k_{1}}Q_{<l-10}\psi_{1} P_{s_{1}}Q_l\nabla^{-1}I^cQ_{\nu j}(P_{k_{2}}\psi_{2}, P_{k_{3}}\psi_{3}))]\nonumber\\&\hspace{8cm}\times P_{r_{1}}\nabla^{-1}IQ_{\mu j}(P_{k_{4}}\psi_{4}, P_{k_{5}}\psi_{5})\big)
\end{align*}
But this we can then estimate via the
$\|.\|_{\dot{X}_0^{-1,-\frac{1}{2},1}}$ -norm of the output, i.e.,
it suffices to bound
\begin{align*}
&\|\nabla_{x,t}P_{[-5,5]}Q_{[l-10, 10]}[\big(P_{0}[\nabla^{-1}(P_{k_{1}}Q_{<l-10}\psi_{1} P_{s_{1}}Q_l\nabla^{-1}I^cQ_{\nu j}(P_{k_{2}}\psi_{2}, P_{k_{3}}\psi_{3}))]\nonumber\\&\hspace{7cm}\times P_{r_{1}}\nabla^{-1}IQ_{\mu j}(P_{k_{4}}\psi_{4}, P_{k_{5}}\psi_{5})\big)\|_{\dot{X}_0^{-1,-\frac{1}{2},1}}\\
&\lesssim 2^{-\frac{l}{2}}\|P_{s_{1}}Q_l\nabla^{-1}I^cQ_{\nu j}(P_{k_{2}}\psi_{2}, P_{k_{3}}\psi_{3})\|_{L_t^2L_x^\infty}
\|P_{k_{1}}Q_{<l-10}\psi_{1}\|_{L_t^\infty L_x^2}\|\\&\hspace{7cm}\times P_{r_{1}}\nabla^{-1}IQ_{\mu j}(P_{k_{4}}\psi_{4}, P_{k_{5}}\psi_{5})\|_{L_{t,x}^\infty}
\end{align*}
From here the estimates are continued in a fashion identical to the ones used to control \eqref{eq:decomp101}. This completes estimating the contribution of $P_{s_{1}}\nabla^{-1}I^cQ_{\nu j}(P_{k_{2}}\psi_{2}, P_{k_{3}}\psi_{3}))$.
\\

{\it{(iib): contribution of the hyperbolic type term.}} Next we consider the contribution of
\[
P_{s_{1}}\nabla^{-1}IQ_{\nu j}(P_{k_{2}}\psi_{2}, P_{k_{3}}\psi_{3})),
\]
which is the expression
\begin{align*}
&\nabla_{x,t}P_{[-5,5]}Q_{<10}[\big(P_{0}[\nabla^{-1}(P_{k_{1}}\psi_{1} P_{s_{1}}\nabla^{-1}IQ_{\nu j}(P_{k_{2}}\psi_{2}, P_{k_{3}}\psi_{3}))]\nonumber\\&\hspace{8cm}\times P_{r_{1}}\nabla^{-1}IQ_{\mu j}(P_{k_{4}}\psi_{4}, P_{k_{5}}\psi_{5})\big)
\end{align*}
We shall again make the reduction $s_1<-10$, the remaining frequency interactions being treated analogously.
This is accomplished using Lemma 4.16. We obtain
\begin{align*}
&\|\nabla_{x,t}P_{[-5,5]}Q_{<10}[\big(P_{0}[\nabla^{-1}(P_{k_{1}}\psi_{1} P_{s_{1}}\nabla^{-1}IQ_{\nu j}(P_{k_{2}}\psi_{2}, P_{k_{3}}\psi_{3}))]\nonumber\\&\hspace{8cm}\times P_{r_{1}}\nabla^{-1}IQ_{\mu j}(P_{k_{4}}\psi_{4}, P_{k_{5}}\psi_{5})\big)\|_{L_t^1\dot{H}^{-1}}\\
&\lesssim \|P_{k_{1}}\psi_{1}\|_{L_t^\infty L_x^2}\|P_{s_{1}}\nabla^{-1}IQ_{\nu j}(P_{k_{2}}\psi_{2}, P_{k_{3}}\psi_{3})\|_{L_t^2L_x^\infty}\|P_{r_{1}}\nabla^{-1}IQ_{\mu j}(P_{k_{4}}\psi_{4}, P_{k_{5}}\psi_{5})\|_{L_t^2L_x^\infty}\\
&\lesssim 2^{\frac{s_1+r_1}{2}}2^{\frac{\min\{s_1, k_{2,3}\}-\max\{s_1, k_{2,3}\}}{2}}2^{\frac{\min\{r_1, k_{4,5}\}-\max\{r_1, k_{4,5}\}}{2}}\prod_{j=1}^{5}\|P_{k_j}\psi_j\|_{S[k_j]}
\end{align*}
This is as desired with $\delta=\frac{1}{2}$.
\end{proof}

\subsubsection{Error terms of order higher than five.} Here we consider the errors generated by repeated
application of Hodge decompositions, which are of higher than quintic degree. We recall that they arise
when we apply repeated Hodge decompositions to the second and third input in
\[
\nabla_{x,t}[\psi\nabla^{-1}(\psi^{2})]
\]
or else to the second and third input in
\[
\nabla_{x,t}[\nabla^{-1}[\psi\nabla^{-1}(\psi^{2})]\nabla^{-1}IQ_{\nu j}(\psi, \psi)]
\]
To simplify the discussion, we shall call terms that arise in the
first situation 'of the first type', while those in that arise in
the second situation will be called of 'second type'. In either
case, we associate a binary graph with each such expression as in
the discussion above. We call expressions whose associated graph has
only directed subgraphs of length at most three 'short', and those
with directed graphs of length at least four 'long'.  For technical
reasons, it will be most convenient to organize the 'short' and
'long' higher order terms into suitable sums, which are easier to
estimate. Specifically, note that each of these higher order terms
consists of nested terms of the form
\begin{align}\label{eq:arrghproblem}
\ldots\nabla^{-1}P_{s_1}[P_{k_1}R_{\nu}\psi P_{r_1}\nabla^{-1}(P_{k_{2}}\psi P_{s_2}\nabla^{-1}P_{s_{2}}[\ldots]],
\end{align}
here the case of a node with one outgoing edge, or alternatively
\begin{align}\label{eq:noproblem}
\ldots\nabla^{-1}P_{s_1}[P_{s_2}\nabla^{-1}[P_{k_1}\psi\nabla^{-1}P_{r_1}[\ldots]]\nabla^{-1}(P_{k_{2}}\psi P_{s_3}\nabla^{-1}[\ldots]]
\end{align}
in case of two outgoing edges.
\\
It is the first type of expression which may cause some mild
difficulties due to the presence of the $R_{\nu}$-operator, which
for $\nu=0$ may be formally unbounded. However, re-combining a term
of type \eqref{eq:arrghproblem} with a suitable term of the form
\eqref{eq:noproblem} and using the relation
\[
R_{\nu}\psi+\chi_\nu=\psi_\nu,\quad \psi=-\sum_{k=1,2}R_{k}\psi_k,
\]
we replace each such 'intermediate' gradient term (i.e., not
contributing to one of the innermost $Q_{\nu j}$ null-forms in case
of 'short' expressions) $R_{\nu}\psi$ by its non-gradient
counterpart $\psi_{\nu}$. We shall call the resulting expressions
'reduced'. Thus for example the (short) quintilinear expression
\[
\nabla_{x,t}[P_{k_{1}}\psi_1\nabla^{-1}P_{r_{1}}[P_{k_{2}}R_\nu\psi_2 P_{r_2}[P_{k_{3}}\psi\nabla^{-1}P_{r_{3}}Q_{\mu j}(P_{k_4}\psi_4, P_{k_5}\psi_5)]]]
\]
has reduced version
\[
\nabla_{x,t}[P_{k_{1}}\psi_1\nabla^{-1}P_{r_{1}}[P_{k_{2}}\psi_{2\nu} P_{r_2}[P_{k_{3}}\psi\nabla^{-1}P_{r_{3}}Q_{\mu j}(P_{k_4}\psi_4, P_{k_5}\psi_5)]]]
\]
Now we can formulate

\begin{prop}\label{prop:HigherOrderTypeI} Let
\[
P_{0} F_{2l+1}(\psi),\quad l=2,3,4
\]
be short reduced higher order expression of first type at frequency $\sim 1$. We can write it in nested form
\begin{align*}
&P_{0}F_{2l+1}(\psi)=\\&\nabla_{x,t}P_0\big[P_{k_1}\psi_1\nabla^{-1}P_{r_1}\big[\ldots \nabla^{-1}P_{r_{j}}\big[P_{k_{j+1}}\psi_{j+1}\nabla^{-1}P_{r_{j+1}}\big[P_{k_{j+2}}\psi_{j+2}\\&\hspace{6cm}\times\nabla^{-1}\big[\ldots \nabla^{-1}P_{r_{2l-1}}Q_{\mu j}(P_{k_{2l}}\psi_{2l}, P_{k_{2l+1}}\psi_{2l+1})\big]\big]\big]\big]\big]
\end{align*}
Then we have the following bounds:
\\

(1) If $r_1\ll -10$, we have
\[
\|P_{0}F_{2l+1}(\psi)\|_{N[0]}\lesssim 2^{\delta[\min\{k_{2},\ldots, k_{2l+1}, r_1,\ldots, r_{2l-1}\}-\max\{k_{2},\ldots, k_{2l+1}, r_1,\ldots, r_{2l-1}\}}\prod_{j=1}^{2l+1}\|P_{k_{j}}\psi_{j}\|_{S[k_{j}]}
\]
for a suitable constant $\delta>0$.
\\

(2) If $r_{1}\in [-10, 10]$, we get the bound
\[
\|P_{0}F_{2l+1}(\psi)\|_{N[0]}\lesssim  2^{\delta k_1}2^{\delta[\min\{k_{2},\ldots, k_{2l+1}, r_1,\ldots, r_{2l-1}\}-\max\{k_{2},\ldots, k_{2l+1}, r_1,\ldots, r_{2l-1}\}}\prod_{j=1}^{2l+1}\|P_{k_{j}}\psi_{j}\|_{S[k_{j}]}
\]

(3) If $r_{1}>10$, we have
\[
\|P_{0}F_{2l+1}(\psi)\|_{N[0]}\lesssim  2^{-\delta k_1}2^{\delta[\min\{k_{2},\ldots, k_{2l+1}, r_1,\ldots, r_{2l-1}\}-\max\{k_{2},\ldots, k_{2l+1}, r_1,\ldots, r_{2l-1}\}}\prod_{j=1}^{2l+1}\|P_{k_{j}}\psi_{j}\|_{S[k_{j}]}
\]
\end{prop}

The proof of this follows the exact same pattern as the one for  Proposition~\ref{quintilinear1}, and is omitted. In fact, for $l>2$, one no longer needs to use the sharp improved Strichartz endpoint as in the case $l=2$.
\\

In a similar vein, we have the analogue of Proposition~\ref{quintilinear2}. A short reduced expression of the second type can be written as
\begin{align}\label{mess2}
P_{0}F_{2l+3}(\psi)=&\nabla_{x,t}P_0\big[\nabla^{-1}P_{r_1}\big[\ldots \nabla^{-1}P_{r_{j}}\big[P_{k_{j+1}}\psi_{j+1}\nabla^{-1}P_{r_{j+1}}\big[P_{k_{j+2}}\psi_{j+2}\\&\hspace{1cm}\times\nabla^{-1}\big[\ldots \nabla^{-1}P_{r_{2l-1}}Q_{\mu j}(P_{k_{2l}}\psi_{2l}, P_{k_{2l+1}}\psi_{2l+1})\big]\big]\big]\big]P_{s_1}Q_{\nu j}\nabla^{-1}(P_{k_{2l+2}}\psi_{2l+2}, P_{k_{2l+3}}\psi_{2l+3} \big]\nonumber,
\end{align}
where $l=1,2,3$. Then we have

\begin{prop}\label{prop:HigherOrderTypeII} Using the representation \eqref{mess2},
let $P_{0}F_{2l+3}(\psi)$ be a short reduced term at frequency $\sim 1$ of the the second type. Then the following hold:
\\

(1) If $s_{1}<-10$, we have
\begin{align*}
\|P_{0}F_{2l+3}(\psi)\|_{N[0]}\lesssim &2^{\delta s_1}2^{\delta[\min\{s_1, k_{2l+2}, k_{2l+3}\}-\max\{s_1, k_{2l+2}, k_{2l+3}\}]}
\\&\times 2^{\delta[\min\{k_{2},\ldots, k_{2l}, r_1,\ldots, r_{2l-1}\}-\max\{k_{2},\ldots, k_{2l}, r_1,\ldots, r_{2l-1}\}}\prod_{j=1}^{2l+3}\|P_{k_{j}}\psi_{j}\|_{S[k_{j}]}
\end{align*}

(2) If $s_{1}\in[-10, 10]$, we have
\begin{align*}
\|P_{0}F_{2l+3}(\psi)\|_{N[0]}\lesssim &2^{\delta r_1}2^{\delta[\min\{s_1, k_{2l+2}, k_{2l+3}\}-\max\{s_1, k_{2l+2}, k_{2l+3}\}]}
\\&\times 2^{\delta[\min\{k_{2},\ldots, k_{2l}, r_1,\ldots, r_{2l-1}\}-\max\{k_{2},\ldots, k_{2l}, r_1,\ldots, r_{2l-1}\}}\prod_{j=1}^{2l+3}\|P_{k_{j}}\psi_{j}\|_{S[k_{j}]}
\end{align*}

(3) If $s_{1}>10$, we have
\begin{align*}
\|P_{0}F_{2l+3}(\psi)\|_{N[0]}\lesssim &2^{-\delta r_1}2^{\delta[\min\{s_1, k_{2l+2}, k_{2l+3}\}-\max\{s_1, k_{2l+2}, k_{2l+3}\}]}
\\&\times 2^{\delta[\min\{k_{2},\ldots, k_{2l}, r_1,\ldots, r_{2l-1}\}-\max\{k_{2},\ldots, k_{2l}, r_1,\ldots, r_{2l-1}\}}\prod_{j=1}^{2l+3}\|P_{k_{j}}\psi_{j}\|_{S[k_{j}]}
\end{align*}
\end{prop}

Again the proof is similar to the one of Proposition~\ref{quintilinear2}.
\\
Note that in order to estimate the expressions of short type, we still need to to use a little bit of null-structure to make
them amenable to estimation by the $S$-spaces. This is no longer the case for 'long' expressions $P_{0}F_11(\psi)$ of reduced type:  write such an expression as
\begin{align*}
&P_{0}F_{11}(\psi)=\\&\nabla_{x,t}P_0\big[P_{k_1}\psi_1\nabla^{-1}P_{r_1}\big[\ldots \nabla^{-1}P_{r_{j}}\big[P_{k_{j+1}}\psi_{j+1}\nabla^{-1}P_{r_{j+1}}\big[P_{k_{j+2}}\psi_{j+2}\\&\hspace{6cm}\times\nabla^{-1}\big[\ldots \nabla^{-1}P_{r_{9}}(P_{k_{10}}\psi_{10} P_{k_{11}}\psi_{11})\big]\big]\big]\big]\big]
\end{align*}
if it is of first type or
\begin{align*}
P_{0}F_{11}(\psi)=&\nabla_{x,t}P_0\big[\nabla^{-1}P_{r_1}\big[\ldots \nabla^{-1}P_{r_{j}}\big[P_{k_{j+1}}\psi_{j+1}\nabla^{-1}P_{r_{j+1}}\big[P_{k_{j+2}}\psi_{j+2}\\&\hspace{1cm}\times\nabla^{-1}\big[\ldots \nabla^{-1}P_{r_{9}}(P_{k_{10}}\psi_{10}P_{k_{11}}\psi_{11})\big]\big]\big]\big]P_{s_1}\nabla^{-1}(P_{k_{12}}\psi_{12}P_{k_{13}}\psi_{13}) \big]\nonumber,
\end{align*}
Note that the innermost bilinear expressions
\[
\nabla^{-1}P_{r_{9}}(P_{k_{10}}\psi_{10} P_{k_{11}}\psi_{11})
\]
are no longer null-forms.
\begin{prop}\label{prop:HigherOrderLong}
Let $P_{0}F_{11}(\psi)$ be a long expression of either first or second type, written as in immediately preceding. The if $P_{0}F_{11}(\psi)$ is of the first type, we have if $r_{1}<-10$
\begin{align*}
\|P_{0}F_{11}(\psi)\|_{N[0]}\lesssim  2^{\delta[\min\{k_{2},\ldots, k_{9}, r_1,\ldots, r_{9}\}-\max\{k_{2},\ldots, k_{9}, r_1,\ldots, r_{9}\}}2^{\delta[\min\{k_{10}, k_{11}\}-\max\{k_{10}, k_{11}\}]}\prod_{j=1}^{11}\|P_{k_{j}}\psi_{j}\|_{S[k_{j}]}
\end{align*}
Thus by contrast to Proposition~\ref{prop:HigherOrderTypeI} case (1), we have an extra factor
\[
2^{\delta[\min\{k_{10}, k_{11}\}-\max\{k_{10}, k_{11}\}]}
\]
whence we cannot gain in case of high-high interactions in the innermost expression
\[
\nabla^{-1}P_{r_{9}}(P_{k_{10}}\psi_{10} P_{k_{11}}\psi_{11})
\]
\end{prop}

\begin{proof} This is purely an application of our available Strichartz norms: indeed, we have for suitable $\delta>0$
\begin{align*}
&\|P_{r_{8}}[P_{k_9}\psi_9 \nabla^{-1}P_{r_9}(P_{k_{10}}\psi_{10}P_{k_{11}}\psi_{11})]\|_{L_t^8L_x^{\frac{4}{3}+}}\\&\lesssim
2^{(\frac{3}{8}+\nu) r_8}2^{\delta[\min\{r_{8,9}, k_9\}-\max\{r_{8,9}, k_9\}]}2^{\delta[\min\{k_{10}, k_{11}\}-\max\{k_{10}, k_{11}\}]}\prod_{j=9}^{11}\|P_{k_{j}}\psi_{j}\|_{S[k_{j}]}
\end{align*}
where we define
\[
\nu=\frac{3}{4}-(\frac{4}{3}+)^{-1}
\]
Further, we have for $p=1,2,\ldots, 7$ and suitable $\delta_p>0$
\[
2^{-\nu_{p+1}r_1}\|P_{r_1}[P_{k}\psi\nabla^{-1}P_{r_2}F]\|_{L_{t}^{\frac{8}{p+1}}L_x^{\frac{4}{3}+}}\lesssim 2^{\delta_p[\min\{r_{1,2},k\}-\max\{r_{1,2},k\}]}[2^{-\nu_p r_2}\|P_{r_2}F\|_{L_{t}^{\frac{8}{p+1}}L_x^{\frac{4}{3}+}}]
 \]
where scaling dictates
\[
\nu_p=2-\frac{p}{8}-2(\frac{4}{3}+)^{-1}
\]
The proposition follows by applying these two inequalities sufficiently often.
\end{proof}

\begin{remark}\label{rem:quinticimprov}
We note that in the estimates above, we have not used wave-packet atoms.
\end{remark}

\section{Some basic perturbative results}
\label{sec:perturb}

This section develops some of the basic perturbative theory required for our work. More
precisely, we introduce a norm locally on some time interval
$(-T_0,T_1)$ which we denoted by $\|\psi\|_{S(-T_0,T_1)}$ with the
property that its finiteness insures that the gauged wave map $\psi$
can be continued outside of that time interval. The second topic we
discuss is the issue of defining wave maps with data which are
merely of energy class. This is accomplished by means of passing to
the limit in energy of smooth wave maps.

\subsection{A blow-up criterion}

Assume we are given a wave map $\bfu: (-T_{0}, T_{1})\times\R^{2} \to \Hyp^{2}$ with Schwartz
data at time $t=0$, by which we mean that the derivative components $\phi^{i}_{\alpha}$, $i=1,2$, $\alpha=0,1,2$,
and thus also the Coulomb components $\psi^{i}_{\alpha}$, are Schwartz functions at time $t=0$. These functions will
 then also be Schwartz on fixed time slices on the maximal interval of existence $(-T_{0}, T_{1})\times\R^{2}$. The
 following norm will provide us with sufficient control for long time existence and scattering.

\begin{defi}
 \label{def:Snorm} For any Schwartz function on $  (-T_{0}, T_{1})\times\R^{2}$ set
\[
\|\psi\|_{S}:=  \Big(\sum_{k\in\Z}\| P_{k}\psi\|_{S[k]((-T_0,
T_1)\times\R^{2})}^{2}\Big)^{\frac{1}{2}}
\]
Here
\[
   \|P_k\psi\|_{S[k]((-T_0,T_1)\times \R^2)} := \sup_{T<T_1,T'<T_0} \|P_k\psi_i\|_{S[k]([-T',T]\times \R^2)}
\]
where the local norms are those from~\eqref{eq:Skloc} {\em using
the} $\trip\cdot\trip$-norm.
\end{defi}

The goal of this section is to prove the following result.

\begin{prop}\label{BlowupCriterion} Let $ (-T_{0}, T_{1})$ be the maximal interval of existence for the wave map~$\bfu$   in the smooth sense. If
$\|\psi\|_{S}<\infty$
then necessarily $T_{0}=T_{1}=\infty$. Moreover, the wave map scatters at infinity, i.e.,  the components $\psi, \phi$
approach free waves in the energy topology as $t\to\pm\infty$.
\end{prop}

The strategy for proving the theorem will be to demonstrate an apriori bound
\[
\sup_{t\in (-T_{0}, T_{1})}\|\phi(t,\cdot)\|_{H^{s}}< \infty
\]
for some $s>0$, using the assumption $\|\psi\|_{S}<\infty$. By the
Klainerman-Machedon local well-posedness theory, this implies that
$u$ may be extended smoothly to some interval $(-T_{0}, T_{1}+\eps)$
for $\eps>0$ provided $T_{1}<\infty$, which contradicts minimality,
and similarly for $T_{0}$. Once we know that $u$ exists for $t\in
(-\infty, \infty)$,  scattering will follow by using a similar
argument. To obtain apriori control over sub-critical norms, we use
Tao's device of {\it{frequency envelope}}: for some $\delta_{1}>0$
depending only on certain apriori parameters specified later,
define
\begin{equation}\label{eq:ck_def}
c_{k}:=(\sum_{\ell\in
\Z}2^{-\delta_{1}|k-\ell|}\|P_{\ell}\phi(0,\cdot)\|_{L_{x}^{2}}^{2})^{\frac{1}{2}}
\end{equation}
Here as always we let $\phi=\{\phi^{i}_{\alpha}\}$ the vector of derivative components.
Proposition~\ref{BlowupCriterion} now follows from the following result.

\begin{prop}\label{BlowupCriterion1}  Let $ (-T_{0}, T_{1})$ be the maximal interval of existence for the wave map~$\bfu$   in the smooth sense. If
$\|\psi\|_{S}<\infty$, then there exists a number $C_{1}=C_{1}(u)<\infty$ (which may depend in a
complicated fashion on the wave map~$\bfu$, and not just its energy)  such that
\[
\|P_{k}\phi \|_{L^\infty_t((-T_{0}, T_{1});L_{x}^{2})}\leq C_{1}c_{k}
\]
In fact,
\[
\|P_{k}\psi\|_{S[k]((-T_0, T_1)\times\R^{2})}\leq C_{1}c_{k}
\]
\end{prop}

\noindent To establish the existence of $C_{1}$, we shall cover the
time interval $(-T_{0}, T_{1})$ by a finite number  of shorter open
intervals $I_{j}$ (which can still be very large): let
$\|\psi\|_{S}<C_0$. Then
\begin{equation}\label{eq:Ij}
(-T_{0}, T_{1})=\cup_{j=1}^{M_{1}}I_{j},\quad M_{1}=M_{1}(C_{0}),
\end{equation}
where $\psi|_{I_{j}}$ will satisfy a suitable smallness property.
The idea then is to bootstrap certain bounds on each $I_{j}$,
beginning with the interval containing the time slice $t=0$. More
precisely, the intervals $I_{j}$ will be chosen so that the wave map
restricted to each~$I_{j}$ is {\it{well approximated by a free
wave}}. While the error can be treated perturbatively, the free wave
has better dispersive properties which we can exploit. All functions
will be smooth in space and time and Schwartz functions on fixed
time slices.

\subsubsection{Splitting the wave map on shorter time intervals}

We first derive a simple estimate on the nonlinearities appearing in~\eqref{eq:psisys1} and~\eqref{eq:psisys2}.
It will be based entirely on the Strichartz estimates, see Lemma~\ref{lem:Strich}.
We will keep the time interval $(-T_0,T_1)$ from above fixed throughout.

\begin{lemma}
 \label{lem:easytrilin}
Let $\max_{i=1,2,3}\|\psi_i\|_{S}<C_0$. Then
\[
 \|P_0(\psi_1 |\nabla|^{-1} (\psi_2 \psi_3))\|_{L^M_t((-T_0,T_1); L^2_x)} \les C_0^2 \sup_{i=1,2,3}\; \sup_{k\in\Z}\;
2^{-\frac{|k|}{M}} \|P_k\psi_i\|_{S[k]((-T_0,T_1)\times \R^2)}
\]
provided $M$ is large and with an absolute implicit constant. Alternatively, one has the bound
\[
 \|P_0(\psi_1 |\nabla|^{-1} (\psi_2 \psi_3))\|_{L^M_t((-T_0,T_1); L^2_x)} \les \|\psi_2\|_{\ener} \|\psi_3\|_{\ener} \;\sup_{k\in\Z}\;
2^{-\frac{|k|}{M}} \|P_k\psi_1\|_{S[k]((-T_0,T_1)\times \R^2)}
\]
with an absolute implicit constant.
\end{lemma}
\begin{proof} Assume to begin with that $\psi_i$ is adapted to $k_i\in\Z$. As in Section~\ref{sec:trilin},
we now consider all possible cases of interactions.  Also, we shall drop the time interval $(-T_0,T_1)$ from
our notation with the understanding that integration in time is to be restricted to this interval. Moreover,
replacing each $\psi_i$ by a globally defined Schwartz function $\tilde \psi_i$ with the property that
\[
 \|\tilde\psi_i\|_{S[k_i]} \le 2 \|\psi_i\|_{S[k_i]([-T',T]\times \R^2)},\quad \tilde\psi_i|_{[-T',T]}=\psi_i|_{[-T',T]}
\]
for some $T',T$ as above, allows us to assume that the $\psi_i$ are globally defined initially.
Finally, fix any $M\ge100$.

\medskip \noindent
 {\em Case 1:} $\mathit{0\le k_1\le k_2+O(1)=k_3+O(1) }$.
Then
\begin{align*}
& \|P_0(\psi_1 |\nabla|^{-1} (\psi_2 \psi_3))\|_{L^M_t L^2_x} \les \|P_0(\psi_1 |\nabla|^{-1} (\psi_2 \psi_3))\|_{L^M_t L^1_x} \\
&\les \|\psi_1\|_{L^M_t L^\infty_x} \, 2^{-k_1}  \|\psi_2\|_{L^\infty_t L^2_x} \|\psi_3\|_{L^\infty_t L^2_x}
\les 2^{-\frac{k_1}{M}} \prod_{i=1}^3 \|\psi_i\|_{S[k_i]}
\end{align*}
where the final estimate is from~\eqref{eq:Strich}.

\medskip
\noindent {\em Case 2:} $\mathit{0\le k_1= k_2+O(1), k_3\le k_2-C }$.
If $k_3\ge0$, one proceeds as in Case~1. Otherwise,
\begin{align*}
& \|P_0(\psi_1 |\nabla|^{-1} (\psi_2 \psi_3))\|_{L^M_t L^2_x} \les \|P_0(\psi_1 |\nabla|^{-1} (\psi_2 \psi_3))\|_{L^M_t L^1_x} \\
&\les \|\psi_1\|_{L^\infty_t L^2_x} \, 2^{-k_1}  \|\psi_2\|_{L^\infty_t L^2_x} \|\psi_3\|_{L^M_t L^\infty_x} \les 2^{k_3-k_1}
\|\psi_1\|_{S[k_1]}    \|\psi_2\|_{S[k_2]} \|\psi_3\|_{L^M_t L^2_x}  \\
&\les  2^{-k_1} 2^{k_3(1-\frac{1}{M})} \prod_{i=1}^3 \|\psi_i\|_{S[k_i]}
\end{align*}
by \eqref{eq:Strich} and Bernstein's inequality.

\medskip
\noindent {\em Case 3:} $\mathit{0\le k_1= k_3+O(1), k_2\le k_3-C .}$ This
case is symmetric to the previous one.

\medskip
\noindent {\em Case 4:} $\mathit{O(1)\le k_2= k_3+O(1), k_1\le -C }$.  Here
\begin{align*}
& \|P_0(\psi_1 |\nabla|^{-1} (\psi_2 \psi_3))\|_{L^M_t L^2_x} \les \|P_0(\psi_1 |\nabla|^{-1} \tilde P_0(\psi_2 \psi_3))\|_{L^M_t L^1_x} \\
&\les \|\psi_1\|_{L^M_t L^\infty_x} \,    \|\psi_2\|_{L^\infty_t L^2_x} \|\psi_3\|_{L^\infty_t L^2_x} \les  2^{k_1}
\|\psi_1\|_{L^M_t L^2_x}    \|\psi_2\|_{L^\infty_t L^2_x} \|\psi_3\|_{L^\infty_t L^2_x}  \\
&\les 2^{(1-\frac{1}{M})k_1} \prod_{i=1}^3 \|\psi_i\|_{S[k_i]}
\end{align*}

\medskip
\noindent {\em Case 5:} $\mathit{k_1=O(1),\; k_2= k_3+O(1)}$. In this case we estimate
\begin{align*}
& \|P_0(\psi_1 |\nabla|^{-1} (\psi_2 \psi_3))\|_{L^M_t L^2_x} \les \sum_{k\le k_2\wedge 0+C}
 \|P_0(\psi_1 |\nabla|^{-1} P_k (\psi_2 \psi_3))\|_{L^M_t L^2_x} \\
&\les \sum_{k\le k_2\wedge 0+C}   \|\psi_1\|_{L^M_t L^2_x} \| |\nabla|^{-1} P_k (\psi_2 \psi_3)\|_{\Linf} \\
&\les \sum_{k\le k_2\wedge 0+C}     \|\psi_1\|_{S[k_1] } \, 2^k \| \psi_2 \psi_3\|_{L^\infty_t L^1_x} \\
&\les  2^{k_2\wedge0} \prod_{i=1}^3 \|\psi_i\|_{S[k_i]}
\end{align*}

\medskip
\noindent {\em Case 6:} $\mathit{O(1)= k_1\ge k_2+O(1)\ge k_3+C}$.
Here one has
\begin{align*}
& \|P_0(\psi_1 |\nabla|^{-1} (\psi_2 \psi_3))\|_{L^M_t L^2_x} \les
 \|P_0(\psi_1 |\nabla|^{-1} \tilde P_{k_1} (\psi_2 \psi_3))\|_{L^M_t L^2_x} \\
&\les     \|\psi_1\|_{L^M_t L^2_x} \| |\nabla|^{-1} \tilde P_{k_1}  (\psi_2 \psi_3)\|_{\Linf} \\
&\les       \|\psi_1\|_{S[k_1] } \, \| \psi_2 \psi_3\|_{L^\infty_t L^2_x} \les  \|\psi_1\|_{S[k_1] }  \| \psi_2\|_{\ener} \| \psi_3\|_{\Linf} \\
&\les  2^{k_3}  \prod_{i=1}^3 \|\psi_i\|_{S[k_i]}
\end{align*}

\medskip
\noindent {\em Case 7:} $\mathit{ k_1=O(1)\ge k_3+O(1)\ge k_2+C}$.
This case is symmetric to the previous one.

\medskip
\noindent {\em Case 8:} $\mathit{k_2=O(1), \max(k_1, k_3)\le -C}$. Finally, in
this case the estimate reads
\begin{align*}
& \|P_0(\psi_1 |\nabla|^{-1} (\psi_2 \psi_3))\|_{L^M_t L^2_x} \les
 \|P_0(\psi_1 |\nabla|^{-1} \tilde P_{0} (\psi_2 \psi_3))\|_{L^M_t L^2_x} \\
&\les     \|\psi_1\|_{L^M_t L^\infty_x} \| \psi_2 \psi_3\|_{\ener}  \les  2^{k_1(1-\frac{1}{M})} 2^{k_3} \prod_{i=1}^3 \|\psi_i\|_{S[k_i]}
\end{align*}

\medskip
\noindent {\em Case 9:} $\mathit{k_3=O(1), \max(k_1, k_2)\le -C}$.
This is symmetric to the preceding case.

\medskip \noindent We now drop the assumption on the frequency support of the inputs. Summing over all these cases yields the bound
\[
 \|P_0(\psi_1 |\nabla|^{-1} (\psi_2 \psi_3))\|_{L^M_t L^2_x} \les
\sup_{i=1,2,3}\; \sup_{k\in\Z}\,\big[
2^{-\frac{|k|}{M}} \|P_k\psi_i\|_{S[k]} \big]\;\max_{j=1,2,3}  \sum_{k\in\Z} \|P_k\psi_j\|_{S[k]}^2
\]
which proves the first bound. The proof of the second estimate is implicit in the preceding
and the lemma is proved.
\end{proof}

\begin{remark}
 \label{rem:ener2} If $\psi_2$, $\psi_3$ are gauged wave maps with energy bounded by~$E$, then the second bound of the preceding lemma becomes
\begin{equation}
 \label{eq:easytrilin2}
\|P_0(\psi_1 |\nabla|^{-1} (\psi_2 \psi_3))\|_{L^M_t((-T_0,T_1); L^2_x)} \les E^2 \;\sup_{k\in\Z}\;
2^{-\frac{|k|}{M}} \|P_k\psi_1\|_{S[k]((-T_0,T_1)\times \R^2)}
\end{equation}
with an absolute implicit constant.
\end{remark}

Our main goal here is to prove the following decomposition of the gauged wave map.

\begin{lemma}\label{lem:LocalSplitting} Let $\|\psi\|_{S}<C_0$.  Given $\eps_0>0$, there exist
$M_{1}=M_{1}(C_{0}, \eps_0)$ many intervals $I_{j}$ as
in~\eqref{eq:Ij} with the following property: for each
$I_{j}=(t_{j}, t_{j+1})$, there is a decomposition
\[
\psi|_{I_{j}}=\psi_{L}^{(j)}+\psi_{NL}^{(j)},\quad
\Box\psi_{L}^{(j)}=0
\]
which satisfies
\begin{align}
\sum_{k\in\Z}\|P_{k}\psi_{NL}^{(j)}\|_{S[k](I_{j}\times\R^{2})}^{2} &<\eps_0 \label{eq:psiNLbd} \\
\|\nabla_{x,t}\psi_{L}^{(j)}\|_{L_{t}^{\infty}\dot{H}^{-1}} &\le
M_2(C_0,\eps_0) \label{eq:psiLbd}
\end{align}
where the constant $M_2=M_2(C_0,\eps_0)$ satisfies $M_2\les C_0^3
\eps_0^{-\frac{1}{M}}$ with some $M\ge100$. Moreover,
$P_k\psi_{NL}^{(j)}$ and $P_k \psi_L^{(j)}$ are Schwartz functions
for each~$k\in \Z$. We also have the bounds
\begin{equation}\label{eq:daniel}
 \|\nabla_{x,t}P_{k}\psi_{L}^{(j)}\|_{\dot{H}_{x}^{-1}}+\|P_{k}\psi_{NL}^{(j)}\|_{S[k](I_{j}\times\R^{2})}\les c_{k}
\end{equation} with implied constant depending on $C_0$, provided
$c_{k}$ is a sufficiently flat frequency envelope with
$\|P_{k}\psi\|_{S[k]}\le c_k$.
\end{lemma}
\begin{proof}
The $\psi_\alpha$ satisfy the  system~\eqref{eq:psisys1}--\eqref{eq:psi_wave}.
Consider the frequency component~$P_0\psi_\alpha$.

\noindent {\em Case 1:} The underlying time interval~$I=(-T_0,T_1)$ is
very small, say $|I|<\eps_1$ with an $\eps_1$ that is to be determined. As explained in Section~\ref{subsec:waveeq}
one uses the div-curl
system~\eqref{eq:psisys1}, \eqref{eq:psisys2} in this case. Schematically,  this system takes
the form
\[
\partial_{t}P_{0}\psi=\nabla_{x}P_{0}\psi+P_{0}[\psi\nabla^{-1}(\psi^{2})]
\]
where we suppress the subscripts and also ignore the null-structure in the nonlinearity. Therefore,
\begin{equation}\label{eq:easystart}
\|P_{0}\psi(t)-P_{0}\psi(0)\|_{L_{x}^{2}}\leq \Big\|\int_{0}^{t}\nabla_{x}P_{0}\psi(s,\cdot)\,ds \Big\|_{L_{x}^{2}}
+\Big\|\int_{0}^{t}P_{0}[\psi\nabla^{-1}(\psi^{2})](s,\cdot)\,ds \Big\|_{L_{x}^{2}}
\end{equation}
For all $j\in\Z$ define
\begin{equation}\label{eq:ajdef}
 a_j:=\sup_{k\in \Z}\; 2^{-\frac{|k-j|}{M}}\|P_{k}\psi\|_{S[k](I\times \R^2)} \les C_0
\end{equation}
Clearly,
\begin{equation} \label{eq:intel1}
 \Big\|\int_{0}^{t}\nabla_{x}P_{0}\psi(s,\cdot)\,ds\Big\|_{L_{x}^{2}}\leq \eps_1 \|P_{0}\psi\|_{S[0](I\times \R^2)} \le a_0\, \eps_1
\end{equation}
Lemma~\ref{lem:easytrilin} implies
 \begin{align*}
\Big\| \int_{0}^{t}P_{0}[\psi\nabla^{-1}(\psi^{2})](s,\cdot)\,ds\Big\|_{L^\infty_t(I;L_{x}^{2})} &\lesssim C_0^2\,a_0\,  \eps_1^{1-\frac{1}{M}}\\
\Big\|
\int_{0}^{t}P_{0}[\psi\nabla^{-1}(\psi^{2})](s,\cdot)\,ds\Big\|_{L^2_t(I;L_{x}^{2})}
&\lesssim C_0^2\, a_0\, \eps_1^{\frac32-\frac{1}{M}}
\end{align*}
From the div-curl system \eqref{eq:psisys1} and \eqref{eq:psisys2},
\begin{align*}
 \|\partial_{t}P_{0}\psi \|_{L^2_t(I;L_{x}^{2})} &\le \|\nabla_{x}P_{0}\psi\|_{L^2_t(I;L_{x}^{2})} + \|P_{0}[\psi\nabla^{-1}(\psi^{2})]\|_{L^2_t(I;L_{x}^{2})}
\les C_0^2\, a_0\, \eps_1^{\frac12-\frac{1}{M}}
\end{align*}
where we assumed without loss of generality that $C_0\ge1$.  We claim that these bounds imply that
\begin{equation}\label{eq:claim63}
\Big\|
\int_{0}^{t}P_{0}[\psi\nabla^{-1}(\psi^{2})]\Big\|_{S[0](I\times
\R^2)}\ll \eps_0 a_0
\end{equation}
provided $\eps_1$ was chosen sufficiently small depending
on~$\eps_0$. To see this, let $I':=[-T',T]\subset I=(-T_0,T_1)$ and
pick any smooth bump function $\chi$ supported in~$I$ so that
$\chi=1$ on~$I'$ and with~$0\le \chi\le1$. Moreover, let
$\tilde\chi$ be any smooth compactly supported function
with~$\tilde\chi=1$ on~$I$ (the choice of this function does not
depend on~$I'$).  Then  define
\[
 \tilde\psi(t):= \tilde\chi(t)\big[ P_0\psi(0)+ \int_0^t \chi(s) \del_s P_0\psi(s)\, ds\big]
\]
By construction, $\tilde\psi$ is a global Schwartz function so that $\tilde\psi=\psi$ on~$I'$. Moreover,
by the preceding bounds,
\[
 \|\tilde\psi\|_{\Ltwotx}+ \|\del_t\tilde\psi\|_{\Ltwotx} \ll \eps_0\,a_0\,
\]
provided $\eps_1$ was chosen small enough (this smallness does not depend on the choice of~$I'$).
This now implies that
\[
 \|\tilde \psi\|_{\dot X_0^{0,\frac12,1}} \ll \eps_0\,a_0
\]
whence~\eqref{eq:claim63}. In view of \eqref{eq:intel1}, \eqref{eq:claim63} and~\eqref{eq:easystart},
\begin{equation}\label{eq:diff15}
 \|P_{0}\psi(t)-P_{0}\psi(0)\|_{S[0](I\times \R^2)}\ll  \eps_0 \, a_0
\end{equation}
We now define $P_{0}\psi_{L}$ to be the free wave with initial data $(P_{0}\psi(0), 0)$ at time $t=0$.
Clearly
\begin{align*}
 \|P_{0}\nabla_{x,t}\psi_{L}\|_{L_{t}^{\infty}\dot{H}^{-1}} &\lesssim \|P_{0}\psi(0)\|_{L_{x}^{2}} \\
 \|P_{0}\psi_L-P_{0}\psi_L(0)\|_{S[0](I\times \R^2)}&\ll \eps_0\,  a_0
\end{align*}
The second inequality here implies that
\[
\|P_{0}\psi-P_{0}\psi_{L}\|_{S[0](I_{j}\times\R^{2})}\ll \eps_0 \, a_0
\]
Thus in the present situation, we approximate $P_{0}\psi$ by the free wave $P_{0}\psi_{L}$ just described and the
bounds which we just obtained should be viewed as versions of~\eqref{eq:psiNLbd} and~\eqref{eq:psiLbd} on a fixed dyadic frequency
block. Several remarks are in order: First, we shall of course need to construct $\psi_L$ and~$\psi_{NL}$ for each such dyadic block $P_k$, and then obtain the
global bounds required by~\eqref{eq:psiNLbd} and~\eqref{eq:psiLbd}. In this regard, any bound depending on~$a_j$ can easily
be square-summed since
\[
 \sum_j a_j^2 \le C(M) \sum_{k\in\Z} \|P_{k}\psi\|_{S[k](I\times \R^2)}^2 \le C(M)\, C_0^2
\]
Second, the construction we just carried out applies to $P_k\psi$
equally well provided $|I|\le 2^{-k}\eps_1$. Moreover, $I$ can be
any time interval on which $\psi$ is defined --- with any $t_0\in I$
playing the role of $t=0$ --- and we shall indeed apply this exact
same procedure to those intervals~$I_j$ which we are about to
construct provided they satisfy this length restriction.

\noindent {\em Case 2:} The underlying time interval~$I=(-T_0,T_1)$
satisfies $|I|>\eps_1$ with $\eps_1$ as in Case~1.  To construct the
$I_j$, we shall use the wave equation~\eqref{eq:psi_wave}
for~$\psi_{\alpha}$. By means of Schwartz extensions and  successive
Hodge type decompositions of the $\psi_{\alpha}$-components as
explained above, the nonlinearity can be written as
\begin{equation}\label{eq:Falpha}
\Box\psi_{\alpha}=F_{\alpha}(\psi)=F_{\alpha}^{3}(\psi)+F_{\alpha}^{5}(\psi)+F_{\alpha}^{7}(\psi)+F_{\alpha}^{9}(\psi)+F_{\alpha}^{11}(\psi),
\end{equation}
where the superscripts denote the degree of multi-linearity, see Section~\ref{sec:hodge}.  The
contribution of the trilinear null-form $F_{\alpha}^{3}(\psi)$ here
is in a sense the principal contribution, and causes the main
technical difficulties. We now make the following claim: {\em There
exists a cover $I=\bigcup_{j=1}^{M_1} I_j$ by open intervals
$I_{j}$, $1\le j\le M_{1}$, $M_{1}=M_{1}(\eps_0)$, such that}
\begin{equation}\label{eq:Ijclaim}
\max_{1\le j\le
M_{1}}\,\sum_{\ell\in\Z}\|P_{\ell}F_{\alpha}(\psi)\|_{N[\ell](I_{j}\times\R^{2})}^{2}<\eps_0
C_0^6
\end{equation}
We verify this for each of the different types of nonlinearities
appearing on the right-hand side of~\eqref{eq:Falpha} starting with
the trilinear ones. Let us schematically write anyone of these
trilinear expressions in the form $\nabla_{t,x} [ \psi_1
|\nabla|^{-1} I^c\calN(\psi_2,\psi_3)]$  or $\nabla_{t,x} [ R\,
\psi_1 |\nabla|^{-1} I\calN(\psi_2,\psi_3)]$, where $\calN$ stands
for the usual bilinear nullforms and~$R$ for a Riesz transform (each
of the $\psi_i=\psi$ but it will be convenient to view these inputs
as independent). Break up to inputs into dyadic frequency pieces:
$\psi_i=\sum_{k_i} P_{k_i} \psi_i$ for $i=1,2,3$. In view of our
discussion in Section~\ref{subsec:improvetrilin}, it suffices to
consider the high-low-low case $|k_2-k_3|<L$, $k_2<k_1+L$ for some
large $L=L(\eps_0)$. In addition, it suffices to restrict attention
to frequencies $k>k_2-L'$ where $P_k$ localized the frequency
of~$\calN$ and $L'=L'(\delta)$ is large.   Finally, one can assume
angular separation between the inputs: there exists
$m_0=m_0(\eps_0)\ll-1$ so that~\eqref{eq:Ijclaim} reduces to the
estimates
\begin{align}
\max_{1\le j\le M_{1}}\,\sum_{\substack{k,\ell,k_1,k_2,k_3\\
\kappa_1,\kappa_2,\kappa_3}} \| \nabla_{t,x} P_{\ell}
[P_{k_1,\kappa_1} \psi_1 |\nabla|^{-1}P_k I^c\calN(P_{k_2,\kappa_2}
\psi_2, P_{k_3,\kappa_3} \psi_3)]
  \|_{N[\ell](I_{j}\times\R^{2})}^{2}
&<\eps_0 C_0^6 \label{eq:cruxA}\\
\max_{1\le j\le M_{1}}\,\sum_{\substack{k,\ell,k_1,k_2,k_3\\
\kappa_1,\kappa_2,\kappa_3}}     \|  \nabla_{t,x} P_{\ell}
[P_{k_1,\kappa_1} R\, \psi_1 |\nabla|^{-1}P_k
I\calN(P_{k_2,\kappa_2} \psi_2, P_{k_3,\kappa_3} \psi_3)]
   \|_{N[\ell](I_{j}\times\R^{2})}^{2} &<\eps_0 C_0^6
   \label{eq:cruxB}
\end{align}
where the sums extend over integers $k,\ell,k_1,k_2,k_3$ as
specified above, as well as over caps
$\kappa_1,\kappa_2,\kappa_3\in\calC_{m_0}$ with
$\dist(\kappa_i,\kappa_j)>2^{m_0}$ for $i\ne j$. Let us first
consider the case where the entire output is restricted by
$Q_{<2m_0+\ell}$ in modulation, and the inputs are in the hyperbolic
regime, i.e., $P_{k_i}\psi_i=Q_{\le k_i+C}P_{k_i}\psi_i$ where $C$
is large depending on~$L$. Then we bound~\eqref{eq:cruxA}
(and~\eqref{eq:cruxB}) as follows, first on the whole time axis~$\R$
(assuming as we may that the inputs have been suitably extended):
\begin{align}
 & \sum_{\substack{k,\ell,k_1,k_2,k_3\\
\kappa_1,\kappa_2,\kappa_3}} \| \nabla_{t,x} Q_{<2m_0+\ell-C}
P_{\ell} [P_{k_1,\kappa_1} \psi_1 |\nabla|^{-1}P_k
I^c\calN(P_{k_2,\kappa_2} \psi_2, P_{k_3,\kappa_3} \psi_3)]
  \|_{N[\ell ]}^{2} \nn \\
& \les \sum_{\substack{k,\ell,k_1,k_2,k_3\\
\kappa_1,\kappa_2,\kappa_3}} \sum_{\kappa\in\calC_{m_0}} \|
P_{\ell,\kappa}Q_{<2m_0+\ell-C} [P_{k_1,\kappa_1} \psi_1
|\nabla|^{-1}P_k I^c\calN(P_{k_2,\kappa_2} \psi_2, P_{k_3,\kappa_3}
\psi_3)]
  \|_{\NF[\kappa]}^{2} \label{eq:cruxC}
\end{align}
Note that the $2^{2\ell}$-factor produced by the output is canceled
against the scaling factor which is part of the~$N[\ell]$-norm, see
Definition~\ref{def:Nk}. By the usual arguments involving disposable
multipliers, we may replace $|\nabla|^{-1}P_k I^c$ by~$2^{-k_2}$
(implicit constants are allowed to depend on~$L$). Since the inputs
are hyperbolic, we may also ignore the null-form~$\calN$. For
any~$\kappa$,
\[
\max_{i=2,3} \dist(\kappa,\kappa_i)\gtrsim 2^{m_0}
\]
Let us assume that this happens for $i=3$. Then
by~\eqref{eq:bilin1}, followed by~\eqref{eq:bilin2},
\begin{align*}
  \eqref{eq:cruxC} &\le C(L,m_0) \sum_{\substack{k,\ell,k_1,k_2,k_3\\
\kappa_1,\kappa_2,\kappa_3}} 2^{-k_2}  \| P_{k_1,\kappa_1} \psi_1
P_{k_2,\kappa_2} \psi_2\|_{\Ltwotx}^2 \| P_{k_3,\kappa_3}
\psi_3\|_{S[k_3]}^2 \\
&\le C(L,m_0) \sum_{\substack{k,\ell,k_1,k_2,k_3\\
\kappa_1,\kappa_2,\kappa_3}}   \| P_{k_1,\kappa_1}
\psi_1\|_{S[k_1,\kappa_1]}^2 \| P_{k_2,\kappa_2}
\psi_2\|_{S[k_2,\kappa_2]}^2 \| P_{k_3,\kappa_3}
\psi_3\|_{S[k_3,\kappa_3]}^2\\
&\le C(L,m_0)  \Big(\sum_{k\in\Z} \|P_{k} \psi\|_{S[k]}^2\Big)^3 \le
C(L,m_0)C_0^6
\end{align*}
Note that we are not assuming that $P_{k_i,\kappa_i} \psi_i$  are
wave-packets, i.e., localized in modulation to $<2^{2m_0+k_i}$ but
only to modulations $<2^{k_i+C}$. Therefore, to pass to the last
line one needs to use Lemma~\ref{lem:square_func} for the
modulations between these two cut-offs. However, this only costs a
factor of~$\les|m_0|$ which is admissible. We now rewrite the first
line in this estimate in the form
\[
\eqref{eq:cruxC} \le C(L,m_0) \int_{\R^3} \sum_{\substack{k,\ell,k_1,k_2,k_3\\
\kappa_1,\kappa_2,\kappa_3}} 2^{-k_2}  | P_{k_1,\kappa_1}
\psi_1(t,x) P_{k_2,\kappa_2} \psi_2(t,x)|^2 \| P_{k_3,\kappa_3}
\psi_3\|_{S[k_3]}^2\, dtdx \le C(L,m_0)C_0^6
\]
By the dominated convergence theorem, we can cover the line (and
especially $(-T_0,T_1)$) into finitely many intervals $I_j$ such
that
\[ C(L,m_0) \int_{I_j\times\R^2} \sum_{\substack{k,\ell,k_1,k_2,k_3\\
\kappa_1,\kappa_2,\kappa_3}} 2^{-k_2}  | P_{k_1,\kappa_1}
\psi_1(t,x) P_{k_2,\kappa_2} \psi_2(t,x)|^2 \| P_{k_3,\kappa_3}
\psi_3\|_{S[k_3]}^2\, dtdx < \eps_0 C_0^6
\]
for each $I_j$. Moreover, the number of these intervals is~$\le
M_1(\eps)$.  Refining the intervals further if necessary, we may
similarly assume that
\begin{equation}
\label{eq:L2dxdt}
 C(L,m_0) \int_{I_j\times\R^2} \sum_{\substack{k,\ell,k_1,k_2,k_3\\
\kappa_1,\kappa_2,\kappa_3}} 2^{-k_2}  | P_{k_1,\kappa_1}
\psi_1(t,x) P_{k_2,\kappa_2} \psi_2(t,x)|^2 \| P_{k_3,\kappa_3}
\psi_3\|_{S[k_3]}^2\, dtdx < \eps_0 C_0^6
\end{equation}
for each $I_j$. Moreover, the number of these intervals can be taken
to depend only on~$C_0$ and~$\eps_0$ (since $L$ and $m_0$ have the
same property). Retracing our steps shows that these intervals have
the desired properties~\eqref{eq:cruxA} and~\eqref{eq:cruxB} under
the modulation assumptions $P_{k_i}\psi_i=Q_{\le
k_i+C}P_{k_i}\psi_i$, and the additional assumption that the output
is limited to size~$\les 2^{2m_0+\ell}$ (the Schwartz extensions
implicit in~\eqref{eq:cruxA} and~\eqref{eq:cruxB} are simply
obtained by multiplying the $L^2_{tx}$ functions by smooth bump
functions). The remaining cases where these modulation assumptions
are violated are handled similarly. For example,
consider~\eqref{eq:cruxB} for outputs of modulations~$\gtrsim
2^{\ell+C}$ but again on the whole time axis
\begin{align*}
& \sum_{\substack{k,\ell,k_1,k_2,k_3\\
\kappa_1,\kappa_2,\kappa_3}}     \|  \nabla_{t,x} P_{\ell} Q_{\ge
\ell+C} [P_{k_1,\kappa_1} R\, \psi_1 |\nabla|^{-1}P_k
I\calN(P_{k_2,\kappa_2} \psi_2, P_{k_3,\kappa_3} \psi_3)]
   \|_{N[\ell]}^{2} \\
&\les  \sum_{\substack{k,\ell,k_1,k_2,k_3\\
\kappa_1,\kappa_2,\kappa_3}}     \|  \nabla_{t,x} P_{\ell} Q_{\ge
\ell+C} [P_{k_1,\kappa_1} R\, \psi_1 |\nabla|^{-1}P_k
I\calN(P_{k_2,\kappa_2} \psi_2, P_{k_3,\kappa_3} \psi_3)]
   \|_{\dot X_{\ell}^{-\frac12+\eps ,-1-\eps,2}}^{2} \\
&\les  \sum_{\substack{k,\ell,k_1,k_2,k_3\\
\kappa_1,\kappa_2,\kappa_3}} \Big(\sum_{m\ge \ell+C}
2^{-(\frac12-\eps)\ell} 2^{(1-\eps)m}  2^{-\ell} \| P_{k_1,\kappa_1}
Q_m \psi_1\|_{\Ltwotx}\Big)^2 \: 2^{-2k} \| P_k
I\calN(P_{k_2,\kappa_2} \psi_2, P_{k_3,\kappa_3} \psi_3)
\|_{\Linf}^2\\
&\les  \sum_{\substack{k,\ell,k_1,k_2,k_3\\
\kappa_1,\kappa_2,\kappa_3}} 2^{-2\ell}  \| P_{k_1,\kappa_1}
\psi_1\|_{\dot X_{\ell}^{-\frac12+\eps ,1-\eps,2}}^2 \: 2^{k} \| P_k
I\calN(P_{k_2,\kappa_2} \psi_2, P_{k_3,\kappa_3} \psi_3)
\|_{\Ltwotx}^2\\
&\les \sum_{\substack{k,\ell,k_1,k_2,k_3\\
\kappa_1,\kappa_2,\kappa_3}} 2^{-2\ell}  \| P_{k_1,\kappa_1}
\psi_1\|_{\dot X_{\ell}^{-\frac12+\eps ,1-\eps,2}}^2 \: 2^{k} \|
\calL(P_{k_2,\kappa_2} \psi_2, P_{k_3,\kappa_3} \psi_3)
\|_{\Ltwotx}^2\\
\\& \le C(L,m_0) C_0^6
\end{align*}
where the final bound again follows from~\eqref{eq:bilin1} ($\calL$
stands for the usual averaged space-time translation operator which
arises via removal of disposable multipliers). Writing out the
$\Ltwotx$-norm explicitly in the previous estimate allows us again
to choose intervals~$I_j$ with the desired properties. The remaining
case of output modulations $Q_m$ with $2m_0+\ell\le m\le \ell+C$ is
similar:
\begin{align*}
& \sum_{\substack{k,\ell,k_1,k_2,k_3\\
\kappa_1,\kappa_2,\kappa_3}}     \|  \nabla_{t,x} P_{\ell} Q_{m}
[P_{k_1,\kappa_1} R\, \psi_1 |\nabla|^{-1}P_k
I\calN(P_{k_2,\kappa_2} \psi_2, P_{k_3,\kappa_3} \psi_3)]
   \|_{N[\ell]}^{2} \\
&\les  \sum_{\substack{k,\ell,k_1,k_2,k_3\\
\kappa_1,\kappa_2,\kappa_3}}     \|  \nabla_{t,x} P_{\ell} Q_{m}
[P_{k_1,\kappa_1} R\, \psi_1 |\nabla|^{-1}P_k
I\calN(P_{k_2,\kappa_2} \psi_2, P_{k_3,\kappa_3} \psi_3)]
   \|_{\dot X_{\ell}^{-\frac12+\eps ,-1-\eps,2}}^{2} \\
&\les  \sum_{\substack{k,\ell,k_1,k_2,k_3\\
\kappa_1,\kappa_2,\kappa_3}} 2^{-\ell}  \| P_{k_1,\kappa_1} Q_m
\psi_1\|_{\ener}^2 \:  \| P_k I\calN(P_{k_2,\kappa_2} \psi_2,
P_{k_3,\kappa_3} \psi_3) \|_{\Ltwotx}^2 \\& \le C(L,m_0) C_0^6
\end{align*}
Due to the $\Ltwotx$-norm one can now proceed as before. Finally,
suppose that the output as well as~$\psi_1$ are hyperbolic, but that
$\psi_2$ and~$\psi_3$ are elliptic. Then
\begin{align*}
& \Big(\sum_{\substack{k,\ell,k_1,k_2,k_3\\
\kappa_1,\kappa_2,\kappa_3}}     \|  \nabla_{t,x} P_{\ell} Q_{\le
\ell+C}  [P_{k_1,\kappa_1} Q_{\le k_1+C} R \psi_1 |\nabla|^{-1}P_k
I\calN(P_{k_2,\kappa_2} \psi_2, P_{k_3,\kappa_3} \psi_3)]
   \|_{N[\ell]}^{2}\Big)^{\frac12} \\
&\les  \sum_{\substack{k,\ell,k_1,k_2,k_3\\
\kappa_1,\kappa_2,\kappa_3}}  \sum_{m\ge k_2+C} 2^{-\ell}  \|
\nabla_{t,x} P_{\ell} Q_{\le \ell+C} [P_{k_1,\kappa_1} Q_{\le k_1+C}
R \psi_1 |\nabla|^{-1}P_k I\calN(P_{k_2,\kappa_2} Q_m \psi_2,
P_{k_3,\kappa_3} \tilde Q_m \psi_3)]
   \|_{\enerN}
   \end{align*}
   which can be further estimated as
   \begin{align*}
   &\les  \sum_{\substack{k,\ell,k_1,k_2,k_3\\
\kappa_1,\kappa_2,\kappa_3}}  \sum_{m\ge k_2+C} 2^{\ell} \|
P_{\ell} Q_{\le \ell+C} [P_{k_1,\kappa_1} Q_{\le k_1+C} R \psi_1
|\nabla|^{-1}P_k I\calN(P_{k_2,\kappa_2} Q_m \psi_2,
P_{k_3,\kappa_3} \tilde Q_m \psi_3)]
   \|_{\Leins} \\
   &\les  \sum_{\substack{k,\ell,k_1,k_2,k_3\\
\kappa_1,\kappa_2,\kappa_3}}  \sum_{m\ge k_2+C}   \|
P_{k_1,\kappa_1}  \psi_1\|_{\ener} \| P_k I\calN(P_{k_2,\kappa_2}
Q_m \psi_2, P_{k_3,\kappa_3} \tilde Q_m \psi_3)]
   \|_{\enerN} \\
     &\les  \sum_{\substack{k,\ell,k_1,k_2,k_3\\
\kappa_1,\kappa_2,\kappa_3}}  \sum_{m\ge k_2+C}   \|
P_{k_1,\kappa_1}  \psi_1\|_{\ener} \: 2^{m-k_2} \| P_{k_2,\kappa_2}
Q_m \psi_2\|_{\Ltwotx} \| P_{k_3,\kappa_3} \tilde Q_m \psi_3
   \|_{L^2_t L^\infty_x} \\
& \le C(L,m_0) C_0^6
\end{align*}
by Bernstein's inequality and the definition of $S[k]$; to pass to
the second line use that $P_{\ell} Q_{\le \ell+C}$ is disposable.
Breaking up the $\Leins$-norm in the third line of this estimate
into disjoint time intervals allows us to obtain the desired
conclusion as before. Alternatively, one can gain smallness here by
taking $C$ in $m\ge k_2+C$ large; this will be important later (see
Remark~\ref{rem:bilineargivesS}). We leave the analogous analysis
of~\eqref{eq:cruxA} to the reader.

\medskip\noindent
The proof of the claim \eqref{eq:Ijclaim} for the higher degree nonlinearities is outlined in the Appendix.

\medskip A crucial feature of the construction of the intervals
$\{I_{j}\}_{1\le j\le M_{1}}$ above is that is {\em universal}, i.e.,
it does not depend on the choice of the underlying frequency scale.
We now conclude the proof of Lemma~\ref{lem:LocalSplitting}. Fix some $I_j$ and localize
$\psi$ to frequency~$2^k$.
If $|I_j|<\eps_1\, 2^{-k}$, then $P_k \psi_{NL}^{(j)}:= P_k\psi$ satisfies the
bound~\eqref{eq:psiNLbd} by the analysis in Case~1. Otherwise,
 one represents the solution via~\eqref{eq:wave_solution}.
The bounds in Case~1 above then imply the estimate
\[
 \|(P_k \psi)\big|_{[t_j-\eps_1\,2^{-k},t_j-\eps_1\,2^{-k}]} \|_{S[k]} \les \|P_k\psi\|_{S[k]}
\]
 The free wave $P_k\psi_{L}$ at
dyadic frequency $2^{k}$ is now defined as the free evolution in~\eqref{eq:wave_solution}, whereas
$P_k \psi_{NL}^{(j)}$ is the sum of the other two terms in that formula. Summing over~$k$ now yields the claimed local splitting of $\psi$ in
Lemma~\ref{lem:LocalSplitting}.

\noindent Finally, the proof of~\eqref{eq:daniel} is implicit in the
preceding and we skip the details.
\end{proof}

\begin{remark}\label{rem:psiLevac}
Later we will apply \eqref{eq:daniel} in the following context. If
\[
 (\sum_{k>k_{0}}\|P_{k}\psi\|_{S[k]}^{2})^{\frac{1}{2}}<\delta_{3}
\]
for some (very small) $\delta_3>0$, we have
\[
 \|\nabla_{x,t}P_{>k_{0}}\psi_{L}^{(j)}\|_{\dot{H}_{x}^{-1}}+(\sum_{k>k_0}\|P_{k}\psi_{NL}^{(j)}\|_{S[k](I_{j}\times\R^{2})}^{2})^{\frac{1}{2}}\les \delta_{3}
\]
where the implied constant depends on $\|\psi\|_{S}$.
\end{remark}

\begin{remark}\label{rem:bilineargivesS}
The preceding proof can be easily modified to give the following
result that will be important later: Let $\psi$ be the gauged
derivative components of an admissible wave map. Assume that we have
an apriori bound of the form
\begin{equation}
 \label{eq:bilinLambda}
 \sum_{k_{1}>k_{2}}\sum_{\substack{\kappa_{1,2}\in \calC_{m_0}\\ \dist(\kappa_1,\kappa_2)\gtrsim 2^{-m_0}}} 2^{-k_2}\|P_{k_1,\kappa_1}\psi P_{k_2,\kappa_2}\psi\|_{L_{t,x}^{2}}^{2}+\sum_{k<l}[2^{(1-\eps)l-(\frac{1}{2}-\eps)k}\|P_kQ_l\psi\|_{L_{t,x}^2}]^{2}<\Lambda
\end{equation}
where $m_0$ is sufficiently large depending on the
energy~$E$ of~$\psi$. Then we can infer a bound of the form
\[
\|\psi\|_{S}\les C_2(E,\,m_0,\,\Lambda)
\]
This is done by a bootstrap, with the desired smallness coming either from the intervals $I_j$ or
the gains from the angular alignment.
Moreover, assume that for each $k_1>k_2$ one has
\[
\sum_{\substack{\kappa_{1,2}\in \calC_{m_0}\\
\dist(\kappa_1,\,\kappa_2)\gtrsim 2^{-m_0}}} P_{k_1,\kappa_1}\psi
P_{k_2,\kappa_2}\psi =f_{k_1,k_2} + g_{k_1,k_2}
\]
\[
P_kQ_{>k}\psi=h_k+i_k
\]
with for some positive integer $\nu$
\begin{equation}\begin{split}
 \label{eq:bilinLambda2}
 \sum_{k_{1}>k_{2}}
2^{-k_2}
\| f_{k_1,k_2}\|_{L_{t,x}^{2}}^{2}+\sum_{k<l}[2^{(1-\eps)l-(\frac{1}{2}-\eps)k}\|Q_l h_k\|_{L_{t,x}^2}]^{2}& <\Lambda\\
 \sum_{k_{1}>k_{2}}
2^{-k_2}
\| g_{k_1,k_2}\|_{L_{t,x}^{2}}^{2} +\sum_{k<l}[2^{(1-\eps)l-(\frac{1}{2}-\eps)k}\|Q_l i_k\|_{L_{t,x}^2}]^{2}&< \delta \|\psi\|^\nu_S
\end{split}\end{equation}
where $\delta>0$ is small depending on the energy and the integer~$\nu$, but independent of~$\|\psi\|_S$. Then one can again conclude
\[
\|\psi\|_{S}\les C_2(E,\,m_0,\,\Lambda)
\]
Note the the time intervals $I_j$ are determined only by means of the $f_{k_1,k_2}$ and not the $g_{k_1,k_2}$.
\end{remark}

\subsubsection{Proof of  Proposition~\ref{BlowupCriterion1}}

Recall that we are making the assumption~$\|\psi\|_{S}<C_0$.  We
first show that the wave map cannot break down in finite time, i.e.,
$T=T'=\infty$. Assume for example that $T<\infty$. For $\eps_0>0$ a
sufficiently small but absolute constant (which will be specified
later),  pick the $M_{1}(C_{0}, \eps_0)$-many intervals $I_{j}$ as
in Lemma~\ref{lem:LocalSplitting}. It will suffice to consider that
interval $I_{j_0}$ which has~$T$ as its endpoints. Alternatively,
starting with that interval $I_{j}$ containing the initial time
slice $t=0$, one can inductively obtain control over the
frequency-localized constituents of $\psi$, the $P_{k}\psi$.

\begin{lemma}\label{bootstrap2}
Let $I_{j}=(t_{j}, t_{j+1})$ be an interval as in
Lemma~\ref{lem:LocalSplitting}. Introduce the frequency envelope
\[
c_{k}:=(\sum_{\ell\in \Z}2^{-\sigma_0|k-\ell|}\|P_{\ell}\psi(t_{j},
\cdot)\|_{L_{x}^{2}}^{2})^{\frac{1}{2}}
\]
where $\sigma_0>0$ is some small constant.
Also, write $\psi|_{I_{j}}=\psi_{L}+\psi_{NL}$. Then there is a
number $C_{1}=C_{1}(\psi_{L})<\infty$ with the property that
\[
\|P_{k}\psi\|_{S[k](I_j\times \R^2)}\leq C_{1}c_{k},\quad\forall\:
k\in\Z
\]
\end{lemma}
\begin{proof}
 We prove this by splitting the interval $I_{j}$ into a finite number of smaller intervals
 depending on $\psi_{L}$. Thus we shall write
\[
I_{j}=\cup_{i}J_{ji}
\]
for a finite number of smaller intervals depending on $\psi_{L}$.
The exact definition of these intervals will be given later in the
proof. On each $J_{ji}$, we now run a bootstrap argument, commencing
with the {\em bootstrap assumption:}
\[ \|P_{k}\psi\|_{S[k](J_{j}\times\R^{2})}\leq A(C_{0})c_{k}\]
Here $A(C_{0})$ is a number that depends purely on the apriori
bound we are making on the wave map. We shall show that provided
$A(C_{0})$ is chosen large enough, the bootstrap assumption implies
the better bound
\[
\|P_{k}\psi\|_{S[k](J_{j}\times\R^{2})}\leq \frac{A(C_{0})}{2}c_{k}
\]
We prove this for each frequency mode. By scaling invariance, we may
assume $k=0$. As before, one needs to distinguish between
$|J_{j}|<\eps_1$ and the opposite case, where $\eps_1$ is chosen
sufficiently small. In the former case, one directly uses the
div-curl system
\[
\partial_{t}\psi=\nabla_{x}\psi+\psi\nabla^{-1}(\psi^{2})
\]
as in the previous section to obtain the desired conclusion for
$P_{0}\psi$. Thus we can assume that the interval satisfies
$|J_{j}|\geq \eps_1$, which means we can control
$\big(P_{0}\psi(t_{0},\cdot), P_{0}\partial_{t}\psi(t_{0},
\cdot)\big)$ for some $t_{0}\in J_{j}$ via
\[
\|\big(P_{0}\psi(t_{0},\cdot), P_{0}\partial_{t}\psi(t_{0},
\cdot)\big)\|_{L_{x}^{2}\times\dot{H}^{-1}}\lesssim
A_{1}(C_{0})c_{0}
\]
for some constant $A_{1}(C_{0})$, which is explicitly computable,
independently of $A(C_{0})$. Passing to the wave equation
\[
\Box
P_{0}\psi=P_{0}F_{\alpha}(\psi)=\sum_{i=1}^{5}P_{0}F_{\alpha}^{2i+1}(\psi)
\]
via Schwartz extensions and Hodge decompositions as before,  we
first consider the principal terms $P_{0}F_{\alpha}^{3}(\psi)$.
These terms  can be schematically  written as
\[
P_{0}\nabla_{x,t}[\psi\nabla^{-1}(\psi^{2})]
\]
More accurately, they are of the form
\begin{align*}
   \nabla_{t,x} P_{\ell}
[P_{k_1,\kappa_1} \psi_1 |\nabla|^{-1}P_k I^c\calN(P_{k_2,\kappa_2}
\psi_2, P_{k_3,\kappa_3} \psi_3)] \\
\nabla_{t,x} P_{\ell} [P_{k_1,\kappa_1} R\, \psi_1 |\nabla|^{-1}P_k
I\calN(P_{k_2,\kappa_2} \psi_2, P_{k_3,\kappa_3} \psi_3)]
\end{align*}
with a Riesz projection $R$ and a nullform~$\calN$.
 Substituting the
decomposition $\psi=\psi_{L}+\psi_{NL}$ into the inner null-form
yields
\begin{equation}\nonumber\begin{split}
&P_{0}\nabla_{x,t}[\psi\nabla^{-1}(\psi^{2})]
=P_{0}\nabla_{x,t}[\psi\nabla^{-1}(\psi_{L}^{2})]
+P_{0}\nabla_{x,t}[\psi\nabla^{-1}(\psi_{L}\psi_{NL})]
+P_{0}\nabla_{x,t}[\psi\nabla^{-1}(\psi_{NL}^{2})]
\end{split}\end{equation}
Note that the last term automatically has the desired smallness
property if we choose $\eps_0$ smaller than some absolute constant.
Indeed, by~\eqref{eq:psiNLbd}, and the trilinear estimates of
Section~\ref{sec:trilin},
\[
\|P_{0}\nabla_{x,t}[\psi\nabla^{-1}(\psi_{NL}^{2})]\|_{N[0]}\lesssim
\eps_0^{2}\sup_{k\in\Z} 2^{-\sigma_0\,
|k|}\|P_{k}\psi\|_{S[k]}\lesssim \eps_0^{2}A(C_{0})c_{0} \ll
A(C_{0})c_{0}
\]
for small~$\eps_0$.  Next, for the mixed term
$P_{0}\nabla_{x,t}[\psi\nabla^{-1}(\psi_{L}\psi_{NL})]$, choosing
$\eps_0$ sufficiently small (depending on $C_{0}$), we can arrange
in light of Lemma~\ref{lem:LocalSplitting} and the trilinear
estimates
\[
\|P_{0}\nabla_{x,t}[\psi\nabla^{-1}(\psi_{L}\psi_{NL})]\|_{N[0]}\lesssim
A(C_0) C_0^3  \eps_0^{1-\frac{1}{M}}\; c_0 \ll A(C_0)\, c_0
\]
The first term
\[
P_{0}\nabla_{x,t}[\psi\nabla^{-1}(\psi_{L}^{2})]
\]
requires a separate argument. In fact, we treat this term by
decomposing the interval $I_{i}$ into smaller ones. In order to
select these intervals, first note that upon localizing the
frequencies of the inputs according to
\[
P_{0}\nabla_{x,t}[P_{k_{1}}\psi\nabla^{-1}P_{k}(P_{k_{2}}\psi_{L}P_{k_{3}}\psi_{L})]
\]
one obtains from the trilinear bounds of Section~\ref{sec:trilin}
\[
\|P_{0}\nabla_{x,t}[P_{k_{1}}\psi\nabla^{-1}P_{k}(P_{k_{2}}\psi_{L}P_{k_{3}}\psi_{L})]\|_{N[0]}\leq
\eps_0\, 2^{-\sigma |k_1|}\|P_{k_{1}}\psi_{1}\|_{S[k_{1}]}\ll
A(C_0)\, c_0
\]
in the following two cases:  $k_1,k_2,k_3$ fall outside the
range~\eqref{eq:outside} (the high-low-low case), or,  if they do
fall in the range~\eqref{eq:outside}, then $k\le k_2-L'$.  Here $L$
and~$L'$ are large constants depending on $C_{0}, \eps_0$, due to
the bounds on $\psi_{L}$ from Lemma~\ref{lem:LocalSplitting}. Thus,
denoting by
\[
P_{0}\nabla_{x,t}[\psi\nabla^{-1}(\psi_{L}^{2})]'
\]
the sum over all frequency interactions described by these
conditions, one then obtains the estimate
\[
\|P_{0}\nabla_{x,t}[\psi\nabla^{-1}(\psi_{L}^{2})]'\|_{N[0]}\ll
A(C_{0})\, c_{0}
\]
Employing the notations of Section~\ref{subsec:improvetrilin}, it
thus suffices to consider the sum of expressions
\[
{\sum}''_{k_1,k_2,k_3\in\Z} \sum_{k=k_2-L'}^{k_2+O(1)}
P_{0}\nabla_{x,t}[P_{k_{1}}\psi\nabla^{-1}P_{k}(P_{k_{2}}\psi_{L}P_{k_{3}}\psi_{L})],
\]
where, of course, the implied constants may be quite large depending
on $C_{0}, \eps_0$. Furthermore, by the results of that section,
we may assume that the inputs have pairwise angular separation on
the Fourier side, and in particular we make this assumption for the
free wave inputs $P_{k_{2}}\psi_{L}$ and $P_{k_{3}}\psi_{L}$. Thus
we have now reduced ourselves to estimating
\[
\sum_{\substack{\kappa_1,\kappa_2,\kappa_3\in\calC_{m_0} \\
\max_{i\ne j}\dist(\kappa_i,\kappa_j)>2^{m_0}}}
{\sum}''_{k_1,k_2,k_3\in\Z} \sum_{k=k_2-L'}^{k_2+O(1)}
P_{0}\nabla_{x,t}[P_{k_{1},\kappa_1}\psi\nabla^{-1}P_{k}(P_{k_{2},\kappa_2}\psi_{L}\:
P_{k_{3},\kappa_3}\psi_{L})],
\]
The next step is to exploit the dispersive properties of the
expression
\[
\nabla^{-1}P_{k}(P_{k_{2},\kappa_{2}}\psi_{L}\:
P_{k_{3},\kappa_{3}}\psi_{L})
\]
First, due to the energy bound for $\psi_{L}$, there exists some
finite set $A\subset\Z$ so that
\begin{align*}
&\sum_{\substack{\kappa_1,\kappa_2,\kappa_3\in\calC_{m_0} \\
\max_{i\ne j}\dist(\kappa_i,\kappa_j)>2^{m_0}}}
{\sum}''_{\substack{k_1,k_2,k_3\in\Z\\ k_2\not\in A}}
\sum_{k=k_2-L'}^{k_2+O(1)} \|
P_{0}\nabla_{x,t}[P_{k_{1},\kappa_1}\psi\nabla^{-1}P_{k}(P_{k_{2},\kappa_2}\psi_{L}\:
P_{k_{3},\kappa_3}\psi_{L})]\|_{N[0]}   \\
&\leq \eps_02^{-\sigma_0\, |k_{1}|} \|P_{k_{1}}\psi\|_{S[k_{1}]}
\end{align*}
On the other hand, assume now that $k_{2}\in A$ and consider
\[
\nabla^{-1}P_{k}(P_{k_{2},\kappa_{2}}\psi_{L}\:P_{k_3,\kappa_{3}}\psi_{L})
\]
where $k,k_3$ are chosen as in~\eqref{eq:outside}. Note that the set
$A$ depends on the dyadic frequency of the output, in this case
frequency $2^{0}$. Changing the  frequency localizations of the
output amounts to a rescaling of~$A$. Nonetheless, one has the
following estimates which are independent under rescaling:
\[
\|\nabla^{-1}P_{k}(P_{k_{2},\kappa_{2}}\psi_{L}P_{k_3,\kappa_{3}}\psi_{L})\|_{L_{t}^{1}L_{x}^{\infty}}<C_{3}(\psi_{L},
k_{2})
\]
In particular,
\[
\sum_{k_{2}\in
A}\|\nabla^{-1}P_{k}(P_{k_{2},\kappa_{2}}\psi_{L}P_{k_3,\kappa_{3}}\psi_{L})\|_{L_{t}^{1}L_{x}^{\infty}}<C_{4}(\psi_{L})<\infty
\]
To prove this bounds, set $k_{2}=0$  by scaling invariance. But then
$P_{k_{2}}\psi_{L}(t_{j}, \cdot)$ is a Schwartz function in the
$x$-variable. Using the angular separation of the inputs it is now
straightforward to see  that
\[
\|\nabla^{-1}P_{k}(P_{k_{2},\kappa_{2}}\psi_{L}\:P_{k_3,\kappa_{3}}\psi_{L})\|_{L_{t}^{1}L_{x}^{\infty}}<C_{3}(\psi_{L},
k_{2})
\]
Indeed, this follows from stationary phase and the angular
separation of the inputs. We now define the intervals $J_{j}$ by
requiring that
\[
\sum_{\substack{\kappa_1,\kappa_2,\kappa_3\in\calC_{m_0} \\
\max_{i\ne j}\dist(\kappa_i,\kappa_j)>2^{m_0}}}
\sum_{\substack{k_2\in A\\ |k_2-k_3|<L}} \sum_{k=k_2-L'}^{k_2+O(1)}
\| \nabla^{-1}P_{k}(P_{k_{2},\kappa_2}\psi_{L}\:
P_{k_{3},\kappa_3}\psi_{L})]\|_{L_{t}^{1}L_{x}^{\infty}(J_{j}\times\R^{2})}<\eps_0
\]
It is furthermore clear that we also obtain
\begin{equation}\nonumber\begin{split}
&\sum_{\substack{\kappa_1,\kappa_2,\kappa_3\in\calC_{m_0} \\
\max_{i\ne j}\dist(\kappa_i,\kappa_j)>2^{m_0}}}
{\sum}''_{\substack{k_1,k_2,k_3\in\Z\\ k_2\in A}}
\sum_{k=k_2-L'}^{k_2+O(1)} \|
P_{0}\nabla_{x,t}[P_{k_{1},\kappa_1}\psi\nabla^{-1}P_{k}(P_{k_{2},\kappa_2}\psi_{L}\:
P_{k_{3},\kappa_3}\psi_{L})]\|_{N[0](J_{j}\times\R^{2})}
\\
&\leq \eps_02^{-\sigma_0\, |k_{1}|}\|P_{k_{1}}\psi\|_{S[k_{1}]}\ll
A(C_{0})c_{0}
\end{split}\end{equation}
which completes the bootstrap for the trilinear source terms.

\medskip
The contribution of the higher order terms is dealt with in the appendix.

\medskip
By applying the above bootstrap argument on each of the finitely
many intervals $J_{ji}$ comprising each $I_{j}$, the proof of
Lemma~\ref{bootstrap2} now follows.
\end{proof}

The proof of Proposition~\ref{BlowupCriterion1} can easily be
concluded. Indeed, one infers from the Klainerman-Machedon criterion
that $T=T'=\infty$. Moreover, we obtain a global apriori bound
\[
\|P_{k}\psi\|_{S[k]}\leq C_{1}c_{k}
\]
were the constant $C_{1}$ depends implicitly on $\psi_{L}$.
Unfortunately we have no apriori way of controlling this number.
Moreover, Lemma~\ref{lem:LocalSplitting} implies the scattering for
large times.

\subsection{Control of wave-maps via a fixed $L^2$-profile}\label{subsec:perturbprofile}

A fundamental issue that we need to address is the {\it{very
definition}} of wave maps with data that are in some sense only of
energy class. To propagate such data under the wave map evolution,
we shall use approximations by smooth wave maps each of  which can
be continued canonically. The following lemma justifies this
procedure.

\begin{lemma}\label{BasicStability}
Let $\phi^{n}_{\alpha}$ be the derivative components of a sequence
of Schwartz class\footnote{In the usual sense that
$\phi^{n}_{\alpha}|_{t=const}$ is Schwartz on $\R^{2}$.} wave maps
$\bfu^{n}\::\: (-T_0^n,T_1^n)\times \R^{2}\to \Hyp^{2}$ on their maximal time interval of existence
and assume that the Coulomb
components $\psi^{n}_{\alpha}(0,\cdot)$ satisfy
\[
\lim_{n\to\infty}\|\psi^{n}_{\alpha}(0,
\cdot)-V_{\alpha}\|_{L_{x}^{2}}=0
\]
for some $V_{\alpha}\in L^2(\R^2)$. Denoting the collection of
components $V_{\alpha}$ by $V$, there is a time $T_0=T_{0}(V)>0$ such
that $\min(T_0^n,T_1^n)>T_0$ for all sufficiently large~$n$ and
\[
\limsup_{n\to\infty}\|\psi^{n}_{\alpha}\|_{S((-T_0,T_0)\times\R^{2})}\leq
C(V)<\infty
\]
Furthermore,
there is a constant $C_{1}(V)$ with the following property: defining
the frequency envelope
\[
c_{k}^{(n)}:=\max_{\alpha=0,1,2} (\sum_{\ell\in
\Z}2^{-\sigma|k-\ell|}\|P_{\ell}\psi^{n}_{\alpha}\|_{L_{x}^{2}}^{2})^{\frac{1}{2}}
\]
for sufficiently small fixed $\sigma >0$, one has for all $k\in \Z$ and all large~$n$
\[
\max_{\alpha=0,1,2} \|P_{k}\psi^{n}_{\alpha}\|_{S[k]((-T_0,T_0)\times\R^{2})}\leq
C_{1}(V)c_{k}^{(n)}
\]
Finally, the wave map propagations of the $\psi^{n}_{\alpha}$
converge on fixed time slices $t=t_{0}\in [-T_{0}, T_{0}]$ in the
$L^{2}$-topology, uniformly in time.
\end{lemma}

The proof of this lemma will occupy the remainder of this section.
Before we begin with the proof, we discuss some related results and
implications of Lemma~\ref{BasicStability}.
Most fundamental is the following stability result:

\begin{prop}
 \label{prop:ener_stable} Let $\bfu:[-T_0,T_1]\times \R^{2}\to\Hyp^2$ be an admissible wave-map with gauged derivative components denoted by~$\psi$.
Assume that $\|\psi\|_{S([-T_0,T_1]\times \R^2)} =A<\infty$. Then
there exists $\eps_1=\eps_1(A)>0$ with the following property: any
other admissible wave-map $\bfv$ defined locally around $t=0$ and
with gauged derivative components $\tilde\psi$ satisfying
$\|\psi(0)-\tilde\psi(0)\|_2< \eps<\eps_1$ extends as an admissible
wave-map to $[-T_0,T_1]$ and satisfies
$\|\tilde\psi\|_{S([-T_0,T_1]\times \R^2)} < A + c(\eps_1)$ where
$c(\eps)\to 0$ as $\eps\to0$.
\end{prop}
\begin{proof}
The proof will be given in Section~\ref{sec:missing proofs}, as it
follows directly from the proof of Proposition~\ref{PsiBootstrap}.
\end{proof}

As a consequence, one has the
following important continuation result.
\begin{cor}
\label{cor:contS} Let $\{\psi^n\}_{n=1}^\infty$ be a sequence of
Coulomb components of admissible wave maps $\bfu^n: I\to\Hyp^2$
where $I$ some fixed nonempty closed interval such that for some
$t_0\in I$ one has
\[
\lim_{n\to\infty}\|\psi^{n}_{\alpha}(t_0,
\cdot)-V_{\alpha}\|_{L_{x}^{2}}=0
\]
with $V_{\alpha}\in L^2(\R^2)$ as well as
\[
\sup_n \|\psi^n\|_{S(I\times\R^2)}<\infty.
\]
Then there exists a true extension $\tilde I$ of $I$ (meaning that
it extends by some positive distance beyond the endpoints of~$I$ in
sofar as they are finite) to which each $\bfu^n$ can be continued as
an admissible wave map provided $n$ is large.
\end{cor}
\begin{proof}
By Proposition~\ref{prop:ener_stable} we can define $\lim_{n\to\infty}\psi^{n}_\alpha(t,.)$ in the $L^{2}$-sense for $t$ an endpoint of $I$.
 By Lemma~\ref{BasicStability} the $\psi^{n}_\alpha$ extend beyond the (finite) endpoints for $n$ large enough.
\end{proof}

We can use the preceding results to define wave maps with
$L^2$ data at the level of the Coulomb gauge.

\begin{defi}\label{CoulombWM}
Assume we are given a family $\{V_\alpha\}$, $\alpha=0, 1, 2$, of $L^{2}(\R^{2})$-functions, to be interpreted as data at time $t=0$. Also, assume we have
\[
 V_\alpha = \lim_{n\to\infty}\psi^{n}_{\alpha}
\]
where $\{\psi^{n}_\alpha\}$ are Coulomb components of admissible wave maps at time $t=0$. Determine $I=(-T_0, T_1)$
to be the maximal open time interval with the property that
\[
 \sup\{\tilde{I}\subset I,\,\text{$\tilde{I}$ closed} \,|\,\liminf_{n\rightarrow\infty}\|\psi^{n}_{\alpha}\|_{S(\tilde{I}\times\R^{2})}<\infty\}
\]
Then we define the Coulomb wave maps propagation of $\{V_\alpha\}$ to be
\[
 \Psi^\infty_\alpha(t, x):=\lim_{n\to\infty}\psi^{n}_\alpha(t,x),\quad t\in I
\]
We call $I\times\R^{2}$ the lifespan of the (Coulomb) wave maps evolution of $\{V_\alpha\}$.
\end{defi}

It is of course important that the life span does not depend on the choice of sequence and, moreover, that the ``solutions'' $V_\alpha$  are unique.
These statements follow from Proposition~\ref{prop:ener_stable}.

The aforementioned uniqueness properties are now immediate -- indeed, simply mix any two sequences which converge to~$V_\alpha$.
Moreover, we can characterize the life-span as follows.

\begin{cor}
 \label{cor:lifespan}
Let $V_\alpha$, $\{\psi_\alpha^n\}$, and $I$ be as in Definition~\ref{CoulombWM}. Assume in addition that $I\neq (-\infty, \infty)$.
Then
\begin{equation}
 \label{eq:Iblow} \sup_{\substack{J\subset I\\ J \text{\ closed}}} \liminf_{n\to\infty} \|\psi_\alpha^n\|_{S(J\times \R^2)} =\infty
\end{equation}
\end{cor}
\begin{proof}
Suppose not. Let $I=(-T_0, T_1)$ where w. l. o. g. we assume $T_1<\infty$. Then there exists a number $M<\infty$ with the property that
for every closed $J\subset I$ with $0\in J$ one has
 \[
  \liminf_{n\to\infty}\|\psi_\alpha^n\|_{S(J\times \R^2)} < M
 \]
Now observe that
\[
 \limsup_{n,\,J\subset I}\|\psi^{n}_\alpha\|_{S(J\times\R^{2})}=\infty,
\]
where $J$ ranges over the closed subsets of $I$. Indeed, if not, we have
\[
\limsup_{n\to\infty}\|\psi^n_\alpha\|_{S([0,T_1]\times\R^2)}<\infty
\]
But then by Corollary~\ref{cor:contS} one can extend
$\psi^{n}_\alpha$ beyond the endpoint $T_1$ of $I$ to some interval
$\tilde{I}$ for $n$ large enough while maintaining the finiteness of
$\|\psi^{n}_\alpha\|_{S(\tilde{I}\times\R^{2})}$, contradicting the
definition of $I$.
\\
Now pick $\epsilon_1$ as in Proposition~\ref{prop:ener_stable}, with $M$ replacing $A$, and pick $J\subset I$, $n_0$ large enough such that
\[
\|\psi^{n_0}_\alpha\|_{S(J\times\R^{2})}\gg M,\quad \sup_{n,m\geq n_0}\|(\psi^{n}_\alpha-\psi^{m}_\alpha)(0,\,\cdot)\|_{L^{2}}<\epsilon_1
\]
But by our definition of $M$ there exists $k_0>n_0$ with the property that
\[
 \|\psi^{k_0}_\alpha\|_{S(J\times\R^{2})}<M
\]
and then applying Proposition~\ref{prop:ener_stable} to $\psi^{k_0}_\alpha(0,\cdot)$ we obtain a contradiction. This proves the corollary.
\end{proof}

Another important property is to be able to ensure the apriori existence of wave maps flows ``at infinity'', i.e.,
the solution of the scattering problem. In this regard, we have the following result.

\begin{prop}\label{prop:waveops}
Assume we are given admissible data at time $t=0$ of the form
\[
 \psi_\alpha=\partial_\alpha \big(S(0-t_0)(\partial_t V, V)\big)+o_{L^{2}}(1),\quad \alpha=0,\,1,\,2
\]
Here $(\partial_t V, V)\in L^{2}\times\dot{H}^{1}$ is a fixed profile. Then for $t_0=t_0(\partial_t V, V)>0$ large enough and
$
o_{L^{2}}(1)
$
small enough, the wave map associated with $\psi_\alpha$ exists on $(-\infty, 0]$, is admissible there, and we have
\[
 \|\psi_\alpha\|_{S((-\infty, 0]\times\R^{2})}<\infty
\]
Moreover, letting $\psi^{n}_\alpha$ be a sequence of admissible Coulomb components (i.e., associated with admissible maps) at time $t=0$ satisfying
\[
 \psi^{n}_\alpha\to \partial_\alpha \big(S(0-t_0)(\partial_t V, V)\big)
\]
for $(\partial_t V, V)$ as before and $t_0$ large enough also as before, the limit
\[
 \lim_{n\to\infty}\psi^{n}_\alpha(t, x)=\Psi^{\infty}_\alpha(t, x),\quad t\in(-\infty, 0]
\]
exists independently of the particular sequence chosen. We call this the Coulomb wave maps evolution of the data
\[
\partial_\alpha \big(S(0-t_0)(\partial_t V, V)\big)
\]
at time $t=-\infty$.  A similar construction applies at time $t=\infty$.
\end{prop}

\begin{cor}\label{temporallyunbounded} Assume that for a sequence of admissible Coulomb components $\psi^{n}_\alpha$ at time $t=0$ we have
\[
 \psi^{n}_\alpha= \partial_\alpha \big(S(t-t^{n})(\partial_{t}V,\,V)\big)+o_{L^{2}}(1)
\]
Then if $t_n\to\infty$, the limits
\[
 \lim_{n\to\infty}\psi^{n}_\alpha(t+t^n,\,x)=\Psi^{\infty}_\alpha(t,\,x)
\]
exist in the $L^{2}$-sense on some interval $(-\infty, -C)$, uniformly on closed subintervals, for $C$ large enough. We have
\[
 \limsup_{n\to\infty}\|\psi^{n}_\alpha(t+t^n,\,x)\|_{S((-\infty, -C_0]\times\R^{2})}<\infty
\]
for $C_0>C$. We call the maximal interval $I=(-\infty, -C)$ for which these statements hold the lifespan of the limiting
object $\Psi^\infty_\alpha$; here $C$ may be negative or $-\infty$.
A similar construction applies when $t_n\to-\infty$.
 \end{cor}

Both Proposition~\ref{prop:waveops} as well as
Corollary~\ref{temporallyunbounded} will be proved in
Section~\ref{subsec: scatteringwellposed}.
Having defined limiting objects $\Psi^{n}_\alpha$ as in Lemma~\ref{CoulombWM} (temporally bounded case) as well as Corollary~\ref{temporallyunbounded},
we can now define in obvious fashion the norms
\[
 \|\Psi^\infty_\alpha\|_{S(J\times\R^{2})}=\lim_{n\to\infty}\|\psi^{n}_\alpha\|_{S(J\times\R^{2})}
\]
for $J\subset I$ closed, with $I$ the lifespan of the limiting object. This is well-defined due to Proposition~\ref{prop:ener_stable}.
We can then also state the following
\begin{lemma}
 \label{lem:lifespan2}
Let $\Psi^\infty_\alpha$ be as before, with lifespan $I$.  Assume in addition that $I\neq (-\infty, \infty)$.
Then
\begin{equation}
 \label{eq:Iblowinfty} \sup_{\substack{J\subset I\\ J \text{\ closed}}}  \|\Psi_\alpha^{\infty}\|_{S(J\times \R^2)} =\infty
\end{equation}
The same conclusion holds for arbitrary $I$ provided the sequence $\psi^{n}_\alpha$ is essentially singular.
\end{lemma}

We now turn to the proof of Lemma~\ref{BasicStability}.
We begin with the the lower bound on the life span of the
$\psi^{n}_{\alpha}$. In essence, this is a consequence of the fact
that $\psi^{n}_{\alpha}\to V_{\alpha}$ in $L^{2}$ implies a uniform
non-concentration property of the energy of the $\psi^{n}_{\alpha}$.
This then allows one to approximate the corresponding wave maps with
derivative components $\phi^{n}_{\alpha}$ on small discs
--- with radii depending only on the limiting ``profile''~$V$ --- by small energy smooth wave
maps;  the small energy theory and finite propagation speed then
imply a uniform lower bound on the life span. Technically speaking,
restricting to small scales requires some care since localizing the
wave map by applying a smooth cutoff does not necessarily decrease
the energy. To see this, let $\chi$ be a cutoff to a small ball~$B$
of size~$r$. Then the first term on the right-hand side of
\begin{equation}\label{eq:leib}
\int_{\R^2} |\nabla (\chi \phi)(x)|^2\, dx \les  r^{-2}\int_{B}
|\phi(x)|^2\, dx + \int_{B} |\nabla \phi(x)|^2\, dx
\end{equation}
does  in general not become small as $r\to0$.

\noindent Let $\eps_0>0$ be the cutoff such that smooth data with
energy less than  $\eps_0$ result in global wave maps. More
precisely, we will rely on the following result by the first
author, see~\cite{Krieger}.

\begin{theorem}
  \label{thm:smallenergy}
Given smooth initial data $({\bf{x}},\,{\bf{y}})[0]:
0\times\R^{2}\to \Hyp^{2}$ which are sufficiently small in the sense
that
\begin{equation}\nonumber
\int_{0\times\R^{2}}\sum_{\alpha=0}^{2}\big[
\big(\frac{\partial_{\alpha}{\bf{x}}}{{\bf{y}}}\big)^{2}+\big(\frac{\partial_{\alpha}{\bf{y}}}{{\bf{y}}}\big)^{2}
\big]\,dx_1 dx_2 <\eps_0^2
\end{equation}
where $\eps_0>0$ is a small absolute constant, there exists a unique
classical wave map from $\R^{2+1}$ to $\Hyp^{2}$ extending these
data  globally in time. Moreover, one has the bound
$\sum_{\alpha=0}^2 \|\psi_\alpha\|_{S(\R^{1+2})} \le C\eps_0$
where~$C$ is an absolute constant.
\end{theorem}

Denoting the actual map at time $t=0$ giving rise (together with the
time derivatives) to $\phi^{n}_{\alpha}, \psi^{n}_{\alpha}$, by
$({\bfx}, {\bfy})(0,\cdot): \R^{2}\to \Hyp^2$, where we
have omitted the superscript~$n$ for simplicity, we now consider a
``re-normalized'' map, subject to a choice of $x_0\in\R^2$ and $r_0>0$,
\begin{equation}\label{eq:x1y1}
({\bfx}_{1}, {\bfy}_{1}):=
\Big(\chi_{[|x-x_{0}|<r_{0}]}\frac{{\bfx}-{\bfx}_{0}}{{\bfy}_{0}},\:e^{\chi_{[|x-x_{0}|<2r_{0}]}\log[\frac{\bfy}{{\bfy}_{0}}(0,\cdot)]}
\Big)
\end{equation}
Here $\chi_{[|x-x_{0}|<r_{0}]}$ is a smooth cutoff to the
disk~$D_{x_{0}, r_{0}}:=\{|x-x_0|<r_0\}$ which equals one on
$|x-x_{0}|<\frac{r_{0}}{2}$, say, and
\begin{align*}
\bfx_0&:= \slashint_{[|x-x_{0}|<r_{0}]} \bfx(x)\, dx_1dx_2,\qquad
 \bfy_0 :=
\exp\big(\slashint_{[|x-x_{0}|<2r_{0}]} \log \bfy(x)\, dx_1dx_2\big)
\end{align*}
with $\slashint_B:=|B|^{-1}\int_B$. Note that we have chosen the
cutoffs on the two components differently -- the one on the second
component is slightly larger than the first. This is merely a
technical convenience which amounts to $\bfy_1=\frac{\bfy}{\bfy_0}$
when $\nabla\chi_{[|x-x_{0}|<r_{0}]}\ne0$.
Lemma~\ref{lem:DataLocalization} below verifies the desired
smallness of energy property for these data. We begin with a basic
imbedding lemma which we shall need in the proof of that lemma. Even
though we only require the case $d=2$, we formulate this lemma in
any dimension.

\begin{lemma}
  \label{lem:embed} $\dot{B}^{\frac{d}{2}}_{2,\infty}(\R^d)\hookrightarrow
  \bmo(\R^d)$.
\end{lemma}
\begin{proof}
  By duality, it suffices to prove that $\calH^1\hookrightarrow
  \dot{B}^{-\frac{d}{2}}_{2,1}(\R^d)$ for the Hardy space~$\calH^1$. Thus, we need to show that
\[
\sum_{j\in\Z} 2^{-j\frac{d}{2}} \|P_j \phi\|_2\le C(d)
\]
for any $\phi$ which is an
  $\calH^1(\R^d)$ atom. Here $P_j$ are the usual Littlewood-Paley projections to frequencies of size~$2^j$.
  By scaling and translation invariance we may
  assume that $\supp(\phi)\subset B(0,1)$, $|\phi|\le 1$ and
  $\int\phi(x)\,dx=0$. If $j\ge0$, then we use that
  \[
\|P_j \phi\|_2\le \|\phi\|_2\le C(d)
  \]
If $j\le 0$, then writing $ P_j \phi = 2^{jd}\psi(2^j\cdot)\ast
\phi$ we conclude that
\[
|P_j\phi(x)|\le C\int_{B(0,1)} 2^{j(d+1)} \int_0^1 |\nabla \psi|(2^j
(x-ty))\,dt|\phi(y)|\,dy
\]
which implies that
\[
\|P_j\phi\|_2 \le C\int_{B(0,1)} 2^{j(d+1)} \|\nabla\psi\|_2
2^{-j\frac{d}{2}}|\phi(y)|\,dy\le C(d)2^{j(\frac{d}{2}+1)}
\]
and we are done.
\end{proof}

The importance of $\bmo$ in this context lies with the fact that
exponentiation maps small balls in~$\bmo$ into the $A_p$-class.
Recall that  $w$ is an $A_p$-weight in the sense of Muckenhaupt (see
Chapter~7 in~\cite{Duo} or Stein~\cite{Stein}  for all this standard material) provided
  \begin{equation}\label{eq:Ap}
 |Q|^{-1}\int_Q w(x)\,dx \; \Big(|Q|^{-1}\int_Q
 w^{1-p'}(x)\,dx\Big)^{p-1} \le A_p(w)
  \end{equation}
uniformly for all cubes~$Q\subset\R^d$ for some constant $A_p(w)$. Here $1<p<\infty$ and
$p'=\frac{p}{p-1}$ as usual. Note that $A_p\subset A_q$ if $p\le q$.
The $A_1$ class is defined as all $w$ with $Mw\le Cw$ a.e., where $M$ is the Hardy-Littlewood
maximal operator.  At the other end one has
$A_\infty:=\bigcup_{1\le p<\infty} A_p$, which  is characterized by the
estimate
\begin{equation}\label{eq:Ainfty}
\frac{w(S)}{w(Q)}\le C \Big(\frac{|S|}{|Q|}\Big)^\delta
\end{equation}
for all $S\subset Q$ (this is deep and requires the ``reverse
H\"older inequality''). Here $C$ and $\delta>0$ only depend on~$w$.
From the John-Nirenberg inequality, $w=e^{\phi}$ is an $A_p$ weight
for some $1<p<\infty$ provided $\|\phi\|_{\bmo}<r_0$ is small enough
and the $A_p$-constant $A_p(w)$ in~\eqref{eq:Ap} only depends
on~$r_0$.

\begin{lemma}
 \label{lem:bmoAp}  Let $\|\varphi\|_{L^2}\le A$ and set $w:= e^{(-\Delta)^{-\frac12} \varphi}$.
Then for any $1<p<\infty$ one has $A_p(w)\le C(p,A)$ where the latter constant only depends on~$p$ and~$A$.
\end{lemma}
\begin{proof}
For any $\delta>0$,
\[
\#\{j\in\Z\:|\: \|P_j \vphi \|_2 \ge \delta\}\le \delta^{-2}A^2
\]
In particular, for any $\delta>0$ there is a decomposition $
\varphi = \wt\vphi + (\vphi-\wt\vphi) $ so that
$\|(-\Delta)^{-\frac12}\wt\vphi\|_\infty\le  C \delta^{-2}A^3$
and, by Lemma~\ref{lem:embed},
\[
\|(-\Delta)^{-\frac12}(\vphi-\wt\vphi)\|_{\bmo}\le C\delta
\]
By the John-Nirenberg inequality one may choose $\delta$ small depending on~$p\in(1,\infty)$ such that
$\exp\big((-\Delta)^{-\frac12}(\vphi-\wt\vphi) \big)\in A_p$ with some absolute $A_p$-constant.
Since
\[
 \big\| e^{(-\Delta)^{-\frac12}\wt\vphi} \big\|_\infty \le e^{C \delta^{-2}A^3}
\]
we are done.
\end{proof}

\noindent The importance of $A_p$ weights lies with the fact that the
Hardy-Littlewood maximal operator $M$ as well as Calderon-Zygmund
operators $T$ are bounded on $L^p(w\,dx)$ with constants that only
depend on the dimension and the $A_p$ constant from~\eqref{eq:Ap}
(and $T$ in case of a singular integral) provided $1<p<\infty$.
In the present context, we will require a version of Poincar\'e's inequality
with~$A_2$ weights.
Now for the small energy lemma.

\begin{lemma}\label{lem:DataLocalization} Let $(\bfx_1^n,\bfy_1^n)$ be
as in~\eqref{eq:x1y1} applied to~$(\bfx^n,\bfy^n)$. Then
 given $\eps_0>0$, we can pick $r_{0}>0$ small enough such that
\[
\big\|\frac{\nabla \bfx^n_{1}}{\bfy^n_{1}}\big\|_{L_{x}^{2}}+\big\|\frac{\nabla \bfy^n_{1}}{\bfy^n_{1}}\big\|_{L_{x}^{2}}\ll\eps_0
\]
Here $\nabla$ is the spatial gradient and $r_{0}$ does not depend on~$n$. Since one can clearly also arrange
\[
\|\chi_{[|x-x_{0}|<r_{0}]}\phi^{n}_{0}\|_{L_{x}^{2}}\ll\eps_0,
\]
we have now achieved smallness of the energy of these data. Moreover, $r_0>0$ can be chosen {\em uniformly} in~$x_0\in\R^2$.
\end{lemma}
\begin{proof} We assume as we may (by rescaling) that
$\|\phi_\alpha^1\|_2+\|\phi_\alpha^2\|_2\le 1$ for $\alpha=0,1,2$.
We shall also suppress the time dependence and drop the superscript~$n$.
 In view of \eqref{eq:leib} it suffices to estimate the
 contributions of those terms in which the derivatives falls on the
 cutoff~$\chi$ in~\eqref{eq:x1y1}. Starting with the component
 ${\bf{y}}_{1}=\chi_{[|x-x_{0}|<r_{0}]}\frac{{\bf{y}}(x)}{{\bf{y}}(x_{0})}$,
note that Poincar\'e's inequality implies that, with
$B:=\{|x-x_{0}|<r_{0}\}$,
\begin{align*}
& r_0^{-2}\int_B
 \Big|\log\big[\frac{{\bf{y}}(x)}{{\bf{y}}_0}\big] \Big|^2\,
 dx_1 dx_2
 \les \int_B
 \Big| \frac{{\nabla\bf{y}}(x)}{{\bf{y}}(x)} \Big|^2\, dx_1dx_2 \\
 &\les \sum_{\alpha=1,2}\; \int_B
  |  \phi_\alpha(x)|^2\, dx_1dx_2  = \sum_{\alpha=1,2}\; \int_B
  |  \psi_\alpha(x)|^2\, dx_1dx_2 \ll\eps_0^2
\end{align*}
uniformly in~$n$ provided $r_{0}$ is small enough. Here we used the
relation~\eqref{eq:xy_fund} and that the gauge change is given by
multiplication by a unimodular factor. For the $\bfx_1$-component,
we make the preliminary observation that $\bfy\in A_p$ for any
$1<p<\infty$. Indeed, by Lemma~\ref{lem:embed}, for any $1\le
M<\infty$  we can find $C=C(M)$ so large that
\[
\|P_{\R\setminus[-C,C]} \log \bfy\|_{\bmo} = \|P_{\R\setminus[-C,C]}
\sum_{j=1,2} \Delta^{-1}\del_j \phi_j^2\|_{\bmo} \ll M^{-1}
\]
which implies that  $\bfy_2:= \exp\big(P_{\R\setminus[-C(p),C(p)]}
\log \bfy \big)\in A_p$ for any $1<p<\infty$ with a suitable $C(p)$.
Since Lemma~\ref{lem:embed} implies that $\|\bfy
\bfy_2^{-1}\|_\infty\le C$, the claim follows. We now use the
following weighted Poincar\'e inequality, see Theorem~1.5 in~\cite{FKS}: for any
$w\in A_2$, and ball $B$ of radius~$r>0$,
\[
\int_B |f(x)-(f)_B|^2\, w(x)dx \le C(w) r^2 \int_B |\nabla f(x)|^2\,
w(x) dx,\qquad (f)_B:= \slashint_B f(x)\, dx
\]
Consequently, with $w=\bfy^{-2}\in A_2$, and in view of our
definition of~$\bfx_0$,
\begin{align*}
  r_0^{-2} \int_B \Big|\frac{\bfx(x)-\bfx_0}{\bfy(x)}\Big|^2\, dx_1dx_2
 &\les  \int_B
\Big|\frac{\nabla\bfx(x)}{\bfy(x)}\Big|^{2}\, dx_1dx_2  \les
 \sum_{j=1,2}\int_B |\phi^1_j(x)|^2\, dx_1 dx_2
\end{align*}
By our  choice of cutoffs in~\eqref{eq:x1y1} we are done. To obtain the final statement of the proof,
simply note that we can always find $r_0>0$ such that
\[
 \sup_{x_0\in\R^2} \sum_{\alpha=0}^2\int_{D_{x_0,r_0}} |V_\alpha(x)|^2\, dx_1dx_2 \ll \eps_0^2
\]
Consequently, for all sufficiently large~$n$,
\[
 \sup_{x_0\in\R^2} \sum_{\alpha=0}^2\int_{D_{x_0,r_0}} |\psi^n_\alpha(x)|^2\, dx_1dx_2 \ll \eps_0^2
\]
and therefore also
\[
 \sup_{x_0\in\R^2} \sum_{\alpha=0}^2\int_{D_{x_0,r_0}} |\phi^n_\alpha(x)|^2\, dx_1dx_2 \ll \eps_0^2
\]
for all large~$n$,
which is all that is needed for the proof.
\end{proof}

We will also require an analogous result on small energy outside of a big ball.
Thus, let $R_0\gg1$ be large and define
\begin{equation}\label{eq:x2y2}
({\bfx}_{2}, {\bfy}_{2}):=
\Big(\chi_{[|x|>R_{0}]}\frac{{\bfx}-{\bfx}_{0}}{{\bfy}_{0}},\:e^{\chi_{[|x|>\frac{R_{0}}{2}]}\log[\frac{\bfy}{{\bfy}_{0}}(0,\cdot)]}
\Big)
\end{equation}
Here $\chi_{[|x|>R_{0}]}$ is a smooth cutoff to the
set~$\{|x|>R_0\}$ which equals one on
$|x|>2R_0$, say, and
\begin{align*}
\bfx_0&:= \slashint_{[R_0<|x|<2R_{0}]} \bfx(x)\, dx_1dx_2,\qquad
 \bfy_0 :=
\exp\big(\slashint_{[\frac{R_0}{2}<|x|<R_{0}]} \log \bfy(x)\, dx_1dx_2\big)
\end{align*}
In analogy to~\eqref{eq:x1y1} the construction here is such that $\bfy_2=\frac{\bfy}{\bfy_0}$
on the set~$\{\nabla \chi_{[|x|>R_{0}]}\ne0\}$.

\begin{lemma}\label{lem:DataLocalization'} Let $(\bfx_2^n,\bfy_2^n)$ be
as in~\eqref{eq:x2y2} applied to~$(\bfx^n,\bfy^n)$. Then
 given $\eps_0>0$, we can pick $R_{0}>0$ large enough such that
\[
\big\|\frac{\nabla \bfx^n_{2}}{\bfy^n_{2}}\big\|_{L_{x}^{2}}+\big\|\frac{\nabla \bfy^n_{2}}{\bfy^n_{2}}\big\|_{L_{x}^{2}}\ll\eps_0
\]
Here $\nabla$ is the spatial gradient and $R_{0}$ does not depend on~$n$. Since one can clearly also arrange
\[
\|\chi_{[|x|>R_{0}]}\phi^{n}_{0}\|_{L_{x}^{2}}\ll\eps_0,
\]
we have now achieved smallness of the energy of these data.
\end{lemma}
\begin{proof}
The argument is completely analogous to the one for Lemma~\ref{lem:DataLocalization}. The only difference is that
ones uses the following Poincar\'e inequalities on annuli instead of disks: for any $R_0>0$,
\begin{align*}
 \int_{R_0<|x|<2R_0} |f(x)-(f)_{R_0}|^2\, w(x)dx \le CR_0^2 \int_{R_0<|x|<2R_0} |\nabla f(x)|^2\, w(x)dx
\end{align*}
for any $A_2$-weight $w$ and a constant $C$ which only depends on
the~$A_2$ constant of~$w$. As usual $(f)_{R_0}$ denotes the average
of~$f$ over the annulus. For $w=1$ this is of course standard, and
for general~$w$ it follows from~\cite{FKS}.
\end{proof}

Next, we wish to establish control over the $\psi^n_\alpha$
in the $S$-norm on a nonempty time interval $(-T_0,T_0)$ uniformly in~$n$. The idea
is to apply Theorem~\ref{thm:smallenergy} to the finitely many small
energy maps given by Lemma~\ref{lem:DataLocalization} and then to
reconstruct and also bound the original sequence~$\psi^n_\alpha$ in
terms of these constituents. The latter of course relies on finite
propagation speed and involves smooth partitions of unity. In order to
handle partitions of unity, we need to derive estimates of the form
\[
 \|\chi\psi\|_S\le C(\chi)\|\psi\|_S
\]
for Schwartz functions~$\chi$ and some constant~$C(\chi)$. Due to issues
 having to do with the slow decay as well as limited regularity of the logarithmic potential
$\Delta^{-1}\del\phi$ which appears in the phase of the gauge change, we will need
to allow for a larger class of functions~$\chi$. The following lemma is tailored to such purposes.

\begin{lemma}
  \label{lem:chiS} Let $\chi\in C^\infty(\R^{2+1})$ satisfy the following properties\footnote{The logarithmic potential in~\eqref{eq:logpot}
decays like $|z|^{-1}$ (but in general no faster) which explains why
we need $p>2$ in the first condition. Since one in fact has
asymptotic equality with $|z|^{-1}$ up to a multiplicative constant,
it follows that the Fourier transform of this potential around zero
exhibits a $|\xi|^{-1}$-singularity, which explains the second
condition. Finally, we cannot control more than one time derivative
of~\eqref{eq:logpot}, and showing that one time derivative can be
controlled in terms of the energy alone is nontrivial and requires
the div-curl system for~$\phi$, see
Corollary~\ref{cor:smalltimeexistence}.}: for some constant $A$
\begin{itemize}
\item  $\max_{k=0,1} \max_{|\alpha|\le 100} \|\del_t^k \nabla_x^\alpha \chi\|_{L^q_t L^p_x} \le A$ for all $2<p\le \infty$, $1\le q\le\infty$
\item $ \|\la\tau\ra \max(|\xi|, |\xi|^{100})  \wh{\chi}(\tau,\xi)\|_{L^q_\tau L^\infty_\xi} \le A$ for all $2\le q\le \infty$
\end{itemize}
Then there exists an absolute  constant $C_0$   such that
$\|\chi\psi\|_S\le C_0 A\|\psi\|_S$ for any Schwartz
function~$\psi$. The $S$-norm here can be defined in terms of {\em
both} the {\em original} $S[k]$-spaces from Definition~\ref{def:Sk}
as well as the stronger~$\trip\cdot\trip$-norm, and it can
be either localized to some interval in time or be defined globally in time.
\end{lemma}
\begin{proof}
It suffices to consider  global in time estimates.  We begin with
the original $S[k]$-norm. We need to prove that
\begin{equation}
\Big(\sum_{k\in\Z} \|P_k(\chi\psi)\|_{S[k]}^2\Big)^{\frac12} \le
C(\chi) \|\psi\|_{S} \label{eq:Schi}
\end{equation}
 We begin with the energy
component of $S[k]$. If $k\le C$, then by Bernstein's inequality
\begin{equation}\label{eq:Pkchipsi}
\|P_k(\chi \psi)\|_{\ener} \les 2^{\frac{2k}{3}} \|\chi\|_{L^\infty_t L^3_x}
\|\psi\|_{\ener} \les 2^{\frac{2k}{3}} \|\chi\|_{L^\infty_t L^3_x}   \|\psi\|_{S} \les 2^{\frac{2k}{3}} A   \|\psi\|_{S}
\end{equation}
Here we used that
\[
 \|\psi\|_{\ener}\le \Big(\sum_{k\in\Z} \|P_k \psi\|_{\ener}^2\Big)^{\frac12} \le \Big(\sum_{k\in\Z} \|P_k \psi\|_{S[k]}^2\Big)^{\frac12}
\le \|\psi\|_{S}
\]
On the other hand, if $k\ge C$, then
\begin{align*}
  \|P_k(\chi\psi)\|_{\ener} &\le \|P_k(P_{\le k-10}\chi\, \tilde P_k
\psi)\|_{\ener} + \|P_k(P_{> k-10}\chi\, \psi)\|_{\ener} \\
&\les \|\chi\|_{\Linf} \|\tilde P_k\psi\|_{\ener} +
\|P_{>k-10} \chi\|_{L^\infty_{t,x}}  \|\psi\|_{\ener}\\
&\les (\|\chi\|_{\Linf}+ \|\nabla \chi\|_{\Linf}) (\|\tilde P_k\psi\|_{S[k+O(1)]} + 2^{-k} \|\psi\|_{S})\\
&\les A(\|\tilde P_k\psi\|_{S[k+O(1)]} + 2^{-k} \|\psi\|_{S})
\end{align*}
where we used the reverse Bernstein inequality
\[
\|P_{>k-10} \chi\|_{\Linf} \les 2^{-k} \| \nabla \chi\|_{\Linf}
\]
Square-summing now implies the desired bound.
To proceed we need to control $\|P_{k} \chi\|_{S[k]}$. If $k\le0$, then
\begin{align} \nn
\|P_{k}\chi\|_{S[k]} &\les \|P_{k} Q_{\le k} \chi\|_{\dot X_{k}^{0,\frac12,1}}
+ \|P_{k} Q_{> k} \chi\|_{\dot X_{k}^{-\frac12+\eps,1-\eps,2}} \\
&\les 2^{\frac{k}{2}} \|P_{k} Q_{\le k} \chi\|_{\Ltwotx}  + 2^{(-\frac12+\eps)k} \|P_{k} Q_{k<\cdot<0} \chi\|_{\Ltwotx} +
2^{(-\frac12+\eps)k} \| P_{k} Q_{>0} \,\del_t\chi\|_{\Ltwotx} \nn \\
&\les 2^{\frac{k}{2}}  \Big(  \int  \sup_{|\eta|\sim 2^k} ( |\eta|^2|\wh{\chi}(\tau,\eta)|^2) \,d\tau \Big)^{\frac12}
+ 2^{(-\frac12+\eps)k}  \Big( \int  \la\tau\ra^2 \sup_{|\eta|\sim 2^k} ( |\eta|^2|\wh{\chi}(\tau,\eta)|^2) \,d\tau
\Big)^{\frac12}\nn \\
&\les A\, 2^{(-\frac12+\eps)k}
\label{eq:chiSk1}\end{align}
whereas if $k\ge0$, then
\begin{align} \nn
\|P_{k}\chi\|_{S[k]} &\les \|P_{k} Q_{\le k} \chi\|_{\dot X_{k}^{0,\frac12,1}}
+ \|P_{k} Q_{> k} \chi\|_{\dot X_{k}^{-\frac12+\eps,1-\eps,2}} \\
&\les 2^{\frac{k}{2}} \|P_{k} Q_{\le k} \chi\|_{\Ltwotx}  +
2^{(-\frac12+\eps)k} \| P_{k} Q_{>k} \,\del_t\chi\|_{\Ltwotx} \nn \\
&\les  2^{-10 k}  \Big( \int  \la\tau\ra^2 \sup_{|\eta|\sim 2^k} ( |\eta|^{22}|\wh{\chi}(\tau,\eta)|^2) \,d\tau
\Big)^{\frac12} \les A\, 2^{-10k}
\label{eq:chiSk2}\end{align}
Next, if $k\le C$,
then
\[
\|P_k Q_{\le k} (\chi\psi)\|_{\dot X_k^{0,\frac12,\infty}} \les
2^{\frac{k}{2}} \|P_k(\chi\psi)\|_{\Ltwotx} \les 2^{\frac{k}{2}}
\|\chi\|_{L^2_t L^\infty_x} \|\psi\|_{S} \les 2^{\frac{k}{2}}
A\, \|\psi\|_{S}
\]
whereas, if $k>C$, then by Lemma~\ref{lem:Sk_prod2}
\begin{equation}\label{eq:G1}
\|P_k Q_{\le k} (\chi\psi)\|_{\dot X_k^{0,\frac12,\infty}} \les
\sum_{k_1,k_2\in\Z} 2^{k_1\wedge k_2} 2^{\frac{k-k_1\vee k_2}{4}}\;
\| P_{k_1}\chi\|_{S[k_1]} \| P_{k_2} \psi\|_{S[k_2]}
\end{equation}
where the sum is only over those $k_1,k_2$ which are admissible by the low-high, high-low, and high-high
trichotomy.
Distinguishing these cases and
using the estimate~\eqref{eq:chiSk1} , \eqref{eq:chiSk2},
one obtains
\begin{align*}
\|P_k(\chi\psi)\|_{\dot X_k^{0,\frac12,\infty}} \les A\big(
 \|\tilde P_{k} \psi\|_{S[k]} + 2^{-k}
\|\psi\|_{S}\big)
\end{align*}
which is sufficient. Next, if $k\le C$, then
\[
\|P_kQ_{k\le \cdot\le C} (\chi\psi)\|_{\dot
X_k^{-\frac12+\eps,1-\eps,2}} \les 2^{(\frac16+\eps)k} \|\chi
\psi\|_{L^2_t L^{\frac65}_x} \les 2^{(\frac16+\eps)k}
\|\chi\|_{L^2_t L^3_x} \|\psi\|_{S} 2^{(\frac16+\eps)k} A
\|\psi\|_{S}
\]
and for $j\ge C$, using that $P_kQ_j$ is disposable,
\begin{align*}
\|P_kQ_{j} (\chi\psi)\|_{\dot X_k^{-\frac12+\eps,1-\eps,2}} &\les
\sum_{j\ge C} 2^{j(1-\eps)} 2^{(\frac16+\eps)k}\big[ \|P_kQ_j(Q_{\le
j-10} \chi\, \psi)\|_{L^2_t L^{\frac65}_x} + \|Q_{> j-10} \chi\,
\psi \|_{L^2_t L^{\frac65}_x}\big]\\
&\les \sum_{j\ge C} 2^{j(1-\eps)} 2^{(\frac16+\eps)k}\big[
\|P_kQ_j(Q_{\le j-10} \chi\, \psi)\|_{L^2_t L^{\frac65}_x} + 2^{-j}\| \del_t \chi\|_{L^2_t L^3_x}
\|\psi\|_{\ener}\big]\\
&\les A\, 2^{(\frac16+\eps)k}
\|\psi\|_{\ener}+\sum_{j\ge C} 2^{j(1-\eps)} 2^{(\frac16+\eps)k}
\|P_kQ_j(Q_{\le j-10} \chi\, \psi)\|_{L^2_t L^{\frac65}_x}
\end{align*}
The first  term on the right-hand side here can be square summed
over~$k\le C$, whereas the second is bounded as follows:
\begin{align}
&\sum_{j\ge C} 2^{j(1-\eps)} 2^{(\frac16+\eps)k}  \|P_kQ_j(Q_{\le
j-10} \chi\, \psi)\|_{L^2_t L^{\frac65}_x} \nn \\
&\les \sum_{j\ge C} 2^{j(1-\eps)} 2^{(\frac16+\eps)k}\Big(
\sum_{\ell\ge C}  \|P_kQ_j(P_\ell Q_{\le j-10} \chi\, \tilde P_\ell
\tilde Q_j\psi)\|_{L^2_t L^{\frac65}_x} + \|P_kQ_j(P_{\le C} Q_{\le j-10}
\chi\,  P_{\le C} \tilde Q_j\psi)\|_{L^2_t L^{\frac65}_x} \Big) \nn \\
&\les \sum_{j\ge C} 2^{j(1-\eps)} 2^{(\frac16+\eps)k} \Big(
\sum_{\ell\ge C} \|P_\ell Q_{\le j-10} \chi\|_{L^\infty_t L^3_x} \|\tilde P_\ell \tilde Q_j\psi\|_{L^2_t
L^2_x} + \| P_{\le C} Q_{\le j-10}
\chi\|_{L^\infty_t L^3_x} \| P_{\le C} \tilde Q_j\psi \|_{L^2_t L^2_x} \Big) \label{eq:wub1}
\end{align}
where we used that $P_kQ_j$ is disposable for $k,j$ in the specified range.
Bernstein's inequality and Lemma~\ref{lem:QLp} imply
\begin{align*}
\eqref{eq:wub1} & \les
\sum_{j\ge C} 2^{j(1-\eps)} 2^{(\frac16+\eps)k}\Big(
\sum_{\ell\ge C} 2^{\frac{\ell}{3}} \|P_\ell  \chi\|_{L^\infty_t L^2_x} \|\tilde P_\ell \tilde Q_j\psi\|_{L^2_t
L^2_x} + \| P_{\le C} Q_{\le j-10}
\chi\|_{L^\infty_t L^2_x} \|P_{\le C} \tilde Q_j\psi \|_{L^2_t L^2_x} \Big) \\
&\les A\, 2^{(\frac16+\eps)k}\Big( \sum_{\ell\ge C} 2^{-\ell}
\|\tilde P_\ell \psi\|_{\dot X_\ell^{-\frac12+\eps,1-\eps,2}} +
\sum_{\ell\le C} 2^{(\frac12-\eps)\ell} \| P_{\ell} \psi \|_{\dot
X_\ell^{-\frac12+\eps,1-\eps,2}} \Big) \\
&\les A \,
2^{(\frac16+\eps)k} \|\psi\|_{S}
\end{align*}
To pass to the final estimate here one uses Cauchy-Schwarz. Since
this bound can be square-summed over~$k\le C$ it is admissible. Now suppose
$k\ge C$ and estimate
\begin{align*}
\sum_{j\ge k} 2^{j(1-\eps)} 2^{(-\frac12+\eps)k}  \|P_kQ_j(P_{\le
k-10} Q_{\le j-10} \chi\, \psi)\|_{L^2_t L^2_x} &\les \sum_{j\ge k}
2^{j(1-\eps)} 2^{(-\frac12+\eps)k}  \|\chi\|_{\Linf} \|\tilde P_k
\tilde Q_j
\psi\|_{L^2_t L^2_x} \\
&\les  A \|\tilde P_k \psi\|_{S[k]}
\end{align*}
where we used that $P_{\le
k-10} Q_{\le j-10} $ is disposable provided $j\ge k$, and similarly,
\begin{align*}
\sum_{j\ge k} 2^{j(1-\eps)} 2^{(-\frac12+\eps)k}  \|P_kQ_j( Q_{>
j-10}\chi\, \psi)\|_{L^2_t L^2_x} &\les \sum_{j\ge k} 2^{j(1-\eps)}
2^{(-\frac12+\eps)k}  2^{-j} \|\del_t \chi\|_{L^2_t L^\infty_x} \|
\psi\|_{L^\infty_t L^2_x} \\
&\les  A\, 2^{-\frac{k}{2}}  \| \psi\|_{S}
\end{align*}
which can be square-summed over $k\ge C$, and finally,
\begin{align*}
&\sum_{j\ge k} 2^{j(1-\eps)} 2^{(-\frac12+\eps)k}  \|P_kQ_j( P_{>
k-10}Q_{\le j-10}\chi\: \psi)\|_{L^2_t L^2_x} \\
&\les  \sum_{\ell \le k-10} \sum_{j\ge k} 2^{j(1-\eps)} 2^{(-\frac12+\eps)k}  \|P_kQ_j( \tilde P_k Q_{\le j-10}\chi\: \tilde P_\ell \tilde Q_j\psi)\|_{L^2_t L^2_x} \\
&\quad +  \sum_{\ell > k-10}\sum_{j\ge k} 2^{j(1-\eps)} 2^{(\frac12+\eps)k}  \|P_kQ_j( P_\ell Q_{\le j-10}\chi\: \tilde P_\ell\tilde Q_j\psi)\|_{L^2_t L^1_x}\\
& \les 2^{-k} (\|\nabla_x \chi\|_{\Linf} \|\psi\|_{S}  + \| \nabla^2 \chi\|_{\Linf}  \|P_\ell \psi\|_{S[\ell]}) \\
&\le A\, 2^{-k} \|
\psi\|_{S}
\end{align*}
In conclusion, for all $k\ge C$,
\[ \|P_kQ_{\ge k}
(\chi\psi)\|_{\dot X_k^{-\frac12+\eps,1-\eps,2}} \les A\,(\|\tilde
P_k \psi\|_{S[k]} + 2^{-k} \|\psi\|_S)
\]
which is admissible.  Next, we need to deal with the Strichartz
norms in~\eqref{eq:Sk2}. If $k\le C$, then the free-wave estimate of
Lemma~\ref{lem:Strich2} implies
\begin{equation}\label{eq:Strichagain}
\sup_{j\in\Z} \sup_{\ell\le0}\; 2^{-(\frac12-\eps)\ell}
2^{-\frac{3k}{4}} \Big(\sum_{c\in \calD_{k,\ell}} \|Q_{<j} P_c (\chi
\psi)\|_{L^4_t L^\infty_x}^2\Big)^{\frac12} \le C
\|P_k(\chi\psi)\|_{\dot X_k^{0,\frac12,1}}
\end{equation}
By the preceding
\begin{align*}
\|P_k(\chi\psi)\|_{\dot X_k^{0,\frac12,1}} &\les 2^{\frac{k}{2}}
\|P_k(\chi\psi)\|_{\Ltwotx}+ \|P_kQ_{\ge k} (\chi\psi)\|_{\dot
X_k^{-\frac12+\eps,1-\eps,2}} \les A\,  2^{\frac{k}{6}} \|\psi\|_S
\end{align*}
Since this bound is square summable over $k\le C$ we are done with
the small frequencies. For $k\ge C$ one needs to argue differently.
The issue here is that we control $\|P_k Q_{\le k} (\chi
\psi)\|_{\dot X_k^{0,\frac12,\infty}}$ but not the stronger $\|P_k
Q_{\le k} (\chi \psi)\|_{\dot X_k^{0,\frac12,1}}$ which is needed
for a reduction to free waves. On the other hand, by
Lemma~\ref{lem:Sk_prod2},
\begin{equation}\label{eq:G2}
\sup_{j\le 0}\|P_k Q_{<j} (\chi\psi)\|_{\dot X_k^{0,\frac12,1}} \les
\sum_{k_1,k_2\in\Z} 2^{\frac{3}{4}k_1\wedge k_2} 2^{\frac{k-k_1\vee
k_2}{4}}\; \| P_{k_1}\chi\|_{S[k_1]} \| P_{k_2} \psi\|_{S[k_2]}
\end{equation}
where the sum is only over admissible $k_1,k_2$, i.e., those which respect the
usual trichotomy.
By the same considerations as before one can see that the right-hand
side here is an admissible bound, i.e., square-summable over~$k$.
This observation shows that~\eqref{eq:Strichagain} can be controlled
provided the supremum in~$j$ is taken over~$j\le 0$. On the other
hand, if $j\ge0$, then we use the commutation relation
\begin{equation}\label{eq:commrel}\begin{aligned}
Q_{<j} (\chi \psi) &= \chi Q_{<j} \psi + 2^{-j} \sum_{i=1,2} \int_{\R^4}
\del_i \chi(\cdot-z_1) \psi(\cdot-z_2)\, \nu_{ij}(dz_1,dz_2) \\
&=: \chi Q_{<j} \psi + 2^{-j} L_j(\nabla \chi,\psi)
\end{aligned}
\end{equation}
where $\nu_{ij}$ is a measure with mass controlled uniformly
in~$j\ge0$. By these considerations, we may ignore the supremum
over~$j$ in~\eqref{eq:Strichagain} altogether. We begin with the
low-high case. With $k\ge C$,
\begin{align}
 & \Big(\sum_{c\in\calD_{k,\ell}} \|P_c(P_{\le k-10}\chi\, \tilde P_k
\psi)\|_{L^4_t L^\infty_x}^2\Big)^{\frac12}\nn \\
&\les \sum_{ j\le
k-10} \Big(\sum_{c\in\calD_{k,\ell}}\big(\sum_{c_1\sim_c c_2}
\|P_c(P_{c_1}P_j\,\chi\: \tilde P_k P_{c_2} \psi)\|_{L^4_t
L^\infty_x}\big)^2\Big)^{\frac12} \nn \\
&\les \sum_{ j\le k-10} \Big(\sum_{c\in\calD_{k,\ell}}
\sum_{c_1\sim_c c_2} \max(1,2^{2(j-k-\ell)})
\|P_c(P_{c_1}P_j\,\chi\: \tilde P_k P_{c_2} \psi)\|_{L^4_t
L^\infty_x}^2\Big)^{\frac12} \label{eq:simc1c2}
\end{align}
Here $c_2$ runs over $\calD_{k,\ell}$, and $c_1$ runs over
$\calD_{j,\ell'}$ where $\ell':=\max(\ell+k-j,-10)$ with the added
property that $\dist(c_1+c_2,c)\les \diam(c)$ (which we denote by $c_1\sim_c c_2$).  The
factor $\max(1,2^{2(j-k-\ell)})$ equals the largest number of disks $c_1$ possible, and it arises due to Cauchy-Schwarz.
By means of
Bernstein's inequality one can bound the final norm
in~\eqref{eq:simc1c2} as follows:
\begin{align*}
 \|P_c(P_{c_1}P_j\,\chi\: \tilde P_k P_{c_2} \psi)\|_{L^4_t
L^\infty_x} &\les \|P_{c_1}P_j\,\chi\: \tilde P_k P_{c_2}
\psi\|_{L^4_t L^\infty_x}\les \|P_{c_1}P_j\,\chi\|_{\Linf} \| P_{c_2} \psi\|_{L^4_t
L^\infty_x} \\
&\les  \min(2^{j}, 2^{k+\ell}) \|\wh{P_{c_1}\chi}\|_{L^1_\tau L^2_\xi} \|
P_{c_2} \psi\|_{L^4_t L^\infty_x}
\end{align*}
Since there are only finitely many choices of~$c$
in~\eqref{eq:simc1c2} one concludes that
\begin{align*}
&\Big(\sum_{c\in\calD_{k,\ell}} \|P_c(P_{\le k-10}\chi\, \tilde P_k
\psi)\|_{L^4_t L^\infty_x}^2\Big)^{\frac12} \\
& \les \sum_{ j\le
k-10} \Big(  \sum_{c_1, c_2}
\max(1,2^{2(j-k-\ell)}) \|P_{c_1}P_j\,\chi\: \tilde P_k P_{c_2}
\psi\|_{L^4_t L^\infty_x}^2\Big)^{\frac12}\\
&\les \sum_{ j\le k-10} \Big( \sum_{c_1, c_2} 2^{2j}
\|\la \tau\ra \wh{P_{c_1}\chi}\|_{L^2_{\tau,\xi}}^2 \| P_{c_2} \psi\|^2_{L^4_t L^\infty_x}
\Big)^{\frac12} \les \sum_{ j\le k-10} 2^{j}
\|\la \tau\ra \wh{P_j\chi}\|_{L^2_{\tau,\xi}} \Big( \sum_{c_2}  \| P_{c_2} \psi\|^2_{L^4_t L^\infty_x}
\Big)^{\frac12} \\
&\les \sum_{j\le k-10} 2^j \| \la\tau\ra \sup_{|\eta|\sim 2^j} |\eta|\, |\wh{\chi}(\eta)|  \|_{L^2_\tau}   \Big( \sum_{c_2} \| P_{c_2} \psi\|^2_{L^4_t
L^\infty_x} \Big)^{\frac12} \\
& \les A\,  2^{(\frac12-\eps)\ell}
2^{\frac{3k}{4}} \|\psi\|_{S[k]}
\end{align*}
Next, the high-low case is estimated similarly. More precisely, with
the roles of $c_1$ and $c_2$ interchanged, one has
\begin{align*}
 & \Big(\sum_{c\in\calD_{k,\ell}} \|P_c(\tilde P_{k}\chi\,
P_{<k-10}
\psi)\|_{L^4_t L^\infty_x}^2\Big)^{\frac12}\nn \\
&\les \sum_{ j\le k-10} \Big( \sum_{c_1, c_2}
\max(1,2^{2(j-k-\ell)}) \|P_{c_1}\,\chi\: P_j P_{c_2} \psi\|_{L^4_t
L^\infty_x}^2\Big)^{\frac12} \\
&\les  \sum_{ j\le k-10} \Big( \sum_{c_1, c_2}
\max(2^{2(k+\ell)},2^{2j})
 \|\la \tau\ra \wh{P_{c_1}\,\chi}\|_{L^2_{\tau,\xi}}^2 \|P_j P_{c_2} \psi\|_{L^4_t
L^\infty_x}^2\Big)^{\frac12}\\
&\les \|\la \tau\ra \la \xi\ra\wh{\tilde P_{k}\,\chi}\|_{L^2_{\tau,\xi}} \sum_{ j\le k-10} \Big( \sum_{c_2\in \calD_{j,\ell'}}
\max(2^{2(k+\ell)},2^{2j})
 2^{-2k}  \|P_{c_2} \psi\|_{L^4_t
L^\infty_x}^2\Big)^{\frac12}\\
&\les A\sum_{ j\le k-10}  \max(2^\ell,2^{j-k})
2^{\frac{3j}{4}}
\min(2^{(\frac12-\eps)(\ell+k-j)},1) \|P_j\psi\|_{S[j]} \\
&\les A\, 2^{\frac{3k}{4}} 2^{(\frac12-\eps)\ell} \sum_{j\le k}
2^{\frac34(j-k)}  \|P_j \psi\|_{S[j]}
\end{align*}
This bound can be square-summed in~$k$.
Finally, in the high-high case one has with $c_1,c_2\in
\calD_{m,\ell}$,
\begin{align*}
 & \sum_{m\ge k-10} \Big(\sum_{c\in\calD_{k,\ell}} \|P_c( P_{m}\chi\,
\tilde P_{m}
\psi)\|_{L^4_t L^\infty_x}^2\Big)^{\frac12}\nn \\
&\les \sum_{m\ge k-10}  \Big( \sum_{c_1, c_2\in\calD_{m,\ell+k-m}}
2^{2(m-k-\ell)} \|P_m P_{c_1}\,\chi\: \tilde P_m P_{c_2}
\psi\|_{L^4_t
L^\infty_x}^2\Big)^{\frac12} \\
&\les \sum_{ m\ge k-10} \Big( \sum_{c_1, c_2\in \calD_{m,\ell+k-m}}
2^{2m}
 \|P_{c_1}\,\chi\|_{\Ltwotx}^2 \|P_m P_{c_2} \psi\|_{L^4_t
L^\infty_x}^2\Big)^{\frac12}\\
&\les A \sum_{ m\ge k-10} 2^{-m} 2^{\frac{3m}{4}}
2^{(\frac12-\eps)(\ell+k-m)} \|P_m\psi\|_{S[m]}
\end{align*}
which is again sufficient due to $k\ge C$.

Lastly, we bound the square function in~\eqref{eq:squarefunc}
for~$\chi\psi$.  If $k\le C$, then by Lemma~\ref{lem:incl_free},
\begin{align*}
&\sup_{\pm}\sup_{\ell\le-100}\;\sup_{\ell\le m\le 0}
  \Big(  \sum_{\kappa\in\caps_\ell} \sum_{R\in\calR_{k,\pm \kappa,m}}
  \|P_{R}
  Q_{\le k+2\ell}^\pm\:(\chi
  \psi)\|_{S[k,\kappa]}^2 \Big)^{\frac12} \\
&\les \sup_{\pm}\sup_{\ell\le-100}\;\sup_{\ell\le m\le 0}
  \Big(  \sum_{\kappa\in\caps_\ell} \sum_{R\in\calR_{k,\pm \kappa,m}}
  \|P_{R}
  Q_{\le k+2\ell}^\pm\:(\chi
  \psi)\|_{\dot X^{0,\frac12,1}_k}^2 \Big)^{\frac12} \\
&\les \sup_{\ell\le 0}\|P_k Q_{k+2\ell} (\chi\psi)\|_{\dot
X^{0,\frac12,1}_k}  \les 2^{\frac{k}{2}} \|\chi\psi\|_{\Ltwotx} \les A\,2^{\frac{k}{2}} \|\psi\|_{S}
\end{align*}
It therefore suffices to consider $k\ge C$. We make a number of
reductions by means of Lemma~\ref{lem:incl_free}: first,
\begin{align*}
 & \sup_{\ell\le0} \|P_k Q_{\le k+2\ell} (P_{<k-10}\chi\, Q_{\ge k+2\ell}\psi)\|_{\dot
X_k^{0,\frac12,1}}  \les \sup_{\ell\le0} 2^{\frac{k}{2}+\ell}
\|\chi\|_{\Linf} \|\tilde P_k Q_{\ge k+2\ell} \psi\|_{\Ltwotx} \les A\|\tilde P_k \psi\|_{S[k]}
\end{align*}
Second,
\begin{align*}
 & \sup_{\ell\le0} \|P_k Q_{\le k+2\ell} (P_{\ge k-10}\chi\, \psi)\|_{\dot
X_k^{0,\frac12,1}}  \les \sup_{\ell\le0} 2^{\frac{k}{2}+\ell}
\|P_{\ge k-10}\chi\|_{L^2_t L^\infty_x} \| \psi\|_{\ener}  \les A\,2^{-k}\|\psi\|_{S}
\end{align*}
and third,
\begin{align*}
 & \sup_{\ell\le0} \|P_k Q_{\le k+2\ell} (P_{<k-10}Q_{\ge k+2\ell}\chi\, \psi)\|_{\dot
X_k^{0,\frac12,1}}  \les \sup_{\ell\le0} 2^{\frac{k}{2}+\ell} \|
Q_{\ge k+2\ell} \chi\|_{L^2_t L^\infty_x}
 \|\tilde P_k \psi\|_{\ener} \les A\|\tilde P_k \psi\|_{S[k]}
\end{align*}
which is again admissible.  To final estimate here requires some justification:
\begin{align*}
 \sup_{m\in\Z} 2^m \|Q_{\ge 2m} \chi\|_{L^2_t L^\infty_x} &\les \sum_{\ell} \sup_{m\in\Z} 2^m \|\calF[Q_{\ge 2m} P_\ell \chi]\|_{L^2_\tau L^1_\xi} \\
&\les \sum_{\ell} 2^{-\ell}  \sup_{m\in\Z} 2^m \|\chi_{[||\xi|-|\tau||\ge 2^{2m}]} \chi_{[|\xi|\sim 2^\ell]}  \sup_{|\eta|\sim 2^\ell} |\eta||\wh{\chi}(\tau,\eta)|\, \|_{L^2_\tau L^1_\xi} \\
 & + \sum_{\ell} 2^{-100\ell}  \sup_{m\in\Z} 2^m \|\chi_{[||\xi|-|\tau||\ge 2^{2m}]} \chi_{[|\xi|\sim 2^\ell]} 2^{-\ell}  \sup_{|\eta|\sim 2^\ell} |\eta|^{100} |\wh{\chi}(\tau,\eta)|\, \|_{L^2_\tau L^1_\xi} \\
&\les \sum_{\ell} \min(2^\ell,2^{-\ell}) \| \la\tau\ra^{\frac12} |\eta|\vee |\eta|^{100} |\wh{\chi}(\tau,\eta)|\, \|_{L^2_\tau L^\infty_\eta} \les A
\end{align*}
Consequently, it suffices to bound
\begin{align}
  & \sup_{\pm}\sup_{\ell\le-100}\;\sup_{\ell\le m\le 0}
  \Big(  \sum_{\kappa\in\caps_\ell} \sum_{R\in\calR_{k,\pm \kappa,m}}
  \|P_{R}
  Q_{\le k+2\ell}^\pm\:(P_{<k-10}Q_{\le k+2\ell}\chi\: \tilde P_k
  Q_{\le k+2\ell}\psi)\|_{S[k,\kappa]}^2 \Big)^{\frac12} \nn \\
& \les \sup_{\pm}\sup_{\ell\le-100}\;\sup_{\ell\le m\le 0}
\sum_{j\le k-10}
  \Big(  \sum_{\kappa\in\caps_\ell} \sum_{R\in\calR_{k,\pm \kappa,m}}
  \|P_{R}
  Q_{\le k+2\ell}^\pm\:(P_{j}Q_{\le k+2\ell}\chi\: \tilde P_k
  Q_{\le k+2\ell}\psi)\|_{S[k,\kappa]}^2 \Big)^{\frac12}
\label{eq:ksumkappa}
\end{align}
We now distinguish two cases: 1) $j\le k+\ell$ and  2) $k+\ell\le
j$.

\noindent In the first case, the projection $P_R$ essentially passes
through $P_j\chi$. More precisely, by elementary geometry one has
\begin{align*} \|P_{R}
  Q_{\le k+2\ell}^\pm\:(P_{j}Q_{\le k+2\ell}\chi\: \tilde P_k
  Q_{\le k+2\ell}\psi)\|_{S[k,\kappa]} &\les \|P_{j}Q_{\le
k+2\ell}\chi\|_{\Linf} \| P_{\tilde R}
  Q_{\le k+2\ell}\psi\|_{S[k,\kappa]} \\
&\les \min(2^{j},2^{-j}) \|\la\tau\ra |\xi| \wh{\chi}(\tau,\xi)\|_{L^2_\tau L^\infty_\xi} \| P_{\tilde R}
  Q_{\le k+2\ell}\psi\|_{S[k,\kappa]}\\
&\les A\, \min(2^{j},2^{-j}) \| P_{\tilde R}
  Q_{\le k+2\ell}\psi\|_{S[k,\kappa]}
\end{align*}
where $\tilde R\in \calR_{k+O(1),\pm \kappa,m}$ has the property
that $\dist(R,\tilde R)\les 2^j$. This bound
reduces~\eqref{eq:ksumkappa} to the square function of~$\psi$ alone
and is therefore sufficient.
 In the second case, one estimates~\eqref{eq:ksumkappa} directly from Lemma~\ref{lem:incl_free}.
In fact, the second case contributes at most
\begin{align*}
 & \sum_{k+\ell \le j\le k-10}    \|P_j Q_{\le k+2\ell} \chi\: \tilde P_k Q_{\le k+2\ell} \psi\|_{\dot X_j^{0,\frac12,1}}  \\
 & \les \sum_{k+\ell \le j\le k-10}  2^{\frac{k}{2}+\ell} \|P_j Q_{\le k+2\ell} \chi\: \tilde P_k Q_{\le k+2\ell} \psi\|_{\Ltwotx}  \\
 &\les \sum_{k+\ell \le j\le k-10}  2^{\frac{k}{2}+\ell} \|P_j Q_{\le k+2\ell} \chi\|_{L^2_t L^\infty_x} \| \tilde P_k Q_{\le k+2\ell} \psi\|_{\ener}   \\
&\les \sum_{k+\ell \le j\le k-10}  2^{\frac{j}{2}} 2^j\wedge 2^{-10j} \||\xi|\vee |\xi|^{100} \wh{\chi}(\tau,\xi)\|_{L^2_\tau L^\infty_\xi}
 \| \tilde P_k   \psi\|_{\ener}   \\
&\les A\,\| \tilde P_k   \psi\|_{S[k]}
\end{align*}
To pass to the fourth line one uses that $\ell\le 0$.
This concludes the proof of~\eqref{eq:Schi}.

It remains to prove the analogous bound for the stronger
norm~$\trip\cdot \trip_{S[k]}$, i.e.,
\begin{equation}
\trip \chi\psi \trip_{S}:= \Big(\sum_{k\in\Z} \trip P_k(\chi\psi)\trip_{S[k]}^2\Big)^{\frac12}
\les A\Big(\sum_{\ell\in\Z} \trip P_\ell
\psi\trip_{S[\ell]}^2\Big)^{\frac12} \label{eq:Schi2}
\end{equation}
Written out, the left-hand side here means
\begin{equation}
 \Big(\sum_{k\in\Z} \sup_{t_0}\| P_k(\chi\psi)(t_0)\|_{L^2_x}^2 +    \|P_k\Box(\chi\psi)\|_{N[k]}^2\Big)^{\frac12}
\label{eq:dreiterm}
\end{equation}
Clearly,
\[
 \sum_{k\in\Z} \| P_k(\chi\psi)(t_0)\|_{L^2_x}^2  \les \| (\chi\psi)(t_0)\|_{L^2_x}^2 \les \|\chi(t_0)\|_{L^\infty_x}^2 \|\psi(t_0)\|_{L^2_x}^2
\les A\trip \psi\trip_{S}
\]
which is admissible.
It remains to bound the second term in~\eqref{eq:dreiterm}.  First,
\begin{align*}
 \sum_{k\le C} \|P_k\Box(\chi\psi)\|_{N[k]}^2 &\les  \sum_{k\le C} \|P_kQ_{\le k}\Box(\chi\psi)\|_{\dot X_k^{-1,-\frac12,1}}^2 + \|P_kQ_{>k}\Box(\chi\psi)\|_{\dot X_k^{-\frac12+\eps,-1-\eps,2}}^2 \\
&\les  \sum_{k\le C} 2^{k} \|P_k(\chi\psi)\|_{L^2_t L^2_x}^2 \les \|\chi\|_{L^2_t L^\infty_x}^2 \|\psi\|_{S}^2 \les A^2 \trip \psi\trip_S^2
\end{align*}
We may therefore assume that henceforth $k\ge C$.
Second, by Lemma~\ref{lem:Sk_prod2}, and using similar arguments as earlier in the proof,
\begin{align}
 \sum_{k> C} \|P_kQ_{\le C}\Box(\chi\psi)\|_{N[k]}^2 &\les  \sum_{k> C} \|P_kQ_{\le C}\Box(\chi\psi)\|_{\dot X_k^{-1,-\frac12,1}}^2
\les  \sum_{k> C} \Big(\sum_{j\le C}\|P_kQ_{j}(\chi\psi)\|_{\dot X_k^{0,\frac12,\infty}}\Big)^2 \nn \\
&\les \sum_{k> C} \Big(\sum_{k_1,k_2} 2^{\frac34 k_1\wedge k_2} 2^{\frac{k-k_1\vee k_2}{4}} \|P_{k_1} \chi\|_{S[k_1]} \|P_{k_2}\psi\|_{S[k_2]}   \Big)^2 \nn  \\
&\le A \|\psi\|_{S}\les A\trip \psi\trip_S \nn
\end{align}
The sum over $k,k_1,k_2$ here respects the usual trichotomy.  As a final reduction, we need to limit the modulation of~$\chi$ (since we
cannot control $\Box\chi$). In fact,
\begin{align*}
 \| P_k  \Box( Q_{>\frac{3k}{4}} \chi\: \psi) \|_{N[k]} &\les \| P_k Q_{\le k}  \Box( Q_{>\frac{3k}{4}} \chi\: \psi) \|_{\dot X_k^{-1,-\frac12,1}} +
\| P_k Q_{> k}  \Box( Q_{>\frac{3k}{4}} \chi\: \psi) \|_{\dot X_k^{-\frac12+\eps,-1-\eps,2} } \\
&\les 2^{\frac{k}{2}} \|P_k ( Q_{>\frac{3k}{4}} \chi\: \psi) \|_{\Ltwotx} \\
&\les  2^{\frac{k}{2}} \|P_k ( P_{>k-10} Q_{>\frac{3k}{4}} \chi\: \psi) \|_{\Ltwotx} + 2^{\frac{k}{2}} \|P_k ( P_{\le k-10} Q_{>\frac{3k}{4}} \chi\: \tilde P_k\psi) \|_{\Ltwotx}  \\
&\les 2^{\frac{k}{2}} \| \chi_{[|\xi|\gtrsim 2^k]} \wht{\chi}(\tau,\xi)\|_{L^2_\tau L^1_\xi} \|\psi\|_{S} + 2^{\frac{k}{2}}
 \|  \chi_{[||\tau|-|\xi||\gtrsim 2^{\frac{3k}{4}}]}  \wht{\chi}(\tau,\xi) \|_{L^2_\tau L^1_\xi} \|\tilde P_k\psi \|_{S[k]}\\
&\les 2^{-k} \| |\xi|^{100} \wht{\chi}(\tau,\xi)\|_{L^2_\tau L^\infty_\xi} \|\psi\|_S + 2^{-\frac{k}{4}}
 \|  \la \tau\ra |\xi|\vee |\xi|^{100}  \wht{\chi}(\tau,\xi) \|_{L^2_\tau L^\infty_\xi} \|\tilde P_k\psi \|_{S[k]}\\
&\les A\, 2^{-\frac{k}{4}} \|\psi\|_S
\end{align*}
which is admissible.
We now estimate each of the three term on the right-hand side of
\begin{equation}\label{eq:Nkdrei} \begin{aligned}
 \| P_kQ_{> C}\Box(Q_{\le \frac{3k}{4}}\chi\:\psi)\|_{N[k]} &\les \| P_kQ_{> C}(\Box Q_{\le \frac{3k}{4}}\chi\: \psi)\|_{N[k]}+ \| P_kQ_{> C} (\del_\alpha Q_{\le \frac{3k}{4}}\chi\:\del^\alpha\psi)\|_{N[k]}\\& +
\| P_kQ_{> C}(Q_{\le\frac{3k}{4}}\chi\:\Box\psi)\|_{N[k]}
\end{aligned}
\end{equation}
First,
\begin{align*}
 \| P_kQ_{> C}(\Box Q_{\le \frac{3k}{4}}\chi\: \psi)\|_{N[k]} &\les 2^{-k} \|\Box Q_{\le \frac{3k}{4}}\chi\:
\psi\|_{L^1_t L^2_x} \les 2^{-k} \|\Box Q_{\le \frac{3k}{4}} \chi\|_{L^1_t L^\infty_x} \|\psi\|_{S} \\
&\les 2^{-\frac{k}{4}} \max_{k=0,1}\max_{|\alpha|\le 1} \|\del_t^k \nabla^\alpha \chi\|_{L^1_t L^\infty_x} \|\psi\|_S\les A\, 2^{-\frac{k}{4}}\|\psi\|_S
\end{align*}
which is admissible. Second,  by estimate (29) in~\cite{T1} as well as~\eqref{eq:chiSk1} and~\eqref{eq:chiSk2},
and with $k,k_1,k_2$ respecting the usual trichotomy,
\begin{align*}
 \| P_kQ_{> C}(\del_\alpha Q_{\le\frac{3k}{4}}\chi\:\del^\alpha\psi)\|_{N[k]} &\les 2^{-k} \sum_{k_1,k_2} 2^{k_1+k_2}
\|P_{k_1}\chi\|_{S[k_1]} \|P_{k_2} \psi\|_{S[k_2]} \\
&\les A 2^{-k} \sum_{k_1,k_2} 2^{k_1+k_2} \min(2^{-(\frac12-\eps)k_1},2^{-10k_1}) \|P_{k_2} \psi\|_{S[k_2]}\\
&\les A( 2^{-k} \|\psi\|_S + \| P_k \psi\|_{S[k]})
\end{align*}
 which is again square-summable in~$k$. As for the third term in~\eqref{eq:Nkdrei} we are reduced to showing the bound
\begin{equation}
 \label{eq:Nkreduc}
\sum_{k\ge C} \|P_k Q_{>C} (Q_{\le\frac{3k}{4}}\chi\: F)\|_{N[k]}^2 \les A^2\sum_{\ell\in\Z} \|P_\ell F\|_{N[\ell]}^2
\end{equation}
This bound in turn  follows via Schur's lemma from the following claim:
\begin{equation}
 \label{eq:Nkclaim}
 \|P_k Q_{>C} (Q_{\le\frac{3k}{4}}\chi\: P_\ell F)\|_{N[k]} \les A 2^{-\frac14|k-\ell|} \|P_\ell F\|_{N[\ell]}
\end{equation}
If $j\le\ell$ and $j\le C$, then by \eqref{eq:Xsb_dom} one always has the bound
\begin{align*}
 \|P_k Q_{>C} (Q_{\le\frac{3k}{4}}\chi\: P_\ell Q_{\le j} F)\|_{N[k]} &\le 2^{-k} \|P_k(Q_{\le\frac{3k}{4}}\chi\:
P_\ell Q_{\le j} F)\|_{\enerN}  \les 2^{-k} \|Q_{\le\frac{3k}{4}}\chi\|_{L^2_t L^\infty_x} \|P_\ell Q_{\le j} F\|_{\Ltwotx} \\
&\les 2^{\ell-k} 2^{\frac{j}{2}} \|\wh{\chi}(\tau,\xi)\|_{L^2_\tau L^1_\xi}  \|P_\ell Q_{\le j} F\|_{N[\ell]}  \les    A\, 2^{\ell-k}  \|F\|_{N[\ell]}
\end{align*}
which agrees with \eqref{eq:Nkclaim} provided $\ell\le k+C$. On the other hand, if $\ell\ge k+C$ but still $j\le C$ the
same estimate holds with an additional high-high gain of $2^{-100\ell}$ coming from~$\chi$ which is of course more than sufficient for~\eqref{eq:Nkclaim}.
Finally, if $C\ge j\ge\ell$, then an additional Bernstein gain yields
\begin{align*}
 \|P_k Q_{>C} (Q_{\le\frac{3k}{4}}\chi\: P_\ell Q_{\le j} F)\|_{N[k]} &\le 2^{-k} \|P_k(Q_{\le\frac{3k}{4}}\chi\:
P_\ell Q_{\le j} F)\|_{\enerN}  \les 2^{-k} \|Q_{\le\frac{3k}{4}}\chi\|_{L^2_t L^4_x} \|P_\ell Q_{\le j} F\|_{L^2_t L^4_x} \\
&\les 2^{\frac{\ell}{2}-k}  \|\wh{\chi}(\tau,\xi)\|_{L^2_\tau L^{\frac43}_\xi}  \|P_\ell Q_{\le j} F\|_{\Ltwotx}  \les    A\, 2^{\ell-k}  \|F\|_{N[\ell]}
\end{align*}
as desired.
Therefore, the claim~\eqref{eq:Nkclaim} holds
provided~$F=Q_{\le C}F$.  Let us now verify~\eqref{eq:Nkclaim} for each of the four types of $N[\ell]$-atoms with the additional assumption that
$F\ne Q_{\le C}F$. If $F$ is an energy atom, then
\[
 \|P_k Q_{>C} (Q_{\le\frac{3k}{4}}\chi P_\ell F)\|_{N[k]}\les 2^{-k} \|Q_{\le\frac{3k}{4}}\chi P_\ell F\|_{\enerN}
\les 2^{\ell-k} \|Q_{\le\frac{3k}{4}}\chi\|_{\Linf} \|P_\ell F\|_{N[\ell]}
\]
which is sufficient if $\ell\le k+C$ and if $\ell>k+C$ then
\[
 \|P_k Q_{>C} (Q_{\le\frac{3k}{4}}\chi P_\ell F)\|_{N[k]}\les 2^{-k} \|Q_{\le\frac{3k}{4}}\chi P_\ell F\|_{\enerN}
\les 2^{-\ell-k} \|Q_{\le\frac{3k}{4}}\chi\|_{\Linf} \|P_\ell F\|_{N[\ell]}
\]
which is more than sufficient.  Here we used the estimate
\begin{equation} \label{eq:Qchi34}
 \|Q_{\le\frac{3k}{4}}\chi\|_{\Linf} \les \|\wh{\chi}(\tau,\xi)\|_{L^1_\tau L^1_\xi} \les \|\la \tau\ra |\xi|\vee|\xi|^{100}
\wh{\chi}(\tau,\xi)\|_{L^2_\tau L^\infty_\xi} \les A
\end{equation}
For the remaining atoms we first make the simplifying assumption that $\wh{\chi}(\tau,\xi)$ is supported on $|\tau|+|\xi|\les 1$.
Now suppose that $P_\ell Q_j F=F$ with $j> C$.
If $\|F\|_{\Ltwotx}\le 2^\ell 2^{\frac{j}{2}}$ and $j\le \ell$, then $\chi$ essentially does not change the Fourier support  of~$F$.
Thus, $\ell=k+O(1)$ and
\begin{equation}\label{eq:Nkwieder}
 \|P_k Q_{>C} (\chi F)\|_{N[k]} \les 2^{-\ell} 2^{-\frac{j}{2}} \|\chi F\|_{\Ltwotx} \les A \|F\|_{\dot X_\ell^{-1,-\frac12,1}}
\end{equation}
as desired. On the other hand, if $j>\ell$ and $\|F\|_{\Ltwotx}\le 2^{\ell(\frac12-\eps)}  2^{j(1+\eps)}$, then
we need to distinguish the case $\ell\le C$ from $\ell>C$. In the latter case, one argues as in~\eqref{eq:Nkwieder}.
In the former case, the modulation of the output is essentially~$2^j$ and $k\le C$ which is excluded.  It remains to consider the null-frame atoms. Thus,
$F=\sum_{\kappa\in\calC_m} F_\kappa$
where $F_\kappa= P_{\ell,\kappa} Q_{\le \ell+2m} F_\kappa$ and $m\le -100$. Due to $F\ne Q_{\le C}F$, one has $\ell+m\ge \ell+2m\ge C$ which implies that the
Fourier support of $\chi F_\kappa$ is essentially that of~$F_\kappa$. Therefore, $\chi F=\sum_{\kappa} \chi F_\kappa$ can be treated as
a wave-packet atom satisfying the bounds
\[
 \sum_{\kappa} \|\chi F_\kappa\|_{\NF[\kappa]}^2 \les A^2 \sum_{\kappa} \|F_\kappa\|_{\NF[\kappa]}^2
\]
Since $k=\ell+O(1)$ we are done.  Next, suppose that $\wh{\chi}(\tau,\xi)$ is supported on~$|\tau|\sim 2^n$ with $n\ge 10$ and~$|\xi|\les 1$.
Then $\|\chi\|_{\Linf}\les 2^{-n}A$. Start with a wave-packet atom~$F$ of the type we just considered. If $n\le k+2m+10$, then $\chi F_\kappa$ has
essentially the same Fourier support as~$F_\kappa$ whence
\[
\|\chi F\|_{N[k]}\les 2^{-k} \Big(\sum_{\kappa} \|\chi F_\kappa\|_{\NF[\kappa]}^2\Big)^{\frac12} \les A 2^{-n}2^{-\ell}
\Big( \sum_{\kappa} \|F_\kappa\|_{\NF[\kappa]}^2\Big)^{\frac12}
\les 2^{-n} A \|F\|_{N[\ell]}
\]
which is summable in $n\ge10$. If $n>k+2m+10$, then $\chi F$ has modulation of size~$2^n$. If $k\ge n$, then
\begin{align*}
 \|\chi F\|_{N[k]} &\les \|\chi F\|_{\dot X_k^{-1,-\frac12,1}} \les 2^{-\frac{n}{2}-k} \|\chi F\|_{\Ltwotx} \\
&\les A\,2^{-\frac{3n}{2}-k} \|F\|_{\Ltwotx} \les A\, 2^{-\frac{3n}{2}-k} 2^{\frac{3\ell}{2}}    \|F\|_{\dot X_\ell^{-1,-\frac12,\infty}} \\
&\les A\, 2^{-n} \|F\|_{N[\ell]}
\end{align*}
where we used \eqref{eq:Xsb_dom} and $\ell=k+O(1)$.
If $k<n$, then
\begin{align*}
 \|\chi F\|_{N[k]} &\les \|\chi F\|_{\dot X_k^{-\frac12+\eps,-1-\eps,2}} \les 2^{-n(1+\eps)} 2^{-k(\frac12-\eps)} \|\chi F\|_{\Ltwotx} \\
&\les A\,2^{-n(2+\eps)} 2^{-k(\frac12-\eps)}  \|F\|_{\Ltwotx} \les A\, 2^{-n(2+\eps)} 2^{k(1+\eps)}     \|F\|_{\dot X_\ell^{-1,-\frac12,\infty}} \\
&\les A\, 2^{-n} \|F\|_{N[\ell]}
\end{align*}
Now suppose that $F$ is a $\dot X_\ell^{-1,-\frac12,1}$-atom with $F=P_\ell Q_j F$. If $j>n+10$, then $\chi F$ is the same kind of atom and one argues as before gaining
a factor of~$2^{-n}$. If $j=n+O(1)$, then
\begin{align*}
 \|\chi F\|_{N[k]} &\les 2^{-k} \|\chi F\|_{\enerN} \les 2^{-\ell} \|\chi\|_{L^2_t L^\infty_x} \|F\|_{\Ltwotx} \\
&\les 2^{-n} A  2^{\frac{j}{2}} \|F\|_{\dot X_\ell^{-1,-\frac12,1}} \les 2^{-\frac{n}{2}} A
\end{align*}
Finally, if $n>j+10$, then $\chi F$ has modulation of size~$2^n$. If $n\le \ell$, then
\[
 \|\chi F\|_{N[k]} \les \|\chi F\|_{\dot X_k^{-1,-\frac12,1}} \les 2^{-\ell-\frac{n}{2}} \|\chi F\|_{\Ltwotx} \les 2^{-n}\, A
\]
whereas in case $n\ge \ell$,  one checks similarly that
\[
 \|\chi F\|_{N[k]} \les \|\chi F\|_{\dot X_k^{-\frac12+\eps,-1-\eps,2}} \les  2^{-n}\, A
\]
as desired. If $F$ is a $\dot X_\ell^{-\frac12+\eps,-1-\eps,2}$-atom
with $F=P_\ell Q_j F$, then analogous arguments lead to a bound of
$\|\chi F\|_{N[\ell]}\les 2^{-\eps n}\, A$ which is again summable
in~$n\ge0$.

Finally, one needs to consider the case where $\wh{\chi}(\tau,\xi)$ is supported on~$|\xi|\sim 2^n$ with $n\ge 10$, say. However, this is easier
due to the rapid decay of $\wh{\chi}$ in~$\xi$. We leave those details to the reader.
\end{proof}

\begin{remark}
 \label{rem:chialso} Lemma~\ref{lem:chiS} of course applies to any space-time Schwartz function~$\chi$. Moreover,
one can check that
the exact same conclusions of Lemma~\ref{lem:chiS} hold for any Schwartz function~$\chi$
which only depends on $t$ and~$x$ alone; the only difference is that $C_0A$
needs to be replaced by $C(\chi)$.
\end{remark}

It is now a simple matter to prove that the $\psi^n_\alpha$ have
uniformly controlled~$S$ norms on some time interval $(-T_0,T_0)$
where $T_0=T_0(V)$.

\begin{cor}
  \label{cor:smalltimeexistence}
Under the assumptions of Lemma~\ref{BasicStability} there exists a
time $T_0=T_0(V)>0$ such that
\begin{equation}
  \label{eq:psinuniform}
\max_{\alpha=0,1,2} \|\psi_\alpha^n\|_{S((-T_0,T_0)\times \R^2)} \le
C(V)<\infty
\end{equation}
uniformly in large~$n$.
\end{cor}
\begin{proof}
Pick $r_{0}>0$
small enough and $R_0$ large enough according to  Lemmas~\ref{lem:DataLocalization} and~\ref{lem:DataLocalization'}, respectively.
In view of~\eqref{eq:x1y1}, Theorem~\ref{thm:smallenergy},
and finite propagation speed, patching up the local evolutions of $(\bfx_1^n,\bfy_1^n)$  shows that the
evolution of  $(\bfx^n,\bfy^n)$ exists on some time interval $(-T_0,T_0)$ uniformly in large~$n$; in fact, one can take $T_0=r_0$.
Note that this part of the argument does not require $(\bfx_2^n,\bfy^n_2)$. These functions are needed to obtain uniform control over $\|\psi_\alpha^n\|_{S((-T_0,T_0)\times \R^2)}$, to which we now turn. The $\phi_\alpha^n$
of the original sequence agree with the $\tilde\phi_\alpha^n$ obtained from~\eqref{eq:x1y1} on the
cone~$K_{x_0,r_0}:=\{(t,x)\:|\: |x-x_0|<r_0-t,\; 0\le t<r_0\}$.
This follows from the construction of~$(\bfx_1^n, \bfy_1^n)$ and finite propagation speed. Note that the $\tilde\phi_\alpha^n$
exist globally in~$\R^{1+2}$ but agree with $\phi_\alpha^n$ only on~$K_{x_0,r_0}$.
A similar observation applies to $(\bfx_2^n,\bfy^n_2)$ on the set $K_{R_0,T_0}:=\{|x|>R_0+t, \:0\le t<T_0\}$.
Cover~$\R^2$ by finitely many $D_j:=D(x_j,r_0)$ as well as the complement of~$D_0:=D(0,R_0)$. This can be done is such a fashion that
there exists a smooth and finite partition of unity $1=\sum_{j=1}^J \chi_j$ on $[0,T_0]\times\R^2$ such that each $\chi_j$ is entirely
supported in either a cone $K_{x_j,r_0}$ or within~$K_{R_0,T_0}$.   Thus
\[
 \psi_\alpha^n = \sum_j \chi_j \psi_\alpha^n = \sum_j \chi_j \wt \psi^{n,j}_\alpha e^{i\Delta^{-1}\del\Re(\tilde\phi^{n,j} - \phi^n)}
\]
Here $\tilde \phi_\alpha^{n,j}$  are the derivative components of the {\em small energy}
wave maps which were constructed by means of Lemmas~\ref{lem:DataLocalization} and~\ref{lem:DataLocalization'}, and  $\tilde \psi_\alpha^{n,j}$
are their gauged counterparts.
If $\chi_j$ has compact support, we now claim that
\[
  \wt \chi_j^n := \chi_j\,  e^{i\Delta^{-1}\del\Re(\tilde\phi^{n,j} - \phi^n)}
\]
satisfies the hypotheses of Lemma~\ref{lem:chiS} with a constant $A$ that can be chosen uniformly in~$n$.
The compact support assumption in time can of course be fulfilled.
Since  for each $j$ and all~$n$
\[
 \chi_j(\tilde\phi_\alpha^{n,j} - \phi_\alpha^n)=0
\]
it follows from the uniform $L^2$ bound on $\tilde\phi_\alpha^{n,j}$ and $\phi_\alpha^n$ that
\begin{equation}
 \label{eq:integ_rep}
\Delta^{-1}\del\Re(\tilde\phi_\alpha^{n,j} - \phi_\alpha^n)(t,x) = \frac{1}{2\pi} \int \frac{x-y}{|x-y|^2} \,\Re(\tilde\phi_\alpha^{n,j} - \phi_\alpha^n)(t,y)\,dy
\end{equation}
is a smooth function relative to~$x$ on the support of~$\chi_j$ with uniform $L^\infty$ bounds on the derivatives (uniform here means relative to large~$n$).
Indeed,
\begin{equation}
 \label{eq:dxest}
 \|\chi_j \nabla_x^\alpha \Delta^{-1}\del\Re(\tilde\phi_\beta^{n,j} - \phi_\beta^n)(t,x)\|_{L^\infty_x} \le C_\alpha \|(\tilde\phi_\beta^{n,j} - \phi_\beta^n)(t,x)\|_{L^2_x} \le C_\alpha \, E
\end{equation}
where $E$ governs the energy uniformly in~$t$. It turns out that we can also incorporate one time derivative into these bounds (but not necessarily any
higher regularity in time). This follows from the div-curl system for~$\phi_\alpha$, see~\eqref{eq:compat1}--\eqref{eq:derwmp2}. Indeed,  if $\alpha\ne0$
then plugging \eqref{eq:compat1} into
\[
 \del_t \Delta^{-1}\del\Re(\tilde\phi_\alpha^{n,j} - \phi_\alpha^n)(x) = \frac{1}{2\pi} \int \frac{x-y}{|x-y|^2} \:\del_0\,\Re(\tilde\phi_\alpha^{n,j} - \phi_\alpha^n)(t,y)\,dy
\]
leads to an expression which is of the schematic form
\[
 \int \frac{x-y}{|x-y|^2} \:\del_\alpha\Re(\tilde\phi_0^{n,j} - \phi_0^n)(t,y)\,dy + \int \frac{x-y}{|x-y|^2} \big[(\tilde \phi^{n,j})^2 - (\phi^n)^2\big](t,y)\, dy
\]
Integrating by parts in the first integral moves the derivative from the $\phi$'s onto the kernel which allows for the same estimate as in~\eqref{eq:dxest}.
As for the second integral on the right-hand side, one has
\[
 \Big\|\nabla_x^\alpha \big[\chi_j \int \frac{x-y}{|x-y|^2} \big[(\tilde \phi^{n,j})^2 - (\phi^n)^2\big](t,y)\, dy\big]\Big\|_{L^\infty_x}\le C_\alpha
\|(\tilde \phi^{n,j})^2 - (\phi^n)^2(t)\|_{L^1_x}  \le C_\alpha \, E
\]
as desired.  If $\alpha=0$, then one uses~\eqref{eq:derwmp2} to arrive at the same conclusion. This establishes our claim concerning the hypotheses of
Lemma~\ref{lem:chiS}; in fact, we  obtained stronger conclusions as far as the conditions for large~$x$ or small $\xi$ are concerned.
Now let us consider the cut-off function $\chi_j$ with unbounded support, which we may assume is $\chi_0$. We can arrange the partition of unity so that $\chi_0(t,x)=\chi_{00}(x)\chi_{01}(t)$
with~$\chi_{01}$ smooth and supported in $(-1,1)$ and with $1-\chi_{00}$ smooth and compactly supported in~$\R^2$.
With $\wt\chi_0^n$ defined as above, we now claim  that
\[
 \chi_{01}(t) - \wt \chi_0^n(t,x)  = \chi_{01}(t)\big(1-  \chi_{00}(x)\,  e^{i\Delta^{-1}\del\Re(\tilde\phi^{n,j} - \phi^n)(t,x)}  \big)
\]
satisfies the requirements of Lemma~\ref{lem:chiS} with a  constant $A$ that  is controlled uniformly in~$n$.
First,
\[
\chi_{01}(t)\chi_{00}(x)\, \Re(\tilde\phi_\alpha^{n,j} - \phi_\alpha^n)(t,x) =0
\]
which shows as before that $\wt \chi_0^n(t,x)$ is smooth in~$x$ with
derivatives that are uniformly bounded in~$L^\infty_x$ relative
to~$n$. In addition, the same arguments involving the div-curl
system allow us to place one~$\del_t$ on $\wt \chi_0^n(t,x)$ without
destroying these conclusions. As for the asymptotic behavior in
$x\to\infty$ and $\xi\to0$, one simply expands
\[
 \frac{x-y}{|x-y|^2} = \frac{x}{|x|^2} + O(\frac{y}{|x|^2})
\]
inside the integral in~\eqref{eq:integ_rep} which is sufficient due to $|y|\les R_0$.
In conclusion, by Lemma~\ref{lem:chiS} and Remark~\ref{rem:chialso}
\begin{align*}
 \|\psi_\alpha^n\|_{S(-T_0,T_0)}& \le \sum_j \|\wt \chi_j^n \wt \psi_\alpha^{n,j} \|_{S(-T_0,T_0)} \le \sum_j  C(\wt \chi_j^n) \| \wt \psi_\alpha^{n,j} \|_{S}
\le \sum_j  C(\wt \chi_j^n) C\eps_0
\end{align*}
is finite uniformly in~$n$.
\end{proof}

The preceding corollary concludes the proof of
Lemma~\ref{BasicStability} up to the assertion about the frequency
envelope at the end. This will be proved in Section~\ref{sec:missing
proofs}.

We close this section with an important strengthening of the bound on~$\psi_L$ from Lemma~\ref{lem:LocalSplitting}.
More specifically, we prove that the intervals $I_j$ can be chosen in such a way that the estimate~\eqref{eq:psiLbd} only
depends on the energy of~$\psi$. This will play an important role later on.
In order to achieve this property, we require an improvement over Lemma~\ref{lem:easytrilin}. We begin with the following technical statements
which allow us to make a better choice of the intervals~$I_j$ in the proof of Lemma~\ref{lem:LocalSplitting}.

\begin{lemma}
 \label{lem:intervals2} Let $\|\psi\|_S<C_0$ and $\eps_0>0$ be arbitrary. Then there exists a partition of~$\R$ into intervals
$\{I_j\}_{j=1}^M$ which depend on~$\psi$ but with $M=M(\eps_0,C_0)$ and which satisfy
\[
 \max_{1\le j\le M} \sum_{k\in\Z} \|P_k (\psi\; |\nabla|^{-1} \psi^2)\|_{L^2_t(I_j; \dot H^{-\frac12})}^2 \le \eps_0
\]
where $\nabla=\nabla_x$ and $\psi\; |\nabla|^{-1} \psi^2$ is  schematic notation which stands for any one of the nonlinearities appearing
on the right-hand side of the div-curl system \eqref{eq:psisys1}, \eqref{eq:psisys2}.
\end{lemma}
\begin{proof}
It of course suffices to show that
\begin{equation}\label{eq:nointervals}
 \sum_{k\in\Z} \|P_k [\psi_1\; |\nabla|^{-1} (\psi_2\psi_3)]\|_{L^2_t(\R; \dot H^{-\frac12})}^2 \les \prod_{i=1}^3 \|\psi_i\|_S^2
\end{equation}
We begin with the case where $\psi_2\psi_3$ is replaced by $ I^c\psi_2\cdot\psi_3$. It is easy to see that
\[
 \|P_k (I^c\psi_2\cdot\psi_3)\|_{L^2_{t,x}} \les 2^{\frac{k}{2}} \|\psi_2\|_S \|\psi_3\|_S
\]
Then by the
usual trichotomy,
\begin{align*}
 \sum_{k\in\Z} 2^{-k} \| P_k(\psi_1 \,|\nabla|^{-1}(I^c\psi_2\cdot \psi_3)) \|_{L^2_{t,x}}^2 &\les
\sum_{k\in\Z} 2^{-k} \| P_k \psi_1\: |\nabla|^{-1} P_{<k-5} (I^c\psi_2\cdot \psi_3) \|_{L^2_{t,x}}^2 \\
&+ \sum_{k\in\Z} 2^{-k} \Big(\sum_{\ell>k} \|P_k [ \tilde P_\ell \psi_1 \:|\nabla|^{-1} P_{\ell} (I^c\psi_2 \cdot \psi_3) ] \|_{L^2_{t,x}}\Big)^2 \\
&+\sum_{k\in\Z} 2^{-k} \| P_{<k-5} \psi_1\: |\nabla|^{-1} P_{k} (I^c\psi_2 \cdot \psi_3) \|_{L^2_{t,x}}^2 \les \prod_{i=1}^3 \|\psi_i\|_S^2
\end{align*}
Hence, we may assume that the two inner inputs are both hyperbolic,
i.e.,  $\psi_i=Q_{\le k_i}\psi_i$ for $i=2,3$.

\noindent Now implement the Hodge decomposition for the inputs of
$|\nabla|^{-1}(\psi^{2})$, i.e., write
\[
\psi_\alpha=R_\alpha\psi+\chi_\alpha
\]
We begin by considering the resulting trilinear expressions, more specifically the one where
the inner null-form is hyperbolic: Suppressing the indices on~$\psi$ for simplicity,
\begin{equation}
 \label{eq:psiNab}
\sum_{k\in\Z} \|P_k (\psi |\nabla|^{-1} I \calN_{\alpha\beta} (\psi,\psi))\|_{L^2_t(\R; \dot H^{-\frac12})}^2 \le C\|\psi\|_S^6
\end{equation}
where $\calN_{\alpha j}$ is the null-form from Definition~\ref{def:Nab}.
As usual, this splits into the high-low, high-high, and low-high cases:
\begin{align*}
  \sum_{k\in\Z} 2^{-k} \|P_k (\psi |\nabla|^{-1} I \calN_{\alpha\beta} (\psi,\psi))\|_{L^2_{t,x}}^2
&\les \sum_{k\in\Z} 2^{-k} \|P_k (\tilde P_k\psi |\nabla|^{-1} P_{<k-5} I \calN_{\alpha\beta} (\psi,\psi))\|_{L^2_{t,x}}^2  \\
&+ \sum_{k\in\Z} 2^{-k} \Big(\sum_{\ell>k}  \|P_k (P_\ell\psi  \tilde P_\ell |\nabla|^{-1} I \calN_{\alpha\beta} (\psi,\psi))\|_{L^2_{t,x}}\Big)^2  \\
&+ \sum_{k\in\Z} 2^{-k} \|P_k (P_{<k-5} \psi \tilde P_k |\nabla|^{-1} I \calN_{\alpha\beta} (\psi,\psi))\|_{L^2_{t,x}}^2  \\
&=: A + B +C
\end{align*}
Next, one writes $A\le A_1+A_2+A_3$ reflecting the high-high, high-low, and low-high decomposition of the $\calN_{\alpha\beta}$-nullform. Thus,
by Lemma~\ref{lem:Nablowmod},
\begin{align*}
 A_1 &\le \sum_{k\in\Z} 2^{-k} \Big(\sum_{k_0<k-5} \;\sum_{k_2=k_3+O(1)>k_0-5} \|P_k (\tilde P_k\psi |\nabla|^{-1} P_{k_0} I \calN_{\alpha\beta} (P_{k_2}\psi,P_{k_3}\psi))\|_{L^2_{t,x}}\Big)^2 \\
&\les \sum_{k\in\Z} 2^{-k} \Big(\sum_{k_0<k-5} \; \sum_{k_2=k_3+O(1)>k_0-5} \|\tilde P_k\psi\|_{\ener} \|  P_{k_0} I \calN_{\alpha\beta} (P_{k_2}\psi,P_{k_3}\psi))\|_{L^2_{t,x}}\Big)^2 \\
&\les \sum_{k\in\Z} 2^{-k} \|\tilde P_k\psi\|_{\ener}^2 \Big(\sum_{k_0<k-5}\; \sum_{k_2=k_3+O(1)>k_0-5} 2^{k_0-\frac{k_2}{2}} \| P_{k_2}\psi\|_{S[k_2]} \|P_{k_3}\psi\|_{S[k_3]} \Big)^2 \les \|\psi\|_S^6
\end{align*}
Similarly, by Lemma~\ref{lem:Nablowmod2},
\begin{align*}
 A_2 &\le \sum_{k\in\Z} 2^{-k} \Big(\sum_{k_2+O(1)=k_0<k-5}\; \sum_{k_3<k_0-5} \|P_k (\tilde P_k\psi |\nabla|^{-1} P_{k_0} I \calN_{\alpha\beta} (P_{k_2}\psi,P_{k_3}\psi))\|_{L^2_{t,x}}\Big)^2 \\
&\les \sum_{k\in\Z} 2^{-k} \Big(\sum_{k_2+O(1)=k_0<k-5}\; \sum_{k_3<k_0-5}   \|\tilde P_k\psi\|_{\ener} \|  P_{k_0} I \calN_{\alpha\beta} (P_{k_2}\psi,P_{k_3}\psi))\|_{L^2_{t,x}}\Big)^2 \\
&\les \sum_{k\in\Z} 2^{-k} \|\tilde P_k\psi\|_{\ener}^2 \Big(\sum_{k_2+O(1)=k_0<k-5}\; \sum_{k_3<k_0-5}  2^{(\frac12-\eps)k_3} 2^{\eps k_0} \| P_{k_2}\psi\|_{S[k_2]} \|P_{k_3}\psi\|_{S[k_3]}    \Big)^2
\les \|\psi\|_S^6
\end{align*}
This concludes the bound on~$A$ since $A_3$ is of course symmetric to~$A_2$. Next, with $B\le B_1+B_2+B_3$ via the same trichotomy,
\begin{align*}
 B_1 &\les \sum_{k\in\Z} 2^{-k} \Big(\sum_{\ell>k} 2^{k} \|P_\ell\psi\|_{\ener} \; 2^{-\ell}\weg  \sum_{k_2=k_3+O(1)>\ell-5} \|  \tilde P_\ell  I \calN_{\alpha\beta} (P_{k_2} \psi, P_{k_3} \psi))\|_{L^2_t L^2_x}\Big)^2\\
&\les  \sum_{k\in\Z} 2^{-k} \Big(\sum_{\ell>k} 2^{k} \|P_\ell\psi\|_{\ener}\; 2^{-\ell}\weg  \sum_{k_2=k_3+O(1)>\ell-5}  2^{\ell-\frac{k_2}{2}} \| P_{k_2}\psi\|_{S[k_2]} \|P_{k_3}\psi\|_{S[k_3]}  \Big)^2 \les \|\psi\|_S^6
\end{align*}
by Lemma~\ref{lem:Nablowmod}, whereas
\begin{align*}
 B_2 &\les \sum_{k\in\Z} 2^{-k} \Big(\sum_{\ell>k} 2^{k} \|P_\ell\psi\|_{\ener}\; 2^{-\ell} \weg \sum_{\ell=k_2+O(1)>k_3-5 } \|  \tilde P_\ell  I \calN_{\alpha\beta} (P_{k_2} \psi, P_{k_3} \psi))\|_{L^2_t L^2_x}\Big)^2\\
&\les  \sum_{k\in\Z} 2^{-k} \Big(\sum_{\ell>k} 2^{k} \|P_\ell\psi\|_{\ener} \; 2^{-\ell}  \sum_{\ell=k_2+O(1)>k_3-5  } 2^{(\frac12-\eps)k_3} 2^{\eps k_2}  \| P_{k_2}\psi\|_{S[k_2]} \|P_{k_3}\psi\|_{S[k_3]}  \Big)^2 \les \|\psi\|_S^6
\end{align*}
by Lemma~\ref{lem:Nablowmod2}.
The low-high case of \eqref{eq:psiNab} is treated in an analogous fashion and we skip it.

Next, we treat the case where the inner null-form is elliptic. Then
the desired bound reads
\begin{equation}
 \label{eq:psiNabhigh}
\sum_{k\in\Z} \|P_k (\psi\, |\nabla|^{-1} I^c \calN_{\alpha\beta} (\psi,\psi))\|_{L^2_t(\R; \dot H^{-\frac12})}^2 \le C\|\psi\|_S^6
\end{equation}
As before,  $A\le A_1+A_2+A_3$ reflecting the high-high, high-low, and low-high decomposition of the $\calN_{\alpha\beta}$-nullform.
{\em We will first excluded the contributions by opposing high-high interactions in the null-form}, cf.~Remark~\ref{rem:log_loss}.
Hence,  by \eqref{eq:sch113} {\em without} the $\la k_1-k\ra^2$ loss,
\begin{align*}
 A_{1} &\les \sum_{k\in\Z} 2^{-k} \Big(\sum_{k_0<k-5} \;\sum_{k_2=k_3+O(1)>k_0-5} \|P_k (\tilde P_k\psi\, |\nabla|^{-1} P_{k_0} Q_{>k_0} \calN_{\alpha\beta} (P_{k_2}\psi,P_{k_3}\psi))\|_{L^2_{t,x}}\Big)^2 \\
&\les \sum_{k\in\Z} 2^{-k} \Big(\sum_{k_0<k-5} \;\sum_{k_2=k_3+O(1)>k_0-5} \|P_k\psi\|_{\ener} \|  P_{k_0} Q_{>k_0} \calN_{\alpha\beta} (P_{k_2}\psi,P_{k_3}\psi)\|_{L^2_{t,x}}\Big)^2 \\
&\les \sum_{k\in\Z} 2^{-k} \Big(\sum_{k_0<k-5} \;\sum_{k_2=k_3+O(1)>k_0-5} \|P_k\psi\|_{\ener} 2^{\frac{k_0}{2}} \,
\| P_{k_2}\psi\|_{S[k_2]}  \|P_{k_3}\psi\|_{S[k_3]}\Big)^2 \les \|\psi\|_S^6
\end{align*}
For $A_2$ one proceeds similarly, using Lemma~\ref{lem:Nabhighmod2} instead. In fact, due to the hyperbolic nature of $\psi_1,\psi_3$,
\begin{align*}
 A_{2} &\les \sum_{k\in\Z} 2^{-k} \Big(\sum_{k_2+O(1)=k_0<k-5}\; \sum_{k_3<k_0-5}  \|P_k (\tilde P_k\psi\, |\nabla|^{-1} P_{k_0} I^c \calN_{\alpha\beta} (P_{k_2}\psi,P_{k_3}\psi))\|_{L^2_{t,x}}\Big)^2 \\
&\les \sum_{k\in\Z} 2^{-k} \Big(\sum_{k_2+O(1)=k_0<k-5}\; \sum_{k_3<k_0-5}   \|P_k\psi\|_{\ener} \|  P_{k_0}  \tilde Q_{k_0} \calN_{\alpha\beta} (P_{k_2}\psi,P_{k_3}\psi)\|_{L^2_{t,x}}\Big)^2 \\
&\les \sum_{k\in\Z} 2^{-k} \Big(\sum_{k_2+O(1)=k_0<k-5}\; \sum_{k_3<k_0-5}  \|P_k\psi\|_{\ener}  2^{\frac{k_3}{2}}
\| P_{k_2}\psi\|_{S[k_2]}  \|P_{k_3}\psi\|_{S[k_3]}\Big)^2 \les \|\psi\|_S^6
\end{align*}
This concludes the high-low case~$A$.
In the high-high case we write $B\le B_1+B_2+B_3$ as before. Therefore,
\begin{align*}
 B_1 &\les \sum_{k\in\Z} 2^{-k} \Big(\sum_{\ell>k} 2^{k} \|  P_\ell\psi\|_{\ener} \; 2^{-\ell}\weg \weg
\sum_{k_2=k_3+O(1)>\ell-5} \|  \tilde P_\ell Q_{>\ell}   \calN_{\alpha\beta} (P_{k_2} \psi, P_{k_3} \psi)\|_{L^2_{t,x}}\Big)^2\\
&\les  \sum_{k\in\Z}  \Big(\sum_{ \ell>k} 2^{\frac{k}{2}}
\|P_\ell\psi\|_{S[\ell]}\; 2^{-\ell}\weg \weg
\sum_{k_2=k_3+O(1)>\ell-5}\weg\weg  2^{\frac{\ell}{2}}  \|
P_{k_2}\psi\|_{S[k_2]} \|P_{k_3}\psi\|_{S[k_3]}  \Big)^2 \les
\|\psi\|_S^6
\end{align*}
by Lemma~\ref{lem:Nabhighmod}, whereas
\begin{align*}
 B_2 &\les \sum_{k\in\Z} 2^{-k} \Big(\sum_{ \ell>k} 2^{k} \|P_\ell Q_j \psi\|_{ L^2_{t,x}}\; 2^{-\ell} \weg\weg
\sum_{\ell=k_2+O(1)>k_3-5 } \|  \tilde P_\ell  \tilde Q_\ell \calN_{\alpha\beta} (P_{k_2} \psi, P_{k_3} \psi)\|_{L^2_t L^2_x}\Big)^2\\
&\les  \sum_{k\in\Z}  \Big(\sum_{\ell>k} 2^{\frac{k}{2}}  \|P_\ell\psi\|_{S[\ell]} \; 2^{-\ell}\weg \weg
\sum_{\ell=k_2+O(1)>k_3-5  }  2^{\frac{k_3}{2}}      \| P_{k_2}\psi\|_{S[k_2]} \|P_{k_3}\psi\|_{S[k_3]}  \Big)^2 \les \|\psi\|_S^6
\end{align*}
by Lemma~\ref{lem:Nabhighmod2} which finishes the analysis of~$B$.
We again leave the low-high case to the reader.

It remains the bound the contributions by the opposing high-high waves in the inner null-form.
Returning to the $\psi_1,\psi_2,\psi_3$ notation, we may assume that $\psi_i=Q_{\le k_i}\psi_i$ for $i=2,3$ and that
there is an angular separation of the Fourier supports of $\psi_1$ and $\psi_2$, say (since the Fourier supports of $\psi_2,\psi_3$ make
a large angle). Hence we may bound the missing contribution to~$A_1$ as follows, where we ignore the nullform and replace the
outer $|\nabla|^{-1}$ with a weight by the usual convolution  logic:
\begin{align*}
 A_{1} &\les \sum_{k\in\Z} 2^{-k} \Big(\sum_{k_0<k-5} \;\sum_{k_2=k_3+O(1)>k_0-5} 2^{-k_0} \big[\sum_{c\in\calD_{k,k_0-k}}
 \|P_{c}\psi_1 P_{k_0}( P_{k_2}\psi_2 \, P_{k_3}\psi_3 )\|_{L^2_{t,x}}^2 \big]^{\frac12} \Big)^2
\end{align*}
We now invoke \eqref{eq:bilin2} to conclude that
\begin{align*}
\big[\sum_{c\in\calD_{k,k_0-k}}  \|P_{c}\psi_1 \, P_{k_0}(  P_{k_2}\psi_2 \, P_{k_3}\psi_3 ) \|_{L^2_{t,x}}^2 \big]^{\frac12}  &\les  2^{k_0}
 \big[ \sum_{c\in\calD_{k,k_0-k}}   \|P_{c}\psi_1  P_{k_2}\psi \, P_{k_3}\psi \|_{L^2_t L^1_x}^2  \big]^{\frac12} \\
& \les 2^{k_0} \big[ \sum_{c\in\calD_{k,k_0-k}}   \|P_{c}\psi_1  P_{k_2}\psi\|_{L^2_{t,x}}^2  \|P_{k_3}\psi \|_{\ener}^2 \big]^{\frac12}  \\
&\les 2^{\frac{3k_0}{2}} \la k_0-k\ra\big[ \sum_{c\in\calD_{k,k_0-k}}    \|P_{c}\psi_1\|_{S[k]}^2   \|  P_{k_2}\psi\|_{S[k_2]}^2  \big]^{\frac12}
 \|P_{k_3}\psi \|_{\ener} \\
&\les  2^{\frac{3k_0}{2}} \la k_0-k\ra \prod_{i=1}^3 \|P_{k_i}\psi_i\|_{S[k_i]}
\end{align*}
The loss of $\la k_0-k\ra$ here is due to the usual issue of wave-packets which are too thick resulting in the need for Lemma~\ref{lem:incl_free}.
Inserting this into the bound on~$A_1$ yields
\begin{align*}
 A_{1} &\les \sum_{k\in\Z}  \Big(\sum_{k_0<k} \; 2^{\frac{k_0-k}{2}}\la k_0-k\ra  \|P_{k} \psi_1\|_{S[k]}    \Big)^2 \|\psi_2\|_S^2 \|\psi_3\|_S^3
\les   \prod_{i=1}^3 \|\psi_i\|_{S}^2
\end{align*}
as desired. The opposing high-high contributions to the other terms are similar and omitted.
We still need to control the contributions from the elliptic terms $\chi$, leading to higher order nonlinearities. This is again done in the appendix.
\end{proof}

We can now state the refined version of
Lemma~\ref{lem:LocalSplitting} which gives better control over the
linear wave $\psi_L$. As in that lemma~$\psi$ are the gauged
components of an admissible wave map locally on some time interval
$[-T_0,T_1]$.

\begin{cor}
 \label{cor:localsplit2}
Let $\|\psi\|_{S}<C_0$.  Given $\eps_0>0$, there exist
$M_{1}=M_{1}(C_{0}, \eps_0)$ many intervals $I_{j}$ as
in~\eqref{eq:Ij} with the following property: for each
$I_{j}=(t_{j}, t_{j+1})$, there is a decomposition
\[
\psi|_{I_{j}}=\psi_{L}^{(j)}+\psi_{NL}^{(j)},\quad
\Box\psi_{L}^{(j)}=0
\]
which satisfies
\begin{align}
\sum_{k\in\Z}\|P_{k}\psi_{NL}^{(j)}\|_{S[k](I_{j}\times\R^{2})}^{2} &<\eps_0 \label{eq:psiNLbd2} \\
\|\nabla_{x,t}\psi_{L}^{(j)}\|_{L_{t}^{\infty}\dot{H}^{-1}} &\les
\eps_0^{-\frac14}(E+1)E  \label{eq:psiLbd2}
\end{align}
where the implied constant in the last inequality is universal and
$E=\|\psi(t)\|_2$ is the conserved energy. In particular,
\[
\|\psi_{NL}^{(j)}\|_{S(I_j\times\R^2)} \,
\|\nabla_{x,t}\psi_L^{(j)}\|_{\dot H^{-1}} \ll 1
\]
by choosing~$\eps_0$ small enough depending on the energy.
\end{cor}
\begin{proof} We first prove \eqref{eq:psiNLbd2}
and~\eqref{eq:psiLbd2} by following the strategy of the proof of
Lemma~\ref{lem:LocalSplitting}; however, we use
Lemma~\ref{lem:intervals2} instead of Lemma~\ref{lem:easytrilin}
when the underlying time interval is small. More precisely, consider
the frequency component~$P_0\psi_\alpha$.

\noindent {\em Case 1:} The underlying time
interval~$I_0:=(-T_0,T_1)$ satisfies $|I_0|<\eps_1$ with an $\eps_1$
that is to be determined. The main property of this parameter is
that it can be chosen to be an absolute constant independently
of~$C_0$.
 The $\psi_\alpha$ satisfy the
system~\eqref{eq:psisys1}--\eqref{eq:psi_wave}. Schematically,  this
system takes the form
\begin{align*}
\partial_{t}P_{0}\psi_j &=\del_j
P_{0}\psi_0+P_{0}[\psi\nabla^{-1}(\psi^{2})],\quad j=1,2\\
\partial_{t}P_{0}\psi_0 &=\sum_{j=1}^2 \del_j
P_{0}\psi_j+P_{0}[\psi\nabla^{-1}(\psi^{2})]
\end{align*}
where the nonlinearity is schematically. Now define the linear wave
$P_0\psi_L$ to be
\begin{align*}
   P_0\psi_{L,j} &:= S(t)(P_0\psi_j(0), \del_j
P_{0}\psi_0),\quad j=1,2\\
P_0\psi_{L,0} &:=  S(t)\big(P_{0}\psi_0, \sum_{j=1}^2
P_0\del_j\psi_j(0)\big)
\end{align*}
whereas $P_0 \psi_{NL,\alpha}:= P_0 \psi_\alpha - P_0
\psi_{L,\alpha}$. Thus, for $j=1,2$,
\begin{align*}
  P_{0} \psi_j(t) &= P_{0}\psi_j(0) + \int_0^t P_0 \del_j
P_{0}\psi_0(s)\, ds + \int_0^t P_{0}[\psi\nabla^{-1}(\psi^{2})](s)\,
ds \\
P_{0}\psi_{L,j}(t) &= P_{0}\psi_j(0) + O_{L^2_x}(t^2)
\end{align*}
and similarly for $\psi_0$,  whence  for all $t\in I_0$,
\begin{align*}
\|P_{0}\psi_{NL}(t)\|_{L^2_x} &\leq t^2 \|P_{0}\psi(0)\|_{L^2_x}
+\Big\|\int_{0}^{t}P_{0}[\psi\nabla^{-1}(\psi^{2})](s,\cdot)\,ds
\Big\|_{L_{x}^{2}} \\
&\les t^2 \|P_{0}\psi(0)\|_{L^2_x} + |t|^{\frac12} \big\|
P_{0}[\psi\nabla^{-1}(\psi^{2})] \big\|_{L_{t,x}^{2}}
\end{align*}
In other words,
\[
\|P_{0}\psi_{NL}\|_{L^\infty_t(I_0;L^2_x)} \les \eps_1^2
\|P_{0}\psi(0)\|_{L^2_x} + \eps_1^{\frac12} \|
P_{0}[\psi\nabla^{-1}(\psi^{2})] \|_{L^2_t(I_0;L^2_x)}
\]
As in the proof of Lemma~\ref{lem:LocalSplitting} one concludes from
this that
\[
\|P_{0}\psi_{NL}\|_{S[0](I_0\times \R^2)} \les \eps_1^2
\|P_{0}\psi(0)\|_{L^2_x} + \eps_1^{\frac12} \|
P_{0}[\psi\nabla^{-1}(\psi^{2})] \|_{L^2_t(I_0;L^2_x)}
\]
Rescaling this bound to general $2^k$ yields the following. Suppose
$|I|\le \eps_1 2^{-k}$. Then
\[
\|P_{k}\psi_{NL}\|_{S[0](I_0\times \R^2)} \les \eps_1^2 \|P_k
\psi(0)\|_{L^2_x} + \eps_1^{\frac12} 2^{-\frac{k}{2}} \|
P_{0}[\psi\nabla^{-1}(\psi^{2})] \|_{L^2_t(I_0;L^2_x)}
\]
Now provided $I_0\subset I_j$ where $\{I_j\}_{j=1}^M$,
$I_j=I_j(\tilde\eps_0,\psi)$, $M=M(\tilde\eps_0,C_0)$ are the
intervals constructed in Lemma~\ref{lem:intervals2}, one concludes
that
\begin{equation}
\label{eq:fall1NL} \sum_{k:|I_0|\le \eps_1 2^{-k}}
\|P_{k}\psi_{NL}\|^2_{S[k](I_0\times \R^2)} \les \eps_1^4
\|\psi(0)\|^2_{L^2_x} + \tilde\eps_0\eps_1
\end{equation}
where $\tilde\eps_0$ is a separate smallness parameter. We now pick
$\eps_1:=\eps_0^{\frac14}(1+E)^{-1}$ and
$\tilde\eps_0:=\eps_0^{\frac34}$ where $E=\|psi(t)\|_{L^2_x}$ is the
conserved energy of $\psi$ (for this one needs to remain on the
interval on which $\psi$ equals the gauged derivative components of
a wave map). This renders the right-hand side of~\eqref{eq:fall1NL}
less than~$\eps_0$.

\noindent  As already explained in the proof of
Lemma~\ref{lem:LocalSplitting}, we will use this analysis also in
the case of large intervals to which we now turn. However, in that
case the estimates obtained here allow one to control the term
$\|\psi\big|_{[-T_0,T_0]}\|_{S}$ in~\eqref{eq:IIc} of
Section~\ref{subsec:waveeq}.

\medskip \noindent {\em Case 2:} The underlying time
interval~$I_0=(-T_0,T_1)$ satisfies $|I_0|>\eps_1$ where $\eps_1$ is
as in Case~1 (again for the $P_0$ frequencies). Here the analysis of
Case~2 of Lemma~\ref{lem:LocalSplitting} applies verbatim, leading
to intervals $\{I_j'\}_{j=1}^{M'}$, with $M'=M'(\eps_0,C_0)$ such
that
\[
\max_{1\le j\le M'} \sum_{k\in\Z} \|P_k F_\alpha\|^2_{N[k](I_j\times
\R^2)} <\eps_0
\]
where $F=\sum_\alpha F_\alpha$ stands for the right-hand side
of~\eqref{eq:psi_wave} as usual.

\medskip\noindent Now we take the intersections of the intervals
$I_j$ and $I_k'$ which appeared in Cases~1 and~2 above. Denote this
collection again by~$\{I_j\}_{j=1}^M$ with $M=M(\tilde
\eps_0,\eps_0,C_0)$. Fix such an~$I_j$. Given $k\in Z$, we define
$P_k \psi_L^{(j)}$ to be the free evolution of $(I\psi)[t_0]$ where
$t_0\in I_j$ is the center of~$I_j$, whereas $P_k \psi_{NL}^{(j)}$
is everything else. By our construction,
\[
\sum_{k:|I_j|\le \eps_1 2^{-k}} \|P_k
\psi_{NL}^{(j)}\|^2_{S[k](I_j\times \R^2)} \le \eps_0
\]
Combining this with \eqref{eq:fall1NL} this bound
implies~\eqref{eq:psiNLbd2}. As for the linear wave~$\psi_L^{(j)}$,
we note that those~$k$ which belong to Case~1 yield
\[
\|P_k\psi_L^{(j)}\|_{S[k](I_j\times\R^2)} \les \|P_k\psi\|_2
\]
with an absolute implicit constant, whereas \eqref{eq:IIc} from
Section~\ref{subsec:waveeq} yields the bound
\[
\|P_k\psi_L^{(j)}\|_{S[k](I_j\times\R^2)} \les \eps_1^{-1}
\|P_k\psi\|_2
\]
These estimates imply~\eqref{eq:psiLbd2}.
\end{proof}
\begin{remark}\label{rem:psiLimprov} Note that if we apriori work on a time interval of length $\geq 1$, say, the statement of the Corollary may be strengthened to
\[
\|\nabla_{x,t}\psi^{(j)}_{L}\|_{L_t^\infty\dot{H}^{-1}}\lesssim E
\]
with universal implied constant. Indeed, in this case, the 'time averaging' around the initial data does not cost a large constant.
\end{remark}

Later we shall need to following corollary which further specifies
the Fourier support of~$\psi_L$.

\begin{cor}
 \label{cor:psiLNLsupp}  Let $\|\psi\|_S<C_0$. Assume that $\psi=\tilde\psi+\breve{\psi}$ where for some $b$
\[
 \|\breve{\psi}\|_{S} + \|P_{(-\infty,b]^c} \tilde\psi\|_S <\delta_1
\]
for some small $\delta_1$. Then there exist intervals $\{I_j\}_{j=1}^{M_1}$ as in Corollary~\ref{cor:localsplit2}
so that on each $I_j$ one has a decomposition
\[
 \psi|_{I_{j}}=\psi_{L}^{(j)}+\psi_{NL}^{(j)},\quad
\Box\psi_{L}^{(j)}=0
\]
where furthermore $\psi_{L}^{(j)} = \tilde\psi_{L}^{(j)} + \breve{\psi}_{L}^{(j)}$ and
$\psi_{NL}^{(j)} = \tilde\psi_{NL}^{(j)} + \breve{\psi}_{NL}^{(j)}$ where both $\tilde\psi_{L}^{(j)}$ and $\breve{\psi}_{L}^{(j)}$
are free waves satisfying~\eqref{eq:psiLbd2} and both $\tilde\psi_{NL}^{(j)}$ and $\breve{\psi}_{NL}^{(j)}$ satisfy~\eqref{eq:psiNLbd2}.
Furthermore,
\begin{align*}
 \|\breve{\psi}_L^{(j)}\|_{S} + \|P_{(-\infty,b]^c} \tilde\psi_L^{(j)}\|_S &\les \delta_1\\
\|\breve{\psi}_{NL}^{(j)}\|_{S} + \|P_{(-\infty,b]^c} \tilde\psi_{NL}^{(j)}\|_S &\les \delta_1
\end{align*}
with an absolute implicit constant.
\end{cor}
\begin{proof}
 The proof of this statement follows the exact same lines as the proof of the previous corollary. The only difference
is that each nonlinearity needs to be split into the contributions made by $\breve{\psi}$ and $\tilde\psi$, respectively.
\end{proof}

\section{$\bmo$, $A_p$, and weighted commutator estimates}
\label{sec:bmo}

In this section we develop some auxiliary tools that will be needed
in the implementation of the Bahouri-Gerard theory for wave maps.
More specifically, due to the lack of an imbedding from energy
to~$L^\infty$ in the critical case we need to invoke methods
involving $\bmo$ and the closely related~$A_p$-classes  in order to
carry out Steps~1 and~2 of the program delineated in Section~\ref{sec:intro}. Lemma~\ref{lem:embed} will play a
crucial role here. Moreover, we require a weighted  version of the Coifman-Meyer commutator theorem, with
the weights belonging to the
$A_p$-class.
Although it does not seem to be widely known, it is an easy
consequence of the standard theory and we sketch the proof for the
sake of completeness. The paper~\cite{ST} contains a more general
form of this result. A Calderon-Zygmund kernel here is defined to be
any linear operator $T$ bounded on~$L^2$ with the additional
property that for any $f\in L^2$ with compact support and all
$x\not\in \supp(f)$,
\[
Tf(x) = \int K(x,y) f(y)\, dy
\]
where $|K(x,y)|\le C|x-y|^{-d}$ and for some $0<\gamma\le 1$,
\begin{align*}
  |K(x,y)-K(x',y)| &\le C\frac{|x-x'|^\gamma}{|x-y|^{d+\gamma}}
  \quad\forall\; |x-y|>2|x-x'|\\
|K(x,y)-K(x,y')| &\le C\frac{|y-y'|^\gamma}{|x-y|^{d+\gamma}}
  \quad\forall\; |x-y|>2|y-y'|
\end{align*}
By the Calderon-Zygmund theorem, any such $T$ is also bounded
on~$L^p(\R^d)$ provided $1<p<\infty$.

\begin{lemma}\label{lem:CoifMeyer}
   Let $1<p<\infty$.   There exists $\delta=\delta(p)>0$ with the following
   property: suppose $\phi=\phi_0+\phi_1$ where
   $\|\phi_0\|_{\bmo(\R^d)}<\delta$ and $\|\phi_1\|_{L^\infty(\R^d)}\le A$. Then
   \begin{equation}
     \label{eq:commute}
     \| e^{-\phi} [T,b] e^\phi\|_{p\to p}\le C(d,A,T,p)\|b\|_{\bmo}
   \end{equation}
   for any Calderon-Zygmund operator $T$ and $b\in\bmo$. Moreover,
   $\inf_{p\in I}\delta(p)>0$ and $\sup_{p\in I} C(d,A,T,p)<\infty$ for any compact $I\subset (1,\infty)$.
\end{lemma}
\begin{proof}
  Since $\phi_1$ contributes at most $e^{2A}$ to the estimate, we
  can assume that $\phi=\phi_0$ with small $\bmo$ norm. In particular, $e^{\phi_0}\in A_p$.  We will require the following inequality involving the
so-called sharp maximal function $M^\sharp f$ which is defined as
\[
(M^\sharp f)(x) =  \sup_{Q:x\in Q} \inf_{c}  |Q|^{-1}\int_Q \big| f(y)- c\big|\,dy
\]
where $c$ is a constant. The optimal choice of $c$ is
$c=f_Q:=|Q|^{-1}\int_Q f(y)\,dy$. The estimate then reads (see
Theorem~7.10 in~\cite{Duo})
\begin{equation}
  \label{eq:sharpM} \int_{\R^d} (M f)^p(x)\,w(x)\, dx\le C \int_{\R^d} (M^\sharp f)^p(x)\,w(x)\, dx
\end{equation}
for any $w\in A_\infty$ with a constant that only depends on the
dimension and the constants in~\eqref{eq:Ainfty}. To avoid
trivialities like $f=\const$ for which~\eqref{eq:sharpM} fails, one
needs to assume $Mf\in L^{p_0}(\R^d)$ for some $1\le p_0\le p$.

The proof of \eqref{eq:commute} combines the standard proof of the
unweighted Coifman-Meyer bound with the sharp function
estimate~\eqref{eq:sharpM}. More precisely, fix a cube~$Q$ and write
\begin{align*}
  [T,b]f &= -(b-b_Q) Tf + T((b-b_Q)\chi_{2Q} f) + T((b-b_Q)\chi_{\R^d\setminus 2Q}
  f)\\
  &=: A_Q + B_Q + C_Q
\end{align*}
To bound $M^\sharp ([T,b]f)$, we simply note that for any $x\in Q$,
and any $1<s<\infty$,
\begin{align*}
  &|Q|^{-1} \int_Q (|A_Q(y)|+|B_Q(y)|)\, dy \le  C(s,d,T)\|b\|_{\bmo} \big(
  (M|Tf|^s)^{\frac1s}(x) + (M|f|^s)^{\frac1s}(x)\big)
\end{align*}
Indeed, for $A$ this follows from H\"older's inequality and the
definition of $\bmo$, whereas for $B$ we also invoke the $L^q$
boundedness of $T$ for some $1<q<s$.  For $C_Q$ we let $y_Q$ be the
center of~$Q$ and estimate for any $y\in Q$,
\begin{align*}
|C_Q(y)-C_Q(y_Q)| &\le \int_{\R^d\setminus 2Q}
|K(y,z)-K(y_Q,z)\|(b-b_Q)(z)\|f(z)|\, dz \\
&\le C \int_{\R^d\setminus 2Q}
\frac{|y-y_Q|^\gamma}{|z-y_Q|^{d+\gamma}}
|(b-b_Q)(z)\|f(z)|\, dz\\
&\le C \|b\|_{\bmo} \inf_{x\in Q} (M|f|^s)^{\frac{1}{s}}(x)
\end{align*}
where $\gamma>0$ is as above.

\noindent In conclusion,
\begin{align*}
M^\sharp ([T,b]f) \le C(s,d,T)\|b\|_{\bmo} \big(
  (M|Tf|^s)^{\frac1s} + (M|f|^s)^{\frac1s}\big)
\end{align*}
The lemma follows from~\eqref{eq:sharpM} and the weighted $L^p$
boundedness of $M$ and~$T$.
\end{proof}

We now apply this to prove the following lemma, which will be
important in the implementation of the Bahouri-Gerard decomposition
for wave maps. Instead of a general Calderon-Zygmund operator, we
restrict ourselves to the subclass of Mikhlin multiplier operators
which are of the form $Tf=(m\hat{f})^{\vee}$ with $m\in
C^3(\R^2\setminus\{0\})$ and with
\begin{equation}\nn
|D^\alpha m(\xi)|\le C(\alpha) |\xi|^{-|\alpha|} \quad\forall\;
\xi\in\R^2\setminus\{0\}
\end{equation}
 for all $|\alpha|\le 3$. For
simplicity, we also limit ourselves to two dimensions.

\begin{lemma}
 \label{lem:geom_dec} Suppose $\{\vphi_n\}_{n=1}^\infty, \{\phi_n\}_{n=1}^\infty$ lie in the unit-ball of~$L^2$. Furthermore, assume that
\[
 \supp(\wht{\vphi_n}),\;  \supp(\wht{\phi_n}) \subset\{\xi\in\R^2\::\: 2^{k_0}\le |\xi|\le 2^{k_1}\}
\]
for arbitrary $k_0<k_1-4$ and let $v_n:= e^{(-\Delta)^{-\frac12}
\vphi_n}$.   Then
\[
 \| P_j (v_n^{-1} T(\phi_n v_n)) \|_{2} \les \min(2^{k_1-j}, 2^{\frac{j-k_0}{3}} )
\]
provided either $j\le k_0$ or $j\ge k_1$.
\end{lemma}
\begin{proof} By Lemma~\ref{lem:bmoAp}, for any $1<p<\infty$ one has $\sup_n \sup_{|t|\le 1} A_p(v_n^t)\le C(p)$. Set $R=2^{k_1}$ and $r=2^{k_0}$.
 If $j>k_1$, then
\begin{align*}
 \| P_j (v_n^{-1} T(\phi_n v_n)) \|_{2} &\les 2^{-j} \| \nabla (v_n^{-1} T(\phi_n v_n)) \|_{2} \\
&\les 2^{-j} (\|\phi_n\|_4 \|\vphi_n\|_4 + \|\nabla \phi_n\|_2) \les
2^{-j}R
\end{align*}
On the other hand, if $j<k_0$, then
\begin{align*}
 \| P_j (v_n^{-1} T(\phi_n v_n)) \|_{2} &\les \int_0^1 \|  P_j (v_n^{-t} [T,(-\Delta)^{-\frac12} \vphi_n] (\phi_n v_n^t)) \|_{2}\, dt \\
 &\les \int_0^1 \|  P_j (v_n^{-t} [T,(-\Delta)^{-\frac12} \vphi_n] (\phi_n v_n^t)) \|_{2}\, dt \\
&\les 2^{\frac{j}{3}}  \int_0^1 \|  P_j (v_n^{-t} [T,(-\Delta)^{-\frac12} \vphi_n] (\phi_n v_n^t)) \|_{\frac32}\, dt \\
&\les 2^{\frac{j}{3}}   \| (-\Delta)^{-\frac12} \vphi_n\|_{6} \|
\phi_n  \|_{2} \les 2^{\frac{j}{3}}  r^{-\frac13}
\end{align*}
In the last line, one interpolates between $\|(-\Delta)^{-\frac12}
\vphi_n\|_{2} \les r^{-1}$ and $\|(-\Delta)^{-\frac12}
\vphi_n\|_{\bmo}\les 1$.
\end{proof}

The following result   allows  us to strip away weights from $T(\phi)$ provided they result from functions with
frequencies which are well-separated from the Fourier support of~$\phi$.  In what follows, we  use the following terminology
from~\cite{BG}: Given a bounded sequence $\uf:=\{f_n\}_{n\ge1}\subset L^2$, and sequence $\ueps:=\{\eps_n\}_{n\ge1}\subset R^+$, we say that
{\em $\uf$ is $\ueps$-oscillatory} iff
\[
 \lim_{R\to\infty} \limsup_{n\to\infty}  \int_{[|\xi|\eps_n\in (0,\infty)\setminus (R^{-1},R)]}  |f_n(\xi)|^2\, d\xi =0
\]
We say that {\em $\uf$ is $\ueps$-singular} iff
\[
  \limsup_{n\to\infty}  \int_{[|\xi|\eps_n\in (a,b)]}   |f_n(\xi)|^2\, d\xi =0
\]
for all $b>a>0$. In what follows, we shall freely use the scale selection algorithm from Section~III.1 from~\cite{BG}, see in particular
Lemma~3.1, Lemma~3.2 part~(iii), and Proposition~3.4 in that section.

\begin{lemma}
  \label{lem:weight_remov} Suppose both $\{\vphi_n\}_{n=1}^\infty\subset
  L^2(\R^2)$ and $\{\phi_n\}_{n=1}^\infty\subset
  L^2(\R^2)$ are $1$-oscillatory,  whereas $\{\psi_n\}_{n=1}^\infty\subset
  L^2(\R^2)$ is $1$-singular. Define
\[
v_n:= \exp((-\Delta)^{-\frac12} \vphi_n),\quad w_n:=
\exp((-\Delta)^{-\frac12} \psi_n)
\]
Then
\begin{equation}   \label{eq:wnweg}
(v_nw_n)^{-1} T(\phi_n \,v_nw_n) = v_n^{-1} T(\phi_n\, v_n) +
o_{L^2}(1)
\end{equation}
as $n\to\infty$.  Moreover, $v_n^{-1} T(\phi_n\, v_n)$ is
$1$-oscillatory\footnote{Note that neither  $(v_nw_n)^{-1} T(\psi_n
\,v_nw_n)$ nor $v_n^{-1} T(\psi_n \,v_n)$ are in general
$1$-singular.}.
\end{lemma}
\begin{proof}
By assumption,
\[
\|\phi_n\|_2+\|\vphi_n\|_2+\|\psi_n\|_2\le A<\infty
\]
for all $n\ge1$.
By Lemma~\ref{lem:bmoAp} one has $v_n\in A_p$ and $v_nw_n\in A_p$ for
all $1<p<\infty$  with $A_p$ constants depending only
on~$A$ and~$p$.
 Now fix $\eps>0$ arbitrarily small. Then
there is $R>1$ so that
\[
\limsup_{n\to\infty} \int\limits_{[|\xi|<R^{-1},\; |\xi|>R]}
|\hat\vphi_n(\xi)|^2\, d\xi <\eps^2
\]
Fix an $R=R(\eps)$ with this property.
Define $\vphi_{1n}:= (\chi_{[R^{-1},R]} \widehat{\vphi_n})^{\vee}$,
$\vphi_{2n}:=\vphi_n-\vphi_{1n}$ and $\phi_{1n}:= (\chi_{[R^{-1},R]}
\widehat{\phi_n})^{\vee}$, $\phi_{2n}:=\phi_n-\phi_{1n}$.
Then
$\|\vphi_{2n}\|_2+\|\phi_{2n}\|_2<\eps$ for large $n$ whence
\begin{align*}
  \|(v_nw_n)^{-1} T(\phi_{2n} \,v_nw_n)\|_2 + \|v_n^{-1} T(\phi_{2n}
  \,v_n)\|_2 \le C(A,T)\eps
\end{align*}
as well as $\|(-\Delta)^{-\frac12}\vphi_{2n}\|_{\bmo}<C\eps$.  Next,
define
\[
v_{jn}:= \exp((-\Delta)^{-\frac12} \vphi_{jn})\qquad j=1,2
\]
By Lemmas~\ref{lem:bmoAp} and~\ref{lem:CoifMeyer},
\begin{align*}
   &\| (v_nw_n)^{-1} T(\phi_{1n} \,v_nw_n) - (v_{1n}w_n)^{-1} T(\phi_{1n}
   \,v_{1n}w_n)\|_2\\
   &\le \int_0^1 \| (w_nv_{1n}v_{2n}^t )^{-1} [\,T,(-\Delta)^{-\frac12} \vphi_{2n}  ](\phi_{1n}
   \,w_n v_{1n}v_{2n}^t)\|_2\, dt\\
   &\le C(A,T) \| (-\Delta)^{-\frac12} \vphi_{2n}\|_{\bmo} \le C(A,T) \eps
\end{align*}
By the same argument,
\begin{align*}
   \| v_n^{-1} T(\phi_{1n} \,v_n) - v_{1n}^{-1} T(\phi_{1n}
   \,v_{1n})\|_2 &\le C(A,T) \eps
\end{align*}
Similarly, set
\[
\psi_{2n}:= (\chi_{[\rho^{-1},\rho]} \widehat{\psi_n})^{\vee},\quad
\psi_{1n}:=\psi_n-\psi_{2n}
\]
where $\rho>1$ will be determined later. By assumption,
$\|\psi_{2n}\|_2\to0$ as $n\to\infty$. In particular,
\[\|(-\Delta)^{-\frac12}\psi_{2n}\|_{\bmo}\to0\] as $n\to\infty$.
Applying Lemma~\ref{lem:CoifMeyer} as before allows one to remove
the weights $w_{2n}$ from~\eqref{eq:wnweg} where
\[
w_{jn}:= \exp((-\Delta)^{-\frac12} \psi_{jn})\qquad j=1,2
\]
Hence, we are reduced to establishing that
\begin{equation}
  \label{eq:reduc}
\| (v_{1n}w_{1n})^{-1} T(\phi_{1n} \,v_{1n}w_{1n})- v_{1n}^{-1}
T(\phi_{1n}\, v_{1n})\|_2 \le C(A,T)\eps
\end{equation}
for sufficiently large $n$.   For ease of notation, we shall now
drop the subscript~$1$ from $\vphi_{1n}$ etc.~with the understanding
that $\widehat{\vphi_n}$ and $\wht{\phi_n}$ are supported on
$[R^{-1},R]$ and that $\widehat{\psi_n}$ is supported off
$[\rho^{-1},\rho]$ where $\rho>1$  is a large number depending
on~$\eps$ to be chosen later. Define
\[
\psi_{n,{\rm low}}:= (\chi_{(0,\rho^{-1}]}
\widehat{\psi_n})^{\vee},\quad \psi_{n,{\rm high}}:=
(\chi_{[\rho,\infty)} \widehat{\psi_n})^{\vee}
\]
and write, correspondingly, $w_n= w_{n,{\rm low}} w_{n,{\rm high}}$.
It is easy to remove $w_{n,{\rm high}}$:
\begin{align}
&\| (v_{n}w_{n})^{-1} T(\phi_{n} \,v_{n}w_{n})- (v_{n}w_{n,{\rm low}})^{-1} T(\phi_{n} \,v_{n}w_{n,{\rm low}})\|_2 \nn \\
&\le \int_0^1 \| (v_{n}w_{n,{\rm low}}w_{n,{\rm high}}^t)^{-1} [\,T,
(-\Delta)^{-\frac12}\psi_{n,{\rm high}}](\phi_{n}
\,v_{n}w_{n,{\rm low}}w_{n,{\rm high}}^t)\|_2 \, dt \nn \\
&\le C(A,T) \| (-\Delta)^{-\frac12}\psi_{n,{\rm high}}\|_{4}\|
\phi_{n}\|_{4}\label{eq:sch12} \\
&\le C(A,T,R) \rho^{-\frac12} \le \eps \nn
\end{align}
provided $\rho$ is sufficiently large.
Here we used that $v_{n}w_{n,{\rm low}}w_{n,{\rm high}}^t$ are $A_2$
weights uniformly in $0\le t\le1$ as well as an interpolation
between $L^2$ and $\bmo$ to pass to the last line. For the final
bound  we need $\rho\gg\eps^{-2}$.

To remove $w_{n,{\rm low}}$ we split
\[
T= P_{<-\lambda}T + P_{-\lambda<\cdot<\lambda}T + P_{>\lambda }T
\]
where $2^{-\lambda} R\ll \eps$ and $P_{<\lambda}$ etc.~denote
Littlewood-Paley projections. Introducing an angular decomposition
into finitely many sectors, we may assume that $|\xi_1|\ge|\xi|/10$
on the support of~$m$. Then for large   $\lambda$, and with
$\mu:=2^{-\lambda}$,
\begin{align*}
  &\|(v_{n}w_{n,{\rm low}})^{-1} P_{>\lambda }T(\phi_{n}
  \,v_{n}w_{n,{\rm low}})\|_2\le C \|(v_{n}w_{n,{\rm low}})^{-1} \partial_1^{-1} P_{>\lambda }T(\partial_1[\phi_{n}
  \,v_{n}w_{n,{\rm low}}])\|_2\\
  &\le C(A)\mu \big(\|\partial_1 (\phi_n v_n)\|_2 + \|\phi_n
  \partial_1 (-\Delta)^{-\frac12}\psi_{n,{\rm low}}\|_2\big) \le \eps
\end{align*}
For the small frequencies $P_{<-\lambda}T$ we first recall the
following standard fact: with $\psi$ a suitable Schwarz function,
\begin{align}
  (P_{<-\lambda} (fg) - gP_{<-\lambda}f)(x) &= -\sum_{j=1}^2 \int_0^1 \int_{\R^2}
  \mu^2\psi(\mu y) y_j f(x-y) \partial_j g(x-sy)\,dy ds \nn \\
  &= \sum_{j=1}^2 L_{j,\lambda} (\mu^{-1}f,\partial_j g) \label{eq:Lj}
\end{align}
where $L_{j,\lambda}$ in the final line denotes a multi-linear expression of
the form
\begin{equation}
\label{eq:muav} L(f,g)(x) = \int_{\R^4} f(x-u) g(x-v)\,\nu(du,dv)
\end{equation}
with a measure $\nu$ of mass bounded by some constant (in this case uniformly in all parameters). Using this
notation, one has (since $\|v_n^{-1}\|_\infty \le C(A,R)$)
\begin{align*}
  \|(v_{n}w_{n,{\rm low}})^{-1} P_{<-\lambda}T(\phi_{n}
  \,v_{n}w_{n,{\rm low}})\|_2
  &\le C(A) \|w_{n,{\rm low}}^{-1} T(P_{<-\lambda}(\phi_{n}
  \,v_{n})w_{n,{\rm low}})\|_2 \\
  &\quad +C(A)\mu^{-1} \sum_{j=1}^2 \| w_{n,{\rm low}}^{-1} T( L_{j,\lambda}(\phi_{n}
  v_{n},w_{n,{\rm low}}\partial_j(-\Delta)^{-\frac12}\psi_{n,{\rm low}} ))\|_2\\
  &=: I_n + II_n
\end{align*}
To bound $I_n$, note that since we may take  $\mu\le R^{-1}$, one
has $P_{<-\lambda}(\phi_{n}
  v_{n})=P_{<-\lambda}(\phi_{n}
  (v_{n}-1))$. Hence, by the boundedness of $T$ relative to the
  weight $w_{n,{\rm low}}$ and Bernstein's inequality,
\begin{align*}
  \| I_n\|_2 &\le \| P_{<-\lambda}(\phi_{n}
  v_{n}) \|_2 = \| P_{<-\lambda}(\phi_{n}
  (v_{n}-1)) \|_2 \le C \mu \| P_{<-\lambda}(\phi_{n}
  (v_{n}-1) \|_1 \\
&\le C \mu \| \phi_{n} \|_2
  \| v_{n}-1 \|_2 \le  \eps
\end{align*}
for  $\lambda\gg 1$. Next, in view of~\eqref{eq:Lj}, $II_n$ is
bounded by (using $\|v_n\|_\infty \le C(A,R)$)
\begin{align}
\nn
   &C(A)\mu^{-1} \int_{\R^2}\int_0^1 \mu^2 |\psi(\mu y)|\, |\mu y|
    \big\| \phi_n(\cdot-y)
   w_{n,{\rm low}}(\cdot- sy)w_{n,{\rm low}}^{-1}(\cdot)
     \nabla(-\Delta)^{-\frac12}\psi_{n,{\rm low}}(\cdot -sy)\big\|_{L^2}\, dy
   ds \\ \nn
&\le C(A) \mu^{-1} \|\phi_{n}
  \|_2\|\nabla(-\Delta)^{-\frac12}\psi_{n,{\rm low}}\|_\infty \int_{\R^2}\int_0^1 \mu^2 |\psi(\mu y)|\,|\mu y|
  \big\|
   w_{n,{\rm low}}(\cdot- sy)w_{n,{\rm low}}^{-1}(\cdot)
     \big\|_{L^\infty}\, dy
   ds \\
  &\le C(A)\mu^{-1} \rho^{-1}  \int_{\R^2} \mu^2 (1+\mu |y|)^{-N}\,
  (1+|y|\mu)^{k(A)}\, dy \label{eq:k(A)}
  \le \eps
\end{align}
for some constant $k(A)>0$  provided we
  choose $\rho$ such that $\mu^{-1}\eps^{-1}\ll \rho$. To pass to
  the bound in~\eqref{eq:k(A)}, assume first that
  $|y|\mu\le1$. Then with
  $h_n:=(-\Delta)^{-\frac12}\psi_{n,{\rm low}}$ so that $w_{n,{\rm low}}=e^{h_n}$,
\begin{align*}
  w_{n,{\rm low}}(x- y)w_{n,{\rm low}}^{-1}(x) &\le   \exp\big( |y|
  \|\nabla h_n\|_\infty \big) \le
  \exp\big(C\mu^{-1}\rho^{-1}\big)\le e^{C\eps}
  \le 2
\end{align*}
where we used that \[\|\nabla
  h_n \|_\infty\le \|\nabla (-\Delta)^{-\frac12}\psi_{n,{\rm low}}\|_\infty \le  C\rho^{-1}\]
This implies that on scales $\le \mu^{-1}$, the weight
$w_{n,{\rm low}}$ is essentially constant (up to multiplicative
constants). Next, observe that for all cubes
\[
|(h_n)_Q- (h_n)_{2^\ell Q}|\le C\|h_n\|_{\bmo}\ell\le C(A) \ell
\qquad \forall\;\ell\ge0
\]
Hence, partitioning $\R^2$ into cubes of side-length $\mu^{-1}$
one obtains that
\begin{equation}\label{eq:hdiff}
|h_n(y)-h_n(y')|\le C(A)\log(2+|y-y'|\mu)
\end{equation}
whence
\[
\sup_x w_{n,{\rm low}}(x- y)w_{n,{\rm low}}^{-1}(x)\le
C(A)(1+|y|\mu)^{k(A)}
\]
as claimed.

\noindent Note that the previous estimates on $P_{<-\lambda} T$ and
$P_{>\lambda}T$ also prove that
\[
\| v_n^{-1} P_{<-\lambda} T(\phi_n v_n)\|_2 + \| v_n^{-1}
P_{>\lambda} T(\phi_n v_n)\|_2 \le \eps
\]
Therefore, it remains to  prove that
  \[
\| (v_{n}w_{n,{\rm low}})^{-1} T_\lambda (\phi_{n} \,v_{n}w_{n,{\rm
low}}) - v_{n}^{-1} T_\lambda (\phi_{n} \,v_{n}) \|_2 \le \eps
  \]
  where
\begin{equation}
  \label{eq:Tlambda}
  T_\lambda:= P_{-\lambda<\cdot<\lambda}T
\end{equation}
is the operator on intermediate frequencies.  Since
$Tf=(m\hat{f})^{\vee}$ with $m\in C^3(\R^2\setminus\{0\})$, we
conclude that
\[
 P_{\lambda<\cdot<\lambda^{-1}}T f(x) = \int K_\lambda (x-y) f(y)\,
 dy
\]
with $|K_\lambda(x)|\le C(\lambda)(1+|x|)^{-3}$. Now, with
$h_n=(-\Delta)^{-\frac12}\psi_{n,{\rm low}}$ as above, and $M$ denoting the Hardy-Littlewood maximal operator,
\begin{align*}
 &\| (v_{n}w_{n,{\rm low}})^{-1} T_\lambda (\phi_{n} \,v_{n}w_{n,{\rm low}}) -
v_{n}^{-1} T_\lambda (\phi_{n} \,v_{n}) \|_2 \\
&\le \int_0^1 \| (v_{n}w_{n,{\rm low}}^t)^{-1} [T_\lambda, h_n ]
(\phi_{n}
\,v_{n}w_{n,{\rm low}}^t) \|_2 \, dt \\
&\le C(A,\lambda)\rho^{-\frac14} \int_0^1 \| (v_{n}w_{n,{\rm
low}}^t)^{-1} M
 (\phi_{n}
\,v_{n}w_{n,{\rm low}}^t) \|_2 \, dt \le C(A,\lambda)\rho^{-\frac14}
\end{align*}
Here we used that the kernel of $[T_\lambda, h_n]$ is of the form $
K_\lambda(x,y)(h_n(x)-h_n(y))$ and satisfies the bounds,
cf.~\eqref{eq:hdiff},
\[
|K_\lambda(x,y)(h_n(x)-h_n(y))|\le C(A,\lambda)\min(\rho^{-1}|x-y|,
|x-y|^{-3}\log(2+|x-y|))
\]
whence
\[
|[T_\lambda,h_n]f(x)|\le C(A,\lambda)\rho^{-\frac14} Mf(x)
\]
Taking $\rho$ sufficiently large (depending on $\eps$, $R$, and~$A$)
finishes the proof of~\eqref{eq:wnweg}. Lemma~\ref{lem:geom_dec} now
implies that $v_n^{-1} T(\phi_n v_n)$ is $1$-oscillatory.
\end{proof}

The following statement will be an essential technical tool for the
Bahouri-Gerard method in the context of wave maps into hyperbolic
space. As before, $T$ is a Mikhlin multiplier operator.

\begin{cor}\label{cor:Jlem}
Let $\{f_n\}_{n=1}^\infty\subset L^2(\R^2)$ satisfy $\sup_{n\ge1}
\|f_n\|_2\le A<\infty$ and define
$y_n=\exp\big((-\Delta)^{-\frac12}f_n\big)$. Let
$\Lambda_j:=\big\{\lambda_{n,j}\big\}_{n=1}^\infty$ be sequences of
positive
  numbers for each $1\le j\le J$ with the property that
  \begin{equation}\label{eq:orth}
\lim_{n\to\infty} \Big\{\frac{\lambda_{n,j}}{\lambda_{n,j'}} +
\frac{\lambda_{n,j'}}{\lambda_{n,j}}\Big\} \to\infty
  \end{equation}
for any $1\le j\ne j'\le J$. Assume further that
\[
f_n = \sum_{j=1}^J \vphi_{n,j} + \omega_n
\]
where $\{\vphi_{n,j}\}_{n=1}^\infty\subset L^2(\R^2)$ is
$\Lambda_j$-oscillatory for each $1\le j\le J$,
$\{\omega_{n}\}_{n=1}^\infty$ is $\Lambda_j$-singular for every
$1\le j\le J$, and $
\sup_{n\ge1}\|\omega_n\|_{\dot{B}^0_{2,\infty}}<\delta $.

\noindent Then $\big\{y_n^{-1}\, T(\vphi_{n,j}\,
y_n)\big\}_{n=1}^\infty$ is $\Lambda_j$-oscillatory,
$\big\{y_n^{-1}\, T(\omega_n\, y_n)\big\}_{n=1}^\infty$ is
$\Lambda_j$-singular for each $1\le j\le J$, and
\begin{equation}\label{eq:Besov}
\limsup_{n\to\infty}\|y_n^{-1}\, T( \omega_n
y_n)\|_{\dot{B}^0_{2,\infty}}< C(A,T)\delta
\end{equation}
where the constant $C(A,T)$ only depends on $A$ and $T$.
\end{cor}
\begin{proof}
Define
\[
v_{n,j}:= \exp((-\Delta)^{-\frac12} \vphi_{n,j}),\qquad w_{n}:=
\exp((-\Delta)^{-\frac12} \omega_{n})
\]
so that $y_n:= w_n \prod_{j=1}^J v_{n,j}$.
  By Lemma~\ref{lem:weight_remov}, both $\big\{v_{n,j}^{-1}\, T(\vphi_{n,j}
v_{n,j})\big\}_{n=1}^\infty$ and $\big\{y_n^{-1}\, T(\vphi_{n,j} \,
y_n)\big\}_{n=1}^\infty$  are $\Lambda_j$-oscillatory. Now suppose
$\{\psi_n\}_{n=1}^\infty$ is an arbitrary $\Lambda_j$-oscillatory
sequence where $1\le j\le J$ is fixed. Then $\widetilde\omega_n :=
y_n^{-1}\, T(\omega_n y_n)$ satisfies
\[
\langle \widetilde\omega_n,\psi_n\rangle = \langle \omega_n,
  y_n T^*(\psi_n y_n^{-1})\rangle
\]
By Lemma~\ref{lem:weight_remov}, $\{y_n T^*(\psi_n
y_n^{-1})\}_{n=1}^\infty$ is $\Lambda_j$-oscillatory whence
\[
\lim_{n\to\infty}\langle \widetilde\omega_n,\psi_n\rangle =0
\]
Therefore, $\big\{\widetilde\omega_n\big\}_{n=1}^\infty$ is
$\Lambda_j$-singular for each $1\le j\le J$.

For the proof of \eqref{eq:Besov}, we first note that passing to a
subsequence if necessary, \eqref{eq:orth} implies that we may assume
that
\[
\lambda_{n,1}>\lambda_{n,2}>\ldots>\lambda_{n,J}
\]
for all large~$n$ whence  for any $1\le j\le J-1$
\begin{equation}
  \label{eq:lamb_rat} \frac{\lambda_{n,j+1}}{\lambda_{n,j}}\to\infty
\end{equation}
as $n\to\infty$. We also note that
\[
\sum_{j=1}^J \|\vphi_{n,j}\|_2^2 + \|w_n\|_2^2 \le A^2 +o(1)
\]
as $n\to\infty$. Now we let $m\ge10$ and $K\ge10$ be integers (to be
determined later) and define
\begin{align*}
  \wt \vphi_{n,j} &:=
  \vphi_{n,j}\, \chi_{[2^{-m}\le|\xi|\lambda_{n,j}\le
  2^m ]} \\
  \wt y_n &:= \prod_{j=1}^J \exp\big((-\Delta)^{-\frac12}
  \wt\vphi_{n,j}\big)\\
  \wt w_n &:= w_n \, \chi_{\R^2\setminus\bigcup_{j=1}^J
  [2^{-Km}\le |\xi|\lambda_{n,j}\le 2^{Km} ]}
\end{align*}
As in the proof of Lemma~\ref{lem:weight_remov},
\begin{equation}\label{eq:delta_m}
\limsup_{n\to\infty} \| y_n^{-1} T(w_n y_n) - \tilde y_n^{-1}
T(\tilde w_n\, \tilde y_n)\|_2 <  C(A,T)\,\delta
\end{equation}
provided $m$ is chosen large enough and irrespective of the choice
of $K\ge1$. We will now fix $m$ so that~\eqref{eq:delta_m} holds. It
therefore suffices to show that
\begin{equation}\label{eq:reduc_main}
\limsup_{n\to\infty}\; \sup_{j\in\Z} \|P_j \,[\tilde y_n^{-1}
T(\tilde w_n\, \tilde y_n)]\|_2 \le C(A,T)\, \delta
\end{equation}
provided $K$ is chosen sufficiently large. The idea
behind~\eqref{eq:reduc_main} is that $\wt w_n$ behaves like a
lacunary series, i.e., each $\wt w_n$ is the sum of functions whose
 Fourier supports consist of disjoint blocks which are
 very strongly separated. In addition, the
$\wt\vphi_{n,j}$ are Fourier supported on intervals which are  well
separated from the Fourier support of~$\wt w_n$. It will turn out
that for each $j$ -- up to negligible errors as $K\to\infty$ -- only
one block of frequencies from~$\wt w_n$ (namely the one
containing~$2^j$) contributes to $ P_j \,[\tilde y_n^{-1} T(\tilde
w_n\, \tilde y_n)] $ and, moreover, only those $\wt\vphi_{n,j}$ with
frequencies much smaller than $2^j$ matter. In this way, we can then
essentially pass $P_j$ onto~$\tilde w_n$.

\noindent To establish~\eqref{eq:reduc_main}, we introduce some more
notation: set \[\psi_n:= \sum_{j=1}^J (-\Delta)^{-\frac12}
  \wt\vphi_{n,j}\] and define $[T,\psi_n]^{(s)}$ iteratively via
  \[[T,\psi_n]^{(1)}:=[T,\psi_n],\qquad [T,\psi_n]^{(s+1)}=[[T,\psi_n]^{(s)},\psi_n]\] Then
\begin{align}
  \tilde y_n^{-1} T(\tilde w_n \tilde y_n) &= \sum_{\ell=0}^s \frac{1}{\ell!}
  [T,\psi_n]^{(\ell)}\,\tilde w_n + \label{eq:main_term} \\
   &\quad +\frac{1}{s!} \int_0^1 (1-t)^s\,
  \tilde y_n^{-t} [T,\psi_n]^{(s+1)} (\tilde w_n\,  \tilde y_n^t)\, dt \label{eq:comm_iter}
\end{align}
Denote the remainder in~\eqref{eq:comm_iter} by~$R_{n,s}$.  To bound
it in~$L^2$, note that $\|\psi_n\|_\infty\le CmJA$ with some
absolute constant~$C$. Therefore, placing $\psi_n$, $\tilde y_n$,
and $\tilde y_n^{-1}$ in $L^\infty$ yields for all $n\ge1$
\begin{align*}
\|R_{n,s}\|_2 &\le \frac{e^{CmJA}}{(s+1)!} (CmJA)^{s+1} \|\tilde
w_n\|_2  \le e^{CmJA} \Big(\frac{CmJA}{s}\Big)^{s}=:\gamma_s
\end{align*}
which clearly goes to zero as $s\to\infty$. In particular,
$\gamma_s\le\delta$ for large~$s$. We now turn to the details of the
analysis of the main terms in~\eqref{eq:main_term}. First, one has $
\tilde w_n = \sum_{j=0}^J \wt w_{n,j} $ where
\begin{align*}
  \wt w_{n,0} &:= \tilde w_n \,\chi_{[|\xi|\,\lambda_{n,1}\le
  2^{-Km}]} \\
\wt w_{n,j} &:= \tilde w_n \,\chi_{[2^{Km}\le |\xi| \lambda_{n,j}\le
  2^{-Km}]} \quad \forall \;1\le j\le J-1 \\
  \wt w_{n,J} &:= \tilde w_n \, \chi_{[2^{Km} \le
|\xi|\lambda_{n,J}]}
\end{align*}
with $n$ large. Then
\begin{align}
\sum_{\ell=0}^s \frac{1}{\ell!}
  [T,\psi_n]^{(\ell)}\,\tilde w_n  &= \sum_{j=0}^J \sum_{\ell=0}^s \frac{1}{\ell!}
  [T,\psi_n]^{(\ell)}\,\wt w_{n,j}
= \sum_{j=0}^J \sum_{\ell=0}^s \frac{1}{\ell!}
  [T,\psi_{n,j}^{(-)}+ \psi_{n,j}^{(+)}]^{(\ell)}\,\wt w_{n,j}
  \label{eq:comm_decomp}
\end{align}
where we have set, for each $0\le j\le J$,
\[
\psi_{n,j}^{(-)} := \sum_{1\le k\le j}  (-\Delta)^{-\frac12}
\wt\vphi_{n,k},\qquad \psi_{n,j}^{(+)} := \sum_{j< k\le J}
(-\Delta)^{-\frac12} \wt\vphi_{n,k}
\]
We shall now show that for a given $\wt w_{n,j}$ only the small
frequency part of $\psi_{n,j}$, i.e., $\psi_{n,j}^{(-)}$,
contributes significantly to the commutators
in~\eqref{eq:comm_decomp} (at least for very large~$K$).
 To this end write
\begin{align} [T,\psi_{n,j}^{(-)}+
\psi_{n,j}^{(+)}]^{(\ell)} &= [T,\psi_{n,j}^{(-)}]^{(\ell)} +
\sum_\eps [[\cdots[T,\psi_{n,j}^{(\eps_1)}],
\;\psi_{n,j}^{(\eps_2)}], \; \ldots, \psi_{n,j}^{(\eps_\ell)}] \nn \\
&= \calK_{n,j,\ell}+\calR_{n,j,\ell} \label{eq:comm_ell}
\end{align}
 where the sum here runs over $\ell$-fold commutators
with each $\eps_k=\pm$, the choice $\eps_k=-$ for all $1\le
k\le\ell$ being excluded (as it is represented by the first -- and
main -- term on the right-hand side). Next, observe that for each
\begin{equation}
  \label{eq:s_def}
  1\le\ell\le s:=\frac{1}{100} 2^{(K-1)m}
\end{equation} one has, for each $1\le
j\le J$ and every $k\in\Z$,
\[
P_k\big( \calK_{n,j,\ell} \wt w_{n,j}\big) = P_k \calK_{n,j,\ell}
P_{k-2<\cdot<k+2}\,\wt w_{n,j}
\]
In fact, this vanishes unless $2^{Km-2} \lambda_{n,j}^{-1}\le 2^k\le
  2^{-Km+2}\,\lambda_{n,j+1}^{-1}$. Writing
  \begin{align}
    P_k \sum_{\ell=0}^s \frac{1}{\ell!}
  [T,\psi_n]^{(\ell)}\,\tilde w_n &= P_k \sum_{j=0}^J\sum_{\ell=0}^s \frac{1}{\ell!}
  \calK_{n,j,\ell}
P_{k-2<\cdot<k+2}\,\wt w_{n,j} \label{eq:first_sum} \\
&\quad + P_k \sum_{j=0}^J\sum_{\ell=0}^s \frac{1}{\ell!}
  \calR_{n,j,\ell}
P_{k-2<\cdot<k+2}\,\wt w_{n,j} \label{eq:second_sum}
  \end{align}
it follows from~\eqref{eq:lamb_rat} that for all sufficiently
large~$n\ge n_0$ depending on~$K,m$, at most one term
in~\eqref{eq:first_sum} can be nonzero for any choice of~$k\in \Z$.
Applying the decomposition~\eqref{eq:main_term}
and~\eqref{eq:comm_iter} with $\psi_{n,j}^{(-)}$ instead of~$\psi_n$
to~\eqref{eq:first_sum} yields
\begin{align*}
&\sup_{k\in\Z}\Big\|P_k  \sum_{j=0}^J\sum_{\ell=0}^s \frac{1}{\ell!}
  \calK_{n,j,\ell}
P_{k-2<\cdot<k+2}\,\wt w_{n,j} \Big\|_2 \le \sup_{k\in\Z}\sup_{0\le
j\le J}\Big\|\sum_{\ell=0}^s \frac{1}{\ell!}
  \calK_{n,j,\ell}
P_{k-2<\cdot<k+2}\,\wt w_{n,j} \Big\|_2 \\
&\le \sup_{k\in\Z}\sup_{0\le j\le J} \Big\| e^{-\psi_{j,n}^{(-)}}
T\big( e^{\psi_{j,n}^{(-)}} P_{k-2<\cdot<k+2}\,\wt w_{n,j} \big)
\Big\|_2 +\gamma_s\\
&\le C(A,T) \|w_n\|_{\dot{B}^0_{2,\infty}} + \gamma_s \le
C(A,T)\delta
\end{align*}
To pass to the final bound, we note that $\gamma_s\le \delta$
provided $K$ is chosen sufficiently large. We also used that the
weights $e^{\psi_{j,n}^{(-)}}\in A_2$ with $A_2$ constant $\le CA$
uniformly in $j,n$, cf.~Lemma~\ref{lem:CoifMeyer}.  As
for~\eqref{eq:second_sum}, we make the following crude estimate for
the $\ell$-fold commutator as in~\eqref{eq:comm_ell}
\[
\big\|[[\cdots[T,\psi_{n,j}^{(\eps_1)}], \;\psi_{n,j}^{(\eps_2)}],
\; \ldots, \psi_{n,j}^{(\eps_\ell)}] \wt w_{n,j}\|_2 \le C(T)
(CmJA)^{\ell-1} \| \psi_{n,j}^{(+)} \|_4 \|\wt w_{n,j}\|_4
\]
It arises by placing one $\psi_{n,j}^{(+)}$ in~$L^4$, all other
$\psi_{n,j}^{(\eps_i)}$ in~$L^\infty$, and $\tilde w_{n,j}$
in~$L^4$. By Bernstein's inequality,
\[
\|\wt w_{n,j}\|_4 \le C (2^{-Km}\lambda_{n,j+1})^{-\frac12}\|\wt
w_{n,j}\|_2 \le CA 2^{-Km/2}\lambda_{n,j+1}^{-\frac12}
\]
whereas by interpolation between the $L^2$ and $\bmo$ bounds,
\[
\| \psi_{n,j}^{(+)} \|_4 \le CA 2^{m/2}\lambda_{n,j+1}^{\frac12}
\]
whence
\[
\big\|[[\cdots[T,\psi_{n,j}^{(\eps_1)}], \;\psi_{n,j}^{(\eps_2)}],
\; \ldots, \psi_{n,j}^{(\eps_\ell)}] \wt w_{n,j}\|_2 \le C(T,A)
(CmJA)^{\ell-1} 2^{(1-K)m/2}
\]
Hence, the error resulting from~\eqref{eq:second_sum} can be made as
small as we wish by taking $K$ large and we are done.
\end{proof}

In what follows, we call sequences $\Lambda_j\subset\R^+$ as in
Lemma~\ref{cor:Jlem} pairwise {\em orthogonal} iff they
satisfy~\eqref{eq:orth}. The following auxiliary
Lemma~\ref{lem:elimhigh} strengthens the result of
Lemma~\ref{lem:weight_remov} by replacing $L^2$ with $\dot
B^{0}_{2,1}$, but under slightly different conditions.  As before,
$T$ is a Mikhlin operator.

\begin{lemma}
 \label{lem:elimhigh}
Suppose $\{\vphi_n\}_{n=1}^\infty, \{\phi_n\}_{n=1}^\infty,
\{\psi_n\}_{n=1}^\infty$ lie in the unit-ball of~$L^2$. Furthermore,
assume that
\[
 \supp(\wht{\vphi_n}),\;  \supp(\wht{\phi_n}) \subset\{|\xi|\le 1\},\quad \supp(\wht{\psi_n}) \subset\{|\xi|\ge1\}
\]
Define
\[
v_n:= \exp((-\Delta)^{-\frac12} \vphi_n),\quad w_n:=
\exp((-\Delta)^{-\frac12} \psi_n)
\]
Then given $\eps>0$ there exists $\delta>0$ such that
\begin{align}   \label{eq:wnBeselim}
\big\| (v_nw_n)^{-1} T(\phi_n \,v_nw_n) - v_n^{-1} T(\phi_n\, v_n) \big\|_{\dot B^0_{2,1}} &< \eps \\
\big\|\nabla^{-1}\big[ (v_nw_n)^{-1} T(\phi_n \,v_nw_n) - v_n^{-1}
T(\phi_n\, v_n) \big] \big\|_{\infty} &<\eps  \label{eq:Linfty_err}
\end{align}
for all sufficiently large~$n$
provided
\begin{equation}
\label{eq:Bcondi}
\limsup_{n\to\infty} \|P_{<k_0} \psi_n\|_{\dot B^0_{2,\infty}}<\delta
\end{equation}
 where $k_0=k_0(T,\eps)$ is some positive integer.
\end{lemma}
\begin{proof} Since \eqref{eq:wnBeselim} implies \eqref{eq:Linfty_err} it suffices to prove the former.
As before,
\begin{equation}
 \label{eq:vnwn_diff}
 (v_nw_n)^{-1} T(\phi_n \,v_nw_n) - v_n^{-1} T(\phi_n\, v_n) = \int_0^1 (v_nw_n^t)^{-1}
[T, (-\Delta)^{-\frac12} \psi_n](\phi_n \,v_nw_n^t)\, dt
\end{equation}
We now estimate the $L^2$-norm of this expression localized to
frequency~$2^j$. First, we consider the case $j\ge0$. Then, with
$y_{n,t}:= v_nw_n^t$, and using Bernstein's inequality, one has the
bound
\begin{align}
& \| P_j \big( y_{n,t}^{-1} [T, (-\Delta)^{-\frac12} \psi_n](\phi_n
\,y_{n,t}) \big) \|_{2} \les 2^{-\frac23 j} \| P_j \nabla \big(
y_{n,t}^{-1}
[T, (-\Delta)^{-\frac12} \psi_n](\phi_n \,y_{n,t}) \big) \|_{\frac32} \nn \\
&\les 2^{-\frac23 j} \|  \nabla \big( y_{n,t}^{-1} \big) [T,
(-\Delta)^{-\frac12} \psi_n](\phi_n \,y_{n,t})  \|_{\frac32} +
2^{-\frac23 j} \|    y_{n,t}^{-1}
[T, (-\Delta)^{-\frac12} \psi_n] \nabla (\phi_n \,y_{n,t})  \|_{\frac32} \label{eq:sch32} \\
&\quad + 2^{-\frac23 j} \|    y_{n,t}^{-1} [T, \nabla
(-\Delta)^{-\frac12} \psi_n]  (\phi_n \,y_{n,t})  \|_{\frac32} \nn
\end{align}
Since uniformly in $0\le t\le1$,
\[
 y_{n,t}\nabla  y_{n,t}^{-1}=-(\nabla (-\Delta)^{-\frac12} \vphi_n+t\nabla (-\Delta)^{-\frac12} \psi_n)=O_{L^2}(1)
\]
we can further estimate
\begin{align*}
 \|  \nabla \big( y_{n,t}^{-1} \big) [T, (-\Delta)^{-\frac12}
\psi_n](\phi_n \,y_{n,t})  \|_{\frac32}
& \les   \|   y_{n,t}^{-1}
[T, (-\Delta)^{-\frac12} \psi_n](\phi_n \,y_{n,t})  \|_{6} \\
&\les   \|  (-\Delta)^{-\frac12} \psi_n \|_{12}
\|\phi_n   \|_{12} \les 1
\end{align*}
To pass to the final bound, the term involving $\psi_n$ is estimated
via an $L^2$-$\bmo$ interpolation, whereas the $\phi_n$ term is
controlled by Bernstein's inequality.  The other two terms on the
right-hand side of~\eqref{eq:sch32} are estimated similarly. As for
the case $j\le0$, Bernstein's inequality yields
\begin{align*}
 \| P_j \big( y_{n,t}^{-1}
[T, (-\Delta)^{-\frac12} \psi_n](\phi_n \,y_{n,t}) \big) \|_{2}
&\les 2^{\frac{j}{3}} \| P_j \big( y_{n,t}^{-1}
[T, (-\Delta)^{-\frac12} \psi_n](\phi_n \,y_{n,t}) \big) \|_{\frac32} \\
&\les 2^{\frac{j}{3}} \|(-\Delta)^{-\frac12} \psi_n\|_6 \|\phi_n\|_2
\les 2^{\frac{j}{3}}
\end{align*}
To obtain~\eqref{eq:wnBeselim}, it suffices to show that for every
$\eps>0$
\[
 \big\| (v_nw_n)^{-1} T(\phi_n \,v_nw_n) - v_n^{-1} T(\phi_n\, v_n) \big\|_{2} < \eps^2
\]
for large~$n$. Indeed, combining this bound with the preceding then implies
\[
 \big\| (v_nw_n)^{-1} T(\phi_n \,v_nw_n) - v_n^{-1} T(\phi_n\, v_n) \big\|_{\dot B^0_{2,1}} \les \eps^2\,  \log\eps
\]
which is more than enough.
To this end, fix a large enough $a$, and let
$w_n=w_{n,{\rm low}} w_{n,{\rm high}}$ where $w_{n,{\rm low}}$
corresponds to~$\calF^{-1}[\chi_{[|\xi|\le a]}\wht{\psi_n}] $ and
$w_{n,{\rm high}}$
 to~$\calF^{-1}[\chi_{[|\xi|\ge a]}\wht{\psi_n}]$ (with sharp cut-offs). By~\eqref{eq:psi_evac} and Lemma~\ref{lem:CoifMeyer},
\[
 \big\| (v_nw_n)^{-1} T(\phi_n \,v_nw_n) - (v_n w_{n,{\rm high}})^{-1} T(\phi_n\, v_n w_{n,{\rm high}}) \big\|_{2} \le C(T) \| \calF^{-1}[\chi_{[|\xi|\le a]}\wht{\psi_n}]  \|_{\dot B^0_{2,\infty}}
\]
whereas
\[
 \sup_{n} \big\| v_n^{-1} T(\phi_n \,v_n) - (v_n w_{n,{\rm high}})^{-1} T(\phi_n\, v_n w_{n,{\rm high}}) \big\|_{2} \le C(T)\,  a^{-\frac12}
\]
by the same argument as in~\eqref{eq:sch12}. Choosing $a$ so that
this final bound is $<\eps$ defines both $k_0(T,\eps)$ and~$\delta$.
\end{proof}

Clearly, one has the following limiting statement.

\begin{cor}
\label{cor:elimhighlimit}
Suppose $\{\vphi_n\}_{n=1}^\infty, \{\phi_n\}_{n=1}^\infty,
\{\psi_n\}_{n=1}^\infty$ lie in the unit-ball of~$L^2$. Furthermore,
assume that
\[
 \supp(\wht{\vphi_n}),\;  \supp(\wht{\phi_n}) \subset\{|\xi|\le1\}
\]
and
\begin{equation}\label{eq:psi_evac}
 \supp(\wht{\psi_n}) \subset\{|\xi|\ge1\},\qquad \lim_{n\to\infty} \int_{|\xi|\le a} |\wht{\psi_n}(\xi)|^2\, d\xi  =0
\end{equation}
for each $a\ge1$. Define
\[
v_n:= \exp((-\Delta)^{-\frac12} \vphi_n),\quad w_n:=
\exp((-\Delta)^{-\frac12} \psi_n)
\]
Then
\begin{align}   \nn
\big\| (v_nw_n)^{-1} T(\phi_n \,v_nw_n) - v_n^{-1} T(\phi_n\, v_n) \big\|_{\dot B^0_{2,1}} &\to 0 \\
\big\|\nabla^{-1}\big[ (v_nw_n)^{-1} T(\phi_n \,v_nw_n) - v_n^{-1}
T(\phi_n\, v_n) \big] \big\|_{\infty} &\to 0 \nn
\end{align}
as $n\to\infty$.
\end{cor}

\section{The Bahouri-Gerard concentration compactness method}
\label{sec:BG}

In this section, we execute the scheme that was sketched in the
introduction. We shall follow the five individual steps which we
outlined there.

\subsection{The precise setup for the Bahouri-Gerard method}

As far as the concentration compactness method is concerned, our
goal is to demonstrate the following main result.

\begin{prop}\label{prop:TheBound}
Let ${\bf{u}}=({\bfx}, {\bfy}): (-T_{0},
T_{1})\times\R^{2}\to{\mathbb{H}}^{2}$ be a Schwartz
class wave map. Then denoting its energy
\[
\sum_{\alpha=0,1,2}\Big(\Big\|\frac{\partial_{\alpha}{\bfx}}{{\bfy}}\Big\|_{L_{x}^{2}}^{2}+
\Big\|\frac{\partial_{\alpha}{\bfy}}{{\bfy}}\Big\|_{L_{x}^{2}}^{2}\Big)=E<\infty,
\]
there is a an increasing function $C(E): \R^+\to \R^+$ with the
property
\[
\|\psi\|_{S((-T_{0}, T_{1})\times\R^{2})}\leq C(E)
\]
 \end{prop}

We refer to the derivative components of ${\bf{u}}$ with respect to
the standard frame $({\bfy}\partial_{{\bfx}},\,{\bfy}
\partial_{{\bfy}})$ as $\phi^{i}_{\alpha}$, $i=1,2$, $\alpha=0,1,2$.   We also use the complex notation
$\phi_{\alpha}:=\phi^{1}_{\alpha}+i\phi^{2}_{\alpha}$.  We shall
refer to a wave map as {\it{admissible}}, provided its derivative
components at time $t=0$, $\phi^{i}_{\alpha}(0,\cdot)$ lie in the
Schwartz class.  Finally, for wave maps of Schwartz class as before,
we denote the Coulomb components by
\[
\psi_{\alpha}:=\phi_{\alpha}\,
e^{-i\sum_{k=1,2}\triangle^{-1}\partial_{k}\phi^{1}_{k}}
\]
The energy is then given by
\[
E({\bf{u}})=\sum_{\alpha=0,1,2}\|\phi_{\alpha}\|_{L_{x}^{2}}^{2}=\sum_{\alpha=0,1,2}\|\psi_{\alpha}\|_{L_{x}^{2}}^{2}
\]
To prove Proposition~\ref{prop:TheBound}, we proceed by
contradiction, assuming that the set of energy levels $E$ for which
it fails is nonempty. Then it has an infimum $\Ecrit>0$ by the small
energy result. We can then find a sequence of wave maps
${\bf{u}}^{n}=({\bfx}^{n}, {\bfy}^{n}): (-T_{0}^{n},
T_{1}^{n})\times\R^{2}\to{\mathbb{H}}^{2}$ with the properties
\begin{itemize}
 \item $\lim_{n\to \infty}E({\bf{u}}^{n})=\Ecrit$ (these energies approach $\Ecrit$ from above)
\item $\lim_{n\to \infty}\|\psi^{n}\|_{S((-T_{0}^{n}, T_{1}^{n})\times\R^{2})}=\infty$.
\end{itemize}
We call such a sequence of wave maps {\it{essentially singular}}. It
is now our goal to apply the Bahouri-Gerard method to the derivative
components of a sequence of essentially singular data
$\phi^{n}_{\alpha}(0,\cdot)$.

\subsection{Step 1: frequency decomposition of initial data}
\label{subsec:BGstep1}

We consider wave maps $\bfu: \R^{2+1}\to \Hyp^{2}$, with Schwartz
initial data. Here $\Hyp^2$ stands for  two-dimensional hyperbolic
space which we identify with the upper half-plane. More precisely,
introducing coordinates $(\bfx, \bfy)$ on $\Hyp^{2}$ in the standard
model as upper half plane, and expressing $u$ in terms of these
coordinates, we assume that $\bfx$,
$\bfy,\,\partial_{t}\bfx,\,\partial_{t}\bfy$ are smooth, decay
toward infinity in the sense that
\[
\lim_{|x|\to\infty}({\mathbf{x}}(x),
{\mathbf{y}}(x))=({\mathbf{x}}_{0},{\mathbf{y}}_{0}) \in {\Hyp}^{2}
\]
and such that the derivative components
\[
\phi^{1}_{\alpha}=\frac{\partial_{\alpha}{\mathbf{x}}}{{\mathbf{y}}},\;\;
\phi^{2}_{\alpha}=\frac{\partial_{\alpha}{\mathbf{y}}}{{\mathbf{y}}},\;\;\alpha=0,1,2,
\]
are Schwartz, all at fixed time $t=0$. We make the following

\begin{defi}
We call initial data $\{{\mathbf{x}}, {\mathbf{y}},
\partial_{t}{\mathbf{x}}, \partial_{t} {\mathbf{y}}\}: \R^{2}\to {\Hyp}^{2}\times T \Hyp^{2}$ {\em{admissible}},
provided the derivative components $\phi^{k}_{\alpha}$ are Schwartz
functions for any $\alpha=0,1,2$ and $k=1,2$.
\end{defi}

We note here that the property of admissibility is propagated along
with the wave map flow on fixed time slices, as long as the wave map
persists and is smooth. This follows from finite propagation speed,
as well as the small-data well-posedness theory. We recall that the
energy associated with given initial data at time $t=0$ is given by
\[
E:=\int_{\R^{2}}\sum_{\alpha=0,1,2}[(\phi^{1}_{\alpha})^{2}+(\phi^{2}_{\alpha})^{2}]\,dx_{1}dx_{2}
\]
We now come to the first step in the Bahouri-Gerard decomposition of
a sequence of initial data, cf.~\cite{BG}. More precisely, we wish
to obtain a decomposition of the derivative initial data which is
analogous to the one of~\cite{BG}. An added feature for wave maps,
which does not appear in~\cite{BG}, consists of the fact that the
decomposition has the be performed in such a way that the individual
summands in it are themselves derivatives of admissible maps. This
requires some care, as the requisite condition is nonlinear, see
Lemma~\ref{lem:comp} below. In what follows we write
$\phi_{\alpha}:=\phi^{1}_{\alpha}+i\phi^{2}_{\alpha}$, any
additional  superscript  referring to the index of a sequence.

\begin{lemma}\label{lem:comp}
The complex-valued Schwartz functions $\phi_{\alpha}$, $\alpha=1,2$,
correspond to the derivative components of admissible data $\bfu:
\R^{2}\to \Hyp^{2}$ iff
\begin{equation}\label{eq:compat}
\partial_{k}\phi_{j}-\partial_{j}\phi_{k}=\phi^{1}_{k}\phi^{2}_{j}-\phi^{2}_{k}\phi^{1}_{j},\;\;k,j=1,2
\end{equation}
are satisfied.
\end{lemma}
\begin{proof} The ``only if'' part follows from~\eqref{eq:compat1}, \eqref{eq:compat2}. For the ``if'' part, note first that we get
\begin{equation}\label{eq:gl1}
\partial_{k}\phi^{2}_{j}-\partial_{j}\phi^{2}_{k}=0
\end{equation}
for the imaginary parts of $\phi_{j}$ and $\phi_k$. This implies
that
\[
\phi^{2}_{j}=\frac{\partial_{j}{\bfy}}{{\bfy}},\;\;j=1,2
\]
for a suitable positive function ${\bfy}: \R^{2}\to \R^+$ which is
unique only up to a multiplicative positive constant. We can
rewrite~\eqref{eq:compat} in the form
\begin{equation}\label{eq:gl2}
\partial_{k}({\bfy}\phi^{1}_{j})-\partial_{j}({\bfy}\phi^{1}_{k})=0,\;\;k,j=1,2
\end{equation}
which in turn implies that
\[
\phi^{1}_{j}=\frac{\partial_{j}{\bfx}}{{\bfy}}
\]
for a suitable function $\bfx: \R^{2}\to \R$.  To understand the
behavior of $(\bfx,\bfy)$ at infinity, we observe the following:
from~\eqref{eq:gl1},
\[
\partial_2\int_{-\infty}^\infty \phi^2_1(x_1,x_2)\,dx_1 =0
\]
which implies that the integral does not depend on~$x_2$ and
therefore is, in fact, zero. Similarly,
\[
\int_{-\infty}^\infty \phi^2_2(x_1,x_2)\,dx_2 =0  \qquad \forall\;
x_1\in \R
\]
It follows that $\bfy$ tends to the same constant at infinity
irrespective of the way in which we approach infinity. Without loss
of generality, we may set this constant equal to~$1$.
From~\eqref{eq:gl2} one further sees that
\[
\int_{-\infty}^\infty \bfy\,
\phi^1_1(x_1,x_2)\,dx_1=\int_{-\infty}^\infty \bfy\,
\phi^1_2(x_1,x_2)\,dx_2=0
\]
whence $\bfx$ approaches a constant $\bfx_0$ at $\infty$.
\end{proof}

Now for the first step in the concentration compactness method, which is the Metivier-Schochet scale selection process,
see~\cite{MS} and Section~III.1
of~\cite{BG}.  As already explained above, the difficulty we face here in contrast to~\cite{BG} is that we need to
make sure that the pieces we decompose the derivative components into are geometric, i.e., they are themselves
derivative components of maps $\R^2\to\Hyp^2$. Section~\ref{sec:bmo} provides us with the tools required for this purpose.

\begin{prop}\label{prop:CCI}
Let $\{{\bfx}_{n},
\,{\bfy}_{n},\,\partial_{t}{\bfx}_{n},\,\partial_{t}{\bfy}_{n}\}_{n\geq
1}$ be any  sequence of admissible data with energy bounded by~$E$
and with associated derivative sequence
$\{\phi^{n}_{\alpha}\}_{n\ge1}$, $\alpha=0,1,2$. Then up to passing
to a subsequence the following holds: given $\delta>0$, there exists
a positive integer $A=A(\delta, E)$ and a decomposition
\[
\phi^{n}_{\alpha}=\sum_{a=1}^{A}\phi^{na}_{\alpha}+w^{nA}_{\alpha}
\]
for $\alpha=0,1,2$ and $n\ge1$. Here the functions
$\phi^{na}_{\alpha}$, $1\le a\le A$ are derivative components of
admissible maps $\bfu^{a}_{n}: \R^{2}\to \Hyp^{2}$, and are
$\lambda_{n}^{a}$-oscillatory for a sequence of pairwise orthogonal
frequency scales $\{\lambda_{n}^{a}\}_{n\ge1}$ while the remainder
$w^{An}_{\alpha}$ is $\lambda_{n}^{a}$-singular  for each $1\le a\le
A$ and satisfies the smallness condition
\[
\sup_{n\ge1}\|w^{nA}_{\alpha}\|_{\dot{B}_{2,\infty}^{0}}<\delta
\]
Finally, given any sequence $R_n\to\infty$ one has the frequency localization
with $\mu_n^a:=-\log \lambda_n^a$,
\begin{equation}\label{eq:phinaloc}
 \sup_{\alpha=0,1,2} \| P_j \phi^{na}_{\alpha} \|_2 \le E\, R_n 2^{-\frac13|j-\mu_n^a|}  \qquad\forall \;j\in\Z
\end{equation}
for all $1\le a\le A$ and all large~$n$.
\end{prop}
\begin{proof}
We omit the time dependence in the notation, keeping in mind that
everything takes place at initial time $t=0$. As in Section~III.1 of~\cite{BG} one
obtains a decomposition
\begin{equation}\label{crudedecomp}
\phi^{n}_{\alpha}=\sum_{a=1}^{A}\tilde{\phi}^{na}_{\alpha}+\tilde{w}^{nA}_{\alpha},\;\;\alpha=0,1,2
\end{equation}
where the functions $\tilde{\phi}^{na}_{\alpha}\in L^{2}(\R^{2})$
are $\lambda^{a}_{n}$-oscillatory for suitable pairwise orthogonal
frequency scales $\{\lambda^{a}_{n}\}_{n\ge1}$ for all $1\le a\le
A$. Moreover, there is the smallness
\[
\|\tilde{w}^{nA}_{\alpha}\|_{\dot{B}^{0}_{2,\infty}}<\delta
\]
We now restrict to Fourier supports of these functions. Pick a
sequence $R_{n}\to \infty$ growing sufficiently slowly such that the
intervals $[(\lambda^{a}_{n})^{-1}R_{n}^{-1},
(\lambda^{a}_{n})^{-1}R_{n}]$ are mutually disjoint for $n$ large
enough and different values of~$a$. Then we replace
$\tilde{w}^{nA}_{\alpha}$ by
\[
P_{\cap_{a=1}^{A}[\mu_n^a-\log R_{n},
\mu_n^a+\log R_{n}]^{c}}
\tilde{w}^{nA}_{\alpha}+\sum_{a=1}^{A}P_{[\mu_n^a-\log
R_{n}, \mu_n^a+\log
R_{n}]^{c}}\tilde{\phi}^{na}_{\alpha}
\]
where $\mu_n^a:=-\log\lambda_n^a$,
while we replace each $\tilde{\phi}^{na}_{\alpha}$, $1\le a\le
A$, by
\[
P_{[\mu_n^a-\log R_{n},
\mu_n^a+\log R_{n}]}\tilde{\phi}^{na}_{\alpha}+
P_{[\mu_n^a-\log R_{n},
\mu_n^a+\log R_{n}]}\tilde{w}^{nA}_{\alpha}
\]
We need to make $R_n$ increase sufficiently slowly so that the
second term here remains $\lambda_n^a$-oscillatory. Then the new
decomposition, which we again refer to as
\[
\phi^{n}_{\alpha}=\sum_{a=1}^{A}\tilde{\phi}^{na}_{\alpha}+\tilde{w}^{nA}_{\alpha}
\]
has the same properties as the original one with the added advantage of the frequency localization around
the scales~$(\lambda_n^a)^{-1}$.  In particular, since the $\phi^n_\alpha$ are Schwartz functions, one concludes
that the $\tilde{\phi}^{na}_{\alpha}$ have the same property which means that
 the components
$\tilde{\phi}^{na}_{\alpha}$ are admissible, and so is~$\tilde{w}^{nA}_{\alpha}$.

In order to prove the proposition we need to show that we can
replace the components $\tilde{\phi}^{na}_{\alpha}$ by components
$\phi^{na}_{\alpha}$ which actually belong to admissible maps
$\bfu^{na}: \R^{2}\to \Hyp^{2}$ up to a small error (which again can be
absorbed into $\tilde{w}^{nA}_{\alpha}$). Note that the $\alpha=0$ component does
not present a problem here. For the $\alpha=1,2$ components, however, we need to
ensure that the compatibility relations~\eqref{eq:compat} hold.
Continuing with the proof of the Proposition~\ref{prop:CCI}, we
notice that
\[
{\bfx}^{n}=\sum_{k=1,2}\triangle^{-1}\partial_{k}[\phi^{1n}_{k}{\bfy}^{n}],\;\;
{\bfy}^{n}=e^{\sum_{k=1,2}\triangle^{-1}\partial_{k}\phi^{2n}_{k}}
\]
for the coordinate functions $({\bfx}^{n},\,{\bfy}^{n})$. In turn,
these identities imply that
\[
\phi^{1n}_{j}=({\bfy}^{n})^{-1}\sum_{k=1,2}\triangle^{-1}\partial_{j}\partial_{k}
[\phi^{1n}_{k}{\bfy}^{n}],\;\;\phi^{2n}_{j}=\sum_{k=1,2}\triangle^{-1}\partial_{j}\partial_{k}\phi^{2n}_{k}
\]
 These relations shall allow us to replace \eqref{crudedecomp}
 by a ``geometric decomposition''. Indeed, we simply substitute the decomposition~\eqref{crudedecomp} to obtain
 \begin{align*}
 \phi^{1n}_{j}&=\sum_{a=1}^{A}({\bfy}^{n})^{-1}
 \sum_{k=1,2}\triangle^{-1}\partial_{k}\partial_{j}[\tilde{\phi}^{1na}_{k}{\bfy}^{n}]
+({\bfy}^{n})^{-1}\sum_{k=1,2}\triangle^{-1}\partial_{j}\partial_{k}[\tilde{w}^{1nA}{\bfy}^{n}]\\
\phi^{2n}_{j}&=\sum_{a=1}^{A}\sum_{k=1,2}\triangle^{-1}\partial_{k}\partial_{j}\tilde{\phi}^{2na}_{k}
+\sum_{k=1,2}\triangle^{-1}\partial_{k}\partial_{j}\tilde{w}^{2nA}_{k}
\end{align*}
This suggests making the following choices:
\begin{align*}
{\bfx}^{na}
&:=\sum_{k=1,2}\triangle^{-1}\partial_{k}[\tilde{\phi}^{1na}{\bfy}^{na}],\quad
{\bfy}^{na}
:=e^{\sum_{k=1,2}\triangle^{-1}\partial_{k}\tilde{\phi}^{2na}_{k}}
\end{align*}
and then defining
\begin{align*}
\phi^{1na}_{j}
&:=({\bfy}^{na})^{-1}\sum_{k=1,2}\triangle^{-1}\partial_{j}\partial_{k}[\tilde{\phi}^{1na}{\bfy}^{na}],\quad
\phi^{2na}_{j}
:=\sum_{k=1,2}\triangle^{-1}\partial_{j}\partial_{k}\tilde{\phi}^{2na}_{k} \\
w^{1nA}_{j}
&:=({\bfy}^{n})^{-1}\sum_{k=1,2}\triangle^{-1}\partial_{j}\partial_{k}[\tilde{w}^{1nA}{\bfy}^{n}],\quad
w^{2nA}_{j}
:=\sum_{k=1,2}\triangle^{-1}\partial_{j}\partial_{k}\tilde{w}^{2nA}_{k}
\end{align*}
as well as $\phi^{na}:=\phi^{1na}+i\phi^{2na}$, $w^{nA}:= w^{1nA}+iw^{2nA}$.
Clearly the components $\phi^{1na}_{j}$, $\phi^{2na}_{j}$ are now
geometric in the sense that they derive from a map into hyperbolic
space; in fact, they are associated with the maps given by the
components $({\bfx}^{na},\,{\bfy}^{na})$. The proof is now concluded
by appealing to Lemma~\ref{lem:geom_dec}, Corollary~\ref{cor:Jlem},
and Lemma~\ref{lem:weight_remov}.
For the final statement, note that by Lemma~\ref{lem:geom_dec}, the ``geometric'' components
$\phi^{na}_{\alpha}$ are also frequency localized to the interval $[\mu_n^a-\log R_{n},
\mu_n^a+\log R_{n}]$ up to
exponentially decaying errors.
\end{proof}

As an immediate consequence of Proposition~\ref{prop:CCI} one
obtains that $\phi_{j}^{kna}, w^{knA}_{j}$, $k=1,2$, are
asymptotically  orthogonal (where $\phi_{j}^{1na}=\Re\,
\phi_{j}^{na}$ and $\phi_{j}^{2na}=\Im\, \phi_{j}^{na}$).

\noindent We now make some preparations for the second stage of the
Bahouri-Gerard procedure. More specifically, we shall have to pass
to the Coulomb gauge components, $\psi_{\alpha}$, and transfer the
above decomposition to the level of these components. One can split
\[
\psi_{\alpha}^{n}=\phi^{n}_{\alpha}e^{-i\sum_{k=1,2}\triangle^{-1}\partial_{k}\phi^{1n}_{k}}=\big[\sum_{a=1}^{A}\phi^{na}_{\alpha}+
w^{nA}_{\alpha}\big]e^{-i\sum_{k=1,2}\triangle^{-1}\partial_{k}\phi^{1n}_{k}}
\]
However,  the components
\[
\phi^{na}_{\alpha}e^{-i\sum_{k=1,2}\triangle^{-1}\partial_{k}\phi^{1n}_{k}}
\]
are not the Coulomb gauge components of a suitable wave map, and
should ideally be replaced by
\[
\phi^{na}_{\alpha}e^{-i\sum_{k=1,2}\triangle^{-1}\partial_{k}\phi^{1na}_{k}}
\]
Due to the lack of $L^{\infty}$ control over the exponent, this
cannot be done without further physical localizations. Nevertheless,
we can state the following fact.

\begin{lemma}\label{lem:Coul}
The components
\[
\phi^{na}_{\alpha}e^{-i\sum_{k=1,2}\triangle^{-1}\partial_{k}\phi^{1n}_{k}},
\qquad
w^{nA}_{\alpha}e^{-i\sum_{k=1,2}\triangle^{-1}\partial_{k}\phi^{1n}_{k}}
\]
are $\lambda^{a}_{n}$-oscillatory and $\lambda^{a}_{n}$-singular,
respectively, for each~$a$ and we have
\[
\big\|w^{nA}_{\alpha}e^{-i\sum_{k=1,2}\triangle^{-1}\partial_{k}\phi^{1n}_{k}}\big\|_{\dot{B}^{0}_{2,\infty}}\lesssim
\delta
\]
where $\delta$ is as in Proposition~\ref{prop:CCI}.
\end{lemma}
\begin{proof}
We may assume $\lambda^{a}_{n}=1$ by scaling invariance. Given any
$\eps>0$, we can choose $k_0$ large enough such that
\[
\limsup_{n\to \infty}\|P_{[-k_0,
k_0]^{c}}\phi^{na}_{\alpha}e^{-i\sum_{k=1,2}\triangle^{-1}\partial_{k}\phi^{1n}_{k}}\|_{L^{2}}<\eps
\]
Next, for $k_1>k_0+C$, consider the expressions
\[
P_{<-k_1}\big[P_{[-k_0,
k_0]}\phi^{na}_{\alpha}e^{-i\sum_{k=1,2}\triangle^{-1}\partial_{k}\phi^{1n}_{k}}\big],\quad
P_{>k_1}\big[P_{[-k_0,
k_0]}\phi^{na}_{\alpha}e^{-i\sum_{k=1,2}\triangle^{-1}\partial_{k}\phi^{1n}_{k}}\big].
\]
Start with the first expression, which we write as
\begin{equation}\nonumber\begin{split}
&P_{<-k_1}\big[P_{[-k_0,
k_0]}\phi^{na}_{\alpha}e^{-i\sum_{k=1,2}\triangle^{-1}\partial_{k}\phi^{1n}_{k}}\big]
\\&=P_{<-k_1}\big[P_{[-k_0, k_0]}\phi^{na}_{\alpha}\sum_{j=1,2}\triangle^{-1}\partial_{j}P_{[-k_0, k_0]}([-i\sum_{k=1,2}
\triangle^{-1}\partial_{j}\partial_{k}\phi^{1n}_{k}]e^{-i\sum_{k=1,2}\triangle^{-1}\partial_{k}\phi^{1n}_{k}})\big]
\end{split}\end{equation}
Using Bernstein's inequality, we can then estimate
\begin{equation}\nonumber\begin{split}
&\|P_{<-k_1}\big[P_{[-k_0,
k_0]}\phi^{na}_{\alpha}\sum_{j=1,2}\triangle^{-1}\partial_{j}P_{[-k_0,
k_0]}([-i\sum_{k=1,2}
\triangle^{-1}\partial_{j}\partial_{k}\phi^{1n}_{k}]e^{-i\sum_{k=1,2}\triangle^{-1}\partial_{k}\phi^{1n}_{k}})\big]\|_{L^{2}}\\
&\lesssim 2^{k_0-k_1}\|P_{[-k_0,
k_0]}\phi^{na}_{\alpha}\|_{L^{2}}\|\phi^{1n}\|_{L^{2}}<\eps
\end{split}\end{equation}
provided we choose $k_1$ sufficiently large in relation to $k_0$.
The estimate for the second term is more of the same.  Next,
consider the ``tail term'' $
 w^{nA}_{\alpha}e^{-i\sum_{k=1,2}\triangle^{-1}\partial_{k}\phi^{1n}_{k}}
$. That this is $\lambda_n^a$-singular for each $1\le a\le A$
follows from the preceding via duality. It therefore remains
 to estimate its $\|\cdot\|_{\dot{B}^{0}_{2,\infty}}$-norm. We localize this term to fixed dyadic frequency $\sim 2^{q}$
 \begin{equation}\nonumber\begin{split}
&
P_{q}\big[w^{nA}_{\alpha}e^{-i\sum_{k=1,2}\triangle^{-1}\partial_{k}\phi^{1n}_{k}}\big]
 =P_{q}\big[w^{nA}_{\alpha}P_{<q-10}e^{-i\sum_{k=1,2}\triangle^{-1}\partial_{k}\phi^{1n}_{k}}\big]\\
 &+P_{q}\big[w^{nA}_{\alpha}P_{[q-10, q+10]}e^{-i\sum_{k=1,2}\triangle^{-1}\partial_{k}\phi^{1n}_{k}}\big]
+P_{q}\big[w^{nA}_{\alpha}P_{>q+10}e^{-i\sum_{k=1,2}\triangle^{-1}\partial_{k}\phi^{1n}_{k}}\big]
\end{split}\end{equation}
and estimate the three terms on the right separately: first, we have
 \begin{equation}\nonumber\begin{split}
&\|P_{q}\big[w^{nA}_{\alpha}P_{<q-10}e^{-i\sum_{k=1,2}\triangle^{-1}\partial_{k}\phi^{1n}_{k}}\big]\|_{L^{2}}
=\|P_{q}\big[P_{[q-10, q+10]}(w^{nA}_{\alpha})P_{<q-10}e^{-i\sum_{k=1,2}\triangle^{-1}\partial_{k}\phi^{1n}_{k}}\big]\|_{L^{2}}\\
&\lesssim \|P_{[q-10, q+10]}(w^{nA}_{\alpha})\|_{L^{2}}\lesssim
\|w^{nA}_{\alpha}\|_{\dot{B}^{0}_{2,\infty}}\lesssim \delta
\end{split}\end{equation}
Next,
 \begin{equation}\nonumber\begin{split}
&\|P_{q}\big[w^{nA}_{\alpha}P_{[q-10, q+10]}e^{-i\sum_{k=1,2}\triangle^{-1}\partial_{k}\phi^{1n}_{k}}\big]\|_{L^{2}}\\
&=\|P_{q}\big[P_{<q+10}(w^{nA}_{\alpha})P_{[q-10, q+10]}e^{-i\sum_{k=1,2}\triangle^{-1}\partial_{k}\phi^{1n}_{k}}\big]\|_{L^{2}}\\
&=\|P_{q}\big[P_{<q+10}w^{nA}_{\alpha}\sum_{j=1,2}\triangle^{-1}\partial_{j}P_{[q-10,
q+10]}([-i\sum_{k=1,2}\triangle^{-1}\partial_{j}
\partial_{k}\phi^{1n}_{k}]e^{-i\sum_{k=1,2}\triangle^{-1}\partial_{k}\phi^{1n}_{k}})\big]\|_{L^{2}}\\
&\lesssim
2^{-q}\|P_{<q+10}w^{nA}_{\alpha}\|_{L^{\infty}}\|\phi^{1n}\|_{L^{2}}
\lesssim \|w^{nA}_{\alpha}\|_{\dot{B}^{0}_{2,\infty}}\lesssim \delta
\end{split}\end{equation}
where Bernstein's inequality was used in the last step. The third
term in the above Littlewood-Paley trichotomy corresponding to
high-high interactions, is treated analogously and omitted.
\end{proof}

By letting $A\to\infty$ the construction of this section yields an
infinite double sequence $\{\phi^{na}_{\alpha}\}_{a\geq 1}$.  For
later reference, it shall be important to  construct ``partial
approximations'' of the components $\phi^{n}_{\alpha}$ in terms of
the $\phi^{na}_{\alpha}$. Specifically, for $I\subset \{1,2,\ldots,
A\}$, we let
\[
\tilde{\phi}^{nI}_{\alpha}:=\sum_{a\in I}\tilde{\phi}^{na}_{\alpha}
\]
Then reasoning exactly as in the preceding, and employing the same
notation as there, one obtains the following statement.

\begin{cor} Let
\[
{\bfy}^{nI}:=e^{\sum_{k=1,2}\sum_{a\in
I}\triangle^{-1}\partial_{k}\tilde{\phi}^{2na}_{k}},\qquad
{\bfx}^{nI}:= \sum_{k=1,2}\sum_{a\in
I}\triangle^{-1}\partial_{k}[\tilde{\phi}_k^{1na}{\bfy}^{nI}]
\]
Then  for $a\in I$
\[
 \phi^{1na}_{j}=({\bfy}^{nI})^{-1}\sum_{k=1,2}\triangle^{-1}\partial_{j}\partial_{k}[\tilde{\phi}_k^{1na}{\bfy}^{nI}]+o_{L^{2}}(1)
 \]
 In particular, we have
 \[
 \sum_{a\in I}\phi^{1na}_{j}=\frac{\partial_{j}{\bfx}^{nI}}{{\bfy}^{nI}}+o_{L^{2}}(1),\qquad \sum_{a\in I}\phi^{2na}_{j}
 =\frac{\partial_{j}{\bfy}^{nI}}{{\bfy}^{nI}}+o_{L^{2}}(1)
 \]
\end{cor}

\subsection{Step 2: frequency localized approximations to the data}
\label{subsec:BGstep2}

Given an essentially singular sequence ${\bf{u}}^{n}$ with
derivatives $\phi^{n}_{\alpha}$, Proposition~\ref{prop:CCI} yields
a new essentially singular sequence $\phi^{n}_{\alpha}$ with the
following property: for any $A\geq 1$ (recall the
$\phi^{na}_{\alpha}$ are defined inductively)
\[
\phi^{n}_{\alpha}=\sum_{a=1}^{A}\phi^{na}_{\alpha}+w^{nA}_{\alpha}
\]
Given $\delta_0>0$, there exists $A\geq 1$  so that $\|w^{nA}_{\alpha}\|_{\dot{B}^{0}_{2,\infty}}<\delta_0$ for large~$n$.
In addition, we assume that for each~$n$, the energy $\sum_{\alpha=0}^2 \|\phi^{na}_\alpha\|_2^2$ is nonincreasing in~$a$.
In what follows, we will use smallness parameters $1\gg \eps_0\gg \delta_1\gg \delta_0>0$, each of which will eventually be
chosen depending only on the energy of the initial data.

\begin{figure*}[ht]
\begin{center}
\centerline{\hbox{\vbox{ \epsfxsize= 8.0 truecm \epsfysize=6.5
truecm \epsfbox{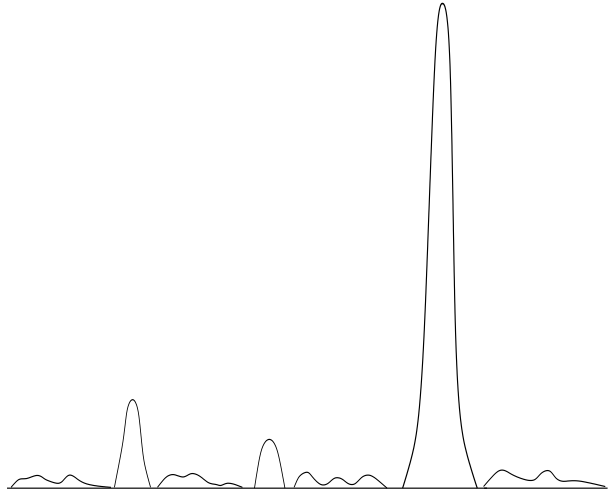}}}} \caption{Atoms and the Besov
error}
\end{center}
\end{figure*}

\noindent Ultimately we wish to show that there can only be a single
frequency block, i.e., $A=1$,  and furthermore, that the energy of
this block converges to the critical energy $\Ecrit$ as $n\to
\infty$.  Thus we now use the following {\em{dichotomy}}:
\begin{itemize}
 \item  We have $A=1$ and $\lim_{n\to \infty}\sum_{\alpha=0,1,2}\|\phi^{na}_{\alpha}\|_{L_{x}^{2}}^{2}=\Ecrit$.
\item  The previous scenario does not occur. Thus, for a suitable subsequence
\[
\limsup_{n\to
\infty}\sum_{\alpha=0,1,2}\|\phi^{na}_{\alpha}\|_{L_{x}^{2}}^{2}<\Ecrit-\delta_{2}
\]
for some $\delta_{2}>0$, and all $a$.
\end{itemize}

If the first alternative occurs, then continue with Step~4 below. Hence we now assume that the first alternative occurs,
in which case we will show that the sequence ${\bf{u}}^{n}$
cannot be essentially singular. We may of course assume that for
each $1\le a\le A$,
\[
\liminf_{n\to
\infty}\sum_{\alpha=0,1,2}\|\phi^{na}_{\alpha}\|_{L_{x}^{2}}>0,
\]
as otherwise we may pass to a subsequence for which the
$\phi^{na}_{\alpha}$ may be absorbed into the error
$w^{nA}_{\alpha}$.
The issue now becomes how to choose the
cutoff~$A$. Due to the asymptotic orthogonality of the
$\phi^{na}_{\alpha}$ as $n\to \infty$, and for each $\alpha=0,1,2$,
\[
\lim_{A_{0}\to \infty}\sum_{a\geq A_{0}}\limsup_{n\to
\infty}\|\phi^{na}_{\alpha}\|_{L_{x}^{2}}^{2}=0
\]
For some absolute  $\eps_0>0$ which is small enough,
in particular smaller than the cutoff for the small energy global
well-posedness theory, we choose $A_{0}$ large enough such that
\[
\sum_{a\geq A_{0}}\limsup_{n\to
\infty}\|\phi^{na}_{\alpha}\|_{L_{x}^{2}}^{2}<\eps_0,
\]
and then put $A=A_{0}$. Thus we now arrive at the decomposition
\[
\phi^{n}_{\alpha}=\sum_{a=1}^{A_{0}}\phi^{na}_{\alpha}+w^{nA_{0}}_{\alpha}
\]
We may further decompose
\[
w^{nA_{0}}_{\alpha}=\sum_{a=A_{0}+1}^{A}\phi^{na}_{\alpha}+w^{nA}_{\alpha},
\]
with the smallness property
\[
\sum_{a\geq
A_{0}+1}\sum_{\alpha=0,1,2}\|\phi^{na}_{\alpha}\|_{L_{x}^{2}}^{2}<\eps_0
\]
By adjusting $A$, we can further achieve
\[
\|w_{\alpha}^{nA}\|_{\dot{B}_{2, \infty}^{0}}<\delta_{0}
\]
for any given $\delta_{0}>0$.

\noindent
Re-ordering the superscripts if
necessary, we may assume that the frequency scales $(\lambda_n^a)^{-1}$ of the
$\phi^{na}_{\alpha}$ are increasing with~$1\le a\le A_0$. The error term
$w^{nA_{0}}_{\alpha}$
may be written as a sum of constituents
\[
w^{nA_{0}}_{\alpha}=w^{nA_{0}^{(0)}}_{\alpha}+w^{nA_{0}^{(1)}}_{\alpha}+\ldots+w^{nA_{0}^{(A_{0})}}_{\alpha} + o_{L^{2}}(1)
\]
which satisfy the property that
\begin{equation}
 \label{eq:w_teil}
w^{nA_{0}^{(k)}}_{\alpha}=P_{\mu_{n}^{k-1}+L_n <\cdot<\mu_{n}^{k}-L_n}w^{nA_{0}^{(k)}}_{\alpha} \text{\ \ as\ \ }n\to\infty
\end{equation}
with $\mu_n^a:=-\log \lambda_n^a$ and a sequence $L_n\to\infty$ which increases very slowly.
This can be done since $w_\alpha^{nA_0}$ is $\lambda_n^a$-singular for each~$1\le a\le A_0$.
Thus the frequency support of $w^{nA_{0}^{(k)}}_{\alpha}$ is
contained in the annulus
\[
(\lambda_{n}^{k})^{-1}e^{L_n} <|\xi|<(\lambda_{n}^{k+1})^{-1}e^{-L_n} ,\quad
(\lambda_{n}^{0})^{-1}:=0,\quad (\lambda_{n}^{A_0^{(A_0)}+1})^{-1}:=\infty
\]
Figure~5 above is a schematic depiction of the situation $A_0=1$
with a unique large atom on the right, but with two smaller atoms on
the left which are too large to be included in the Besov error (the
three bumpy curves between the atoms). More precisely,
$w_\alpha^{nA_0}$ consists of the four small curves between the
atoms, and $w_\alpha^{n A_0^{(0)}}$ is the sum of the three curves
to the left of the big atom {\em together with the two small atoms},
and $w_\alpha^{n A_0^{(1)}}$ the one to the right of the big atom.

\noindent  Note that if we refine the frequency decomposition, i.e.,
increase $A_{0}$ to $A^{(k)}\geq A_{0}$, then the components
$w^{nA_{0}^{(k)}}_{\alpha}$ are decomposed into
\[
w^{nA_{0}^{(k)}}_{\alpha}=\sum_{j}\phi^{na_{j}^{k}}_{\alpha}+w^{nA^{(k)}}_{\alpha}
\]
for suitable $a_{j}^{k}\in [A_{0}+1, A^{(k)}]$. In Figure~5 one has $j=1,2$ for $k=0$ corresponding
to the two small atoms to the left of the large one. We may again assume
that the $a_{j}^{k}$ are increasing in $j$ and have frequency
support with increasing value of $|\xi|$, for each $k$. Furthermore,
we have
\[
\limsup_{n\to \infty}\sum_{\alpha=0, 1,
2}\sum_{j}\|\phi^{na_{j}^{k}}_{\alpha}\|_{L_{x}^{2}}^{2}<\eps_0
\]
by asymptotic orthogonality and the choice of~$A_0$.
Our first goal, to be dealt with in the following section, is to
control the nonlinear evolution of the minimum frequency components
$w^{nA_{0}^{(0)}}_{\alpha}$. The idea behind this is as follows: due
to the energy constraint
\[
\limsup_{n\to \infty}\|w^{nA_{0}^{(0)}}_{\alpha}\|_2\leq \Ecrit,
\]
we may subdivide $w^{nA_{0}^{(0)}}_{\alpha}$ into finitely many
pieces be means of frequency localizations $\big\{P_{J_{\ell}}
w^{nA_{0}^{0}}_{\alpha}\big\}_{1\le \ell\le
\frac{1000\Ecrit}{\eps_0}}$ such that the dyadic intervals
$J_{\ell}$ are disjoint, with $\cup_{\ell}J_{\ell}=(-\infty,
(\lambda_{n}^{1})^{-1}e^{-L_n})$, and furthermore
\[
\big\|P_{J_{\ell}}w^{nA_{0}^{(0)}}_{\alpha}\big\|_{L_{x}^{2}}< \eps_0 \quad
\forall \ell
\]
Recall that $(\lambda_{n}^{1})^{-1}$ is the frequency scale of the
first frequency atom $\phi^{n1}_{\alpha}$. In particular, this means
that the frequency localized pieces
$P_{J_{\ell}}w^{nA_{0}^{0}}_{\alpha}$ should be treatable via a
perturbative argument. More precisely, we shall
 run an induction in~$\ell$ on a sequence of  approximating maps with
data
\[
\sum_{1\leq j\leq
\ell}P_{J_{j}}w^{nA_{0}^{(0)}}_{\alpha}  \:e^{-i\Re\sum_{k=1,2}\triangle^{-1}\partial_{k}\sum_{1\leq
j\leq \ell}P_{J_{j}}w^{nA_{0}^{(0)}}_{k}}
\]
As always, we face the issue at this point that these gauged components are not necessarily
admissible, i.e., they are not given by derivative components of maps~$\R^2\to\Hyp^2$.
In order to apply the perturbative theory we shall need to show that they are close to such
admissible data.
This in turn follows from Lemma~\ref{lem:elimhigh} provided we chose the intervals $J_\ell$ carefully; for this it is
essential that the endpoints of these intervals do not fall onto one of the 'small' atoms $\phi^{na}_\alpha$.
Otherwise, condition~\eqref{eq:Bcondi} would be violated. In detail, this is done as follows.
Recall that the $J_{j}$ are chosen to
be disjoint and such that
\[
w^{nA_{0}^{(0)}}_{\alpha}=\sum_{j}P_{J_{j}}w^{nA_{0}^{(0)}}_{\alpha},\qquad
\sup_j
\big\|P_{J_{j}}w^{nA_{0}^{(0)}}_{\alpha}\big\|_{L_{x}^{2}}\les \eps_0
\]
On the other hand, upon refining the Bahouri-Gerard frequency
decomposition applied to $w^{nA_{0}^{(0)}}_{\alpha}$, we can also
write
\begin{equation}
 \label{eq:wfurtherrefine}
w^{nA_{0}^{(0)}}_{\alpha}=\sum_{j\geq
1}\phi^{na_{j}^{(0)}}_{\alpha}+w^{nA^{(0)}}_{\alpha}
\end{equation}
 Here $A^{(0)}>A_{0}^{(0)}$ is chosen such that $\|w^{nA^{(0)}}_{\alpha}\|_{\dot{B}_{2,\infty}^{0}}\ll\delta_{0}$ for some
 constant $\delta_{0}>0$ which is to be determined, while the $a_{j}^{(0)}$ are certain indices in the interval $[A_{0}^{(0)}, A^{(0)}]$.
Our choice of $A_{0}$ ensures that
\[
\limsup_{n\to \infty}\sum_{j\geq
1}\big\|\phi^{na_{j}^{(0)}}_{\alpha}\big\|_{L_{x}^{2}}^{2}<\eps_0
\]
Now, to choose the $J_{j}$, pick for each of the
$\phi^{na_{j}^{(0)}}_{\alpha}$ (which are finite in number) a
frequency interval
\[
\big[(\lambda_n^{a_{j}^{(0)}}R_{j}^{(0)})^{-1},
(\lambda_n^{a_{j}^{(0)}})^{-1}R_{j}^{(0)}\big]
\]
with $R_{j}^{(0)}$ large enough such that
\begin{equation} \label{eq:freq_sep}
\limsup_{n\to
\infty}\big\|P_{[\log(\lambda_n^{a_{j}^{(0)}})^{-1}-\log R_{j}^{(0)},
\log(\lambda_n^{a_{j}^{(0)}})^{-1} +\log
R_{j}^{(0)}]^{c}}\phi^{na_{j}^{(0)}}_{\alpha}\big\|_{L_{x}^{2}}\ll\delta_{0},
\end{equation}
which is possible due to the frequency localization of the atoms
$\phi^{na_{j}^{(0)}}_{\alpha}$. Here $\delta_{0}>0$ is a
sufficiently small constant  such that
$\delta_{0}=\delta_{0}(\Ecrit, \eps_0)$. Picking $n$ large enough,
we may assume that the intervals
\[
\big[(\lambda_n^{a_{j}^{(0)}}R_{j}^{(0)})^{-1},
(\lambda_n^{a_{j}^{(0)}})^{-1}R_{j}^{(0)} \big]
\]
are disjoint. We can now exactly specify how to select the $J_{j}$:
inductively, assume that \[J_{1}=[a_{1}, b_{1}], \ldots,
J_{k-1}=[a_{k-1}, b_{k-1}]\] have been chosen. Then pick
$\tilde{J}_{k}=[\tilde{a}_{k}, \tilde{b}_{k}]$ such that
$\tilde{a}_{k}=b_{k-1}$ and such that the integer $\tilde{b}_{k}$ is
maximal with the property that
\[
\sum_{\alpha=0,1,2}\big\|P_{[\tilde{a}_{k},
\tilde{b}_{k}]}w^{nA_{0}^{(0)}}_{\alpha}\big\|_{L_{x}^{2}}^{2}<
\eps_0
\]
Then if $\tilde{b}_{k}\in [\log(\lambda_n^{a_{j}^{(0)}})^{-1}-\log
R_{j}^{(0)}, \log(\lambda_n^{a_{j}^{(0)}})^{-1}+\log R_{j}^{(0)}]$
for some $j$, we let
\[
b_{k}=\log(\lambda_n^{a_{j}^{(0)}})^{-1}+\log R_{j}^{(0)}
\]
Otherwise, we let $b_{k}=\tilde{b}_{k}$.
The point of this construction is that
if the endpoint of $\tilde{J}_{k}$ happens to fall on a ``small
atom'' which may still be too large in $\dot{B}_{2,\infty}^{0}$ for
our later purposes, we simply absorb this atom into~$J_{k}$.

\noindent We can now state the approximate admissibility fact alluded to above.   Recall that $\Re\, w^{nA_{0}^{(0)}}_{k} =
w^{1nA_{0}^{(0)}}_{k}$. Moreover, the constant $\delta_0$ controls the Besov norm
of the tails and is kept fixed.  We begin with a statement which does not involve the~$J_\ell$.

\begin{lemma}\label{lem:wapproximation} There is an admissible map $\R^{2}\to {\mathbb{H}}^{2}$ with
derivative components $\Phi^{nA_{0}^{(0)}}_{\alpha}$ such that
\[
\Big\| w^{nA_{0}^{(0)}}_{\alpha}e^{-i\sum_{k=1,2}\triangle^{-1}\partial_{k}w^{1nA_{0}^{(0)}}_{k}}
 -\Phi^{nA_{0}^{(0)}}_{\alpha}e^{-i\sum_{k=1,2}\triangle^{-1}\partial_{k}\Phi^{1nA_{0}^{(0)}}_{k}}\Big\|_{L^2_x} \to 0
\]
as $n\to \infty$. The same applies to the difference  $w^{nA_{0}^{(0)}}_{\alpha}
 -\Phi^{nA_{0}^{(0)}}_{\alpha}$.
\end{lemma}
\begin{proof}
Recall the relation that defines
$w^{nA_{0}^{(0)}}_{j}=w^{1nA_{0}^{(0)}}_{j}+iw^{2nA_{0}^{(0)}}_{j}$:
\[
w^{1nA_{0}^{(0)}}_{j}=({\bfy}^{n})^{-1}\sum_{k=1,2}\triangle^{-1}\partial_{k}\partial_{j}[\tilde{w}_{k}^{1nA_{0}^{(0)}}{\bfy}^{n}],
\qquad
w^{2nA_{0}^{(0)}}_{j}=\sum_{k=1,2}\triangle^{-1}\partial_{k}\partial_{j}\tilde{w}^{2nA_{0}^{(0)}}_{j}
\]
We now claim that the components $w^{1nA_{0}^{(0)}}_{j}$,
$w^{2nA_{0}^{(0)}}_{j}$ are $o_{L^{2}}(1)$- close to the derivative components $\Phi^{1,2nA_{0}^{(0)}}_{j}$ of a map, when $n\to \infty$.
 Moreover, the error satisfies $\nabla^{-1}o_{L^{2}}(1)=o_{L^{\infty}}(1)$.
First, observe that  by Corollary~\ref{cor:elimhighlimit},  the component
$
w^{1nA_{0}^{(0)}}_{j}
$
is close in the above sense to
 \[
 \Phi_j^{1nA_0^{(0)}} := ({\bfy}^{nA_{0}^{(0)}})^{-1}\sum_{k=1,2}\triangle^{-1}\partial_{k}\partial_{j}[\tilde{w}_{k}^{1nA_{0}^{(0)}}{\bfy}^{nA_{0}^{(0)}}],
\qquad
{\bfy}^{nA_{0}^{(0)}}:=e^{\sum_{k=1,2}\triangle^{-1}\partial_{k}\tilde{w}_k^{2nA_{0}^{(0)}}}
 \]
Next,  introduce the
auxiliary map $({\bfx}^{nA_{0}^{(0)}},\, {\bfy}^{nA_{0}^{(0)}}):
\R^{2}\to{\mathbb{H}}^{2}$, with components defined by
\[
{\bfx}^{nA_{0}^{(0)}}
:=\sum_{k=1,2}\triangle^{-1}\partial_{k}[\tilde{w}_k^{1nA_{0}^{(0)}}{\bfy}^{nA_{0}^{(0)}}],
\qquad
{\bfy}^{nA_{0}^{(0)}}=e^{\sum_{k=1,2}\triangle^{-1}\partial_{k}\tilde{w}_k^{2nA_{0}^{(0)}}}
\]
Furthermore, as before we have
\[
w^{1nA_{0}^{(0)}}_{j}=({\bfy}^{n})^{-1}\sum_{k=1,2}\triangle^{-1}\partial_{k}\partial_{j}[\tilde{w}_{k}^{1nA_{0}^{(0)}}{\bfy}^{n}],
\qquad
w^{2nA_{0}^{(0)}}_{j}=\sum_{k=1,2}\triangle^{-1}\partial_{j}\partial_{k}\tilde{w}_{k}^{2nA_{0}^{(0)}}
\]
and we set
$ \Phi^{2nA_{0}^{(0)}}_{j}:= w^{2nA_{0}^{(0)}}_{j}$, $w^{1,2nA_{0}^{(0)}}_{0}=\tilde{w}^{1,2nA_{0}^{(0)}}_{0}$. In view of the
preceding,
\[
w^{nA_{0}^{(0)}}_{\alpha}e^{-i\sum_{k=1,2}\triangle^{-1}\partial_{k}w^{1nA_{0}^{(0)}}_{k}}=\Phi^{nA_{0}^{(0)}}_{\alpha}
e^{-i\sum_{k=1,2}\triangle^{-1}\partial_{k}\Phi^{1nA_{0}^{(0)}}_{k}}+o_{L^{2}}(1)
\]
as $n\to \infty$.
\end{proof}

A similar result now applies to the frequency localized pieces. This time one has to use Lemma~\ref{lem:elimhigh}.

\begin{lemma}\label{lem:Truncatedwapproximation} Given any $\delta_{1}>0$ one can choose $\delta_0\ll \delta_1$ such that
for all large~$n$
\[
\sum_{j\leq
\ell}P_{J_{j}}w^{nA_{0}^{(0)}}_{\alpha}e^{-i\sum_{k=1,2}\triangle^{-1}\partial_{k}\sum_{j\leq
\ell}P_{J_{j}}w^{1nA_{0}^{(0)}}_{k}},\quad \ell\geq 1
\]
may be approximated within $\delta_{1}$ in the energy topology by
Coulomb components
\begin{equation}
 \label{eq:Psielldef}
\Psi^{\ell nA_{0}^{(0)}}_{\alpha}:=\Phi^{\ell nA_{0}^{(0)}}_{\alpha}  \:e^{-i\Re \sum_{k=1,2}\triangle^{-1}\partial_{k}\Phi^{\ell nA_{0}^{(0)}}_{k}}
\end{equation}
of actual maps from $\R^{2}\to{\mathbb{H}}^{2}$,
uniformly in $\ell$. The same statement holds for the functions without any exponential phases.
\end{lemma}
\begin{proof}
This follows exactly along the lines of the proof of
Lemma~\ref{lem:wapproximation}: for the components
$
\sum_{j\leq \ell}P_{J_{j}}w^{nA_{0}^{(0)}}_{\alpha}$
we use the approximating maps
\[
{\bfx}^{\ell
nA_{0}^{(0)}}:=({\bfy}^{\ell nA_{0}^{(0)}})^{-1}\sum_{k=1,2}\sum_{j\leq
\ell}\triangle^{-1}
\partial_{k}\big[P_{J_{j}}\tilde{w}_k^{1nA_{0}^{(0)}}{\bfy}^{\ell nA_{0}^{(0)}}],\qquad {\bfy}^{\ell nA_{0}^{(0)}}:=e^{\sum_{k=1,2}
\triangle^{-1}\partial_{k}\sum _{j\leq
\ell}\tilde{w}^{2nA_{0}^{(0)}}_{k}}
\]
However, this time, the smallness of the error is contingent on the
$\|\cdot\|_{\dot{B}_{2,\infty}^{0}}$-norm of the non-atomic
 part of $\tilde{w}^{2nA_{0}^{(0)}}_{\alpha}$, while the contribution of the atomic part can be made small by choosing $n$ large enough.
More precisely, \eqref{eq:Bcondi} holds for all large~$n$ due to the frequency separation properties which we have imposed on the various
components, see~\eqref{eq:freq_sep} and~\eqref{eq:w_teil}. These separations become effective for large~$n$ due to the orthogonality of
the scales involved.
\end{proof}

As a general comment, we would like to remind the reader that all constructions here are not unique; moreover, they are subject to
errors of the form~$o_{L^2}(1)$ as $n\to\infty$.

\subsection{Step 3: Evolving the lowest-frequency nonatomic part}\label{subsec:BGstep3}

As far as the evolution of $w^{nA_{0}^{(0)}}_{\alpha}$ is concerned,
we claim the following result.  Note that we phrase it in terms of
the derivative components that we just constructed. Once we have
evolved {\em all constituents} of the decomposition from Step~1, the
perturbative theory of Section~\ref{sec:perturb} will then allow us
to conclude that the representation that we obtain is accurate up to
a small energy error {\em globally in time}.

\begin{prop}\label{ControlNonatomicComponent1}
Let $\Phi^{nA_{0}^{(0)}}_{\alpha}$ be as in Lemma~\ref{lem:wapproximation} and set
\[
 \Psi^{nA_{0}^{(0)}}_{\alpha} :=  \Phi^{nA_{0}^{(0)}}_{\alpha}e^{-i\sum_{k=1,2}\triangle^{-1}\partial_{k}\Phi^{1nA_{0}^{(0)}}_{k}}
\]
Then provided $\eps_0\gg \delta_1 \gg \delta_0>0$ above are chosen
sufficiently small, and provided $n$ is large enough, the
$\Phi^{nA_{0}^{(0)}}_{\alpha}$ exist globally in time as derivative
components of an admissible wave map. Moreover, there is a constant
$C_{1}(\Ecrit)$ such that the solution of the  gauged counterparts
of these components, i.e.,
 $\Psi^{nA_{0}^{(0)}}_{\alpha}$   satisfy the bound
\[
\sup_{T_{0,1}>0}\|\Psi^{nA_{0}^{(0)}}_{\alpha}\|_{S([-T_{0},
T_{1}]\times\R^{2})}\leq C_{1}(\Ecrit)
\]
Finally, $\Psi^{nA_{0}^{(0)}}_{\alpha}$ has essential Fourier support contained in $(0,(\lambda_n^1)^{-1})$.
More precisely, for some sequence $\{R_n\}_{n=1}^\infty$ going to $\infty$ sufficiently slowly, one has
\begin{equation}\label{eq:firstPsitails}
 \| P_k \Psi^{nA_{0}^{(0)}}_{\alpha} \|_{S[k]} \le R_n^{-1} e^{-\sigma|k-\mu_n^1|}
\end{equation}
for all $k>\mu_n^1=-\log\lambda_n^1$ and some absolute constant
$\sigma$. As usual, all functions belong to the Schwartz class on
fixed time slices.
\end{prop}

The proof of this result will occupy this entire section.  The idea
is to run an induction in~$\ell$ on a sequence of approximating maps
with data $\Psi_\alpha^{\ell nA_0^{(0)}}$, see~\eqref{eq:Psielldef}.
As we start from the low frequencies, it will turn out that the
differences between two consecutive such approximating components is
of small energy (provided $\delta_1\gg\delta_0$ are both
sufficiently small). This allows us to pass from one approximation
to the next better one by applying a perturbative argument, albeit
with a linear operator involving a magnetic potential. Moreover, we
need to divide the time-axis into a number of intervals which is
controlled by the total energy. A key fact here which prevents
energy build-up as we pass from one time interval to the next, is
that these approximating components essentially preserve their
energy, see Corollary~\ref{cor:epsenergyconservation}. The
approximate energy conservation, in turn, comes from the fact that
these components are all essentially the Coulomb components of
suitable maps as demonstrated in
Lemma~\ref{lem:Truncatedwapproximation}.

\noindent For the remainder of this section we drop the superscript~$A_0^{(0)}$ from our notation since
we will limit ourselves entirely to the low frequency part. We begin by showing that (still at time $t=0$) the step from
$
\Psi^{\ell-1,n}_{\alpha}$
to
$
\Psi^{\ell,n}_{\alpha}
$
amounts to adding on a term of much larger frequency, up to small errors in energy.

\begin{lemma}\label{lem:WinducPrep}
One has
\begin{align*}
\Psi^{\ell, n}_{\alpha}-\Psi^{\ell-1,n}_{\alpha} &=\epsilon_{\alpha}^{\ell,n}=P_{J_{\ell}}\epsilon_{\alpha}^{\ell, n}+\tilde{\epsilon}_{\alpha}^{\ell, n}
\end{align*}
with
$
\|\tilde{\epsilon}_{\alpha}^{\ell,n}\|_{L_{x}^{2}}\les \delta_{1}$.
Furthermore,
\[
\Psi^{\ell-1,n}_{\alpha}=P_{\cup_{j\leq \ell-1}J_{j}}\Psi^{\ell-1,n}_{\alpha}+\tilde{\Psi}^{\ell-1,n}_{\alpha}
\]
with
$
 \|\tilde{\Psi}^{\ell-1,n}_{\alpha}\|_{L_{x}^{2}}\les \delta_{1}
$.
 Similar statements hold on the level of the $\Phi$-components.
\end{lemma}
\begin{proof} In view of Lemma~\ref{lem:Truncatedwapproximation} we may switch from $\Psi^{\ell,n}$ to the corresponding expressions involving~$w^n$.
For simplicity, write
\[
\sum_{j\leq \ell}P_{J_{j}}w^{n}_{\alpha}  \:e^{-i\Re\sum_{k=1,2}\triangle^{-1}\partial_{k}\sum_{j\leq \ell}P_{J_{j}}w^{n}_{k}} =: f_\ell e^{ig_\ell}
\]
with $g_\ell$ real-valued. Since the Fourier support of $f_\ell$ is contained in $\cup_{j\leq \ell}J_{j}=(-\infty,b_\ell]$,
for any $k\ge b_\ell+10$ one has
\begin{align*}
 \|P_k (f_\ell \, e^{ig_\ell})\|_2 &\les \|f_\ell P_{k+O(1)} \, e^{ig_\ell}\|_2 \les 2^{-k} \|f_\ell\|_2 \|\nabla P_{k+O(1)} e^{ig_\ell}\|_\infty\\
&\les 2^{-k} \|f_\ell\|_2 \|\Delta^{-1} D^2 f_\ell \|_\infty \les 2^{-k} \|f_\ell\|_2 \|\Delta^{-1} D^3 f_\ell \|_2 \les 2^{b_\ell-k} \|w^n\|_2^2\\
&\les \Ecrit 2^{b_\ell-k}
\end{align*}
where $\Ecrit$ controls the total energy, and thus also the $L^2$-norm of~$w^n$. By construction of $w^n_\alpha$, one has for any $L>0$
\[
 \limsup_{n\to\infty} \big\| P_{[b_\ell-L,b_\ell]} w^{n}_{\alpha} \big\|_2\les L\delta_0
\]
Together with the preceding bound this implies that
\begin{align*}
 &\Big\| \sum_{j\leq \ell}P_{J_{j}}w^{n}_{\alpha}  \:e^{-i\Re\sum_{k=1,2}\triangle^{-1}\partial_{k}\sum_{j\leq \ell}P_{J_{j}}w^{n}_{k}}
- P_{\cup_{j\leq \ell}J_{j}} \sum_{j\leq \ell}P_{J_{j}}w^{n}_{\alpha}  \:e^{-i\Re\sum_{k=1,2}\triangle^{-1}\partial_{k}\sum_{j\leq \ell}P_{J_{j}}w^{n}_{k}}\Big\|_2 \\
&\les \log[(\Ecrit+1)\delta_0^{-1}]\delta_0 \ll \delta_1
\end{align*}
for small $\delta_0$. Next, observe that
\begin{equation}\nonumber\begin{split}
&\sum_{j\leq \ell}P_{J_{j}}w^{n}_{\alpha}  \:e^{-i\Re\sum_{k=1,2}\triangle^{-1}\partial_{k}\sum_{j\leq \ell}P_{J_{j}}w^{n}_{k}}\\
&=\sum_{j\leq \ell-1}P_{J_{j}}w^{n}_{\alpha}  \:e^{-i\Re\sum_{k=1,2}\triangle^{-1}\partial_{k}\sum_{j\leq \ell}P_{J_{j}}w^{n}_{k}}+P_{J_{\ell}}w^{n}_{\alpha}  \:e^{-i\Re\sum_{k=1,2}\triangle^{-1}\partial_{k}\sum_{j\leq \ell}P_{J_{j}}w^{n}_{k}}\\
\end{split}\end{equation}
The first assertion of the lemma therefore follows from the following claims:
\begin{itemize}
 \item The function
\[
P_{J_{\ell}}w^{n}_{\alpha}  \:e^{-i\Re\sum_{k=1,2}\triangle^{-1}\partial_{k}\sum_{j\leq \ell}P_{J_{j}}w^{n }_{k}}
\]
has frequency support in $J_{\ell}=[a_{\ell}, b_{\ell}]$ up to exponentially decaying errors, and we also have
\[
\limsup_{n\to\infty} \big\|P_{J_{\ell}^{c}}\big[P_{J_{\ell}}w^{n }_{\alpha}  \:e^{-i\Re\sum_{k=1,2}\triangle^{-1}\partial_{k}\sum_{j\leq \ell}P_{J_{j}}w^{n }_{k}}\big] \big\|_{L_{x}^{2}}<\delta_{1}
\]
\item Furthermore, we have
\begin{equation}\nonumber\begin{split}
&\Big\|\sum_{j\leq \ell-1}P_{J_{j}}w^{n }_{\alpha}e^{-i\Re\sum_{k=1,2}\triangle^{-1}\partial_{k}\sum_{j\leq \ell}P_{J_{j}}w^{n }_{k}}-\sum_{j\leq \ell-1}P_{J_{j}}w^{n }_{\alpha}e^{-i\Re\sum_{k=1,2}\triangle^{-1}\partial_{k}\sum_{j\leq \ell-1}P_{J_{j}}w^{n }_{k}} \Big\|_{L_{x}^{2}}<\delta_{1}
\end{split}\end{equation}
for $n$ large enough.
\end{itemize}
As for the first claim, note that we have already dealt with the case of frequencies larger than $b_\ell$. Thus, assume that $j\le a_\ell-10$ and estimate
\begin{align*}
 \|P_j (P_{J_\ell} w^n \: e^{ig_\ell})\|_2 &\les 2^{\frac{j}{3}} \| P_{J_\ell} w^n\: e^{ig_\ell}\|_{\frac32}
\les 2^{\frac{j}{3}} \| P_{J_\ell} w^n\|_2 \| P_{J_\ell+O(1)} e^{ig_\ell}\|_{6}  \\
&\les   2^{\frac{j}{3}} \| P_{J_\ell} w^n\|_2 \sum_{k\in J_\ell+O(1)} 2^{-\frac{k}{3}} \| P_{k} \nabla e^{ig_\ell}\|_{2}
\les \Ecrit\, 2^{\frac{j-a_\ell}{3}}
\end{align*}
Furthermore, as before one can ``fudge at the edges'' meaning
\[
 \limsup_{n\to\infty} \big\| P_{[a_\ell,a_\ell+L]} w^{n}_{\alpha} \big\|_2\les L\delta_0
\]
which concludes the first claim.
For the second claim we need to show
\[
\Big\|\sum_{j\leq \ell-1}P_{J_{j}}w^{n }_{\alpha}  \:e^{-i\Re\sum_{k=1,2}\triangle^{-1}\partial_{k}\sum_{j\leq \ell-1}P_{J_{j}}w^{n }_{k}}\big(1-  \:e^{-i\Re\sum_{k=1,2}\triangle^{-1}\partial_{k}P_{J_{\ell}}w^{n }_{k}}\big) \Big\|_{L_{x}^{2}}\lesssim\delta_{1}
\]
where the implied constant is absolute (not depending on any of the other parameters).
However, this follows easily from the frequency localization up to exponentially decaying errors of
\[
\sum_{j\leq \ell-1}P_{J_{j}}w^{n }_{\alpha}  \:e^{-i\Re\sum_{k=1,2}\triangle^{-1}\partial_{k}\sum_{j\leq \ell-1}P_{J_{j}}w^{n }_{k}}
\]
as well as the fact that
\[
\limsup_{n\to\infty}\|P_{[a_\ell,a_{\ell}+L]\cup[b_\ell-L, b_{\ell}]}P_{J_{\ell}}w^{n }_{k}\|_{L_{x}^{2}}\les L\delta_{0}
\]
and we are done. The claim of the lemma about $\Phi$ is easier since it does not involve any phases, cf.~Lemma~\ref{lem:Truncatedwapproximation}
and Lemma~\ref{lem:wapproximation}.
\end{proof}

Our strategy now is to inductively control the nonlinear evolution of the
$\Psi^{\ell,n}_{\alpha}$, the Coulomb components of the approximation maps, starting with $\ell=1$ . At
each induction step we add a term $\epsilon_{\alpha}^{\ell, n}$ of energy less than~$\eps_{0}$.
The key then is the following perturbative result.  Recall that $\eps_0>0$ is a small constant which determines the perturbative
energy-cutoff (it depends on~$\Ecrit$).

\begin{prop}\label{PsiInduction} Let $\Psi^{\ell,n  }_{\alpha}$, $\epsilon_{\alpha}^{\ell, n  }$, be as before, with
$1\le \ell \le C_{1}(\Ecrit, \eps_{0})$. Also, let
\[
c_{k}^{(\ell-1)}:=\max_\alpha (\sum_{r\in \Z}2^{-\sigma
|r-k|}\|P_{r}P_{\cup_{j\le \ell-1} J_j}\Psi^{\ell-1,n
}_{\alpha}\|_{L_{x}^{2}}^{2})^{\frac{1}{2}}
\]
for some small enough constant $\sigma>0$ (an apriori constant). We
now make the following induction hypotheses, valid for all large
$n$: there is a decomposition $\Psi^{\ell-1,n  }_{\alpha} =
\tilde\Psi^{\ell-1,n  }_{\alpha} + \breve\Psi^{\ell-1,n }_{\alpha}$
so that
\begin{align}
\max_\alpha \|P_{k}\tilde\Psi^{\ell-1,n  }_{\alpha}\|_{S[k]([-T_{0},
T_{1}]\times\R^{2})}&< C_{2}\,c_{k}^{(\ell-1)}\label{eq:C2envelope} \\
\|\breve\Psi^{\ell-1,n }_{\alpha}\|_{S} &< C_2\, \delta_1
\label{eq:Snorm}
\end{align}
for some positive number $C_{2}$.

\noindent Then there exists a partition $\Psi^{\ell,n  }_{\alpha} =
\tilde\Psi^{\ell,n  }_{\alpha} + \breve\Psi^{\ell,n }_{\alpha}$ so
that
\begin{align}
\max_\alpha \|P_{k}\tilde\Psi^{\ell,n  }_{\alpha}\|_{S[k]([-T_{0},
T_{1}]\times\R^{2})}&< C_{3}\,c_{k}^{(\ell)}\label{eq:C2envelope2} \\
\|\breve\Psi^{\ell,n }_{\alpha}\|_{S} &< C_3 \,\delta_1
\label{eq:Snorm2}
\end{align}
provided $\delta_1<\delta_{1}^{0}=\delta_{1}^{0}(C_{2})$ and
provided  $n$ is sufficiently large. Here $C_3=C_3(C_2,\Ecrit)$.
\end{prop}

 It is important to note that we iterate
Proposition~\ref{PsiInduction} $O(\frac{C_{1}(\Ecrit,
\eps_{0})}{\eps_{0}})$ many times, obtaining the induction start
from the small data result of~\cite{Krieger}. It is  clear that
there is some  constant $\delta_{11}>0$ (depending only on~$\Ecrit$)
such that choosing $\delta_{1}<\delta_{11}$ in each step, this
proposition can be applied. This $\delta_{11}>0$ dictates our choice
of $A^{(0)}$ in the decomposition
\[
w_{\alpha}^{n  }=\sum_{j}\phi_{\alpha}^{na_{j}^{(0)}}+w_{\alpha}^{nA^{(0)}}
\]
from before, see~\eqref{eq:wfurtherrefine}.  Another essential feature of the construction is that
\begin{equation}
 \label{eq:PsiellS}
\|\Psi^{\ell,n  }\|_S \le K(\Ecrit)
\end{equation}
where $K$ is some rapidly growing function of the energy. This follows immediately
from the inductive nature of the proof and the fact that the number of steps
is controlled by the energy alone. However, it is crucial to the argument that
we do not have to make $\eps_0$ small depending on the function~$K(\Ecrit)$ as we
go through the inductive process. In other words, we have to make sure that one can fix~$\eps_0$
throughout.

\noindent
 The idea of the proof of Proposition~\ref{PsiInduction}
is as follows: under the assumptions~\eqref{eq:C2envelope}
and~\eqref{eq:Snorm} we can find time intervals $I_{1}, I_{2},
\ldots, I_{M_{1}}$, $M_{1}=M_{1}(\tilde{C}_{2})$ as in
Section~\ref{sec:perturb}, such that locally on~$I_j$,
\[
\Psi^{\ell-1,n}=\Psi^{\ell-1,n}_{L}+\Psi^{\ell-1,n}_{NL}
\]
Here $\psi_L$ is a linear wave and $\psi_{NL}$ is small in a
suitable sense,  see Lemma~\ref{lem:LocalSplitting} and
Corollary~\ref{cor:localsplit2}. In order to control the evolution
of $\Psi^{\ell,n}_{\alpha}$, we need to control the evolution of
\[
\epsilon_{\alpha}^{\ell,n}=\Psi^{\ell,n}_{\alpha}-\Psi^{\ell-1,n}_{\alpha}
\]
This we do inductively, over each interval $I_{j}$, starting with the one containing the initial time slice $t=0$.

\noindent
At this point one encounters the danger that the energy of $\epsilon_\alpha^{\ell,n}$ keeps growing as we move to later (or earlier)
intervals $I_{j}$, thereby effectively leaving the perturbative regime. The idea here is that we have {\em{ apriori energy conservation
for the components $\Psi^{\ell-1,n}$, $\Psi^{\ell,n}$}}, while at the same time, due to our assumptions on the frequency distribution of
energy for $\Psi^{\ell-1,n}$, $\epsilon_\alpha^{\ell,n}$, {\em{there cannot be much energy transfer between the latter two types of components}};
more precisely, we can enforce this by choosing $\delta_{1}$ small enough.
This means that we have effectively {\em{approximate energy conservation for $\epsilon_\alpha^{\ell,n}$}},
whence the induction can be continued to all the $I_{j}$.  We can now begin the proof in earnest.

\begin{figure*}[ht]
\begin{center}
\centerline{\hbox{\vbox{ \epsfxsize= 10.0 truecm \epsfysize=6.5
truecm \epsfbox{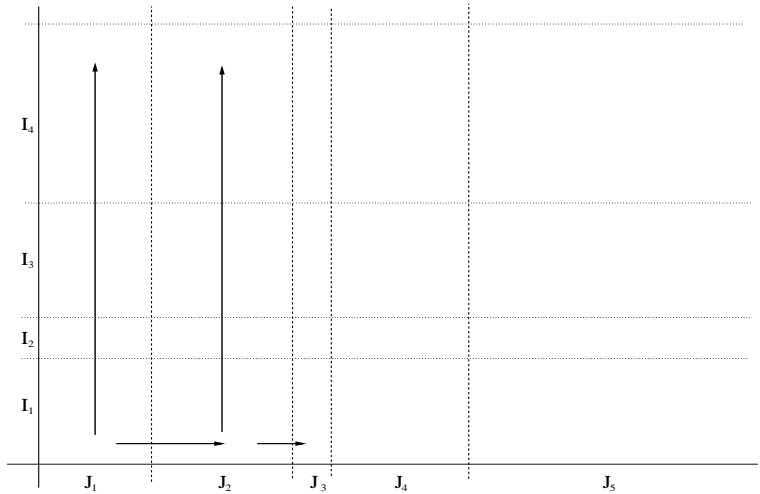}}}} \caption{The two directions of the
induction}
\end{center}
\end{figure*}

\begin{proof}(Proposition~\ref{PsiInduction}) We inductively control the nonlinear evolution of $\epsilon_\alpha^{\ell,n}$.
For ease of notation, we set
$\epsilon_\alpha:=\epsilon_\alpha^{\ell,n}$ and
$\psi_\alpha:=\Psi^{\ell-1,n}_{\alpha}$ and for the most part we
also ignore the~$\alpha$ subscript. Note that while $\psi$ exists
globally in time, $\epsilon$ exists only locally in time but we will
of course need to prove global existence and bounds for~$\epsilon$.
But for now, any statement  we make for~$\epsilon$ will be locally
in time on some interval~$I_0$ around $t=0$. Applying the
fungibility statements Lemma~\ref{lem:LocalSplitting} and
Corollary~\ref{cor:localsplit2} to~$\psi$ generates a decomposition
of~$\R$ into intervals~$\{I_j\}_{j=1}^M$ where
$M=M(\eps_0,\|\psi\|_S)$. We may of course intersect these intervals with~$I_0$ which we will
tacitly assume.  Fix one of these intervals, say $I_1$,
which contains~$t=0$. It will of course be necessary for us to pass to later intervals
in the temporal sense until we have exhausted the entire existence interval~$I_0$.
In other words, our induction has two direction, namely a temporal one (referring to the interval~$I_j$),
as well as a frequental one (referring to the interval $J_\ell$). These two directions are indicated
as vertical and horizontal ones, respectively, in Figure~6.

\noindent
By construction, there is a decomposition
\begin{equation}\label{eq:psiLNL}
\psi=\psi_{L}+\psi_{NL}
\end{equation}
where $\|\psi_L\|_S\le \eps_2^{-\frac14} E^2_{\mathrm{crit}}$ and
such that $\|\psi_{NL}\|_S <\eps_2$,
\begin{equation}\label{eq:prod_small}
\|\psi_{NL}\|_S\|\psi_L\|_S < \sqrt{\eps_2} \end{equation} Here
$\eps_2$ is small depending on~$\Ecrit$ and with $\eps_0\ll
\eps_2\ll1$. We note the following important improvement
over~\eqref{eq:PsiellS}:
\begin{equation}
 \label{eq:PsiellS2}
\max_{j}\|\Psi^{\ell,n  }\|_{S(I_j)} \les \eps_2^{-\frac14} E^2_{\mathrm{crit}}
\end{equation}
Proposition~\ref{PsiInduction} will follow from a bootstrap
argument, which is based on the following crucial result. Recall
that $J_\ell$ is the Fourier support of $\epsilon(0)$ up to errors
which can be made arbitrarily small in energy.

 \begin{prop}\label{PsiBootstrap} Let $\psi$ satisfy the inductive
assumptions~\eqref{eq:C2envelope} and~\eqref{eq:Snorm} and let
$\epsilon$ be  defined as above.
 Suppose there is a decomposition
 $
 \epsilon=\epsilon_{1}+\epsilon_{2}
 $
which satisfies the bounds
 \begin{equation}\label{eq:bootass}\begin{aligned}
 \|\epsilon_{2}\|_{S(I_1\times\R^{2})}&<C_2C_{4}\,\delta_{1} \\
 \|P_{k}\epsilon_{1}\|_{S[k](I_1\times\R^{2})} &\leq C_{4} d_{k}\quad
\forall\;k\in\Z
 \end{aligned}
\end{equation}
 where we define
 \[
 d_{k}:=(\sum_{r\in\Z}2^{-\sigma|r-k|}\|P_{r}P_{J_\ell}\epsilon(0,\cdot)\|_{L_{x}^{2}}^{2})^{\frac{1}{2}}
  \]
 for some $C_{4}=C_{4}(\Ecrit)$ sufficiently large, and some small absolute constant $\sigma>0$.
 Then we can improve this to a similar decomposition with
 \begin{equation}
  \label{eq:improvboot}
 \|\epsilon_{2}\|_{S(I_1\times\R^{2})}<\frac{C_{4}}{2}\,C_2\,\delta_{1},\qquad  \|P_{k}\epsilon_1\|_{S[k](I_1\times\R^{2})}\leq \frac{C_{4}}{2} d_{k}
 \end{equation}
for all $k\in\Z$.
 \end{prop}

This proposition is the key ingredient in the proof.  It asserts
that the frequency profile of~$\epsilon$ at time $t=0$ is
essentially preserved under the evolution up to some frequency
leakage, which however is controlled by the size of the {\em
underlying Besov error.} What allows us to prevent energy
of~$\epsilon$ moving from high to low frequencies (which is the main
difficulty here) are gains in the high-high-low interactions in the
nonlinearities. Without these gains, there could indeed be this kind
of energy transfer and the argument would break down. It is essential in Proposition~\ref{PsiBootstrap}
that~$C_4$ is a constant that {\em does not} change throughout the induction, whereas~$C_2$
{\em does} change.

\smallskip
\noindent
If we accept Proposition~\ref{PsiBootstrap} for now, then it is an easy matter to derive the aforementioned
approximate energy conservation.

\begin{cor}\label{cor:epsenergyconservation}
Under the induction hypothesis of Proposition~\ref{PsiInduction} and assuming the validity of Proposition~\ref{PsiBootstrap}, one
has the following:
For sufficiently small $\delta_{1}$ (depending on $C_2$ and $C_4$) and large $n$, we have
\[
\sum_{\alpha=0,1,2}\|\epsilon_{\alpha}\|_{L_{t}^{\infty}L_{x}^{2}(I_1\times\R^{2})}^{2}<\eps_{0}
\]
where $I_1$ is as above.
\end{cor}
\begin{proof}(Corollary~\ref{cor:epsenergyconservation}) Due to energy conservation for the evolution of $\psi+\epsilon$, we have
\[
\sum_{\alpha=0,1,2}\|\psi_{\alpha}+\epsilon_{\alpha}\|_{L_{x}^{2}}^{2}=\text{constant}
\]
Similarly, we have
\[
\sum_{\alpha=0,1,2}\|\psi_{\alpha}\|_{L_{x}^{2}}^{2}=\text{constant}
\]
The crucial observation now is that
\[
\|\psi+\epsilon\|^2_{L_{x}^{2}}=\|\psi\|_{L_{x}^{2}}^{2}+\|\epsilon\|_{L_{x}^{2}}^{2}+2\Re\sum_{k\in\Z}\int_{\R^{2}}P_{k}\psi\overline{P_{k}\epsilon}\, dx
\]
on fixed time slices $t=t_{0}\in I_{1}$, and we can split
\[
\sum_{k\in\Z}\int_{\R^{2}}P_{k}\psi\overline{P_{k}\epsilon}\, dx
=\sum_{k\in\cup_{j\leq l-1}J_{j}}\int_{\R^{2}}P_{k}\psi\overline{P_{k}\epsilon}\, dx+\sum_{k\in J_{l}}\int_{\R^{2}}P_{k}\psi\overline{P_{k}\epsilon}\, dx
\]
Both contributions on the right are $\lesssim C_{4}C_{2}^2\delta_{1}$, which can be made arbitrarily small  by choosing $\delta_{1}$ small enough.
To obtain this bound, observe that the induction hypothesis and Proposition~\ref{PsiBootstrap}
allow one to transfer Lemma~\ref{lem:WinducPrep} to all times in the interval~$I_1$.
Cauchy-Schwarz then implies the bound of~$\lesssim C_{4}C_{2}^2\delta_{1}$.
\end{proof}

Corollary~\ref{cor:epsenergyconservation} allows us to keep the energy under control as we inductively pass from $I_1$ to its successor~$I_2$ and so forth
by restarting the procedure.
Indeed, since the number of the ``fungibility'' intervals is bounded by $M(\eps_0,\Ecrit)$, we can make $\delta_1$ in the corollary so small
(depending on this number) and $n$ so large that even the energy of the final $\epsilon$ is no bigger than $2\eps_0$, say. Even though
we will now work on $I_1$, all arguments carried out below apply to any of the later intervals $I_2, I_3,\ldots$ as well.

\medskip
\begin{proof}[Proof of Proposition~\ref{PsiBootstrap}]
We may reduce ourselves to proving the statement for frequency
$2^{0}$, i.e.,  $k=0$, by scaling invariance. Recall that we have
chosen the intervals $I_{j}$ in such fashion that~\eqref{eq:psiLNL} holds with the stated bounds.
In order to obtain the desired estimates on~$\epsilon$, we distinguish between two cases, depending on the size
of the underlying time interval. If it is short, we use the div-curl system. Otherwise we use the wave equation.

\medskip
\noindent
{\em Case 1:}  $|I_{1}|<T_{1}$ where $T_{1}>0$ is some absolute small constant (to be specified).
We shall use the div-curl system  {\em linearized} around~$\psi$,
see~\eqref{eq:psisys1}, \eqref{eq:psisys2}, which takes the schematic form
 \[
 \partial_{t}\epsilon=\nabla_{x}\epsilon+\epsilon\nabla^{-1}(\psi^{2})+\psi\nabla^{-1}(\psi\epsilon)
+\epsilon\nabla^{-1}(\psi\epsilon)+\psi\nabla^{-1}(\epsilon^{2})+\epsilon\nabla^{-1}(\epsilon^{2})
 \]
The first linear term $\nabla_x\epsilon$ on the right-hand side is
estimated by bootstrap, choosing $T_1$ smaller than some absolute
constant. For each of the five nonlinear terms on the right-hand
side one needs to consider two cases, depending on whether
$\epsilon$ gets replaced by~$\epsilon_1$ or~$\epsilon_2$.

\noindent{(a)} The  term $\epsilon_1\nabla^{-1}(\psi^{2})$;  we
cannot just use Lemma~\ref{lem:easytrilin}
of~Section~\ref{sec:perturb}, since smallness there can only be
enforced by choosing $T_1$ very small, which is counter productive
in Case~2, when we work on a larger interval. Hence we have to
exploit the fungibility of the expression, which  forces us to
exploit the hidden null-structure. However, we can easily conclude
from the proof of Lemma~\ref{lem:easytrilin} that
\[
\|P_{0}[\epsilon_1\nabla^{-1}P_{<-C}(\psi^{2})]\|_{L_t^M L_x^2}\ll d_0
\]
provided we pick $C=C(\Ecrit)$ sufficiently large, and thence
\[
\|\int_{0}^{t}P_{0}[\epsilon_1\nabla^{-1}P_{<-C}(\psi^{2})]\,ds\|_{L_t^\infty L_x^2}\ll d_0
\]
\[
\|\int_{0}^{t}P_{0}[\epsilon_1\nabla^{-1}P_{<-C}(\psi^{2})]\,ds\|_{L_t^2 L_x^2}\ll d_0
\]
for $t\in [-T_1, T_1]$, and from there
\[
\|\int_{0}^{t}P_{0}[\epsilon_1\nabla^{-1}P_{<-C}(\psi^{2})]\,ds\|_{S[0]}\ll d_0,
\]
compare \eqref{eq:claim63} (provided $T_1<1$, say). Similarly, one checks that the contribution of
\[
 P_{0}[\epsilon_1\nabla^{-1}P_{>C}(\psi^{2})]
\]
is acceptable, and so we now need to force smallness for
\[
P_{0}[\epsilon_1\nabla^{-1}P_{[-C, C]}(\psi^{2})],
\]
which we do by subdivision into small time intervals (whose number depends on $\|\psi\|_{S}$).
First, we observe that choosing $C_1$ large enough depending on $C$ and $\Ecrit$, we can force that
\[
\|P_{0}[\epsilon_1\nabla^{-1}P_{[-C, C]}(Q_{>C_1}\psi\psi)]\|_{L_{t,x}^{2}}\ll d_0,
\]
and from here one can again infer that
\[
\int_{0}^{t}P_{0}[\epsilon_1\nabla^{-1}P_{[-C, C]}(Q_{>C_1}\psi\psi)]\,ds\|_{S[0]}\ll d_0
\]
for $t\in [-T_1, T_1]$, $T_1<1$, say. The same applies to
\[
P_{0}[Q_{>C_1}\epsilon_1\nabla^{-1}P_{[-C, C]}(\psi^{2})]
\]
Hence we may reduce to considering
\[
P_{0}[\epsilon_1\nabla^{-1}P_{[-C, C]}(\psi^{2})]
\]
where we automatically assume that $\psi=Q_{<C_1}\psi$, $\epsilon_1=Q_{<C_1}\epsilon_1$. Now we implement the customary Hodge decomposition
\[
\psi_\nu=R_\nu\psi+\chi_\nu
\]
First, substitute the gradient term for either factor $\psi$, which results in the expression
\[
P_{0}[\epsilon_1\nabla^{-1}P_{[-C, C]}\calN_{\nu j}(\psi, \psi)]
\]
Now due to Lemma~\ref{lem:Nablowmod} etc that in case of high-low or low-high interactions inside $\calN_{\nu j}(\psi, \psi)$  we can estimate
\[
\|P_{[-C, C]}\calN_{\nu j}(\psi, \psi)\|_{L_{t,x}^{2}}\lesssim \|\psi\|_{S}^{2}
\]
and one may then pick time intervals $I_j$ with the property that
\[
\sum_{k\in\Z}\|\chi_{I_j}P_{[k-C, k+C]}\calN_{\nu j}(\psi, \psi)\|_{L_{t,x}^{2}}^2\ll 1
\]
which ensures ``fungibility''. Thus it remains to deal with the
expression
\[
P_{0}[\epsilon_1\nabla^{-1}P_{[-C, C]}\calN_{\nu j}(P_{>C_2}\psi, P_{>C_2}\psi)]
\]
and indeed in light of Lemma~\ref{lem:Nablowmod} etc only the case
when $\nu=0$ needs to be considered. We choose $C_2\gg \max\{C,
C_1\}$.  Note that in this case the inner null-form may have very
large modulation (comparable to the frequency of the inputs), in
which case we cannot take advantage of the null-structure. The idea
then is to use the smoothing effect of integration over time.
Specifically, we write schematically
\begin{align}
&P_{0}[\epsilon_1\nabla^{-1}P_{[-C, C]}\calN_{\nu j}(P_{>C_2}\psi, P_{>C_2}\psi)]\nonumber\\
&=P_{0}[\epsilon_1\nabla^{-1}P_{[-C, C]}\partial_t(P_{>C_2}|\nabla|^{-1}\psi P_{>C_2}R_j\psi)]\nonumber\\
&-P_{0}[\epsilon_1\nabla^{-1}P_{[-C, C]}\partial_j(P_{>C_2}|\nabla|^{-1}\psi P_{>C_2}R_0\psi)]\nonumber\\
&=P_{0}\partial_t[\epsilon_1\nabla^{-1}P_{[-C, C]}(P_{>C_2}|\nabla|^{-1}\psi P_{>C_2}R_j\psi)]\label{eq:decomp401}\\
&-P_{0}[\partial_t\epsilon_1\nabla^{-1}P_{[-C, C]}(P_{>C_2}|\nabla|^{-1}\psi P_{>C_2}R_j\psi)]\label{eq:decomp402}\\
&-P_{0}[\epsilon_1\nabla^{-1}P_{[-C, C]}\partial_j(P_{>C_2}|\nabla|^{-1}\psi P_{>C_2}R_0\psi)]\label{eq:decomp403}
\end{align}
Now it is straightforward to analyze the contribution of each term, keeping in mind our assumptions about hyperbolicity of each input.
For the contribution of \eqref{eq:decomp401}, note that we have
\begin{align*}
&\int_0^t P_{0}\partial_s[\epsilon_1\nabla^{-1}P_{[-C, C]}(P_{>C_2}|\nabla|^{-1}\psi P_{>C_2}R_j\psi)]\,ds\\
&=P_{0}[\epsilon_1\nabla^{-1}P_{[-C, C]}(P_{>C_2}|\nabla|^{-1}\psi P_{>C_2}R_j\psi)](t,\,\cdot)\\
&-P_{0}[\epsilon_1\nabla^{-1}P_{[-C, C]}(P_{>C_2}|\nabla|^{-1}\psi P_{>C_2}R_j\psi)](0,\,\cdot)\\
\end{align*}
and we can then crudely bound (assuming $T_1<1$, say)
\begin{align*}
&\|\chi_{[-T_1,T_1]}\big[P_{0}[\epsilon_1\nabla^{-1}P_{[-C, C]}(P_{>C_2}|\nabla|^{-1}\psi P_{>C_2}R_j\psi)](t,\,\cdot)\\
&-P_{0}[\epsilon_1\nabla^{-1}P_{[-C, C]}(P_{>C_2}|\nabla|^{-1}\psi P_{>C_2}R_j\psi)](0,\,\cdot)\big]\|_{L_{t,x}^{2}}\ll d_0
\end{align*}
This again suffices for the bootstrapping.
\\
Next, for the expression \eqref{eq:decomp402}, we estimate it by
\begin{align*}
&\|\chi_{[-T_1,T_1]}P_{0}[\partial_t\epsilon_1\nabla^{-1}P_{[-C, C]}(P_{>C_2}|\nabla|^{-1}\psi P_{>C_2}R_j\psi)]\|_{L_{t,x}^{2}}\\
&\lesssim \|\partial_t\epsilon_1\|_{L_t^\infty L_x^2}\|\nabla^{-1}P_{[-C, C]}(P_{>C_2}|\nabla|^{-1}\psi P_{>C_2}R_j\psi)\|_{L_t^\infty L_x^2}\\
&\ll d_0
\end{align*}
Finally, expression \eqref{eq:decomp403} is more of the same (due to the hyperbolicity of the inputs) and omitted.\\
The corresponding estimate for $\epsilon_2\nabla^{-1}(\psi^{2})$ is essentially the same, the only difference
being that one square-sums over the frequencies at the end.

\smallskip
\noindent{(b)} The  term $\psi\nabla^{-1}(\psi\epsilon_1)$ as well as $\psi\nabla^{-1}(\psi\epsilon_2)$ both will be placed
in the $\epsilon_2$ component, meaning that we will prove that they have small~$S$-norm.
We start with $\epsilon_1$. We claim that
\begin{equation}
 \label{eq:fallb}
\| P_0 [ \psi\nabla^{-1}(P_{k_2} \psi P_{k_3}\epsilon_1)] \|_{L^M_t L^2_x} \les  2^{-\sigma_0 |k_2-k_3|} \|P_{k_2}\psi\|_{S[k_2]}
 \| P_{k_3}\epsilon_1\|_{S[k_3]} \;\sup_{k_1\in\Z} 2^{-\sigma_0 |k_1|} \|P_{k_1}\psi\|_{S[k_1]}
\end{equation}
for some $\sigma_0>0$. This follows by inspecting the proof of Lemma~\ref{lem:easytrilin}. If $|k_2-k_3|>B|\log\delta_1|$ where $B$ is large,
one concludes from~\eqref{eq:fallb} that
\[
 \sum_{ |k_2-k_3|>C|\log\delta_1| } \| P_0 [ \psi\nabla^{-1}(P_{k_2} \psi P_{k_3}\epsilon_1)] \|_{L^M_t L^2_x} \les C_4\,
 \delta_1^{B\sigma_0} \|\psi\|_S \|\epsilon_1(0)\|_2 \; \sup_{k_1\in\Z} 2^{-\sigma_0 |k_1|} \|P_{k_1}\psi\|_{S[k_1]}
\]
Replacing $P_0$ by $P_k$ and square summing in~$k$ yields a bound of
\[
 C_4\, \delta_1^{B\sigma_0} \|\psi\|_S^2 \|\epsilon_1(0)\|_2 \ll C_4C_2\, \delta_1
\]
for the contribution of this case. This can be done by choosing $B$ large depending on $\Ecrit$, see~\eqref{eq:PsiellS2}.
On the other hand, if $|k_2-k_3|\le B|\log\delta_1|$, then we exploit that
the Fourier supports of $\psi$ and $\epsilon$ are essentially disjoint up to small errors (bounded by $\les\delta_1$ in the~$S$-norm)
and exponentially decaying tails.
Now we sum~\eqref{eq:fallb} over this range to obtain
\begin{align}
\sum_{|k_2-k_3|\le B|\log\delta_1|} \| P_0 [ \psi\nabla^{-1}(P_{k_2} \psi P_{k_3}\epsilon_1)] \|_{L^M_t L^2_x} &\les
\sum_{|k_2-k_3|\le B|\log\delta_1|} \| P_0 [ \psi\nabla^{-1}(P_{k_2} \breve{\psi} P_{k_3}\epsilon_1)] \|_{L^M_t L^2_x} \label{eq:Psibreve} \\
&+\quad \sum_{|k_2-k_3|\le B|\log\delta_1|} \| P_0 [ \psi\nabla^{-1}(P_{k_2} \tilde \psi P_{k_3}\epsilon_1)] \|_{L^M_t L^2_x} \label{eq:Psitilde}
\end{align}
For \eqref{eq:Psibreve} one obtains as above
\[
 \eqref{eq:Psibreve} \les    \|\breve{\psi}\|_S \|\epsilon_1(0)\|_2 \; \sup_{k_1\in\Z} 2^{-\sigma_0 |k_1|} \|P_{k_1}\psi\|_{S[k_1]}
\]
with an absolute implicit constant. Replacing $P_0$ with $P_k$ and summing over all scales yields the bound
\[
 \les \|\psi\|_S \|\breve{\psi}\|_S  \|\epsilon_1(0)\|_2 \les \eps_2^{-\frac14}E^2_{\mathrm{crit}}\, C_2\,\delta_1\,\eps_0 \ll C_2C_4\,\delta_1
\]
provided we choose $\eps_2^{-\frac14}E^2_{\mathrm{crit}}\eps_0\ll C_4$.
Next, by the definition of the frequency envelopes $c_k$ and $d_k$,
\begin{align*}
 \eqref{eq:Psitilde} &\les \sup_{k_1\in\Z} 2^{-\sigma_0 |k_1|} \|P_{k_1}\psi\|_{S[k_1]} \;
 \sum_{|k_2-k_3|\le B|\log\delta_1|} 2^{-\sigma_0 |k_2-k_3|} \|P_{k_2} \tilde\psi\|_{S[k_2]} \|P_{k_3}\epsilon_1\|_{S[k_3]} \\
&\les \sup_{k_1\in\Z} 2^{-\sigma_0 |k_1|} \|P_{k_1}\psi\|_{S[k_1]} \;
 \sum_{k_2-k_3|\le B|\log\delta_1|} 2^{-\sigma_0 |k_2-k_3|} C_2 \, c_{k_2}^{(\ell-1)} \, C_4 d_{k_3}\\
&\les C_2C_4\, \delta_0\; \sup_{k_1\in\Z} 2^{-\sigma_0 |k_1|} \|P_{k_1}\psi\|_{S[k_1]} \;
\end{align*}
This follows from the fact that $\delta_0$ was chosen to control the Besov norm of $w_\alpha^{nA_0^{(0)}}$, as well as the fact that the intervals
$J_i$ where chosen in such a way that any of the smaller atoms contained within  $w_\alpha^{nA_0^{(0)}}$ are arbitrarily far away from the endpoints
of~$J_i$ as $n\to\infty$.
Rescaling this bound to~$P_k$ from $P_0$ and square summing yields a bound of $C_2 C_4\, \|\psi\|_S\delta_0 \ll C_2 C_4 \delta_1$ by taking~$\delta_0$
small enough, cf.~\eqref{eq:PsiellS2}.

Next, we turn to $\psi\nabla^{-1}(\psi\epsilon_2)$. Here the
smallness comes from ``fungibility'' again as in case (a). More
precisely, reasoning as in (a), we may reduce this expression to the
form
\[
P_0[\psi\nabla^{-1}P_{[-C,C]}(\psi\epsilon_2)]
\]
where we moreover have $\psi=Q_{<C_1}\psi$, $\epsilon_2=Q_{<C_1}\epsilon_2$. Again the argument from (a) shows that we may assume
both inputs of $P_{[-C,C]}(\psi\epsilon_2)$ to have frequency $O(1)$ (implied constant depending on $C, C_1$, and $\Ecrit$).
Furthermore, it is straightforward to check that if the two factors $\psi$ have closely aligned Fourier supports, we obtain the desired
smallness via Bernstein's inequality. But if the Fourier supports of the two $\psi$ have some angular separation, interpreting the operator
$\nabla^{-1}P_{[-C,C]}$ as convolution with a kernel $K(x)$ of bounded (although possibly large) $L^1$-mass, we may write
\[
P_0[\psi\nabla^{-1}P_{[-C,C]}(\psi\epsilon_2)]=\int_{\R^{2}}P_0[\psi(\cdot, x)K(y)(\psi(\cdot, x-y)\epsilon_2)(\cdot, x-y)]\,dy
\]
and then
\[
\|\psi(\cdot, x)\psi(\cdot, x-y)\|_{L_{t,x}^{2}}\lesssim \|\psi\|_{S}^{2},
\]
which follows from our assumption about the Fourier supports, as well as the fact that both frequencies here are $<O(1)$.
But then we can again smallness by picking the $I_j$ suitably, such that
\[
\sum_{k\in\Z}\|\chi_{I_j}\int_{\R^{2}}P_k[P_{<k+O(1)}\psi(\cdot, x)|K(y)|P_{<k+O(1)}\psi(\cdot, x-y)]\,dy\|_{L_t^2\dot{H}^{-\frac{1}{2}}}^{2}\ll 1
\]

\smallskip
\noindent{(c)}  The term $\epsilon_1 \nabla^{-1} (\psi\epsilon_1)$ is easy, since it inherits the frequency profile of~$\epsilon_1$.
More precisely, using the same type of trilinear estimates as in (a) and (b) one obtains
\[
 \|P_k (\epsilon_1 \nabla^{-1} (\psi\epsilon_1))\|_{S[k](I_1\times\R^2)} \les C_4\,d_k \|\psi\|_{\ener} \|\epsilon_1\|_{\ener}  \ll C_4 \, d_k
\]
using \eqref{eq:PsiellS2} and the fact that $\|\epsilon_1\|_{\ener} \le 2\eps_0$ (taking $\delta_1$ small).
The other cases are easier due to the presence of~$\delta_1$ coming from~$\epsilon_2$.

\smallskip
\noindent{(d)} The term $\psi \nabla^{-1}(\epsilon^2)$ splits into the terms $\psi \nabla^{-1}(\epsilon_1^2)$,
$\psi \nabla^{-1}(\epsilon_1\,\epsilon_2)$, and $\psi \nabla^{-1}(\epsilon_2^2)$. The last two are easier due to the smallness of~$\epsilon_2$.
The first one is harder, as it inherits the frequency profile of~$\psi$ and therefore needs to be incorporated in~$\epsilon_2$.
This means that we need to gain the very small $\delta_0$, which is only possible if there are high-high gains in the
inner term of $\psi \nabla^{-1}(\epsilon_1^2)$ resulting from~$\epsilon_1$. Of course, this requires that we expand this inner expression
into a null-form via the usual Hodge decomposition.

{\it{(i): High-High-Low interactions in $\nabla^{-1}(\epsilon^{2})$.}} This is the following (schematic) type of term:
\[
 \sum_{k,\,k_{1,2,3},\,k\ll k_{2}}P_{0} [P_{k_{1}}\psi\nabla^{-1}P_{k}(P_{k_{2}}\epsilon P_{k_{3}}\epsilon)].
\]
It is straightforward to see that we may assume $|k|<\sigma_3 k_{2}$
for some $\sigma_3>0$ (absolute constant independent of the other
smallness parameters), and furthermore
$k_{2}=k_{3}+O(1)>B|\log\delta_{1}|$, since otherwise the desired
smallness follows as in the preceding Case~(b). We may thus
essentially assume $k_{1}=O(1),\,k=O(1)$, and reduce to the
simplified expression
\[
 \sum_{k_{1}=O(1)=k,\,k_{2}>B|\log\delta_{1}|}P_{0} [P_{k_{1}}\psi\nabla^{-1}P_{k}(P_{k_{2}}\epsilon P_{k_{3}}\epsilon)]
\]
Suppressing the frequency localizations for now,
we use the schematic relation
\begin{equation}\nonumber\begin{split}
P_{0}\big[\psi\nabla^{-1}(\epsilon^{2})\big]=&P_{0}\big[\psi\nabla^{-1}(R_{\nu}\epsilon^{1}R_{j}\epsilon^{2}-R_{j}\epsilon^{1}
R_{\nu}\epsilon^{2})+\psi\nabla^{-1}(\nabla^{-1}(\epsilon\nabla^{-1}(\epsilon^{2}))R_{\nu}\epsilon)+\ldots\\
&+\psi\nabla^{-1}(\nabla^{-1}\big([\epsilon\nabla^{-1}(\epsilon^{2})]^{2}\big)\big]+\ldots
\end{split}\end{equation}
where we omit the remaining quintilinear and septilinear terms.
More precisely, we shall use this provided both inputs $\epsilon$ have relatively small modulation, i.e., are of hyperbolic type.
Thus for $k_{2}=k_{3}+O(1)>B|\log\delta_{1}|$, we write
\begin{align}
P_{0}\big[ P_{k_{1}} \psi\nabla^{-1}(P_{k_{2}}\epsilon P_{k_{3}}\epsilon)\big]=&P_{0}\big[P_{k_{1}}\psi\nabla^{-1}
(P_{k_{2}}Q_{>k_{2}}\epsilon P_{k_{3}}\epsilon)\big]+P_{0}\big[P_{k_{1}}\psi\nabla^{-1}(P_{k_{2}}Q_{<k_{2}}\epsilon
P_{k_{3}}Q_{>k_{3}}\epsilon)\big] \label{eq:entw1}\\
&+P_{0}\big[P_{k_{1}}\psi\nabla^{-1}(R_{\nu}P_{k_{2}}Q_{<k_{2}}\epsilon R_{j}P_{k_{3}}Q_{<k_{3}}\epsilon-R_{j}P_{k_{2}}Q_{<k_{2}}\epsilon R_{\nu}P_{k_{3}}Q_{<k_{3}}\epsilon)\label{eq:entw2} \\
&+P_{k_{1}}\psi\nabla^{-1}(\nabla^{-1}P_{k_{2}}Q_{<k_{2}}(\epsilon\nabla^{-1}(\epsilon^{2}))R_{\nu}P_{k_{3}}Q_{<k_{3}}\epsilon) \label{eq:entw3} \\
&+P_{k_{1}}\psi\nabla^{-1}(\nabla^{-1}P_{k_{2}}Q_{<k_{2}}(\psi\nabla^{-1}(\epsilon^{2}))R_{\nu}P_{k_{3}}Q_{<k_{3}}\epsilon) \label{eq:entw4}\\
& +P_{k_{1}}\psi\nabla^{-1}(\nabla^{-1}P_{k_{2}}Q_{<k_{2}}(\epsilon\nabla^{-1}(\psi\epsilon))R_{\nu}P_{k_{3}}Q_{<k_{3}}\epsilon)\label{eq:entw5} \\
&+P_{k_{1}}\psi\nabla^{-1}(\nabla^{-1}P_{k_{2}}Q_{<k_{2}}(\epsilon\nabla^{-1}(\psi^{2}))R_{\nu}P_{k_{3}}Q_{<k_{3}}\epsilon)\label{eq:entw6}\\
&+P_{k_{1}}\psi\nabla^{-1}(\nabla^{-1}P_{k_{2}}Q_{<k_{2}}(\psi\nabla^{-1}(\psi\epsilon))R_{\nu}P_{k_{3}}Q_{<k_{3}}\epsilon)\label{eq:entw7} \\
&+P_{k_{1}}\psi\nabla^{-1}(\nabla^{-1}P_{k_{2}}Q_{<k_{2}}[\epsilon\nabla^{-1}(\epsilon^{2})]\nabla^{-1}P_{k_{3}}Q_{<k_{3}}
[\epsilon\nabla^{-1}(\epsilon^{2})])\big]+\ldots\label{eq:entw8}
\end{align}
where $\ldots$ denotes the remaining septilinear terms containing mixed $\psi$-$\epsilon$-interactions. Again we
may substitute $\epsilon_{1}$ everywhere for $\epsilon$, the contributions from $\epsilon_{2}$ leading to much smaller
contributions. The first two terms on the right are straightforward to estimate: using Bernstein's inequality, one obtains for \eqref{eq:entw1} the bound
\begin{align*}
&\|P_{0}\big[P_{k_{1}}\psi\nabla^{-1}(P_{k_{2}}Q_{>k_{2}}\epsilon_1 P_{k_{3}}\epsilon_1)\big]\|_{L_{t,x}^{2}}\lesssim
\min\{\|P_{k_{1}}\psi\|_{L_{t}^{\infty}L_{x}^{2}}, \|P_{k_{1}}\psi \|_{L_{t}^{\infty}L_{x}^{\infty}}\}
\|P_{k_{2}}Q_{>k_{2}}\epsilon_1\|_{L_{t,x}^{2}}\|P_{k_{3}}\epsilon_1\|_{L_{t}^{\infty}L_{x}^{2}}\\&\lesssim
2^{-\frac{k_{2}}{2}}\min\{\|P_{k_{1}}\psi_{1}\|_{L_{t}^{\infty}L_{x}^{2}}, \|P_{k_{1}}\psi_{1}\|_{L_{t}^{\infty}L_{x}^{\infty}}\}
\|P_{k_{2}}\epsilon_1\|_{S[k_{2}]}\|P_{k_{3}}\epsilon_1\|_{S[k_{3}]}
\end{align*}
Keep in mind here we assume $k_{1}=O(1)$. Then by an argument similar to the one used to estimate \eqref{eq:Psitilde},
replacing the output frequency by $2^{k}$ and square summing over $k=k_{1}+O(1)$ while also summing over $|k_{1}-k_{2}|>B|\log\delta_{1}|$,
one can bound this contribution by $\les C_{2}C_{4}^{2}\eps_0^{2}\delta_0$, which is enough to incorporate this term into $\epsilon_2$.
The second term in the expansion is of course handled identically, and so we now turn to the third term \eqref{eq:entw2},
which is the most delicate one. The potential difficulty comes when $\nu=0$, as the $Q_{\nu j}$-null-form allows us to pull
out one derivative otherwise; indeed, assume first that $\{\nu, j\}=\{1,2\}$. Then using the identity (and omitting the
subscript from $\epsilon$ for simplicity)
\[
R_{1}\epsilon^{1}R_{2}\epsilon^{2}-R_{1}\epsilon^{2}R_{2}\epsilon^{1}=\partial_{1}[\nabla^{-1}\epsilon^{1}R_{2}\epsilon^{2}]-
\partial_{2}[\nabla^{-1}\epsilon^{1}R_{1}\epsilon^{2}],
\]
we can estimate (always under the assumption $k_{1}=O(1)=k$)
\begin{align*}
&\|P_{0}\big[P_{k_{1}}\psi\nabla^{-1}P_{k}(R_{1}P_{k_{2}}Q_{<k_{2}}\epsilon R_{2}P_{k_{3}}Q_{<k_{3}}\epsilon-R_{2}
P_{k_{2}}Q_{<k_{2}}\epsilon R_{1}P_{k_{3}}Q_{<k_{3}}\epsilon)\|_{L_{t,x}^{2}}\\
&\les\|P_{k_{1}}\psi\|_{L_{t}^{\infty}L_{x}^{2}}\|P_{k}[\nabla^{-1}P_{k_{2}}Q_{<k_{2}}\epsilon R_{1,2}P_{k_{3}}Q_{<k_{3}}\epsilon]\|_{L_{t}^{2}L_{x}^{\infty}}
\end{align*}
In order to estimate the right-hand factor, we use the improved Strichartz estimates: we have
\begin{align*}
P_{k}[\nabla^{-1}P_{k_{2}}Q_{<k_{2}}\epsilon R_{1,2}P_{k_{3}}Q_{<k_{3}}\epsilon]=\sum_{\substack{c_{1,2}\in \calD_{k_{2}, -k_{2}}
\\ \text{dist}(c_{1},-c_{2})=O(1)}}P_{k}[\nabla^{-1}P_{c_1}Q_{<k_{2}}\epsilon R_{1,2}P_{c_2}Q_{<k_{3}}\epsilon]
\end{align*}
whence we get
\begin{align*}
 \|P_{k}[\nabla^{-1}P_{k_{2}}Q_{<k_{2}}\epsilon R_{1,2}P_{k_{3}}Q_{<k_{3}}\epsilon]\|_{L_{t}^{2}L_{x}^{\infty}}\lesssim
 &2^{-k_{2}}\big(\sum_{c\in \calD_{k_{2}, -k_{2}}}\|P_{k_{2}}Q_{<k_{2}}\epsilon\|_{L_{t}^{4}L_{x}^{\infty}}^{2}\big)^{\frac{1}{2}}\big(
 \sum_{c\in \calD_{k_{3}, -k_{3}}}\|P_{k_{3}}Q_{<k_{3}}\epsilon\|_{L_{t}^{4}L_{x}^{\infty}}^{2}\big)^{\frac{1}{2}}\\
&\les 2^{-\frac{k_2}{2+}}\prod_{j=2,3}\|P_{k_{j}}\epsilon\|_{S[k_j]},
\end{align*}
whence we now have
\begin{align*}
&\|P_{0}\big[P_{k_{1}}\psi\nabla^{-1}P_{k}(R_{1}P_{k_{2}}Q_{<k_{2}}\epsilon R_{2}P_{k_{3}}Q_{<k_{3}}\epsilon-R_{2}P_{k_{2}}Q_{<k_{2}}
\epsilon R_{1}P_{k_{3}}Q_{<k_{3}}\epsilon)\|_{L_{t,x}^{2}}\\
&\les \|P_{k_{1}}\psi\|_{L_{t}^{\infty}L_{x}^{2}}2^{-\frac{k_2}{2+}}\prod_{j=2,3}\|P_{k_{j}}\epsilon\|_{S[k_j]}
\end{align*}
From here one can again conclude as in case (b).
\\
Hence we now consider now the more difficult case where $\nu=0$. First, it is straightforward to check that we may reduce the first input
 $P_{k_{1}}\psi$ to modulation $<2^{\sigma_4 k_{2}}$, where for example we may put $\sigma_4=\frac{1}{2}$. Then we use the schematic representation
\begin{align*}
&P_{0}\big[P_{k_{1}}Q_{<\frac{k_{2}}{2}}\psi\nabla^{-1}P_{k}(R_{0}P_{k_{2}}Q_{<k_{2}}\epsilon R_{1}P_{k_{3}}Q_{<k_{3}}\epsilon-R_{1}
P_{k_{2}}Q_{<k_{2}}\epsilon R_{0}P_{k_{3}}Q_{<k_{3}}\epsilon)\\
&=P_{0}\partial_{t}[P_{k_{1}}Q_{<\frac{k_{2}}{2}}\psi\nabla^{-1}P_k[\nabla^{-1}P_{k_{2}}Q_{<k_{2}}\epsilon R_{1}P_{k_{3}}Q_{<k_{3}}\epsilon]]-
P_{0}[P_{k_{1}}Q_{<\frac{k_{2}}{2}}\partial_{t}\psi\nabla^{-1}P_k[\nabla^{-1}P_{k_{2}}Q_{<k_{2}}\epsilon R_{1}P_{k_{3}}Q_{<k_{3}}\epsilon]]\\
&-P_{0}[P_{k_{1}}Q_{<\frac{k_{2}}{2}}\psi\nabla^{-1}P_k R_1[\nabla^{-1}P_{k_{2}}Q_{<k_{2}}\epsilon R_{0}P_{k_{3}}Q_{<k_{3}}\epsilon]]
\end{align*}
If one then integrates the transport equation for $\epsilon$, the contribution from the above terms is
\begin{align*}
 P_{0}\epsilon(t, \cdot)=P_0\epsilon(0,\cdot)+&\int_{0}^{t}\nabla_{x}\epsilon(s,\cdot)\, dsP_{0}[P_{k_{1}}Q_{<\frac{k_{2}}{2}}
 \psi\nabla^{-1}P_k[\nabla^{-1}P_{k_{2}}Q_{<k_{2}}\epsilon R_{1}P_{k_{3}}Q_{<k_{3}}\epsilon]](t, \cdot)\\&-P_0[P_{k_{1}}
 Q_{<\frac{k_{2}}{2}}\psi\nabla^{-1}P_k[\nabla^{-1}P_{k_{2}}Q_{<k_{2}}\epsilon R_{1}P_{k_{3}}Q_{<k_{3}}\epsilon]](0,\cdot)\\
&-\int_{0}^{t}P_{0}[P_{k_{1}}Q_{<\frac{k_{2}}{2}}\partial_{t}\psi\nabla^{-1}P_k[\nabla^{-1}P_{k_{2}}Q_{<k_{2}}\epsilon R_{1}
P_{k_{3}}Q_{<k_{3}}\epsilon]](s, \cdot)\,ds\\
&-\int_{0}^{t}P_{0}[P_{k_{1}}Q_{<\frac{k_{2}}{2}}\psi\nabla^{-1}P_k R_1[\nabla^{-1}P_{k_{2}}Q_{<k_{2}}\epsilon R_{0}P_{k_{3}}Q_{<k_{3}}\epsilon]](s,\cdot)\,ds
\end{align*}
But under our current assumption $k_{1}=O(1)$, $k=O(1)$, we have the estimate (using Bernstein's inequality)
\begin{align*}
 &\|P_{0}[P_{k_{1}}Q_{<\frac{k_{2}}{2}}\psi\nabla^{-1}P_k[\nabla^{-1}P_{k_{2}}Q_{<k_{2}}\epsilon R_{1}P_{k_{3}}Q_{<k_{3}}\epsilon]](t, \cdot)\\
 &-P_0[P_{k_{1}}Q_{<\frac{k_{2}}{2}}\psi\nabla^{-1}P_k[\nabla^{-1}P_{k_{2}}Q_{<k_{2}}\epsilon R_{1}P_{k_{3}}Q_{<k_{3}}\epsilon]](0,\cdot)
 \|_{L_{t}^{\infty}L_{x}^{2}}\\
&\les 2^{-k_{2}}\|P_{k_{1}}\psi\|_{\ener}\|P_{k_{2}}\epsilon\|_{\ener}\|P_{k_{3}}\epsilon\|_{\ener}
\end{align*}
and the remaining integral expressions on the right also easily lead to exponential gains in $-k_{2}$ due to the extra $\nabla^{-1}$ applied to
 $P_{k_{2}}Q_{<k_{2}}\epsilon $. Our assumption $k_{2}>B|\log\delta_{1}|$ then allows us to incorporate the contribution of all these source
 terms into $\epsilon_{2}$. Note that the cutoff $Q_{<\frac{k_{2}}{2}}$ in front of $\partial_{t}\psi$ allows us to control the effect of the $\partial_{t}$.
\\
The remaining terms \eqref{eq:entw3}-\eqref{eq:entw8} no longer require an integration by parts trick and can be directly placed into $L_{t,x}^{2}$
 with the requisite gain in $k_{2}$. We treat here the term \eqref{eq:entw4} given by
\[
 P_{k_{1}}\psi\nabla^{-1}(\nabla^{-1}P_{k_{2}}Q_{<k_{2}}(\psi\nabla^{-1}(\epsilon^{2}))R_{\nu}P_{k_{3}}Q_{<k_{3}}\epsilon)
\]
where we always keep in mind the localizations $k_{1}=O(1)=k$, $k_{2}=k_{3}+O(1)>B|\log\delta_{1}|$. The key here is as before the improved
Strichartz estimates. Write
\[
 P_{k_{2}}Q_{<k_{2}}(\psi\nabla^{-1}(\epsilon^{2}))=P_{k_{2}}Q_{<k_{2}}(\psi\nabla^{-1}P_{<0}(\epsilon^{2}))+\sum_{s\geq 0}P_{k_{2}}
 Q_{<k_{2}}(\psi\nabla^{-1}P_s(\epsilon^{2}))
\]
We treat here the contribution of the second term on the right, the first being treated in the same vein. Now if $s<k_{2}-10$, we get
\[
\big(\sum_{c\in \calD_{k_{2}, s-k_{2}}}\|P_{c}Q_{<k_{2}}(\psi\nabla^{-1}P_s(\epsilon^{2}))\|_{L_{t}^{4}L_{x}^{1}}^{2})^{\frac{1}{2}}\les
2^{\frac{3 k_{2}}{4}}2^{\frac{s-k_{2}}{2+}}2^{-s}\|P_{k_{2}}\psi\|_{S[k_{2}]}\|\epsilon\|_{\ener}^{2}
\]
Thus in the case $s<k_{2}-10$ from Bernstein's inequality we get
\begin{align*}
&\|P_{k_{1}}\psi\nabla^{-1}(\nabla^{-1}P_{k_{2}}Q_{<k_{2}}(\psi\nabla^{-1}P_{s}(\epsilon^{2}))R_{\nu}P_{k_{3}}Q_{<k_{3}}\epsilon)\|_{L_{t,x }^{2}}\\
&=\sum_{\substack{c_{1,2}\in \calD_{k_{2}, s-k_{2}}\\\text{dist}(c_{1}, -c_{2})\les 2^{s}}}\|P_{k_{1}}\psi\nabla^{-1}(\nabla^{-1}
P_{c_{1}}Q_{<k_{2}}(\psi\nabla^{-1}P_{s}(\epsilon^{2}))R_{\nu}P_{c_{2}}Q_{<k_{3}}\epsilon)\|_{L_{t,x }^{2}}\\
&\les \|P_{k_{1}}\psi\|_{\ener}\big(\sum_{c\in \calD_{k_{2}, s-k_{2}}}\|P_{c_{1}}Q_{<k_{2}}(\psi\nabla^{-1}P_{s}(\epsilon^{2}))
\|_{L_{t}^{4}L_{x}^{1}}^{2}\big)^{\frac{1}{2}}\big(\sum_{c\in \calD_{k_{2}, s-k_{2}}}\|R_{\nu}P_{c_{2}}Q_{<k_{3}}\epsilon
\|_{L_{t}^{4}L_{x}^{\infty}}^{2})^{\frac{1}{2}}\\
&\les 2^{-k_{2}}2^{\frac{3 k_{2}}{2}}2^{2(\frac{s-k_{2}}{2+})}2^{-s}\|P_{k_{1}}\psi_{1}\|_{S[k_{1}]}\|P_{k_{2}}\psi\|_{S[k_{2}]}\|\epsilon\|_{\ener}^{2}
\end{align*}
Summing over $0<s<k_{2}$ results in the bound
\[
 \les 2^{-\frac{k_{2}}{2+}}\|P_{k_{1}}\psi_{1}\|_{S[k_{1}]}\|P_{k_{2}}\psi\|_{S[k_{2}]}\|\epsilon\|_{\ener}^{2}
\]
On the other hand, when $s\geq k_{2}-10$, we simply bound
\[
 \|P_{k_{2}}Q_{<k_{2}}(\psi\nabla^{-1}P_s(\epsilon^{2}))\|_{L_{t}^{4}L_{x}^{1}}\les 2^{-\frac{k_{1}}{4}}\|\psi\|_{S}\|\epsilon\|_{\ener}^{2}
\]
and from here one estimates the $L_{t,x}^{2}$-norm of the output as before but without using the improved Strichartz, just the
standard $L_{t}^{4}L_{x}^{\infty}$-bound. The remaining terms \eqref{eq:entw5} are handled similarly.
\\
{\it{(ii) : High-Low/ Low-High interactions within
$\nabla^{-1}(\epsilon^{2})$}} In this case one  gains exponentially
in the maximum frequency occurring among the two factors $\epsilon$,
provided this is much larger than $1$. In this case one can argue as
in case (b) to include this contribution into $\epsilon_{2}$.

\smallskip
\noindent{(e)} The cubic term $\epsilon\nabla^{-1}(\epsilon^2)$ is easy, and can be treated as in (a) and (b) above.

\smallskip\noindent
The bootstrap argument for $\epsilon$ in the small time case is now completed as in the proof of Lemma~\ref{lem:LocalSplitting}, cf.~\eqref{eq:claim63}.

\medskip
{\em Case 2:} $|I_{1}|\geq T_{1}$, where $T_{1}>0$ is a small
constant depending on $\Ecrit$. Here we have to work with the wave equation satisfied by $P_{0}\epsilon$. We start by recording this equation
schematically in its original trilinear form, to which we apply various Hodge type decompositions as well as localizations in frequency space.
The goal is to write the equation in the form of a nonlinear wave equation with a low-frequency magnetic potential term, which we
will treat as part of the linear operator.
To begin with, we have the schematic equation (here we suppress the fact that $\epsilon$ really stands for the system of
variables $\{\epsilon_{\alpha}\}$,\,$\alpha=0,1,2$)
\begin{align}
\label{eq:st1} \Box P_{0}\epsilon&=\nabla_{x,t}P_{0}\big[(\psi+\epsilon)\nabla^{-1}([\psi+\epsilon]^{2})\big]
-\nabla_{x,t}P_{0}\big[(\psi)\nabla^{-1}(\psi^{2})\big]\\
&=P_{0}\nabla_{x,t}[\epsilon\nabla^{-1}(\psi^{2})]+P_{0}\nabla_{x,t}[\psi\nabla^{-1}(\psi\epsilon)]+
P_{0}\nabla_{x,t}[\epsilon\nabla^{-1}(\psi\epsilon)]+P_{0}\nabla_{x,t}[\psi\nabla^{-1}(\epsilon^{2})]+
P_{0}\nabla_{x,t}[\epsilon\nabla^{-1}(\epsilon^{2})]\nn
\end{align}
More precisely,  the terms on the right-hand side of~\eqref{eq:st1}
are exactly those given by~\eqref{eq:psi_wave}. It is precisely the
first term on the last line which causes technical difficulties for
the bootstrap argument, and we shall have to include parts of it
into the linear operator. However, this will only be made specific
once we have localized the terms suitably in frequency space. To
begin with, note that we will implement a bootstrap argument in
order to deduce bounds on $\epsilon$. For this we substitute
Schwartz extensions $\tilde{\epsilon}_{\alpha}$ for each
$\epsilon_{\alpha}$ on the right-hand side (these extensions
agreeing with $\epsilon_{\alpha}$ on the time interval
$I_{1}\times\R^{2}$ we are working on), and then solve the
inhomogeneous wave equation for $\epsilon_{\alpha}$, improving the
bounds we used for $\tilde{\epsilon}_{\alpha}$. Denoting the
right-hand source term above --- with $\tilde{\epsilon}_{\alpha}$
instead of $\epsilon_{\alpha}$ --- by $\tilde{F}_{\alpha}$, what we
really do is solving the problem
\[
 \Box P_{0}\epsilon_{\alpha}=P_{0}\tilde{F}_{\alpha}
\]
In order to deduce the $S$-bounds on $P_{0}\epsilon_{\alpha}$, we split this variable into two parts
\[
 P_{0}\epsilon_{\alpha}=P_{0}Q_{\geq D}\epsilon_{\alpha}+P_{0}Q_{<D}\epsilon_{\alpha}
\]
Here the parameter $D$ is chosen sufficiently large depending on
$T_{1}$ from Case~1 and thus depends on $\Ecrit$ (but is independent
of the induction stage). Then we solve the preceding wave equation
by setting
\begin{align*}
P_{0}Q_{\geq D}\epsilon_{\alpha} &=\Box^{-1}Q_{\geq
D}P_{0}\tilde{F}_{\alpha}\\
 P_{0}Q_{<D}\epsilon_{\alpha}
&=S(t)\big(P_{0}Q_{<D}\epsilon_{\alpha}\big)[0]+\int_{0}^{t}U(t-s)P_{0}Q_{<D}\tilde{F}_{\alpha}(s)\,ds
\end{align*}
In other words, $P_{0}Q_{<D}\epsilon_{\alpha}$ solves the following inhomogeneous wave equation:
\begin{align}\label{eq:epslowmod}
 \Box P_{0}Q_{<D}\epsilon=P_{0}Q_{<D}\nabla_{x,t}\big[\tilde{\epsilon}\nabla^{-1}(\psi^{2})
 +\psi\nabla^{-1}(\psi\tilde{\epsilon})+\tilde{\epsilon}\nabla^{-1}(\psi\tilde{\epsilon})
 +\psi\nabla^{-1}(\tilde{\epsilon}^{2})+\tilde{\epsilon}\nabla^{-1}(\tilde{\epsilon}^{2})\big]
\end{align}
First, we identify the terms which can be included in the right-hand
side as source terms since they gain smallness, which is achieved in
part by introducing suitable Fourier localizations. To begin with,
recall that the basic version of the wave maps equation at the level
of the Coulomb gauge is of the schematic form
\[
\Box\psi_{\alpha}=i\partial^{\beta}[\psi_{\alpha} A_{\beta}]-i\partial^{\beta}[\psi_{\beta} A_{\alpha}]+i\partial_{\alpha}[\psi^{\nu} A_{\nu}]
\]
The estimates of Section~\ref{sec:trilin} will be seen to imply that
the middle term here can  be included entirely in the right-hand
side, and the immediately ensuing discussion is only applied to the
first and third terms. Split the first term on the right in
\eqref{eq:epslowmod} (which is understood to be of the first or
third type) into
\[
 P_{0}Q_{<D}\nabla_{x,t}\big[\tilde{\epsilon}\nabla^{-1}(\psi^{2})\big]=
P_{0}Q_{<D}\nabla_{x,t}\big[\tilde{\epsilon}\nabla^{-1}P_{<-D_{1}}(\psi^{2})\big]+
P_{0}Q_{<D}\nabla_{x,t}\big[\tilde{\epsilon}\nabla^{-1}P_{\geq -D_{1}}(\psi^{2})\big]
\]
Here $D_{1}$ is a large constant depending like $D$ on the energy in
a ``mild'' way, i.e., independently of the stage of the induction we
are at, as will be seen shortly. Recalling that on
$I_{1}\times\R^{2}$ we have the decomposition
\[
 \psi=\psi_{L}+\psi_{NL},
\]
we further decompose (schematically)
\begin{align*}
 &P_{0}Q_{<D}\nabla_{x,t}\big[\tilde{\epsilon}\nabla^{-1}P_{<-D_{1}}(\psi^{2})\big]\\
&=P_{0}Q_{<D}\nabla_{x,t}\big[\tilde{\epsilon}\nabla^{-1}P_{<-D_{1}}(\psi_L^{2})\big]
+P_{0}Q_{<D}\nabla_{x,t}\big[\tilde{\epsilon}\nabla^{-1}P_{<-D_{1}}(\psi_{NL}^{2})\big]
+P_{0}Q_{<D}\nabla_{x,t}\big[\tilde{\epsilon}\nabla^{-1}P_{<-D_{1}}(\psi_L\psi_{NL})\big]
\end{align*}
Due to the smallness of $\psi_{NL}$ and~\eqref{eq:prod_small}, it is
only the first term on the right which we need to incorporate in
part into the linear operator. Of course this requires replacing
$\tilde{\epsilon}$ by $\epsilon$, which requires some care due to
the non-local operator $Q_{<D}$ interfering with our aim. First,
write
\[
P_{0}Q_{<D}\nabla_{x,t}\big[\tilde{\epsilon}\nabla^{-1}P_{<-D_{1}}(\psi_L^{2})\big]=
P_{0}\nabla_{x,t}\big[\tilde{\epsilon}\nabla^{-1}P_{<-D_{1}}(\psi_L^{2})\big]-
P_{0}Q_{\geq D}\nabla_{x,t}\big[\tilde{\epsilon}\nabla^{-1}P_{<-D_{1}}(\psi_L^{2})\big]
\]
Since we only need to solve the equation on $I_{1}\times\R^{2}$, where $\tilde{\epsilon}$ and $\epsilon$ agree, we may replace the right-hand side by
\begin{align*}
 &P_{0}\nabla_{x,t}\big[\epsilon\nabla^{-1}P_{<-D_{1}}(\psi_L^{2})\big]-
P_{0}Q_{\geq
D}\nabla_{x,t}\big[\tilde{\epsilon}\nabla^{-1}P_{<-D_{1}}(\psi_L^{2})\big]\\
&=
P_{0}\nabla_{x,t}\big[Q_{<D}\epsilon\nabla^{-1}P_{<-D_{1}}(\psi_L^{2})\big]+P_{0}\nabla_{x,t}\big[Q_{\geq
D}\epsilon\nabla^{-1}P_{<-D_{1}}(\psi_L^{2})\big]-P_{0}Q_{\geq
D}\nabla_{x,t}\big[\tilde{\epsilon}\nabla^{-1}P_{<-D_{1}}(\psi_L^{2})\big]
\end{align*}
Now we introduce null-structure by performing Hodge decompositions
as in Section~\ref{sec:hodge}, for all the trilinear terms. In
particular, the preceding discussion yields that we replace the
schematic term
\[
P_{0}Q_{<D}\nabla_{x,t}\big[\tilde{\epsilon}\nabla^{-1}P_{<-D_{1}}(\psi_L^{2})\big]
\]
by
\begin{align*}
 &\sum_{j=1,3}F_\alpha^{3j}(P_{0}Q_{<D}\epsilon; P_{<-D_{1}};\psi_L, \psi_L)+P_{0}Q_{<D}F_\alpha^{32}(\tilde{\epsilon}; P_{<-D_{1}};\psi_L, \psi_L)\\
&+\big[\sum_{j=1,3}P_{0}F_\alpha^{3j}(Q_{<D}\epsilon; P_{<-D_{1}};\psi_L, \psi_L)-\sum_{j=1,3}F_\alpha^{3j}(P_{0}Q_{<D}\epsilon; P_{<-D_{1}};\psi_L, \psi_L)\big]\\
&+\big[\sum_{j=1,3}P_{0}F_\alpha^{3j}(Q_{\geq D}\epsilon; P_{<-D_{1}};\psi_L, \psi_L)-\sum_{j=1,3}P_0 Q_{\geq D}F_\alpha^{3j}(\tilde{\epsilon}; P_{<-D_{1}};\psi_L, \psi_L)\big]\\
&+\sum_{k=2}^{4}P_{0}Q_{<D}F_\alpha^{2k+1}(\epsilon; P_{<-D_{1}}; \psi_L, \psi_{L})
\end{align*}
We can now write the wave equation that we use to solve for $P_{0}Q_{<D}\epsilon$ as follows:
\begin{equation}\label{eq:epsilonwave}
\begin{aligned}
\Box(P_{0}Q_{<D}\epsilon)= &\sum_{j=1,3}F_\alpha^{3j}(P_{0}Q_{<D}\epsilon; P_{<-D_{1}};\psi_L, \psi_L)+P_{0}Q_{<D}F_\alpha^{32}(\tilde{\epsilon}; P_{<-D_{1}};\psi_L, \psi_L)\\
&+\big[\sum_{j=1,3}P_{0}F_\alpha^{3j}(Q_{<D}\epsilon; P_{<-D_{1}};\psi_L, \psi_L)-\sum_{j=1,3}F_\alpha^{3j}(P_{0}Q_{<D}\epsilon; P_{<-D_{1}};\psi_L, \psi_L)\big]\\
&+\big[\sum_{j=1,3}P_{0}F_\alpha^{3j}(Q_{\geq D}\epsilon; P_{<-D_{1}};\psi_L, \psi_L)-\sum_{j=1,3}P_0 Q_{\geq D}F_\alpha^{3j}(\tilde{\epsilon}; P_{<-D_{1}};\psi_L, \psi_L)\big]\\
&+\sum_{k=2}^{4}P_{0}Q_{<D}F_\alpha^{2k+1}(\psi+\tilde{\epsilon}, (\psi+\tilde{\epsilon}), (\psi+\tilde{\epsilon}))-P_{0}Q_{<D}F_\alpha^{2k+1}(\psi, \psi, \psi)\\
&+P_{0}Q_{<D}F_\alpha^{3}(\tilde{\epsilon}; P_{<-D_{1}};\psi_{NL}, \psi_L)+P_{0}Q_{<D}F_\alpha^{3}(\tilde{\epsilon}; P_{<-D_{1}};\psi_L, \psi_{NL})\\
&+P_{0}Q_{<D}F_\alpha^{3}(\tilde{\epsilon},\psi_{NL}, \psi_{NL})+P_{0}Q_{<D}F_\alpha^{3}(\tilde{\epsilon}; P_{\geq-D_{1}}; \psi_L, \psi_L)\\
&+P_{0}Q_{<D}F_\alpha^{3}(\psi,\tilde{\epsilon}, \psi)+P_{0}Q_{<D}F_\alpha^{3}(\psi,\psi,\tilde{\epsilon})+P_{0}Q_{<D}F_\alpha^{3}(\tilde{\epsilon},\tilde{\epsilon}, \psi)+P_{0}Q_{<D}F_\alpha^{3}(\psi,\tilde{\epsilon}, \tilde{\epsilon})\\&+P_{0}Q_{<D}F_\alpha^{3}(\tilde{\epsilon},\tilde{\epsilon}, \tilde{\epsilon})
\end{aligned}
\end{equation}
The significance of the first term on the right, i.e., the expression
\[
 \sum_{j=1,3}F_\alpha^{3j}(P_{0}Q_{<D}\epsilon,\psi_L, \psi_L),
\]
is that it implicitly contains a magnetic potential interaction
term, see the discussion at the end of Section~\ref{sec:hodge}. In
order to deduce estimates, we shall re-arrange terms and move the
magnetic interaction term contained in the above term
\begin{equation}\label{eq:magn_def}
 2i\partial^{\beta}(P_{0}Q_{<D}\epsilon)A_\beta,\qquad A_\beta:=-P_{<-D_{1}}\del_j^{-1}I\calN_{\beta j}(\psi_L, \psi_L)
\end{equation}
to the left, thereby obtaining an equation of the schematic type
\begin{equation}\label{twistedwave}
\Box(P_{0}Q_{<D}\epsilon)+2i\partial^{\beta}(P_{0}Q_{<D}\epsilon)A_\beta = F
\end{equation}
The next issue occupying us is the derivation of apriori estimates for this type of equation, at first treating $F$ as a function with good
 Fourier localization properties and bounded with respect to $\|\cdot\|_{N}$.

\subsubsection{Solving the wave equation with a magnetic potential in the Coulomb gauge}

For simplicity's sake, replace $(P_{0}Q_{<D}\epsilon)$ at the end of the preceding section by $\epsilon$ for this subsection. The
key fact that is proven here is the following:

\begin{prop}\label{TwistedWaveEquation} Assume that $F$ is a function at frequency $\sim 1$ satisfying the bound
\[
\|F\|_{N[0]}\leq \alpha
\]
Also, assume  the solution to \eqref{twistedwave} with data $(\epsilon(0,\cdot), \partial_{t}\epsilon(0,\cdot))=(f, g)$, all supported at frequency $\sim 1$, to be supported at frequency $\sim 1$ and modulation $\lesssim 1$.  Finally, assume that
\[
D_{1}>D_{1}(\Ecrit),\quad
\alpha<\alpha_{0}(\Ecrit)
\]
Then $\epsilon$ satisfies the bound
\[
\|\epsilon\|_{S[0]}\lesssim  \|F\|_{N[0]}+\|(f, g)\|_{L_{x}^{2}\times\dot{H}^{-1}}
\]
with implied constant only depending on $\Ecrit$. Furthermore, if
$F=0$, there is approximate energy conservation:
\[
\|\partial_{t}\epsilon(t,\cdot)\|_{L_{x}^{2}}^{2}+\|\nabla_{x}\epsilon(t,\cdot)\|_{L_{x}^{2}}=\|\partial_{t}\epsilon(0,\cdot)\|_{L_{x}^{2}}^{2}+\|\nabla_{x}\epsilon(0,\cdot)\|_{L_{x}^{2}}+c(D_1)
\]
with $c(D_1)\to 0$
as $D_1\to \infty$, independently of $t$.
\end{prop}
\begin{proof}
Recall that  $0\le \beta\le 2$,
\[
A_{\beta}=-\triangle^{-1}\sum_{j=1,2}\partial_{j}P_{<-D_1}I[R_{\beta}\psi^{1}_{L}R_{j}\psi^{2}_{L}-R_{\beta}\psi^{2}_{L}R_{j}\psi^{1}_{L}]
\]
and observe that these functions are real-valued and Schwartz for fixed times.
The key difficulty comes from the fact that there appears no obvious way to obtain smallness for the linear interaction term $2i\partial^{\beta}\epsilon A_{\beta}$, even when restricting to small time intervals.
 The easiest way out of this impasse is to use an approximate apriori bound resulting from energy conservation.
This will allow us to split the bad interaction term into two, one of which is small due to angular alignment of the inputs,
 the other of which is controlled due to the apriori bound. Moreover, we note that we may always move parts of the expression $2i\partial^{\beta}\epsilon A_{\beta}$ with additional smallness properties, such as extreme frequency discrepancies inside $A_{\beta}$ or special angular alignments, to the right-hand side, since we gain smallness for them as shown in~Section~\ref{sec:trilin}.  More precisely,
one writes the underlying equation~\eqref{twistedwave} as
\[
\Box\epsilon+2i\partial^{\beta}\epsilon \tilde{A}_{\beta}=F-2i\partial^{\beta}\epsilon {A}^\dagger_{\beta},\quad \epsilon[0]=(f,g)
\]
Here we define
\begin{equation}\nonumber\begin{split}
A^\dagger_{\beta}&:= -\sum_{\kappa_{1,2}\in K_{-C_{6}(\Ecrit)},\text{dist}(\kappa_{1,2})<c_{6}(\Ecrit)}\quad \sum_{\substack{\max\{k_{1,2,3}\}\le \min\{k_{1,2,3}\}+C_6(\Ecrit)\\k_1<-D_1}} \\
&I\triangle^{-1}\sum_{j=1,2}\partial_{j}P_{k_{1}+O(1)}[R_{\beta}P_{k_{1},\kappa_{1}}\psi^{1}_{L}R_{j}P_{k_2,\kappa_{2}}\psi^{2}_{L}-R_{\beta}P_{k_2,\kappa_{2}}\psi^{2}_{L}R_{j}P_{k_{1},\kappa_{1}}\psi^{1}_{L}]\\
&-\sum_{\substack{\max\{k_{1,2,3}\}>\min\{k_{1,2,3}\}+C_6(\Ecrit)\\k_1<-D_1}}I\triangle^{-1}\sum_{j=1,2}\partial_{j}P_{k_{1}}[R_{\beta}P_{k_{2}}\psi^{1}_{L}R_{j}P_{k_{3}}\psi^{2}_{L}-R_{\beta}P_{k_{2}}\psi^{2}_{L}R_{j}P_{k_{3}}\psi^{1}_{L}]
\end{split}\end{equation}
and furthermore $\tilde{A}_{\beta}:=A_{\beta}-A^\dagger_{\beta}$. The
angular separation $c_{6}(\Ecrit)$ is chosen in such fashion that,
for some sufficiently large $C_{7}=C_{7}(\Ecrit)$,
\[
\|-i\partial^{\beta}\epsilon A^\dagger_{\beta}\|_{N[0]}\leq \frac{1}{C_{7}}\|\epsilon\|_{S[0]}
\]
We  then use the preceding reasoning with $A_{\beta}$ replaced by $\tilde{A}_{\beta}$, assuming $D_1(\Ecrit)$ to be large enough.

\noindent After these preparations, we commence by establishing the aforementioned apriori bound: Specifically, consider the {\it{covariant  energy density}}
\[
\frac{1}{2}[|\partial_{t}\epsilon+iA_{0}\epsilon|^{2}+\sum_{j=1,2}|\partial_{x_{j}}\epsilon+iA_{j}\epsilon|^{2}]
\]
where we write $A_\beta$ instead of $\tilde A_\beta$ for simplicity. Now compute
\begin{equation}\nonumber\begin{split}
&\partial_{t}\big[\frac{1}{2}|\partial_{t}\epsilon+iA_{0}\epsilon|^{2}+\sum_{j=1,2}\frac{1}{2}|\partial_{x_{j}}\epsilon+iA_{j}\epsilon|^{2}\big]\\
& =\Re\big[(\partial_{t}\epsilon+iA_{0}\epsilon)\overline{(\partial_{t}+iA_{0})^2\epsilon}+\sum_{j=1,2}(\partial_{x_{j}}\epsilon+iA_{j}\epsilon)\overline{\partial_{t}(\partial_{x_{j}}\epsilon+iA_{j}\epsilon)}\big]
\end{split}\end{equation}
 The second term on the right satisfies
\begin{equation}\nonumber\begin{split}
\Re[(\partial_{x_{j}}\epsilon+iA_{j}\epsilon)\overline{\partial_{t}(\partial_{x_{j}}\epsilon+iA_{j}\epsilon)}]
&=\Re[(\partial_{x_{j}}\epsilon+iA_{j}\epsilon)\overline{(\partial_{x_{j}}+iA_{j})(\partial_{t}+iA_{0})\epsilon}] \\
&\qquad +\Re[(\partial_{x_{j}}\epsilon+iA_{j}\epsilon)\overline{i(\partial_{t}A_{j}-\partial_{x_{j}}A_{0})\epsilon}]\\
&=\partial_{x_{j}}\Re[(\partial_{x_{j}}\epsilon+iA_{j}\epsilon)\overline{(\partial_{t}+iA_{0})\epsilon}]-\Re[(\partial_{x_{j}}+iA_{j})^2\epsilon\,  \overline{(\partial_{t}+iA_{0})\epsilon}]\\
& +\Re[(\partial_{x_{j}}\epsilon+iA_{j}\epsilon)\overline{i(\partial_{t}A_{j}-\partial_{x_{j}}A_{0})\epsilon}]
\end{split}\end{equation}
In summary, one obtains the following {\em local form of energy conservation}:
\begin{equation}\label{DivIdentity}\begin{split}
&\partial_{t}\Big[\frac{1}{2}|\partial_{t}\epsilon+iA_{0}\epsilon|^{2}+\sum_{j=1,2}\frac{1}{2}|\partial_{x_{j}}\epsilon+iA_{j}\epsilon|^{2}\Big]
-\sum_{j=1}^2\partial_{x_{j}}\Re\big[(\partial_{x_{j}}\epsilon+iA_{j}\epsilon)\overline{(\partial_{t}+iA_{0})\epsilon}\big]\\
&=\Re\big[\big[(\partial_{t}+iA_{0})^2-\sum_{j=1,2}(\partial_{x_{j}}+iA_{j})^2\big]\epsilon\,
\overline{(\partial_{t}\epsilon+iA_{0}\epsilon) }
\big]+\Re[(\partial_{x_{j}}\epsilon+iA_{j}\epsilon)\overline{i(\partial_{t}A_{j}-\partial_{x_{j}}A_{0})\epsilon}]\\
\end{split}\end{equation}
We furthermore observe that any solution of $\Box\eps+2iA^\alpha\,\del_\alpha\eps=F$ satisfies
\[
\big[(\partial_{t}+iA_{0})^{2}-\sum_{j=1,2}(\partial_{x_{j}}+iA_{j})^{2}\big]\epsilon=F+i(\partial_{t}A_{0}-\sum_{j=1,2}\partial_{x_{j}}A_{j})\epsilon+(\sum_{j=1,2}A_{j}^{2}-A_{0}^{2})\epsilon
\]
We now integrate the above relation over a time slice $[0,t_0]\times\R^{2}$, which gives
\begin{equation}\nonumber\begin{split}
&\int_{\R^{2}}\big[\frac{1}{2}|\partial_{t}\epsilon+iA_{0}\epsilon|^{2}+\sum_{j=1,2}\frac{1}{2}|\partial_{x_{j}}\epsilon+iA_{j}\epsilon|^{2}\big](t_0, x)\,dx\\
&=\int_{\R^{2}}\big[\frac{1}{2}|\partial_{t}\epsilon+iA_{0}\epsilon|^{2}+\sum_{j=1,2}\frac{1}{2}|\partial_{x_{j}}\epsilon+iA_{j}\epsilon|^{2}\big](0, x)\,dx\\
&+\int_{[0, t_0]\times\R^{2}}\Re\big[\big(F+i(\partial_{t}A_{0}-\sum_{j=1,2}\partial_{x_{j}}A_{j})\epsilon+(\sum_{j=1,2}A_{j}^{2}-A_{0}^{2})\big)\epsilon \,
\overline{(\partial_{t}\epsilon+iA_{0}\epsilon)}\big]\, dtdx\\
&+\int_{[0, t_0]\times\R^{2}}\Re[(\partial_{x_{j}}\epsilon+iA_{j}\epsilon)\overline{i(\partial_{t}A_{j}-\partial_{x_{j}}A_{0})\epsilon}]\, dtdx\\
\end{split}\end{equation}
In order to proceed, we now make the following bootstrap assumption:
\begin{equation}\label{bootstr}
\|\epsilon\|_{S[0]}\leq C [\|F\|_{N[0]}+\|\epsilon[0]\|_{L_{x}^{2}\times\dot{H}^{-1}}]
\end{equation}
We shall then show that if $C>C_{0}(\Ecrit)$, then one may replace
$C$ by $\frac{C}{2}$ in~\eqref{bootstr}. To accomplish this, we
first estimate the expression
\begin{equation}\nonumber\begin{split}
&\int_{[0, t_0]\times\R^{2}}\Re\big[(\partial_{t}\epsilon+iA_{0}\epsilon)\overline{\big(F+i(\partial_{t}A_{0}-\sum_{j=1,2}\partial_{x_{j}}A_{j})\epsilon+(\sum_{j=1,2}A_{j}^{2}-A_{0}^{2})\epsilon}\big)\big]\, dtdx\\
&+\int_{[0, t_0]\times\R^{2}}\Re[(\partial_{x_{j}}\epsilon+iA_{j}\epsilon)\overline{i(\partial_{t}A_{j}-\partial_{x_{j}}A_{0})\epsilon}]\, dtdx,\\
\end{split}\end{equation}
For each of these terms we must gain a smallness constant.
One can classify three types of terms.

\smallskip
\noindent (1) The term $\int_{[0, t_0]\times\R^{2}}\Re\big[(\partial_{t}\epsilon+iA_{0}\epsilon)\overline{F}\big]\, dtdx$.
Here one uses the duality of $N$ and $S$, Lemma~\ref{lem:Ncut}, as well as the space-time frequency localization of $\epsilon$:
\[
\Big|\int_{[0, t_0]\times\R^{2}}\Re\big[(\partial_{t}\epsilon+iA_{0}\epsilon)\overline{F}\big]\, dtdx\Big|\lesssim \|F\|_{N[0]}\|\epsilon\|_{S}\leq \alpha \|\epsilon\|_{S}
\]
Application of Lemma~\ref{lem:Ncut} is justified due to our assumptions on the modulation of~$\epsilon$, which in turn
restrict the modulation of~$F$ to the hyperbolic regime via the equation.

\smallskip\noindent
(2) The terms of the form $\int_{[0, t_0]\times\R^{2}}\nabla_{x,t}A\; \epsilon\nabla_{x,t}\epsilon\, dtdx$.
These are quite delicate and we can just barely control them.  Note the schematic identity
\[
\nabla_{x,t}A_{\beta}=\sum_{k_{1,2}<-D_1}\nabla_{x,t}\nabla^{-1}P_{k_{1}}[P_{k_{2}}\psi_{L}P_{k_{23}}\psi_{L}]
\]
Here our reductions for $A_{\beta}$ (which is $\tilde{A}_{\beta}$ in the discussion above) imply $k_{1}=k_{2}+O(1)=k_{3}+O(1)$ (where the implied constant
may be quite large depending on $\Ecrit$) and furthermore the inputs $P_{k_{2,3}}\psi_L$
have some angular separation between their Fourier supports. Using the mixed-norm Wolff-type endpoint result
established by Tataru~\cite{TatWolff}, one obtains
\[
\|P_{k_{1}}\psi_{L}P_{k_{2}}\psi_{L}\|_{L_{t}^{\frac{4}{3}}L_{x}^{2}}\lesssim 2^{\frac{k_{1}}{4}}
\]
with implied constant depending on $\|P_{k_{1,2}}\psi_{L}\|_{L_{x}^{2}}\les \Ecrit+1$. But then using the fact that
\[
\|\epsilon\|_{L_{t}^{8}L_{x}^{4}}\lesssim \|\epsilon\|_{S[0]}
\]
see Lemma~\ref{lem:Strich},
one infers that
\[
\Big |\int_{[0, t_0]\times\R^{2}}\nabla_{x,t}P_{<k}A\; \epsilon\nabla_{x,t}\epsilon\, dtdx \Big|\lesssim
2^{\frac{k}{4}}\|\epsilon\|_{S[0]}^{2}
\]

\smallskip\noindent
(3) The terms \[ \int_{[0, t_0]\times\R^{2}}A^{2}\: \epsilon\nabla_{x,t}\epsilon\, dtdx \]
which are handled similarly. Here, one uses again that
\[
\|[P_{k}A]^{2}\|_{L_{t}^{\frac{4}{3}}L_{x}^{2}}\lesssim E^4_{\mathrm{crit}}\, 2^{\frac{k}{4}},
\]
which  follows from the usual Strichartz estimates, cf.~Lemma~\ref{lem:Strich}:
\[
 \|\nabla^{-1} P_k(\psi_L^2)\|_{L_t^{\frac83} L^4_x} \les 2^{-k} \|\psi_L\|_{L_t^{\frac{16}{3}} L^8_x}^2 \les 2^{\frac{k}{8}} \|\psi_L\|_2^2
\]
Summation over small $k$ now yields the desired smallness provided $k$ is small enough, which is ensured if $D_1$ is large enough.

\medskip
\noindent In view of the preceding, we may conclude that
\begin{equation}
 \label{eq:cov_ener}
\|\nabla_{x,t}\epsilon(t, \cdot)\|_{L_{x}^{2}}^{2}=\|\nabla_{x,t}\epsilon(0, \cdot)\|_{L_{x}^{2}}^{2}+O(\gamma^2 \|\epsilon\|_{S[0]}^{2}+\alpha\|\epsilon\|_{S[0]}),
\end{equation}
where $\gamma$ may be made arbitrarily small by choosing
$D_1(\Ecrit)$ in the statement of the proposition large enough.
Note that we eliminated the magnetic potential here from the
covariant energy by means of the estimate $\|A_\beta\|_{\Linf}\les
\gamma\ll1$.

\noindent
Unfortunately, this a apriori bound is still insufficient to estimate the magnetic potential interaction term
$2i\partial^{\beta}\epsilon A_{\beta}$, and instead we need to gain apriori control over one of the null-frame ingredients of~$\|\cdot\|_{S[0]}$.
Indeed, a very slight modification of the above essentially allows us to control
\[
\sum_{\kappa\in K_{-C_{1}}}\|P_{0,\kappa}\epsilon\|_{NF[\kappa]^{*}}^{2}
\]
for any $C_{1}=C_{1}(\Ecrit)$ (with implied constant depending on
$\Ecrit$), which will then be good enough to move the entire
magnetic potential interaction to the right-hand side. To get this
extra control we argue almost exactly as before, but
integrating~\eqref{DivIdentity} over a region
$A^{\omega}_{t,c}:=[0,t]\times\R^{2}\cap \{t_{\omega}>c\}$ for
arbitrary~$c$, with $\omega\in S^{1}$ being a fixed direction. This
yields
\begin{equation}\nonumber\begin{split}
&\int_{A^{\omega}_{t,c}}\Re\big[(\partial_{t}\epsilon+iA_{0}\epsilon)\overline{\big[(\partial_{t}+iA_{0})^2-\sum_{j=1,2}(\partial_{x_{j}}
+iA_{j})^2\big]\epsilon}\big]+\Re[(\partial_{x_{j}}\epsilon+iA_{j}\epsilon)\overline{i(\partial_{t}A_{j}-\partial_{x_{j}}A_{0})\epsilon}]\, dtdx\\
&=\int_{0\times\R^{2}\cap A^{\omega}_{t,c}}\big[\frac{1}{2}|\partial_{t}\epsilon+iA_{0}\epsilon|^{2}+\sum_{j=1,2}\frac{1}{2}|\partial_{x_{j}}\epsilon+iA_{j}\epsilon|^{2}\big]\, dtdx\\
&-\int_{t\times\R^{2}\cap A^{\omega}_{t,c}}\big[\frac{1}{2}|\partial_{t}\epsilon+iA_{0}\epsilon|^{2}+\sum_{j=1,2}\frac{1}{2}|\partial_{x_{j}}\epsilon+iA_{j}\epsilon|^{2}\big]\, dtdx\\
&+\int_{\{t_{\omega}=c\}\cap A^{\omega}_{t,c}}\big[\frac{1}{2}|\partial_{t}\epsilon+iA_{0}\epsilon|^{2}+
\sum_{j=1,2}\big(\frac{1}{2}|\partial_{x_{j}}\epsilon+iA_{j}\epsilon|^{2}-\omega_{j}
\Re[(\partial_{x_{j}}\epsilon+iA_{j}\epsilon)\overline{(\partial_{t}+iA_{0})\epsilon}]\big)\big] \, dx_{\omega}
\end{split}\end{equation}
It is the latter integral expression that gives us the additional
information we need: to see this, note that we may localize the
entire equation \eqref{twistedwave} to an angular sector
$\kappa\subset S^{1}$ by applying $P_{0,\kappa}$ to both sides,
where $|\kappa|\sim 2^{-C_{1}}$ and with $C_{1}=C_{1}(\Ecrit)$ being
a fixed large constant. Replacing $\epsilon$ by
$P_{0,\kappa}\epsilon$ generates an error which can be incorporated
into the right-hand $F$, as is easily seen. Hence in the above we
may replace $\epsilon$ by $P_{0,\kappa}\epsilon$. Then use the
decomposition
\[
\sum_{j=1,2}  |\partial_{x_{j}}P_{0,\kappa}\epsilon+iA_{j}P_{0,\kappa}\epsilon|^{2}
= \big|\sum_{j=1,2}(\omega_{j}\partial_{x_{j}}P_{0,\kappa}\epsilon+i\omega_{j}A_{j}P_{0,\kappa}\epsilon)\big|^{2}
+ \big|\sum_{j=1,2}(\omega^{\perp}_{j}\partial_{x_{j}}P_{0,\kappa}\epsilon+i\omega^{\perp}_{j}A_{j}P_{0,\kappa}\epsilon) \big|^{2}
\]
If now we have $\omega\notin 2\kappa$, then we can conclude that
\[
\sup_{c}\sup_{t}\int_{\{t_{\omega}=c\}\cap A^{\omega}_{t,c}} \Big| \sum_{j=1,2}(\omega^{\perp}_{j}\partial_{x_{j}}P_{0,\kappa}\epsilon+i\omega^{\perp}_{j}A_{j}P_{0,\kappa}\epsilon) \Big|^{2}\gtrsim \|P_{0,\kappa}\epsilon\|_{L_{t_{\omega}}^{\infty}L_{x_{\omega}}^{2}}^{2},
\]
since the magnetic potential is small in~$\Linf$. Using~\eqref{eq:cov_ener} as well as
\[
\frac{1}{2}\Big|\sum_{j=1,2}(\omega_{j}\partial_{x_{j}}P_{0,\kappa}\epsilon+i\omega_{j}A_{j}P_{0,\kappa}\epsilon)\Big|^{2}
+\frac{1}{2} |\partial_{t}\epsilon+iA_{0}\epsilon|^{2}-\sum_{j=1,2}\omega_{j}
\Re[(\partial_{x_{j}}\epsilon+iA_{j}\epsilon)\overline{(\partial_{t}+iA_{0})\epsilon}]\geq 0,
\]
we infer that
\[
\Big(\sum_{\kappa\in K_{-C_{1}}}\|P_{0,\kappa}\epsilon\|_{NF[\kappa]^{*}}^{2}\Big)^{\frac{1}{2}}\lesssim \|\nabla_{x,t}\epsilon(0,\cdot)\|_{L_{x}^{2}}+
\frac{\alpha}{\gamma}+\gamma \|\epsilon\|_{S[0]},
\]
where $\gamma$ can be made small independently of $\Ecrit$, and with
an implicit constant depending only on~$\Ecrit$ as well as $C_1=2C_6$.

\noindent
We can finally complete the proof of Proposition~\ref{TwistedWaveEquation}:
By choosing $D_1$ large enough, in relation to $C_{6}$ in the definition of $A^\dagger_{\beta}$ and then making  $\sup_{k\in\Z}\|P_{k}\psi_{L}\|_{L_{x}^{2}}$ small enough, one obtains
\begin{equation}\label{encon}
(\sum_{\kappa\in K_{-C_{1}}}\|P_{0,\kappa}\epsilon\|_{NF[\kappa]^{*}}^{2})^{\frac{1}{2}}\lesssim \|\nabla_{x,t}\epsilon(0,\cdot)\|_{L_{x}^{2}}+
\frac{\alpha}{\gamma}+\gamma \|\epsilon\|_{S[0]},
\end{equation}
as before, where the implicit constant only depends on $\Ecrit,
C_{1}$, while $\gamma$ may be made small independently of $\Ecrit,
C_{1}$. Now we {\it{feed this information back into the magnetic
interaction term}}:  write \eqref{twistedwave} in the form
\[
\Box \epsilon = -2iA_{\beta}\epsilon+F
\]
and decompose
\begin{equation}\nonumber\begin{split}
&2iA_{\beta}\partial^{\beta}\epsilon=\\
&\sum_{\kappa_{1,2,3}\in K_{-C_{1}},\,\min_{i\neq j}\{\text{dist}(\kappa_{i}, \kappa_{j})\}\geq 2^{-C_{1}+10}}
2i\partial^{\beta}P_{0,\kappa_{1}}\epsilon\; I\triangle^{-1}\sum_{j=1,2}\partial_{j}[R_{\beta}P_{\kappa_{2}}\psi^{1}_{L}R_{j}P_{\kappa_{3}}\psi^{2}_{L}-R_{\beta}P_{\kappa_{2}}\psi^{2}_{L}R_{j}P_{\kappa_{3}}\psi^{1}_{L}]\\
&+ \sum_{\kappa_{1,2,3}\in K_{-C_{1}},\,\min_{i\neq j}\{\text{dist}(\kappa_{i}, \kappa_{j})\}<2^{-C_{1}+10}}  2i\partial^{\beta}P_{0,\kappa_{1}}\epsilon\; I\triangle^{-1}\sum_{j=1,2}\partial_{j}[R_{\beta}P_{\kappa_{2}}\psi^{1}_{L}R_{j}P_{\kappa_{3}}\psi^{2}_{L}-R_{\beta}P_{\kappa_{2}}\psi^{2}_{L}R_{j}P_{\kappa_{3}}\psi^{1}_{L}]\\
\end{split}\end{equation}
Picking $C_{1}$ large enough in relation to $\Ecrit$ as well as
$\nu$, one obtains
\begin{equation}\nonumber\begin{split}
&\| \sum_{\kappa_{1,2,3}\in K_{-C_{1}},\,\min_{i\neq j}\{\text{dist}(\kappa_{i}, \kappa_{j})\}<2^{-C_{1}+10}} \partial^{\beta}P_{0,\kappa_{1}}\epsilon\; I\triangle^{-1}\sum_{j=1,2}\partial_{j}[R_{\beta}P_{\kappa_{2}}\psi^{1}_{L}R_{j}P_{\kappa_{3}}\psi^{2}_{L}-R_{\beta}P_{\kappa_{2}}\psi^{2}_{L}R_{j}P_{\kappa_{3}}\psi^{1}_{L}]\|_{N[0]}\\
&\leq \nu \|\epsilon\|_{S[0]}
\end{split}\end{equation}
for any $\nu\in (0,1)$. Furthermore, for the first term above, using \eqref{encon} we directly arrive at the estimate
\begin{equation}\nonumber\begin{split}
 &\|\sum_{\kappa_{1,2,3}\in K_{-C_{1}},\,\min_{i\neq j}\{\text{dist}(\kappa_{i}, \kappa_{j})\}\geq 2^{-C_{1}+10}} \partial^{\beta}P_{0,\kappa_{1}}\epsilon\; I\triangle^{-1}\sum_{j=1,2}\partial_{j}[R_{\beta}P_{\kappa_{2}}\psi^{1}_{L}R_{j}P_{\kappa_{3}}\psi^{2}_{L}-R_{\beta}P_{\kappa_{2}}\psi^{2}_{L}R_{j}P_{\kappa_{3}}\psi^{1}_{L}]\|_{N[0]}\\
 &\lesssim  \|\nabla_{x,t}\epsilon(0,\cdot)\|_{L_{x}^{2}}+C_{2}(\nu)\alpha+\nu \|\epsilon\|_{S[0]},
 \end{split}\end{equation}
 with implied constant depending on $C_{1}, \Ecrit, \nu$. Here $\nu$ may be made small independently of the latter parameters.
Summarizing, using the fundamental energy estimate for $\Box u = F$, we then obtain (for a universal constant $C_{0}$)
\[
\|\epsilon\|_{S[0]}\leq
C_{0}[\|\epsilon[0]\|_{L_{x}^{2}\times\dot{H}^{-1}}+\alpha+\nu\|\epsilon\|_{S[0]}+C_{11}(C_{1},
\Ecrit, \nu)[
\|\nabla_{x,t}\epsilon(0,\cdot)\|_{L_{x}^{2}}+C_{2}(\nu)\alpha+\nu\|\epsilon\|_{S[0]}]
\]
If  we now select $C$ in \eqref{bootstr} large enough in relation to
$C_{0}, \Ecrit$, then $\nu$ small enough, and finally $\gamma$
and then $\alpha$ small enough, we conclude that
\[
\|\epsilon\|_{S[0]}\leq \frac{C}{2} [\|F\|_{N[0]}+\|\epsilon[0]\|_{L_{x}^{2}\times\dot{H}^{-1}}],
\]
which is the desired bootstrap.
 \end{proof}

Due to frequency leakage coming from the magnetic term we shall also require energy estimates that take $N \to S$, or alternatively, preservation
of frequency envelope.
For the following lemma, we allow more general frequency support of~$A_\alpha$. Hence consider the following equation
\begin{equation}
\label{eq:uFfg}
 \Box u + 2i\del^\alpha u\, A_\alpha=F, \quad u[0]=(f,g)
\end{equation}
where $F$ has the property that $F=F_1+F_2$ where $\|F_1\|_{N}:=
\Big(\sum_{k\in\Z} \|F\|_{N[k]}^2\Big)^{\frac{1}{2}}$ is finite and
with~$F_2$ controlled by a frequency envelope, i.e., $\|P_k
F_2\|_{N[k]}\le c_k$ and $\{c_k\}_{k\in\Z}$ is sufficiently flat (as
defined above). Furthermore, $(f,g)=(f_1,g_1)+(f_2,g_2)$ with
$\|(f_1,g_1)\|_{L^2\times \dot H^{-1}}$ finite and
$\|P_k(f_2,g_2)\|_{L^2\times \dot H^{-1}}\le d_k$ where $d_k$ is
again a sufficiently flat envelope. Finally,
\[
 A_\alpha = \Delta^{-1}\sum_{j=1,2}\partial_{j}[R_{\beta}\psi^{1}_{L}R_{j}\psi^{2}_{L}-R_{\beta}\psi^{2}_{L}R_{j}\psi^{1}_{L}]
\]
is more general than in the previous proposition. Here $\psi^1_L$ and $\psi^2_L$ are finite energy free waves (with energy bounded by~$\Ecrit$).
Now one has the following result.

\begin{lemma}
 \label{lem:F1F2} Let $u$ be a solution of \eqref{eq:uFfg} with $F$ and $f,g$ as above. Then $u=u_1+u_2+u_3$
where $\|u_1\|_S\les \|F_1\|_N+ \|(f_1,g_1)\|_{L^2\times \dot H^{-1}}$, and $\|P_k u_2\|_{S[k]} \les c_k$, $\|P_k u_3\|_{S[k]}\les d_k$.
The implied constants only depend on the energy of~$\psi_L^{1,2}$.
\end{lemma}
\begin{proof}
 We restrict ourselves to $P_j u$. By scaling $j=0$. Now split $A_\alpha=\sum_{i=1}^3 A_\alpha^{(i)}$ where
\begin{equation}\label{eq:Asplit}
 A_\alpha^{(1)}=\sum_{k<-C} P_k A_\alpha, \quad A_\alpha^{(2)}=\sum_{-C<k<C} P_k A_\alpha
\end{equation}
The constant $C$ in~\eqref{eq:Asplit} is chosen such that the proof of Proposition~\ref{TwistedWaveEquation}
applies to the low frequency part of~$A_\alpha$. Then we write
\begin{align*}
 \Box P_0 u + 2i\del^\alpha P_0 u\, A_\alpha^{(1)} &=P_0 F - 2iP_0[\del^\alpha u\, A_\alpha^{(2)} + \del^\alpha u\, A_\alpha^{(3)}] +
L(\del^\alpha  u, \nabla A_\alpha^{(1)}) \\
 P_0 u[0] &= P_0(f,g)
\end{align*}
Here $L(\cdot,\cdot)$ stands for the commutator in $P_0(uv)=vP_0u+ L(u,\nabla v)$.
We now divide $\R=\bigcup_{\ell=1}^M I_\ell$ into disjoint intervals with the property that
\[
 \max_\ell \|P_0[\del^\alpha u\, A_\alpha^{(2)} + \del^\alpha u\, A_\alpha^{(3)}]\|_{N[0](I_\ell\times\R^2)} \le \gamma\,
\|\psi_L\|_{S}^2\:\sum_{k\in\Z} 2^{-\sigma_0|k|} \|P_k u\|_{S[k]}
\]
where $\gamma>0$ can be made arbitrarily small, and $M=M(\Ecrit,\gamma)$. To see this, one argues as in several previous instances.
First consider~$A^{(1)}$.  In case of angular alignment of the Fourier supports of any two of the inputs, one
obtains a gain as shown in Section~\ref{sec:trilin}. On the other hand, in case of angular separation
the desired smallness is achieved by a careful choice of the~$I_\ell$, see Section~\ref{sec:perturb}. Finally, for $A^{(2)}$ one uses the high-high-low gains
in the trilinear estimates (see the form of the weights $w(j_1,j_2,j_3)$ in Section~\ref{sec:trilin} when $\max(j_2,j_3)>C$).
The commutator terms satisfies the bound
\[
 \| L(\del^\alpha  u, \nabla A_\alpha^{(1)})\|_{N[0]} \les 2^{-C} \, \|\psi_L\|_{S}^2\:\sum_{k\in\Z} 2^{-\sigma_0|k|} \|P_k u\|_{S[k]}
\]
since $\nabla A_\alpha^{(1)}$ gains a factor of $2^{-C}$. We now apply the covariant energy bound of Proposition~\ref{TwistedWaveEquation} to conclude
that (with $P_j$ instead of $P_0$)
\[
 \|P_j u\|_{S[j]} \le C(\Ecrit)(\|P_j(f,g)\|_{L^2\times \dot H^{-1}} + (\gamma+2^{-C})\sum_{k\in\Z} 2^{-\sigma_0|k-j|} \|P_k u\|_{S[k]} + \|P_j F\|_{N[j]})
\]
The lemma now follows from this estimate provided the frequency envelopes are flat enough compared to~$\sigma_0$.
\end{proof}

\subsubsection{Controlling the error terms}

In this section, we complete the proof of Proposition~\ref{PsiBootstrap}. This amounts to
bounding each of the terms on the right-hand side of~\eqref{eq:epsilonwave} one by one using
the covariant energy estimate of the previous section.

We begin with the first term in~\eqref{eq:epsilonwave}, i.e., $\sum_{j=1,3}F_\alpha^{3j}(P_{0}Q_{<D}\epsilon; P_{<-D_1}; \psi_L, \psi_L)$
from which we which we have subtracted the magnetic potential term. Thus, we claim that we can
decompose, with $A_\beta$ as in \eqref{eq:magn_def},
\[
 \sum_{j=1,3}F_\alpha^{3j}(P_{0}Q_{<D}\epsilon; P_{<-D_1};\psi_L, \psi_L) - 2i\del^\beta\epsilon\, A_\beta
\]
into the sum of two terms, one of which has controlled frequency envelope and the other small $S$-norm as in~\eqref{eq:improvboot}.
By~\eqref{eq:Fthree}, this difference equals
\begin{align*}
& iP_0Q_{<D}\epsilon_\alpha \, I \partial^{\beta}
\del_j^{-1}P_{<-D_{1}}\calN_{\beta j}(\psi_L,\psi_L) +iP_0 Q_{<D} R^\beta \epsilon\,
\,\del_j^{-1} I P_{<-D_1}\partial_{\alpha}\calN_{\beta j}(\psi_L,\psi_L)
\end{align*}
Denoting these terms by $\term_{11}$ and $\term_{12}$, respectively, we now proceed to estimate
them by means of Section~\ref{subsec:hyp_trilin}. Let us now assume that $\epsilon$ is of the envelope type, see~$\epsilon_1$ in \eqref{eq:bootass}.
Then by~\eqref{eq:dergain}
\[
 \|\term_{11}\|_{N[0]} \les 2^{-\sigma D_1} \|P_0\epsilon_\alpha\|_{S[0]} \|\psi_L\|_S^2 \ll C_4\, d_0
\]
for $D_1$ large depending on~$\Ecrit$. The contribution of $\epsilon_2$ is estimated similarly.
For $\term_{12}$ one  uses that
\[
 \|P_0 Q_{<D}  R^\beta \epsilon\,
\,\del_j^{-1} I P_{k}\partial_{\alpha}\calN_{\beta j}(\psi_L,\psi_L)\|_{N[0]} \les 2^{k} \|P_0\epsilon\|_{S[0]} \|\psi_L\|_S^2
\]
which is sufficient for both $\epsilon_{1,2}$ since $k\le -D_1$.

The second term in~\eqref{eq:epsilonwave} is bounded by
\[
 \| P_0 Q_{<D} [R_\beta \tilde\epsilon\,
\,\del_j^{-1} I P_{<-D_1}\partial^{\beta}\calN_{\alpha j}(\psi_L,\psi_L)]\|_{N[0]} \les 2^{-D_1} \|P_0\epsilon\|_{S[0]} \|\psi_L\|_S^2
\]
which is sufficient.
This follows from Lemma~\ref{lem:tri_hyp}.

The third term in~\eqref{eq:epsilonwave} is the commutator
\begin{align*}
& \sum_{j=1,3}P_{0}F_\alpha^{3j}(Q_{<D}\epsilon; P_{<-D_{1}};\psi_L, \psi_L)-\sum_{j=1,3}F_\alpha^{3j}(P_{0}Q_{<D}\epsilon; P_{<-D_{1}};\psi_L, \psi_L)\\
&= \tilde P_0 \del^\beta[Q_{<D}\epsilon_\alpha\, \nabla A_\beta] + \tilde P_0 \del^\alpha[Q_{<D}\epsilon^\beta\,\nabla A_\beta]
\end{align*}
where the second line is schematic. Hence, the smallness for this term is obtained just as in the preceding term.

Next, as the fourth term we face the commutator
\begin{align}
& \sum_{j=1,3}P_{0}F_\alpha^{3j}(Q_{\geq D}\epsilon; P_{<-D_{1}};\psi_L, \psi_L)-\sum_{j=1,3}P_0 Q_{\geq D}F_\alpha^{3j}(\tilde{\epsilon}; P_{<-D_{1}};\psi_L, \psi_L) \nn \\
& =\sum_{j=1,3}P_{0}F_\alpha^{3j}(Q_{\geq D}\epsilon; P_{<-D_{1}};\psi_L, \psi_L)-\sum_{j=1,3}P_0 Q_{\geq D}F_\alpha^{3j}(\epsilon; P_{<-D_{1}};\psi_L, \psi_L) \label{eq:diffeins} \\
& +\sum_{j=1,3}P_0 Q_{\geq D}F_\alpha^{3j}(\epsilon; P_{<-D_{1}};\psi_L, \psi_L)-\sum_{j=1,3}P_0 Q_{\geq D}F_\alpha^{3j}(\tilde{\epsilon}; P_{<-D_{1}};\psi_L, \psi_L)  \label{eq:diffzwei}
\end{align}
First,  \eqref{eq:diffeins} is bounded by
\[
 2^{-D_1-D} \|P_0\epsilon\|_{S[0]} \|\psi_L\|_S^2
\]
by the same commutator logic as before and Lemma~\ref{lem:tri_hyp}. Second, since $\epsilon=\tilde\epsilon$ on $I_1$, the length of which is bounded below by an absolute constant by Case~1 above,  one obtains that
\begin{align*}
& \big\|\sum_{j=1,3}P_0 Q_{\geq D}F_\alpha^{3j}(\epsilon-\tilde\epsilon; P_{<-D_{1}};\psi_L, \psi_L) \big\|_{N[0](I_1\times\R^2)} \\
& \les \big\|\sum_{j=1,3}P_0 Q_{\geq D}F_\alpha^{3j}(\epsilon-\tilde\epsilon; P_{<-D_{1}};\psi_L, \psi_L) \big\|_{L^1_t(I_1;L^2_x)} \\
& \les 2^{-D} \|P_0(\epsilon-\tilde\epsilon)\|_{\ener} \| P_{<-D_{1}} \nabla^{-1}(\psi_L^2)\|_{\Linf} \\
& \les 2^{-D-D_1} \|P_0(\epsilon-\tilde\epsilon)\|_{\ener} \| \psi_L\|_S^2
\end{align*}
where we used Bernstein's inequality in the last line.

The fifth term is a collection of quintilinear and higher oder terms, and we deal with it in the appendix.

The terms six through eight are easy by Lemma~\ref{lem:tri_hyp} and \eqref{eq:prod_small}. More precisely, they
inherit the frequency profile of $\tilde\epsilon$ times a factor of~$\sqrt{\eps_2}$; this is good enough to bootstrap both~$\epsilon_1$
and~$\epsilon_2$.

The ninth term in~\eqref{eq:epsilonwave} is split as follows, see \eqref{eq:Fthree}:
\begin{align*}
& P_{0}Q_{<D}F_\alpha^{3}(\tilde{\epsilon}; P_{\geq-D_{1}}; \psi_L, \psi_L) \\
&= i\partial^{\beta}P_0Q_{<D}[\tilde{\epsilon}_\alpha \, I
\del_j^{-1}P_{\geq-D_{1}}\calN_{\beta j}(\psi_L,\psi_L)] -i P_{0}Q_{<D}[P_0 R_\beta \tilde{\epsilon}\,
\,\del_j^{-1} I  \partial^{\beta}P_{-D_{1}\leq\cdot<-5}\calN_{\alpha j}(\psi_L,\psi_L)]\\
&-i \partial^{\beta}P_0Q_{<D}[P_{>0}R_\beta \tilde{\epsilon}\, \del_j^{-1} I
P_{>0} \calN_{\alpha j}(\psi_L,\psi_L)]- i
\partial^{\beta} P_0Q_{<D}[P_{<-5}R_\beta \tilde{\epsilon}\, \del_j^{-1} I
P_{0} \calN_{\alpha j}(\psi_L,\psi_L)] \\
&
+i P_{0}Q_{<D}[
 P_0\partial^{\beta}\tilde{\epsilon}_\alpha \,\del_j^{-1} I
P_{-D_{1}\leq \cdot<-5}\calN_{\beta j}(\psi_L,\psi_L)] +i P_{0}Q_{<D}[ P_0R^\beta \tilde{\epsilon}\,
\,\del_j^{-1} I \partial_{\alpha}P_{-D_{1}\leq\cdot <-5}\calN_{\beta j}(\psi_L,\psi_L)] \\
&+i \partial_{\alpha}P_0Q_{<D}[P_{>0}R^\beta \tilde{\epsilon}\, \del_j^{-1} I
P_{>0}\calN_{\beta j}(\psi_L,\psi_L)]+i
\partial_\alpha P_0 Q_{<D}[ P_{<-5}R^\beta \tilde{\epsilon}\, \del_j^{-1} I
P_{0} \calN_{\beta j}(\psi_L,\psi_L)]
\end{align*}
Denote these terms in this order by $\term_{91}$ through $\term_{98}$. First, we rewrite $\term_{91}$
in terms of the usual trichotomy:
\begin{equation}\label{eq:term91split}
\begin{aligned}
 & i\partial^{\beta}P_0Q_{<D}[\tilde{\epsilon}_\alpha \, I
\del_j^{-1}P_{\geq-D_{1}}\calN_{\beta j}(\psi_L,\psi_L)]
=  i\partial^{\beta}P_0Q_{<D}[P_0\tilde{\epsilon}_\alpha \, I
\del_j^{-1}P_{-D_{1}\le \cdot<-5}\calN_{\beta j}(\psi_L,\psi_L)] \\
&+  i\partial^{\beta}P_0Q_{<D}[P_{<-5}\tilde{\epsilon}_\alpha \, I
\del_j^{-1}P_{0}\calN_{\beta j}(\psi_L,\psi_L)] +   i\partial^{\beta}P_0Q_{<D}[P_{>0}\tilde{\epsilon}_\alpha \, I
\del_j^{-1}P_{>0}\calN_{\beta j}(\psi_L,\psi_L)]
\end{aligned}
\end{equation}
The first term in~\eqref{eq:term91split} is rewritten as  the sum
\begin{align*}
 &\partial^{\beta}P_0Q_{<D}[P_0\tilde{\epsilon}_\alpha \, I
\del_j^{-1}P_{-D_{1}\le \cdot<-5}\calN_{\beta j}(\tilde\psi_L,\tilde\psi_L)] + \partial^{\beta}P_0Q_{<D}[P_0\tilde\epsilon_\alpha \, I
\del_j^{-1}P_{-D_{1}\le \cdot<-5}\calN_{\beta j}(\tilde\psi_L,\breve\psi_L)] \\
&+ \partial^{\beta}P_0Q_{<D}[P_0\tilde\epsilon_\alpha \, I
\del_j^{-1}P_{-D_{1}\le \cdot<-5}\calN_{\beta j}(\breve\psi_L,\tilde\psi_L)] + \partial^{\beta}P_0Q_{<D}[P_0\tilde\epsilon_\alpha \, I
\del_j^{-1}P_{-D_{1}\le \cdot<-5}\calN_{\beta j}(\breve\psi_L,\breve\psi_L)]
\end{align*}
where we followed the notation of Corollary~\ref{cor:psiLNLsupp}.  Each of the terms containing $\breve\psi$ is bootstrapped
easily, using the smallness of~$\tilde\epsilon_\alpha$ and Lemma~\ref{lem:tri_hyp}. Rescaling and square-summing these contributions are placed in~$\epsilon_2$;
alternatively, one can recover the frequency envelope using the smallness of~$\delta_1$ for the bootstrap. For the first term,
we proceed as in (b) of Case~1. More precisely, using the smallness of~$\delta_0$ (and the fact that the Besov smallness of $\psi$
at the edges of the intervals $J_\ell$ inherits itself to~$\tilde\psi$) as well as the frequency evacuation property for large~$n$, one obtains
that
\[
 \|\partial^{\beta}P_0Q_{<D}[P_0\tilde\epsilon_\alpha \, I
\del_j^{-1}P_{-D_{1}\le \cdot<-5}\calN_{\beta j}(\tilde\psi_L,\tilde\psi_L)] \|_{N[0]} \les C_4^3 C_2^2 \delta_0 \, d_0
\]
As usual, this gets turned into an $S$ bound, leading to an $\epsilon_2$ contribution.  If $\tilde\epsilon=\tilde\epsilon_2$, then it again
suffices to consider $\tilde\psi$. In this case, one needs to gain extra smallness by partitioning $I_1$ further; however, the number of intervals
needed for this partition only depends on the energy in an absolute way (i.e., not on the stage of the induction). First, we may assume
that there is angular separation between the Fourier supports of the two $\tilde \psi_L$ inputs due to the bound
\[
 \sum_{\substack{\kappa_1,\kappa_2\in \calC_{-m_0}\\ \dist(\kappa_1,\kappa_2)\les 2^{-m_0}}} \| \partial^{\beta}P_0Q_{<D}[P_0\tilde\epsilon_\alpha \, I
\del_j^{-1}P_{-D_{1}\le \cdot<-5}\calN_{\beta j}(P_{\kappa_1}\tilde\psi_L, P_{\kappa_2}\tilde\psi_L)] \|_{N[0]} \ll  \|P_0\tilde\epsilon\|_{S[0]}
\]
see Section~\ref{subsec:improvetrilin}. Here we used that $\|\psi_L\|_S^2$ is bounded by the energy in an absolute way, which allows
us to chose~$m_0$ in the same fashion.  On the other hand, the remaining term
\begin{align*}
 & \sum_{\substack{\kappa_1,\kappa_2\in \calC_{-m_0}\\ \dist(\kappa_1,\kappa_2)> 2^{-m_0}}} \| \partial^{\beta}P_0Q_{<D}[P_0\tilde\epsilon_\alpha \, I
\del_j^{-1}P_{-D_{1}\le \cdot<-5}\calN_{\beta j}(P_{\kappa_1}\tilde\psi_L, P_{\kappa_2}\tilde\psi_L)] \|_{N[0]}
\end{align*}
is estimated by placing $\calN_{\beta j}(P_{\kappa_1}\tilde\psi_L, P_{\kappa_2}\tilde\psi_L)$ into~$L^2_{t,x}$, see the reasoning leading
up to~\eqref{eq:L2dxdt}, followed by a decomposition of the interval of integration. Here is important to note that $D_1$ only depends on
the energy.

For the second term in~\eqref{eq:term91split} consider first $\tilde\epsilon_1$; then the frequency envelope of $\tilde\epsilon_1$ is inherited by this expression.
More precisely, for $P_k \tilde\epsilon$ one gains a weight $2^{-\sigma k}$ from Lemma~\ref{lem:tri_hyp} which is sufficient for
the bootstrap provided $k$ is sufficiently large and negative; if not, then one applies the same fungibility as for the previous terms.
the same reasoning applies to~$\tilde\epsilon_2$.

Finally, for the third term in~\eqref{eq:term91split}  consider first the contribution by $\tilde\epsilon_1$. In that case one has
\begin{equation}
 \label{eq:91HH}
 \| \partial_{\alpha}P_0Q_{<D}[P_{>0}R^\beta\tilde\epsilon\, \del_j^{-1} I
P_{>0}\calN_{\beta j}(\psi_L,\psi_L)] \|_{N[0]} \les \sup_{k_1,k_2>0} 2^{-\sigma_0|k_1-k_2|}\|P_{k_1} \tilde\epsilon\|_{S[k_1]}  \|P_{k_2}\psi_L\|_{S[k_2]}
\end{equation}
which can be made $\ll C_4C_2\,\delta_1$ by choosing $\delta_0$ small and $n$ large.
On the other hand, if $\tilde\epsilon=\tilde\epsilon_2$, then one gains smallness in two ways: if  any one of $\tilde\epsilon$, or the two $\psi_L$ inputs
has large frequency, then ones gains smallness from the weight $w$ in Lemma~\ref{lem:tri_hyp}. If the three inputs have frequency of size~$O(1)$,
then one gains smallness by fungibility as before.

Next, we note that $\term_{92}$ is treated in the same fashion as the first term on the right-hand side of~\eqref{eq:term91split}.

The terms~$\term_{93}$ and~$\term_{97}$ are of the high-high type and are estimated exactly as in~\eqref{eq:91HH}, and $\term_{94}$, $\term_{99}$ are essentially the same
as the low-high term on the righ-hand side of~\eqref{eq:term91split}. To bound~$\term_{95}$ and~$\term_{96}$ one applies the same fungibility considerations
as in the high-low case of $\term_{91}$.

The tenth and eleventh terms in~\eqref{eq:epsilonwave} are essentially the same so it suffices to estimate the former. Since the details
are quite similar to the preceding arguments, we will proceed schematically. Beginning with $\tilde\epsilon=\tilde\epsilon_1$, we split
\begin{equation}
 \label{eq:term10a}\begin{aligned}
P_{0}Q_{<D}F_\alpha^{3}(\psi,\psi,\tilde{\epsilon})& = P_{0}Q_{<D}F_\alpha^{3}(\psi_L,\psi_L,\tilde{\epsilon}) + P_{0}Q_{<D}F_\alpha^{3}(\psi_L,\psi_{NL},\tilde{\epsilon})\\
&+ P_{0}Q_{<D}F_\alpha^{3}(\psi_{NL},\psi_L,\tilde{\epsilon}) + P_{0}Q_{<D}F_\alpha^{3}(\psi_{NL},\psi_{NL},\tilde{\epsilon})
\end{aligned}\end{equation}
and furthermore, using Corollary~\ref{cor:psiLNLsupp},
\begin{equation}
 \label{eq:term10b}\begin{aligned}
 P_{0}Q_{<D}F_\alpha^{3}(\psi_L,\psi_L,\tilde{\epsilon}) & = P_{0}Q_{<D}F_\alpha^{3}(\tilde\psi_L,\tilde\psi_L,\tilde{\epsilon}) + P_{0}Q_{<D}F_\alpha^{3}(\tilde\psi_L,\breve\psi_L,\tilde{\epsilon}) \\
& +P_{0}Q_{<D}F_\alpha^{3}(\breve\psi_L,\tilde\psi_L,\tilde{\epsilon}) + P_{0}Q_{<D}F_\alpha^{3}(\breve\psi_L,\breve\psi_L,\tilde{\epsilon})
\end{aligned}\end{equation}
All terms here are going to be placed into the $S$ error~$\epsilon_2$ since they inherit the frequency envelope of~$\psi$. The
trilinear estimates of Lemma~\ref{lem:tri_hyp} allows for this, with the required smallness for the terms containing~$\psi_{NL}$ is
gained by the smallness of $\|\psi_L\|_S \|\psi_{NL}\|_S$.  Furthermore, the terms containing~$\breve\psi$ are easy due to
the smallness $\|\breve \psi\|_S< C_2\,\delta_1$ and the bootstrap assumption on~$\tilde\epsilon$ (one then chooses $\eps_0$ small enough).
The most interesting term here is $P_{0}Q_{<D}F_\alpha^{3}(\tilde\psi_L,\tilde\psi_L,\tilde{\epsilon})$. To place it in~$\epsilon_2$ one
uses the same small Besov/frequency evacuation logic that we have used several times before.

The twelfth term in~\eqref{eq:epsilonwave} is easy since it inherits
the  frequency envelope of~$\tilde\epsilon$ and basically bootstraps
itself.

The thirteenth term has to be placed entirely into the $S$-error
$\epsilon_2$. This can be done using the high-high gain in
Lemma~\ref{lem:tri_hyp} and in~\eqref{eq:JoachimIc} as demonstrated
several times before.

Finally, the fourteenth term is the cubic one which is again easy. This concludes obtaining the bootstrap for $Q_{<D}\epsilon$.
We now need to do the same thing for
\[
Q_{\geq D}\epsilon
\]
Since this is a technical repetition of similar reasoning, we again defer this to the appendix.
This concludes the proof of Proposition~\ref{PsiBootstrap}.
\end{proof}

It is now easy to conclude the proof of~\ref{PsiInduction}. More precisely, as indicated in Figure~6, one  proceeds in the direction of increasing time by passing
from $I_1$ to $I_2$ and so on. Corollary~\ref{cor:epsenergyconservation} guarantees that the energy does not increase
without bounds in the proccess, in fact, it remains always $<2\eps_0$ provided $\delta_1$ is sufficiently small.
Even though $\epsilon$ is initially only defined locally, Proposition~\ref{BlowupCriterion} and the $\|\cdot\|_S$-norm
bound of Proposition~\ref{PsiBootstrap} imply that $\epsilon$ exists globally with the bounds stated in Proposition~\ref{PsiInduction},
see~\eqref{eq:C2envelope2} and~\eqref{eq:Snorm2}.  \end{proof}

\begin{proof}[Proof of Proposition~\ref{ControlNonatomicComponent1}]
This follows simply by iterating Proposition~\ref{PsiBootstrap}, i.e., by passing from $J_1$ to $J_2$ and so forth in Figure~6.
Even though the constant~$C_2$ increases with~$\ell$,
in the end one obtains a bound of the form~\eqref{eq:PsiellS}.  The final statement~\eqref{eq:firstPsitails} is a consequence
of our proof of Proposition~\ref{PsiBootstrap} due to the frequency evacuation of the first Besov error from the atom $\phi_n^{a}$.
In fact, our estimates are based on control of the frequency envelope which therefore implies~\eqref{eq:firstPsitails} at all
stages of the induction.
\end{proof}

\subsection{Completion of the proofs of Lemma~\ref{BasicStability} and Proposition~\ref{prop:ener_stable}}\label{sec:missing proofs}

We commence by proving the final assertion of Lemma~\ref{BasicStability}. This follows exactly as in the proof of
Proposition~\ref{PsiBootstrap}, with $\epsilon_2=0$ and $\psi=\epsilon_1$. Note that we never need to place the nonlinear
source terms into the 'small term' $\epsilon_2$ (which is not present in this situation), since the outputs always inherit the frequency envelope of $\psi$.
\\
Next, the proof of Proposition~\ref{prop:ener_stable} follows again as in the proof of Proposition~\ref{PsiBootstrap},
but this time with $\epsilon_1=0$, $\epsilon_2=\epsilon$.

\subsection{Step 4: Adding the first large atomic component; preparing the second stage of Bahouri Gerard}

Recall  from Section~\ref{subsec:BGstep1} that we wrote the data
$\phi^{n}_{\alpha}$ of the essentially singular sequence (at time
$t=0$) in the form
\[
\phi^{n}_{\alpha}=\sum_{a=1}^{A_{0}}\phi^{na}_{\alpha}+w^{nA_{0}}_{\alpha},
\]
where $A_{0}$ was chosen such that the sum
\[
\limsup_{n\to\infty}\sum_{a\geq
A_{0}+1}\|\phi^{na}_{\alpha}\|_{L_{x}^{2}}^{2}\ll \eps_{0}
\]
As before $\eps_{0}(\Ecrit)>0$ is an absolute constant that depends
only on the energy. Then recall from Section~\ref{subsec:BGstep2}
that the atoms $\phi^{na}_{\alpha}$ ``split'' the error term
$w^{nA_{0}}_{\alpha}$ into finitely many pieces
$w^{nA_{0}^{(i)}}_{\alpha}$, $0\leq i\leq A_{0}$, ordered by the
size of $|\xi|$ in their Fourier support. Of course our eventual
goal is to describe the evolution of the Coulomb components (with
$\phi^{n}_{\alpha}=\phi^{n1}_{\alpha}+i\phi^{n2}_{\alpha}$)
\[
\psi^{n}_{\alpha}=\phi^{n}_{\alpha}e^{-i\sum_{k=1,2}\triangle^{-1}\partial_{k}\phi^{n1}_{k}}
\]
Our strategy then is to construct ``intermediate wave maps''
bootstrapping the bounds from one to the next, starting with the low
frequency ones to the higher frequency ones. In the previous
section, we have shown that we can derive an apriori bound
\[
\|\Psi_\alpha^{nA_{0}^{(0)}}\|_S=\big\|\Phi^{nA_{0}^{(0)}}_{\alpha}e^{-i\sum_{k=1,2}\triangle^{-1}\partial_{k}\Phi^{nA_{0}^{(0)}1}_{k}}\big\|_{S}<C_{10}(\Ecrit)
\]
provided we choose $A_{0}$ above large enough and also pick $n$
large enough. Moreover, we can then prove frequency localized bounds
of the form
\[
\big\|P_{k}\big[\Phi^{nA_{0}^{(0)}}_{\alpha}e^{-i\sum_{k=1,2}\triangle^{-1}\partial_{k}\Phi^{nA_{0}^{(0)}1}_{k}}\big]\big\|_{S[k](\R^{2+1})}\leq
C_{11}(\Ecrit)c_{k}
\]
for a suitable frequency envelope $c_{k}$ with
$\sum_{k\in\Z}c_{k}^{2}\leq 1$, say, and $c_{k}$ rapidly decaying
for $k\notin (-\infty, \log(\lambda_{n}^{1})^{-1})$, where the
frequency scales of the $\phi^{na}$ are given by
$(\lambda_{n}^{a})^{-1}$.

\noindent
 We now pass to the next approximating map, with data given by
\begin{align*}
[w^{nA_{0}^{(0)}}_{\alpha}+\phi^{n1}_{\alpha}]e^{-i\sum_{k=1,2}\triangle^{-1}\partial_{k}[w^{nA_{0}^{(0)}1}_{k}+\phi^{n1}_{k}]}&=
w^{nA_{0}^{(0)}}_{\alpha}e^{-i\sum_{k=1,2}\triangle^{-1}\partial_{k}[w^{nA_{0}^{(0)}1}_{k}+\phi^{n1}_{k}]}\\
&\quad
+\phi^{n1}_{\alpha}e^{-i\sum_{k=1,2}\triangle^{-1}\partial_{k}[w^{nA_{0}^{(0)}1}_{k}+\phi^{n1}_{k}]}
\end{align*}
Here the first component satisfies
\[
w^{nA_{0}^{(0)}}_{\alpha}e^{-i\sum_{k=1,2}\triangle^{-1}\partial_{k}[w^{nA_{0}^{(0)}1}_{k}+\phi^{n1}_{k}]}
=w^{nA_{0}^{(0)}}_{\alpha}e^{-i\sum_{k=1,2}\triangle^{-1}\partial_{k}w^{nA_{0}^{(0)}1}_{k}}+o_{L^{2}}(1)
\]
as $n\to\infty$ since $w^{nA_{0}^{(0)}}_{\alpha}$ is singular with
respect to the scale of~$\phi^{nA_{0}^{(0)}1}_{k}$. Technically
speaking, this follows by means of  the usual trichotomy
considerations.
 We now need to understand
the lack of compactness of the large added term
\[
\tilde{\psi}^{na}_{\alpha}=\tilde{\psi}^{n1}_{\alpha}:=\phi^{n1}_{\alpha}e^{-i\sum_{k=1,2}\triangle^{-1}\partial_{k}[w^{nA_{0}^{(0)}1}_{k}+\phi^{n1}_{k}]},
\]
which is where the second phase of Bahouri-Gerard needs to come in.

We now normalize via re-scaling to $\lambda_{n}^{1}=1$. This means
now that the frequency support of $ \tilde{\psi}^{na}_{\alpha} $ is
uniformly concentrated around frequency $|\xi|\sim 1$. Observe that
here {\em{we cannot get rid of the phase
$e^{-i\sum_{k=1,2}\triangle^{-1}\partial_{k}w^{nA_{0}^{(0)}1}_{k}}$}},
which may indeed ``twist'' the Coulomb components additionally. This
will have a negligible effect, however, since the $\psi$-system
\eqref{eq:psisys1}--\eqref{eq:psi_wave} is invariant with respect to
the modulation symmetry $\psi\mapsto e^{i\gamma}\psi$.

\noindent  For technical reasons\footnote{This has to do with the
fact that the energy of the free wave equation involves a
derivative.}, we now apply a Hodge type decomposition to the
components $\psi^{na}_{1,2}$ (here $1,2$ refer to the derivatives on
$\R^{2}$ with respect to the two coordinate directions), as well as
for $\tilde{\psi}^{na}_{1,2}$. Thus write
\begin{align}
\phi^{na}_{1}
&=\partial_{1}\tilde{\phi}^{na}+\partial_{2}{\mathring{\phi}}^{na}
\label{eq:yetanotherHodge1}\\
\phi^{na}_{2}&=\partial_{2}\tilde{\phi}^{na}-\partial_{1}{\mathring{\phi}}^{na}
\label{eq:yetanotherHodge2}\\
\tilde{\psi}^{na}_{1}&=\partial_{1}\zeta^{na}+\partial_{2}\eta^{na}
\label{eq:yetanotherHodge3}\\
\tilde{\psi}^{na}_{2}&=\partial_{2}\zeta^{na}-\partial_{1}\eta^{na}\label{eq:yetanotherHodge4}
\end{align}
More precisely, we define the components $\tilde{\phi}^{na},
{\mathring{\phi}}^{na}$, $\zeta^{na}, \eta^{na}$. All of this is at
time $t=0$, of course. Now following the procedure of the preceding
section, using the bound
\[
\|\Psi^{nA_{0}^{(0)}}_{\alpha}\|_{S}<C_{10}(\Ecrit),
\]
we can select finitely many intervals $I_{j}$ (whose number depends
on $C_{10}(\Ecrit)$) such that
\begin{equation}
\label{eq:psindecomp}
\Psi^{nA_{0}^{(0)}}|_{I_{j}}=\Psi^{nA_{0}^{(0)}}_{jL}+\Psi^{nA_{0}^{(0)}}_{jNL}
\end{equation}
for each interval $j$, see Corollary~\ref{cor:localsplit2}.
Moreover, it is straightforward to verify that our normalization
$\lambda_{n}^{1}=1$ implies that $|I_{j}|\to\infty$ as
$n\to \infty$; indeed, this follows from
$L^{\infty}$-bounds.

\noindent Next, pick the interval $I_{1}$ containing the initial
time slice $t=0$.  Consider the magnetic potential  (note that we do
not use the Hodge decomposition here)
\[
A^{n}_{\nu}:=\sum_{j=1,2}\triangle^{-1}\partial_{j}[\Psi^{1nA_{0}^{(0)}}_{\nu}\Psi^{2nA_{0}^{(0)}}_{j}-\Psi^{1nA_{0}^{(0)}}_{j}\Psi^{2nA_{0}^{(0)}}_{\nu}]
\]
Here we restrict everything to a non-resonant situation, i.e., we
shall replace the above by
\begin{equation}\label{eq:Anudef}
\begin{split}
A^{n}_{\nu}=&\sum_{\substack{\kappa_{1,2}\in
K_{-\Lambda_{n}}\\\text{dist}(\kappa_{1,2})\gtrsim
2^{-\Lambda_{n}}}}\;\sum_{|k-k_{1}|<\Lambda_{n},
|k_{1}-k_{2}|<\Lambda_{n}}\\
&\sum_{j=1,2}\triangle^{-1}\partial_{j}P_{k}[P_{k_{1},\kappa_{1}}\Psi^{1nA_{0}^{(0)}}_{\nu}P_{k_{2},\kappa_{2}}\Psi^{2nA_{0}^{(0)}}_{j}
-P_{k_{1},\kappa_{1}}\Psi^{1nA_{0}^{(0)}}_{j}P_{k_{2},\kappa_{2}}\Psi^{2nA_{0}^{(0)}}_{\nu}]
\end{split}
\end{equation}
Here we shall let $\Lambda_{n}\to\infty$ as
$n\to\infty$ sufficiently slowly. The errors thereby
generated shall be treatable as perturbative errors. {\it{This time
we use the full~$\Psi$, and not just the free wave part.}} Our
notation is somewhat inconsistent, since we do not include
$A_{0}^{(0)}$. Since we keep this parameter fixed throughout
this section, this omission will be inconsequential. From now on we shall denote $\Psi^{1nA_{0}^{(0)}}_{\nu}=\Psi^{1n}_{\nu}$ etc.\ to simplify the notation.

\begin{defi}
\label{def:cov_box}
 The  covariant wave operator $\Box_{A^n}$ is defined via
\[
\Box_{A^{n}}u:=\Box
u+2i\partial^{\nu}uA^{n}_{\nu}+iu\partial^{\nu}A^{n}_{\nu}
\]
\end{defi}

The fundamental fact about this operator is that solutions obeying
$\Box_{A^{n}}u=0$ preserve the energy in the limit
$n\to\infty$. This will allow us to modify the second stage of
the Bahouri-Gerard method to the covariant d'Alembertian instead of
the ``flat'' d'Alembertian. We state this rigorously as follows.

\begin{lemma}\label{BoxAEnergy}  Assume that $u$ is essentially supported at frequency $1$,
and that $A_{\nu}$ is essentially supported at frequencies $\ll 1$. By this we mean that
\begin{equation}
 \lim_{R\to\infty}\|P_{[-R, R]^{c}}u[0]\|_{\dot{H}^{-1}}=0
\label{eq:oscill1}
\end{equation}
as well as
\begin{equation}
\label{eq:Aevac}
 \lim_{n\to\infty}\|P_{>-R}\Psi^{1n}\|_{L_{x}^{2}}=0
\end{equation}
for any $R>0$.
If $u$ solves
\[
\Box_{A^{n}}u=0,\,u[0]=(\partial_{t}u,\nabla_{x}u)=(u_{0}, u_{1})\in
L^{2}\times L^{2}
\]
then one obtains a global bound (uniformly in the implicit
$\Lambda_{n}$)
\[
\|u\|_{S(\R^{2+1})}\lesssim \|u[0]\|_{L_{x}^{2}}
\]
with implied constant depending on $\Ecrit$ (which controls $A^{n}$),
and we can conclude that
\[
\|\partial_{t}u(t,\cdot)\|_{L_{x}^{2}}^{2}+\|\nabla_{x}u(t,
\cdot)\|_{L_{x}^{2}}^{2}=\|u_{0}(t, \cdot)\|_{L_{x}^{2}}^{2}+\|u_{1}(t,
\cdot)\|_{L_{x}^{2}}^{2}+o_{L^{2}}(1)
\]
as $n\to\infty$, uniformly in $t\in\R$, provided
$\Lambda_{n}\to\infty$ sufficiently slowly.
\end{lemma}
\begin{proof}
 This follows by the same argument that we used to prove Proposition~\ref{TwistedWaveEquation}.
\end{proof}

In our applications of Lemma~\ref{BoxAEnergy}, \eqref{eq:oscill1} will hold due to the frequency
localization inherent in our construction of the atoms; in other words, $u$ will be $1$-oscillatory after rescaling.
The other condition~\eqref{eq:Aevac} will hold due to~\eqref{eq:firstPsitails}, at least at the first stage of
the construction (i.e., when adding the first atom as we are doing here). For $a=2$ etc.~we will use the exact
same frequency evacuation property which gave rise to~\eqref{eq:firstPsitails} in the first place.

\subsubsection{Dispersion for the covariant wave equation}

In this section we prove a weak form of dispersion for the initial value
problem
\begin{equation}\label{eq:BoxAnCP}
 \Box_{A^n} u =0,\quad u[0]:=(f,g)
\end{equation}
where $\Box_{A^n}$ is as in Definition~\ref{def:cov_box}. For simplicity, we first
consider the case where $A_\nu^n$ is defined as in~\eqref{eq:Anudef} but with {\em free waves}~$\Psi_L$.
We shall assume that $(f, g)$, whence
also $u$ by Lemma~\ref{lem:F1F2}, are essentially supported at frequency~$1$, see Lemma~\ref{BoxAEnergy}.
Generally speaking,  $u$  depends  on $n$ away from the time $t_{0}=0$, but the above limit is uniform in~$n$ and holds on any time-slice.
We assume that the free waves $\psi_L$ satisfy
\begin{equation}
\label{eq:Aevac2}
 \lim_{n\to\infty}\|P_{>-R}\psi_L\|_{L_{x}^{2}}=0
\end{equation}
for any $R>0$.
We now claim the following main result of this subsection for the covariant wave equation~\eqref{eq:BoxAnCP}. For simplicity, we drop~$n$ as
a superscript.

\begin{prop}
 \label{prop:covdisp} Let $u$ be a solution of~\eqref{eq:BoxAnCP}.
Given $\gamma>0$,  there exists a decomposition
\[
 u=u_1+u_2
\]
with the following properties:
\begin{itemize}
 \item $u_{1,2}$ satisfy the same apriori estimates which were proved for~$u$ in Lemma~\ref{lem:F1F2}
\item $\|u_2\|_S<\gamma $
\item there exists
$t_0=t_0(\gamma,f,g,\Ecrit)$ (but $t_0$ does not depend otherwise on~$\psi_L$) such that
for $|t|>t_0$ one has that
\begin{equation}
 \label{eq:covdisp} \|u_1(t,\cdot)\|_{L^\infty_x} <\gamma,
\end{equation}
uniformly for large enough $n$.
\end{itemize}
\end{prop}

The proof of this result will be split into several pieces. The idea is to
first obtain a ``parametrix'' for $u$, which is established by restricting to suitable time intervals (this is done via ``fungibility''). Once we
have such a parametrix (more precisely, a representation of $u$ as a sum of Volterra iterates starting with the free wave), we can use the dispersion
of the wave equation  to prove the desired result. First, we follow Tao to establish
the following fungibility lemma.

\begin{lemma}
 \label{lem:Taofung} For any $\eps_1>0$ there exist a partition of $\R$ into intervals $\{I_j\}_{j=1}^M$
where $M\les (\Ecrit\eps_1^{-1})^C$ for some absolute constant~$C$
with the property that for any~$u$
\[
\max_{1\le j\le M} \| \del^\alpha u\:A_\alpha\|_{N(I_j\times\R^2)}\le  \eps_1\|u\|_{S}
\]
Note that the intervals depend on~$\psi_L$ (but not on~$u$), but
their number does not (other than through the energy).
\end{lemma}
\begin{proof}
According to the trilinear estimates of Section~\ref{sec:trilin}, we may assume that there is angular
separation between $\hat{u}$ and the waves in~$A_\alpha$. Otherwise there is the desired gain.
The amount of angular separation is very small and
depends on~$\Ecrit$ and~$\eps_1$. We shall now implicitly assume that $\del^\alpha u\:A_\alpha$ respects this
type of angular separation. Note that we may restrict ourselves to the case of high-low interactions between $u$ and $A_{\alpha}$, since
for the other cases, the fungibility follows by using the same argument as in the proof of Lemma~\ref{lem:F1F2}.

\noindent  By \eqref{eq:bilin1},
\begin{align}
\nn \| \del^\alpha u\:A_\alpha \| &\le C(\Ecrit,\eps_1)  \sum_{k_1<-C} 2^{-\frac{k_1}{2}} \|P_0 u\, P_{k_1} \psi_L\|_{L^2_{t,x}} \| P_{k_2}\psi_L \|_{S}\\
&\le  C(\Ecrit,\eps_1)  \Big( \sum_{k_1<-C} 2^{-k_1}  \|P_0 u\, P_{k_1} \psi_L\|_{L^2_{t,x}}^2 \Big)^{\frac12} \| \psi_L \|_{S}
\label{eq:upsiL}
\end{align}
Next, by Theorem~1.11 of~\cite{T8}, assuming $u$ to be a free wave,
for each $k\in\Z$ there exists a collection~$\calT_k $ of
tubes~$\tau_k^i$
  of size $\infty\times 2^k\times 2^k$ centered along a light-ray and aligned with the Fourier
support of~$u$  such that $\#\calT_k \le (\Ecrit\eps_2^{-1})^C$ and so that, where $\eps>0$ is small and will be determined,
\begin{equation}
 \label{eq:Taotubes}
 \| P_0 u\, P_{k_1} \psi_L\|_{L^2_{t,x}\setminus \Omega_{k_1}} \le \eps \, 2^{\frac{k_1}{2}} \|u\|_2 \|P_{k_1}\psi_L\|_2
\end{equation}
where $\Omega_k:= \bigcup_{\tau\in\calT_{k_1}} \tau$. In our case $u$ is of course not a free wave; however, by
Remark~\ref{rem:slightimprov} as well as Remark~\ref{rem:quinticimprov} in conjunction with Lemma~\ref{lem:Null_rep}, we conclude that we can write
\[
 u=u_{1}+u_{2}
\]
where
\[
 \|u_2\|_S<\eps_2\|u\|_S
\]
while
\[
 u_1= \int  f_a u_a\, \nu( da)
\]
is a superposition of free waves $u_a$ with the same frequency support properties as $u$ and
\[
 \int \| f_a u_{a}\|_{L_{x}^{2}}\,\nu(da)\le C(\eps_2)\|u\|_{S},\quad f_{a}\in L_{t,x}^{\infty},\quad\|f_{a}\|_{L_{t,x}^{\infty}}\le C
\]
Thus for $u$ in the original sense, choosing $\eps$ in~\eqref{eq:Taotubes} of the form $C(\Ecrit,\eps_1)^{-1}\eps_2$, we get
\[
 \| P_0 u\, P_{k_1} \psi_L\|_{L^2_{t,x}\setminus \Omega_{k_1}} \le \eps_2\, 2^{\frac{k_1}{2}} \|u\|_{S[0]}  \|P_{k_1}\psi_L\|_2
\]
Inserting this bound in~\eqref{eq:upsiL} yields
\begin{align*}
  \| \del^\alpha u\:A_\alpha \|_{L^2_{t,x}} &\le C(\Ecrit,\eps_1) \Big( \sum_{k_1<-C} 2^{-k_1}
\|P_0 u\, P_{k_1} \psi_L\|_{L^2_{t,x}\setminus \Omega_{k_1}}^2 \Big)^{\frac12} \| \psi_L \|_{S} \\
& \quad + C(\Ecrit,\eps_1) \Big( \sum_{k_1<-C} 2^{-k_1} \sum_{\tau_k^i\in\calT_k} \|\chi_{\tau_k^i} P_0 u\, P_{k_1} \psi_L\|_{L^2_{t,x}}^2 \Big)^{\frac12}
\| \psi_L \|_{S} \\
&\le C(\Ecrit,\eps_1) \eps_2    \|u\|_{S[0]} \| \psi_L \|^2_{S} \\
& \quad + C(\Ecrit,\eps_1) \Big( \sum_{k_1<-C} 2^{-k_1} \sum_{\tau_{k_1}^i\in\calT_{k_1}} \|P_0 u\|_{\ener}^2
\| \chi_{\tau_k^i} \, P_{k_1} \psi_L\|_{L^2_t L^\infty_x}^2 \Big)^{\frac12} \| \psi_L \|_{S}
\end{align*}
Next, by a standard $TT^*$ estimate, and for all $k_1\in\Z$,
\[
 \| \chi_{\tau_k^i} \, P_{k_1} \psi_L\|_{L^2_t L^\infty_x} \les 2^{\frac{k_1}{2}} \| P_{k_1} \psi_L\|_2
\]
whence
\[
 \Big( \sum_{k_1\in\Z} 2^{-k_1} \sum_{\tau_{k_1}^i\in\calT_{k_1}}
\| \chi_{\tau_k^i} \, P_{k_1} \psi_L\|_{L^2_t L^\infty_x}^2 \Big)^{\frac12}  \les (\Ecrit\eps^{-1})^C \|\psi_L\|_2
\]
Therefore, the exist intervals $\{I_j\}_{j=1}^M$ as claimed.
Since the constants $C(\eps_2)$ and $C(\Ecrit,\eps_1)$ depend polynomially on the parameters, we are done.
\end{proof}

We can now prove Proposition~\ref{prop:covdisp}. We will assume that the energy of the data $(f,g)$ is also controlled by~$\Ecrit$ although
this is only a notational  convenience.

\begin{proof}[Proof of Proposition~\ref{prop:covdisp}]  For simplicity, we drop the zero order term from~$\Box_A$. This is admissible,
since it only presents a notational inconvenience and is amenable to the exact same arguments that we apply to the first order term.
 With $\{I_j\}_{1\le j\le M}$ as in the lemma, we relabel them as follows: with initial time $0\in I_{j_0}$, we set $J_0:= I_{j_0}$.
At the next step, we define $J_1=I_{j_1}$ and $J_{-1}:=I_{j_2}$
where $I_{j_1}$ is the successor of $I_{j_0}$ (with respect to
positive orientation of time), whereas $I_{j_2}$ is the predecessor.
In this fashion one obtains a sequence $J_{i}$ with $0\le i\le M'$
and $M'\le (\Ecrit\eps_1^{-1})^C$ as in Lemma~\ref{lem:Taofung}
where $\eps_1$ is small depending only on~$\Ecrit$. Next, let $u$ be
the solution of
\[
 \Box u+ 2i \del^\alpha u\, A_\alpha =0,\quad u[0]=(f,g)
\]
We claim that  $u^{(0)}:= u\big|_{J_0}$ can be written as an infinite Duhamel expansion in the form
\begin{align*}
 u^{(0)} &:= \sum_{\ell=0}^\infty u^{(J_0,\ell)}, \quad  u^{(J_0,0)}(t):= S(t)u[0],\\
 u^{(J_0,\ell)} &:= - 2i\int_0^t U(t-s)  \del^\alpha u^{(J_0,\ell-1)}\,
A_\alpha(s)\, ds
\end{align*}
where $S(t)=(U,V)(t)$ is the free wave evolution, and $U(t)=\frac{\sin(t|\nabla|)}{|\nabla|}$, $V(t)=\cos(t|\nabla|)$. Of course, $t\in J_0$ in this
equation.
Due to the energy estimate of Section~\ref{subsec:energy} and Lemma~\ref{lem:Taofung}, this series converges with respect to the $S$-norm.
In a similar fashion, we can pass to later times: $u^{(i)}:= u\big|_{J_i}$ satisfies
\begin{align}
 u^{(i)} &:= \sum_{\ell=0}^\infty u^{(J_i,\ell)}, \quad  u^{(J_i,0)}(t):= S(t-t_i)u^{(i-1)}[t_i], \nn \\
 u^{(J_i,\ell)} &= - 2i\int_{t_i}^t U(t-s)  \del^\alpha u^{(J_i,\ell-1)}\,
A_\alpha(s)\, ds   \label{eq:uJiell}
\end{align}
where $t\in J_i$ and $t_i:=\max J_{i-1}=\min J_i$ for $i\ge1$ and $t_0:=0$.  Observe that
\begin{equation}
 \label{eq:Duhi}
\begin{split}
u^{(J_i,0)}(t) &:= S(t-t_{i-1}) u^{(i-1)}[t_{i-1}]  \\
&\quad -2i \sum_{\ell=1}^\infty \int_{t_{i-1}}^{t_i} U(t-s) \chi_{J_{i-1}}(s) \del^\alpha u^{(J_{i-1},\ell)}(s) \, A_\alpha(s)\, ds
\end{split}
\end{equation}
for all $t\in J_i$. If $i\ge2$, we expand further to obtain
\begin{align*}
 S(t-t_{i-1}) u^{(i-1)}[t_{i-1}]  &:= S(t-t_{i-2}) u^{(i-2)}[t_{i-2}]  \\
&\quad -2i \sum_{\ell=1}^\infty \int_{t_{i-2}}^{t_{i-1}} U(t-s) \chi_{J_{i-2}}(s) \del^\alpha u^{(J_{i-2},\ell)}(s) \, A_\alpha(s)\, ds
\end{align*}
This procedure can be continued all the way back to $t_0=0$ and yields
\begin{equation}
\label{eq:Jiell0}
 u^{(J_i,0)}(t) := S(t)(f,g) - \sum_{k=0}^{i-1} 2i \sum_{\ell=1}^\infty \int_{t_{k}}^{t_{k+1}} U(t-s)
\chi_{J_{k}}(s) \del^\alpha u^{(J_{k},\ell)}(s) \, A_\alpha(s)\, ds
\end{equation}
for all $t\in J_i$. Inductively, one passes from this term to $ u^{(J_i,\ell)}$ for all $\ell\ge0$ by means of~\eqref{eq:uJiell}.
We next claim that for each $j$, the functions $u^{(J_{i},\ell)}$ become small with respect to $\|\cdot\|_S$ provided $\ell$ is large enough.
 This is a direct  consequence of applying Lemma~\ref{lem:Taofung} to the above iterative definition of $u^{(J_i,\ell)}$ as well as the basic energy estimate.

\noindent
Now fix a number $\gamma>0$. We will show that there exist $t_{0}=t_{0}(\gamma)$ and $n_{0}(\gamma)$ with the property that if $|t|>t_{0}(\gamma)$
and $n>n_{0}(\gamma)$, then we can write
\[
u=u_1+u_2
\]
where
\[
 \|u_2\|_S<\gamma
\]
and
\[
 |u_1(t, x)|<\gamma
\]
for $|t|>t_{0}$, uniformly in $n>n_{0}(\gamma)$.
We start by reducing ourselves to a double light cone. Indeed, pick a large enough disc $D_{\gamma}$ in the time slice $\{0\}\times\R^{2}$
with the property that \[
 \|\chi_{D_{\gamma}^{c}}u[0]\|_{L_{x}^{2}}\ll\gamma                                                                                                             \]
Here $\chi_{D_{\gamma}^{c}}$ is a smooth cutoff localizing to a large dilate of $D_{\gamma}$.
If we denote the covariant propagation of $\chi_{D_{\gamma}^{c}}u[0]$ by $\tilde u_2$, then we can achieve that
\[
 \|\tilde u_2\|_{S}\ll\gamma
\]
by means of Lemma~\ref{lem:F1F2}.
We are thus reduced to estimating $\tilde u_1 = u-\tilde u_2$, which by construction is supported in a (large) double
cone whose base depends only on $\gamma$. We can then expand $\tilde u_1$ in terms of Volterra iterates just as before,
and there exists $\ell_\gamma$ with the property that
\begin{equation}
\label{eq:summell}
 \sum_i\sum_{\ell >\ell_{\gamma}}\|\tilde u_1^{J_{i}, \ell}\|_{S}\ll\gamma
\end{equation}
Furthermore, note that all the iterates $\tilde u_1^{J_{i}, \ell}$ are supported in the same double light cone with base $D_{\gamma}$.
We now show that $\tilde u_1=u_1+u_2^\dagger$ where $\|u_2^\dagger\|_S\ll\gamma$ and $u_1$ has the desired dispersive property. Setting
$u_2:=\tilde u_2+ u_2^\dagger$ then concludes the argument. First, in view of~\eqref{eq:summell} and the fact that the total number of $J_i$
is controlled by the energy, we may include the contributions of $\ell>\ell_\gamma$ in~$u_2^\dagger$.

\noindent By Huyghens principle, $ \tilde u_1=\chi(t,x) \tilde u_1$ where for the remainder of the proof $\chi(t,x)$
is a smooth cut-off to the region $|x|\le |t|+\rho$ with $\rho$ being the radius of~$D_\gamma$. Then we can write
\begin{align}
 \tilde u_1^{(J_i,\ell)}(t) &= - 2i\chi(t,x)\int_{t_i}^t U(t-s) P_{[-k_0<\cdot<k_0]}\big[ \del^\alpha \tilde u_1^{(J_i,\ell-1)}\,
A_\alpha(s)\big]\, ds \label{eq:chiUt} \\
& - 2i\chi(t,x)\int_{t_i}^t U(t-s) P_{[-k_0<\cdot<k_0]^c}\big[ \del^\alpha \tilde u_1^{(J_i,\ell-1)}\,
A_\alpha(s)\big]\, ds \nn
\end{align}
We now show that the second integral splits into a term of small $\Linf$-norm, and one of small~$S$ norm.
First, consider $P_{[<-k_0]}$. Then by Bernstein's inequality, and the energy estimate
\begin{align}
 & \Big\| \chi(t,x)\int_{t_i}^t U(t-s) P_{[<-k_0]}\big[ \del^\alpha \tilde u_1^{(J_i,\ell-1)}\,
A_\alpha(s)\big]\, ds \Big\|_{L^\infty_x} \nn \\
&\les 2^{-k_0} \Big\| \int_{t_i}^t U(t-s) P_{[<-k_0]}\big[ \del^\alpha \tilde u_1^{(J_i,\ell-1)}\,
A_\alpha(s)\big]\, ds \Big\|_{\ener} \\
&\les 2^{-k_0} \Ecrit \|\tilde u_1^{(J_i,\ell-1)}\|_S \le C(\Ecrit) 2^{-k_0}
\end{align}
whereas for $P_{[>k_0]}$ one can essentially (up to tails which are handled by Lemma~\ref{lem:chiS}, for example) remove the exterior~$\chi$
since the interior~$\del^\alpha \tilde u_1^{(J_i,\ell-1)}$ obeys that very localization. In conclusion, the resulting term is placed in~$u_2^\dagger$.
Now consider the main term~\eqref{eq:chiUt}. Decompose $S^1$ into caps~$\kappa$ of size $c(\Ecrit,\gamma)$ which is a small constant. Denote
the corresponding decomposition of the double light-cone $\{|x|\le |t|+\rho\}$ into angular sectors by~$\{S_\kappa\}_\kappa$. Associated
with the $S_\kappa$ there is a smooth partition of unity~$\sum_\kappa\chi_\kappa=\chi$. Write \eqref{eq:chiUt}  as the sum
\begin{align}
 & 2i\sum_\kappa \chi_{\kappa}(t,x)\int_{t_i}^t U(t-s) P_{[-k_0<\cdot<k_0]} P_{[\hat\xi\in \mp 2\kappa]} \big[ \del^\alpha \tilde u_1^{(J_i,\ell-1)}\,
A_\alpha(s)\big]\, ds \label{eq:chikappa}\\
&+ 2i\sum_\kappa \chi_{\kappa}(t,x)\int_{t_i}^t U(t-s) P_{[-k_0<\cdot<k_0]} P_{[\hat\xi\not\in \mp 2\kappa]} \big[ \del^\alpha \tilde u_1^{(J_i,\ell-1)}\,
A_\alpha(s)\big]\, ds \label{eq:chikappa'}
\end{align}
Here $\hat\xi=\frac{\xi}{|\xi|}$ and the sign is selected according to the decomposition into incoming and
outgoing propagator:
\[
 U(t) = \frac{1}{2i|\nabla|}[e^{i t|\nabla|} - e^{-i t|\nabla|}]
\]
By Bernstein's inequality the first term~\eqref{eq:chikappa} satisfies
\[
 \| \eqref{eq:chikappa} \|_{\Linf} \le C(\Ecrit) |\kappa|^{\frac12} \| \tilde u_1^{(J_i,\ell-1)} \|_S
\]
which can be made small for small~$\kappa$.
Also note that
\[
 |t|\xi| \pm x\cdot\xi|\gtrsim |t| \quad \forall\; (t,x,\xi) \text{\ such that\ } \chi(t,x)\ne0, \; 2^{-k_0}<|\xi|< 2^{k_0}
\]
 where the choice of $\pm$ depends on whether the propagator $U$ is incoming or outgoing.
Now we make the inductive assumption (relative to $\ell$ and~$i$) that
\begin{equation}\label{eq:uinduc}
 \| [\chi_\kappa P_{[\pm\hat\xi\not\in -2\kappa]} \tilde u_1^{(J_i,\ell-1)} ](t)\|_{L^\infty_x}+
\|\chi_{[||t|-|x||\gtrsim |t|]}\tilde u_1^{(J_i,\ell-1)}(t,x)\|_{L_{x}^{\infty}} \le C_N(i,\ell,\Ecrit,\gamma) |t|^{-N}
\end{equation}
where the $\pm$ sign is according to whether the function has space-time Fourier support in the upper or lower half-spaces, i.e., whether $\tau>0$ or~$\tau<0$.
Strictly speaking, the cap size here depends on $(i,\ell)$ with the size $c(\Ecrit,\gamma)$ from above being the size at the end of the induction (recall
that there are only finitely many choices for these parameters). But for simplicity of notation, we suppress this dependence from the notation.
Note that we only have finitely many values of $\ell,i$. Now to estimate the second integral term \eqref{eq:chikappa'}, we distinguish
between a number of cases: first if $|s|\ll|t|$ (where the implicit small constant depends on~$|\kappa|$), due to the apriori support conditions satisfied by $\tilde u_1^{(J_i,\ell-1)}$ which forces $|y|<|s|+\rho$, we obtain the desired gain in~$t$ by integrating by parts with respect to~$|\xi|$.
Next, assume that $|s|\sim |t|$ (where the implicit small constant again depends on~$|\kappa|$ - this will be tacitly understood for the remainder
of the proof). Then we first reduce to $||s|-|y||\ll|s|$. For this consider the term
\[
 2i\sum_\kappa \chi_{\kappa}(t,x)\int_{t_i}^t U(t-s) P_{[-k_0<\cdot<k_0]} P_{[\hat\xi\not\in \mp 2\kappa]} \big[ \chi_{[||s|-|y||\gtrsim |s|]}\del^\alpha \tilde u_1^{(J_i,\ell-1)}\,
A_\alpha(s)\big]\, ds
\]
Since we assume $|s|\sim |t|$, the desired gain $t^{-N}$ here follows by using the induction hypothesis. Hence we now reduce to estimating
\[
 2i\sum_\kappa \chi_{\kappa}(t,x)\int_{t_i}^t U(t-s) P_{[-k_0<\cdot<k_0]} P_{[\hat\xi\not\in \mp 2\kappa]} \big[ \chi_{[||s|-|y||\ll|s|]}\del^\alpha \tilde u_1^{(J_i,\ell-1)}\,
A_\alpha(s)\big]\, ds
\]
Here we apply a further decomposition
\begin{align*}
\chi_{[||s|-|y||\ll|s|]}\del^\alpha \tilde u_1^{(J_i,\ell-1)}&=\chi_{[||s|-|y||\ll|s|]}\sum_{\kappa'}\chi_{\kappa'}(s,y)\del^\alpha \tilde u_1^{(J_i,\ell-1)}\\
&=\chi_{[||s|-|y||\ll|s|]}\sum_{\kappa'}\chi_{\kappa'}(s,y)P_{[\pm\hat{\xi}\in-2\kappa']}\del^\alpha \tilde u_1^{(J_i,\ell-1)}\\
&\quad +\chi_{[||s|-|y||\ll|s|]}\sum_{\kappa'}\chi_{\kappa'}(s,y)P_{[\pm\hat{\xi}\in-(2\kappa')^{c}]}\del^\alpha \tilde u_1^{(J_i,\ell-1)}
\end{align*}
The contribution of the second term here is again rapidly decaying due to the induction assumption. Hence we have now reduced to estimating
\begin{align*}
 2i\sum_\kappa \chi_{\kappa}(t,x)\int_{t_i}^t U(t-s) P_{[-k_0<\cdot<k_0]} &P_{[\hat\xi\not\in \mp 2\kappa]} \big[\chi_{[||s|-|y||\ll  |s|]}\\
&\sum_{\kappa'}\chi_{\kappa'}(s,y)P_{[\pm\hat{\xi}\in-2\kappa']}\del^\alpha \tilde u_1^{(J_i,\ell-1)}\,
A_\alpha(s)\big]\, ds
\end{align*}
Now writing out the free wave parametrix, we see that on the support of the resulting integral in the variables $\xi,\,y,\,s$, we have that
\[
 |\pm|\xi|s+y\cdot\xi|\ll |t|,
\]
and choosing $\kappa'$ as well as the implied constant in $||s|-|y||\ll|s|$ suitably small, we can ensure that
\[
 |\pm t|\xi|+x\cdot\xi|\sim |t|\gg |\pm |\xi|s+y\cdot\xi|
\]
on the support of the integrand. Integrations by parts in~$|\xi|$ yield the desired rapid decay with respect to $|t|$.
This recovers the first part of the inductive assumption,  and the second follows identically, since if $||t|-|x||\gtrsim |t|$, then we necessarily have
\[
 |\pm t|\xi| + x\cdot\xi|\gtrsim |t|
\]
The inductive procedure is now completed by means of \eqref{eq:Jiell0} which takes account of the changes in the level~$i$.
\end{proof}

Recall that we restricted $\Psi$ to be a free wave
in~\eqref{eq:Anudef}. In order to treat the general case, we apply
the usual decomposition~\eqref{eq:psindecomp}.  As usual, the
smallness of the $\Psi_{NL}$ allows one to iterate these terms away.
Furthermore, the proof of Proposition~\ref{prop:covdisp} applies to
these terms equally well since we do not rely on any specific
structure of the $u_1^{(J_i,\ell-1)}$ other than the inductive
assumption~\eqref{eq:uinduc}, and the formalism of the Volterra
iteration by which we represented these solutions.

\subsubsection{The second stage of Bahouri Gerard, applied to the first large atomic component}

Recall that we are considering only $a=1$. Nevertheless, we keep the parameter ``$a$'' in our notation general.
We now need to quantify the lack of compactness for the functions
$\tilde{\phi}^{na}$, ${\mathring{\phi}}^{na}$, $\zeta^{na}$,
$\eta^{na}$, all at time $t=0$. We
evolve each of these using the covariant wave flow from before and
select a number of concentration profiles. The method for this
follows exactly the Bahouri-Gerard template, but using Lemma~\ref{BoxAEnergy}
 instead of standard energy conservation for the free wave
flow. In order to define the temporal flow for each component, we
need to impose time derivatives at time $t=0$. We do this by
defining
\begin{align*}
\partial_{t}\tilde{\phi}^{na}(0, \cdot):=\phi^{na}_{0}(0,\cdot),\quad \partial_{t}{\mathring{\phi}}(0,\cdot):=0\\
\partial_{t}\zeta^{na}(0, \cdot):=\tilde{\psi}^{na}_{0},\quad \partial_{t}\eta^{na}(0,\cdot):=0
\end{align*}
Introduce the following terminology:

\begin{defi} Given data $u[0]=(u_{0}, u_{1})$ at time $t=0$, we denote by
\[
S_{A^{n}}(t)(u[0])
\]
the solution of $\Box_{A^{n}}(u)=0$ with the given data, evaluated
at time $t$.
\end{defi}

We now describe the important process of {\em {extraction of
concentration profiles:}} Consider $S_{A^{n}}(\zeta^{na}[0])$, with
$\zeta^{na}[0]=(\tilde{\psi}^{na}_{0}, \zeta^{na})$.
Following~\cite{BG} introduce the family
$\mathcal{V}_{A^{n}}(\underline{\zeta^{a}})$, consisting of all
functions on $V_{\zeta}(t, x)\in L_{t,\text{loc}}^{2}H_{x}^{1}\cap
C^{1}L_{x}^{2}$ such that
\[
\big(S_{A^{n}}(\zeta^{na}[0])\big)(t+t_{n}, x+x_{n})\rightharpoonup
V_{\zeta}(t, x)
\]
as $n\to\infty$ for some sequence $\{(t_{n},
x_{n})\}_{n=1}^\infty\in \R\times\R^{2}$. Here, the weak limit is in
the sense of $L_{t,\text{loc}}^{2}H_{x}^{1}$. Observe that such a
function $V_{\zeta}(t,x)$ solves $\Box V_{\zeta}=0$ in the sense of
distributions. Thus it makes sense to introduce the quantity
\[
\eta_{A^{n}}(\underline{\zeta^{a}}):=\sup\{
E(V_{\zeta}),\,V_{\zeta}\in
\mathcal{V}_{A^{n}}(\underline{\zeta^{a}})\},
\]
where
\[
E(V_{\zeta}):=\int_{\R^{2}}|\nabla_{x,t}V_{\zeta}|^{2} dx
\]
We can now state the following lemma that is at the core of the
second stage of the Bahouri-Gerard process for wave maps. Recall
that $a=1$ here.

\begin{lemma}\label{CoreBGII} There exists a collection of sequences $\{(t_{n}^{ab}, x_{n}^{ab})\}\subset \R\times\R^{2}$, $b\ge 1$, as
well as a family of concentration profiles $V_{\zeta}^{ab}[0]:=(V_{\zeta 0}^{ab}(x), V^{ab}_{\zeta 1}(x))\in L^{2}(\R^{2})\times \dot{H}^{1}(\R^{2})$,
 with the following properties: introducing the shifted gauge potentials
\begin{equation}\label{eq:Ashift}
\tilde{A}^{nab}:=A^{n}(t+t_{n}^{ab}, x+x_{n}^{ab}),
\end{equation}
one has \begin{itemize}
              \item   For any $B\geq 1$, one can write
\[
\big(S_{A^{n}}(\zeta^{na}[0])\big)(t,
x)=\sum_{b=1}^{B}\big(S_{\tilde{A}^{nab}}(V_{\zeta}^{ab}[0])\big)(t-t_{n}^{ab},
x-x_{n}^{ab})+W_{\zeta}^{naB}(t,x)
\]
Here each function
$\big(S_{\tilde{A}^{nab}}(V_{\zeta}^{ab}[0])\big)(t-t_{n}^{ab},
x-x_{n}^{ab})$, $W_{\zeta}^{naB}(t,x) $, solves the equation
$\Box_{A^{n}}u=0$, and we have
\begin{equation}
\lim_{B\to\infty}\big[\limsup_{n\to\infty}\eta_{A^{n}}(\underline{W}^{aB})\big]=0
\label{eq:Weta}
\end{equation}
\item One has  the divergence relations
\[
\lim_{n\to\infty}[|t_{n}^{ab}-t_{n}^{ab'}|+|x_{n}^{ab}-x_{n}^{ab'}|]=\infty
\]
for $b\neq b'$.
\item There is the asymptotic orthogonality relation
\[
E(\zeta^{na}[0])=\sum_{b=1}^{B}E(V_{\zeta}^{ab}[0])+E(W_{\zeta}^{naB}(t,\cdot))+o(1)
\]
Here $E$ refers to the standard (flat) energy and the $o$-term
satisfies $\lim_{B\to\infty}\limsup_{n\to\infty} o(1)=0$.
\item All $V^{ab}_\zeta[0]$, as well as their evolutions $S_{\tilde
A^{nab}}(V^{ab}_\zeta[0])$ and the $W_{\zeta}^{naB}$ are
$1$-oscillatory.
             \end{itemize}
\end{lemma}
\begin{proof} We follow \cite{BG}: There is nothing to do provided $\eta_{A^{n}}(\underline{\zeta^{a}})=0$.
 Hence assume this quantity is $>0$. Then pick a profile $V^{a1}_{\zeta}(t,x)\in  L_{t,\text{loc}}^{2}H_{x}^{1}\cap C^{1}L_{x}^{2}$
and associated sequence $\{(t^{a1}_{n}, x^{a1}_{n})\}_{n\geq 1}$ such that
\[
\big(S_{A^{n}}(\zeta^{na}[0])\big)(t+t^{a1}_{n},
x+x^{a1}_{n})\rightharpoonup V^{a1}_{\zeta}(t, x)
\]
with
\[
E(V^{a1}_{\zeta})>\frac12 \eta_{A^{n}}(\underline{\zeta^{a}})
\]
 Using the notation of the lemma, consider then
\begin{equation}\nonumber\begin{split}
&\big(S_{A^{n}}(\zeta^{na}[0])\big)(t+t^{a1}_{n}, x+x^{a1}_{n})-\big[S_{A^{n}}
\big(S_{\tilde{A}^{na1}}(V^{a1}_{\zeta}[0])(0-t_{n}^{a1}, \cdot-x_{n}^{a1})\big)\big](t+t_{n}^{a1}, x+x_{n}^{a1})\\
&=\big(S_{A^{n}}(\zeta^{na}[0])\big)(t+t^{a1}_{n},
x+x^{a1}_{n})-\big(S_{\tilde{A}^{na1}}(V^{a1}_{\zeta}[0])\big)(t,x)
\end{split}\end{equation}
But by our construction, this expression converges weakly to $0$.
\\
Furthermore, due to Lemma~\ref{BoxAEnergy}, we have that
\[
E(S_{\tilde{A}^{na1}}(V^{a1}_{\zeta}[0])(0-t_{n}^{a1},
\cdot-x_{n}^{a1}))=E(V^{a1}_{\zeta}[0])+o_{L^{2}}(1)
\]
Now we repeat the preceding step, but replace $\zeta^{na}[0]$ by
\[
\zeta^{na}[0]-S_{\tilde{A}^{na1}}(V^{a1}_{\zeta}[0])(0-t_{n}^{a1},
\cdot-x_{n}^{a1})
\]
Thus select a sequence $\{(t_{n}^{a2}, x_{n}^{a2})\}_{n\geq 1}$ and
a concentration profile $V_{\zeta}^{a2}(t,x)$ such that
\[
E(V_{\zeta}^{a2})\geq
\frac{1}{2}\eta(\zeta^{na}-S_{\tilde{A}^{na1}}(V^{a1}_{\zeta}[0])(0-t_{n}^{a1},
\cdot-x_{n}^{a1}))
\]
and furthermore
\[
\big[S_{A^{n}}\big(\zeta^{na}-S_{\tilde{A}^{na1}}(V^{a1}_{\zeta}[0])(0-t_{n}^{a1},
\cdot-x_{n}^{a1})\big)\big](t+t_{n}^{a2}, x+x_{n}^{a2})\rightharpoonup
V_{\zeta}^{a2}(t,x)
\]
We obtain that
necessarily
\[
\lim_{n\to\infty}|t_{n}^{a1}-t_{n}^{a2}|+|x_{n}^{a1}-x_{n}^{a2}|=\infty
\]
Furthermore, we claim that
\[
E(V^{a1}_{\zeta}[0])+E(\zeta^{na}[0]-S_{\tilde{A}^{na1}}(V^{a1}_{\zeta}[0])[0-t_{n}^{a1}])=E(\zeta^{na}[0])+o_{L^{2}}(1)
\]
This follows again just as in the free case, using
Lemma~\ref{BoxAEnergy}: We need to show that
\[
\int_{\R^{2}}\nabla_{x,t}S_{\tilde{A}^{na1}}(V^{a1}_{\zeta}[0])(0-t_{n}^{a1},
\cdot)\cdot
\nabla_{x,t}[\zeta^{na}[0]-S_{\tilde{A}^{na1}}(V^{a1}_{\zeta}[0])(0-t_{n}^{a1},
\cdot)]\, dx =o_{L^{2}}(1)
\]
Due to Lemma~\ref{BoxAEnergy}, up to $o_{L^{2}}(1)$, the left-hand
side equals
\begin{align*}
&\int_{0}^{1}\int_{\R^{2}}\nabla_{x,t}S_{A^{n}}\big(S_{\tilde{A}^{na1}}(V^{a1}_{\zeta}[0])[0-t_{n}^{a1}]\big)(t+t_{n}^{a1},
\cdot+x_{n}^{a1})\cdot\\
&\qquad\cdot \nabla_{x,t}S_{A^{n}}\big([\zeta^{na}[0]-S_{\tilde{A}^{na1}}(V^{a1}_{\zeta}[0])[0-t_{n}^{a1}]]\big)(t+t_{n}^{a1},
\cdot+x^{a1}_{n})\, dx
\end{align*}
But here we can again use that
\[
S_{A^{n}}\big(S_{\tilde{A}^{na1}}(V^{a1}_{\zeta}[0])[0-t_{n}^{a1}]\big)(t+t_{n}^{a1},
\cdot+x_{n}^{a1})= V^{a1}_{\zeta}(t,\cdot)+o_{L^{2}}(1)
\]
provided $t\in [0,1]$, while by construction
\[
S_{A^{n}}\big([\zeta^{na}[0]-S_{\tilde{A}^{na1}}(V^{a1}_{\zeta}[0])[0-t_{n}^{a1}]]\big)(t+t_{n}^{a1},
\cdot+x^{a1}_{n})\rightharpoonup 0
\]
The conclusion is that
\begin{equation}\nonumber\begin{split}
&\int_{0}^{1}\int_{\R^{2}}\nabla_{x,t}S_{A^{n}}\big(S_{\tilde{A}^{na1}}(V^{a1}_{\zeta}[0])[0-t_{n}^{a1}]\big)
(t+t_{n}^{a1}, \cdot+x_{n}^{a1})\cdot \\
&\qquad\qquad \cdot\nabla_{x,t}S_{A^{n}}\big([\zeta^{na}[0]-S_{\tilde{A}^{na1}}(V^{a1}_{\zeta}[0])
[0-t_{n}^{a1}]]\big)(t+t_{n}^{a1}, \cdot+x^{a1}_{n}) \,dx\\
&=o_{L^{2}}(1),\\
\end{split}\end{equation}
from which the asymptotic orthogonality follows.
All assertions of the lemma now follow by applying the preceding
considerations inductively $B$ times.
\end{proof}

\begin{figure*}[ht]
\begin{center}
\centerline{\hbox{\vbox{ \epsfxsize= 8.0 truecm \epsfysize=8.0
truecm \epsfbox{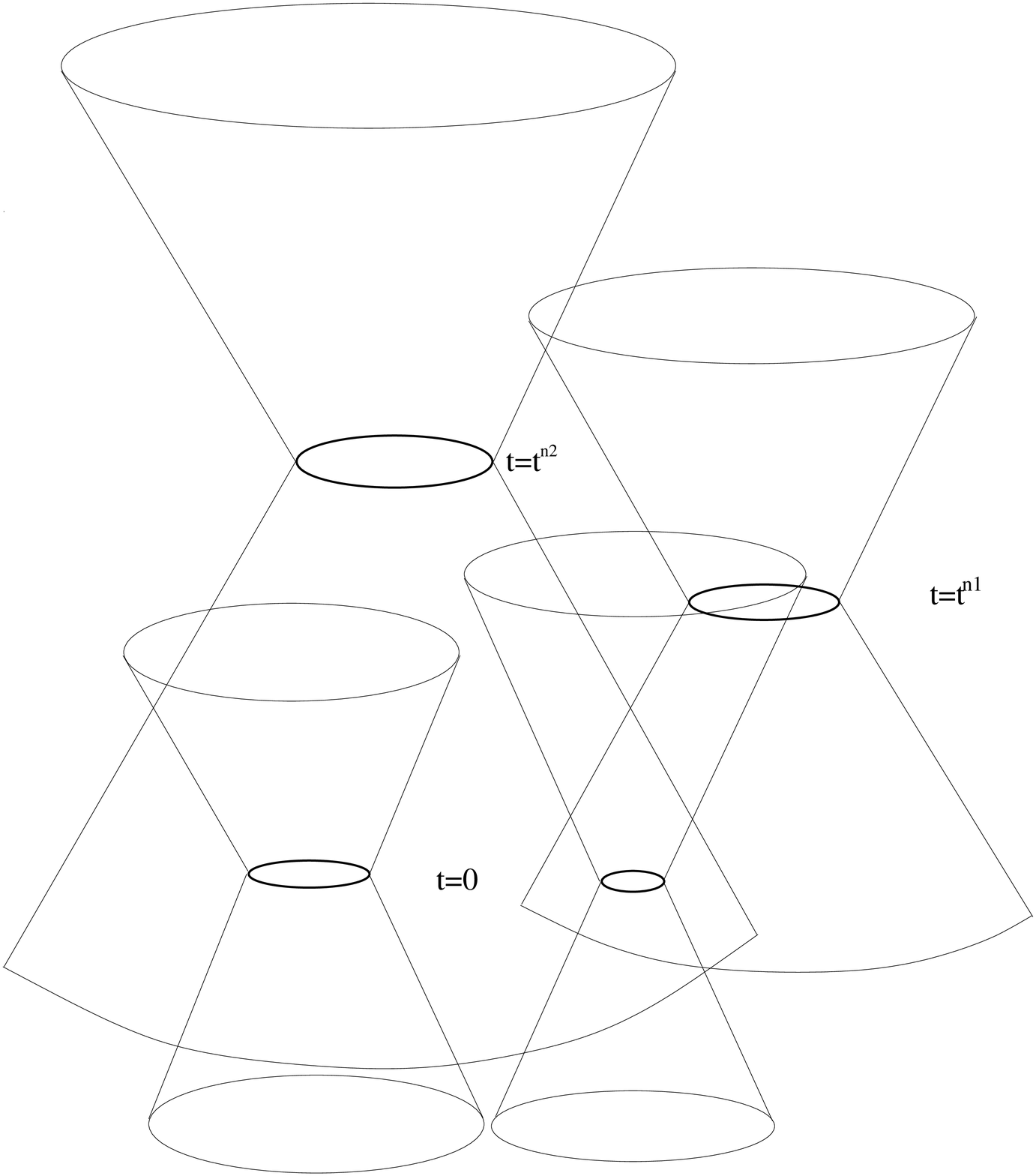}}}} \caption{The dependence domains of
various concentration profiles}
\end{center}
\end{figure*}

Figure~7 depicts various concentration profiles. More precisely, one should view these
profiles as being well-localized in physical space centered at their cores in space-time.
The figure then show the approximate support of the wave evolutions of these profiles.

Generally speaking, $a$ will always refer to an atom, whereas $b$
refers to the concentration profile generated by an atom. We shall
now apply Lemma~\ref{CoreBGII} to the covariant evolution of
$\zeta^{na}[0]$, as well as the remaining components $\eta^{na}$,
$\tilde{\psi}^{na}$, ${\mathring{\psi}}^{na}$.

\subsubsection{Selecting geometric concentration profiles}\label{subsec:geomprof}

At this stage, we face the same issue as in Step~1 above: we have a
sequence of component functions $V_{\alpha}^{ab}$ associated with
the essentially singular sequence $\phi^{n}_{\alpha}$, but in order
to apply the ``energy induction hypothesis'', i.e., the assumption
that $\Ecrit$ is the minimal energy for which uniform control fails,
we need to show that the $V_{\alpha}^{ab}$ can be assembled to form
the Coulomb components of actual maps from $\R^{2}\to \Hyp^2$. We
now address this task. To begin with, we may assume that
$\Psi^{nA_{0}^{(0)}}\ne o_{L^{2}}(1)$ since otherwise
$\Psi^{nA_{0}^{(0)}}$ is a perturbative error.

To summarize our construction of the concentration profiles: we
started with the Hodge decomposition (all at time $t=0$)
\begin{align*}
\phi^{na}_{1}&=\partial_{1}\tilde{\phi}^{na}+\partial_{2}{\mathring{\phi}}^{na},\,\partial_{t}\tilde{\phi}^{na}:=\phi^{na}_{0},\quad
\partial_{t}{\mathring{\phi}}^{na}=0
\\
\phi^{na}_{2}&=\partial_{2}\tilde{\phi}^{na}-\partial_{1}{\mathring{\phi}}^{na}
\\
\tilde{\psi}^{na}_{1}&=\partial_{1}\zeta^{na}+\partial_{2}\eta^{na},\quad
\partial_{t}\zeta^{na}:=\tilde{\psi}^{na}_{0},\quad \partial_{t}\eta^{na}=0
\\
\tilde{\psi}^{na}_{2}&=\partial_{2}\zeta^{na}-\partial_{1}\eta^{na}
\end{align*}
From here it is immediate that
\[
E(\phi^{na})=\sum_{\alpha=0}^{2}\|\phi^{na}_{\alpha}\|_{L_{x}^{2}}^{2}=\sum_{\alpha=0}^{2}
\|\partial_{\alpha}\tilde{\phi}^{na}\|_{L_{x}^{2}}^{2}+\sum_{\alpha=0}^{2}\|\partial_{\alpha}{\mathring{\phi}}^{na}\|_{L_{x}^{2}}^{2}
=\sum_{\alpha=0}^{2}\|\partial_{\alpha}\zeta^{na}\|_{L_{x}^{2}}^{2}+\sum_{\alpha=0}^{2}\|\partial_{\alpha}\eta^{na}\|_{L_{x}^{2}}^{2}
\]
Now, we evolve each of the $\tilde{\phi}^{na}$ etc.\ in time using
the covariant flow, and apply Lemma~\ref{CoreBGII}. Changing the
notation from that lemma, one obtains the decompositions (with
$\tilde A$ as in~\eqref{eq:Ashift})
\[
\nabla_{x,t}\tilde{\phi}^{na}=\sum_{b=1}^{B}\nabla_{x,t}\big[S_{\tilde{A}^{nab}}(\tilde{V}_{1}^{ab}[0])\big](0-t^{nab},
x-x^{nab})+\nabla_{x,t}\tilde{W}_{1}^{naB}
\]
\[
\nabla_{x,t}{\mathring{\phi}}^{na}=\sum_{b=1}^{B}\nabla_{x,t}\big[S_{\tilde{A}^{nab}}(\tilde{V}_{2}^{ab}[0])\big](0-t^{nab},
x-x^{nab})+\nabla_{x,t}\tilde{W}_{2}^{naB}
\]
\[
\nabla_{x,t}\zeta^{na}=\sum_{b=1}^{B}\nabla_{x,t}\big[S_{\tilde{A}^{nab}}(V_{1}^{ab}[0])\big](0-t^{nab},
x-x^{nab})+\nabla_{x,t}W_{1}^{naB}
\]
\[
\nabla_{x,t}\eta^{na}=\sum_{b=1}^{B}\nabla_{x,t}\big[S_{\tilde{A}^{nab}}(V_{2}^{ab}[0])\big](0-t^{nab},
x-x^{nab})+\nabla_{x,t}W_{2}^{naB}
\]
where the $W$ errors are small in the $\eta$-sense when $B$ is
large, see~\eqref{eq:Weta}. Here we use the same sequences $t^{nab},
x^{nab}$ for all decompositions, which of course we can by passing
to suitable subsequences. Note that we are working with {\em both}
the $\phi$ and $\psi$ components here, which is needed for the
following result.

\begin{prop}\label{ConcentrarionProfileApprox} For any $1\le b\le B$, and any $\delta_{2}>0$, there exists
an admissible (derivative components are Schwartz) map from $\R^{2}$ into $\Hyp^{2}$, with
derivative components $\phi^{nab}_{j\delta_{2}}$, $j=1,2$, and a number $\gamma_{\delta_{2}nab}\in \R$, such that
\begin{equation}\nonumber\begin{split}
&\Big\|\big(\partial_{1}\big[S_{\tilde{A}^{nab}}(V_{1}^{ab}[0])\big](0-t^{nab}, x)+\partial_{2}\big[S_{\tilde{A}^{nab}}(V_{2}^{ab}[0])\big](0-t^{nab}, x),\,\partial_{2}\big[S_{\tilde{A}^{nab}}(V_{1}^{ab}[0])\big](0-t^{nab}, x) \\
&-\partial_{1}\big[S_{\tilde{A}^{nab}}(V_{2}^{ab}[0])\big](0-t^{nab}, x)\big)-e^{i\gamma_{\delta_{2}nab}}\big(\phi_{1\delta_{2}}^{nab}e^{-i\sum_{k=1,2}\triangle^{-1}\partial_{k}\phi_{k\delta_{2}}^{nab}},
\phi_{2\delta_{2}}^{nab}e^{-i\sum_{k=1,2}\triangle^{-1}\partial_{k}\phi_{k\delta_{2}}^{nab}}\big)\Big\|_{L_{x}^{2}}<\delta_{2}
\end{split}\end{equation}
for large $n$.
\end{prop}
\begin{proof} Due to the asymptotic orthogonality relation of
Lemma~\ref{CoreBGII}, given $\delta_2$ there exists $B_0$ so that
for all $b>B_0$ one can simply take the derivative components to
equal zero. In other words, it suffices to consider $1\le b\le B_0$.

Fix a $b$, we shall pick $B$ larger, if necessary, and also pick $n$
large enough later. For simplicity introduce the notation
\[
V_{3}^{nab}:=\partial_{1}S_{\tilde{A}^{nab}}(V_{1}^{ab}[0])+\partial_{2}S_{\tilde{A}^{nab}}(V_{2}^{ab}[0])
\]
\[
V_{4}^{nab}:=\partial_{2}S_{\tilde{A}^{nab}}(V_{1}^{ab}[0])-\partial_{1}S_{\tilde{A}^{nab}}(V_{2}^{ab}[0])
\]
and similarly for $W_{3,4}^{nab}$. Note that we here introduce
dependence on~$n$ again.
\\
Thus at time $t=0$, we have the identities
\begin{align*}
\phi_{1}^{na}e^{-i\sum_{k=1,2}\triangle^{-1}\partial_{k}\phi^{1n}_{k}} &=\sum_{b=1}^{B}V_{3}^{nab}(0-t^{nab},
x-x^{nab})+W_{3}^{naB}\\
\phi_{2}^{na}e^{-i\sum_{k=1,2}\triangle^{-1}\partial_{k}\phi^{1n}_{k}}
&=\sum_{b=1}^{B}V_{4}^{nab}(0-t^{nab}, x-x^{nab})+W_{4}^{naB}
\end{align*}
where the $W$'s satisfy the smallness property~\eqref{eq:Weta}.
Then we distinguish between the following two cases:

\medskip\noindent
{\bf{(A)}}: {\em $V_{3,4}^{nab}(\cdot-t^{nab}, x-x^{nab})$ is of temporally
bounded type}. By this we mean that
\[
\liminf_{n\to\infty}|t^{nab}|<\infty
\]
By passing to a subsequence, we may then assume that
\[
\limsup_{n\to\infty}|t^{nab}|<\infty
\]
or in fact, that $\lim_{n\to\infty}t^{nab}$ exists.

\medskip\noindent
{\bf{(B)}}: {\em $V_{3,4}^{nab}(\cdot-t^{nab}, x-x^{nab})$ of temporally
unbounded type}. By this we mean that
\[
\lim_{n\to\infty}|t^{nab}|=\infty
\]
Observe that in this latter case, due to
Proposition~\ref{prop:covdisp}, we can conclude that
\[
V_{3,4}^{nab}(\cdot-t^{nab}, x-x^{nab})=o_{L^{\infty}}(1)+o_{L^{2}}(1)
\]
as $n\to\infty$.
\\

We treat these cases separately, commencing with the temporally
bounded Case~A. We need to show that $V_{3,4}^{nab}(\cdot-t^{nab},
x-x^{nab})$ can be approximated arbitrarily well by the Coulomb
components of admissible maps. We shall do this by physical
localization:  Note that for $b'\neq b$, we have either
\[
\lim_{n\to\infty}|t^{nab'}|=\infty
\]
or else
\[
\lim_{n\to\infty}|x^{nab}-x^{nab'}|=\infty
\]
We conclude that if $\chi^{nab}_{R}$ is a smooth spatial cutoff
localizing to a disc of radius $R$, $R<\infty$, centered at
$x^{nab}$, then  we have
\[
\lim_{n\to\infty}\|\chi_{R}^{nab}V_{3,4}^{nab'}(0-t^{nab'},
x-x^{nab'})\|_{L_{x}^{2}}=0,
\]
using Proposition~\ref{prop:covdisp}.  We also claim

\begin{lemma}\label{TailDispersion} Choosing $B=B(\delta_{2}, R)$ large enough, we get (here $\delta_{2}>0$ can be prescribed arbitrarily)
\[
\limsup_{n\to\infty}\|\chi_{R}^{nab}W_{3,4}^{naB}\|_{L_{x}^{2}}\ll
\delta_{2}
\]
for all $1\le b\le B_0$.
\end{lemma}
\begin{proof} Recall that
\[
W_{3,4}^{naB}=\partial_{1,2}W_{1}^{naB}\pm\partial_{2,1}W_{2}^{naB},
\]
where $W_{1,2}^{naB}$ solve the covariant wave equation
$\Box_{A^{n}}u=0$ (where as before $A^{n}$ is defined using
$\Psi^{nA_{0}^{(0)}}$). But then it is straightforward to check that
the space-time Fourier support of $u$ is contained in a small
neighborhood of the light cone intersected with the set $|\xi|\sim
1$, up to arbitrarily small errors. One can then reason exactly as
in~\cite{BG}, see Lemma~3.8 in that paper.
\end{proof}

Therefore, given $\delta_{2}>0$, we can pick $R=R(\delta_{2},
V_{1,2}^{ab})$ with the property that
\begin{equation}\nonumber\begin{split}
\limsup_{n\to\infty}&\Big\|\chi_{R}^{nab}\big(\phi_{1}^{na}e^{-i\sum_{k=1,2}\triangle^{-1}\partial_{k}[w_{k}^{nA_{0}^{(0)}}+\phi_{k}^{1na}]},
\phi_{2}^{na}e^{-i\sum_{k=1,2}\triangle^{-1}\partial_{k}[w_{k}^{nA_{0}^{(0)}}+\phi_{k}^{1na}]}\big)\\&
-\big(V_{3}^{ab}(0-t^{nab}, x-x^{nab}), V_{4}^{ab}(0-t^{nab},
x-x^{nab})\big)\Big\|_{L_{x}^{2}}\ll \delta_{2}
\end{split}\end{equation}
We now need to show that the components
\[
\big(\chi_{R}^{nab}\phi_{1}^{na}e^{-i\sum_{k=1,2}\triangle^{-1}\partial_{k}[w_{k}^{nA_{0}^{(0)}}+\phi_{k}^{1na}]},
\chi_{R}^{nab}\phi_{2}^{na}e^{-i\sum_{k=1,2}\triangle^{-1}\partial_{k}[w_{k}^{nA_{0}^{(0)}}+\phi_{k}^{1na}]}\big)
\]
are close to the Coulomb components of an admissible map, up to a
constant phase shift. To achieve this, we insert the profile
decomposition we obtained above for $\phi^{na}_{1,2}$, i.e., write
\[
\chi_{R}^{nab}\phi^{na}_{1}=\chi_{R}^{nab}[\sum_{b'=1}^{B}\sum_{j=1,2}\partial_{j}\big(S_{\tilde{A}^{nab'}}(\tilde{V}_{j}^{ab'}[0])(0-t^{nab'},
x-x^{nab'})+\tilde{W}_{j}^{naB}\big)]
\]
\[
\chi_{R}^{nab}\phi^{na}_{2}=\chi_{R}^{nab}[\sum_{b'=1}^{B}\sum_{j=1,2}(-1)^{j+1}\partial_{j+1}\big(S_{\tilde{A}^{nab'}}(\tilde{V}_{j}^{ab'}[0])(0-t^{nab'},
x-x^{nab'})+\tilde{W}_{j}^{naB}\big)]
\]
where $\del_{j+1}$ has to be interpreted modulo~$2$.  Now if we
choose $B$ large enough (depending on $R$, chosen further above),
and then choose $n$ large enough, we can ensure that
\[
\|\chi_{R}^{nab}\phi^{na}_{1}-\chi_{R}^{nab}\sum_{j=1,2}\partial_{j}\big(S_{\tilde{A}^{nab}}(\tilde{V}_{j}^{ab}[0])(0-t^{nab},
x-x^{nab})\|_{L_{x}^{2}}\ll \delta_{2}
\]
\[
\|\chi_{R}^{nab}\phi^{na}_{1}-\chi_{R}^{nab}\sum_{j=1,2}(-1)^{j+1}\partial_{j+1}\big(S_{\tilde{A}^{nab}}(\tilde{V}_{j}^{ab}[0])(0-t^{nab},
x-x^{nab})\|_{L_{x}^{2}}\ll \delta_{2}
\]
We continue by approximating the truncated components
$\chi_{R}^{nab}\phi^{na}_{1}$ by the derivative components of an
admissible map $(\tilde{\bf{x}}^{nab}, \tilde{\bf{y}}^{nab}):
\R^{2}\to\Hyp^2$.
\\
For this purpose we recall the identity
\[
\phi^{1na}_{j}=({\bf{y}}^{na})^{-1}\sum_{k=1,2}\triangle^{-1}\partial_{k}\partial_{j}[\phi^{1na}_{k}{\bf{y}}^{na}],\quad
j=1,2
\]
Inserting the above decomposition for the $\phi^{1na}_{k}$, we
obtain
\begin{equation}\nonumber\begin{split}
\phi^{1na}_{j}=&({\bf{y}}^{na})^{-1}\sum_{k=1,2}\triangle^{-1}\partial_{1}\partial_{j}[[\sum_{b'=1}^{B}\sum_{\tilde{j}=1,2}\partial_{\tilde{j}}\big(\Re S_{\tilde{A}^{nab'}}(\tilde{V}_{\tilde{j}}^{ab'}[0])(0-t^{nab'}, x-x^{nab'})+\Re\tilde{W}_{\tilde{j}}^{naB}\big)]{\bf{y}}^{na}]\\
&+({\bf{y}}^{na})^{-1}\sum_{k=1,2}\triangle^{-1}\partial_{2}\partial_{j}[[\sum_{b'=1}^{B}\sum_{\tilde{j}=1,2}(-1)^{\tilde{j}+1}\partial_{\tilde{j}+1}\big(\Re
S_{\tilde{A}^{nab'}}(\tilde{V}_{\tilde{j}}^{ab'}[0])(0-t^{nab'},
x-x^{nab'})+\Re\tilde{W}_{\tilde{j}}^{naB}\big)]{\bf{y}}^{na}]
\end{split}\end{equation}
Using the frequency localization of all functions involved, and
increasing $R$ if necessary (independently of $B$), we can then
achieve that for $n$ large enough
\begin{equation}\nonumber\begin{split}
\big\|&\phi^{1na}_{j}-({\bf{y}}^{na})^{-1}\sum_{k=1,2}\sum_b\triangle^{-1}\partial_{1}\partial_{j}
[\chi_{R}^{nab}[\sum_{\tilde{j}=1,2}\partial_{\tilde{j}}\big(\Re S_{\tilde{A}^{nab}}(\tilde{V}_{\tilde{j}}^{ab}[0])
(0-t^{nab}, x-x^{nab})+\Re\tilde{W}_{\tilde{j}}^{naB}\big)]{\bf{y}}^{na}]\\
&+({\bf{y}}^{na})^{-1}\sum_{k=1,2}\sum_b\triangle^{-1}\partial_{2}\partial_{j}[\chi_{R}^{nab}
[\sum_{\tilde{j}=1,2}(-1)^{\tilde{j}+1}\partial_{\tilde{j}+1}\big(\Re
S_{\tilde{A}^{nab}}(\tilde{V}_{\tilde{j}}^{ab}[0])(0-t^{nab},
x-x^{nab})+\Re\tilde{W}_{\tilde{j}}^{naB}\big)]{\bf{y}}^{na}]\big\|_{L_{x}^{2}}\ll
\delta_{2}
\end{split}\end{equation}
From  here we infer that
\[
\limsup_{n\to\infty}\big\|\chi_{R}^{nab}\phi^{1na}_{j}-({\bf{y}}^{na})^{-1}
\sum_{k=1,2}\triangle^{-1}\partial_{k}\partial_{j}[\chi_{R}^{nab}\phi^{1na}_{k}{\bf{y}}^{na}]\big\|_{L_{x}^{2}}\ll
\delta_{2}
\]
Now modify ${\bf{y}}$ to a function $\tilde{{\bf{y}}}^{nab}$ by
picking numbers $R'', R'$ with
\[
R\ll R''\ll R',
\]
to be specified shortly, and setting
\[
\tilde{{\bf{y}}}^{nab}=e^{\sum_{k=1,2}\chi_{R'}^{nab}\triangle^{-1}\partial_{k}P_{[-R'',
R'']}\phi^{2na}_{k}}
\]
whence
\[
\frac{\partial_{j}\tilde{{\bf{y}}}^{nab}}{\tilde{{\bf{y}}}^{nab}}=\partial_{j}\sum_{k=1,2}\chi_{R'}^{nab}\triangle^{-1}\partial_{k}P_{[-R'',
R'']}\phi^{2na}_{k}=\chi_{R'}^{nab}\phi^{2na}_{j}+\text{error},
\]
where we can achieve that $\|\text{error}\|_{L_{x}^{2}}\ll
\delta_{2}$ by choosing $R''$ large enough depending on $\delta_{2}$
and the localization of $\phi^{na}$ in frequency space, and then
$R'$ large enough in relation to $R''$. Increasing $B$ if necessary
and then choosing $n$ large enough, we can then also achieve that
\[
\Big\|\chi_{R'}^{nab}\phi^{na}_{2}-\chi_{R'}^{nab}\sum_{j=1,2}(-1)^{j+1}\partial_{j+1}S_{\tilde{A}^{nab}}(\tilde{V}_{\tilde{j}}^{ab}[0])(0-t^{nab},
x-x^{nab})\Big\|_{L_{x}^{2}}\ll \delta_{2}
\]
and then
\[
\|\chi_{R'}^{nab}\phi^{na}_{2}-\chi_{R}^{nab}\phi^{na}_{2}\|_{L_{x}^{2}}\ll
\delta_{2}
\]
We next show that the expression
\[
({\bf{y}}^{na})^{-1}\sum_{k=1,2}\triangle^{-1}\partial_{k}\partial_{j}[\chi_{R}^{nab}\phi^{1na}_{k}{\bf{y}}^{na}]
\]
is well approximated by
\[
(\tilde{{\bf{y}}}^{nab})^{-1}\sum_{k=1,2}\triangle^{-1}\partial_{k}\partial_{j}[\chi_{R}^{nab}\phi^{1na}_{k}\tilde{{\bf{y}}}^{nab}]
\]
To see this, write
\[
({\bf{y}}^{na})^{-1}\sum_{k=1,2}\triangle^{-1}\partial_{k}\partial_{j}
[\chi_{R}^{nab}\phi^{1na}_{k}{\bf{y}}^{na}]=e^{-\sum_{k=1,2}\chi_{R'}^{nab}\triangle^{-1}\partial_{k}\phi^{2na}_{k}}
\sum_{k=1,2}\triangle^{-1}\partial_{k}\partial_{j}[\chi_{R}^{nab}\phi^{1na}_{k}e^{\sum_{k=1,2}\chi_{R'}^{nab}\triangle^{-1}
\partial_{k}\phi^{2na}_{k}}]+\text{error},
\]
where we can achieve $\|\text{error}\|_{L_{x}^{2}}\ll \delta_{2}$ by
choosing $R'$ large enough in relation to $R$ and the intrinsic
Fourier localization properties of $\phi^{na}$.
\\

Split the phase into the product
\[
e^{\sum_{k=1,2}\chi_{R'}^{nab}\del_k^{-1}\phi^{2na}_{k}}=e^{\sum_{k=1,2}\chi_{R'}^{nab}\del_k^{-1}P_{[-R'',
R'']}\phi^{2na}_{k}}e^{\sum_{k=1,2}\chi_{R'}^{nab}\del_k^{-1}P_{[-R'',
R'']^{c}}\phi^{2na}_{k}}
\]
We need to show that we can eliminate the factor
$e^{\sum_{k=1,2}\chi_{R'}^{nab}\del_k^{-1}P_{[-R'',
R'']^{c}}\phi^{2na}_{k}}$. Using similar arguments as in Step~1,
choosing $R''$ large enough in relation to $R$, it is
straightforward to show that, with $\del_k^{-1}:=\Delta^{-1}\del_k$,
\begin{equation}\nonumber\begin{split}
&\|e^{-\sum_{k=1,2}\chi_{R'}^{nab}\del_k^{-1}\phi^{2na}_{k}}\sum_{k=1,2}\del_k^{-1}\partial_{j}
[\chi_{R}^{nab}\phi^{1na}_{k}e^{\sum_{k=1,2}\chi_{R'}^{nab}P_{(-\infty,
R'']}\del_k^{-1}\phi^{2na}_{k}}(e^{\sum_{k=1,2}\chi_{R'}^{nab}\del_k^{-1}P_{>R''}\phi^{2na}_{k}}-1)]\|_{L_{x}^{2}}\\&\ll
\delta_{2}
\end{split}\end{equation}
\begin{equation}\nonumber\begin{split}
&\|(e^{-\sum_{k=1,2}\chi_{R'}^{nab}\del_k^{-1}P_{>R''}\phi^{2na}_{k}}-1)e^{-\sum_{k=1,2}\chi_{R'}^{nab}\del_k^{-1}P_{(-\infty,
R'']}\phi^{2na}_{k}}\sum_{k=1,2}\del_k^{-1}\partial_{j}[\chi_{R}^{nab}\phi^{1na}_{k}e^{\sum_{k=1,2}\chi_{R'}^{nab}\triangle^{-1}P_{(-\infty,
R'']}\partial_{k}\phi^{2na}_{k}}]\|_{L_{x}^{2}}\\&\ll \delta_{2}
\end{split}\end{equation}
We next show that we can also eliminate
$e^{\sum_{k=1,2}\chi_{R'}^{nab}\del_k^{-1}P_{<-R''}\phi^{2na}_{k}}$.
Indeed, proceeding as in the first section, write
\begin{equation}\nonumber\begin{split}
&e^{-\sum_{k=1,2}\chi_{R'}^{nab}\del_k^{-1}P_{<-R''}\phi^{2na}_{k}}(\tilde{{\bf{y}}}^{nab})^{-1}
\sum_{k=1,2}\del_k^{-1}\partial_{j}[\chi_{R}^{nab}\phi^{1na}_{k}\tilde{{\bf{y}}}^{nab}e^{\sum_{k=1,2}\chi_{R'}^{nab}\del_k^{-1}P_{<-R''}\phi^{2na}_{k}}]\\
&=\sum_{l\geq 2}[\chi_{lR}^{nab}-\chi_{(l-1)R}^{nab}]e^{-\sum_{k=1,2}\chi_{R'}^{nab}\del_k^{-1}P_{<-R''}\phi^{2na}_{k}}
(\tilde{{\bf{y}}}^{nab})^{-1}\sum_{k=1,2}\del_k^{-1}\partial_{j}[\chi_{R}^{nab}\phi^{1na}_{k}\tilde{{\bf{y}}}^{nab}
e^{\sum_{k=1,2}\chi_{R'}^{nab}\del_k^{-1}P_{<-R''}\phi^{2na}_{k}}]\\
&+\chi_{R}^{nab}e^{-\sum_{k=1,2}\chi_{R'}^{nab}\del_k^{-1}P_{<-R''}\phi^{2na}_{k}}(\tilde{{\bf{y}}}^{nab})^{-1}
\sum_{k=1,2}\del_k^{-1}\partial_{j}[\chi_{R}^{nab}\phi^{1na}_{k}\tilde{{\bf{y}}}^{nab}e^{\sum_{k=1,2}\chi_{R'}^{nab}\del_k^{-1}P_{<-R''}\phi^{2na}_{k}}]\\
\end{split}\end{equation}
Here the cutoff $\chi_{lR}^{nab}$ localizes to a disc of radius $lR$
around $x^{nab}$. Then pick a point $p_{lnab}$ in this disc, for
each $l$, and write for fixed $l\geq 2$
\begin{equation}\nonumber\begin{split}
&e^{-\sum_{k=1,2}\chi_{R'}^{nab}\del_k^{-1}P_{<-R''}\phi^{2na}_{k}}(\tilde{{\bf{y}}}^{nab})^{-1}\sum_{k=1,2}\del_k^{-1}\partial_{j}[\chi_{R}^{nab}\phi^{1na}_{k}\tilde{{\bf{y}}}^{nab}e^{\sum_{k=1,2}\chi_{R'}^{nab}\del_k^{-1}P_{<-R''}\phi^{2na}_{k}}]\\
&=\big(\frac{e^{\sum_{k=1,2}\chi_{R'}^{nab}\del_k^{-1}P_{<-R''}\phi^{2na}_{k}}}{e^{\sum_{k=1,2}\chi_{R'}^{nab}\del_k^{-1}P_{<-R''}\phi^{2na}_{k}(p_{lnab})}}\big)^{-1}(\tilde{{\bf{y}}}^{nab})^{-1}\sum_{k=1,2}\del_k^{-1}\partial_{j}[\chi_{R}^{nab}\phi^{1na}_{k}\tilde{{\bf{y}}}^{nab}\frac{e^{\sum_{k=1,2}\chi_{R'}^{nab}\del_k^{-1}P_{<-R''}\phi^{2na}_{k}}}{e^{\sum_{k=1,2}\chi_{R'}^{nab}\del_k^{-1}P_{<-R''}\phi^{2na}_{k}(p_{lnab})}}]\\
\end{split}\end{equation}
But then we can estimate
\[
[\chi_{lR}^{nab}-\chi_{(l-1)R}^{nab}]\big(\frac{e^{\sum_{k=1,2}\chi_{R'}^{nab}\del_k^{-1}P_{<-R''}\phi^{2na}_{k}}}{e^{\sum_{k=1,2}\chi_{R'}^{nab}\del_k^{-1}P_{<-R''}\phi^{2na}_{k}(p_{lnab})}}-1\big)=O(\frac{lR}{R''}),
\]
and then using the machinery from  Step~1  (which yields a $l^{-N}$
gain), and choosing $R''\gg R$, we can achieve that
\begin{equation}\nonumber\begin{split}
&\Big\|[\big(\frac{e^{\sum_{k=1,2}\chi_{R'}^{nab}\del_k^{-1}P_{<-R''}\phi^{2na}_{k}}}{e^{\sum_{k=1,2}\chi_{R'}^{nab}\del_k^{-1}P_{<-R''}\phi^{2na}_{k}(p_{lnab})}}\big)^{-1}-1](\tilde{{\bf{y}}}^{nab})^{-1}\sum_{k=1,2}\del_k^{-1}\partial_{j}[\chi_{R}^{nab}\phi^{1na}_{k}\tilde{{\bf{y}}}^{nab}\frac{e^{\sum_{k=1,2}\chi_{R'}^{nab}\del_k^{-1}P_{<-R''}\phi^{2na}_{k}}}{e^{\sum_{k=1,2}\chi_{R'}^{nab}\del_k^{-1}P_{<-R''}\phi^{2na}_{k}(p_{lnab})}}]\Big\|_{L_{x}^{2}}\\
&\ll \delta_{2}l^{-N}
\end{split}\end{equation}
Similarly, one can eliminate the second instance of
\[
\frac{e^{\sum_{k=1,2}\chi_{R'}^{nab}\del_k^{-1}P_{<-R''}\phi^{2na}_{k}}}{e^{\sum_{k=1,2}\chi_{R'}^{nab}\del_k^{-1}P_{<-R''}\phi^{2na}_{k}(p_{lnab})}}
\]
Hence we have now shown that for $B$ and then $n$ large enough, we
have that the functions
\[
\tilde{\phi}^{1nab}_{j}:=(\tilde{{\bf{y}}}^{nab})^{-1}\sum_{k=1,2}\del_k^{-1}\partial_{j}[\chi_{R}^{nab}\phi^{1na}_{k}\tilde{{\bf{y}}}^{nab}],\,\tilde{\phi}^{2nab}_{j}:=\frac{\partial_{j}\tilde{{\bf{y}}}^{nab}}{\tilde{{\bf{y}}}^{nab}},
\]
which of course are the derivative components of admissible maps,
satisfy the inequalities
\[
\|\tilde{\phi}^{1nab}_{j}-\chi_{R}^{nab}\phi^{1na}_{j}\|_{L_{x}^{2}}\ll
\delta_{2},\quad
\|\tilde{\phi}^{2nab}_{j}-\chi_{R}^{nab}\phi^{2na}_{j}\|_{L_{x}^{2}}\ll
\delta_{2}
\]
Our next task is to approximate the {\em{Coulomb components}}. For
this consider
\[
\tilde{\phi}^{nab}_{j}e^{-i\sum_{k=1,2}\del_k^{-1}[w_{k}^{nA_{0}^{(0)}}+\phi^{1na}_{k}]},\quad\tilde{\phi}^{nab}_{j}
=\tilde{\phi}^{1nab}_{j}+i\tilde{\phi}^{2nab}_{j}
\]
From the preceding, we can arrange that
\[
\|\tilde{\phi}^{nab}_{j}e^{-i\sum_{k=1,2}\del_k^{-1}\phi^{1n}_{k}}-\chi_{R}^{nab}\phi^{na}_{j}
e^{-i\sum_{k=1,2}\del_k^{-1}[w_{k}^{nA_{0}^{(0)}}+\phi^{1na}_{k}]}\|_{L_{x}^{2}}\ll
\delta_{2}
\]
We need to show that we can also arrange (i.e., upon choosing $B$,
$n$ large enough) that
\[
\|\tilde{\phi}^{nab}_{j}e^{-i\sum_{k=1,2}\del_k^{-1}[w_{k}^{nA_{0}^{(0)}}+\phi^{1na}_{k}]}-\tilde{\phi}^{nab}_{j}e^{-i\sum_{k=1,2}\del_k^{-1}\tilde{\phi}^{1nab}_{k}}e^{i\gamma_{nab}}\|_{L_{x}^{2}}\ll
\delta_{2}
\]
for a suitable constant $\gamma_{nab}$.
\\
We first get rid of the phase
$e^{-i\sum_{k=1,2}\del_k^{-1}w_{k}^{nA_{0}^{(0)}}}$:
simply pick a point $p_{nab}$ in the support of $\chi_{R}^{nab}$,
and replace
$e^{-i\sum_{k=1,2}\del_k^{-1}w_{k}^{nA_{0}^{(0)}}}$
by
$e^{-i\sum_{k=1,2}\del_k^{-1}w_{k}^{nA_{0}^{(0)}}(p_{nab})}$.

\noindent
Next, we need to show that
$\tilde{\phi}^{nab}_{j}e^{-i\sum_{k=1,2}\del_k^{-1}\phi^{1na}_{k}}$
is close to
$\tilde{\phi}^{nab}_{j}e^{-i\sum_{k=1,2}\del_k^{-1}\tilde{\phi}^{1nab}_{k}}$,
up to a constant phase shift. First, pick $R_{1}\gg R$ such that
\[
\|\tilde{\phi}^{nab}_{j}e^{-i\sum_{k=1,2}\del_k^{-1}\phi^{1na}_{k}}-\chi_{R_{1}}^{nab}\tilde{\phi}^{nab}_{j}e^{-i\sum_{k=1,2}\del_k^{-1}\phi^{1na}_{k}}\|_{L_{x}^{2}}\ll
\delta_{2}
\]
Next, pick $R'\gg R_{1}$ such that
\[
\|\chi_{R_{1}}^{nab}\tilde{\phi}^{nab}_{j}e^{-i\sum_{k=1,2}\del_k^{-1}\phi^{1na}_{k}}-e^{i\gamma_{1nab}}\chi_{R_{1}}^{nab}\tilde{\phi}^{nab}_{j}e^{-i\sum_{k=1,2}\del_k^{-1}P_{[-R',
R']}\phi^{1na}_{k}}\|_{L_{x}^{2}}\ll \delta_{2}
\]
for suitable $\gamma_{1nab}$. Next, we claim that picking
$R_{2}\gg R'$, we can arrange that
\[
\|\chi_{R_{1}}^{nab}\tilde{\phi}^{nab}_{j}e^{-i\sum_{k=1,2}\del_k^{-1}P_{[-R',
R']}\phi^{1na}_{k}}-\chi_{R_{1}}^{nab}\tilde{\phi}^{nab}_{j}e^{-i\sum_{k=1,2}\del_k^{-1}P_{[-R',
R']}[\chi_{R_{2}}^{nab}\phi^{1na}_{k}]}\|_{L_{x}^{2}}\ll \delta_{2}
\]
This is a consequence of the fact that
\[
\|\chi_{R_{1}}^{nab}e^{-i\sum_{k=1,2}\del_k^{-1}P_{[-R',
R']}[(1-\chi_{R_{2}}^{nab})\phi^{1na}_{k}]}\|_{L_{x}^{\infty}}\ll
\delta_{2}
\]
Finally, we claim that
\[
\chi_{R_{1}}^{nab}\tilde{\phi}^{nab}_{j}e^{-i\sum_{k=1,2}\del_k^{-1}P_{[-R',
R']}[\chi_{R_{2}}^{nab}\phi^{1na}_{k}]}
\]
is very close to
$\tilde{\phi}^{nab}_{j}e^{-i\sum_{k=1,2}\del_k^{-1}\tilde{\phi}^{nab}_{k}}$,
which is what we need to finish case {\bf{(A)}}. To see this, note
that by choosing $R_{2}$ large enough in relation to $R'$, we get
from Bernstein's inequality
\[
\|\sum_{k=1,2}\del_k^{-1}P_{[-R',
R']}[\chi_{R_{2}}^{nab}\phi^{1na}_{k}]-\sum_{k=1,2}\del_k^{-1}P_{[-R',
R']}\phi^{1na}_{k}\|_{L_{x}^{\infty}}\ll \delta_{2}
\]
This immediately implies
\[
\|\chi_{R_{1}}^{nab}\tilde{\phi}^{nab}_{j}e^{-i\sum_{k=1,2}\del_k^{-1}P_{[-R',
R']}[\chi_{R_{2}}^{nab}\phi^{1na}_{k}]}-\chi_{R_{1}}^{nab}\tilde{\phi}^{nab}_{j}e^{-i\sum_{k=1,2}\del_k^{-1}P_{[-R',
R']}\phi^{1na}_{k}}\|_{L_{x}^{2}}\ll \delta_{2}
\]
To conclude, picking $R'$ large enough in relation to $R_{1}$ allows
us to find a phase $e^{i\gamma_{2nab}}$ such that
\[
\|\chi_{R_{1}}^{nab}\tilde{\phi}^{nab}_{j}e^{-i\sum_{k=1,2}\del_k^{-1}P_{[-R',
R']}\phi^{1na}_{k}}-\chi_{R_{1}}^{nab}\tilde{\phi}^{nab}_{j}e^{-i\sum_{k=1,2}\del_k^{-1}\phi^{1na}_{k}}e^{i\gamma_{2nab}}\|_{L_{x}^{2}}\ll
\delta_{2}
\]
Since we also have, as mentioned before, that
\[
\|\chi_{R_{1}}^{nab}\tilde{\phi}^{nab}_{j}-\tilde{\phi}^{nab}_{j}\|_{L_{x}^{2}}\ll
\delta_{2}
\]
Combining all of the preceding steps, we infer the existence of a
phase $e^{i\gamma_{nab}}$ such that
\[
\|\chi_{R}^{nab}\phi^{na}_{j}e^{-i\sum_{k=1,2}\del_k^{-1}\phi^{1na}_{k}}-
\tilde{\phi}^{nab}_{j}e^{-i\sum_{k=1,2}\del_k^{-1}\phi^{1na}_{k}}e^{i\gamma_{nab}}\|_{L_{x}^{2}}\ll
\delta_{2}
\]
We then get for suitable $\gamma'_{nab}$
\[
\|e^{i\gamma'_{nab}}\tilde{\phi}^{nab}_{1,2}e^{-i\sum_{k=1,2}\del_k^{-1}\phi^{1na}_{k}}e^{i\gamma_{nab}}-V_{3,4}^{ab}(0-t^{nab},
x-x^{nab})\|_{L_{x}^{2}}\ll \delta_{2}
\]
This finally concludes case {\bf{(A)}}, i.e., the temporally bounded
case.
\\

{\bf{(B)}}: {\em temporally unbounded case.} Here we have
$\lim_{n\to\infty}|t^{nab}|=\infty$, whence using
Proposition~\ref{prop:covdisp},  we get that
\[
V_{3,4}^{nab}(0-t^{nab}, x-x^{nab})=o_{L^{\infty}}(1)+o_{L^{2}}(1),
\]
where we recall the notation
\[
V_{3}^{nab}=\partial_{1}S_{\tilde{A}^{nab}}(V_{1}^{ab}[0])+\partial_{2}S_{\tilde{A}^{nab}}(V_{2}^{ab}[0])
\]
\[
V_{4}^{nab}:=\partial_{2}S_{\tilde{A}^{nab}}(V_{1}^{ab}[0])-\partial_{1}S_{\tilde{A}^{nab}}(V_{2}^{ab}[0])
\]
We now make the following
\\

{\bf{Claim}}: {\it{Choosing $n$ large enough, we have
\[
\|\nabla_{x,t}S_{\tilde{A}^{nab}}(V_{2}^{ab}[0])\|_{L_{x}^{2}}\ll
\delta_{3}
\]
for any given $\delta_{3}>0$.}}
\\

Thus in Case (B), the components $V_{3,4}^{nab}(0-t^{nab},
x-x^{nab})$ are approximately given by the gradient of a suitable
(complex valued) function. Once the Claim is established, Case (B)
will be straightforward to conclude.
\\

In order to prove the Claim, we shall use the curl equations
satisfied by the components $\phi^{na}_{j}$. To begin with, pick $R$
large enough such that
\[
\|P_{[-R,
R]}(\phi^{na}_{j})e^{-i\sum_{k=1,2}\del_k^{-1}\phi^{na}_{k}}-\sum_{b=1}^{B}(V_{j+2}(0-t^{nab},
x-x^{nab})+W_{j+2}^{naB})\|_{L_{x}^{2}}\ll \delta_{2},\quad j=1,2
\]
Then using the Littlewood-Paley trichotomy, and choosing $R$ larger
if necessary, we can arrange that
\[
\|P_{[-10R, 10R]}\big(P_{[-R,
R]}(\phi^{na}_{j})e^{-i\sum_{k=1,2}\del_k^{-1}\phi^{na}_{k}}\big)-\sum_{b=1}^{B}(V_{j+2}(0-t^{nab},
x-x^{nab})+W_{j+2}^{naB})\|_{L_{x}^{2}}\ll \delta_{2},\quad j=1,2
\]
Now fix a cutoff $\chi^{nab}$ which localizes to a large annulus of
radius $|t^{nab}|$ around $x^{nab}$ and thickness $R_{n}$ large
enough, such that
\[
\limsup_{n\to\infty}\|\chi^{nab}V_{3,4}^{nab}(0-t^{nab},
x-x^{nab})-V_{3,4}^{nab}(0-t^{nab}, x-x^{nab})\|_{L_{x}^{2}}\ll
\delta_{2}
\]
By removing finitely many 'holes' from this annulus and adjusting
$\chi^{nab}$ correspondingly, we can ensure that
\[
\lim_{n\to\infty}\|\chi^{nab}V_{3,4ab'}(0-t^{nab'},
x-x^{nab'})\|_{L_{x}^{2}}=0,\quad b\neq b',\; 1\leq b'\leq B
\]
We cannot simply arrange that
\[
\lim_{n\to\infty}\|\chi^{nab}W_{3,4}^{naB}\|_{L_{x}^{2}}=0
\]
and it will be more complicated to disentangle $W_{3,4}^{naB}$,
$V_{3,4}^{nab}$. From the preceding, choosing $R$ and then $n$ large
enough, we can arrange that
\[
\|\chi^{nab}P_{[-10R, 10R]}\big(P_{[-R,
R]}(\phi^{na}_{j})e^{-i\sum_{k=1,2}\del_k^{-1}\phi^{na}_{k}}\big)-\chi^{nab}[V_{j+2}(0-t^{nab},
x-x^{nab})+W_{j+2}^{naB}]\|_{L_{x}^{2}}\ll \delta_{2},\quad j=1,2
\]
Here $R$ only depends on the frequency concentration of $\phi^{na}$.
We now analyze the curl expression
\begin{equation}\nonumber\begin{split}
&\nabla^{-1}\partial_{1}\big[\chi^{nab}P_{[-10R, 10R]}\big(P_{[-R, R]}(\phi^{na}_{2})e^{-i\sum_{k=1,2}\del_k^{-1}\phi^{na}_{k}}\big)\big]\\
&-\nabla^{-1}\partial_{2}\big[\chi^{nab}P_{[-10R, 10R]}\big(P_{[-R, R]}(\phi^{na}_{1})e^{-i\sum_{k=1,2}\del_k^{-1}\phi^{na}_{k}}\big)\big]\\
\end{split}\end{equation}
We shall show that this expression becomes arbitrarily small when
$n$ is sufficiently large. Decompose the above expression into
\begin{equation}\nonumber\begin{split}
&\nabla^{-1}\big[\partial_{1}\chi^{nab}P_{[-10R, 10R]}\big(P_{[-R, R]}(\phi^{na}_{2})e^{-i\sum_{k=1,2}\del_k^{-1}\phi^{na}_{k}}\big)\big]\\
&-\nabla^{-1}\big[\partial_{2}\chi^{nab}P_{[-10R, 10R]}\big(P_{[-R, R]}(\phi^{na}_{1})e^{-i\sum_{k=1,2}\del_k^{-1}\phi^{na}_{k}}\big)\big]\\
&+\nabla^{-1}\big[\chi^{nab}P_{[-10R, 10R]}\big(P_{[-R, R]}(\partial_{1}\phi^{na}_{2}-\partial_{2}\phi^{na}_{1})e^{-i\sum_{k=1,2}\del_k^{-1}\phi^{na}_{k}}\big)\big]\\
&+\nabla^{-1}\big[\chi^{nab}P_{[-10R, 10R]}\big(P_{[-R, R]}(\phi^{na}_{2})\partial_{1}(e^{-i\sum_{k=1,2}\del_k^{-1}\phi^{na}_{k}})\big)\big]\\
&-\nabla^{-1}\big[\chi^{nab}P_{[-10R, 10R]}\big(P_{[-R, R]}(\phi^{na}_{1})\partial_{2}(e^{-i\sum_{k=1,2}\del_k^{-1}\phi^{na}_{k}})\big)\big]\\
\end{split}\end{equation}
For the first two terms, choosing the cutoff $\chi^{nab}$ suitably,
it is clear that for $n$ large enough we have
\begin{equation}\nonumber\begin{split}
&\|\nabla^{-1}\big[\partial_{1}\chi^{nab}P_{[-10R, 10R]}\big(P_{[-R, R]}(\phi^{na}_{2})e^{-i\sum_{k=1,2}\del_k^{-1}\phi^{na}_{k}}\big)\big]\\
&-\nabla^{-1}\big[\partial_{2}\chi^{nab}P_{[-10R, 10R]}\big(P_{[-R, R]}(\phi^{na}_{1})e^{-i\sum_{k=1,2}\del_k^{-1}\phi^{na}_{k}}\big)\big]\|_{L_{x}^{2}}\\
&\ll \delta_{2}
\end{split}\end{equation}
For the third term, we use the schematic curl relation
$\partial_{1}\phi_{2}^{na}-\partial_{2}\phi_{1}^{na}="(\phi^{na})^{2}"$.
Note that by including a suitable cutoff $\tilde{\chi}^{nab}$ having
similar characteristics as $\chi^{nab}$, we get
\begin{equation}\nonumber\begin{split}
&\|\nabla^{-1}\big[\chi^{nab}P_{[-10R, 10R]}\big(P_{[-R, R]}(\partial_{1}\phi^{na}_{2}-\partial_{2}\phi^{na}_{1})e^{-i\sum_{k=1,2}\del_k^{-1}\phi^{na}_{k}}\big)\big]\\
&-\nabla^{-1}\big[\chi^{nab}P_{[-10R, 10R]}\big(P_{[-R, R]}(\tilde{\chi}^{nab}"(\phi^{na})^{2}")e^{-i\sum_{k=1,2}\del_k^{-1}\phi^{na}_{k}}\big)\big]\|_{L_{x}^{2}}\ll \delta_{2}\\
\end{split}\end{equation}
Now we insert the decomposition
\[
\phi^{na}_{j}=\sum_{b'=1}^{B}\tilde{V}_{j+2}^{nab}(0-t^{nab},
x-x^{nab})+\tilde{W}_{j+2}^{naB},\quad j=1,2
\]
For any chosen $B$, by picking $n$ large enough, we can achieve that
\[
\Big\|\tilde{\chi}^{nab}\sum_{b'=1,\quad b'\neq
b}^{B}\tilde{V}_{j+2}^{nab}(0-t^{nab},
x-x^{nab})\Big\|_{L_{x}^{2}}\ll \delta_{2},
\]
and hence we reduce to estimating
\[
\|\nabla^{-1}\big[\chi^{nab}P_{[-10R, 10R]}\big(P_{[-R, R]}(\tilde{\chi}^{nab}[\tilde{V}_{3,4}^{nab}+\tilde{W}_{3,4}^{naB}]^{2})e^{-i\sum_{k=1,2}\del_k^{-1}\phi^{na}_{k}}\big)\big]\|_{L_{x}^{2}}\\
\]
Now recall from Proposition~\ref{prop:covdisp} that
\[
\tilde{\chi}^{nab}\tilde{V}_{3,4}^{nab}=o_{L^{\infty}}(1)+o_{L^{2}}(1)
\]
Hence we obtain
\begin{equation}\nonumber\begin{split}
&\|\nabla^{-1}\big[\chi^{nab}P_{[-10R, 10R]}\big(P_{[-R, R]}(\tilde{\chi}^{nab}[\tilde{V}_{3,4}^{nab}]^{2})e^{-i\sum_{k=1,2}\del_k^{-1}\phi^{na}_{k}}\big)\big]\|_{L_{x}^{2}}\\
&+\|\nabla^{-1}\big[\chi^{nab}P_{[-10R, 10R]}\big(P_{[-R,
R]}(\tilde{\chi}^{nab}\tilde{V}_{3,4}^{nab}\tilde{W}_{3,4}^{naB})e^{-i\sum_{k=1,2}\del_k^{-1}\phi^{na}_{k}}\big)\big]\|_{L_{x}^{2}}\ll
\delta_{2}
\end{split}\end{equation}
for $n$ large enough. Finally, consider the term
\[
\nabla^{-1}\big[\chi^{nab}P_{[-10R, 10R]}\big(P_{[-R,
R]}(\tilde{\chi}^{nab}[\tilde{W}_{3,4}^{naB}]^{2})e^{-i\sum_{k=1,2}\del_k^{-1}\phi^{na}_{k}}\big)\big]
\]
Here we split
\[
\tilde{W}_{3,4}^{naB}=P_{[-R_{1},
R_{1}]^{c}}\tilde{W}_{3,4}^{naB}+P_{[-R_{1},
R_{1}]}\tilde{W}_{3,4}^{naB}
\]
Then if $B$ is chosen large enough in relation to $R_{1}$, we obtain
both
\[
\|P_{[-R_{1}, R_{1}]^{c}}\tilde{W}_{3,4}^{naB}\|_{L_{x}^{2}}\ll
\delta_{2}, \|P_{[-R_{1},
R_{1}]}\tilde{W}_{3,4}^{naB}\|_{L_{x}^{\infty}}\ll \delta_{2}
\]
Here the first inequality holds of course uniformly ion in $n, B$
due to the frequency localization. From here we infer that for $B$
and then $n$ large enough, we get
\[
\|\nabla^{-1}\big[\chi^{nab}P_{[-10R, 10R]}\big(P_{[-R,
R]}(\tilde{\chi}^{nab}[\tilde{W}_{3,4}^{naB}]^{2})e^{-i\sum_{k=1,2}\del_k^{-1}\phi^{na}_{k}}\big)\big]\|_{L_{x}^{2}}\ll
\delta_{2}
\]
The argument for showing
\begin{equation}\nonumber\begin{split}
&\|\nabla^{-1}\big[\chi^{nab}P_{[-10R, 10R]}\big(P_{[-R, R]}(\phi^{na}_{2})\partial_{1}(e^{-i\sum_{k=1,2}\del_k^{-1}\phi^{na}_{k}})\big)\big]\\
&-\nabla^{-1}\big[\chi^{nab}P_{[-10R, 10R]}\big(P_{[-R, R]}(\phi^{na}_{1})\partial_{2}(e^{-i\sum_{k=1,2}\del_k^{-1}\phi^{na}_{k}})\big)\big]\|_{L_{x}^{2}}\ll \delta_{2}
\end{split}\end{equation}
of course proceeds in identical fashion.
\\

Summarizing what we have achieved thus far in Case (B), we have
shown that for $n$ large enough, we get
\begin{equation}\nonumber\begin{split}
&\|\nabla^{-1}\partial_{1}\big[\chi^{nab}P_{[-10R, 10R]}\big(P_{[-R, R]}(\phi^{na}_{2})e^{-i\sum_{k=1,2}\del_k^{-1}\phi^{na}_{k}}\big)\big]\\
&-\nabla^{-1}\partial_{2}\big[\chi^{nab}P_{[-10R, 10R]}\big(P_{[-R,
R]}(\phi^{na}_{1})e^{-i\sum_{k=1,2}\del_k^{-1}\phi^{na}_{k}}\big)\big]\|_{L_{x}^{2}}\ll
\delta_{2}
\end{split}\end{equation}
In light of the fact pointed out earlier that ($j=1,2$)
\[
\big[\chi^{nab}P_{[-10R, 10R]}\big(P_{[-R,
R]}(\phi^{na}_{j})e^{-i\sum_{k=1,2}\del_k^{-1}\phi^{na}_{k}}\big)\big]
\]
is well approximated by
\[
\chi^{nab}[V_{j+2}^{nab}+W_{j+2}^{naB}],
\]
we then infer that (recalling the definition of $V_{3,4}$,
$W_{3,4}$)
\begin{equation}\nonumber\begin{split}
&\|\nabla^{-1}\partial_{1}\big[\chi^{nab}\partial_{2}[((S_{\tilde{A}^{nab}}(V_{1}^{ab}[0]))(0-t^{nab},
x-x^{nab})+W_{1}^{naB})]\\&
-\partial_{1}[(S_{\tilde{A}^{nab}}(V_{2}^{ab}[0]))(0-t^{nab}, x-x^{nab})+W_{2}^{naB}]\big]\\
&-\nabla^{-1}\partial_{2}\big[\chi^{nab}\partial_{1}[((S_{\tilde{A}^{nab}}(V_{1}^{ab}[0]))(0-t^{nab},
x-x^{nab})+W_{1}^{naB})]\\&
+\partial_{2}[(S_{\tilde{A}^{nab}}(V_{2}^{ab}[0]))(0-t^{nab}, x-x^{nab})+W_{2}^{naB}]\big]\|_{L_{x}^{2}}\\
&\ll \delta_{2}
\end{split}\end{equation}
But then choosing the cutoff $\chi^{nab}$ as above and picking $n$
large enough, we conclude (noting cancelations in the preceding
expression) that
\[
\|\nabla^{-1}\triangle[(S_{\tilde{A}^{nab}}(V_{2}^{ab}[0]))(0-t^{nab},
x-x^{nab})+\chi^{nab}W_{2}^{naB}]\|_{L_{x}^{2}}\ll \delta_{2}
\]
This inequality, together with the approximate orthogonality of the
two summands involved, then gives the smallness of either summand
separately: recall from Lemma~\ref{CoreBGII} and its proof  that we
have
\[
\Big|\int_{\R^{2}}\nabla_{x,t}(S_{\tilde{A}^{nab}}(V_{2}^{ab}[0]))(0-t^{nab},
x-x^{nab})\cdot\overline{\nabla_{x,t}W_{2}^{naB}(0,\cdot)}\, dx\Big|\ll
\delta_{2}
\]
Now recall the vanishing condition at time $t=0$
\[
\sum_{b'=1}^{B}\partial_{t}S_{\tilde{A}^{nab'}}(V_{2}^{ab'}[0]))(0-t^{nab'},
x-x^{nab'})+\partial_{t}W_{2}^{naB}=0
\]
which we used to define the linear covariant evolution of
$\eta^{na}$. Applying the cutoff $\chi^{nab}$, and choosing $n$
large enough, we get that
\[
\|\partial_{t}S_{\tilde{A}^{nab}}(V_{2}^{ab}[0]))(0-t^{nab},
x-x^{nab})+\chi^{nab}\partial_{t}W_{2}^{naB}\|_{L_{x}^{2}}\ll
\delta_{2}
\]
However, this inequality, together with the two preceding ones,
implies that
\[
\|\nabla_{x,t}(S_{\tilde{A}^{nab}}(V_{2}^{ab}[0]))(0-t^{nab},
x-x^{nab})\|_{L_{x}^{2}}+\|\chi^{nab}\nabla_{x,t}W_{2}^{naB}\|_{L_{x}^{2}}\ll
\delta_{2}
\]

Summarizing the state of affairs in Case (B), we have shown thus far
that the Claim holds. But this then says that the 'diluted
concentration profile' given by
\[
V_{3}^{nab}=[\partial_{1}S_{\tilde{A}^{nab}}(V_{1}^{ab}[0])+\partial_{2}S_{\tilde{A}^{nab}}(V_{2}^{ab}[0])](0-t^{nab},
x-x^{nab})
\]
\[
V_{4}^{nab}:=[\partial_{2}S_{\tilde{A}^{nab}}(V_{1}^{ab}[0])-\partial_{1}S_{\tilde{A}^{nab}}(V_{2}^{ab}[0])](0-t^{nab},
x-x^{nab})
\]
is given, up to an $L^{2}$-error of size $\delta_{2}$, by the
{\it{pure gradient term}}
\[
V_{3}^{nab}=\partial_{1}S_{\tilde{A}^{nab}}(V_{1}^{ab}[0])](0-t^{nab},
x-x^{nab})
\]
\[
V_{4}^{nab}=\partial_{2}S_{\tilde{A}^{nab}}(V_{1}^{ab}[0])](0-t^{nab},
x-x^{nab})
\]
We shall now use this to construct a map from $\R^{2}\to \Hyp^2$ whose Coulomb derivative components are close to
$V_{3,4}^{nab}$.
\\

Indeed, picking $R$ large enough and then $n$ sufficiently large
depending on $R$, it is straightforward to check that
\begin{equation}\nonumber\begin{split}
&\partial_{j}P_{[-R,
R]}\big(S_{\tilde{A}^{nab}}(V_{1}^{ab}[0])](0-t^{nab},
x-x^{nab})\big)
\\&=\partial_{j}P_{[-R, R]}\big(S_{\tilde{A}^{nab}}(V_{1}^{ab}[0])](0-t^{nab}, x-x^{nab})\big)
e^{-i\sum_{k=1,2}\del_k^{-1}\del_k P_{[-R,
R]}\big(S_{\tilde{A}^{nab}}(V_{1}^{ab}[0])](0-t^{nab},
x-x^{nab})\big)}+\text{error},
\end{split}\end{equation}
where we have $\|\text{error}\|_{L_{x}^{2}}\ll \delta_{2}$. Then we
define a map $({\bf{x}}, {\bf{y}}): \R^{2}\to
\Hyp^2$ (here we abuse notation heavily, this map of
course depends on $n, a, b$) via
\[
{\bf{x}}:=\Re P_{[-R,
R]}\big(S_{\tilde{A}^{nab}}(V_{1}^{ab}[0])](0-t^{nab},
x-x^{nab})\big),\quad {\bf{y}}:=\Im P_{[-R,
R]}\big(S_{\tilde{A}^{nab}}(V_{1}^{ab}[0])](0-t^{nab},
x-x^{nab})\big),
\]
These then satisfy
\[
\Big\|\frac{\partial_{j}{\bf{x}}}{\bf{y}}+i\frac{\partial_{j}{\bf{y}}}{{\bf{y}}}-\partial_{j}P_{[-R,
R]}\big(S_{\tilde{A}^{nab}}(V_{1}^{ab}[0])](0-t^{nab},
x-x^{nab})\big)\Big\|_{L_{x}^{2}}\ll \delta_{2},
\]
and the associated Coulomb derivative components are the desired
approximations. This concludes the proof of
Proposition~\ref{ConcentrarionProfileApprox}.
\end{proof}

Summary thus far, for both (A), (B): {\it{we have shown that
we have the ``covariant Bahouri Gerard decompositions''
\[
\phi_{1}^{na}e^{-i\sum_{k=1,2}\del_k^{-1}[w_{k}^{nA_{0}^{(0)}}+\phi^{1na}_{k}]}=\sum_{b=1}^{B}V_{3}^{nab}(0-t^{nab},
x-x^{nab})+W_{3}^{nab}
\]
\[
\phi_{2}^{na}e^{-i\sum_{k=1,2}\del_k^{-1}[w_{k}^{nA_{0}^{(0)}}+\phi^{1na}_{k}]}=\sum_{b=1}^{B}V_{4}^{nab}(0-t^{nab},
x-x^{nab})+W_{4}^{nab},
\]
where we have
\[
V_{3}^{nab}:=\partial_{1}S_{\tilde{A}^{nab}}(V_{1}^{ab}[0])+\partial_{2}S_{\tilde{A}^{nab}}(V_{2}^{ab}[0])
\]
\[
V_{4}^{nab}:=\partial_{2}S_{\tilde{A}^{nab}}(V_{1}^{ab}[0])-\partial_{1}S_{\tilde{A}^{nab}}(V_{2}^{ab}[0])
\]
and similarly for $W_{3,4}$. Furthermore, for $n$ large enough and
any given $\delta_{2}>0$, we can find maps
$({\bf{x}}^{\delta_{2}nab}, {\bf{y}}^{\delta_{2}nab}) :
\R^{2}\to\Hyp^2$, with the property that their
(spatial) Coulomb derivative components are $\delta_{2}$ close
(within the $L^{2}$-metric) to constant phase shifts of the
$V_{3,4}^{nab}(0-t^{nab}, x-x^{nab})$.}}
\\

We shall now refine the information we have by proving the following

\begin{lemma}\label{BGIIErrorControl} Given $\delta_{2}>0$, we can pick $B$ and then $n$ large enough such that
\[
\|\nabla_{x,t}W_{2}^{naB}\|_{L_{x}^{2}}\ll \delta_{2}
\]
\end{lemma}

\begin{remark}Recalling the identities
\[
W_{3}^{naB}=\partial_{1}W_{1}^{naB}+\partial_{2}W_{2}^{naB}
\]
\[
W_{4}^{naB}=\partial_{2}W_{1}^{naB}-\partial_{1}W_{2}^{naB}
\]
We see that this says that $W_{3,4}^{naB}$ are essentially pure
gradient terms, like in Case (B).
\end{remark}

\begin{proof}(Lemma~\ref{BGIIErrorControl}) The proof is quite similar to the Case (B) above. Given $\delta_{2}>0$,
first choose an index $B_{1}$ such that we have
\[
\limsup_{n\to\infty}\|\sum_{b=B_{1}}^{B}V_{3,4}^{nab}(0-t^{nab},
x-x^{nab})\|_{L_{x}^{2}}\ll \delta_{2},
\]
for any $B\geq B_{1}$. Further, pick $R=R(\delta_{2})$ with the
property that
\[
\limsup_{n\to\infty}\|P_{[-R,
R]^{c}}(\phi^{na}_{j})e^{-i\sum_{k=1,2}\del_k^{-1}[\phi^{1na}_{k}+w^{1nA_{0}^{(0)}}_{k}]}\|_{L_{x}^{2}}\ll
\delta_{2},\quad j=1,2
\]
Increasing $R$ if necessary, we can then also achieve that (for $n$
large enough)
\begin{equation}\nonumber\begin{split}
&\|P_{[-10R, 10R]}\big[P_{[-R,
R]}(\phi^{na}_{1,2})e^{-i\sum_{k=1,2}\del_k^{-1}[\phi^{1na}_{k}+w^{1nA_{0}^{(0)}}_{k}]}\big]-\sum_{b=1}^{B_{1}}V_{3,4}^{nab}(0-t^{nab},
x-x^{nab})-W_{3,4}^{naB}\|_{L_{x}^{2}}\ll \delta_{2}
\end{split}\end{equation}
Here we will choose $B$ sufficiently large in relation to $B_{1},
\delta_{2}$. Now pick a cutoff $\chi$ which localizes to the union
of large discs covering most of the support (in the $L^{2}$-sense)
of the atoms $V_{3,4}^{nab}(0-t^{nab}, x-x^{nab})$ of bounded type,
i.e., for which $\limsup|t^{nab}|<\infty$, $1\leq b\leq B_{1}$. Of
course $\chi$ then depends on $a, B_{1}, n$, but we suppress this
dependence here. Picking $\chi$ suitably and then choosing $n$ large
enough, we can then ensure that
\[
\|(1-\chi)[\sum_{b=1}^{B_{1}}V_{3,4}^{nab}(0-t^{nab},
x-x^{nab})-\sum_{b=1}^{'B_{1}}V_{3,4}^{nab}(0-t^{nab},
x-x^{nab})]\|_{L_{x}^{2}}\ll \delta_{2},
\]
where $\sum_{b=1}^{'B_{1}}$ indicates that we only sum over the
atoms of ``unbounded type''. Summarizing the above steps, we now have
\[
\|(1-\chi)P_{[-10R, 10R]}\big[P_{[-R,
R]}(\phi^{na}_{1,2})e^{-i\sum_{k=1,2}\del_k^{-1}[\phi^{1na}_{k}+w^{1nA_{0}^{(0)}}_{k}]}\big]-\sum_{b=1}^{'B_{1}}V_{3,4}^{nab}(0-t^{nab},
x-x^{nab})-(1-\chi)W_{3,4}^{naB}\|_{L_{x}^{2}}\ll \delta_{2}
\]
By picking $B$ large enough (recall that we can do so independently
of $B_{1}$), we may also assume that
\[
\|(1-\chi)W_{3,4}^{naB}-W_{3,4}^{naB}\|_{L_{x}^{2}}\ll \delta_{2}
\]
Here we use Lemma~\ref{TailDispersion}.
\\
Next, we calculate the curl of the Coulomb components, localized as
above, and with an extra cutoff $(1-\chi)$.  Thus we want to
estimate the expression
\begin{equation}\nonumber\begin{split}
&\partial_{2}\big((1-\chi)P_{[-10R, 10R]}\big[P_{[-R,
R]}(\phi^{na}_{1})e^{-i\sum_{k=1,2}\del_k^{-1}[\phi^{1na}_{k}+w^{1nA_{0}^{(0)}}_{k}]}\big]\big)\\&-\partial_{1}\big((1-\chi)P_{[-10R,
10R]}\big[P_{[-R,
R]}(\phi^{na}_{2})e^{-i\sum_{k=1,2}\del_k^{-1}[\phi^{1na}_{k}+w^{1nA_{0}^{(0)}}_{k}]}\big]\big)
\end{split}\end{equation}
This we can estimate as in Case~(B):  Of course the case when a
derivative falls on $(1-\chi)$ is negligible. Then repeating the
arguments in Case~(B) above, we need to estimate the schematic
expression
\[
\big((1-\chi)P_{[-10R, 10R]}\big[P_{[-R,
R]}([\phi^{na}]^{2})e^{-i\sum_{k=1,2}\del_k^{-1}[\phi^{1na}_{k}+w^{1nA_{0}^{(0)}}_{k}]}\big]\big)
\]
Here we use the Bahouri Gerard decomposition of the $\phi^{na}$,
i.e.,
\[
\phi^{na}_{1,2}=\sum_{b=1}^{B}V_{3,4}(0-t^{nab},
x-x^{nab})+W_{3,4}^{naB}
\]
It is clear that we then reduce to estimating
 \[
\big((1-\chi)P_{[-10R, 10R]}\big[P_{[-R,
R]}((1-\tilde{\chi})[\sum_{b=1}^{'B_{1}}V_{3,4}(0-t^{nab},
x-x^{nab})+W_{3,4}^{naB}
]^{2})e^{-i\sum_{k=1,2}\del_k^{-1}[\phi^{1na}_{k}+w^{1nA_{0}^{(0)}}_{k}]}\big]\big)
\]
But the contribution of the terms $V_{3,4}(0-t^{nab}, x-x^{nab})$,
$1\leq b\leq B_{1}$ can be made arbitrarily small by choosing $n$
large enough, while the contribution of $W_{3,4}^{naB}$ is handled
by placing one factor into $L_{x}^{2}$ and the other into
$L_{x}^{\infty}$.
\\

Summarizing, we have now shown that
\begin{equation}\nonumber\begin{split}
&\Big\|\nabla^{-1}\partial_{2}\big[\sum_{b=1}^{'B_{1}}V_{3}^{nab}(0-t^{nab},
x-x^{nab})+W_{3}^{naB}\big]
-\nabla^{-1}\partial_{2}\big[\sum_{b=1}^{'B_{1}}V_{4}^{nab}(0-t^{nab},
x-x^{nab})+W_{4}^{naB}\big]\Big\|_{L_{x}^{2}}\ll \delta_{2}
\end{split}\end{equation}
But then recalling the defining relations for $V_{3, 4}^{nab},
W_{3,4}^{naB}$, we can repeat the argument from part (B) in the
preceding proof to conclude that for $B$ and then $n$ large enough,
we have
\[
\|\nabla_{x,t}W_{2}^{naB}\|_{L_{x}^{2}}\ll \delta_{2},
\]
as desired.
\end{proof}

Proposition~\ref{ConcentrarionProfileApprox} together with
Lemma~\ref{BGIIErrorControl} are key technical tools we shall use in
the next section when bounding the wave maps with data
\[
w^{nA_{0}^{0}}+\phi^{na},
\]
where $a=1$.

\subsection{Step 4: Adding the first large atomic component and invoking the induction hypothesis}

In Step~3, we constructed a wave map with data corresponding to the
lowest frequency ``non-atomic'' part, whose Coulomb components are
\[
\Phi_{\alpha}^{nA_{0}^{(0)}}e^{-i\sum_{k=1,2}\triangle^{-1}\partial_{k}\Phi_{k}^{nA_{0}^{(0)}}}= w_{\alpha}^{nA_{0}^{(0)}}
e^{-i\sum_{k=1,2}\triangle^{-1}\partial_{k}w_{k}^{nA_{0}^{(0)}}}  +o_{L^{2}}(1)
\]
Our next step now is to prove bounds for the wave map whose Coulomb
components are given by
\[
\psi^{n(<1)}_{\alpha}:=[w_{\alpha}^{nA_{0}^{(0)}}+\phi_{\alpha}^{n1}]e^{-i\sum_{k=1,2}\triangle^{-1}\partial_{k}[w_{k}^{nA_{0}^{(0)}}
+\phi^{n1}_{k}]}+o_{L^{2}}(1),
\]
provided we make the following key
\\
{\bf{Energy Assumption:}} {\em All concentration profiles have energy $<\Ecrit$. Thus}
\begin{align}\label{eq:physicalnonconc}
E(V^{ab})<\Ecrit\,\forall\,b
\end{align}
As before, in order to avoid confusion, we shall denote the
superscript $1$ here instead by $a$, it being understood that $a=1$.
Thus we now intend to prove global bounds for the evolution of the
Coulomb data
\begin{equation}\nonumber\begin{split}
&\psi^{n(<a)}_{\alpha}:=[w_{\alpha}^{nA_{0}^{(0)}}+\phi_{\alpha}^{na}]e^{-i\sum_{k=1,2}\triangle^{-1}\partial_{k}[w_{k}^{nA_{0}^{(0)}}+\phi^{na}_{k}]}+o_{L^{2}}(1)\\
&=w_{\alpha}^{nA_{0}^{(0)}}e^{-i\sum_{k=1,2}\triangle^{-1}\partial_{k}w_{k}^{nA_{0}^{(0)}}}+\phi_{\alpha}^{na}e^{-i\sum_{k=1,2}\triangle^{-1}\partial_{k}[w_{k}^{nA_{0}^{(0)}}+\phi^{na}_{k}]}+o_{L^{2}}(1)
\end{split}\end{equation}
From the preceding section, we obtain a decomposition of the added
term
\[
\tilde{\psi}^{na}_{\alpha}:=\phi_{\alpha}^{na}e^{-i\sum_{k=1,2}\triangle^{-1}\partial_{k}[w_{k}^{nA_{0}^{(0)}}+\phi^{na}_{k}]}
\]
as a sum of concentration profiles {\em{at time $t=0$}}. Note that
this time in principle plays no distinguished role, other than that
we are guaranteed existence of the evolution of the wave maps with
above data on some small time interval centered at $t=0$. Recall the
decompositions (for any $B\geq 1$)
\[
\phi_{1}^{na}e^{-i\sum_{k=1,2}\triangle^{-1}\partial_{k}[w_{k}^{nA_{0}^{(0)}}+\phi^{na}_{k}]}
=\sum_{b=1}^{B}V_{3}^{nab}(0-t^{nab}, x-x^{nab})+W_{3}^{naB}
\]
\[
\phi_{2}^{na}e^{-i\sum_{k=1,2}\triangle^{-1}\partial_{k}[w_{k}^{nA_{0}^{(0)}}+\phi^{na}_{k}]}
=\sum_{b=1}^{B}V_{4}^{nab}(0-t^{nab}, x-x^{nab})+W_{4}^{naB}
\]
where
\[
V_{3}^{nab}:=\partial_{1}S_{\tilde{A}^{nab}}(V_{1}^{ab}[0])+\partial_{2}S_{\tilde{A}^{nab}}(V_{2}^{ab}[0])
\]
\[
V_{4}^{nab}:=\partial_{2}S_{\tilde{A}^{nab}}(V_{1}^{ab}[0])-\partial_{1}S_{\tilde{A}^{nab}}(V_{2}^{ab}[0])
\]
and similarly for $W_{3,4}$, while we also have
\[
\phi_{0}^{na}e^{-i\sum_{k=1,2}\triangle^{-1}\partial_{k}[w_{k}^{nA_{0}^{(0)}}+\phi^{na}_{k}]}
=\sum_{b=1}^{B}\partial_{t}S_{\tilde{A}^{nab}}(V_{1}^{ab}[0])(0-t^{nab},
x-x^{nab})+\partial_{t}W_{1}^{naB},
\]
\[
\sum_{b=1}^{B}\partial_{t}S_{\tilde{A}^{nab}}(V_{2}^{ab}[0])(0-t^{nab},
x-x^{nab})+\partial_{t}W_{2}^{naB}=0
\]
{\em{These decompositions are understood to hold at time $t=0$, of
course.}} Now the fact that for $B$ large enough (and then $n$ large
enough) we can arrange that
$\|\nabla_{x,t}W_{2}^{naB}\|_{L_{x}^{2}}\ll \delta_{2}$ implies that
\[
\sum_{b=1}^{B}\|\partial_{t}S_{\tilde{A}^{nab}}(V_{2}^{ab}[0])(0-t^{nab},
x-x^{nab})\|_{L_{x}^{2}}^{2}\ll \delta_{2}
\]
by increasing $n$ if necessary, due to the approximate orthogonality
of these functions.
\\

Recall that we have {\it{temporally bounded concentration
profiles}}, as well as {\it{temporally unbounded ones}}. Then it is
intuitively clear that the evolution of $\tilde{\psi}^{na}$ (this is
not well-defined strictly speaking, we can only evolve Coulomb
components of actual maps; however, we can think of
$\tilde{\psi}^{na}$ as the difference between the components of
maps) will be dominated for a large time interval around $t=0$ by
the evolution of the temporally bounded concentration profiles,
which will exhibit nonlinear behavior, while the temporally
unbounded ones will behave like free waves for a long time. In order
to make things precise, we introduce a {\it{hierarchy of temporal
scales}}, which means we order the times $t^{nab}$ according to
whether they are positive or negative and then whether
\[
\lim_{n\to\infty}(t^{nab}-t^{nab'})=\pm\infty
\]
Assume that this way, we arrive at the list of representative time
scales, $M=M(B)$,
\[
0=t^{nab_{1}}, t^{nab_{2}},\ldots, t^{nab_{M}}
\]
where we have $t^{nab_{i}}>0$, say, and
\[
\lim_{n\to\infty}(t^{nab_{j}}-t^{nab_{j-1}})=\infty,
\]
and furthermore for each $b\in\{1,2,\ldots, B\}$, we have
$t^{nab}=t^{nab_{j}}$ for some $j$ as above. Note that we have
chosen to equate those times here that do not diverge from each
other. This can of course be done by passing to a subsequence such
that the difference of these times converges, and the redefining the
concentration profiles accordingly.

 We then implement an
inductive procedure, controlling the evolution of
$\psi^{n(<a)}_{\alpha}$ on the interval $[0, t^{nab_{2}}-C]$ for
some huge $C$ (such that we are guaranteed that all the
concentration profiles focussing at times $t^{nab_{j}}, j\geq 2$,
will not display any nonlinear behavior there yet), while the
temporally bounded ones start to disperse and behave linearly around
time $t^{nab_{2}}-C$, for sufficiently large $n$. {\it{This then
guarantees that there is essentially no nonlinear interactions going
on between evolutions of concentration profiles at different time
scales.}}

\subsubsection{Proving apriori bounds for the evolution of $\psi^{n(<a)}_{\alpha}$; the lowest time scale.}

Here we prove apriori bounds on the (wave map) evolution of the
Coulomb components $\psi_{\alpha}^{n(<a)}$. Recall that at time
$t=0$, we have the decomposition
\[
\psi_{1,2}^{n(<a)}(0,\cdot)=w_{1,2}^{nA_{0}^{(0)}}e^{-i\sum_{k=1,2}\triangle^{-1}\partial_{k}w_{k}^{nA_{0}^{(0)}}}+\sum_{b=1}^{B}V_{3,4}^{nab}(0-t^{nab},
x-x^{nab})+W_{3,4}^{naB}+o_{L^{2}}(1)
\]
\[
\psi_{0}^{n(<a)}(0,\cdot)=w_{0}^{nA_{0}^{(0)}}e^{-i\sum_{k=1,2}\triangle^{-1}\partial_{k}w_{k}^{nA_{0}^{(0)}}}+\sum_{b=1}^{B}\partial_{t}S_{\tilde{A}^{nab}}(V_{1}^{ab}[0])(0-t^{nab},
x-x^{nab})+\partial_{t}W_{1}^{naB}
\]
The ``tail ends'' $W_{3,4}^{naB}$, $\partial_{t}W_{1}^{naB}$ here
satisfy the smallness condition
\[
W_{3,4}^{naB},\,\partial_{t}W_{1}^{naB}=o_{L^{2}}(1)+o_{L^{\infty}}(1)
\]
where $o(\cdot)$ here is meant in case $B, n\to\infty$. Observe
from  Lemma~\ref{BGIIErrorControl} that we actually have
\[
W_{3,4}^{naB}=\partial_{1,2}W_{1}^{naB}+\text{error},\quad
\|\text{error}\|_{L_{x}^{2}}\ll \delta_{2},
\]
provided we choose $B$ and then $n$ large enough. Furthermore, the
proof of Proposition~\ref{ConcentrarionProfileApprox}, case (B),
reveals that for concentration profiles which are {\it{temporally
unbounded}}, we have
\[
V_{3,4}^{nab}(0-t^{nab},
x-x^{nab})=\partial_{1,2}S_{\tilde{A}^{nab}}(V_{1}^{ab}[0])(0-t^{nab},
x-x^{nab})
\]
We shall now build the evolution of $\psi_{\alpha}^{n(<a)}$ as the
sum of well-known pieces, namely the evolutions of the atomic
profiles, plus an error term, which we will show will remain small.
To make things precise, we now use the following construction: We
shall use $\delta_{2}>0$ as a smallness parameter which will
ultimately hinge on intrinsic properties of the concentration
profiles as well as the $S$-bound on the already constructed low frequency part $\Psi_\alpha^{nA_0^(0)}$,
and be specified at the end of the construction. Thinking
of $\delta_{2}>0$ as fixed for now, we first pick a large cutoff
$B_{1}$ with the property that
\[
\limsup_{n\to\infty}\|\sum_{B\geq b\geq
B_{1}}V_{4,3}^{nab}\|_{L_{x}^{2}}+\|\sum_{B\geq b\geq
B_{1}}\partial_{t}S_{\tilde{A}^{nab}}(V_{1}^{ab}[0])\|_{L_{x}^{2}}
+\|\partial_{\alpha}W_1^{naB}\|_{L_{t,x}^\infty}\ll
\delta_{2}
\]
for any $B\geq B_{1}$. Then we evolve the concentration profiles
corresponding to a $b\in \{1,2,\ldots, B_{1}\}$ as follows:

\smallskip

{\bf{(I): Evolution of temporally bounded concentration profile.}}

\smallskip
Here, by passing to a subsequence, we may assume that $t^{nab}$
converges as $n\to\infty$, and we may then set $t^{nab}=0$
by time translation. Also, it is apparent that then
\begin{align*}
V_{3}^{nab}&=\partial_{1}V_{1}^{ab}(0,\cdot))+\partial_{2}V_{2}^{ab}(0,\cdot)+o_{L^{2}}(1)\\
V_{4}^{nab}
&=\partial_{2}V_{1}^{ab}(0,\cdot))-\partial_{1}V_{2}^{ab}(0,\cdot)+o_{L^{2}}(1),
\partial_{t}S_{\tilde{A}^{nab}}(V_{1}^{ab}[0])\\
&=\partial_{t}V_{1}^{ab}(0,\cdot)+o_{L^{2}}(1)
 \end{align*}
 are all
essentially independent of $n$. Now according to
Proposition~\ref{ConcentrarionProfileApprox}, we can find, for each
$\delta_{2}>0$, a constant phase $\gamma_{\delta_{2}ab}$ and an
admissible  map from $\R^{2}\to \Hyp^2$ whose Coulomb components
$\psi^{ab\delta_{2}}_{\alpha}$ satisfy
\[
\|e^{i\gamma_{\delta_{2}ab}}\psi^{ab\delta_{2}}_{1,2}-V_{3,4}^{nab}\|_{L_{x}^{2}}\ll
\delta_{2},\quad
\|e^{i\gamma_{\delta_{2}ab}}\psi^{ab\delta_{2}}_{1,2}-\partial_{t}V_{1}^{ab}(0,\cdot)\|_{L_{x}^{2}}\ll
\delta_{2},
\]
For the sake of simplicity, we now refer to the Coulomb components
of such a map, which we choose for $\delta_{2}$ extremely small
(depending on $B_{1}$ etc.~and to be specified later), simply as
$\psi^{ab}_{\alpha}$.

\noindent  First, we evolve the components of $\psi^{ab}_{\alpha}$
on a large time interval $I^{ab}$ centered at $t=0$, using the wave
maps flow for the Coulomb components. This yields an apriori bound
\[
\|\psi^{ab}\|_{S}<C_{ab}
\]
due to our energy assumption \eqref{eq:physicalnonconc}.
Furthermore, due to Corollary~\ref{cor:localsplit2} as well as
Remark~\ref{rem:psiLimprov}, given $\delta_{2}>0$, one can then
choose time intervals
\[
I_{1}, I_{2}, \ldots, I_{M_{ab}(\delta_{2})},
\]
where the final one is of the form $[t_{1}^{ab\delta_2},\infty)$,
say, such that
\[
\psi^{ab}|_{I_{j}}=\psi^{ab}_{jL}+\psi^{ab}_{jNL}
\]
with\footnote{The implied constant in the second inequality here is
universal, independent of $\delta_2$.}
\[
\|\psi^{ab}_{jNL}\|_{S(I_{j}\times\R^{2})}\ll \delta_{2},\quad
\|\nabla_{x,t}\psi^{M_{ab}ab}_{L}\|_{L_t^\infty\dot{H}^{-1}}\lesssim
\Ecrit
\]
Here of course $\Box \psi^{ab}_{jL}=0$. Note that the intervals
$I_{j}$ here only depend on $a, b$ as well as the smallness
parameter $\delta_{2}$. By the Huyghen's principle, one may assume
that the support of $\psi^{ab}_{jL}$ is contained in the set
$|x|\leq |t|+D_{ab}(\delta_{2})$ for some (possibly very large
number $D_{ab}(\delta_{2})$). But then by choosing a much larger
time $T^{ab\delta_{2}}$, we can arrange that
\[
\|\psi^{ab}_{M_{ab}L}([T^{ab\delta_{2}},\infty), \cdot)\|_{L_{t,x}^{\infty}+L_t^\infty L_x^2}\ll
\delta_{2}
\]
These considerations reveal that pursuing the wave maps evolution of
the components $\psi^{ab}$ long enough, we eventually find that
\[
\psi^{ab}(t, \cdot)=o_{L^{2}}(1)+o_{L^{\infty}}(1),
\]
where $o(\cdot)$ is in the sense as $|t|\to\infty$.
\\
This is conclusion is of critical importance: note that thus far we
have not taken  the low frequency contribution from
$w^{nA_{0}^{(0)}}$ (from Step~3) into account, which starts to play
an important role for extremely large times. The above asymptotic
description allows us to incorporate this low-frequency effect by
adjusting the linear evolution of $\psi^{ab}_{M_{ab}L}$ from flat to
covariant. In more precise terms, we now make the following choice
of an extension $\tilde{\psi}^{ab}_{\alpha}$ of the data
$\psi^{ab}_{\alpha} $:

\begin{itemize}
\item
 On the interval $[0, T^{ab\delta_{2}}]$, we let
$\tilde{\psi}^{ab}=\psi^{ab}$.
\item
 On the interval $[T^{ab\delta_{2}}, \infty)$, we let
$\tilde{\psi}^{ab}$ be the covariant extension of
$\psi^{ab}[T^{ab\delta_{2}}]$, i.e., we have
\[
\Box_{A^{n}}\tilde{\psi}^{ab}=0
\]
on $[T^{ab\delta_{2}}, \infty)$, where $A^{n}$ is defined with
inputs
$w^{nA_{0}^{(0)}}e^{-i\sum_{k=1,2}\triangle^{-1}w_{k}^{nA_{0}^{(0)}}}$.
More precisely, we apply a Hodge decomposition to the data
$\tilde{\psi}^{ab}_{\alpha}$ as in \eqref{eq:yetanotherHodge3},
\eqref{eq:yetanotherHodge4}, and evolve these components as in
Step~3. In order to avoid a ``kink'' at the juncture of these two
regimes, we define
\begin{equation}\label{glue}
\tilde{\psi}^{ab}=\chi_{(-\infty, T^{ab\delta_{2}}+10)}(t)\psi^{ab}+(1-\chi_{(-\infty, T^{ab\delta_{2}}+10)})(t)S_{A^n}\psi^{ab}[T^{ab\delta_{2}}]
\end{equation}
where the notation for the second term is schematic, and $\chi_{(-\infty, T^{ab\delta_{2}}+10)}(t)$ smoothly localizes to
the indicated interval and satisfies
\[
\chi_{(-\infty, T^{ab\delta_{2}}+10)}|_{[0, T^{ab\delta_{2}}]}=1
\]
\end{itemize}

With these definitions, one can prove the following bound.

\begin{prop}\label{BGIITemporallyBounded} We have a bound of the form
\[
\|\tilde{\psi}^{ab}\|_{S}<C(C_{ab}),
\]
where we recall the assumption $\|\psi^{ab}\|_{S}<C_{ab}$ from
above. Furthermore, denoting by $c_{k}$, $k\in\Z$, a frequency
envelope controlling the data at time $t=0$, i.e.,
\[
c_{k}=(\sum_{l\in\Z}2^{-\sigma|l-k|}\|P_{l}\psi^{ab}(0,\cdot)\|_{L_{x}^{2}}^{2})^{\frac{1}{2}}
\]
for sufficiently small apriori constant $\sigma>0$, one has
\[
\|P_{k}\tilde{\psi}^{ab}\|_{S[k]}\leq C(C_{ab})c_{k}
\]
\end{prop}
\begin{proof}
The proof of this follows from
 Proposition~\ref{TwistedWaveEquation}, as well as Lemma~\ref{lem:chiS} and its proof.
\end{proof}

The idea now in Case (I) is to use $\tilde{\psi}^{ab}$ as
approximate evolution of the data $\psi^{ab}$ globally in time, for
$n$ large enough. Thus $\tilde{\psi}^{ab}|_{[0, T^{ab\delta_{2}}]}$
is the actual wave maps flow, while beyond time $T^{ab\delta_{2}}$,
we use the covariant linear evolution.

\smallskip

{\bf{(II): Evolution of temporally unbounded concentration
profile.}}

\smallskip
Here we have $\lim_{n\to\infty}|t^{nab}|=\infty$, and as before,
$1\leq b\leq B_{1}$, where we have chosen $B_{1}$ above. In this
case, using the argument from Case~(B) in the proof of
Proposition~\ref{ConcentrarionProfileApprox} and  arguing as at the
beginning of the preceding Case~(I) (we again write $\psi^{nab}$
instead of $\psi^{nab\delta_{2}}$),
\[
\psi^{nab}_{\alpha}=\partial_{\alpha}S_{\tilde{A}^{nab}}(V_{1}^{ab}[0])(0-t^{nab},
x-x^{nab})+\text{error},\quad \alpha=0,1,2,
\]
with $\|\text{error}\|_{L_{x}^{2}}\ll \delta_{2}$. In this case we
set
\[
\tilde{\psi}^{nab}_{\alpha}=\partial_{\alpha}S_{\tilde{A}^{nab}}(V_{1}^{ab}[0])(t-t^{nab},
x-x^{nab}),
\]
the covariant linear evolution. Of course this becomes inaccurate
when $t\to t^{nab}$ and the nonlinear effects start to become
relevant, but we recall that we are on the lowest time scale in this
subsection, i.e., $t\ll t^{nab_2}$. Then we have the following
bound.

\begin{prop}\label{BGIITemporallyUnBounded} There is a bound of the form
\[
\|\tilde{\psi}^{nab}\|_{S}<C(\Ecrit),
\]
Furthermore, denoting by $c_{k}$, $k\in\Z$, a frequency envelope
controlling the data at time $t=0$, i.e.,
\[
c_{k}=(\sum_{l\in\Z}2^{-\sigma|l-k|}\|P_{l}\psi^{nab}(0,\cdot)\|_{L_{x}^{2}}^{2})^{\frac{1}{2}}
\]
for sufficiently small apriori constant $\sigma>0$,
\[
\|P_{k}\tilde{\psi}^{ab}\|_{S[k]}\leq C(\Ecrit)c_{k}
\]
\end{prop}

\smallskip

{\bf{(III): Evolution of the weakly small error.}}

\smallskip
These are the components $W_{3,4}^{naB}, \partial_{t}W_{1}^{naB}$.
From Lemma~\ref{BGIIErrorControl}, we know that
\[
W_{3,4}^{naB}=\partial_{1,2}W_{1}^{naB}+\text{error},
\]
where we can force $\|\text{error}\|_{L_{x}^{2}}\ll \delta_{2}$ by
choosing $B$ and then $n$ large enough. We then evolve $W_{1}^{naB}$
using the covariant linear evolution, i.e.,
\[
\Box _{A^{n}}W_{1}^{naB}(t,x)=0,\quad W_{1}^{naB}[0]=(W_{1}^{naB},
\partial_{t}W_{1}^{naB}),
\]
and then  define $W_{3,4}^{naB}(t, x)=\partial_{1,2}W_{1}^{naB}(t,
x)$.

\smallskip

We have now defined the evolutions of all the ingredients of
$\tilde{\psi}^{na}_{\alpha}$. We claim that by choosing $\delta_{2}$
small enough and then $B$ and $n$ large enough, the sum of all these
constituents gives the correct evolution of
$\tilde{\psi}^{na}_{\alpha}$ up to a small error. This is clarified
the following  {\em{Core Proposition for Bahouri Gerard~II}} which
ties it all together.

\begin{prop}\label{BGIIHard} There is a cutoff $\delta_{2}>0$ sufficiently small, depending on the
 profiles $V_{1,2}^{ab}[0]$, $1\leq b\leq B_{1}$, as well as the apriori bound we have established for $\Psi^{nA_{0}^{(0)}}$,
 such that the following holds: picking $B_1$ and then $n$ large enough, we can write (with $B_{1}$ chosen as above)
 on $[0, t^{nab_{2}}-C]\times\R^{2}$ for $C$ sufficiently large and depending on the $\tilde{\psi}^{nab}_\alpha$ of unbounded type, $b=1,2,\ldots, B_1$,
 \[
\psi_{\alpha}^{n(<a)}(t,
x)=\Phi_{\alpha}^{nA_{0}^{(0)}}e^{-i\sum_{k=1,2}\triangle^{-1}\partial_{k}\Phi_{k}^{nA_{0}^{(0)}}}(t,
x)+\sum_{b=1}^{B_{1}}\tilde{\psi}^{nab}_{\alpha}(t,
x)+\partial_{\alpha}W_{1}^{naB_1}(t, x)+\epsilon_{\alpha}(t,
x),\quad \alpha=0,1,2
\]
where the components $\tilde{\psi}^{nab}_{\alpha}(t, x)$,
$\partial_{\alpha}W_{1}^{naB_1}(t, x)$, are constructed as in
(I)--(III) above, and with
\[
\|\epsilon \|_{S([0, t^{nab_{2}}-C]\times\R^{2})}\ll \delta_{2}
\]
Moreover, $\|P_k\epsilon\|_{S([0, t^{nab_{2}}-C]\times\R^{2})}$ is
exponentially decaying for frequencies $k>-\log(\lambda^a_n)$ Thus
the inequality above implies uniform smallness of $\epsilon(t, x)$
for $t\in [0, t^{nab_{2}}-C]$.
\end{prop}

\begin{remark} There appears to be circular reasoning in the statement of this result: we need to choose $\delta_{2,3}>0$
extremely small depending on the profiles $V_{1,2}^{ab}[0]$, $1\leq
b\leq B_{1}$, but here $B_{1}$ itself was defined based on
$\delta_{2}$. This is clarified by noting that all the profiles
$V_{1,2}^{ab}[0]$ are small (more precisely, the square sum of their
energies is small) for $b$ sufficiently large, and this implies that
enlarging $B_{1}$ past a certain cutoff will not affect the
condition on $\delta_{2}$; for more clarification see the
``important technical observation'' below.
\end{remark}

\begin{proof}(Proposition~\ref{BGIIHard}) We will prove the inequality for $P_{k}\epsilon$ using a bootstrap argument.
The challenge consists in careful book-keeping of all the possible
interactions. The idea is to essentially replicate the proof of
Proposition~\ref{PsiBootstrap} with $\epsilon=\epsilon_2$. The main
novel feature here is that we now have to deal with a large number
of additional source terms stemming from the nonlinear interactions
of the various constituents in the decomposition of
$\psi^{n(<a)}_\alpha$. To begin with, we split the (large) time
interval $[0, t^{nab_{2}}-C]$ into finitely many intervals
\[
[0, t^{nab_{2}}-C]=\cup_{j=1}^{M_{1}}I_{j},
\]
where we have a decomposition (with
$\Psi^{nA_{0}^{(0)}}_{\alpha}=\Phi_{\alpha}^{nA_{0}^{(0)}}e^{-i\sum_{k=1,2}\triangle^{-1}\partial_{k}\Phi_{k}^{1nA_{0}^{(0)}}}$)
\[
\Psi^{nA_{0}^{(0)}}|_{I_{j}}=\Psi^{nA_{0}^{(0)}}_{jL}+\Psi^{nA_{0}^{(0)}}_{jNL}
\]
with, see Corollary~\ref{cor:localsplit2},
\[
\|\Psi^{nA_{0}^{(0)}}_{jNL}\|_{S(I_{j}\times\R^{2})}<\eps_2,\quad \|\nabla_{x,t}\psi^{nA_{0}^{(0)}}_{jL}\|_{L_{t}^{\infty}\dot{H}^{-1}}\lesssim
\eps_2^{-\frac12}\, E^2_{\mathrm{crit}}
\]
We then run a bootstrap argument inductively
on each of these intervals, where of course $M_{1}=M_{1}(\Ecrit)$ is
not too large. We shall now work on the interval $I_{1}$, say. This
enables us to use the covariant energy estimate from Step~3.

Clearly,  the evolved concentration profiles also interact with
$\epsilon$; we then further subdivide the intervals $I_{j}$ into
smaller ones, which by abuse of notation we again label as $I_{j}$,
such that
\[
\big(\sum_{b\in\{1,2,\ldots
B_{1}\}}\tilde{\psi}^{nab}\big)|_{I_{j}}= \sum_{b\in\{1,2,\ldots
B_{1}\}}\tilde{\psi}^{nab}_{jL}+\sum_{b\in\{1,2,\ldots
B_{1}\}}\tilde{\psi}^{nab}_{jNL}
\]
Note that now the number of intervals is of the form
$M_{1}=M_{1}(\Ecrit,
\{V_{1,2}^{ab}[0]\}_{b\in\{1,2,\ldots,B_{1}\}})$. Furthermore, one
has
\[
\|\big(\sum_{b\in\{1,2,\ldots
B_{1}\}}\tilde{\psi}^{nab}\big)_{jNL}\|_{S(I_{j}\times\R^{2})}< \eps_2
\]
while also
\[
\|\big(\sum_{b\in\{1,2,\ldots
B_{1}\}}\nabla_{x,t}\tilde{\psi}^{nab}\big)_{jL}\|_{L_t^\infty\dot{H}^{-1}}\lesssim \eps_2^{-\frac{1}{2}}E^2_{\mathrm{crit}}
\]
where $\eps_2$ is a universal constant  depending only on~$\Ecrit$.
The fact that we get the last inequality with universal implied
constant hinges on the approximate orthogonality of the
$\tilde{\psi}^{nab}$ for $n$ large enough. One may object at this
point that the choice of $B_1$ was dictated by~$\delta_2$, and hence
may be extremely large, which in turn means that the number $M_1$ of
intervals above depends on $\delta_2$ and may also become extremely
large. The following observation, however, shows that $M_1$ only
depends on a {\em fixed number} of concentration profiles
independent of $\delta_2$:

\smallskip

{\bf{Important technical observation:}}

\smallskip

Here we note that $M_{1}$ really only depends on
$\{V_{1,2}^{ab}[0]\}_{b\in\{1,2,\ldots,B_{0}\}}$, for some $B_{0}$
with the property that
\[
\sum_{b\geq B_{0}}\|V_{1,2}^{ab}[0]\|_{L_{x}^{2}}^{2}<\epsilon_{0}
\]
where $\epsilon_{0}$ is the small-energy global well-posedness
cutoff. Thus we can make $\delta_{2}$ small without increasing
$M_{1}$ concurrently. To see this, write
\[
\{V_{1,2}^{ab}[0]\}_{b\in\{B_{0},B_{0}+1,\ldots,B_{1}\}}=\{V_{1,2}^{ab}[0]\}_{b\in\Lambda_{1}}\cup\{V_{1,2}^{ab}[0]\}_{b\in\Lambda_{2}},\quad
\Lambda_{1}\cup\Lambda_{2}=\{B_{0},B_{0}+1,\ldots,B_{1}\},
\]
so that $\{S_{\tilde{A}^{nab}}(V_{1,2}^{ab}[0])(0-t^{nab},
x-x^{nab})\}_{b\in\Lambda_{1}}$ is the collection of temporally
bounded concentration profiles with $b\in
\{B_{0},B_{0}+1,\ldots,B_{1}\}$. Then the argument that was used for
Case~(A) in the proof of
Proposition~\ref{ConcentrarionProfileApprox} reveals that we can
approximate
\[
\sum_{b\in\Lambda_{2}}V_{3,4}^{nab}(0-t^{nab}, x-x^{nab})
\]
up to a constant phase shift arbitrarily well by the Coulomb
components of an admissible map, and then Proposition~\ref{TwistedWaveEquation}
allows us to evolve the data
\[
\sum_{b\in\Lambda_{2}}V_{3,4}^{nab}(0-t^{nab},
x-x^{nab})+o_{L^{2}}(1)
\]
using the covariant linear flow on $[0, t^{nab_1}]$. This leads to
bounds that are uniform in $B_{0}, n$ only involving $\epsilon_{0}$.
Handling the contribution of
\[
\sum_{b\in\Lambda_{1}}V_{3,4}^{nab}(0-t^{nab},
x-x^{nab})+o_{L^{2}}(1)
\]
i.e., the ``tail'' of bounded concentration profiles, is  more
complicated since we may no longer necessarily approximate this sum
by Coulomb components of admissible maps, but only the individual
summands $V_{3,4}^{nab}(0-t^{nab}, x-x^{nab})$. Thus the correct
evolution of this term has to consist of the evolution of the
individual ingredients, and one then needs to bound the $S$-norm of
this (very large) sum in terms of an apriori bound, provided $n$ is
large enough. In this regard we have the following result.

\begin{lemma}\label{lem:boundedconcproftail}
For each $b\in\Lambda_1$ and $t\in[0, t^{nab_1}]$, denote by $V_{3,4}^{nab}(t-t^{nab},
x-x^{nab})+o_{L^{2}}(1)$ the (nonlinear) wave maps evolution of the Coulomb components
of an admissible sufficiently good approximation to the data $V_{3,4}^{nab}(t-t^{nab},
x-x^{nab})$, as in the preceding discussion. Then for $n$ large enough, we have
\[
\|\sum_{b\in\Lambda_{1}}V_{3,4}^{nab}(t-t^{nab},
x-x^{nab})+o_{L^{2}}(1)\|_{S[0, t^{nab_1}]}\lesssim \eps_0
\]
for a suitable universal implied constant.
\end{lemma}

\begin{proof}(Lemma~\ref{lem:boundedconcproftail}) For each $V_{3,4}^{nab}(t-t^{nab},
x-x^{nab})+o_{L^{2}}(1)$, pick an interval $[0, \tilde{t}^{ab}]$ with the property that we can write
\[
\big[V_{3,4}^{nab}(t-t^{nab},
x-x^{nab})+o_{L^{2}}(1)\big]|_{[0, \tilde{t}^{ab}]^{c}}=\big[V_{3,4}^{nab}(t-t^{nab},
x-x^{nab})+o_{L^{2}}(1)\big]_L+\big[V_{3,4}^{nab}(t-t^{nab},
x-x^{nab})+o_{L^{2}}(1)\big]_{NL}
\]
where we impose the condition
\[
\|\nabla_{x,t}\big[V_{3,4}^{nab}(t-t^{nab},
x-x^{nab})+o_{L^{2}}(1)\big]_L\|_{L_{t}^{\infty}\dot{H}_x^{-1}}\les \|V_{3,4}^{nab}(0-t^{nab},
x-x^{nab})\|_{L_{x}^{2}}
\]
\[
\|\big[V_{3,4}^{nab}(t-t^{nab},
x-x^{nab})+o_{L^{2}}(1)\big]_{NL}\|_{S([0,
\tilde{t}^{ab}]^{c}\times\R^{2})}\ll \frac{\eps_0}{B_{1}}
\]
where the implied constant in the first inequality is universal. That this is possible follows from
Corollary~\ref{cor:localsplit2} and Remark~\ref{rem:psiLimprov}.  Choosing $n$ large enough and exploiting
essential disjointness of the supports at time $t^{nab_1}$, we can arrange that
\[
\big(\sum_{b\in\Lambda_1}\|\nabla_{x,t}\big[V_{3,4}^{nab}(t^{nab_1}-t^{nab},
x-x^{nab})+o_{L^{2}}(1)\big]_L\|_{\dot{H}_x^{-1}}^{2}\big)^{\frac{1}{2}}\les \eps_0
\]
which then implies
\[
\|\sum_{b\in\Lambda_1}\big[V_{3,4}^{nab}(t-t^{nab},
x-x^{nab})+o_{L^{2}}(1)\big]_L\|_{S([0, \tilde{t}^{ab}]^{c}\times\R^{2})}\les \eps_0
\]
In order to complete the proof of the lemma, we need to also control
\[
\|\sum_{b\in\Lambda_1}\big[V_{3,4}^{nab}(t-t^{nab},
x-x^{nab})+o_{L^{2}}(1)\big]\|_{S([0, \tilde{t}^{ab}]\times\R^{2})}
\]
Here we exploit the fact that  for $n$ large, the functions  $\big[V_{3,4}^{nab}(t-t^{nab},
x-x^{nab})+o_{L^{2}}(1)\big]$ are supported on disjoint light cones up to small errors with respect to $S$.
One then concludes by means of a simple orthogonality type argument that for $n$ large enough
\[
\|\sum_{b\in\Lambda_1}\big[V_{3,4}^{nab}(t-t^{nab},
x-x^{nab})+o_{L^{2}}(1)\big]\|_{S([0, \tilde{t}^{ab}]\times\R^{2})}\les \eps_0
\]
where the implied constant is universal.
\end{proof}

Now assume the bound
\[
\|P_{k}\epsilon_{1}\|_{S}\leq
C_{5}\delta_2
\]
We show that provided we choose
$C_{5}=C_{5}(\Ecrit)$ large enough, we can bootstrap $C_{5}$ to
$\frac{C_{5}}{2}$, whence we get the bound on all of $I_{1}$. Then
we continue the argument to $I_{2}$ etc.  Note that by choosing
$\delta_{2}$ small enough in relation to $M_{1}$ as well as the
other apriori data $\Ecrit$, $V_{1,2}^{ab}[0]$, $b=1,2,\ldots,
B_{0}$, the error term will then remain small.
\\

By scaling invariance, it suffices to bootstrap the estimate for
$P_{0}\epsilon$. We now bootstrap the bound for $P_{0}\epsilon$.
Here we essentially proceed as in step (3), the apriori bound for
the first non-atomic component
$\psi^{nA_{0}^{(0)}}=w^{nA_{0}^{(0)}}e^{-i\sum_{k=1,2}\triangle^{-1}\partial_{k}w_{k}^{nA_{0}^{(0)}}}$.
Thus we distinguish as there between the small time case, when the
div-curl system suffices, and the large time case, when the wave
equations are important: we shall work here on the interval $I_{1}$
containing the initial time slice $t=0$.
\\

{\it{(i): small time case}} $|I_{1}|<T_{1}$. Here
$T_{1}$ is a sufficiently small absolute constant. Write the equation for
$\epsilon$, using the div-curl system, schematically as follows:
\begin{equation}\nonumber\begin{split}
\partial_{t}P_{0}\epsilon=&\nabla_{x}P_{0}\epsilon+P_{0}\big[\epsilon\nabla^{-1}([\psi^{nA_{0}^{(0)}}+\sum_{b=1}^{B_{0}}\tilde{\psi}^{nab}+\tilde{\psi}^{a(>B_{0})}+W^{naB}]^{2})\big]\\
&+P_{0}\big[[\psi^{nA_{0}^{(0)}}+\sum_{b=1}^{B_{0}}\tilde{\psi}^{nab}+\tilde{\psi}^{a(>B_{0})}+W^{naB}]\nabla^{-1}(\epsilon [\psi^{nA_{0}^{(0)}}+\sum_{b=1}^{B_{0}}\tilde{\psi}^{nab}+\tilde{\psi}^{a(>B_{0})}+W^{naB}])\big]\\
&+P_{0}\big[\epsilon\nabla^{-1}(\epsilon [\psi^{nA_{0}^{(0)}}+\sum_{b=1}^{B_{0}}\tilde{\psi}^{nab}+\tilde{\psi}^{na(>B_{0})}+W^{naB}])\big]+P_{0}\big[[\psi^{nA_{0}^{(0)}}+\sum_{b=1}^{B_{0}}\tilde{\psi}^{nab}+\tilde{\psi}^{a(>B_{0})}+W^{naB}]\nabla^{-1}(\epsilon^{2})\big]\\
&+P_{0}[\epsilon\nabla^{-1}(\epsilon^{2})]+\text{interactions terms}
\end{split}\end{equation}
 ``Interactions terms'' here refers to all possible expressions which do
not involve the radiation term~$\epsilon$ such as
\[
P_{0}[\psi^{nA_{0}^{(0)}}\nabla^{-1}[(\tilde{\psi}^{nab})^{2}]]
\]
Indeed, the complete list of the error terms included under this
heading is  complicated, due to our construction of the evolutions
$\tilde{\psi}^{nab}$ in (I)-(III) above. Recall that for the
{\it{temporally bounded type components}}, we use the nonlinear wave
maps flow on a large time interval $T^{ab\delta_{2}}$, but we then
use the covariant linear evolution past that time. This means that
on $[0, T^{ab\delta_{2}}]$, we generate error interaction terms like
the preceding one coming from the interactions with the low
frequency part $\psi^{nA_{0}^{(0)}}$, while on the interval
$[T^{ab\delta_{2}}, t^{nab_{2}}-C ]$ generate errors due to the
{\it{nonlinear self-interactions}} of~$\tilde{\psi}^{nab}$.

On the other hand, for the {\it{temporally unbounded type
components}}, we use the linear covariant evolution on $[0,
t^{nab_{2}}-C]$, which means that we generate errors due to the
nonlinear self-interactions.

In addition to all these, we generate errors due to different
concentration profiles interacting with each other, as well with the
small frequency component $w^{nA_{0}^{(0}}$, or the weakly small
error, and the latter also generates nonlinear errors due to
interactions with itself. We will deal with this rather large
collection of errors later, showing that we can make its $N[0]$-norm
arbitrarily small by choosing $B_1$ large enough, and then $n$ large
enough.

We also use the notation $\tilde{\psi}^{a(>B_{0})}$ for the
evolution of
\[
\sum_{b\in\Lambda_{1}\cup\Lambda_{2}}V_{3,4}^{nab}(0-t^{nab},
x-x^{nab})+o_{L^{2}}(1),
\sum_{b\in\Lambda_{1}\cup\Lambda_{2}}\partial_{t}S_{\tilde{A}^{nab}}(V_1^{ab}[0])(0-t^{nab},
x-x^{nab})+o_{L^{2}}(1),
\]
as explained in the ``important technical observation'' above.

We first deal with the terms involving $\epsilon$. Our task is to
gain a smallness constant that allows us to improve the apriori
bound we are assuming about $\epsilon$.

\medskip
 {\it{(i.1): Terms involving $\epsilon$}}. These can be
handled exactly as Case~1 in the proof of
Proposition~\ref{PsiBootstrap}, in light of the bound
\[
\|[\psi^{nA_{0}^{(0)}}+\sum_{b=1}^{B_{0}}\tilde{\psi}^{nab}+\tilde{\psi}^{a(>B_{0})}+W^{naB}]\|_{S}\leq C(\psi^{nA_{0}^{(0)}},
 \{\tilde{\psi}^{nab}\}_{b=1}^{B_0}, \Ecrit)
\]
Thus for example paralleling Case~1~(a) in the proof of
Proposition~\ref{PsiBootstrap}, one obtains a bound
\[
\sum_{k\in\Z}\|\chi_{I_j}P_{k}\big[\epsilon\nabla^{-1}([\psi^{nA_{0}^{(0)}}+\sum_{b=1}^{B_{0}}\tilde{\psi}^{nab}+\tilde{\psi}^{a(>B_{0})]^{2})}
\|_{L_t^2\dot{H}^{-\frac{1}{2}}}^{2}\ll \|\epsilon\|_{S}
\]
provided we choose the time cutoffs suitably (such that the number $M_1$ of such time intervals is as above).

\smallskip
{\it{(i.2) Errors due to nonlinear (self)interactions of the
$\tilde{\psi}^{nab}$, $\psi^{nA_{0}^{(0)}}$, $W^{naB}$.}} Note that
these errors serve as source terms for $\epsilon$, and hence we need
to show that they are extremely small (of order controlled by
$\delta_{2}$). The mechanism for this is first choosing $B_1$
sufficiently large (for the contributions involving $W^{naB}$), and
then choosing $n$ large enough.

\smallskip
{\it{(i.2.a) Errors generated by the temporally bounded type
$\tilde{\psi}^{nab}$.}} If $\tilde{\psi}^{nab}$ is the evolution of
a temporally bounded concentration profile, then recall that we let
$\tilde{\psi}^{nab}$ be the wave maps evolution on the interval $[0,
T^{ab\delta_{2}}]$, {\it{provided $\tilde{\psi}^{nab}$ is supported
at frequency scale $\sim 1$}}. Now we want to track the evolution of
an arbitrary frequency mode $P_{k}\epsilon$, which we have scaled to
$k=0$. But then  we have also re-scaled all the source terms. Now
the source terms generated by $\tilde{\psi}^{nab}$ itself come from
a number of sources:  first, the ``gluing definition'' of
\eqref{glue} implies that we generate errors of the form (before
frequency localization)
\[
\chi_{(-\infty, T^{ab\delta_{2}}+10)}'(t)\psi^{ab}-\chi_{(-\infty, T^{ab\delta_{2}}+10)}'(t)S_{A^n}\psi^{ab}[T^{ab\delta_{2}}]
\]
The only way for this term to contribute in the Case~(i) for a fixed
frequency (which we assume equals one after scaling) is when the
original frequency (which gets scaled to one) is extremely large.
But this contribution is then easily seen to be very small in
$L_t^\infty L_x^2$, say, due to the frequency localization of
$\tilde{\psi}^{ab}$. Next, the self-interaction errors generated
from the usual div-curl system are (schematically)
\[
P_{0}[\partial_{t}\tilde{\psi}^{nab}-\nabla_{x}\tilde{\psi}^{nab}-\tilde{\psi}^{nab}\nabla^{-1}[(\tilde{\psi}^{nab})^{2}]],
\]
which vanishes provided the $I_{1}$ fits into the {\it{re-scaled
interval}} $[0, T^{ab\delta_{2}}]$. Otherwise, one obtains a
contribution of the above form on the complement of  the
{\it{re-scaled interval}} $[0, T^{ab\delta_{2}}]$ inside $I_{1}$,
and which is of the above form. We need to show that picking $n$
large enough, this can be made arbitrarily small. For this purpose
we use the following observation.

\begin{lemma}\label{AsymptoticStructure} Let $\tilde{\psi}^{nab}$ be the evolved Coulomb components
of a temporally bounded type concentration profile, concentrated at
frequency $\sim 1$. Then letting $T^{ab\delta_{2}}$ be the time
indicating transition from nonlinear to linear evolution (as
explained in the preceding discussion), we have
\[
\tilde{\psi}^{nab}_{\alpha}(T^{ab\delta_{2}},
\cdot)=\partial_{\alpha}\tilde{\psi}^{nab}(T^{ab\delta_{2}},
\cdot)+\text{error},
\]
where
\[
\|\text{error}\|_{L_{x}^{2}}\to 0
\]
as $\delta_{2}\to 0$, and furthermore
\[
\tilde{\psi}^{nab}=\sum_{k=1,2}\triangle^{-1}\partial_{k}\tilde{\psi}^{nab}_{k}
\]
\end{lemma}

The proof of this lemma follows exactly as in the proof of Case~(B)
of Proposition~\ref{ConcentrarionProfileApprox}. It then follows
that in case we are on the complement of the re-scaled interval $[0,
T^{ab\delta_{2}}]$ inside $I_{1}$, we generate errors of the form
\[
P_{0}[\tilde{\psi}^{nab}\nabla^{-1}[(\tilde{\psi}^{nab})^{2}]]+\text{error}
\]
with $\text{error}$ as in the preceding lemma, in addition to errors
stemming from interactions of $\tilde{\psi}^{nab}$ with the other
components $\psi^{nA_{0}^{(0)}}$ etc.~to be considered later. But
then, using the $L_{t,x}^{\infty}$-dispersion for the
$\tilde{\psi}^{nab}(t,\cdot)$ as $|t|\to\infty$, it is seen that
\[
\|P_{0}[\tilde{\psi}^{nab}\nabla^{-1}[(\tilde{\psi}^{nab})^{2}]]\|_{L_{t}^{\infty}L_{x}^{2}(I_{1}\cap[0,
T^{ab\delta_{2}}]^{c}\times\R^{2})}\ll\delta_2
\]
if we choose $T^{ab\delta_2}$ large enough in relation
to~$\delta_2$. Next, we need to analyze the errors generated by
$\tilde{\psi}^{nab}$ through interaction with the other ingredients
$\psi^{nA_{0}^{(0)}}$, $\tilde{\psi}^{a(>B_{0})}$, and~$W^{naB}$. We
begin  with the interactions between two {\em distinct} terms
$\tilde{\psi}^{nab}$, $b=1,2,\ldots, B_1$. Thus we are considering
\[
P_0[\tilde{\psi}^{nab_1}\nabla^{-1}(\tilde{\psi}^{nab_2}\tilde{\psi}^{nab_3})]
\]
where $b_i\neq b_j$ for some $i, j$.  By the frequency localization
of all these factors, we may assume that, up to negligible errors,
each of them satisfies $\tilde{\psi}^{nab_j}=P_{[-C_6,
C_6]}\tilde{\psi}^{nab_j}$ where $C_6$ is a potentially extremely
large constant depending on the frequency localizations (i.e., how
well-localized the factors are in frequency space), as well as
$\delta_2$ and $B_1$, and that $\log\big[(\lambda_n^a)^{-1}\big]\in
[-C_6, C_6]$. Now assume first that $I_1\subset [0,
T^{ab_j\delta_2}]$ for all $j$, i.e., our time interval is such that
we are in the ``nonlinear regime'' for each of these factors. But
then choosing $n$ large enough, we can force
\[
\|P_0[\tilde{\psi}^{nab_1}\nabla^{-1}(\tilde{\psi}^{nab_2}\tilde{\psi}^{nab_3})]\|_{L_{t}^{\infty}L_{x}^{2}}
\ll \delta_2\frac{1}{C_{6}^{100}}
\]
by the essential disjointness of the supports of the factors, and
this suffices to handle Case~1, see the proof of
Proposition~\ref{PsiBootstrap}. Indeed, the expression
\[
[\tilde{\psi}^{nab_1}\nabla^{-1}(\tilde{\psi}^{nab_2}\tilde{\psi}^{nab_3})]
\]
 is essentially supported in a frequency interval $[-10C_6, 10C_6]$, and repeating the above estimate
 for each of these frequencies and square summing easily yields the bound
 \[
 \sum_{k}\|\chi_{[-T_1, T_1]}(2^{k}t)P_k[\tilde{\psi}^{nab_1}\nabla^{-1}(\tilde{\psi}^{nab_2}\tilde{\psi}^{nab_3})]
 \|_{L_t^{2}\dot{H}^{-\frac{1}{2}}}\ll \delta_2
 \]
 If, on the other hand, at least one of the factors $\tilde{\psi}^{nab_j}$ is in the ``linear regime'' (i.e.,
  satisfies the covariant wave equation), then smallness follows from the $L^\infty$-decay.

Next we consider the term
\[
P_{0}[\tilde{\psi}^{nab}\nabla^{-1}[(\psi^{nA_{0}^{(0)}})^{2}]]
\]
Here of course it is essential that we are in Case~(i) and so it
suffices to estimate the $L_{t}^{\infty}L_{x}^{2}$ or also
$L_{t,x}^{2}$ norm of this term, see Case~1 of the proof of
Proposition~\ref{PsiBootstrap}. Due to the essential disjointness of
the Fourier supports of $\tilde{\psi}^{nab}$ and
$\psi^{nA_{0}^{(0)}}$, see
Proposition~\ref{ControlNonatomicComponent1}, we may assume that the
first input $\tilde{\psi}^{nab}$ has frequency of size one,  while
the second input has extremely small frequency (controlled by
picking $n$ large enough). But then we may estimate this
contribution by placing $\nabla^{-1}[(\psi^{nA_{0}^{(0)}})^{2}]$
into $L_{t,x}^\infty$, and re-scaling and square-summing over the
output frequencies results in the desired small bound.

\noindent Finally, the interactions of temporally bounded
$\tilde{\psi}^{nab}$ with the remaining weakly small errors
$W^{naB_1}$ are handled similarly by exploiting the smallness of the
latter with respect to $L_{t,x}^{\infty}$. Here the ``Important
Technical  Observation'' from before becomes important again.
\\

{\it{(i.2.b) Errors generated by temporally unbounded
$\tilde{\psi}^{nab}$.}} Again the errors generated are of the form
\[
P_{0}[\partial_{t}\tilde{\psi}^{nab}-\nabla_{x}\tilde{\psi}^{nab}-\tilde{\psi}^{nab}\nabla^{-1}[(\tilde{\psi}^{nab})^{2}]],
\]
as well as terms involving interactions of $\tilde{\psi}^{nab}$ with
$\psi^{nA_{0}^{(0)}}$, $\tilde{\psi}^{a(>B_{0})}$, as well as
$W^{naB}$. From Part~(B) of the proof of
Proposition~\ref{ConcentrarionProfileApprox}, we know that
$\tilde{\psi}^{nab}$ is of gradient form up to an error which can be
made arbitrarily small. Hence the above simplifies, up to a
negligible error, to the nonlinear term
\[
-P_{0}[\tilde{\psi}^{nab}\nabla^{-1}[(\tilde{\psi}^{nab})^{2}]],
\]
To estimate this, we can first reduce this to
\[
-P_{0}[\tilde{\psi}^{nab}\nabla^{-1}P_{[-C_{6},C_{6}]}[(P_{[-C_{6},
C_{6}]}\tilde{\psi}^{nab})^{2}]],
\]
arguing as in Case~1 of the proof of Proposition~\ref{PsiBootstrap},
and then by using the $L_{t,x}^{\infty}$-dispersion, i.e.,
Lemma~\ref{prop:covdisp}, to write
\[
P_{[-C_{6}, C_{6}]}\tilde{\psi}^{nab}(0,\cdot)=o_{L^{\infty}}(1)
\]
from which the desired smallness follows easily. The interaction
terms of temporally unbounded $\tilde{\psi}^{nab}$ with the
remaining components $\psi^{nA_{0}^{(0)}}$,
$\tilde{\psi}^{a(>B_{0})}$, $W^{naB}$, are handled as before and are
omitted.
\\

{\it{(i.2.c) Errors generated by the weakly small remainder
$W^{naB}$.}} Again recalling Part~(B) of the proof of
Proposition~\ref{ConcentrarionProfileApprox}, and
Lemma~\ref{BGIIErrorControl}, we know that $W^{naB}_{\alpha}$ is of
pure gradient form up to a negligible error (provided $B$ and $n$
are large enough). The conclusion is that the error of the form
\[
P_{0}[\partial_{t}W^{naB}-\nabla_{x}W^{naB}-W^{naB}\nabla^{-1}([W^{naB}]^{2})
\]
reduces up to a negligible error to the nonlinear self-interaction
term
\[
P_{0}[-W^{naB}\nabla^{-1}([W^{naB}]^{2})],
\]
which can be estimated as in the preceding case, using the smallness
of $\|W^{naB}\|_{L_{t,x}^{\infty}}$.
\\

(ii) $|I_{1}|>T_{1}$, $T_{1}$ as in Case 1. Proceeding as in Case~2
of the proof of Proposition~\ref{PsiBootstrap}, we decompose
$P_{0}\epsilon$ into
\[
P_{0}\epsilon=P_{0}Q_{<D}\epsilon+P_{0}Q_{\geq D}\epsilon,
\]
where $D=D(\Ecrit)$ is a sufficiently large constant. Then arguing as
in the proof of Proposition~\ref{PsiBootstrap},  we obtain two equations
\[
\Box_A P_{0}Q_{<D}\epsilon_\alpha =F^1_\alpha
\]
\[
\Box P_{0}Q_{\geq D}\epsilon=F^2_\alpha
\]
Here the magnetic potential $A$ in the first equation  is defined as
in the proof of Proposition~\ref{PsiBootstrap} but with $\psi_L$
replaced by
\[
\psi^{nA_{0}^{(0)}}+\sum_{b=1}^{B_{0}}\tilde{\psi}_L^{nab}+W^{naB}
\]
The source terms $F^1_\alpha$ are obtained as in
Section~\ref{sec:hodge}, and here we of course linearize around the
above expression. Then we re-iterate the estimates in the proof of
Proposition~\ref{PsiBootstrap}, with $\epsilon$ replacing
$\epsilon_2$ and $\epsilon_1=0$. As in Case~1 above, the only new
feature are the source terms coming from nonlinear interactions
between the various $\tilde{\psi}^{nab}$, $W^{naB_1}$. Fortunately,
the fact that each of these functions is essentially frequency
localized to the same interval, the mechanisms that force smallness
reduce as before to either physical separation or dispersive decay.
We explain here how to obtain smallness for the trilinear null-form
source terms, which we write schematically in the form
$\nabla_{x,t}[\rho_1\nabla^{-1}\calN_{\nu j}(\rho_2, \rho_3)]$, were
$\rho$ represents one of the functions $\psi^{nA_{0}^{(0)}}$,
$\sum_{b=1}^{B_{0}}\tilde{\psi}_L^{nab}$, $W^{naB}$. We consider the
following cases:

\smallskip

(ii.0) {\it{Errors due to the gluing construction \eqref{glue}.}}
These errors are of the form
\[
\chi_{(-\infty, T^{ab\delta_{2}}+10)}''(t)\psi^{ab}-\chi_{(-\infty, T^{ab\delta_{2}}+10)}''(t)S_{A^n}\psi^{ab}[T^{ab\delta_{2}}]
\]
\[
\chi_{(-\infty, T^{ab\delta_{2}}+10)}'(t)\partial_t\psi^{ab}-\chi_{(-\infty, T^{ab\delta_{2}}+10)}'(t)\partial_t S_{A^n}\psi^{ab}[T^{ab\delta_{2}}]
\]
To show the smallness of these, note that $\psi^{ab}$ solves the schematic div-curl system
\[
\nabla_t\psi-\nabla_x\psi=\psi\nabla^{-1}(\psi^2)
\]
Now since we have $\psi^{ab}=o_{L^{\infty}}(1)+o_{L^2}(1)$, choosing
$T^{ab\delta_2}$ large enough, we see that (with $o(1)$ in case
$T\to\infty$)
\[
\chi_{[T, T+10]}\big[\nabla_t\psi^{ab}-\nabla_x\psi^{ab}\big]=o_{L_t^M\dot{H}^{-(1-\frac{1}{M})}}(1)
\]
Similarly, by construction, the extensions $S_{A^n}\psi^{ab}[T^{ab\delta_{2}}]$ also satisfy the (schematic) relations
\[
\chi_{[T, T+10]}\big[\partial_t\big(S_{A^n}\psi^{ab}[T^{ab\delta_{2}}])-\nabla_x\big(S_{A^n}\psi^{ab}[T^{ab\delta_{2}}])\big]=o_{L_T^\infty L_x^2}(1),
\]
see Lemma~\ref{AsymptoticStructure}. But then it easily follows that
\[
\|\chi_{(-\infty, T^{ab\delta_{2}}+10)}''(t)\psi^{ab}-\chi_{(-\infty, T^{ab\delta_{2}}+10)}''(t)S_{A^n}\psi^{ab}[T^{ab\delta_{2}}]\|_{N}\ll \delta_2
\]
\[
\|\chi_{(-\infty, T^{ab\delta_{2}}+10)}'(t)\partial_t\psi^{ab}-\chi_{(-\infty, T^{ab\delta_{2}}+10)}'(t)\partial_t S_{A^n}\psi^{ab}[T^{ab\delta_{2}}]\|_{N}\ll \delta_2
\]

(ii.1) {\it{Self-interactions of temporally bounded
$\tilde{\psi}^{nab}$}}. These only occur provided $I_1\cap [0,
T^{ab\delta_2}]^c\neq \emptyset$.  Thus assume the latter is the
case, and consider
\[
\nabla_{x,t}[\tilde{\psi}^{nab}\nabla^{-1}\calN_{\nu j}(\tilde{\psi}^{nab}, \tilde{\psi}^{nab})]
\]
Now the estimates of Section~\ref{subsec:improvetrilin} imply that
we obtain
\[
\|\nabla_{x,t}[\tilde{\psi}^{nab}\nabla^{-1}\calN_{\nu j}(\tilde{\psi}^{nab}, \tilde{\psi}^{nab})]\|_{N(I_1\times\R^2)}\ll \delta_2
\]
provided at least two of the inputs have Fourier support with very close angular alignment, depending on
$\|\psi^{nab}\|_{S}$, $\delta_2$. Thus we may assume that these inputs have Fourier supports with some amount
(albeit very small) of angular separation. Similarly, localizing the Fourier support to frequency $\sim 1$, say, we may reduce to the expression
\[
\sum_{k_{1,2,3}=O(1)}\nabla_{x,t}P_0[P_{k_1}\tilde{\psi}^{nab}\nabla^{-1}\calN_{\nu j}(P_{k_2}\tilde{\psi}^{nab},P_{k_3} \tilde{\psi}^{nab})],
\]
where the implied constant $O(1)$ is of course potentially extremely
large, depending on $\|\psi^{nab}\|_{S}$, $\delta_2$. We may
similarly assume that all the modulations present are of size $O(1)$
at most (which may again be quite large, depending on
$\|\psi^{nab}\|_{S}$, $\delta_2$). But then the assumed angular
separation between all factors allows us to bound this expression
(for fixed frequencies) by
\[
\|\nabla_{x,t}P_0[P_{k_1}\tilde{\psi}^{nab}\nabla^{-1}\calN_{\nu j}(P_{k_2}\tilde{\psi}^{nab},P_{k_3} \tilde{\psi}^{nab})]
\|_{N[0]}\lesssim \|P_{k_1}\tilde{\psi}^{nab}\|_{S[k_1]}\|\nabla^{-1}\calN_{\nu j}(P_{k_2}\tilde{\psi}^{nab},P_{k_3} \tilde{\psi}^{nab})\|_{L_{t,x}^{2}}
\]
But then the desired smallness follows by interpolating the improved bilinear Strichartz type bound
\[
\|\nabla^{-1}\calN_{\nu j}(P_{k_2}\tilde{\psi}^{nab},P_{k_3} \tilde{\psi}^{nab})\|_{L_{t,x}^{p}}\lesssim \prod_{j=1,2}\|P_{k_j}\tilde{\psi}^{nab}\|_{S[k_j]}
\]
for some $p<2$ following from a result due to Bourgain\footnote{Of
course one also has the endpoint result due to Wolff and Tao, but
this is not really needed here.}~\cite{B} as well as
Lemma~\ref{lem:Null_rep}, and the smallness bound
\[
\|\nabla^{-1}\calN_{\nu j}(P_{k_2}\tilde{\psi}^{nab},P_{k_3} \tilde{\psi}^{nab})\|_{L_{t,x}^{\infty}}\ll 1
\]
which we obtain by letting $T^{ab\delta_2}$ be large enough in
relation to $\delta_2$. Replacing $0$ by $k$ and square summing over
the output frequencies, the desired bound follows easily.
\\

(ii.2): {\it{ Interactions of two different temporally bounded
$\tilde{\psi}^{nab}$}}. Here the mechanism at work is the physical
separation of the centers of mass for $n$ large. Thus consider
\[
\nabla_{x,t}P_0[P_{k_1}\tilde{\psi}^{nab_1}\nabla^{-1}\calN_{\nu j}(P_{k_2}\tilde{\psi}^{nab_2},P_{k_3} \tilde{\psi}^{nab_3})]
\]
where we have $b_i\neq b_j$ for at least one pair $i, j$. Now if we have $I_1\subset [0, T^{ab_j\delta_2}]$ for both $i, j$,
then $\tilde{\psi}^{nab_{i, j}}$ are essentially supported in disjoint light cones for $n$ large enough. Specifically, due to
Lemma~\ref{lem:DataLocalization'} as well as Lemma~\ref{lem:chiS}, given $\delta_2>0$, we may write
\[
\tilde{\psi}^{nab_{i}}=\tilde{\psi}^{nab_{i}}_{\text{cone}}+\tilde{\psi}^{nab_{i}}_{\text{cone}^{c}}
\]
\[
\tilde{\psi}^{nab_{j}}=\tilde{\psi}^{nab_{j}}_{\text{cone}}+\tilde{\psi}^{nab_{j}}_{\text{cone}^{c}}
\]
where we have
\[
\|\tilde{\psi}^{nab_{i,j}}_{\text{cone}^{c}}\|_{S}\ll \delta_2
\]
while the functions $\tilde{\psi}^{nab_{i,j}}_{\text{cone}}$ are supported in disjoint double cones while still satisfying
\[
\|\tilde{\psi}^{nab_{i,j}}_{\text{cone}}\|_{S}\lesssim C(\tilde{\psi}^{nab})
\]
It is then straightforward to conclude that by choosing $n$ large enough, we may force
\[
\|\nabla_{x,t}[\tilde{\psi}^{nab_1}\nabla^{-1}\calN_{\nu j}(\tilde{\psi}^{nab_2}, \tilde{\psi}^{nab_3})]\|_{N}\ll\delta_2
\]

(ii.3) {\it{Interactions of temporally bounded $\tilde{\psi}^{nab}$
and $\psi^{nA_0^{(0)}}$}}.  We distinguish between $I_1\subset [0,
T^{ab\delta_2}]$ and $I_{1}\cap [0, T^{ab\delta_2}]^{c}\neq
\emptyset$. In the former case, where $\tilde{\psi}^{nab}$ is given
by the actual wave map propagation, we generate error terms of the
form
\[
\nabla_{x,t}[\tilde{\psi}^{nab}\nabla^{-1}\calN_{\nu j}(\psi^{nA_0^{(0)}}, \psi^{nA_0^{(0)}})]
\]
As in case (i.2.a) above, localizing the output to frequency $\sim
1$, we may reduce to the case when $\nabla^{-1}\calN_{\nu
j}(\psi^{nA_0^{(0)}}, \psi^{nA_0^{(0)}})$ has extremely small
frequency. But then one obtains
\[
\|\nabla_{x,t}P_0[\tilde{\psi}^{nab}\nabla^{-1}\calN_{\nu j}(\psi^{nA_0^{(0)}}, \psi^{nA_0^{(0)}})]\|_{L_t^1\dot{H}^{-1}}
\ll \|P_{[-5,5]}[\tilde{\psi}^{nab}\|_{L_t^\infty L_x^2}
\]
and re-scaling and square summing over the output frequencies, we
can force an upper bound $\ll\delta_2$ by choosing $n$ large enough.
We further generate interaction terms of the form
\[
\nabla_{x,t}[\psi^{nA_0^{(0)}}\nabla^{-1}\calN_{\nu
j}(\tilde{\psi}^{nab}, \psi^{nA_0^{(0)}})],\quad
\nabla_{x,t}[\psi^{nA_0^{(0)}}\nabla^{-1}\calN_{\nu
j}(\tilde{\psi}^{nab}, \tilde{\psi}^{nab})],\,
\]
However, the trilinear estimates in Section~\ref{sec:trilin} in
addition to the frequency support properties of these inputs reveal
that choosing $n$ large enough, we can force
\[
\|\nabla_{x,t}[\psi^{nA_0^{(0)}}\nabla^{-1}\calN_{\nu
j}(\tilde{\psi}^{nab}, \psi^{nA_0^{(0)}})]\|_{N}\ll \delta_2,\quad
\|\nabla_{x,t}[\psi^{nA_0^{(0)}}\nabla^{-1}\calN_{\nu
j}(\tilde{\psi}^{nab}, \tilde{\psi}^{nab})]\|_{N}\ll \delta_2
\]
Note that in the second situation above, i.e., $I_{1}\cap [0,
T^{ab\delta_2}]^{c}\neq \emptyset$, we no longer generate errors of
the form
\[
\nabla_{x,t}[\tilde{\psi}^{nab}\nabla^{-1}\calN_{\nu j}(\psi^{nA_0^{(0)}}, \psi^{nA_0^{(0)}})],
\]
within $I_{1}\cap [0, T^{ab\delta_2}]^{c}$ since now $\tilde{\psi}^{nab}$ is given by the linear covariant evolution.
\\

(ii.4) {\it{The remaining interactions the $\tilde{\psi}^{nab}$ of
bounded or unbounded type, $\psi^{nA_0^{(0)}}$, as well as
$\partial_{\alpha}W^{naB_1}$.}} These offer nothing new: note that
both the components $\tilde{\psi}^{nab}$ of unbounded type as well
as the covariant linear waves $W^{naB_1}$ have extremely small
$L_{t,x}^{\infty}$-norm, but enjoy the same frequency localization
properties as $\tilde{\psi}^{nab}$; indeed, for unbounded type
$tilde{\psi}^{nab}$, this follows by choosing the $C$ in the
interval we work on $[0, t^{nab_2}-C]$ sufficiently large. Thus any
trilinear interactions involving them can be handled as in case
(ii.1) in the asymptotic regime. Also, not that interactions of
 $\partial_{\alpha}W^{naB_1}$ with $\psi^{nA_0^{(0)}}$ are of schematic type
\[
\nabla_{x,t}[\psi^{nA_0^{(0)}}\nabla^{-1}\calN_{\nu j}(\partial_{\alpha}W^{naB_1}, \psi^{nA_0^{(0)}})]
\]
\[
\nabla_{x,t}[\psi^{nA_0^{(0)}}\nabla^{-1}\calN_{\nu j}(\partial_{\alpha}W^{naB_1}, \partial_{\alpha}W^{naB_1})]
\]
\[
\nabla_{x,t}[\partial_{\alpha}W^{naB_1}\nabla^{-1}\calN_{\nu j}(\psi^{nA_0^{(0)}}, \partial_{\alpha}W^{naB_1})],
\]
and hence can be made arbitrarily small with respect to $\|\cdot\|_N$ by choosing $n$ large enough.
\\

We omit the treatment of the higher order interactions between the $\tilde{\psi}^{nab}$  as this offers nothing
qualitatively new. Applying the arguments from the proof of Proposition~\ref{PsiBootstrap}, we now conclude the
proof of Proposition~\ref{BGIIHard}.
\end{proof}

Proposition~\ref{BGIIHard} allows us to extend the Coulomb
components $\psi^{n(<a)}_\alpha$ to the interval $[0, t^{nab_2}-C]$.
But now the profiles $\tilde{\psi}^{nab}$ which were temporally
bounded {\it{with respect to $t=0$}} become temporally unbounded
with respect to the new starting time $t^{nab_2}-C$ as $n\to\infty$.
Now by repeating the arguments in Section~\ref{subsec:geomprof}, we
see that for those concentration profiles $\tilde{\psi}^{nab}$ for
which (see the discussion in Section~\ref{subsec:geomprof})
$\limsup_{n\to\infty}|t^{nab_2}-t^{nab}|<\infty$, i.e., they
concentrate at time $t^{nab_2}$ or alternatively time $t^{nab_2}-C$
, the exact same arguments as in that subsection imply that they can
be approximated arbitrarily well in the $L^{2}$-sense by Coulomb
components of admissible maps (but for this we have to know that the
Components $\psi^{n(<a)}_{\alpha}$ and the associated wave maps
actually extend to time $t^{nab_2}-C$). But then we have an exact
analogue for Proposition~\ref{BGIIHard} on the interval
$[t^{nab_2}-C, t^{nab_3}-\tilde{C}]$. Repeating this process
finitely many times, we extend $\psi^{n(<a)}_\alpha$ to $\R^{2+1}$,
and obtain an apriori bound
\[
 \|\psi^{n(<a)}_{\alpha}\|_{S}<C_a
 \]
 as well as exponential decay of the $ \|P_k\psi^{n(<a)}_{\alpha}\|_{S[k]}$ for $k\gg \log[(\lambda_{n}^a)^{-1}]$.

\subsection{Completion of the proof of  Proposition~\ref{prop:waveops} as well as of Corollary~\ref{temporallyunbounded}}\label{subsec: scatteringwellposed}

Both of these can be deduced by a simpler version of the proof of Proposition~\ref{BGIIHard}. For Proposition~\ref{prop:waveops}, one makes the ansatz
\[
\psi_\alpha=\partial_\alpha \big(S(0-t_0)(\partial_t V, V)\big)+\epsilon_\alpha
\]
and performs a bootstrap argument for
$\|\epsilon_\alpha\|_{S((-\infty, 0]\times\R^2)}$ for $t_0$ large
enough. This is as in the proof of Proposition~\ref{BGIIHard} where
the free linear evolution of $\partial_\alpha
\big(S(0-t_0)(\partial_t V, V)\big)$ replaces one of the temporally
unbounded $\tilde{\phi}^{nab}$, say, while all the other components
$\tilde{\phi}^{nab}$, $\psi^{nA_0^{(0)}}$, $\partial_\alpha
W^{naB_1}$ vanish. If we pick $t_0$ large enough, all the error
terms due to nonlinear self-interactions of $\partial_\alpha
\big(S(t-t_0)(\partial_t V, V)\big)$ become arbitrarily small due to
the reasoning in case (ii.1) of the proof of
Proposition~\ref{BGIIHard}. As there, one then obtains the estimates
for $\epsilon$ via the technique used in the proof of
Proposition~\ref{PsiBootstrap}. We conclude that for given
$\delta_3>0$, if $t_0$ is chosen large enough, we obtain the apriori
bound
\[
\|\epsilon_\alpha\|_{S((-\infty, 0]\times\R^2)}\ll \delta_3
\]
and from here the smoothness of the solution follows, see Proposition~\ref{BlowupCriterion1}.
\\

Next, we prove Corollary~\ref{temporallyunbounded}: from  Proposition~\ref{prop:waveops}, we know that we can
construct admissible Coulomb components of the form
\[
\psi^n_\alpha:= \partial_\alpha \big(S(t-t_n)(\partial_t V, V)\big)+\epsilon_\alpha
\]
for $t\in (-\infty, t_n-C]$ for some large enough absolute constant $C$, with
\[
\limsup_{n\to\infty}\|\epsilon_\alpha\|_{S((-\infty, t_n-C]\times\R^2)}\ll 1.
\]

Now we claim that the functions $\psi_\alpha^n(t_n-10C,\cdot)$ form a Cauchy sequence in the $L_x^2$-sense. To see this, note that for $n>m$
\[
\psi_\alpha^n(t_n-t_m, \cdot)=\psi_\alpha^{m}(0,\cdot)+o_{L^2}(1)
\]
as $n, m\to\infty$, whence by Proposition~\ref{prop:ener_stable} one
has
\[
\psi_\alpha^n(t_n-10C, \cdot)=\psi_\alpha^m(t_m-10C, \cdot)+o_{L^2}(1)
\]
But then also
\[
\psi_\alpha^n(t+t_n, \cdot)=\psi_\alpha^m(t+t_m,
\cdot)+o_{L^2}(1),\quad t\in(-\infty, -10C)
\]
again by Proposition~\ref{prop:ener_stable} , and furthermore, due to the uniform bounds
\[
\limsup_{n\to\infty}\|\psi^n_\alpha\|_{S((-\infty, t_n-C]\times\R^2)}<M<\infty
\]
for suitable $M\in\R$, we conclude upon denoting
\[
\Psi_\alpha^\infty(t, \cdot):=\lim_n\psi_\alpha^n(t+t_n, \cdot)
\]
that
\[
\|\Psi_\alpha^\infty\|_{S((-\infty, \tilde{C}]\times\R^2)}\leq M
\]
for any $\tilde{C}<10C$, as desired.

\subsection{Step 5 of the Bahouri Gerard process; adding all atoms}\label{subsec:BGEnd}

In the preceding subsection we  derived apriori bounds for the wave
maps evolution of the (admissible) Coulomb components
\[
 (w^{nA_0^{(0)}}_\alpha+\phi^{n1}_\alpha)e^{-i\sum_{k=1,2}\triangle^{-1}\partial_k(w^{nA_0^{(0)}}_k+\phi^{n1}_k)} + o_{L^2}(1)
\]
under the assumption that either
\[
 \liminf_{n\to\infty}\|w^{nA_0^{(0)}}\|_{L_{x}^{2}}>0
\]
or else, applying the second stage Bahouri Gerard decomposition to the large atom $\phi^{n1}$, that at all the
concentration profiles have energy $<\Ecrit$. We shall henceforth make this assumption. Now we continue the
process by extending the data at time $t=0$ for the Coulomb components to
\[
 (w^{nA_0^{(0)}}_\alpha+\phi^{n1}_\alpha+w^{nA_0^{(1)}}_\alpha)e^{-i\sum_{k=1,2}\triangle^{-1}\partial_k(w^{nA_0^{(0)}}_k
 +\phi^{n1}_k+w^{nA_0^{(1)}}_k)}+o_{L^{2}}(1),
\]
where we recall that the error term $o_{L^{2}}(1)$ is necessary in order to ensure that the data correspond to exact
Coulomb components of an admissible map. Denote the wave maps evolution of
\[
 (w^{nA_0^{(0)}}_\alpha+\phi^{n1}_\alpha+w^{nA_0^{(1)}}_\alpha)e^{-i\sum_{k=1,2}\triangle^{-1}\partial_k(w^{nA_0^{(0)}}_k
 +\phi^{n1}_k+w^{nA_0^{(1)}}_k)}+o_{L^{2}}(1),
\]
which is defined at time $t=0$, by the same symbol.
We state the result:
\begin{prop}
Under the preceding assumptions, the evolution of the preceding Coulomb components exists globally in time. For $n$ large enough, we have an apriori bound
\[
 \|(w^{nA_0^{(0)}}_\alpha+\phi^{n1}_\alpha+w^{nA_0^{(1)}}_\alpha)e^{-i\sum_{k=1,2}\triangle^{-1}
 \partial_k(w^{nA_0^{(0)}}_k+\phi^{n1}_k+w^{nA_0^{(1)}}_k)}+o_{L^{2}}(1)\|_{S(\R^{2+1})}<\infty
\]
The bound here depends on $\Ecrit$ as well as the apriori bounds for
the evolution of the concentration profiles extracted by adding
$\phi^{n1}$. Furthermore, we have the same bounds as in
Proposition~\ref{PsiInduction} (applied to the union of all
$J_{j}$), where the implied constants depend on $\Ecrit $ as well as
the apriori bounds for the evolution of the concentration profiles
extracted by adding $\phi^{n1}$.
\end{prop}

The proof of this is a precise replica of the one given in Step~3. The difference consists in the fact that in the decomposition (see Step~2)
\[
 w^{nA_0^{(1)}}=\sum_{j}\phi^{na_j^{k}}+w^{nA^{(1)}}
\]
we now need to ensure that $\|w^{nA^{(1)}}\|_{\dot{B}^{0}_{2,
\infty}}$ is small enough depending on both $\Ecrit$ as well as the
apriori bounds for the concentration profiles from Step~4.

 Next, one
extends the data at time $t=0$ to
\[
(w^{nA_0^{(0)}}_\alpha+\phi^{n1}_\alpha+w^{nA_0^{(1)}}_\alpha+\phi^{n2}_\alpha)e^{-i\sum_{k=1,2}\triangle^{-1}\partial_k(w^{nA_0^{(0)}}_k
+\phi^{n1}_k+w^{nA_0^{(1)}}_k+\phi^{n2}_k)}+o_{L^{2}}(1)
\]
Repeating the procedure of Step~4 but with magnetic potential defined in terms of the $\psi$-evolution of
\[
 (w^{nA_0^{(0)}}_\alpha+\phi^{n1}_\alpha+w^{nA_0^{(1)}}_\alpha)e^{-i\sum_{k=1,2}\triangle^{-1}\partial_k(w^{nA_0^{(0)}}_k+\phi^{n1}_k
 +w^{nA_0^{(1)}}_k)}+o_{L^{2}}(1),
\]
one again derives the same types of bounds as in Proposition~\ref{BGIIHard} and the process continues $A_{0}$
many times, as we recall from the discussion at the beginning of Step~2. We have finally arrived at the following grand conclusion to this section.

\begin{theorem}
\label{thm:BGcore}
 Let $\psi^n$ be a sequence of gauged derivative components of admissible wave
 maps $\bfu^n:[-T_0^n, T_1^n]\times\R^2\to\Hyp^2$. The hypothesis
\begin{equation}
 \label{eq:inf_hypo}
 \lim_{n\to\infty}\|\psi^{n}\|_{S([-T_0^n, T_1^n]\times\R^2)}=\infty
\end{equation}
implies that two possible cases occur: up to rescaling and spatial translations, either we have
\[
 \psi^{n}_\alpha(0,\,\cdot)=V_\alpha+o_L^{2}(1)
\]
for some fixed $L^{2}$-profile $V_\alpha$, or else we have for some sequence $t^n\to\infty$ (or $t^n\to-\infty$) and suitable
$(\partial_t V, V)\in L^{2}\times\dot{H}^{1}$,
\[
 \psi^{n}_\alpha(0,\,\cdot)=\partial_\alpha\big(S(0-t^{n})[\partial_t V, V]\big)+o_L^{2}(1),\,
\]
where $S(t)$ refers to the standard free wave propagator. In the former case
\[
\sum_{\alpha=0}^{2}\|V_\alpha\|_{L^{2}}^{2}=\Ecrit,
\]
while in the latter case, one has
\begin{equation}\label{eq:V2}
\sum_{\alpha=0}^{2}\|\partial_\alpha V\|_{L^{2}}^{2}=\Ecrit
\end{equation}
\end{theorem}

Note that due to Lemma~\ref{BasicStability}, in the first case, there exist $T_{0}>0$, $T_1>0$ with the property that
\[
\sup_n\|\psi^{n}\|_{S([-T_0, T_1]\times\R^{2})}< \infty,
\]
and we can then define
\begin{equation}
 \lim_{n\to\infty}\psi_\alpha^{n}(t,x)=:\Psi^\infty_\alpha(t,x)
\label{eq:aa}
\end{equation}
where the limit is in the sense of $L^\infty_{\mathrm{loc}}([-T_0,
T_1];L^2(\R^{2}))$. Similarly, in the second case, due to
Corollary~\ref{temporallyunbounded}, we have the corresponding
statements on some semi-infinite interval $I=(-\infty, T_0)$
respectively $(T_0, \infty)$. We call the maximal such open interval
$(-T_0, T_1)$ (respectively $(-\infty, T_0)$ or $(T_0, \infty)$) the
{\em lifespan} of the asymptotic object $\Psi^\infty_\alpha(t,x)$.
Finally, in order to apply the Kenig-Merle type argument, we need
the following essential {\em compactness property}:

\begin{cor}\label{cor:compactV}
There exist continuous functions $\bar{x}:I\to\R^2$ and
$\lambda\::\:I\to \R^+$ so that the family of functions
$\{\lambda(t)^{-1}\,
\Psi^\infty_\alpha(t,(\cdot-\bar{x}(t))\lambda(t)^{-1})\}_{t\in
I}\subset L_{x}^{2}$ is pre-compact.
\end{cor}
\begin{proof} We may assume that
\begin{equation}
 \sup_{0<T_2<T_1} \|\Psi_\alpha^\infty\|_{S([0,T_2)\times\R^2)} = \infty
\label{eq:T1blow}
\end{equation}
see Lemma~\ref{lem:lifespan2}.
The proof follows~\cite{KeM1}, \cite{KeM2} and amounts to an argument by contradiction.
More precisely, we begin by showing that one can find functions $\lambda(t)$, $\bar{x}(t)$ not necessarily continuous with
the desired compactness property. Suppose this fails. Then there exists $\eps>0$ and a sequence of times $\{t_n\}\subset I$
so that
\begin{equation}
 \label{eq:compfail}
\inf_{\lambda>0,\bar{x}\in \R^2} \| \lambda^{-1}\, \Psi^\infty_\alpha(t_n,(\cdot-\bar{x})\lambda^{-1}) -  \Psi^\infty_\alpha(t_m,\cdot) \|_2\ge \eps
\end{equation}
for any $n\ne m$. Necessarily $t_n\rightarrow T_1$.
 Now apply Theorem~\ref{thm:BGcore} to the sequence $\{\Psi_\alpha^\infty(t_n,\cdot)\}_{n=1}^\infty$, which satisfies \eqref{eq:inf_hypo},
but on a shifted time-interval. Note that $\Psi_\alpha^\infty(t_n,\cdot)$ are not admissible in the sense that they are not necessarily
given as the Coulomb derivative components of admissible wave maps. However, by approximation by the original sequence $\psi_\alpha^n$ (up
to symmetries)
one concludes that either for some $V_\alpha\in L^2(\R^2)$,
\begin{equation}
 \Psi^\infty_\alpha(t_n,x)=\lambda_n^{-1} V_\alpha((x-x_n)\lambda_n^{-1})+o_L^{2}(1)
\label{eq:Veins}
\end{equation}
for some sequence $\lambda_n, x_n$, or that for some $s^n\to\infty$ or $s_n\to -\infty$,
\begin{equation}
 \Psi^{\infty}_\alpha(t_n,x)=\lambda_n^{-1} \partial_\alpha\big(S(-s^{n})[\partial_t V, V]\big)((x-x_n)\lambda_n^{-1})+o_L^{2}(1),
\label{eq:Vzwei}
\end{equation}
where $V$ is as in~\eqref{eq:V2}. Clearly, \eqref{eq:Veins} contradicts \eqref{eq:compfail}. For~\eqref{eq:Vzwei},
we first show that $\{s_n\}_{n=1}^\infty$ has to be bounded.
Assume that $s_n \to-\infty$. Then Proposition~\ref{prop:waveops} implies for large~$n$ that $\Psi^\infty_\alpha$ exists on $[0,\infty)\times\R^2$
and \[  \|\Psi^\infty_\alpha\|_{S([0,\infty)\times\R^2)}<\infty\] which contradicts our assumption~\eqref{eq:T1blow}. If on the other hand
$s_n \to \infty$, then this implies by the same proposition that
\[
 \sup_n  \|\Psi^\infty_\alpha\|_{S((-\infty,t_n])\times\R^2)}<\infty
\]
This again contradicts our assumption \eqref{eq:T1blow} and we are done. As in~\cite{KeM1} one proves by approximation  that $\lambda$ and $\bar{x}$
can be taken to be continuous.
\end{proof}

\section{The proof of the main theorem}\label{sec:conclusion}

For the purposes of this section, it is sometimes preferable pass to the {\em extrinsic point of view}.  Specifically,
let $\calS$ be a compact Riemann surface of the hyperbolic type, i.e., it is uniformized by the hyperbolic plane.
Given a covering map $\pi: \Hyp^2\to \calS$, we obtain a Riemannian structure on $\calS$ which makes $\pi$ a local isometry.
By Nash's theorem, we may isometrically embed $\calS\hookrightarrow\R^{N}$ into an ambient Euclidean space. Now denote the compositions
\[
 U^{n}:=\pi\circ\bfu^n\: : \:I\times \R^{2}\to\calS
\]
defined on $I\times \R^{2}$, see the above discussion. We can express these maps in terms of the ambient coordinates.
Our first task is to identify an actual {\em map} $U$ from $I\times \R^{2}$ into $\calS\hookrightarrow\R^{N}$ which in
some sense corresponds to the limiting object $\Psi^{\infty}_\alpha(t,x)$. The fact that this can be done follows again
from the compactness property of the $\Psi^{\infty}_\alpha(t,x)$. We have the following

\begin{prop}\label{WeakWM} Under the above assumptions, there exists a subsequence of $\{U^n,\,\phi^n,\,\psi^n\}$ which
we denote in the same fashion as well as a function $U(t,\cdot)\in C^0(I; \dot{H}^1)\cap C^1(I;L^{2})$ , such that
\[
 \lim_{n\to\infty}U^n(t,x)=:U(t,x),\quad \lim_{n\to\infty}\nabla_{x,t}U^n(t,x)=\nabla_{x,t}U(t,x)
\]
 where the former limit is the a.e.~pointwise sense and the latter limit is in the $L_{x}^{2}$-sense on fixed time intervals.
 The map $U$ is a weak wave map (in the distributional sense). Also, the second limit is uniform on compact intervals $J\subset I$.
 Finally, the family of functions
\[
 \{\nabla_{x,t}U(t,\,\cdot)\}_{t\in I}\subset L_{x}^{2}
\]
is compact up to rescaling and translational symmetries (which may depend on time).
 \end{prop}
\begin{proof}
We may assume that for times $t\in I$ we have
\[
 \psi^{n}_\alpha(t,\,\cdot)=\Psi^\infty_\alpha(t,\,\cdot)+o_{L^{2}}(1)
\]
But then it follows that for each such $t\in I$, there is a subsequence (depending on $t$) such that also $\phi^n_\alpha(t,\,\cdot)$
converges in the $L^{2}$-sense. To see this, note that
\[
 \phi^n_{\alpha}(t,\cdot) = (\Psi^\infty_\alpha(t,\,\cdot)+o_{L^{2}}(1)) e^{i\Delta^{-1}\sum_{j=1}^2 \phi^n_{j} }
\]
inherits both the physical $L^2$-localization coming from $\Psi^\infty_\alpha(t,\cdot)$ as well as the
Fourier localization of this profile\footnote{This follows
as usual from a Littlewood-Paley trichotomy argument and the energy conservation of the $\phi^n$.}  whence
it is compact and a subsequence converges as claimed.
Picking a dense subset of times $\{t_{i}\}_{i=1}^{\infty}\subset I$ and using the Cantor diagonal argument, one obtains  a subsequence which
we again denote by $\psi^{n}$ etc.~such that $\phi^{n}(t_{i},\,\cdot)$ converges for each $i$ in the $L^{2}$ sense.
By Corollary~\ref{cor:compactV}, it then follows that $\phi^{n}(t,\,\cdot)$ converges in the $L^{2}$ sense, uniformly
on compact sub-intervals of~$I$. In particular, the limit $\phi^\infty$ satisfies $\phi^\infty\in C^0(I;L^2(\R^2))$.
We now use this to infer the existence of $U(t,x)$.  First, introduce a global frame $\{e_{1,2}\}$ on the pull-back bundle of
$T\calS$ under the wave map $U^{n}$ by projecting down the standard frame $\{\bf{e}_1,\,\bf{e}_2\}$, i.e., $e_{j}(t,x):=\pi_*({\bf e}_j)(\bfu^n(t,x))$. Thus
\begin{equation}
 \label{eq:Unconv}
 \partial_\alpha U^n(t, x)=\sum_{k=1,2}e_k^n(t,x)\phi_\alpha^{kn}(t,x)
\end{equation}
Fix some $I'\subset I$ which is compactly contained in~$I$. We now use that the pull-back frame is bounded.
By the preceding, given $\eps>0$ there exists $R$ so large that
\[
 \limsup_{n\to\infty} \big\| \nabla_{t,x} U^n \chi_{[|x|>R]} \big\|_{L^\infty(I'; L^2(|x|>R))} <\eps
\]
On the other hand, it is clear that
\[
 \limsup_{n\to\infty} \big\| \nabla_{t,x} U^n \|_{L^\infty(I'; L^2)} <\infty
\]
By Rellich's theorem we now conclude that up to passing to a subsequence,
$\del_\alpha U^{n}\rightharpoonup X_\alpha $ in $L^\infty(I'; L^2)$ (in the weak-* sense),
as well as  $U^n\to U$ in $L^\infty_{\mathrm{loc}}(I; L^2)$ strongly.
Necessarily then $U\in L^\infty(I',\dot H^1(\R^2))$, see~\eqref{eq:Unconv} as well as $X_\alpha=\del_\alpha U$.
One immediately obtains the stronger statement that $U\in C^0(I,L^2)$ by integrating in time.
One in fact has stronger convergence: first note that
\[
 \del_\alpha e_k^n(t,x) = d(\pi_*)(d{\bf e}_j)(\bfu^n(t,x)) \del_\alpha \bfu^n(t,x)
\]
which implies that $\{e_k^n\}_{n=1}^\infty$ is compact in~$\dot H^1(\R^2)$.
It now follows from~\eqref{eq:Unconv} and Rellich's theorem as before that up to a subsequence one has
\[
 \partial_\alpha U^n(t_i,\cdot)\to \del_\alpha U(t_i,\cdot)
\]
 strongly in $L^2$. By compactness, one therefore also has strongly in $L^2$
\[
 \partial_\alpha U^n(t,\cdot)\to \del_\alpha U(t,\cdot)
\]
uniformly on compact subsets of $I$.  This implies all the convergence and regularity statements of the proposition.
The fact that $U$ is a weak wave map follows from this, as well as from~\cite{FMS}.
\end{proof}

Note that we do not claim that we have uniqueness for the limiting object $U$, and indeed we only have a well-posedness
theory at the level of the $\psi_\alpha$. Thus we cannot purely work at the level of wave maps with compact target $\calS$.
Nevertheless, the latter will play an important role when ruling out certain pathological behaviors, or also to formulate the conservation laws.

For example, we have the following
\begin{cor}
 \label{cor:conservationlaw} Let $U$ be the weak wave map as in Proposition~\ref{WeakWM}. Then one has the
following conservation laws: with $|\cdot|^2=\la\cdot,\cdot\ra$ being the metric on~$\calS$,
\begin{itemize}
 \item $\frac{d}{dt}\sum_{\alpha=0}^2 \int_{\R^2} |\del_\alpha U(t,x)|^2\, dx = 0$
\item $\frac{d}{dt} \int_{\R^2} \la \del_t U(t,x), \del_i U(t,x)\ra\, dx
=0$\quad $i=1,2$
\item  $\frac{d}{dt} \sum_{i=1}^2 \int_{\R^2} x_i \phi(x/R) \la \del_t U(t,x), \del_i U(t,x)\ra\, dx = -\int_{\R^2} |\del_t U(t,x)|^2\, dx + O(r(R))$
\item $\frac{d}{dt}   \sum_{\alpha=0}^2 \int_{\R^2} x_i \phi(x/R)  \frac12 |\del_\alpha U(t,x)|^2\, dx = -\int_{\R^2} \la \del_i U,\del_t U\ra\, dx + O(r(R))$
\end{itemize}
where $\phi$ is a fixed bump function which is equal to one on $|x|\le1$ and
\[
r(R) := \int_{[|x|\ge R]} \sum_{\alpha=0}^2  |\del_\alpha U(t,x)|^2\, dx
\]
\end{cor}
\begin{proof}
These are standard calculations for smooth wave maps. By Proposition~\ref{WeakWM}
one can then pass to the limit.
\end{proof}

Note that one could alternatively express these in terms of $\Psi^\infty_\alpha$.
We will now closely  follow the arguments in~\cite{KeM1}.

\subsubsection{Some preliminary properties of the limiting profiles}

We begin with the following
consequence of finite propagation speed.  Let $I^+:=I\cap [0,\infty)$ where $I$ is the life span of~$\Psi_\alpha^\infty$.

\begin{lemma}
 \label{lem:finitespeed}
Let $M>0$ have the property that
\begin{equation}
 \label{eq:enerM0}
\int_{|x|>\frac{M}{2}} \sum_{\alpha=0}^2 |\Psi_\alpha^\infty(0,x)|^2\, dx<\eps
\end{equation}
Then
\begin{equation}
 \label{eq:enerMt}
\int_{|x|>2 M+t} \sum_{\alpha=0}^2 |\Psi_\alpha^\infty(t,x)|^2\, dx<C\eps
\end{equation}
for all $t\in I^+$. Here $C$ is an absolute constant.
\end{lemma}
\begin{proof}
 By definition, there exist $\bfu^n=(\bfx^n,\bfy^n)\::\: I^+\to \Hyp^2$ which are admissible wave maps such that~\eqref{eq:aa} holds.
Now define
\[
({\bfx}^n_{2}, {\bfy}^n_{2})(0,\cdot):=
\Big(\chi_{[|x|>M]}\frac{{\bfx^n}(0,\cdot)-{\bfx^n}_{0}}{{\bfy}^n_{0}},\:e^{\chi_{[|x|>\frac{M}{2}]}\log[\frac{\bfy^n}{{\bfy}^n_{0}}(0,\cdot)]}
\Big)
\]
where  $\chi_{[|x|>M]}$ is a smooth cutoff to the
set~$\{|x|>M\}$ which equals one on
$\{|x|>\frac54 M\}$, say, and
\begin{align*}
\bfx_0^n&:= \slashint_{[M<|x|<\frac54 M]} \bfx^n(x)\, dx_1dx_2,\qquad
 \bfy_0^n :=
\exp\big(\slashint_{[\frac{M}{2}<|x|< \frac58 M]} \log \bfy^n(x)\, dx_1dx_2\big)
\end{align*}
The construction here is such that $\bfy^n_2=\frac{\bfy^n}{\bfy^n_0}$
on the set~$\{\nabla \chi_{[|x|> M]}\ne0\}$.  Let $\tilde \bfu^n$ be the wave map evolution of the data
\[ \Big( ({\bfx}^n_{2}, {\bfy}^n_{2})(0,\cdot), \big(\frac{\del_t \bfx^n(0,\cdot)}{{\bfy}^n_{0}}, \frac{\del_t \bfy^n(0,\cdot)}{{\bfy}^n_{0}}\big) \Big)\]
 By construction, the energy of $\tilde\bfu^n$
does not exceed~$C\eps$. This requires the use of Poincar\'e's inequality as in the proof of Lemma~\ref{lem:DataLocalization'}.
One now concludes by means of finite propagation speed for classical wave maps, and by passing to the limit $n\to\infty$.
\end{proof}

Next, one has the following lower bound on $\lambda(t)$ in Corollary~\ref{cor:compactV}.

\begin{lemma}
 \label{lem:4.7} Assume $I^+$ is finite. After rescaling, we may assume that $I^+=[0,1)$.
There exists a constant $C_0(K)$ depending on the compact set~$K$ in  Corollary~\ref{cor:compactV},
such that
\begin{equation}
 \label{eq:lambdalower}
0<\frac{C_0(K)}{1-t}\le \lambda(t)
\end{equation}
for all $0\le t<1$.
\end{lemma}
\begin{proof}
 Take any sequence $t_j\to 1$. Consider the limiting profile $\{\tilde\Psi^\infty_{\alpha,j}\}_{\alpha=0}^2$ with data
$\lambda(t_j)^{-1}\Psi^\infty_{\alpha}
(t_j,(\cdot-\bar{x}(t_j))\lambda(t_j)^{-1})\}_{\alpha=0}^2$. By the
well posedness theory of the limiting profiles in
Section~\ref{subsec:perturbprofile}, one infers that the
$\{\tilde\Psi^\infty_{\alpha,j}\}_{\alpha=0}^2$ have a fixed life
span  independent of~$j$ which depends only on the compact set~$K$.
By the uniqueness property of the solutions and rescaling,
$(1-t_j)\lambda(t_j)\ge C_0(K)$ as claimed.
\end{proof}

Next, combining this with Lemma~\ref{lem:finitespeed} one concludes the following support property of the $\Psi_\alpha^\infty$
with finite life span.

\begin{lemma}
 \label{lem:4.8} Let $\Psi_\alpha^\infty$ be as in the previous lemma. Then there exists $x_0\in\R^2$ such that
\[
 \supp(\Psi^\infty_\alpha(t,\cdot)) \subset B(x_0,1-t)
\]
 for all $0\le t<1$, $\alpha=0,1,2$.
\end{lemma}
\begin{proof}
This follows the exact same reasoning as in Lemma~4.8 of~\cite{KeM1}. One uses
Lemma~\ref{lem:finitespeed} instead of their Lemma~2.17 and Lemma~\ref{lem:4.7} instead of their
Lemma~4.7.
\end{proof}

Next, we turn to the vanishing moment condition of Propositions 4.10 and 4.11 in~\cite{KeM1}.

\begin{prop}
 \label{prop:4.10} Let $\Psi_\alpha^\infty$ be as above and assume that $I^+$ is finite. Then for $i=1,2$,
\[
 \int_{\R^2} \la \del_i U, \del_t U\ra\, dx = \Re\int_{\R^2} \Psi^\infty_i \bar{\Psi}_0^\infty\, dx = 0
\]
for all times in~$I^+$.
\end{prop}
\begin{proof}
Assume that
\[
 \Re\int_{\R^2} \Psi^\infty_1 \bar{\Psi}_0^\infty\, dx >\gamma>0
\]
This implies that the approximating sequence $\bfu^n$ satisfies
\[
 \int_{\R^2} \la \del_1 \bfu^n, \del_t \bfu^n\ra\, dx > \gamma>0
\]
for large $n$. Following~\cite{KeM1} we apply a Lorentz transformation
\[
 L_d(t,x):= \Big(\frac{t-dx_1}{\sqrt{1-d^2}},\frac{x_1-dt}{\sqrt{1-d^2}},x_2\Big)
\]
to the $u^n$. Note that for any $\eps>0$ one has from Lemma~\ref{lem:4.8} that
\[
 \sum_{\alpha=0}^2\int_{|x|\ge 1-t} |\del_\alpha \bfu^n(t,x)|^2\, dx <\eps
\]
for all $t\in I^+=[0,1)$ and sufficiently large~$n$. Then the argument in~\cite{KeM1} implies that
there exists $d$ small with the property that
\[
 \limsup_{n\to\infty} E(\bfu^n\circ L_d)<\Ecrit
\]
By our induction hypothesis, $\| \psi^{n,d}\|_{S(I^+\times\R^2)} < M<\infty$ for all sufficiently large $n$. Here
$\psi^{n,d}$ are the Coulomb components of the admissible wave maps $\bfu^n\circ L_d$. Note that the Coulomb components
 $\psi^{n,d}$ do not obey a simple transformation law relative to the Coulomb components $\psi^n$  of~$\bfu^n$.
Nonetheless, it is possible to conclude from this that \[\limsup_{n\to\infty}\| \psi^{n}\|_{S(I^+\times\R^2)} <M_1<\infty\]
via Remark~\ref{rem:bilineargivesS} which gives us the desired contradiction.
Thus, we need to prove that for each $k_{1}>k_{2}$
\begin{equation}\label{eq:psingoal}
\sum_{\substack{\kappa_{1,2}\in \calC_{m_0}\\ \dist(\kappa_1,\,\kappa_2)\gtrsim 2^{-m_0}}}
2^{-k_2}
P_{k_1,\kappa_1}\psi^n P_{k_2,\kappa_2}\psi^n = f_{k_1,k_2} + g_{k_1,k_2}
\end{equation}
where $m_0$ is a large depending on~$E_C$,
where we have the bounds~\eqref{eq:bilinLambda2} for $f_{k_1,k_2}$ and $g_{k_1,k_2}$. Furthermore, we need to show that
\[
P_kQ_{>k}\psi^n=h_k+i_k
\]
with the bounds stated in Remark~\ref{rem:bilineargivesS}. We establish this for the bilinear expression,
the corresponding computations for $P_kQ_{>k}\psi^n$ being similar.
First, we claim the following bound for $\psi^{n,d}$:
\begin{equation}\label{eq:psind}
  \sum_{k_{1}>k_{2}}\sum_{\substack{\kappa_{1,2}\in \calC_{m_0}\\ \dist(\kappa_1,\,\kappa_2)\gtrsim 2^{-m_0}}}
2^{-k_2}
\|P_{k_1,\kappa_1}\psi^{n,d} P_{k_2,\kappa_2}\psi^{n,d}\|_{L_{t,x}^{2}}^{2}<\Lambda'
\end{equation}
 This, however, is immediate from the angular separation
and~\eqref{eq:bilin2} with a constant $\Lambda'$ which depends on~$M$ and $\Ecrit$.
In fact, we need something slightly stronger due to the usual tail issues:
\begin{equation}\label{eq:psind'}
 \sup_y \sum_{k_{1}>k_{2}}\sum_{\substack{\kappa_{1,2}\in \calC_{m_0}\\ \dist(\kappa_1,\,\kappa_2)\gtrsim 2^{-m_0}}}
2^{-k_2}
\|P_{k_1,\kappa_1}\psi^{n,d} \tau_y P_{k_2,\kappa_2}\psi^{n,d}\|_{L_{t,x}^{2}}^{2}<\Lambda'
\end{equation}
where $\tau_y$ is a translation by~$y\in\R^2$.
Next, we claim the following estimate:
\begin{equation}\label{eq:phindcrux}
  \sum_{k_{1}>k_{2}}\sum_{\substack{\kappa_{1,2}\in \calC_{m_0}\\ \dist(\kappa_1,\,\kappa_2)\gtrsim 2^{-m_0}}}
2^{-k_2}
\|P_{k_1,\kappa_1}\phi^{n,d} P_{k_2,\kappa_2}\phi^{n,d}\|_{L_{t,x}^{2}}^{2}<\Lambda'
\end{equation}
where $\phi^{n,d}$ are the derivative components of the $\bfu^n\circ L_d$. This is the same as
\begin{equation}
 \nn
  \sum_{k_{1}>k_{2}}\sum_{\substack{\kappa_{1,2}\in \calC_{m_0}\\ \dist(\kappa_1,\,\kappa_2)\gtrsim 2^{-m_0}}}
2^{-k_2}
\|P_{k_1,\kappa_1}\big(\psi^{n,d}e^{-i\del^{-1}\phi^{n,d}}\big)\cdot P_{k_2,\kappa_2}\big( \psi^{n,d} e^{-i\del^{-1}\phi^{n,d}} \big)  \|_{L_{t,x}^{2}}^{2}<\Lambda'
\end{equation}
where we wrote the phase $-i\del^{-1}\phi^{n,d} = -i \Re\sum_{j=1}^2 (-\Delta)^{-1}\del_j \phi^{n,d}$ schematically.
This follows from~\eqref{eq:psind} and the Strichartz estimate
\begin{equation}
 \label{eq:phiStrich}
\Big( \sum_{k\in\Z}  2^{-\frac32 k} \sup_{j\ge10} \sum_{c\in\calD_{k,j}} 2^{-(1-2\eps)j} \|P_c \phi^{n,d}\|_{L^4_t L^\infty_x}^2  \Big)^{\frac12} \les M
\end{equation}
To prove \eqref{eq:phiStrich}, one uses the corresponding bound on~$\psi^{n,d}$ (which is part of the~$S$-norm), energy conservation,
and a simple Littlewood-Paley trichotomy. To prove~\eqref{eq:phindcrux}, one argues as follows. Split
\begin{align}
& \sum_{k_{1}>k_{2}}\sum_{\substack{\kappa_{1,2}\in \calC_{m_0}\\ \dist(\kappa_1,\,\kappa_2)\gtrsim 2^{-m_0}}}
2^{-k_2}
\|P_{k_1,\kappa_1}\big(\psi^{n,d}e^{-i\del^{-1}\phi^{n,d}}\big)\cdot P_{k_2,\kappa_2}\big( \psi^{n,d} e^{-i\del^{-1}\phi^{n,d}} \big)  \|_{L_{t,x}^{2}}^{2} \nn    \\
&  \les
  \sum_{k_{1}>k_{2}}\weg\sum_{\substack{\kappa_{1,2}\in \calC_{m_0}\\ \dist(\kappa_1,\,\kappa_2)\gtrsim 2^{-m_0}}}
\weg \weg
2^{-k_2}
\|P_{k_1,\kappa_1}\big(\psi^{n,d} P_{<k_1-m_0} e^{-i\del^{-1}\phi^{n,d}}\big)\cdot P_{k_2,\kappa_2}\big( \psi^{n,d} P_{<k_2-m_0} e^{-i\del^{-1}\phi^{n,d}} \big)  \|_{L_{t,x}^{2}}^{2}
\label{eq:S1} \\
& +  \sum_{k_{1}>k_{2}}\weg \sum_{\substack{\kappa_{1,2}\in \calC_{m_0}\\ \dist(\kappa_1,\,\kappa_2)\gtrsim 2^{-m_0}}}\weg\weg 2^{-k_2} \|P_{k_1,\kappa_1}\big(\psi^{n,d} P_{<k_1-m_0} e^{-i\del^{-1}\phi^{n,d}}\big)\cdot P_{k_2,\kappa_2}\big( \psi^{n,d} P_{>k_2-m_0} e^{-i\del^{-1}\phi^{n,d}} \big)  \|_{L_{t,x}^{2}}^{2} \label{eq:S2}\\
&  + \sum_{k_{1}>k_{2}}\weg\sum_{\substack{\kappa_{1,2}\in \calC_{m_0}\\ \dist(\kappa_1,\,\kappa_2)\gtrsim 2^{-m_0}}}\weg\weg 2^{-k_2}  \|P_{k_1,\kappa_1}\big(\psi^{n,d} P_{>k_1-m_0} e^{-i\del^{-1}\phi^{n,d}}\big)\cdot P_{k_2,\kappa_2}\big( \psi^{n,d} P_{<k_2-m_0} e^{-i\del^{-1}\phi^{n,d}} \big)  \|_{L_{t,x}^{2}}^{2} \label{eq:S3} \\
& +  \sum_{k_{1}>k_{2}}\weg\sum_{\substack{\kappa_{1,2}\in \calC_{m_0}\\ \dist(\kappa_1,\,\kappa_2)\gtrsim 2^{-m_0}}}\weg\weg 2^{-k_2} \|P_{k_1,\kappa_1}\big(\psi^{n,d} P_{>k_1-m_0} e^{-i\del^{-1}\phi^{n,d}}\big)\cdot P_{k_2,\kappa_2}\big( \psi^{n,d} P_{>k_2-m_0} e^{-i\del^{-1}\phi^{n,d}} \big)  \|_{L_{t,x}^{2}}^{2}
\label{eq:S4}
\end{align}
 In \eqref{eq:S1} one reduces matters to \eqref{eq:psind'} by placing the exponential in~$\Linf$. Next, to bound~\eqref{eq:S2} one notes that
\[
 \big\| P_{>k_2-m_0} e^{-i\del^{-1}\phi^{n,d}} \big\|_{L^4_t L^\infty_x} \les 2^{-\frac{k_2}{4}} M
\]
where the implicit constant depends on~$\Ecrit$. Therefore,
\begin{align*}
 & \sum_{k_{1}>k_{2}}\weg \sum_{\substack{\kappa_{1,2}\in \calC_{m_0}\\ \dist(\kappa_1,\,\kappa_2)\gtrsim 2^{-m_0}}}\weg\weg 2^{-k_2} \|P_{k_1,\kappa_1}\big(\psi^{n,d} P_{<k_1-m_0} e^{-i\del^{-1}\phi^{n,d}}\big)\cdot P_{k_2,\kappa_2}\big( \psi^{n,d} P_{>k_2-m_0} e^{-i\del^{-1}\phi^{n,d}} \big)  \|_{L_{t,x}^{2}}^{2} \\
 &\les  \sum_{k_{1}>k_{2}}\weg \sum_{\substack{\kappa_{1,2}\in \calC_{m_0}\\ \dist(\kappa_1,\,\kappa_2)\gtrsim 2^{-m_0}}}\weg\weg\weg 2^{-k_2} \|P_{k_1} \psi^{n,d}\|_{\ener}^2  \big( \sum_{\ell>k_2}
\sum_{\substack{c,c'\in\calD_{\ell,k_2-\ell} \\ \dist(c,c')\les 2^{k_2}}} \weg\weg
 \| P_c P_\ell \psi^{n,d}\|_{L^4_t L^\infty_x} \|P_{c'} P_{\ell+O(m_0) } (e^{-i\del^{-1}\phi^{n,d}}-1) \|_{L^4_t L^\infty_x}\big)^{2} \\
& \les M^6
\end{align*}
using the Strichartz estimate from above. The remaining terms are the same.
This concludes the proof of~\eqref{eq:phindcrux}. By the same logic, one also obtains
\[
   \sum_{k_{1}>k_{2}}\sum_{\substack{\kappa_{1,2}\in \calC_{m_0}\\ \dist(\kappa_1,\,\kappa_2)\gtrsim 2^{-m_0}}}
2^{-k_2}
\|P_{k_1,\kappa_1}Q_{\le k_1+C_2}\phi^{n,d} P_{k_2,\kappa_2}Q_{\le k_2+C_2}\phi^{n,d}\|_{L_{t,x}^{2}}^{2}<\Lambda'
\]
where $C_2$ is a large constant depending only on the energy which will be determined later.
This then implies the following version {\em without} the Lorentz transforms
\[
   \sum_{k_{1}>k_{2}}\sum_{\substack{\kappa_{1,2}\in \calC_{m_0'}\\ \dist(\kappa_1,\,\kappa_2)\gtrsim 2^{-m_0'}}}
2^{-k_2}
\|P_{k_1,\kappa_1}Q_{\le k_1+C_2}\phi^{n} P_{k_2,\kappa_2}Q_{\le k_2+C_2}\phi^{n}\|_{L_{t,x}^{2}}^{2}<\Lambda'
\]
provided $d$ is chosen small enough, but depending only on~$\Ecrit$ (so that $m_0'$ is close to $m_0$). Finally we claim that
\begin{align}\label{eq:hypphi1}
 &P_{k_1}Q_{> k_1+C_2}\phi^{n} P_{k_2}\phi^{n} \\
&P_{k_1}\phi^{n} P_{k_2}Q_{> k_2+C_2}\phi^{n} \label{eq:hypphi2}
\end{align}
can  both   be included in $g_{k_1, k_2}$.
To see this, one first expands
\begin{align}
P_{k_1}Q_{> k_1+C_2}\phi^{n} &= P_{k_1}Q_{> k_1+C_2} \big[ \psi^{n} e^{-i\del^{-1}\phi^{n}} \big] \nn \\
&  = \sum_{\ell > k_1+C_2-10} P_{k_1}Q_{> k_1+C_2} \big[ P_\ell \psi^{n} P_\ell e^{-i\del^{-1}\phi^{n}} \big] \label{eq:hhphind1} \\
& + \sum_{k_1<\ell \le k_1+C_2-10} P_{k_1}Q_{> k_1+C_2} \big[ P_\ell \psi^{n} P_\ell e^{-i\del^{-1}\phi^{n}} \big] \label{eq:hhphind2} \\
& + P_{k_1}Q_{> k_1+C_2} \big[ P_{<k_1-5}  \psi^{n} P_{k_1} e^{-i\del^{-1}\phi^{n}} \big] \label{eq:lhphind} \\
& + P_{k_1}Q_{> k_1+C_2} \big[ P_{k_1}  \psi^{n} P_{<k_1-5} e^{-i\del^{-1}\phi^{n}} \big] \label{eq:hlphind}
\end{align}
and then inserts these decompositions into~\eqref{eq:hypphi1}. For \eqref{eq:hhphind1} one places  $P_{k_2}\phi^{n}$ into~$L^4_t L^\infty_x$, and
$P_{k_1}Q_{> k_1+100}\phi^{n}$ into~$L^4_t L^2_x$ followed by an application of~\eqref{eq:phiStrich} with
caps of size~$2^{k_1}$; more precisely, $P_\ell e^{-i\del^{-1}\phi^{n}}$ goes into~$L^4_t L^\infty_x$ as before, and $P_\ell \psi^{n}$
gets placed into~$\ener$ (see Lemma~\ref{lem:enersquaresum} for the issue of square-summing the $\ener$-norm of $\psi^{n}$ over caps of size~$2^{k_1}$).
Note that one gains a smallness factor of the form~$2^{-\frac{C_2}{10}}$ due to the improved Strichartz bounds. Next, we consider~\eqref{eq:hlphind}
and the remaining terms~\eqref{eq:hhphind2} and~\eqref{eq:hlphind} will follow similar arguments.
 Now we decompose further:
\begin{align}
 & P_{k_1}Q_{> k_1+C_2} \big[ P_{k_1}  \psi^{n} P_{<k_1-5} e^{-i\del^{-1}\phi^{n}} \big]  \nn \\
& = P_{k_1}Q_{> k_1+C_2} \big[ Q_{>k_1+C_2-10} P_{k_1}  \psi^{n} P_{<k_1-5} e^{-i\del^{-1}\phi^{n}} \big] \label{eq:hlphind1} \\
&+ P_{k_1}Q_{> k_1+C_2} \big[ Q_{\le k_1+C_2-10} P_{k_1}  \psi^{n} P_{<k_1-5}Q_{> k_1+C_2-10} e^{-i\del^{-1}\phi^{n}} \big] \label{eq:hlphind2}
\end{align}
For~\eqref{eq:hlphind1} one estimates
\begin{align*}
 & \| P_{k_1}Q_{> k_1+C_2} \big[ Q_{>k_1+C_2-10} P_{k_1}  \psi^{n} P_{<k_1-5} e^{-i\del^{-1}\phi^{n}} \big]\;P_{k_2}\phi^{n}\|_{L^2_{t,x}} \\
&\les 2^{k_2} \| Q_{>k_1+C_2-10} P_{k_1}  \psi^{n}\|_{L^2_{t,x}} \|P_{k_2}\phi^{n}\|_{\ener} \les 2^{k_2} 2^{-\frac{k_1+C_2}{2}} \|P_{k_1}\psi^{n}\|_{S[k_1]}
\|P_{k_2}\phi^{n}\|_{\ener}
\end{align*}
which is sufficient since it gains the smallness $2^{-\frac{C_2}{2}}$. Finally, we use~\eqref{eq:compat1} for the case when we
substitute \eqref{eq:hlphind2} for $P_{k_1}Q_{>k_1+C_2}\phi^n$; one can then write
\begin{align*}
 & P_{k_1}Q_{> k_1+C_2}\phi^{n} P_{k_2}\phi^{n} \\
& = P_{k_1}Q_{> k_1+C_2} \big[ Q_{\le k_1+C_2-10} P_{k_1}  \psi^{n} \del_t^{-1} P_{<k_1-5}Q_{> k_1+C_2-10} ((\phi^{n}+\nabla^{-1}(\phi^{n}\phi^{n})) e^{-i\del^{-1}\phi^{n}}) \big] P_{k_2}\phi^{n}
\end{align*}
where we have written \eqref{eq:compat1} schematically in the form
\[
 \del_t \del^{-1}\phi^{n}= \phi^{n}+\nabla^{-1}(\phi^{n}\phi^{n})
\]
The contribution of $\phi^{n}$ is easy, it is placed again in $L^4_t L^\infty_x$ (of course after applying the usual trichotomy to
$\phi^{n}e^{-i\del^{-1}\phi^{n}}$). On the other hand, due to the determinant structure of $\nabla^{-1}(\phi^{n}\phi^{n})$ we have
\[
\nabla^{-1}(\phi^{n}\phi^{n})=\nabla^{-1}(\psi^{n}\psi^{n})
\]
By using a further Hodge decomposition of the inputs on the right, we have for each $k\in\Z$
\[
\|P_{k}\nabla^{-1}(\psi^{n}\psi^{n})\|_{L_t^2\dot{H}^{\frac{1}{2}}}\lesssim \|\psi^n\|_{S}^{2},
\]
and from here we get
\begin{equation}\label{eq:tedious}
\|\del_t^{-1} P_{<k_1-5}Q_{> k_1+C_2-10}(\nabla^{-1}(\psi^{n}\psi^{n})e^{-i\del^{-1}\phi^{n}})\|_{L_t^2L_x^{\infty}}\ll 2^{-\frac{k_1}{2}}\|\psi\|_{S}^{2}
\end{equation}
and from here we get
\begin{align*}
&\|P_{k_1}Q_{> k_1+C_2} \big[ Q_{\le k_1+C_2-10} P_{k_1}  \psi^{n} \del_t^{-1} P_{<k_1-5}Q_{> k_1+C_2-10}
(\nabla^{-1}(\phi^{n}\phi^{n})e^{-i\del^{-1}\phi^{n}}) \big] P_{k_2}\phi^{n}\|_{L_{t,x}^{2}}
\\&\ll 2^{k_2-k_1}\|P_{k_2}\phi^{n}\|_{L_{t}^{\infty}L_{x}^{2}}\|\psi\|_{S}^{2}
\end{align*}
This concludes the proof of~\eqref{eq:hypphi1}.
For~\eqref{eq:hypphi2} one argues similary. By following the same Littlewood-Paley trichotomies, one is eventually lead to the most difficult case
\begin{align*}
 & P_{k_1}\phi^{n} P_{k_2}Q_{> k_2+C_2}\phi^{n} \\
& = P_{k_1}\phi^{n} P_{k_2} Q_{> k_2+C_2} \big[ Q_{\le k_2+C_2-10} P_{k_2}  \psi^{n} \del_t^{-1}
P_{<k_2-5}Q_{> k_2+C_2-10} ((\phi^{n}+\nabla^{-1}(\phi^{n}\phi^{n})) e^{-i\del^{-1}\phi^{n}}) \big]
\end{align*}
where we again used the curl equation~\eqref{eq:compat1}. The $\phi^{n}$ term is again easier, whereas for the nonlinear term we again use
\[
 \nabla^{-1}(\phi^{n}\phi^{n}) = \nabla^{-1}(\psi^{n}\psi^{n})
\]
Then as before we use \eqref{eq:tedious}, in order to infer that
\begin{align*}
&\|P_{k_1}\phi^{n} P_{k_2} Q_{> k_2+C_2} \big[ Q_{\le k_2+C_2-10} P_{k_2}  \psi^{n} \del_t^{-1} P_{<k_2-5}Q_{> k_2+C_2-10}
(\nabla^{-1}(\phi^{n}\phi^{n})) e^{-i\del^{-1}\phi^{n}}) \big]\|_{L_{t,x}^{2}}\\
&\lesssim \|P_{k_1}\phi^{n}\|_{L_{t}^\infty L_x^2}\|Q_{\le k_2+C_2-10} P_{k_2}  \psi^{n}\|_{L_{t,x}^\infty}
\| \del_t^{-1} P_{<k_2-5}Q_{> k_2+C_2-10} (\nabla^{-1}(\phi^{n}\phi^{n}) e^{-i\del^{-1}\phi^{n}})\|_{L_{t}^2L_x^\infty}\\
&\ll 2^{\frac{k_2}{2}}\|\psi\|_{S}^2 \|P_{k_1}\phi^{n}\|_{L_t^\infty L_x^2},
\end{align*}
which justifies us in $g_{k_1, k_2}$.
In conclusion, we have now shown that we can write
 \begin{equation}\label{eq:phibil}
 \sum_{\substack{\kappa_{1,2}\in \calC_{m_0'}\\ \dist(\kappa_1,\,\kappa_2)\gtrsim 2^{-m_0'}}}
P_{k_1,\kappa_1}\phi^{n} P_{k_2,\kappa_2}\phi^{n}=\tilde{f}_{k_1, k_2}+\tilde{g}_{k_1, k_2}
\end{equation}
with bounds as in \eqref{eq:bilinLambda2}.
The goal is now to deduce \eqref{eq:psingoal} from this estimate.  For this purpose, fix $k_1>k_2+C_1$ and caps $\kappa_1,\kappa_2\in \calC_{m_0'}$
as above. We now describe how to break up
\[
 P_{k_1,\kappa_1} \psi^n\cdot P_{k_2,\kappa_2} \psi^n = P_{k_1,\kappa_1} (\phi^n e^{-i\del^{-1}\phi^{n}}) \cdot P_{k_2,\kappa_2} (\phi^n
e^{-i\del^{-1}\phi^{n}})
\]
into various pieces which then constitute $f_{k_1,k_2}$ and $g_{k_1,k_2}$, respectively when summed over the caps.
First, write
\begin{align}
 & P_{k_1,\kappa_1} (\phi^n e^{-i\del^{-1}\phi^{n}}) \cdot P_{k_2,\kappa_2} (\phi^n
e^{-i\del^{-1}\phi^{n}}) \nn \\
&= \sum_{i_1,i_2=1}^3 P_{k_1,\kappa_1} (\phi^n  A_{i_1} e^{-i\del^{-1}\phi^{n}}) \cdot P_{k_2,\kappa_2} (\phi^n B_{i_2}
e^{-i\del^{-1}\phi^{n}})  \label{eq:gack}
\end{align}
where
\[
 A_1 = P_{<k_1-m_0'-10},\quad A_2 = P_{k_1-m_0'-10\le \cdot <k_1+C_2},\quad A_3 = P_{\ge k_1+C_2}
\]
and similarly for $B_i$. Here $C$ is large depending on $\Ecrit$.  If $i_2=3$, then one estimates
\begin{align*}
 & \| P_{k_1,\kappa_1} (\phi^n  A_{i_1} e^{-i\del^{-1}\phi^{n}}) \cdot P_{k_2,\kappa_2} (\phi^n B_{i_2}
e^{-i\del^{-1}\phi^{n}})  \|_{L^2_{t,x}} \\
&\les \sum_{m} 2^{-\sigma|k_1-m|} \|P_{m}\phi^n\|_{\ener} \sum_{\ell\ge k_2+C}\sum_{\substack{c_1,c_2\in\calD_{\ell,k_2-\ell} \\\dist(c_1,c_2)\les 2^{k_2}}} 2^{-\ell}
\|P_{c_1} \phi^n\|_{L^4_t L^\infty_x} \| P_{c_2} [\phi^n e^{-i\del^{-1}\phi^{n}}]  \|_{L^4_t L^\infty_x} \\
&\les 2^{-\frac{C_2}{10}}  2^{\frac{k_2}{2}} \sum_{m} 2^{-\sigma|k_1-m|} \|P_{m}\phi^n\|_{\ener}  \Big( \sum_{\ell>k_2} 2^{-\sigma(\ell-k_2)} \|P_\ell \psi^n\|_{S[\ell]}\Big)^2
\end{align*}
with an implicit constant which is allowed to depend on the energy. Therefore, this is placed in~$g_{k_1,k_2}$.
The case where $i_1=3$ is similar. Next, suppose that $i_1=1$ {\em and} $i_2=1$.  Then  the cap localization passes on to the $\phi^n$
and due to~\eqref{eq:phibil} one places the resulting expression into~$f_{k_1,k_2}$. We are left with three cases: $i_1=1, i_2=2$, and $i_1=2,i_2=1$,
and $i_1=i_2=2$. Next, observe that we may assume that
\[
 P_{k_1,\kappa_1} (\phi^n  A_{i_1} e^{-i\del^{-1}\phi^{n}}) = P_{k_1,\kappa_1} (P_{>k_1-C_2}\phi^n  A_{i_1} e^{-i\del^{-1}\phi^{n}})
\]
and
\[ P_{k_2,\kappa_2} (\phi^n B_{i_2}
e^{-i\del^{-1}\phi^{n}})= P_{k_2,\kappa_2} (P_{>k_2-C_2}\phi^n B_{i_2}
e^{-i\del^{-1}\phi^{n}})
\]
for otherwise one obtains smallness from Bernstein's inequality. For example, consider now $i_1=1$, $i_2=2$ which is
\begin{align*}
& P_{k_1,\kappa_1} (P_{>k_1-C_2}\phi^n  P_{<k_1-m_0'-10} e^{-i\del^{-1}\phi^{n}}) \cdot P_{k_2,\kappa_2} (P_{>k_2-C_2}\phi^n  P_{k_2-m_0'-10\le \cdot <k_2+C_2}
e^{-i\del^{-1}\phi^{n}})   \\
& = P_{k_1,\kappa_1} (P_{>k_1-C_2}\phi^n  P_{<k_1-m_0'-10} e^{-i\del^{-1}\phi^{n}}) \cdot P_{k_2,\kappa_2} (P_{>k_2-C_2}\phi^n  P_{k_2-m_0'-10\le \cdot <k_2+C_2}
\del^{-1} [
\phi^n
e^{-i\del^{-1}\phi^{n}}]   )
\end{align*}
Now we distinguish two more cases: either the exponential in the second factor has frequency $<2^{ k_2-m_0'-20}$ or not. In the former case,
one obtains
\begin{align*}
& P_{k_1,\kappa_1} (P_{>k_1-C_2}\phi^n  P_{<k_1-m_0'-10} e^{-i\del^{-1}\phi^{n}}) \cdot P_{k_2,\kappa_2} (P_{>k_2-C_2}\phi^n  P_{k_2-m_0'-10\le \cdot <k_2+C_2}
\del^{-1} [
\phi^n
e^{-i\del^{-1}\phi^{n}}]   ) \\
&= P_{k_1,\kappa_1} (P_{>k_1-C_2}\phi^n  P_{<k_1-m_0'-10} e^{-i\del^{-1}\phi^{n}}) \cdot P_{k_2,\kappa_2} (P_{>k_2-C_2}\phi^n  P_{k_2-m_0'-10\le \cdot <k_2+C_2}
\del^{-1} [
\phi^n P_{<k_2-m_0'-20}
e^{-i\del^{-1}\phi^{n}}]   )
\end{align*}
Now perform a cap decomposition of the first and second $\phi^n$ factors inside the $P_{k_2,\kappa_2}$ term. Observe that due to the fact
that the frequencies of these factors are approximately $2^{k_2}$ at least one of them has to have angular separation with the cap $\kappa_1$
from the first factor by an amount comparable to~$2^{-m_0'}$. We may therefore place this expression into~$f_{k_1,k_2}+g_{k_1,k_2}$ in view of~\eqref{eq:phibil}.
If, on the other hand, the exponential in the second factor has frequency $>2^{ k_2-m_0'-20}$, then one writes
\begin{align*}
 &= P_{k_1,\kappa_1} (P_{>k_1-C_2}\phi^n  P_{<k_1-m_0'-10} e^{-i\del^{-1}\phi^{n}}) \cdot P_{k_2,\kappa_2} (P_{>k_2-C_2}\phi^n  P_{k_2-m_0'-10\le \cdot <k_2+C_2}
\del^{-1} [
\phi^n P_{\ge k_2-m_0'-20}
e^{-i\del^{-1}\phi^{n}}]   )  \\
 &= P_{k_1,\kappa_1} (P_{>k_1-C_2}\phi^n  P_{<k_1-m_0'-10} e^{-i\del^{-1}\phi^{n}}) \cdot\\
&\qquad\qquad \cdot P_{k_2,\kappa_2}
(P_{>k_2-C_2}\phi^n  P_{k_2-m_0'-10\le \cdot <k_2+C_2}
\del^{-1} [
\phi^n P_{\ge k_2-m_0'-20}\del^{-1}[ \phi^n
e^{-i\del^{-1}\phi^{n}}]]   )
\end{align*}
The idea here is to place the entire expression into~$L^2_{t,x}$ by putting the
first factor into $\ener$, i.e., estimating
\[ \|P_{k_1,\kappa_1} (P_{>k_1-C_2}\phi^n  P_{<k_1-m_0'-10} e^{-i\del^{-1}\phi^{n}}) \|_{\ener} \les \|P_{k_1} \phi^n\|_{\ener}
\]
followed by the estimate
\begin{align}
 &  \| P_{k_2,\kappa_2}
(P_{>k_2-C_2}\phi^n  P_{k_2-m_0'-10\le \cdot <k_2+C_2}
\del^{-1} [
\phi^n P_{\ge k_2-m_0'-20}\del^{-1}[ \phi^n
e^{-i\del^{-1}\phi^{n}}]   )  \|_{L^2_t L^\infty_x} \nn \\
&\les 2^{k_2} \| P_{k_2,\kappa_2}
(P_{>k_2-C_2}\phi^n  P_{k_2-m_0'-10\le \cdot <k_2+C_2}
\del^{-1} [
\phi^n P_{\ge k_2-m_0'-20}\del^{-1}[ \phi^n
e^{-i\del^{-1}\phi^{n}}] ]  )  \|_{L^2_t L^2_x}  \nn \\
&\les 2^{k_2} \| P_{k_2,\kappa_2}
(P_{>k_2-C_2}\phi^n  P_{k_2-m_0'-10\le \cdot <k_2+C_2}
\del^{-1} [ P_{<k_2-m_0'- C_4}
\phi^n P_{\ge k_2-m_0'-20}\del^{-1}[ \phi^n
e^{-i\del^{-1}\phi^{n}}] ]  )  \|_{L^2_t L^2_x} \label{eq:threephis1} \\
& + 2^{k_2} \sum_{\substack{|k-k'|\le m_0'+C_4\\ k\ge k_2-m_0'- C_4 }} \| P_{k_2,\kappa_2}
(P_{>k_2-C_2}\phi^n  P_{k_2-m_0'- C_4\le \cdot <k_2+C_2}
\del^{-1} [ P_k
\phi^n P_{k'}\del^{-1}[ \phi^n
e^{-i\del^{-1}\phi^{n}}] ]  )  \|_{L^2_t L^2_x}  \label{eq:threephis2}
\end{align}
Note that we may reduce~\eqref{eq:threephis2} to (with possibly very large $O(1)$ but only depending on the energy)
\begin{equation}\label{eq:threephis3}
 2^{k_2}  \| P_{k_2,\kappa_2}
(P_{k_2+O(1)}\phi^n
\del^{-1} P_{k_2+O(1)}[ P_{k_2+O(1)}
\phi^n P_{k_2+O(1)}\del^{-1}[ P_{k_2+O(1)}\phi^n\:
e^{-i\del^{-1}\phi^{n}}] ]  )  \|_{L^2_t L^2_x}
\end{equation}
since the extremely large frequencies give a gain of a smallness factor whence that case can be place entirely into
the bootstrap term~$g_{k_1,k_2}$.
We chose $C_4$ here so large that the entire expression~\eqref{eq:threephis1} is placed in the bootstrap term $g_{k_1,k_2}$.
To see this, one estimates
\begin{align*}
 \eqref{eq:threephis1} &\les 2^{-k_2}   \|P_{k_2+O(1)}  \phi^n\|_{L^6_{t,x}} \| P_{<k_2-m_0'-25}
\phi^n \|_{L^6_t L^6_x} \|P_{k_2+O(1)} [   \phi^n
e^{-i\del^{-1}\phi^{n}}] \|_{L^6_t L^{6}_x} \\
&\les 2^{-k_2} \|P_{k_2+O(1)} \phi^n\|_{L^6_{t,x}} \sum_{\ell<k_2-m_0'-C_4}  \| P_\ell
\phi^n \|_{L^6_t L^6_x} \|P_{k_2+O(1)} [ \phi^n
e^{-i\del^{-1}\phi^{n}}] \|_{L^6_t L^6_x} \\
&\les 2^{-k_2} \|P_{k_2+O(1)} \phi^n\|_{L^6_{t,x}} \sum_{\ell<k_2-m_0'-C_4} 2^{\frac{\ell}{2}} \| P_\ell
\psi^n \|_{S[\ell]} \|P_{k_2+O(1)} [ \phi^n
e^{-i\del^{-1}\phi^{n}}] \|_{L^6_t L^6_x} \\
&\les 2^{-\frac{C_4}{2}} 2^{\frac{k_2}{2}}  \Big(
\sum_{\ell} 2^{-\frac14|\ell-k_2|} \|P_\ell \psi^n\|_{S[\ell]}\Big)^3
\end{align*}
Second, with each  $S_0:=\sum_j P_jQ_{\le j+C_3}$ and $S_1:=\sum_j P_jQ_{> j+C_3}$ where $C_3$ is a large constant depending only on the energy,
\begin{align*}
 \eqref{eq:threephis3} &\les \sum_{i_1,i_2,i_3=0,1} 2^{-k} \|P_{k_2+O(1)} S_{i_1}\phi^n\|_{L^6_{t,x}} \| P_{k_2+O(1)} S_{i_2}
\phi^n \|_{L^6_t L^6_x} \|P_{k_2+O(1)} S_{i_3} \phi^n
 \|_{L^6_t L^{6}_x}
\end{align*}
Now note the following:
\begin{align*}
 \sum_{k\in\Z} 2^{-k}\|P_k Q_{\le k+C_2} \phi^n\|^2_{L^6_{t,x}} &\les \sum_{k\in\Z}2^{-k}   \|P_k Q_{\le k+C_2'}\phi^{n,d} \|^2_{L^{6}_{t,x}} \les \sum_{k\in\Z}
 2^{-k} \|P_k \psi^{n,d} \|^2_{L^{6}_{t,x}} \\
& \les \sum_{k\in\Z}
  \|P_k \psi^{n,d} \|^2_{S[k]} = \|\psi^{n,d}\|_{S}^2\les M^2
\end{align*}
see above. On the other hand, the elliptic piece satisfies
\begin{align*}
   \|P_k Q_{> k+C_2} \phi^n\|_{L^6_{t,x}} &\les  2^{-{C_2}{10}} 2^{\frac{k}{2}} \sum_{\ell\in\Z} 2^{-\frac14|\ell-k|} \|P_\ell\psi^n\|_{S[\ell]}
\end{align*}
via the same arguments we used in the elliptic case earlier in this proof.
The remaining cases $i_1=2, i_2=1$, $i_{1,2}=$, are treated similarly.
This now concludes the proof of~\eqref{eq:psingoal}, and therefore of the proposition.
\end{proof}

Next, we formulate the analogue of Proposition~4.11 in our context.

\begin{prop}
 \label{prop:4.11} Let $I^+=[0,\infty)$ and assume that $\lambda(t)>\lambda_0>0$ for all $t\ge0$. Then
for $i=1,2$,
\[
 \int_{\R^2} \la \del_i U, \del_t U\ra\, dx = \Re\int_{\R^2} \Psi^\infty_i \bar{\Psi}_0^\infty\, dx = 0
\]
for all times in~$I^+$.
\end{prop}
\begin{proof}
 In view of Proposition~\ref{prop:4.10} we may also assume that $I^-=(-\infty,0]$. For a contradiction,
assume that
\[
 \Re\int_{\R^2} \Psi^\infty_1 \bar{\Psi}_0^\infty\, dx =\gamma>0
\]
As in~\cite{KeM1} one now obtains the following statements,
cf.~(4.10) and~(4.11) in~\cite{KeM1}:
\begin{itemize}
  \item Given $\eps>0$ there exists $R_0(\eps)>0$ so that for all
  $t\ge0$ one has
  \begin{equation}
    \label{eq:me1} \int_{\big|x+\frac{\bar{x}(t)}{\lambda(t)}\big|\ge
    R_0(\eps)} |\del_\alpha \Psi_\alpha^\infty(t,x)|^2\, dx \le \eps
  \end{equation}
  \item There exists $M>0$ so that for all $t\ge0$, one has
  $\Big|\frac{\bar{x}(t)}{\lambda(t)}\Big|\le t+M$
\end{itemize}
These are a consequence of the compactness in
Corollary~\ref{cor:compactV} and Lemma~\ref{lem:finitespeed}. Recall from the proof of Proposition~\ref{WeakWM} that upon passing to
a suitable subsequence of the approximating maps $\bfu^n$, we may extract an $L^2$-limit for the standard derivative components $\phi^n_\alpha$; denote this by $\Phi^\infty_\alpha$ (which, in contrast to $\Psi^\infty_\alpha$, we do not claim to be canonical).
Now
define for each $d>0$, $R>0$,
\[
Z_\alpha^{d,R}(t,x):=
\Phi_\alpha^{\infty,R}\Big(\frac{t-dx_1}{\sqrt{1-d^2}},
\frac{x_1-dt}{\sqrt{1-d^2}}, x_2\Big)
\]
where
\[
\Phi_\alpha^{\infty,R}(s,y):= R\Phi_\alpha^{\infty}(Rs, Ry)
\]
These rescaled limiting profiles again have energy~$\Ecrit$. Now
define $\theta$ to be a smooth cutoff function supported on $|x|\le
2$ and $\theta=1$ on~$|x|\le1$. The main calculation in the proof of
Proposition~4.11 of~\cite{KeM1} now reveals that, see (4.20) there,
uniformly in $t_0\in[1,2]$,
\begin{equation}
  \label{eq:energydecrease}
\sum_{\alpha=0}^2 \int_{\R^2} \theta^2(x)
|Z_\alpha^{d,R}(t_0,x)|^2\, dx = \Ecrit -\gamma d + d\eta(R,d) +
\tilde \eta(R,d) + O(d^2)
\end{equation}
with $\eta(R,d)$ and $\tilde\eta(R,d)\to0$ as $R\to\infty$,
uniformly in $0<d<d_0$ and with $O(d^2)$ uniform in~$R$. Furthermore, the argument in ~\cite{KeM1} yields that for fixed $\eps>0$, $R>0$, $d>0$ as above, one may find $t_0\in [1,2]$ such that
\[
\int_{\frac{1}{2}\le |x| \le 2}|Z_\alpha^{d,R}(t_0,x)|^2\, dx \le \eps
\]
We shall later pick $\epsilon$, $R$ depending on $\gamma, d$ and $d$ depending on $\gamma, \Ecrit$. Now for fixed choices of these parameters, pick $n$ large enough such that for $\bfu^n=(\bfx^n, \bfy^n)$ an element of the approximating sequence of wave maps from $\R^{2+1}\to\Hyp^2$, denoting by $\psi^{n,d, R}_\alpha$ the Coulomb components of $\bfu^n\circ L_d$ dilated by factor $R$ as above, and similarly by $\phi^{n,d, R}_\alpha$ the standard derivative components, an averaging argument over different timelike foliations yields that we may also assume
\[
\int_{\R^{2}}|\phi^{n,d, R}_\alpha(t_0, x)-Z_\alpha^{d,R}(t_0, x)|^{2}\,dx<\eps.
\]
Note that now $t_0$ may depend on $n$, but this does not affect the
argument. The idea now is to truncate the data \[(\bfu^n\circ L_d(R
t_0, R x),\quad R\partial_t\bfu^n\circ L_d(R t_0, R x))\] solve the
Cauchy problem backwards, and undo the Lorentz transform. We thereby
obtain a good approximation to the original essentially singular
sequence $\psi^n_\alpha$, but which satisfies good $S$-estimates,
which gives us the desired  contradiction. Thus, write $\bfu^n\circ
L_d(R t, R x)= (\bfx^{n, d, R}, \bfy^{n, d, R})$. To do this, we
consider data
\[
h^{n, d, R}(t_0,\cdot):=\Big(\chi_{[|x|<\frac{1}{2}]}\frac{{\bfx^{n, d, R}}(t_0,\cdot)-\bfx^{n, d, R}_{0}}{{\bfy}^{n, d, R}_{0}},\:
e^{\chi_{[|x|<1]}\log[\frac{\bfy^{n, d, R}}{{\bfy}^{n, d, R}_{0}}(t_0,\cdot)]}
\Big),
\]
where  $\chi_{[|x|>M]}$ is a smooth cutoff to the
set~$\{|x|>M\}$ which equals one on
$\{|x|>\frac54 M\}$, say, and $\chi_{[|x|<M]}:=1-\chi_{[|x|>M]}$. Moreover,
\begin{align*}
\bfx_0^{n,d,R}&:= \slashint_{[\frac14<|x|<\frac12]} \bfx^{n,d,R}(x)\, dx_1dx_2,\qquad
 \bfy_0^{n,d,R} :=
\exp\big(\slashint_{[\frac{1}{2}<|x|< 1]} \log \bfy^{n,d,R}(x)\, dx_1dx_2\big)
\end{align*}
Also, denote by $h^{n, d, R}(t,\cdot)$ the above expressions with $t_0$ replaced by $t$.
As in the proof of Lemma~\ref{lem:finitespeed}, one then checks that for these data we have
\[
\int e(h^{n, d, R})(t_0,\cdot)\,dx<\Ecrit-\frac{\gamma d}{2}
\]
where $e$ is the energy density,
provided we choose $R$ large enough, $\eps$ and $d$ small enough, and then $n$ large enough. Now consider the wave maps evolution of the data
\[
H^{n, d, R}(t_0,\cdot):= \big(h^{n, d, R}(t_0,\cdot),\quad
\partial_th^{n, d, R}(t_0,\cdot)\big)
\]
Our energy induction hypothesis implies that this evolution is defined globally in time, and  upon denoting the corresponding Coulomb derivative components by
\[
 \psi^{n,d, R}_{\chi, \alpha},
\]
we obtain a global bound
\[
 \|\psi^{n,d, R}_{\chi, \alpha}\|_{S(\R^{2+1})}\le \Lambda(\Ecrit,\,d,\,\gamma)<\infty
\]
Denote the time evolution of the data $H^{n, d, R}(t_0,\cdot)$ by
$H^{n, d, R}(t,\cdot)$, and the corresponding derivative components
(not in the Coulomb Gauge) by
\[
 \phi^{n,d, R}_{\chi, \alpha}
\]
We now undo the Lorentz transformation $L_d$, i.e., consider
\[
 h^{n,d, -d, R}(t,\cdot):=h^{n, d, R}(t,\cdot)\circ L_{-d}
\]
The argument in the proof of the preceding proposition then yields that we also can conclude that
the Coulomb derivative components of $h^{n,d,-d, R}(t,\cdot)$, which we denote by $\psi^{n,d,-d, R}_{\chi, \alpha}$, also satisfy a bound of the form
\[
 \|\psi^{n,d,-d,R}_{\chi, \alpha}\|_{S(\R^{2+1})}\le \Lambda'(\Ecrit,\,d,\,\gamma)<\infty
\]
Furthermore, denoting the standard derivative components of $h^{n,d,-d, R}(t,\cdot)$ by $\phi^{n,d,-d,R}_{\chi, \alpha}$, by finite propagation speed we have
\[
 \phi^{n,d,-d,R}_{\chi, \alpha}(0,x)=\phi^{n, R}_\alpha(0,x),\quad \alpha=0,\,1,\,2,
\]
provided $|x|<\frac{1}{10}$, say, where $\phi^{n,
R}_\alpha(0,\cdot)$ are the standard derivative components of
$\bfu^n(R t,\,R x)$ at time $t=0$. To conclude the proof of the
proposition, we note  that by the convergence of the $\psi^n_\alpha$
at time $t=0$ in the $L^{2}$-sense, picking $R$ large enough and
then also $n$ large enough, we may arrange that (for suitable
constants $\gamma_{nm}\in\R$)
\[
 \|\phi^{n,d,-d,R}_{\chi, \alpha}(0,\cdot)-e^{i\gamma_{nm}}\psi^{m}_\alpha(0,\cdot)\|_{L_x^2}\le \eps_1,\quad m\geq n
\]
where $\eps_1$ is as in Proposition~\ref{prop:ener_stable}, with $A=\Lambda'(\Ecrit,\,d,\,\gamma)$. But this then yields the contradiction
\[
 \limsup_{m\to\infty}\|\psi^{m}_\alpha\|_{S(\R^{2+1})}<\infty
\]
 and we are done.
\end{proof}

\subsubsection{Rigidity I: harmonic maps and reduction to the self-similar case}

As in \cite{KeM1} one now has the following rigidity theorem.

\begin{prop}
 \label{prop:rigidI}
With $\{\Psi_\alpha^\infty\}_{\alpha=0}^2$ as above, and with life span $(-T_0,T_1)$ one cannot have $T_1$ or $T_0$ finite.
Moreover, if $\lambda(t)\ge \lambda_0>0$ for all $t\in\R$, one necessarily has $\Psi_\alpha^\infty=0$ for $\alpha=0,1,2$.
\end{prop}

The proof of it will follow from a sequence of lemmas, and only be
completed after Proposition~\ref{prop:nonexist}. We begin with the
case where $T_1=\infty$ and $\lambda(t)\ge\lambda_0>0$ on
$[0,\infty)$. Assuming that $\Psi^\infty_\alpha$ do not all vanish,
the logic then is to extract a nonconstant harmonic map of finite
energy into the compact Riemann surface~$\calS$, leading to a
contradiction. The following lemma is the analogue of Lemma~5.4
in~\cite{KeM1}. While the statement is identical with that
in~\cite{KeM1}, its proof is slightly different and invokes in a
crucial way the geometry of the target.

\begin{lemma}
 \label{lem:5.4}
There exists $\eps_1>0$, $C>0$ such that if $\eps\in (0,\eps_1)$ there exists $R_0(\eps)$ so that if $R>2R_0(\eps)$ then there
exists $t_0=t_0(R,\eps)$, $0\le t_0\le CR$ with the property that for all $0<t<t_0$ one has
\[
 \Big| \frac{\bar{x}(t)}{\lambda(t)}\Big| < R-R_0(\eps),\qquad \Big| \frac{\bar{x}(t_0)}{\lambda(t_0)}\Big| = R-R_0(\eps)
\]
\end{lemma}
\begin{proof}
 As a preliminary argument, we show that there exists $\alpha\in\R$ with
\begin{equation}
 \label{eq:Iut}
\int_{I}\int_{\R^2} |\Psi_0^\infty|^2(t,x)\, dxdt \ge \alpha>0
\end{equation}
for all intervals $I$ of length one. If not,
there exists  a sequence of intervals $J_n:=[t_n, t_n+1]$ with the property
that $t_n\to\infty$ and
\begin{equation}
 \label{eq:Iutfail}
\int_{J_n}\int_{\R^2} |\Psi_0^\infty|^2(t,x)\, dxdt \le \frac{1}{n}
\end{equation}
Then there exist times $s_n\in J_n$ with the property that $\|\Psi_0^\infty(s_n,\cdot)\|_2\to0$ as $n\to\infty$.
By Corollary~\ref{cor:compactV} one has that
\[
 \Big\{ \lambda(s_n)^{-1} \Psi_\alpha^\infty\big(s_n,(\cdot-\bar{x}(s_n))\lambda(s_n)^{-1}\big) \Big\}_{n=0}^\infty
\]
forms a compact set for $\alpha=0,1,2$. Passing to a subsequence, we may assume that strongly in~$L^2$
\[
  \lambda(s_n)^{-1} \Psi_\alpha^\infty\big(s_n,(\cdot-\bar{x}(s_n))\lambda(s_n)^{-1}\big)\to \Psi_\alpha^*(\cdot)
\]
By Lemma~\ref{BasicStability} there exists some nonempty time interval $I^*$ around zero such that
\[
 \lambda(s_n)^{-1} \Psi_\alpha^\infty\big(s_n+t\lambda(s_n)^{-1},(\cdot-\bar{x}(s_n))\lambda(s_n)^{-1}\big)\to \Psi_\alpha^*(t,\cdot)
\]
in $L^\infty_{\mathrm{loc}}(I^*;L^2(\R^2))$.
Distinguish two cases: $\{\lambda(s_n)\}$ is bounded or not. In the former case, note that   $\lambda(t)\ge\lambda_0>0$ implies that
there exists a nonempty $I^\dagger\subset I^*$ such that $s_n+\lambda(s_n)^{-1} I^\dagger\subset I^*$ for each $n$. Therefore,
\eqref{eq:Iutfail} implies that
\begin{align*}
 \int_{I^\dagger} \int_{\R^2}  |\Psi_0^*|^2(t,x)\, dxdt =0
\end{align*}
This implies that $\Psi_0^*(t,\cdot)=0$ for all $t\in I^\dagger$. On the other hand, if $\{\lambda(s_n)\}$
is unbounded for every sequence $\{s_n\}$ with $s_n\in J_n$, we invoke the covering argument from \cite{Struwe1}. Thus write for each $n$
\[
 J_n=\bigcup_{s\in J_n}[s-\lambda^{-1}(s), s+\lambda^{-1}(s)]
\]
By the Vitali covering lemma, we may pick a disjoint subcollection of intervals
$\{I_s\}_{s\in A^n}$, $I_s:=[s-\lambda^{-1}(s), s+\lambda^{-1}(s)]$ for some subset $A^n\subset J_n$ with the property that
\[
 \bigcup_{s\in A^{n}}|I_s|\geq \frac{1}{5}
\]
But then the defining property of the $J_n$ implies that for each $J_n$, we may pick times $s_n\in J_n$ with the property that
\[
 \int_{I_s\cap J_n}\|\Psi^{\infty}_0(t,\,\cdot)\|_{L_x^{2}}^{2}\,dt=o(\lambda^{-1}(s_n))
\]
Alternatively, this implies that as $n\to\infty$
\[
\int_{-1}^{1}\|\big(\chi_{J_n}\Psi^{\infty}_0\big)(s_n+t\lambda^{-1}(s_n),\,\cdot)\|_{L_x^{2}}^{2}\,dt=o(1)
\]
Now pick a converging subsequence of
\[
 \lambda(s_n)^{-1}\Psi^{\infty}_0(s_n+t\lambda^{-1}(s_n),(\cdot-\bar{x}(s_n))\lambda(s_n)^{-1})
\]
to again obtain a limiting object $\Psi^*_\alpha$ with the property that
\[
\Psi^*_0(t,\cdot)=0
\]
provided $t\in I^*$, the latter its lifespan interval.
\\
We now deduce the desired contradiction from this situation: as in
Proposition~\ref{WeakWM}, we can associate a weak wave map $U^*$
from $\R^{2+1}\to\calS$ with the limiting object $\Psi^*_\alpha$,
and this wave map has the property that
\[
 \partial_t U^*=0,\quad t\in I^*
\]
Moreover, we have
\[
 \sum_{\alpha=1}^{2}\|\partial_\alpha U^*\|_{L_{x}^{2}}^{2}=\sum_{\alpha=1,2}\|\Psi^{*}_\alpha\|_{L_{x}^{2}}^{2}\neq 0
\]
We have thus obtained a nonvanishing finite energy harmonic map $U^*: \R^{2}\to \calS$, which is impossible, see \cite{SchoenYau}.

We therefore conclude that \eqref{eq:Iut} holds. The remainder of the argument is essentially the same as that in Lemma~5.4
of~\cite{KeM1}: by Corollary~\ref{cor:conservationlaw},
\begin{equation}
 \label{eq:Morawetz}
 \frac{d}{dt} \sum_{i=1}^2 \int_{\R^2} x_i \phi(x/R) \la \del_t U(t,x), \del_i U(t,x)\ra\, dx = -\int_{\R^2} |\del_t U(t,x)|^2\, dx + O(r(R))
\end{equation}
where
\[
r(R) := \int_{[|x|\ge R]} \sum_{\alpha=0}^2  |\del_\alpha U(t,x)|^2\, dx
\]
Furthermore, by compactness, for every $\eps>0$ there exists $R_0(\eps)>0$ such that for all $t\ge0$ one has
\[
 \int_{\big|x+\frac{\bar{x}(t)}{\lambda(t)}\big|\ge R_0(\eps)} |\del_\alpha U(t,x)|^2\, dx \le \eps
\]
since $\lambda(t)\ge \lambda_0>0$ for all $t\ge0$. Therefore, if the lemma were to fail, then (assuming $\bar{x}(0)=0$ as we may) one would have
\[
 \big|\frac{\bar{x}(t)}{\lambda(t)}\big|\le R-R_0(\eps)
\]
for all $0\le t < CR$. In view of the preceding, one concludes that $r(R)\le C_5 \eps $ for some absolute constant~$C_5$.
Now choose $\eps>0$ so small that
\[
 \int_{I} \Big(-\int_{\R^2} |\del_t U(t,x)|^2\, dx + O(r(R))\Big)\, dt \le -\frac{\alpha}{2}
\]
for all $I$ of unit length. In view of the apriori bound
\[
 \sup_t\Big| \int_{\R^2} x_i \phi(x/R) \la \del_t U(t,x), \del_i U(t,x)\ra\, dx \Big|\le C_6\, R\Ecrit
\]
one obtains a contradiction by integrating \eqref{eq:Morawetz} over a sufficiently large time interval.
\end{proof}

Next, we obtain a contradiction to Lemma~\ref{lem:5.4} by means of Proposition~\ref{prop:4.11}. This is
completely analogous to Lemma~5.5 in~\cite{KeM1}.

\begin{lemma}
 \label{lem:5.5}
There exists $\eps_2>0$, $R_1(\eps)>0$, $C_0>0$ such that if $R>R_1(\eps)$, $t_0=t_0(R,\eps)$ are as in Lemma~\ref{lem:5.4},
then for $0<\eps<\eps_2$ one has
\[
 t_0(R,\eps)>\frac{C_0 R}{\eps}
\]
\end{lemma}
\begin{proof}
 This follows from Proposition~\ref{prop:4.11} by the same argument as in~\cite{KeM1}.
\end{proof}

\begin{proof}[Proof of Proposition~\ref{prop:rigidI} for $T_1=\infty$] Choosing $\eps$ small in Lemma~\ref{lem:5.4}
and Lemma~\ref{lem:5.5} leads to a contradiction.
\end{proof}

It remains to prove Proposition~\ref{prop:rigidI} in case $T_1<\infty$. This will be lead to a contradiction as in~\cite{KeM1},
by a reduction to the case of a self-similar blow-up scenario.  More precisely, recall from Lemma~\ref{lem:4.7} above that
\[
 \lambda(t) \ge \frac{C_0(K)}{1-t}, \qquad 0<t<1
\]
where we assumed that $T_1=1$ as we may. Recall also that in this case
\[
 \supp(\Psi_\alpha^\infty(t,\cdot)) \subset B(0,1-t), \qquad 0<t<1
\]
see Lemma~\ref{lem:4.8}. Next, we prove an upper bound on $\lambda(t)$ which places us in the self-similar context.

\begin{lemma}
 \label{lem:5.6} Assuming that $T_1=1$ there exists a constant $C_1(K)$ such that
\[
\frac{C_1(K)}{1-t}\ge\lambda(t), \qquad 0<t<1
\]
\end{lemma}
\begin{proof} Suppose this fails.
 Let
\[
 z(t):= \sum_{j=1}^2 \int x_j \Psi_j^\infty(t,x)\bar{\Psi}_0^\infty(t,x)\, dx , \qquad 0<t<1
\]
Note that $z(t)\to 0$ as $t\to1$. Moreover, by Corollary~\ref{cor:compactV} one has
\[
 z'(t) = -\int |\Psi_0^\infty(t,x)|^2\, dx
\]
Hence,
\[
 z(t) = \int_{t}^1 \int |\Psi_0^\infty(s,x)|^2\, dx ds
\]
We now distinguish two cases: either there exists $\alpha>0$ such that
\begin{equation}
 \label{eq:caseKeM}
\int_{t}^1 \int |\Psi_0^\infty(s,x)|^2\, dx ds \ge \alpha(1-t),  \qquad 0<t<1
\end{equation}
or not, i.e., there exists a sequence $J_n=(t_n,1)$ with $t_n\to1$ such that
\begin{equation}
 \label{eq:notcaseKeM}
|J_n|^{-1}\int_{J_n} \int |\Psi_0^\infty(s,x)|^2\, dx ds \to 0\quad \text{\ as\ }n\to\infty
\end{equation}
If the first alternative~\eqref{eq:caseKeM} holds, then one is lead
to a contradiction as in~\cite{KeM1}. On the other hand, we will now
reduce the second alternative~\eqref{eq:notcaseKeM} to the existence
of a nontrivial harmonic map into~$\calS$ by a similar argument as
in the proof of Lemma~\ref{lem:5.4}, see also Struwe~\cite{Struwe1}.
By the Vitali argument from above, one selects intervals
$J_n':=(s_n-\lambda(s_n)^{-1},s_n+\lambda(s_n)^{-1})$ with $s_n\in
J_n$ such that
\[
 |J_n'|^{-1} \int_{J_n'} \int |\Psi_0^\infty(s,x)|^2\, dx ds \to 0\quad \text{\ as\ }n\to\infty
\]
Now one uses compactness as in the proof of Lemma~\ref{lem:5.4} to
conclude that there exists a limiting wave map $\Psi_\alpha^*$ on
some nonempty interval $I^*$ with $\Psi^*_0=0$ on~$I^*$. Therefore,
$\Psi^*$ leads to a  a harmonic map $U^*$ of energy~$\Ecrit$
into~$\calS$, which gives the desired contradiction.
\end{proof}

This now allows us to reduce to the exactly self-similar case.

\begin{cor}
 \label{cor:self-similar} If $T_1=1$, then the set
\[
 \Big\{ (1-t)\Psi^\infty_\alpha(t,(1-t)x)\::\: 0<t<1,\quad \alpha=0,1,2\Big\}
\]
is compact in $L^2$.
\end{cor}
\begin{proof}
 This is as in Proposition~5.7 of \cite{KeM1}.
\end{proof}

\subsubsection{Rigidity II: the self-similar case}

We now turn to the last step in the Kenig-Merle program (modulo the issue of removing the assumption $\lambda(t)>\lambda_0$
for infinite times) which consists of excluding the possibility of self-similar blow-up. As in~\cite{MerleZaag1}, \cite{MerleZaag2} we set
\[
 y=\frac{x}{1-t},\quad s=-\log(1-t),\quad 0<t<1
\]
and
\[
  W(y,s,0):=  U(x,t)= U(e^{-s}y,1-e^{-s}),\quad 0\le s<\infty
\]
where $U$ is a weak wave map as constructed in
Proposition~\ref{WeakWM}. By construction, $\nabla_{s,y}W$ is
supported in $\{|y|\le 1\}$. Next, for $\delta>0$, introduce
\[
y=\frac{x}{1-t+\delta},\quad s=-\log(1-t+\delta),\quad 0<t<1
\]

\begin{equation}\label{eq:Wdef}
W(y,s,\delta):=U(e^{-s}y,1+\delta-e^{-s})
\end{equation}
Then we have that $W(y,s,\delta)$ is defined for $0\le s<-\log\delta$ and
\[
 \supp(\del_\alpha W(\cdot,\delta))\subset\{ |y|\le 1-\delta\}
\]
The $W$ solve the equation in the distributional sense
\begin{equation}
 \label{eq:Wmain}
\del_s^2 W = \frac{1}{\rho}\mathrm{div}(\rho\nabla W-\rho(y\cdot \nabla W)y) - A(W)((\partial_{s}+y\cdot\nabla_y)W,\nabla_y W)
\end{equation}
where the nonlinearity stands for the second fundamental form on the Riemann surface~$\calS$ relative to its embedding into $\R^N$.

We now state the following properties of $W$. Henceforth, $|\cdot|$
when applied to derivatives of~$W$ will denote the metric on~$\calS$
and $W=W(\cdot,\delta)$.

\begin{lemma}
 \label{lem:6.1}  For $\delta>0$ fixed,
\begin{itemize}
 \item  $\supp(\del_\alpha W(\cdot,\delta))\subset\{ |y|\le 1-\delta\}\quad\alpha=0,1,2$
\item  $\int (|\nabla_y W|^2+|\del_s W|^2)\,dy \le C$
\item $\sum_{\alpha=0}^2 \int |\del_\alpha W(s,y)|^2\, |\log(1-|y|^2)|\,dy \le C|\log\delta|$
\item $\sum_{\alpha=0}^2 \int |\del_\alpha W(s,y)|^2\, (1-|y|^2)^{-\frac12}\,dy \le C\delta^{-\frac12}$
\end{itemize}
\end{lemma}
\begin{proof}
 By direct calculation.
\end{proof}

As in \cite{KeM1} one now introduces a Lyapunov functional
\[
 \tilde E(W(s)):= \frac12\int_{\disk} [|\del_s W|^2 + |\nabla_y W|^2 - |y\cdot \nabla_y W|^2]\, (1-|y|^2)^{-\frac12}\,dy
\]
This quantity satisfies

\begin{prop}\label{p:6.7}
For $0<s_1<s_2<\log(\frac{1}{\delta})$,
the following identities hold:
\begin{enumerate}
\item $\tilde{E}(W(s_2))-\tilde{E}(W(s_1)) = \int_{s_1}^{s_2} \int_{\disk}
\frac{|\del_s W|^2}{(1-|y|^2)^{3/2}} \, dy\, ds$
\item $\lim\limits_{s\to\log(\frac{1}{\delta})} \tilde{E}(W(s))\leq \Ecrit$.
\end{enumerate}
\end{prop}
\begin{proof}
 This is proved as in \cite{MerleZaag1}, see Lemma~2.1 there. The difference is of course that
we have a different equation, namely~\eqref{eq:Wmain}. However, the point is that the second fundamental
form is perpendicular to $\del_s W$ and $\nabla_y W$ whence it drops out of the calculation needed for the first identity.

The second property is verified as in \cite{KeM1}.
\end{proof}

As a corollary, one now has the following:

\begin{lemma}
 \label{lem:logdel} For each $\delta>0$ there exists $\bar{s}_\delta\in (1,|\log\delta|)$ such that
\[
 \int_{\bar{s}_\delta}^{\bar{s}_\delta + |\log\delta|^{\frac12}} \int_{\disk}
\frac{|\del_s W|^2}{(1-|y|^2)^{\frac32}}\, dyds \le \frac{\Ecrit}{|\log\delta|^{\frac12}}
\]
\end{lemma}
\begin{proof}
 By Proposition~\ref{p:6.7},
\[
\int_0^{|\log\delta|} \int_{\disk}
\frac{|\del_s W|^2}{(1-|y|^2)^{\frac32}}\, dyds \le \Ecrit
\]
whence the claim.
\end{proof}

The goal is now to obtain a limit $W^*$ as $\delta\to0$ and to show that $W^*$ is a stationary solution of~\eqref{eq:Wmain}.
To this end, select $\delta_j\to0$ such that for each $\alpha=0,1,2$,
\[
 (1-\bar{t}_{\delta_j})\Psi^\infty_\alpha(\bar{t}_{\delta_j},(1-\bar{t}_{\delta_j})x) \to \Psi^*_\alpha(x)
\]
strongly in $L^2$, see Corollary~\ref{cor:self-similar}. In fact, we may arrange also that
\begin{equation}
\label{eq:deltaj}
 (1+\delta_j -\bar{t}_{\delta_j})\Psi^\infty_\alpha(\bar{t}_{\delta_j},(1+\delta_j-\bar{t}_{\delta_j})x) \to \Psi^*_\alpha(x)
\end{equation}
in $L^2$.  Now consider the evolution on the level of the $\Psi$ with data given by the left-hand side of~\eqref{eq:deltaj}, see
Section~\ref{subsec:perturbprofile}. By our perturbation theory of Section~\ref{subsec:perturbprofile} we conclude from~\eqref{eq:deltaj}
 that these
evolutions $\Psi_\alpha^{j*}(t,x)$ exist on some fixed lifespan, and moreover,
\[
 \Psi_\alpha^{j*}(t,x) =  (1+\delta_j -\bar{t}_{\delta_j})\Psi^\infty_\alpha(\bar{t}_{\delta_j}+(1+\delta_j-\bar{t}_{\delta_j})t ,
 (1+\delta_j-\bar{t}_{\delta_j})x)
\]
on that lifespan $[0,T^*)$ where we may assume that $T^*<1$. Note that on account of this identity,
\[
 \supp(\Psi_\alpha^{j*}(t,\cdot))\subset \Big\{|y|\le \frac{1-\bar{t}_{\delta_j}}{1-\delta_j+\bar{t}_{\delta_j}}< 1-t\Big\}
\]
for each $\alpha=0,1,2$ and $0<t<T^*$. Now note that by the construction in the proof of Proposition~\ref{WeakWM} we may
arrange that the weak wave maps $U^{j*}$ associated with $\Psi_\alpha^{j*}$ and $U$ associated with~$\Psi_\alpha^\infty$ satisfy
\[
 U^{j*}(t, x)= U(\bar{t}_{\delta_j}+(1+\delta_j-\bar{t}_{\delta_j})t ,(1+\delta_j-\bar{t}_{\delta_j})x)
\]
Note that for fixed times $t\in(0,T^*)$ one has that $\{\del_\alpha U^{j*}(t, \cdot)\}$ form a compact set in~$L^2$ whence
the argument in the proof of Proposition~\ref{WeakWM} implies that up to passing to a subsequence
\[
 \del_\alpha U^{j*}(t, \cdot) \to \del_\alpha U^*(t,\cdot)
\]
strongly in $L^2$ uniformly on compact subintervals of time. Moreover, $U^*$ is a weak wave map and satisfies the conservation laws.
Next, we switch to the $(s,y)$ variables. Define
\[
 W_j^* (y,s) := U(\bar{t}_{\delta_j}+(1+\delta_j-\bar{t}_{\delta_j})t ,(1+\delta_j-\bar{t}_{\delta_j})x)
\]
with the same relation between $(s,y)$ and $(t,x)$ as above. Similarly, define
\[
 W^*(y,s) = U^*(y,s)
\]
Then by the preceding, uniformly in  $0\le s\le -\log(1-T^*/2)=:\tilde T$ and for $\alpha=0,1,2$,
\[
  \del_\alpha W_j^* (\cdot,s)  \to  \del_\alpha W^*(\cdot,s)
\]
in the strong $L^2$ sense. Moreover, with $W$ as in~\eqref{eq:Wdef}, one has  with $\bar{s}_{\delta_j}=-\log(1+\delta_j-\bar{t}_{\delta_j})$,
\[
 W_j^*(s,y)=W(y,\bar{s}_{\delta_j}+s,\delta_j)
\]
and therefore also
\begin{equation}
 \label{eq:6.9}
\del_\alpha W(y,\bar{s}_{\delta_j}+s,\delta_j)\to \del_\alpha W^*(\cdot,s)
\end{equation}
strongly in $L^2$ uniformly in $0\le s\le \tilde T$.  Moreover, $W^*$ is a solution of~\eqref{eq:Wmain} and
\[
 \supp(\del_\alpha W^*(s,\cdot)) \subset \{|y|\le 1\}
\]
as well as
\[
 {\mathrm{trace}}(W^*(s,\cdot)) = \const
\]
where ${\mathrm{trace}}$ is the $L^2$-trace.

\begin{lemma}
 \label{lem:6.9}
Let $W^*$ be as above. Then, $$W^*(y,s)=W^*(y)  \ \mbox{and} \ W^*\not\equiv\const.$$
\end{lemma}

\begin{proof}
With $S=-\log(1-\tilde{T})$ and   $j$ large one has
\[
\int_0^S \int_{\disk} \frac{|\del_{s}W^*(y,s)|^2}{(1-|y|^2)^{3/2}}\,
dyds \leq \varliminf\limits_{j\to\infty} \int_0^S \int_{\disk}
\frac{|\del_{s}W(y,\overline{s}_{\delta_j}+s,\delta_j)|^2}{(1-|y|^2)^{3/2}}\,
dyds
\]
by \eqref{eq:6.9}. The right-hand side is bounded by
\[
\varliminf\limits_{j\to\infty} \int_{\overline{s}_{\delta_j}}^{S+\overline{s}_{\delta_j}}
\int_{\disk} \frac{|\del_{s}W(y,s,\delta_j)|^2}{(1-|y|^2)^{3/2}}\,
dyds \les \lim\limits_{j\to\infty} |\log \delta_j |^{-1/2}=0,
\]
by Lemma \ref{lem:logdel}. This shows that $W^*(y,s)=W^*(y)$
as claimed. The fact that $W^*\not\equiv \const$ follows as in~\cite{KeM1}.
\end{proof}

In other words, we have now obtained a stationary, nonconstant, distributional solution to~\eqref{eq:Wmain}
with finite energy (relative to the $y$ variable) (as well as finite $\tilde E(W^*)$). The following proposition
now leads to the desired contradiction.

\begin{prop}
 \label{prop:nonexist}
Let $W^*$ be a distributional stationary solution to
\eqref{eq:Wmain} of finite energy
\[
 \int_{\disk} |\nabla W^*(y)|^2\,dy <\infty 
\]
Then $W^*=\const$.
This thus contradicts the preceding construction of $W^*$ and completes the proof of Proposition~\ref{prop:rigidI}.
\end{prop}
\begin{proof}
We follow the argument of Shatah-Struwe, see \cite{SStruwe}: first, $W^*$ is a
weakly harmonic map from $\disk\to \calS$ where $\disk$ is equipped with the hyperbolic metric
\[
 \frac{d\rho^2}{(1-\rho^2)^2} + \frac{\rho^2}{1-\rho^2}\,d\omega^2
\]
where $(\rho,\omega)$ are polar coordinates on~$\disk$. This means that
\[
 -(\rho\sqrt{1-\rho^2}\, W^*_\rho)_\rho + \frac{\Delta_\omega W^*}{\rho\sqrt{1-\rho^2}} \perp T_{W^*} \calS
\]
Note that by Helein's theorem, this holds in the classical sense in the interior.  Integrating by parts against
$\rho\sqrt{1-\rho^2}\,W^*_\rho$ implies that
\[
 \frac{d}{d\rho}\Big( \int_{S^1} \rho^2(1-\rho^2)|W^*_\rho|^2\,d\omega - \int_{S^1} |W^*_\omega|^2\,d\omega\Big)=0
\]
and thus
\[
 \int_{S^1} \rho^2(1-\rho^2)|W^*_\rho|^2\,d\omega - \int_{S^1} |W^*_\omega|^2\,d\omega = C_0
\]
Setting $\rho=0$ one concludes that $C_0=0$ and sending $\rho\to1$ along a suitable subsequence $\rho_j$ implies that
\[
 \lim_{\rho_j\to1} \int_{S^1} |W^*_\omega(\rho_j\omega)|^2\,d\omega \to 0
\]
On the other hand, by the trace theorem, $\sup_{\frac12<\rho<1}\|W^*_\omega(\rho\omega)\|_{\dot H^{\frac12}(S^1)}\le C\|W^*\|_{H^1(\disk)}$.
Since clearly also $\sup_{\frac12<\rho<1}\|W^*(\rho\omega)\|_{L^2(S^1)}<\infty$, one concludes via interpolation that
${\mathrm{trace}}(W^*)=\const$ as the $L^2$ trace on~$S^1$.
The change of variables
\[
 \sigma(\rho)=\exp\Big(-\int_{\rho}^1 \frac{du}{u\sqrt{1-u^2}}\Big)
\]
provides a conformal equivalence between the hyperbolic disk and the disk~$\disk$ with the Euclidean metric.
In fact,
\[
 d\sigma^2+\sigma^2 d\omega^2 = \Big(\frac{\sigma}{\rho}\Big)^2 (1-\rho^2) \Big(\frac{d\rho^2}{(1-\rho^2)^2} + \frac{\rho^2}{1-\rho^2}\,d\omega^2\Big)
\]
By the conformal invariance of the Dirichlet energy in two dimensions, it follows
that $v(\sigma,\omega):=W^*(\rho,\omega)$ is a weakly harmonic map $\disk\to\calS$ with the Euclidean disk~$\disk$. Moreover, one checks that
$v$ has finite $\dot H^1$ energy relative to the $(\sigma,\omega)$-coordinates and that ${\mathrm{trace}}(v)=\const$ in this setting as well.
By a result of Qing~\cite{Quing}, it follows that $v$ is $C^\infty$ on~$\bar\disk$.
And then the result of Lemaire~\cite{Lemaire} gives the desired conclusion that $W^*=\const$.
\end{proof}

The only remaining case is to show that $\lambda(t)$ does not approach zero along some subsequence.
This case is handled as in~\cite{KeM1} or \cite{Merle}. We follow the argument~\cite{KeM1} essential verbatim.

\begin{lemma} Let $\Psi_\alpha^\infty$ be the limiting object as above and suppose that $T_1=\infty$.
Then $\lambda(t)>\lambda_0>0$ for all $t\ge0$.
\end{lemma}
\begin{proof}
Suppose this fails. Then there exist  $t_n\to\infty$ so that $\lambda(t_n)\to
0$; in fact, one may assume even that
\[
\lambda(t_n) \leq \inf\limits_{t\in[0,t_n]}\lambda(t).
\]
From Corollary~\ref{cor:compactV} one has
\[
 \Psi_\alpha^{n}:= \lambda(t_n)^{-1}\Psi_\alpha^\infty\big(t_n,(\cdot-\bar{x}(t_n))\lambda(t_n)^{-1}\big) \to \Psi_\alpha^\dagger
\]
strongly in~$L^2$.  Then $E(\Psi^\dagger)=\Ecrit$ and we may assume that the lifespan $(-T_0^\dagger,T_1^\dagger)$ of $\Psi_\alpha^\dagger$ has the
property that $T_0^\dagger<\infty$. Otherwise one obtains a contradiction from Proposition~\ref{prop:rigidI}. Now define
$\Psi_\alpha^{n}(\tau,x)$ and $\Psi_\alpha^\dagger(\tau,x)$ to be the evolutions of $\Psi_\alpha^{n}$ and~$\Psi_\alpha^\dagger$.
By the perturbation theory of Section~\ref{subsec:perturbprofile} we conclude that $\liminf_{n\to\infty} T_0(\Psi_\alpha^n)=\infty$ and
\[
 \Psi_\alpha^n(\tau,x) \to \Psi_\alpha^\dagger(\tau,x)
\]
in $L^\infty_{\mathrm{loc}}((-\infty,0]\times L^2)$. By uniqueness of the $\Psi$-evolutions
\[
 \Psi_\alpha^n(\tau,x) = \lambda(t_n)^{-1} \Psi_\alpha^\infty(t_n+\tau\lambda(t_n)^{-1},(x-\bar{x}(t_n))\lambda(t_n)^{-1} )
\]
for all $0\le t_n+\frac{\tau}{\lambda(t_n)}$.
We claim that $\tau_n:= -t_n\lambda(t_n)$ satisfies
\[
\varliminf\limits_{n} (-\tau_n) = \infty
\]
so that for all $\tau\in(-\infty,0]$, for $n$ large, $0\leq
t_n + \frac{\tau}{\lambda(t_n)} \leq t_n$. In fact, if $-\tau_n\to-\tau_0<\infty$,
then
\[
\Psi^n_\alpha(x,-\tau_n) =  \lambda(\tau_n)^{-1}  \Psi_\alpha^\infty \big(
\frac{x-x(t_n)}{\lambda(t_n)}, 0 \big)
\]
would converge to $\Psi_\alpha^\dagger (x,-\tau_0)$ in $L^2$, with
$\lambda(t_n)\to 0$, which   contradicts  $\Psi_\alpha^\dagger\ne0$.

We now make the further claim that $\|\Psi_\alpha^\dagger\|_{S(-\infty,0)}=+\infty$.
Otherwise, by the perturbation theory of Section~\ref{subsec:perturbprofile}
 for $n$ large, $T_0(\Psi_\alpha^n)=\infty$
and $\|\Psi_\alpha^n\|_{S(-\infty,0)}\leq M$, uniformly in $n$, which
 contradicts our assumption that $\|\Psi_\alpha^\infty\|_{S(0,+\infty)}=+\infty$. This is on account of Corollary~\ref{cor:lifespan},
 since for every interval $[0,\tilde{\tau}]$, one may find $[-\tilde{\tau}_1, 0]$ with the property that the map
 $\tau\to t_n+\frac{\tau}{\lambda(t_n)}$ takes the
latter interval into the former.

Now fix $\tau\in (-\infty,0]$. Then for $n$ sufficiently large,
$t_n+\frac{\tau}{\lambda(t_n)}\geq 0$ and $\lambda(t_n+\frac{\tau}{\lambda(t_n)})$ is defined. Let
\begin{align*}
& \lambda(t_n+\frac{\tau}{\lambda(t_n)})^{-1}
 \Psi_\alpha^\infty\Big( \frac{x-x(t_n+\frac{\tau}{\lambda(t_n)})}{\lambda(t_n+\frac{\tau}{\lambda(t_n)})},
t_n+\frac{\tau}{\lambda(t_n)} \Big)  = \widetilde{\lambda}_n(\tau)^{-1}
 \Psi_\alpha^n \left( \frac{x-\widetilde{x}_n(\tau)}{\widetilde{\lambda}_n(\tau)},
\tau \right) \in K,
\end{align*}
with
 \begin{equation}\label{eq:7.2}
\widetilde{\lambda}_n(\tau) = \frac{\lambda(t_n+\frac{\tau}{\lambda(t_n)})}{\lambda(t_n)} \geq 1, \
\widetilde{x}_n(\tau) = x(t_n+\frac{\tau}{\lambda(t_n)}) -
\frac{x(t_n)}{\widetilde{\lambda}_n(\tau)}.
 \end{equation}
Now, since
$
\lambda_n^{-1}f\left( \frac{x-x_n}{\lambda_n} \right)
\underset{n\to\infty}{\xrightarrow{\ \ \ \ \ }}f$ strongly in~$L^2$
with either $\lambda_n\to 0$ or $+\infty$, or $|x_n|\to\infty$
implies that $f\equiv 0$, we see that we can assume, after passing to a subsequence, that
$\widetilde{\lambda}_n(\tau)\to\widetilde{\lambda}(\tau)$,
$1\leq\widetilde{\lambda}(\tau)<\infty$ and $\widetilde{x}_n(\tau)
\to\widetilde{x}(\tau)\in\R^2$. This implies that
\[
\widetilde{\lambda}(\tau)^{-1}
 \Psi_\alpha^\dagger \left( \frac{x-\widetilde{x}(\tau)}{\widetilde{\lambda}(\tau)},
\tau \right) \in \overline{K}.
\]
Hence, by Proposition~\ref{prop:4.11} and~\ref{prop:rigidI}, $\Psi_\alpha^\dagger=0$,
which is  a contradiction.
\end{proof}

\begin{proof}[Proof of Theorem~\ref{thm:main}]
We first address global existence and regularity and the global control of the $S$-norms. In fact,
instead of~\eqref{eq:Sbound} we of course require the stronger
\[
 \|\Psi_\alpha\|_{S}\le K(\Ecrit)
\]
from which \eqref{eq:Sbound} then follows by standard Littlewood-Paley calculus and the Strichartz
component of the~$S$-norm. Assume that this strengthened assertion of the theorem fails.
Recall that $\Ecrit$ was defined as the smallest energy with the property that there exists an essentially singular
sequence of admissible maps at energy $\Ecrit$. In other words, there exists a sequence $\{\bfu^n\}_{n=1}^\infty$ of
admissible wave maps $(-T_0^n,T_1^n)\times\R^2\to \Hyp^2$ with associated gauged derivative components $\{\psi_\alpha^n\}_{n=1}^\infty$
and such that
\begin{itemize}
 \item $E(\bfu^n)\to \Ecrit$
\item $\max_{\alpha=0,1,2}\|\psi_\alpha^n\|_{S((-T_0^n,T_1^n)\times\R^2)}\to\infty$
\end{itemize}
as $n\to\infty$.  The Bahouri-Gerard decomposition of Section~\ref{sec:BG} together with the Kenig-Merle argument
of this section now lead to a contradiction whence such an essentially singular sequence cannot exist.
This now gives the result, at least up to the scattering statement. As for the latter, we argue as follows.
It suffices to carry this out for~$\Hyp^2$. Then by applying Lemma~\ref{lem:LocalSplitting} we may represent
the gauged derivative components $\psi$
for any $\delta>0$ in the form
\[
 \psi = \psi_{L}^{(\delta)}+\psi_{NL}^{(\delta)}
\]
on a time interval of the form $(T_0,\infty)$ where
$\|\psi_{NL}^{(\delta)}\|_{\ener}<\delta$ and $\psi_L^{(\delta)}$ is
a free wave. The scattering for the free wave is automatic, and the
$\psi_{NL}^{(\delta)}$ error can be iterated away.
\end{proof}

\section{Appendix}

\subsection{Completing the proof of Lemma~\ref{lem:LocalSplitting}}

We need to show, see \eqref{eq:Ijclaim}, that there exist time intervals $I_j$, $j=1,2,\ldots, M_1$,
with $M_1$ only depending on $\|\psi\|_{S}, \eps_0$,
with the property that
\begin{equation}\label{eq:Ijclaim2}
\max_{1\le j\le
M_{1}}\,\sum_{\ell\in\Z}\|P_{\ell}F_{\alpha}(\psi)\|_{N[\ell](I_{j}\times\R^{2})}^{2}<\eps_0
C_0^6
\end{equation}
Here we need to verify this for $F_\alpha$ of at least quintic
degree. In fact, the verification of this is more or less the same
for all the higher order terms, and we explain it in detail for a
quintic term of first type.  From the discussion at the end of
Section~\ref{sec:quintic} we see that we may assume the expression
to be reduced. Thus consider for example the expression
\[
\nabla_{x,t}[P_{k_0}\psi_{0}\nabla^{-1}P_{r_{1}}(P_{k_{1}}\psi_{1}
\nabla^{-1}P_{r_{2}}(P_{k_{2}}\psi_{2}\nabla^{-1}P_{r_{3}}Q_{\nu k}(P_{k_{3}}\psi_{3},P_{k_{4}}\psi_{4})))]
\]
From Lemma~\ref{quintilinear1} we infer that
\begin{align*}
\|\nabla_{x,t}P_k[P_{k_0}\psi_{0}\nabla^{-1}P_{r_{1}}(P_{k_{1}}\psi_{1} \nabla^{-1}
P_{r_{2}}(P_{k_{2}}\psi_{2}\nabla^{-1}P_{r_{3}}Q_{\nu k}(P_{k_{3}}\psi_{3},P_{k_{4}}\psi_{4})))]\|_{N[k]}\lesssim 2^{-\delta k_0}\|P_{k_0}\psi_0\|_{S[k_0]}
\end{align*}
It then follows upon square summing over all $k\in \Z$ that the
contribution from those expressions with $k_{0}\gg k$ in the sense
that $k_0-k>C(\|\psi\|_{S}, \eps_0)$ may be bounded by $\ll
\eps_0\|\psi_0\|_{S}$. In fact, similar reasoning allows us to
reduce to the case when $r_1<k_0+O(1)$, $k=k_0+O(1)$, where the
implied constant $O(1)$ may of course be quite large depending on
$\|\psi\|_{S}$ and $\eps_0$, and furthermore we may assume that
$k_{i}=r_j+O(1)$, $i=1,2,3,4$, $j=1,2,3$.   The proof of
Lemma~\ref{quintilinear1} also implies  that we may assume all
inputs other than the ones of the null-form $Q_{\nu
k}(P_{k_{3}}\psi_{3},P_{k_{4}}\psi_{4})$ to be essentially in the
hyperbolic regime, i.e., we may replace $P_{k_{j}}\psi_{j}$ by
$P_{k_{j}}Q_{<k_j+O(1)}\psi_{j}$, $j=1,2,3$, with $O(1)$ as before.
Now assume at least one of the inputs of $Q_{\nu
k}(P_{k_{3}}\psi_{3},P_{k_{4}}\psi_{4})$ is of elliptic type, in the
sense that the difference between its modulation and frequency is
large enough. W. l. o. g. write this as
\[
\nabla_{x,t}P_k[P_{k_0}\psi_{0}\nabla^{-1}P_{r_{1}}(P_{k_{1}}\psi_{1} \nabla^{-1}P_{r_{2}}(P_{k_{2}}\psi_{2}\nabla^{-1}
P_{r_{3}}Q_{\nu k}(P_{k_{3}}Q_{>k_3+C}\psi_{3},P_{k_{4}}\psi_{4})))]
\]
where the implied constant $C$ is large enough, depending on $\|\psi\|_{S}, \eps_0$. Then if we write
\[
P_{k_{3}}Q_{>k_3+C}\psi_{3}=P_{k_{3}}Q_{[k_3+C, k_0+10]}\psi_{3}+P_{k_{3}}Q_{>k_0+10}\psi_{3},
\]
the contribution of the first term on the right is seen to be very
small, by placing the output into either
$\dot{X}_{k_0}^{-1,-\frac{1}{2},1}$ or $L_t^1\dot{H}^{-1}$. On the
other hand, consider now the contribution of the second term on the
right. Here one places the output into
$\dot{X}_{k_0}^{-\frac{1}{2}+\eps, -1-\eps, 2}$ provided the output
is in the elliptic regime, or else into $L_{1}^1\dot{H}^{-1}$. In
either case, one verifies that provided $r_1<-C$ is sufficiently
negative, the contribution is small in the above sense. Hence assume
now that $r_1=O(1)$ (which again means an interval depending on
$\|\psi\|_{S}$ as well as $\eps_0$), and as before
$P_{k_3}\psi_3=P_{k_{3}}Q_{>k_0+10}\psi_{3}$. Then we may replace
$P_{k_4}\psi_4$ by $P_{k_4}Q_{<k_4+\frac{C}{2}}\psi_4$, as otherwise
it is again straightforward to see that we gain smallness. Hence we
have now reduced to estimating
\[
\nabla_{x,t}P_k[P_{k_0}\psi_{0}\nabla^{-1}P_{r_{1}}(P_{k_{1}}\psi_{1} \nabla^{-1}P_{r_{2}}(P_{k_{2}}\psi_{2}\nabla^{-1}
P_{r_{3}}Q_{\nu k}(P_{k_{3}}Q_{>k_3+C}\psi_{3},P_{k_{4}}Q_{<k_4+\frac{C}{2}}\psi_{4})))],
\]
but where now $k_{j}=r_i+O(1)$ for all $i, j$, and the output
inherits the modulation from the large modulation term
$P_{k_{3}}Q_{>k_3+C}\psi_{3}$, provided we dyadically localize the
latter. But then a straightforward argument using the
``fungibility'' of $L_{t,x}^{2}$ reveals that we may pick intervals
$\{I_{j}\}_{j=1}^{M_1}$ with $M_1=M_1(\|\psi\|_{S}, \eps_0)$ such
that
\[
\sum_{k_0\in\Z}\|\nabla_{x,t}\chi_{I_j}\big(P_k[P_{k_0}\psi_{0}\nabla^{-1}P_{r_{1}}(P_{k_{1}}\psi_{1} \nabla^{-1}P_{r_{2}}
(P_{k_{2}}\psi_{2}\nabla^{-1}P_{r_{3}}Q_{\nu k}(P_{k_{3}}Q_{>k_3+C}\psi_{3},P_{k_{4}}Q_{<k_4+\frac{C}{2}}\psi_{4})))]\big)\|_{N[k_0]}^{2}<\eps_0
\]

Hence we have now reduced to establishing ``fungibility'' for the
space-time frequency reduced expression (with $k_{j}=r_i+O(1)$ for
all $i, j$)
\[
\nabla_{x,t}P_k[P_{k_0}\psi_{0}\nabla^{-1}P_{r_{1}}(P_{k_{1}}\psi_{1} \nabla^{-1}P_{r_{2}}(P_{k_{2}}\psi_{2}\nabla^{-1}
P_{r_{3}}Q_{\nu k}(P_{k_{3}}Q_{<k_3+C}\psi_{3},P_{k_{4}}Q_{<k_4+C}\psi_{4})))]
\]
But since we may estimate this by
\begin{align*}
&\|\nabla_{x,t}P_k[P_{k_0}\psi_{0}\nabla^{-1}P_{r_{1}}(P_{k_{1}}\psi_{1} \nabla^{-1}P_{r_{2}}(P_{k_{2}}\psi_{2}\nabla^{-1}
P_{r_{3}}Q_{\nu k}(P_{k_{3}}Q_{<k_3+C}\psi_{3},P_{k_{4}}Q_{<k_4+C}\psi_{4})))]\|_{N[k_0]}\\
&\lesssim \|P_{k_0}\psi_{0}\|_{S[k_0]}\|P_{k_{1}}\psi_{1}\|_{L_t^4L_x^\infty}\|P_{k_{2}}\psi_{2}\|_{L_t^4L_x^\infty}
\|\nabla^{-1}P_{r_{3}}Q_{\nu k}(P_{k_{3}}Q_{<k_3+C}\psi_{3},P_{k_{4}}Q_{<k_4+C}\psi_{4})\|_{L_{t}^{2}L_{x}^{\infty}}
\end{align*}
Then use the bound
\begin{align*}
&\|\nabla^{-1}P_{r_{3}}Q_{\nu k}(P_{k_{3}}Q_{<k_3+C}\psi_{3},P_{k_{4}}Q_{<k_4+C}\psi_{4})\|_{L_{t}^{2}L_{x}^{\infty}}\\
&\lesssim \|P_{r_{3}}Q_{\nu k}(P_{k_{3}}Q_{<k_3+C}\psi_{3},P_{k_{4}}Q_{<k_4+C}\psi_{4})\|_{L_{t}^{2}L_{x}^{2}}\\&
\lesssim 2^{\frac{r_3}{2}}\prod_{j=3,4}\|P_{k_{j}}\psi_{j}\|_{S[k_j]}
\end{align*}
which follows from Lemma 4. 16, as well as Bernstein's inequality
and our assumptions on the frequencies/modulations. But then again
using the ``fungibility'' of the space $L_{t,x}^2$, we may pick time
intervals $\{I_j\}$ as before such that
\begin{align*}
&\sum_{k_0\in \Z}\|\nabla_{x,t}\chi_{I_j}P_k[P_{k_0}\psi_{0}\nabla^{-1}P_{r_{1}}(P_{k_{1}}\psi_{1} \nabla^{-1}P_{r_{2}}
(P_{k_{2}}\psi_{2}\nabla^{-1}P_{r_{3}}Q_{\nu k}(P_{k_{3}}Q_{<k_3+C}\psi_{3},P_{k_{4}}Q_{<k_4+C}\psi_{4})))]\|_{N[k_0]}^{2}<\eps_0
\end{align*}
for all $I_j$. This furnishes the proof of claim \eqref{eq:Ijclaim2}
for the first type of quintilinear null-form. The remaining short
error terms of either first or second type are treated similarly.
For the higher order errors of long type (see the discussion at the
end of Section~\ref{sec:quintic} for the terminology), the claim
follows from Proposition~\ref{prop:HigherOrderLong} as well as the
fungibility of $L_{t,x}^{8}$.

\subsection{Completing the proof of Lemma~\ref{bootstrap2}.}

Recall the setup in the proof of Lemma~~\ref{bootstrap2}: we have a frequency envelope $c_k$ controlling the data
$\psi$ at time $t_j$. We then make the bootstrapping assumption
\[
\|P_{k}\psi\|_{S[k](I_j\times\R^2)}\leq A(C_0)c_k
\]
The time intervals $I_j$ have been chosen such that we have a clean separation
\[
\psi|_{I_j}=\psi_L^{(j)}+\psi_{NL}^{(j)}
\]
where we
\[
\sum_{k\in\Z}\|P_k\psi^{(j)}_{NL}\|_{S[k](I_j\times\R^2)}^{2}<\eps_0
\]
\[
\|\nabla_{x,t}\psi^{(j)}_L\|_{L_t^\infty\dot{H}^{-1}}\lesssim \|\psi\|_{S}^3\eps_0^{-\frac{1}{M}}
\]
for large $M$, say $M=100$. We need to check that by refining each $I_{j}$ if necessary into finitely many subintervals $J_{ji}$ such that we have
\[
\|P_{k}F^{2l+1}_\alpha(\psi)\|_{N[k](J_{ji}\times\R^2)}\ll c_k
\]
where now $l=2,3,4,5$. We outline the argument for the quintic errors of first type, the remaining ones following a similar
pattern. Thus consider the expression
\[
\sum_{k_j, r_i}\nabla_{x,t}P_k[P_{k_0}\psi\nabla^{-1}P_{r_{1}}(P_{k_{1}}\psi \nabla^{-1}P_{r_{2}}(P_{k_{2}}\psi\nabla^{-1}
P_{r_{3}}Q_{\nu k}(P_{k_{3}}\psi,P_{k_{4}}\psi)))]
\]
By picking $M$ large enough, it is clear that the only contribution
that matters is when we replace each factor $P_{k_j}\psi$,
$j=1,2,3,4$, by $P_{k_j}\psi_L$. However, we note here in passing
that one can also handle interactions of $\psi_L$ and $\psi_NL$
terms with at least  factors $\psi_L$ present by means of the type
of ``fungibility'' argument to follow. Hence consider now
\[
\sum_{k_j, r_i}\nabla_{x,t}P_k[P_{k_0}\psi\nabla^{-1}P_{r_{1}}(P_{k_{1}}\psi_L \nabla^{-1}P_{r_{2}}(P_{k_{2}}
\psi_L\nabla^{-1}P_{r_{3}}Q_{\nu k}(P_{k_{3}}\psi_L,P_{k_{4}}\psi_L)))]
\]
Due to Proposition~\ref{quintilinear1}, it is clear that we obtain the desired bound
\[
\|\sum_{k_j, r_i}\nabla_{x,t}P_k[P_{k_0}\psi\nabla^{-1}P_{r_{1}}(P_{k_{1}}\psi_L \nabla^{-1}P_{r_{2}}
(P_{k_{2}}\psi_L\nabla^{-1}P_{r_{3}}Q_{\nu k}(P_{k_{3}}\psi_L,P_{k_{4}}\psi_L)))]\|_{N[k]}\ll c_k
\]
provided either $|k_0-k|\gg 1$, and similarly we may assume that $k_j=r_i+O(1)$ for $j=1,2,3,4$, $i=1,2,3$.
Thus we now reduce to estimating the expression where the summation is reduced to $k_0=k+O(1),  k_j=r_i+O(1)$
for $j=1,2,3,4$, $i=1,2,3$. But in this case, the same type of fungibility argument used in the immediately
preceding proof reveals that we may pick intervals $J_{ji}$ whose number depends only on $\Ecrit$ and which are independent of $k$ such that
\[
\sum_{r_3=k_3+O(1)=k_4+O(1)}\|P_{r_{3}}Q_{\nu k}(P_{k_{3}}\psi_L,P_{k_{4}}\psi_L)\|_{L_{t}^{2}\dot{H}^{-\frac{1}{2}}}^{2}\ll 1
\]
and then the same estimates as in the preceding proof reveal that
\[
\|\sum_{k_j, r_i}\nabla_{x,t}P_k[P_{k_0}\psi\nabla^{-1}P_{r_{1}}(P_{k_{1}}\psi_L \nabla^{-1}P_{r_{2}}
(P_{k_{2}}\psi_L\nabla^{-1}P_{r_{3}}Q_{\nu k}(P_{k_{3}}\psi_L,P_{k_{4}}\psi_L)))]\|_{N[k]}\ll c_k,
\]
as desired. The argument for the remaining error terms is similar.

\subsection{Completion of the proof of Lemma~\ref{lem:intervals2}}

To complete the proof, we need to show that the contributions of the $\chi$-factors when implementing the
Hodge decomposition for the factors of $|\nabla|^{-1}(\psi^{2})$ in
\[
\sum_{k\in\Z}\|P_k(\psi|\nabla|^{-1}(\psi^{2}))\|_{L_t^2\dot{H}^{-\frac{1}{2}}}^2
\]
is also controllable in terms of $\|\psi\|_{S}$. Using the schematic relation
\[
\chi=|\nabla|^{-1}[\psi|\nabla|^{-1}(\psi^{2})],
\]
we need to bound
\[
\sum_{k\in\Z}\|P_k(\psi|\nabla|^{-1}(|\nabla|^{-1}[\psi|\nabla|^{-1}(\psi^{2})]\psi))\|_{L_t^2\dot{H}^{-\frac{1}{2}}}^2
\]
\[
\sum_{k\in\Z}\|P_k(\psi|\nabla|^{-1}(|\nabla|^{-1}[\psi|\nabla|^{-1}(\psi^{2})]|\nabla|^{-1}[\psi|\nabla|^{-1}(\psi^{2})]))\|_{L_t^2\dot{H}^{-\frac{1}{2}}}^2
\]
We deal with the first expression, the second being treated along similar lines. Thus consider
\begin{align*}
P_k(\psi|\nabla|^{-1}(|\nabla|^{-1}[\psi|\nabla|^{-1}(\psi^{2})]\psi))=&P_k(P_{[k-10,k+10]}\psi|\nabla|^{-1}(|\nabla|^{-1}[\psi|\nabla|^{-1}(\psi^{2})]\psi))\\
&+P_k(P_{>k+10}\psi|\nabla|^{-1}(|\nabla|^{-1}[\psi|\nabla|^{-1}(\psi^{2})]\psi))\\
&+P_k(P_{<k-10}\psi|\nabla|^{-1}(|\nabla|^{-1}[\psi|\nabla|^{-1}(\psi^{2})]\psi))\\
\end{align*}
Start with the first term on the right, the high-low interactions, which we further express as
\begin{align*}
&P_k(P_{[k-10,k+10]}\psi|\nabla|^{-1}(|\nabla|^{-1}[\psi|\nabla|^{-1}(\psi^{2})]\psi))\\
&=\sum_{r<k+15}P_k(P_{[k-10,k+10]}\psi|\nabla|^{-1}P_r(|\nabla|^{-1}[\psi|\nabla|^{-1}(\psi^{2})]\psi))\\
\end{align*}
Now assume the most delicate case, in which we have a high-high-low scenario inside the expression
\[
|\nabla|^{-1}P_r(|\nabla|^{-1}[\psi|\nabla|^{-1}(\psi^{2})]\psi)
\]
with respect to the factors $|\nabla|^{-1}[\psi|\nabla|^{-1}(\psi^{2})]$, $\psi$. Thus in this case we can write
\begin{align*}
&|\nabla|^{-1}P_r(|\nabla|^{-1}[\psi|\nabla|^{-1}(\psi^{2})]\psi)\\
&=\sum_{r_1=r_2+O(1)>r+O(1)}|\nabla|^{-1}P_r(|\nabla|^{-1}P_{r_1}[\psi|\nabla|^{-1}(\psi^{2})] P_{r_2}\psi)\\
&=\sum_{r_1=r_2+O(1)>r+O(1)}|\nabla|^{-1}P_r(|\nabla|^{-1}P_{r_1}[\psi|\nabla|^{-1}P_{<r}(\psi^{2})] P_{r_2}\psi)\\
&+\sum_{r_1=r_2+O(1)>r+O(1)}|\nabla|^{-1}P_r(|\nabla|^{-1}P_{r_1}[\psi|\nabla|^{-1}P_{\geq r}(\psi^{2})] P_{r_2}\psi)\\
\end{align*}
Now observe that for the first factor on the right we have the estimate
\begin{align*}
&\|\sum_{r_1=r_2+O(1)>r+O(1)}|\nabla|^{-1}P_r(|\nabla|^{-1}P_{r_1}[\psi|\nabla|^{-1}P_{<r}(\psi^{2})] P_{r_2}\psi)\|_{L_t^2L_x^\infty}\\
&=\|\sum_{r_1=r_2+O(1)>r+O(1)}\sum_{\substack{c_{1,2}\in \calD_{r_1, r-r_1}}{\text{dist}(c_1, -c_2)\lesssim 2^{r}}}
\nabla|^{-1}P_r(|\nabla|^{-1}P_{c_1}[\psi|\nabla|^{-1}P_{<r}(\psi^{2})] P_{c_2}\psi)\|_{L_t^2L_x^\infty}\\
&\lesssim 2^{-r}2^{(1-\eps)(r-r_1)}2^{\frac{r_1}{2}}\|P_{r_1}\psi\|_{S[r_1]}\|P_{r_2}\psi\|_{S[r_2]}\|P_{<r}(\psi^{2})]\|_{L_{t,x}^\infty}\\
&\lesssim  2^{(1-\eps)(r-r_1)}2^{\frac{r_1}{2}}\|P_{r_1}\psi\|_{S[r_1]}\|P_{r_2}\psi\|_{S[r_2]}\|\psi\|_{E}^{2}
\end{align*}
Hence we obtain the bound
\begin{align*}
&\|P_k(P_{[k-10,k+10]}\psi\sum_{r_1=r_2+O(1)>r+O(1)}|\nabla|^{-1}P_r(|\nabla|^{-1}P_{r_1}[\psi|\nabla|^{-1}
P_{<r}(\psi^{2})] P_{r_2}\psi))\|_{L_t^2\dot{H}^{-\frac{1}{2}}}\\
&\lesssim \sum_{r_1=r_2+O(1)>r+O(1)}2^{\frac{r_1-k}{2}} 2^{(1-\eps)(r-r_1)}2^{\frac{r_1}{2}}\|P_{[k-10,k+10]}
\psi\|_{L_t^\infty L_x^2}\|P_{r_1}\psi\|_{S[r_1]}\|P_{r_2}\psi\|_{S[r_2]}\|\psi\|_{E}^{2}
\end{align*}
If we now square this expression and sum over $k\in\Z$, it is straightforward to check that we get the upper bound
\[
\lesssim \|\psi\|_{S}^{6}\|\psi\|_{E}^{4}
\]
Next, consider the contribution of the expression
\begin{align*}
&\sum_{r_1=r_2+O(1)>r+O(1)}|\nabla|^{-1}P_r(|\nabla|^{-1}P_{r_1}[\psi|\nabla|^{-1}P_{\geq r}(\psi^{2})] P_{r_2}\psi)\\
&=\sum_{r_1=r_2+O(1)>r+O(1)}\sum_{\tilde{r}\geq r}|\nabla|^{-1}P_r(|\nabla|^{-1}P_{r_1}[\psi|\nabla|^{-1}P_{\tilde{r}}(\psi^{2})] P_{r_2}\psi)\\
&=\sum_{r_1=r_2+O(1)>r+O(1)}\sum_{\tilde{r}\geq r}\sum_{\substack{c_{1,2}\in \calD_{r_1, \tilde{r}-r_1}}
{\text{dist}(c_1, -c_2)\lesssim 2^{r}}}|\nabla|^{-1}P_r(|\nabla|^{-1}P_{c_1}[\psi|\nabla|^{-1}P_{\tilde{r}}(\psi^{2})] P_{c_2}\psi)\\
\end{align*}
Now for fixed $r, \tilde{r}, r_{1,2}$, we can estimate, using as before the improved Strichartz estimates as well as Bernstein's inequality
\begin{align*}
&\|\sum_{\substack{c_{1,2}\in \calD_{r_1, \tilde{r}-r_1}}{\text{dist}(c_1, -c_2)\lesssim 2^{r}}}|\nabla|^{-1}
P_r(|\nabla|^{-1}P_{c_1}[\psi|\nabla|^{-1}P_{\tilde{r}}(\psi^{2})] P_{c_2}\psi)\|_{L_{t}^{2}L_{x}^\infty}\\
&\lesssim 2^{(1-\eps)r}\|\sum_{\substack{c_{1,2}\in \calD_{r_1, \tilde{r}-r_1}}{\text{dist}(c_1, -c_2)
\lesssim 2^{r}}}P_r(|\nabla|^{-1}P_{c_1}[\psi|\nabla|^{-1}P_{\tilde{r}}(\psi^{2})] P_{c_2}\psi)\|_{L_{t}^{2}L_{x}^{1+}}\\
&\lesssim 2^{\frac{r_1}{2}}2^{(1-\epsilon)(\tilde{r}-r_1)}2^{-(1-\epsilon)\tilde{r}}\prod_{j=1,2}\|P_{r_j}
\psi\|_{S[r_j]}\|\psi\|_{E}^{2}\lesssim 2^{\frac{r}{2}}2^{(\frac{1}{2}-\eps)(r-\tilde{r})}\prod_{j=1,2}\|P_{r_j}\psi\|_{S[r_j]}\|\psi\|_{E}^{2},
\end{align*}
and from here the estimate continues as before. The remaining frequency interactions inside
\[
|\nabla|^{-1}P_r(|\nabla|^{-1}[\psi|\nabla|^{-1}(\psi^{2})]\psi)
\]
are handled similarly and omitted.
\\
Next, consider the case of high-high interactions, i.e.,
\begin{align*}
&P_k(P_{>k+10}\psi|\nabla|^{-1}(|\nabla|^{-1}[\psi|\nabla|^{-1}(\psi^{2})]\psi))\\
&=\sum_{k_1=r+O(1)>k+10}P_k(P_{k_1}\psi|\nabla|^{-1}P_r(|\nabla|^{-1}[\psi|\nabla|^{-1}(\psi^{2})]\psi))\\
\end{align*}
We shall again consider the most delicate case when there are high-high interactions within
\begin{align*}
&\nabla|^{-1}P_r(|\nabla|^{-1}[\psi|\nabla|^{-1}(\psi^{2})]\psi)\\
&=\sum_{r_1=r_2+O(1)>r+O(1)}\nabla|^{-1}P_r(|\nabla|^{-1}P_{r_1}[\psi|\nabla|^{-1}(\psi^{2})]P_{r_2}\psi)\\
\end{align*}
But then arguing just as above one obtains the bound
\[
\|\nabla|^{-1}P_r(|\nabla|^{-1}P_{r_1}[\psi|\nabla|^{-1}(\psi^{2})]P_{r_2}\psi)\|_{L_{t}^{2}L_{x}^{2+}}
\lesssim 2^{-(\frac{1}{2}-\eps)r}2^{(\frac{1}{2}-\eps)(r-r_1)}\prod_{j=1,2}\|P_{r_j}\psi\|_{S[r_j]}\|\psi\|_{E}^{2},
\]
and from here one obtains
\begin{align*}
&\|\sum_{k_1=r+O(1)>k+10}P_k(P_{k_1}\psi|\nabla|^{-1}P_r(|\nabla|^{-1}[\psi|\nabla|^{-1}(\psi^{2})]\psi))
\|_{L_{t}^{2}\dot{H}^{-\frac{1}{2}}}\\
&\lesssim \sum_{k_1=r+O(1)>k+10}2^{(\frac{1}{2}-\eps)(k-r)}2^{(\frac{1}{2}-\eps)(r-r_1)}\|P_{k_1}\psi\|_{S[k_1]}
\prod_{j=1,2}\|P_{r_j}\psi\|_{S[r_j]}\|\psi\|_{E}^{2}
\end{align*}
Squaring and summing over $k$ again results in the same bound as before.
\\
The case of low-high interactions, i.e.,
\[
P_k(P_{<k-10}\psi|\nabla|^{-1}(|\nabla|^{-1}[\psi|\nabla|^{-1}(\psi^{2})]\psi)),
\]
is more of the same and omitted.
\subsection{Completion of the proof of Proposition~\ref{PsiBootstrap}, part I}

Here we show how to deal with the higher order terms encountered in
the decomposition \eqref{eq:epsilonwave}, i.e., the fifth term
there. We shall again explain the method for the quintilinear terms
of first type, the remaining higher order terms being treated
similarly. Thus consider the expression
\[
\nabla_{x,t}P_k[P_{k_0}\rho_{0}\nabla^{-1}P_{r_{1}}(P_{k_{1}}\rho_{1} \nabla^{-1}P_{r_{2}}(P_{k_{2}}
\rho_{2}\nabla^{-1}P_{r_{3}}Q_{\nu k}(P_{k_{3}}\rho_{3},P_{k_{4}}\rho_{4})))]
\]
We use the letter $\rho$ here to imply either a $\psi$-factor or one of $\epsilon_{1,2}$, the the setup in the proof of
Proposition~\ref{PsiBootstrap}.  Now we distinguish between a number of cases:
\\

(1) At least one factor of both $\epsilon_1$ and $\epsilon_2$ is present. In this case, the entire expression contributes
to $\epsilon_2$, as follows from Proposition~\ref{quintilinear1}. Indeed, we can sum over all $k_j, r_j$ and then square sum
over $k\in Z$ and bound the entire expression by
\[
\lesssim \|\epsilon_2\|_{S}\|\epsilon_1\|_{S}
\]
where the implied constant only depends on $\Ecrit$. By choosing $\eps_0$, which controls $\|\epsilon_1\|_{S}$, small enough, we can bootstrap.
\\

(2) Only $\epsilon_1$ factors in addition to $\psi$-factors. First, assume that there are at least two $\epsilon_1$ factors.
If one of them is $\rho_0$, then the output inherits the frequency envelope of $\epsilon_1$ from Proposition~\ref{quintilinear1},
and the smallness follows from the presence of the extra factor $\epsilon_1$. If the first factor $\rho_0$ is a $\psi$, then we need
to show that the expression contributes to $\epsilon_2$. But this again follows from Proposition~\ref{quintilinear1}, essentially as in
 Case (1) (d) of the proof of Proposition~\ref{PsiBootstrap}.\\
Next, assume that there is only one $\epsilon_1$ factor present. If
this factor is not $\rho_0$, then the expression contributes to
$\epsilon_2$, following the same reasoning as in Case~(1), (b). Thus
assume now that we have $\rho_0=\epsilon_1$, which is the expression
\[
\nabla_{x,t}P_k[P_{k_0}\epsilon_1\nabla^{-1}P_{r_{1}}(P_{k_{1}}\psi \nabla^{-1}P_{r_{2}}(P_{k_{2}}\psi\nabla^{-1}
P_{r_{3}}Q_{\nu k}(P_{k_{3}}\psi,P_{k_{4}}\psi))]
\]
Recall from the proof of Proposition~\ref{quintilinear1} that here $\psi$ really stands for $\psi_L$ or $\psi_NL$,
but we suppress this here. What matters is that $\|\psi\|_{S}$ depends on $\Ecrit$ in a universal way independent of the
stage of the iteration in the proof. As usual we may reduce to $k_j=r_i+O(1)$, $j=1,2,3,4$, $i=1,2,3$, and $k_0=k+O(1)>r_1+O(1)$.
Furthermore, all inputs may be assumed to be in the hyperbolic regime (up to large constants only depending on $\Ecrit$). But then
the smallness can be forced by shrinking $I_j$ suitably and forcing that
\[
\sum_{r\in\Z}\|\chi_{I_j}Q_{\nu k}P_r(P_{r+O(1)}\psi,P_{r+O(1)}\psi)\|_{L_t^2\dot{H}^{-\frac{1}{2}}}^2\ll 1,
\]
see the proof of Proposition~\ref{quintilinear1}. For the higher
order errors of long type (recall the discussion in
Section~\ref{sec:quintic}), the smallness is achieved by exploiting
the ``fungibility'' of the norms $L_{t,x}^{8}$.
\\

(3) Only $\epsilon_2$ factors present in addition to factors $\psi$.
All of these terms contribute to $\epsilon_2$. If at least two
factors $\epsilon_2$ are present, we clearly obtain the desired
smallness from Proposition~\ref{quintilinear1}. Hence now assume
that only one such factor is present. If this factor is in the
position of $\rho_0$, then we obtain smallness via ``fungibility''
or $L_{t,x}^2$ as in Case~(2). If this factor is in the position of
some $\rho_j$ with $j=1,2,3,4$, one obtains smallness
 via a slightly different fungibility argument: first, reduce to the case when $\rho_0$ and one of the $\rho_j$ which represents a $\psi$
  have angular separation between their Fourier supports: to do this, consider for example
\[
\nabla_{x,t}P_k[P_{k_0}\psi\nabla^{-1}P_{r_{1}}(P_{k_{1}}\psi \nabla^{-1}P_{r_{2}}
(P_{k_{2}}\epsilon_2\nabla^{-1}P_{r_{3}}Q_{\nu k}(P_{k_{3}}\psi,P_{k_{4}}\psi))]
\]
Again we may assume that $k_j=r_i+O(1)$ for $j=1,2,3,4$, $i=1,2,3$,
and $k_0=k+O(1)>r_1+O(1)$. Here we can use the fungibility of
$L_{t,x}^{2}$ by placing $P_{r_{3}}Q_{\nu
k}(P_{k_{3}}\psi,P_{k_{4}}\psi)$ into
$L_{t}^2\dot{H}^{-\frac{1}{2}}$, see the proof of
Proposition~\ref{quintilinear1}. On the other hand, for the
expression
\[
\nabla_{x,t}P_k[P_{k_0}\psi\nabla^{-1}P_{r_{1}}(P_{k_{1}}\psi \nabla^{-1}P_{r_{2}}(P_{k_{2}}\psi\nabla^{-1}
P_{r_{3}}Q_{\nu k}(P_{k_{3}}\epsilon_2,P_{k_{4}}\psi))],
\]
one obtains smallness from the fungibility of $L_{t}^{4}L_{x}^{\infty}$, more precisely, that of
\[
\sum_{k\in\Z}\|P_k\psi\|_{L_{t}^{4}L_{x}^{\infty}}^{4}
\]

\subsection{Completion of the proof of Proposition~\ref{PsiBootstrap}, part II}.
Here we show how to obtain the bootstrap for the elliptic part of
$\epsilon$, i.e., $Q_{\geq D}\epsilon$. Recall that we solve for
$Q_{\geq D}\epsilon$ via the equation
\[
\Box Q_{\geq D}\epsilon = Q_{\geq D}[\sum_{i=1}^{5}F_\alpha^{2i+1}(\psi+\epsilon)]-Q_{\geq D}[\sum_{i=1}^{5}F_\alpha^{2i+1}(\psi)]
\]
where the $F_{\alpha}^{2i+1}$ are obtained as described in
Section~\ref{sec:hodge}. In particular, $F_\alpha^{3}(\psi)$
constitutes the trilinear null-forms. Of course the proper
interpretation of the right-hand side is that we substitute suitable
Schwartz extensions for $\psi$ and $\epsilon$ but which agree with
the actual dynamic variables on the time interval that we work on.
We start by considering the trilinear null-forms, which with the
appropriate localizations we schematically write as
\[
\nabla_{x,t}P_0Q_{\geq D}[(\psi+\epsilon)\nabla^{-1}\calN_{\nu j}(\psi+\epsilon, \psi+\epsilon)]-\nabla_{x,t}P_0Q_{\geq D}
[\psi\nabla^{-1}\calN_{\nu j}(\psi, \psi)]
\]
We need to show that we can write the above expression as the sum of two terms, which, when evaluated with respect to
$\|\cdot\|_{N[0]}$, improve the bootstrap assumption \eqref{eq:bootass}.
Now we distinguish between various cases:
\\

(1) Here we consider the trilinear terms which are schematically of the form
\[
\nabla_{x,t}P_0Q_{\geq D}[\epsilon\nabla^{-1}\calN_{\nu j}(\psi, \psi)]
\]
We decompose this into two further terms according to the type of $\epsilon$:
\\

(1a): This is the expression $\nabla_{x,t}P_0Q_{\geq
D}[\epsilon_1\nabla^{-1}\calN_{\nu j}(\psi, \psi)]$. Recalling the
fine structure of the trilinear terms described in
Section~\ref{sec:hodge}, we see that this can be decomposed into two
types of terms
\begin{align}
&\nabla_{x,t}P_0Q_{\geq D}[\epsilon_1\nabla^{-1}\calN_{\nu j}(\psi, \psi)]\nonumber\\
&=\nabla_{x,t}P_0Q_{\geq D}[\epsilon_1\nabla^{-1}\calN_{\nu j}I^{c}(\psi, \psi)]\label{eq:1a1}\\
&+\nabla_{x,t}P_0Q_{\geq D}[(R_\mu)\epsilon_1\nabla^{-1}\calN_{\nu j}I(\psi, \psi)]\label{eq:1a2}
\end{align}
where in the last term an operator $R_\mu$ may be present or not. Start with the first term on the right, which we write as
\begin{align*}
&\nabla_{x,t}P_0Q_{\geq D}[\epsilon_1\nabla^{-1}\calN_{\nu j}I^{c}(\psi, \psi)]\\
&=\sum_{k_{1,2,3}, r}\nabla_{x,t}P_0Q_{\geq D}[P_{k_1}\epsilon_1\nabla^{-1}P_r\calN_{\nu j}I^{c}(P_{k_2}\psi, P_{k_3}\psi)]\\
\end{align*}
Now the fundamental trilinear estimates in Section~\ref{sec:trilin},
see in particular~\eqref{eq:JoachimIc}, imply that under the
bootstrap assumption
\[
\|P_k\epsilon_1\|_{S[k]}\leq C_4 d_k
\]
with some $C_4=C_4(\Ecrit)$, we have
\begin{align*}
\|\sum_{|k_{1}|\gg 1, k_{2,3}, r}\nabla_{x,t}P_0Q_{\geq D}[P_{k_1}\epsilon_1\nabla^{-1}P_r\calN_{\nu j}I^{c}(P_{k_2}\psi, P_{k_3}\psi)]\|_{N[0]}\ll C_4 d_0,
\end{align*}
which is as desired. In fact, the proof of~\eqref{eq:JoachimIc}
cited above implies that one also obtains
\begin{align*}
\|\sum_{k_{1,2,3}, |r|\gg 1}\nabla_{x,t}P_0Q_{\geq D}[P_{k_1}\epsilon_1\nabla^{-1}P_r\calN_{\nu j}I^{c}(P_{k_2}\psi, P_{k_3}\psi)]\|_{N[0]}\ll C_4 d_0,
\end{align*}
and finally, again the trilinear estimates from
Section~\ref{sec:trilin} imply that we may also assume
$k_{2,3}=O(1)$ (implied constant depending on $\Ecrit$). Hence we
may assume for the present term that all frequencies are $O(1)$.
Thus we may now reduce to considering
\begin{align*}
\sum_{k_{1,2,3}+O(1)=r=O(1)}\nabla_{x,t}P_0Q_{\geq D}[P_{k_1}\epsilon_1\nabla^{-1}P_r\calN_{\nu j}I^{c}(P_{k_2}\psi, P_{k_3}\psi)]\\
\end{align*}
Now if one of the inputs of the null-form $\calN_{\nu j}I^{c}(P_{k_2}\psi, P_{k_3}\psi)$ is of elliptic type, either at
least one of $\epsilon_1$ and the other input has at least comparable modulation, or else the output inherits the modulation
from the large modulation input. In the former case, it is straightforward to obtain smallness: indeed, consider for example
\begin{align*}
\sum_{k_{1,2,3}+O(1)=r=O(1)}\sum_{l\gg 1}\nabla_{x,t}P_0Q_{\geq D}[P_{k_1}Q_{<l-10}\epsilon_1\nabla^{-1}
P_r\calN_{\nu j}I^{c}(P_{k_2}Q_{l}\psi, P_{k_3}Q_{l+O(1)}\psi)]\\
\end{align*}
We can estimate this by (using Bernstein's inequality)
\begin{align*}
&\|\nabla_{x,t}P_0Q_{\geq D}[P_{k_1}Q_{<l-10}\epsilon_1\nabla^{-1}P_r\calN_{\nu j}I^{c}(P_{k_2}Q_{l}\psi, P_{k_3}Q_{l+O(1)}\psi)]\|_{N[0]}\\
&=\|\nabla_{x,t}P_0Q_{[D, l+O(1)]}[P_{k_1}Q_{<l-10}\epsilon_1\nabla^{-1}P_r\calN_{\nu j}I^{c}(P_{k_2}Q_{l}\psi, P_{k_3}Q_{l+O(1)}\psi)]\|_{N[0]}\\
&\lesssim \sum_{D\leq j\leq l+O(1)}2^{-\eps j}2^{\frac{j}{2}}\|R_{\nu}P_{k_2}Q_{l}\psi\|_{L_{t,x}^{2}}
\|P_{k_3}Q_{l+O(1)}\psi\|_{L_{t,x}^{2}}\|P_{k_1}Q_{<l-10}\epsilon_1\|_{L_t^\infty L_x^2}\\
&\ll C_4d_0
\end{align*}
The case when $\epsilon_1$ has comparable modulation is of course similar. Hence we may assume that if one
of the inputs $P_{k_{2,3}}\psi$ is of elliptic type, the output inherits its modulation. In order to obtain smallness
in this case, we can form example use fungibility of $L_{t,x}^{2}$ by applying suitable cutoffs $\chi_{I_j}$ for which
\[
\sum_{k_{2}\in \Z}\|\chi_{I_j}R_{\nu}P_{k_{2}}Q_{\gg k_{2}}\psi\|_{L_{t,x}^{2}}^{2}\ll 1
\]
Next, assume that both inputs $P_{k_2,3}\psi$ of the null-form are
of hyperbolic type. Then using the bilinear estimates of
Section~\ref{sec:bilin}, we can estimate
\begin{align*}
&\|\sum_{k_{1,2,3}+O(1)=r=O(1)}\nabla_{x,t}P_0Q_{\geq D}[P_{k_1}\epsilon_1\nabla^{-1}
P_r\calN_{\nu j}I^{c}(P_{k_2}Q_{<k_2+O(1)}\psi, P_{k_3}Q_{<k_3+O(1)}\psi)]\|_{N[0]}\\
&\leq \|\sum_{k_{1,2,3}+O(1)=r=O(1)}\nabla_{x,t}P_0Q_{\geq D}[P_{k_1}\epsilon_1\nabla^{-1}
P_r\calN_{\nu j}I^{c}(P_{k_2}Q_{<k_2+O(1)}\psi, P_{k_3}Q_{<k_3+O(1)}\psi)]\|_{\dot{X}_0^{-\frac{1}{2}+\eps, -1-\eps, 2}}\\
&\lesssim \|P_{k_1}\epsilon_1\|_{L_t^\infty L_x^2}\|\calN_{\nu j}I^{c}(P_{k_2}Q_{<k_2+O(1)}\psi, P_{k_3}Q_{<k_3+O(1)}\psi)\|_{L_{t,x}^{2}}
\end{align*}
In this case, smallness is again forced by subdividing into suitable time intervals $I_j$ with the property that
\[
\|\chi_{I_j}\calN_{\nu j}I^{c}(P_{k_2}Q_{<k_2+O(1)}\psi, P_{k_3}Q_{<k_3+O(1)}\psi)\|_{L_{t,x}^{2}}\ll 1
\]
This completes treatment of \eqref{eq:1a1}. Next we turn to
\eqref{eq:1a2}. The same reasoning as for \eqref{eq:1a1} shows that
we may assume all frequencies $k_{1,2,3}, r$ (which we introduce in
the same fashion as before) to be of size $O(1)$. Now a technical
issue arises when the operator $R_\nu=R_0$. Indeed, in this case, it
may happen that the output inherits the modulation of the first
input $P_{k_1}R_{\mu}\epsilon_1$, and the remaining inputs
necessarily need to be placed into the energy space which is ``not
fungible''. However, this problem is somewhat artificial, since of
course the Hodge decomposition for the temporal components becomes
counterproductive for in the large modulation (elliptic) case. Thus
for the expression
\[
\sum_{k_{1,2,3}+O(1)=r=O(1)}\nabla_{x,t}P_0Q_{\geq
D}[R_0P_{k_1}Q_{\gg 1}\epsilon_1\nabla^{-1}P_r\calN_{\nu
j}I(P_{k_2}\psi, P_{k_3}\psi)],
\]
it is best to re-combine it with the term
\[
\sum_{k_{1,2,3}+O(1)=r=O(1)}\nabla_{x,t}P_0Q_{\geq
D}[R_0P_{k_1}Q_{\gg 1}\epsilon_2\nabla^{-1}P_r\calN_{\nu
j}I(P_{k_2}\psi, P_{k_3}\psi)],
\]
as well as the ``elliptic'' error $\chi_0$ coming from
\[
\epsilon_0=R_0\epsilon+\chi_0
\]
and replace it by $P_{k_1}Q_{\gg 1}\epsilon =P_{k_1}Q_{\gg
1}\epsilon_1+P_{k_1}Q_{\gg 1}\epsilon_2$. Unfortunately, we
encounter here the technical issue that the inputs $\epsilon_{1,2}$,
$\psi$ on the right-hand side are really Schwartz extensions of the
actual components beyond the time interval $I$ we work on, and hence
do not exactly satisfy the div-curl system. The way around this is
to work on a slightly smaller time interval $\tilde{I}$ obtained by
removing small intervals $I_{1,2}$ from the endpoints  of $I$ with
$I_{1,2}$ of length $\sim T_1$ with $T_1$ as in case 1 of the roof
of Proposition~\ref{PsiBootstrap}. When we restrict the source terms
to $I$, we may invoke the div-curl system for extremely elliptic
(i.e., difference of modulation and frequency very large) terms up
to negligible errors. This allows us to obtain bootstrapped bounds
for $\epsilon_{1,2}$ on $\tilde{I}$, and at the endpoints, we can
re-iterate the argument of Case 1. Then the $\epsilon_{1,2}$ on the
full interval $I$ can be re-assembled from these pieces via
partition of unity with respect to time.
\\
The preceding discussion reveals that we may as well suppress the operator $R_{\mu}$. But once this is done,
the fungibility argument used for \eqref{eq:1a1} may be repeated to give the desired smallness upon suitably restricting the time intervals.
\\

(1b): The argument for $\nabla_{x,t}P_0Q_{\geq D}[\epsilon_1\nabla^{-1}\calN_{\nu j}(\psi, \psi)]$ is exactly
the same, one square sums over the output frequencies instead.
\\

(2): Next we consider the schematically written terms of type $\nabla_{x,t}P_0Q_{\geq D}[\psi\nabla^{-1}\calN_{\nu j}(\epsilon, \psi)]$.
Again these split into two sub-types:
\\

(2a): Terms of type $\nabla_{x,t}P_0Q_{\geq
D}[\psi\nabla^{-1}\calN_{\nu j}(\epsilon_1, \psi)]$. These
contribute to $\epsilon_2$, and indeed apart from the fact that one
uses trilinear estimates from Section~\ref{sec:trilin} for elliptic
outputs, the smallness follows formally just as in Case 1 (b) (of
the proof of Proposition~\ref{PsiBootstrap}).
\\

(2b): Terms of type $\nabla_{x,t}P_0Q_{\geq
D}[\psi\nabla^{-1}\calN_{\nu j}(\epsilon_2, \psi)]$. Here one
encounters again the issue with the terms containing $R_0\epsilon_2$
and of extremely large modulation. As in (1a) above this is handled
by undoing the Hodge decomposition for these terms by restricting to
a smaller time interval, up to negligible errors. This, as well as
arguments as in (1) above, allow one to reduce to an expression of
the form
\[
\nabla_{x,t}P_0Q_{\geq D}[P_{k_1}\psi\nabla^{-1}\calN_{\nu j}(P_{k_2}\epsilon_2, P_{k_3}\psi)]
\]
where all inputs are of hyperbolic type (up to large constants
depending on the energy alone). But then the smallness can be forced
by reducing to frequency-separated inputs $P_{k_{1,3}}\psi$ (via the
estimates of Section~\ref{subsec:improvetrilin}. But then
fungibility is obtained by grouping the inputs $P_{k_{1,3}}\psi$
together and placing their product into $L_{t,x}^{2}$.
\\

(3) The remaining trilinear null-forms with elliptic output are
easier to handle, since they contain at least two factors of type
$\epsilon_{1,2}$, and hence the smallness follows simply by the
smallness assumptions on these factors (bootstrap assumptions), as
well as the trilinear estimates of Section~\ref{sec:trilin}. We omit
the details.
\\

The higher order contributions from the
\[
Q_{\geq D}[\sum_{i=2}^{5}F_\alpha^{2i+1}(\psi+\epsilon)]-Q_{\geq D}[\sum_{i=2}^{5}F_\alpha^{2i+1}(\psi)]
\]
are estimated in a similar vein and omitted.


\begin{thebibliography}{99}

\bibitem{BG} Bahouri, H., G\'erard, P. {\em High frequency approximation of
solutions to critical nonlinear wave equations.}  Amer.\ J.~Math.\
121 (1999),  no.~1, 131--175.



\bibitem{B} Bourgain, J. {\em Estimates for cone multipliers.}     Geometric aspects of functional analysis (Israel, 1992--1994),  41--60,
Oper.\ Theory Adv.\ Appl., 77, Birkh\"auser, Basel, 1995.

\bibitem{Caz} Cazenave, T.,  Shatah, J.,  Tahvildar-Zadeh, A. {\em Shadi Harmonic maps of the hyperbolic space and development
 of singularities in wave maps and Yang-Mills fields.}  Ann.\ Inst.\ H.~Poincar\'e Phys. Th\'eor.~68  (1998),  no.~3, 315--349.

\bibitem{CT1} Christodoulou, D., Tahvildar-Zadeh, A. {\em On the
asymptotic behavior of spherically symmetric wave maps.}  Duke
Math.\  J.~71  (1993),  no.~1, 31--69.

\bibitem{CT2} Christodoulou, D., Tahvildar-Zadeh, A. {\em
  On the regularity of spherically symmetric
wave maps.} Comm.\ Pure Appl.\ Math.~46 (1993), no.~7, 1041--1091.



\bibitem{CKM}
C\^{o}te, R.,  Kenig, C.,  Merle, F. {\em Scattering below critical
energy for the radial 4D Yang-Mills equation and for the 2D
corotational wave map system.}  Comm.\ Math.~Phys.\ 284 (2008),
no.~1, 203--225.

\bibitem{Duo} Duoandikoetxea, J. {\em Fourier analysis.}
Graduate Studies in Mathematics, 29. American Mathematical Society,
Providence, RI, 2001.

\bibitem{FKS} Fabes, E., Kenig, C., Serapioni, R. {\em The local regularity of solutions of degenerate elliptic equations.}
  Comm.\ Partial Differential Equations  7  (1982), no.~1, 77--116.

\bibitem{FMS} Freire, A., M\"{u}ller, S.,  Struwe, M. {\em Weak convergence of wave maps from $(1+2)$-dimensional Minkowski space to Riemannian manifolds.}
  Invent.\ Math.\  130  (1997),  no.~3, 589--617.

\bibitem{FK} Foschi, D.,  Klainerman, S.  {\em Bilinear space-time estimates for homogeneous wave equations.},
 Ann.\ Sci.\ \'Ecole Norm.\ Sup.~(4)  33  (2000), no.~2, 211--274.

\bibitem{Ger} G\'erard, P. {\em Oscillations and concentration effects in semilinear dispersive wave equations.}\
  J.\
   Funct.\ Anal.~141  (1996),  no.~1, 60--98.

\bibitem{GrotShatah} Grotowski, J., Shatah, J. {\em Geometric evolution equations in critical dimensions.}
  Calc.\ Var.\ Partial Differential Equations  30  (2007),  no.~4, 499--512.

\bibitem{KeM1} Kenig, C., Merle, F. {\em Global well-posedness, scattering and blow-up for the energy-critical focusing non-linear wave equation.}
  Acta Math.~201  (2008),  no.~2, 147--212.

\bibitem{KeM2} Kenig, C., Merle, F. {\em
Global well-posedness, scattering and blow-up for the
energy-critical, focusing, non-linear Schr\"odinger equation in the
radial case.}  Invent.\ Math.~166  (2006),  no.~3, 645--675.

\bibitem{KM3} Klainerman, S., Machedon, M. {\em On the optimal local regularity for gauge field theories.}  Differential Integral Equations  10  (1997),  no.~6, 1019--1030.

\bibitem{KM2} Klainerman, S., Machedon, M. {\em On the regularity properties of a model problem related to wave maps.}  Duke Math.\ J.~87  (1997),  no.~3, 553--589.

\bibitem{KM1} Klainerman, S., Machedon, M. {\em Smoothing estimates for null forms and applications. A celebration of John F.~Nash, Jr.}  Duke Math.\ J.~81  (1995),  no.~1, 99--133.

\bibitem{KlRod}
Klainerman, S.,  Rodnianski, I. {\em On the global regularity of
wave maps in the critical Sobolev norm.}  Internat.\ Math.\ Res.
Notices 2001, no.~13, 655--677.

\bibitem{KS1} Klainerman, S.,  Selberg, S. {\em Remark on the optimal
regularity for equations of the wave maps type.}   Comm.\ Partial
Differential Equations  22  (1997),  no.~5-6, 901--918.

\bibitem{KS2} Klainerman, S., Selberg, S.  {\em Bilinear estimates and
applications to nonlinear wave equations.}   Commun.\ Contemp.\
Math.~4  (2002),  no.~2, 223--295.

\bibitem{KlainTat} Klainerman, S.,  Tataru, D.  {\em On the optimal local regularity for Yang-Mills equations in $R\sp {4+1}$.}
  J.\ Amer.\ Math.\ Soc.~12  (1999),  no.~1, 93--116.

\bibitem{Krieger} Krieger, J. {\em  Global regularity of wave maps from $\mathbf{R}\sp {2+1}$ to $H\sp 2$. Small energy.}
  Comm.\ Math.\ Phys.\  250  (2004),  no.~3, 507--580.

\bibitem{Krieger2} Krieger, J. {\em Null-form estimates and nonlinear waves.}
  Adv.\ Differential Equations~8  (2003),  no.~10, 1193--1236.

\bibitem{Krieger3} Krieger, J. {\em Global regularity of wave maps from $R\sp {3+1}$ to surfaces.}
  Comm.\ Math.\ Phys.~238  (2003),  no.~1-2, 333--366.

\bibitem{KST} Krieger, J., Schlag, W., Tataru, D. {\em Renormalization and blow up for charge one equivariant critical wave maps.}
  Invent.\ Math.~171  (2008),  no.~3, 543--615.

\bibitem{Lemaire} Lemaire, L. {\em Applications harmoniques de surfaces riemanniennes.}   J.\ Differential Geom.~13  (1978), no.~1, 51--78.

\bibitem{Lions} Lions, P.\ L. {\em The compensated compactness
principle in the calculus of variations.  The locally compact case
I.} Ann.\ Inst.\ H.\ Poincar\'e Anal.\ Non Lineaire 1 (1984),
109--145.

\bibitem{Merle} Merle, F. {\em Existence of blow-up solutions in the
energy space for the critical generalized KdV equation},
J.\ Amer.\ Math.\ Soc.~\textbf{14} (2001), 555--578.

\bibitem{MerleZaag1} Merle, F., Zaag, H. {\em Determination of the blow-up rate for the semilinear wave equation. } Amer.\ J.\ Math.~125  (2003),  no.~5, 1147--1164.

\bibitem{MerleZaag2} Merle, F.,  Zaag, H. {\em A Liouville theorem for vector-valued nonlinear heat equations and applications.}  Math.\ Ann.\  316  (2000),
 no.~1, 103--137.

\bibitem{MS} M\'etivier, G., Schochet, S. {\em Trilinear resonant interactions of semilinear hyperbolic waves.}
  Duke Math.\ J.~95  (1998),  no.~2, 241--304.

\bibitem{Gerd} Mockenhaupt, G. {\em  A note on the cone multiplier.}  Proc.\ Amer.\ Math.\ Soc.~117  (1993),  no.~1, 145--152.

\bibitem{MSS} Mockenhaupt, G., Seeger, A., Sogge, C.\ D. {\em Wave front sets, local smoothing and Bourgain's circular maximal theorem.}
  Ann.\ of Math.~(2)  136  (1992),  no.~1, 207--218.

\bibitem{NSU}
Nahmod, A.,  Stefanov, A.,  Uhlenbeck, K. {\em On the well-posedness
of the wave map problem in high dimensions.}  Comm.\ Anal.\ Geom.\
11 (2003),  no.~1, 49--83.

\bibitem{Quing} Qing, J. {\em Boundary regularity of weakly harmonic maps from surfaces.}
  J.\ Funct.\ Anal.~114  (1993),  no.~2, 458--466.

\bibitem{RodSterb} Rodnianski, I., Sterbenz, J. {\em On the Formation of Singularities in the Critical $O(3)$ Sigma-Model.} To appear in
 Annals of Math.


\bibitem{SchoenYau} Schoen, R., Yau, S.\ T. {\em Lectures on harmonic maps.}
Conference Proceedings and Lecture Notes in Geometry and Topology, II. International Press, Cambridge, MA, 1997.


\bibitem{ST} Segovia, C., Torrea, J. {\em Weighted inequalities for
commutators of fractional and singular integrals.} Conference on
Mathematical Analysis (El Escorial, 1989).  Publ.\ Mat.~35  (1991),
no.~1, 209--235.

\bibitem{Shatah} Shatah, J. {\em Weak solutions and development of singularities of the ${\rm SU}(2)$ $\sigma$-model.}
  Comm.\ Pure Appl.\ Math.~41  (1988),  no.~4, 459--469.

\bibitem{SStruwe} Shatah, J.,  Struwe, M. {\em Geometric wave equations.} Courant Lecture Notes in Mathematics, 2. New York University,
Courant Institute of Mathematical Sciences, New York; American Mathematical Society, Providence, RI, 1998.

\bibitem{ShaT1} Shatah, J.,  Tahvildar-Zadeh, A. {\em On the Cauchy problem
for equivariant wave maps.} Comm.\ Pure Appl.\ Math.~47  (1994),
no.~5, 719--754.

\bibitem{ShaT2} Shatah, J., Tahvildar-Zadeh, A. {\em Regularity of harmonic maps from the
Minkowski space into rotationally symmetric manifolds.}  Comm.\ Pure
Appl.\ Math.~45  (1992),  no.~8, 947--971.

\bibitem{Stein} Stein, E. {\em Harmonic Analysis.} Princeton, 1994.

\bibitem{SterbTat1} Sterbenz, J.,
Tataru, D. {\em Regularity of Wave-Maps in dimension $2+1$.}
Preprint 2009.

\bibitem{SterbTat2} Sterbenz, J.,
Tataru, D. {\em Energy dispersed large data wave maps in $2+1$
dimensions.} Preprint 2009.

\bibitem{Struwe} Struwe, M. {\em Variational methods, Applications
to Nonlinear PDEs and Hamiltonian Systems.} Second edition, Springer
Verlag, New York 1996.

\bibitem{Struwe1}   Struwe, M. {\em Equivariant wave maps in two space dimensions.}   Comm.\ Pure Appl.\ Math.~56  (2003),  no.~7, 815--823.


\bibitem{T9} Tao, T. {\em Endpoint bilinear restriction theorems for the cone, and some sharp null form estimates.}
  Math.\ Z.\  238  (2001),  no.~2, 215--268.

\bibitem{T8} Tao, T. {\em An inverse theorem for the bilinear $L^2$ Strichartz estimate for the wave equation.} Preprint 2009.

\bibitem{T7} Tao, T. {\em  Global regularity of wave maps VII. Control of
delocalised or dispersed solutions.} Preprint 2009.

\bibitem{T6} Tao, T. {\em  Global regularity of wave maps VI. Abstract theory of
minimal-energy blowup solutions.} Preprint 2009.

\bibitem{T5} Tao, T. {\em Global regularity of wave maps V. Large data local
wellposedness and perturbation theory in the energy class.} Preprint
2008.

\bibitem{T4} Tao, T. {\em Global regularity of wave maps IV. Absence of
stationary or self-similar solutions in the energy class.} Preprint
2008.

\bibitem{T3} Tao, T. {\em Global regularity of wave maps III. Large energy from $R^{1+2}$
to hyperbolic spaces.} Preprint 2008.

\bibitem{T2} Tao, T. {\em Global regularity of wave maps II. Small
energy in two dimensions.}   Comm.\ Math.\ Phys.~224  (2001), no.~2,
443--544.

\bibitem{T1} Tao, T. {\em  Global regularity of wave maps. I. Small critical Sobolev norm in high dimension.}
  Internat.\ Math.\ Res.\ Notices  2001,  no.~6, 299--328.

\bibitem{Tat} Tataru, D. {\em Rough solutions for the wave maps equation.}  Amer.\ J.\ Math.~127  (2005),  no.~2, 293--377.

\bibitem{Tat1} Tataru, D. {\em The wave maps equation.}  Bull.\ Amer.\ Math.\ Soc.~41  (2004),  no.~2, 185--204.

\bibitem{TatWolff} Tataru, D. {\em Null form estimates for second order hyperbolic operators with rough coefficients.}
Harmonic analysis at Mount Holyoke (South Hadley, MA, 2001),  383--409,
Contemp. Math., 320, Amer. Math. Soc., Providence, RI, 2003.

\bibitem{Tat2} Tataru, D. {\em On global existence and scattering for the wave maps equation.}
  Amer.\ J.\ Math.~123  (2001),  no.~1, 37--77.

\bibitem{Tat3} Tataru, D. {\em Local and global results for wave maps. I.}  Comm.\ Partial Differential Equations  23  (1998),  no.~9-10, 1781--1793.

\bibitem{W} Wolff, T. {\em A sharp bilinear cone restriction estimate.}  Ann.\ of Math.\ (2)  153  (2001),  no.~3, 661--698.



\end{thebibliography}
\end{document}